\documentclass[12pt]{article}
\usepackage{multind}
\makeindex{A}
\makeindex{B}

\usepackage{color}

\usepackage[all]{xy}
\usepackage{amscd} 
\usepackage{amsmath,amssymb}
\usepackage{amsthm}
\usepackage{graphicx}
\usepackage{enumerate}
\usepackage{epstopdf}
\usepackage{mathtools}
\usepackage{xr}

\numberwithin{equation}{section}
\numberwithin{figure}{section}

\newtheorem{theorem}{Theorem}[section]

\newtheorem{lemma}[theorem]{Lemma}
\newtheorem{proposition}[theorem]{Proposition}
\newtheorem{corollary}[theorem]{Corollary}
\newtheorem{fact}[theorem]{Fact}
\newtheorem{problem}[theorem]{Problem}
\theoremstyle{definition}
\newtheorem{definition}[theorem]{Definition}
\newtheorem{definitiontheorem}[theorem]{Definition-Theorem}
\newtheorem{definitionlemma}[theorem]{Definition-Lemma}
\newtheorem{example}[theorem]{Example}
\newtheorem{observation}[theorem]{Observation}
\theoremstyle{remark}
\newtheorem{remark}[theorem]{Remark}
\theoremstyle{conjecture}
\newtheorem{conjecture}[theorem]{Conjecture}

\newtheorem{convention}[theorem]{Convention}
\renewcommand{\setminus}{-}

\newcommand{\Abb}[4]{{\mathbb{A}}_{{#1},{#2},{#3}}^{{#4}}}
\newcommand{\Atbb}[4]
{{\widetilde{\mathbb{A}}}_{{#1},{#2},{#3}}^{{#4}}}

\newcommand{\Attbb}[4]
{\tilde{\tilde{{\mathbb{A}}}}_{{#1},{#2},{#3}}^{{#4}}}
\newcommand{\pr}[2]{\operatorname{pr}_{{#1}\to{#2}}}

\renewcommand{\Bbb}[3]{{\mathbb{B}}_{{#1},{#2}}^{{#3}}}

\newcommand{\Kirredrep}[2]{F^{{#1}}({#2})}

\newcommand{\Atcal}[4]{\widetilde{{\mathcal{A}}}_{{#1},{#2},{#3}}^{{#4}}}
\newcommand{\Attcal}[4]
{\tilde{\tilde{{\mathcal{A}}}}_{{#1},{#2},{#3}}^{{#4}}}
\newcommand{\Cbb}[3]{{\mathbb{C}}_{{#1},{#2}}^{{#3}}}
\newcommand{\Ctbb}[3]{\widetilde{{\mathbb{C}}}_{{#1},{#2}}^{{#3}}}

\newcommand{\Ctcal}[3]{\widetilde{{\mathcal{C}}}_{{#1},{#2}}^{{#3}}}

\newcommand{\Exterior}{\mathchoice{{\textstyle\bigwedge}}%
    {{\bigwedge}}%
    {{\textstyle\wedge}}%
    {{\scriptstyle\wedge}}} 
\newcommand{\Reducible}{{\mathcal{R}}ed}
\newcommand{\RedI}{\Reducible_{\operatorname{I}}}
\newcommand{\RedJ}{\Reducible_{\operatorname{I I}}}
\newcommand{\RInt}{{\mathcal{RI}}nt}
\newcommand{\Rij}[2]{R^{{#1},{#2}}}

\newcommand{\Rq}[3]{({#1}\,||\,{#2})_{#3}}

\newcommand{\Ttbb}[3]{\widetilde{{\mathbb{T}}}_{{#1},{#2}}^{{#3}}}
\newcommand{\Tttbb}[3]{\tilde{\tilde{{\mathbb{T}}}}_{{#1},{#2}}^{{#3}}}

\newcommand{\Ttcal}[3]{\widetilde{\mathcal{T}}_{{#1},{#2}}^{{#3}}}

\newcommand{\mub}{\mu^{\flat}}
\newcommand{\mus}{\mu^{\sharp}}
\newcommand{\Psising}{\Psi_{{\operatorname{sp}}}}
\newcommand{\redO}{\mbox{\color{red}O}}
\newcommand{\blueO}{\mbox{\color{blue}O}}
\newcommand{\greenO}{\mbox{\color{green}O}}
\newcommand{\magentaO}{\mbox{\color{magenta}O}}
\newcommand{\blackO}{\mbox{\color{black}O}}

\usepackage{tikz}
\tikzset{>=latex}

\begin{document}
\title
{Symmetry breaking for representations of rank one orthogonal 
groups II}
\author{Toshiyuki Kobayashi, 
\\
\normalsize{Graduate School of Mathematical Sciences,
 and Kavli IPMU,}
\\
\normalsize{the University of Tokyo, 3-8-1 Komaba, Tokyo, 153-8914 Japan,}
\\
\\
Birgit Speh
\\
\normalsize{Department of Mathematics, Cornell University,} 
\\
\normalsize{Ithaca, 14853-4201, NY, USA}
}
\maketitle%

\newpage
\tableofcontents
\newpage
\begin{abstract}

For a pair $(G,G')=(O(n+1,1), O(n,1))$
 of reductive groups, 
 we investigate intertwining operators ({\textit{symmetry breaking operators}})  between principal series representations
 $I_\delta(V,\lambda )$ of $G$, 
 and $J_\varepsilon(W,\nu)$ of the subgroup $G'$.  
The representations are parametrized by finite-dimensional representations $V$,  $W$ of $O(n)$ respectively of $O(n-1)$,
 characters $\delta$, $\varepsilon$ of $O(1)$, 
 and $\lambda, \ \nu \in {\mathbb C}$. 
Denote by $[V:W]$  the multiplicity of $W$ 
 occurring in the restriction $V|_{O(n-1)}$, 
 which is either 0 or 1.  
If $[V:W] \ne 0$
 then 
 we construct a holomorphic family of symmetry breaking operators
 and prove that 
$
   \dim_{\mathbb{C}} \operatorname{Hom}_{G'}(I_{\delta}(V, \lambda)|_{G'}, J_{\varepsilon}(W, \nu))
$
 is nonzero
 for all the parameters $\lambda$, $\nu$ and $\delta$, $\varepsilon$, 
 whereas if $[V:W] = 0$ there may exist {\textit{sporadic}} differential symmetry breaking operators. 

We propose a {\textit{classification scheme}} to find all matrix-valued symmetry breaking operators explicitly,
and carry out this program completely in the case
 $ (V,W)=(\bigwedge^i({\mathbb{C}}^{n}), \bigwedge^j({\mathbb{C}}^{n-1}))$.  
In conformal geometry, our results yield the complete classification
of conformal covariant operators from differential forms on a Riemannian manifold $X$
to those on a submanifold $Y$ in the model space $(X, Y) = (S^n, S^{n-1})$.

We use this information to determine the space of symmetry breaking operators 
 for any pair of irreducible representations
 of $G$ and the subgroup $G'$  with trivial infinitesimal character.  
Furthermore we prove the multiplicity conjecture by B.~Gross and D.~Prasad for  tempered principal series representations of $(SO(n+1,1), SO(n,1))$ 
 and also for 3 tempered representations $\Pi, \pi, \varpi$ of 
$SO(2m+2,1)$, $SO(2m+1,1)$ and $SO(2m,1)$
 with trivial infinitesimal character.
In connection to automorphic form theory,
 we apply our main results to find {\textit{periods}}
 of irreducible representations of the Lorentz group 
having nonzero $({\mathfrak{g}}, K)$-cohomologies.

This book is an extension of the recent work in the two research monographs:
Kobayashi--Speh [Memoirs Amer.~Math.~Soc., 2015]
for {\textit{spherical}} principal series representations 
and
Kobayashi--Kubo--Pevzner [Lecture Notes in Math., 2016]
for conformally covariant {\textit{differential}} symmetry breaking operators.

\end{abstract}
\newpage
\section{Introduction}
\label{sec:Intro}

\newcommand{\bZ}{{\mathbb Z}}
\newcommand{\bR}{{\mathbb R}}
\newcommand{\bC}{{\mathbb C}}
\newcommand{\bN}{{\mathbb N}}

A representation $\Pi$ of a group $G$
 defines a representation
 of a subgroup $G'$
 by restriction.  
In general irreducibility 
 is not preserved by the restriction.  
If $G$ is compact
 then the restriction $\Pi|_{G'}$ is isomorphic
 to a direct sum of irreducible finite-dimensional representations $\pi$ of $G'$ with multiplicities  $m(\Pi,\pi)$. 
These multiplicities are studied by using combinatorial techniques.  
We are interested in the case
 where $G$ and $G'$ are (noncompact) real reductive Lie groups.  
Then most irreducible representations $\Pi$ of $G$ are
 infinite-dimensional, 
 and generically the restriction $\Pi|_{G'}$
 is not a direct sum of irreducible representations
 \cite{KInvent98}.  
So we have to consider another notion of multiplicity.

For a continuous representation $\Pi$ of $G$
 on a complete, 
 locally convex topological vector space ${\mathcal{H}}$,
 the space ${\mathcal{H}}^\infty $ of $C^\infty$-vectors of ${\mathcal{H}}$
 is naturally endowed with a Fr{\'e}chet topology,
 and $(\Pi,{\mathcal{H}})$ induces a continuous representation $\Pi^{\infty}$ of $G$
 on ${\mathcal{H}}^\infty$.  
If $\Pi$ is an admissible representation
 of finite length on a Banach space ${\mathcal{H}}$, 
 then the Fr{\'e}chet representation 
 $(\Pi^{\infty}, {\mathcal{H}}^{\infty})$, 
 which we refer to as an 
\index{B}{admissiblesmoothrepresentation@admissible smooth representation}
{\it{admissible smooth representation}}, 
 depends only 
 on the underlying $({\mathfrak {g}}, K)$-module
 ${\mathcal{H}}_K$.  
In the context of asymptotic behaviour 
 of matrix coefficients, 
 these representations
 are also referred to as an admissible representations
 of moderate growth \cite[Chap.~11]{W}.  
We shall work with these representations
 and write simply $\Pi$ for $\Pi^{\infty}$.  
We denote by 
\index{A}{Irredrep@${\operatorname{Irr}}(G)$}
${\operatorname{Irr}}(G)$
 the set of equivalence classes
 of irreducible admissible smooth representations.  
We also sometimes call these representations
 \lq\lq{irreducible admissible representations}\rq\rq\
 for simplicity.  

Given another admissible smooth representation $\pi$
 of a reductive subgroup $G'$, 
 we consider the space of continuous $G'$-intertwining operators 
\index{B}{symmetrybreakingoperators@symmetry breaking operators}
({\it{symmetry breaking operators}})
\[ \operatorname{Hom}_{G'} ({\Pi}|_{G'}, {\pi}) .\] 
If $G=G'$ then these operators include the Knapp--Stein operators \cite{KS}
 and the differential intertwining operators 
 studied by B.~Kostant \cite{Kos}.
If $G\not = G'$ the dimension 
\index{A}{mPipi@$m(\Pi, \pi)$, multiplicity|textbf}
\[
   m(\Pi,\pi) :=\mbox{dim}_{\mathbb{C}} \operatorname{Hom}_{G'} ({\Pi}|_{G'}, {\pi}) 
\]
 yields important information of the restriction of $\Pi $ to $G'$
 and is called the
\index{B}{multiplicity@multiplicity|textbf}
 {\it{multiplicity}} of $\pi$
 occurring in the restriction $\Pi|_{G'}$.  
In general,
 $m(\Pi,\pi)$ may be infinite.  
The finiteness criterion in \cite{xKOfm} asserts that the
 multiplicity $m(\Pi,\pi)$ is finite
 for all $\Pi \in {\operatorname{Irr}}(G)$
 and for all $\pi \in {\operatorname{Irr}}(G')$
 if and only if a minimal parabolic subgroup $P'$
 of $G'$ has an open orbit 
 on the real flag variety $G/P$, 
 and that the multiplicity is uniformly bounded
 with respect to $\Pi$ and $\pi$
 if and only if a Borel subgroup of $G_{\mathbb{C}}'$
 has an open orbit
 on the complex flag variety of $G_{\mathbb{C}}$.  
 
\medskip
The latter condition depends only
 on the complexified pairs
 $({\mathfrak {g}}_{\mathbb{C}}, {\mathfrak {g}}_{\mathbb{C}}')$, 
 of which the classification
 was already known in 1970s
 by Kr{\"a}mer 
 \cite{Kr1} and Kostant.  
In particular,
 the multiplicity $m(\Pi,\pi)$ is uniformly bounded
 if the Lie algebras $({\mathfrak {g}}, {\mathfrak {g}}')$ of $(G,G')$ 
 are real forms
 of $({\mathfrak {sl}}(N+1,{\mathbb{C}}),
 {\mathfrak {gl}}(N,{\mathbb{C}}))$
 or $({\mathfrak {o}}(N+1,{\mathbb{C}}),
 {\mathfrak {o}}(N,{\mathbb{C}}))$.  
On the other hand, 
 the former condition depends on real forms
 $({\mathfrak {g}}, {\mathfrak {g}}')$, 
 and the classification
 of such symmetric pairs was recently accomplished
 in \cite{xKMt}.  
For instance,
 let $(G,G')=(O(n+1,1), O(n+1-k,1))$.  
Then the classification theory \cite{xKMt}
 and the finiteness criterion \cite{xKOfm}
 imply the following upper and lower estimates
 of the multiplicity
$
  m(\Pi,\pi)
$:
\begin{enumerate}
\item[(1)]
For $2 \le k \le n+1$, 
\begin{enumerate}
\item[]
$m(\Pi,\pi)<\infty$
 for every pair $(\Pi,\pi) \in {\operatorname{Irr}}(G)
 \times {\operatorname{Irr}}(G')$;
\item[]
$\underset{\Pi \in {\operatorname{Irr}}(G)}\sup
\,\,
 \underset{\pi \in {\operatorname{Irr}}(G')}\sup
 m(\Pi,\pi)=\infty.
$
\end{enumerate}
\item[(2)]
For $k=1$, 
there exists $C>0$ such that 
\begin{equation}
\label{eqn:mbdd}
m(\Pi,\pi) \le C
\text{ for all $\Pi \in {\operatorname{Irr}}(G)$
 and for all $\pi \in {\operatorname{Irr}}(G')$.}
\end{equation}
\end{enumerate}

B.~Sun  and C.-B.~Zhu \cite{SunZhu} showed
 that one can take $C$ to be one in \eqref{eqn:mbdd}, 
 namely,
 the  multiplicity
$
    m(\Pi,\pi) \in \{0,1\} 
$
 in this case.  
Thus one of the open problems is to determine
 when $m(\Pi, \pi) \ne 0$
 for irreducible representations $\Pi$ and $\pi$.

In the previous publication \cite{sbon} we initiated a thorough study of symmetry breaking operators
 between {\it{spherical}} principal series representations
of
\begin{equation}
\label{eqn:GG}
     (G,G')=(O(n+1,1),O(n,1)).  
\end{equation}
In particular,
 we determined the multiplicities $m(\Pi, \pi)$
 when both $\Pi$ and $\pi$ are irreducible composition factors
 of the spherical principal series representations.

In this article we will determine
 the multiplicities $m(\Pi,\pi)$ for all irreducible representations
 $\Pi$ and $\pi$
 with trivial infinitesimal character $\rho$ of $G=O(n+1,1)$ and $G'=O(n,1)$, 
 respectively, 
 and also for irreducible principal series representations.

More than just determining the dimension $m(\Pi,\pi)$
 of the space of symmetry breaking operators, 
 we investigate these operators of their own 
 for {\it{general}} principal series representations of $G$ and the subgroup $G'$, 
{\it{i.e.}}, 
 for representations induced from irreducible finite-dimensional representations of a parabolic subgroup.  
We construct a holomorphic family
 of symmetry breaking operators, 
 and present a 
\index{B}{classificationscheme@classication scheme, symmetry breaking operators}
 classification scheme of {\it{all}} symmetry breaking operators $T$
 in Theorem \ref{thm:VWSBO}
 through an analysis of their distribution kernels $K_T$.  
In particular, 
 we prove that any symmetry breaking operators
 in this case is either a 
\index{B}{sporadicsymmetrybreakingoperator@sporadic symmetry breaking operator}
 sporadic differential symmetry breaking operator
 (cf. \cite{KKP})
 or the analytic continuation
 of integral symmetry breaking operators
 and their renormalization
 in Theorem \ref{thm:VWSBO}.

The proof for the explicit formula
 of the multiplicity $m(\Pi, \pi)$ is built on 
 the 
\index{B}{functionalequation@functional equation}
 functional equations 
 (Theorems \ref{thm:TAA} and \ref{thm:ATA})
 satisfied by the regular symmetry breaking operators.

\vskip 2pc
A principal series representation
 $I_\delta(V,\lambda)$ of $G=O(n+1,1)$
 is an (unnormalized) induced representation from an irreducible
 finite-dimensional representation
 $V\otimes \delta \otimes {\mathbb{C}}_\lambda$
 of a minimal parabolic subgroup $P=MAN_+$. 
In our setting, 
 $M \simeq O(n) \times  {\mathbb{Z}}/2{\mathbb{Z}}$ and $A \simeq \bR_+$. 
We assume that $V$ is a representation of $O(n+1)$, 
 $\delta \in \{\pm\}$,
 and $\lambda \in \bC$.
In what follows,
 we identify the representation space
 of $I_\delta(V, \lambda)$ 
 with the space
 of $C^{\infty}$-sections
 of the $G$-equivariant  bundle 
 $G \times_P {\mathcal V}_{\delta,\lambda} \to G/P$, 
so that 
 $I_\delta(V, \lambda)^{\infty} = I_\delta (V,\lambda)$
 is the Fr{\'e}chet globalization
 having moderate growth
 in the sense 
 of Casselman--Wallach \cite{W}.  
The parametrization is chosen
 so that the representation $ I_\delta (V,\frac n2)$ is a unitary tempered representation. 
The representations $I_\delta(V,\lambda)$ are either irreducible
 or of composition series of length 2, 
 see Corollary \ref{cor:length2} in Appendix I.

The group $P'= G' \cap P = M'AN_+'$ is a minimal parabolic subgroup
 of $G'=O(n,1)$.  
For an irreducible representation $(\tau, W)$ of $O(n-1)$, a character 
 $\varepsilon\in \{\pm\}$ of $O(1)$, 
 and $\nu \in{\mathbb{C}}$ we define the principal series representation $J_\varepsilon(W,\nu)$ of $G'$.  

\medskip
We set 
\[
   [V:W]:= \dim_{\mathbb{C}}\operatorname{Hom}_{O(n-1)} (W,V|_{O(n-1)})
         = \dim_{\mathbb{C}}\operatorname{Hom}_{O(n-1)} (V|_{O(n-1)},W).  
\]
For principal series representations $I_\delta(V,\lambda )$ of $G$
 and $J_\varepsilon(W,\nu)$ of the subgroup $G'$, 
 we consider the cases $[V:W] \ne 0 $ and $[V:W] = 0$ separately. In the first case we obtain a lower bound 
for the multiplicity.

In what follows, 
 it is convenient to introduce the set of 
 \lq\lq{special parameters}\rq\rq:
\index{A}{1psi@$\Psising$,
          special parameter in ${\mathbb{C}}^2 \times \{\pm\}^2$|textbf}
\begin{alignat}{2}
  \Psising
  :=
  \left\{(\lambda,\nu, \delta, \varepsilon) \in {\mathbb{C}}^2 \times \{\pm\}^2
   :\,\, \right.
&   \nu-\lambda \in 2 {\mathbb{N}}
&&\text{when $\delta \varepsilon =+$}
\notag
\\
\label{eqn:singset}
   \text{ or }\quad
& \nu-\lambda \in 2 {\mathbb{N}}+1
 \qquad
&&\text{when $\delta \varepsilon =-$}
  \left. \right\}.  
\end{alignat}

\begin{theorem}
[see Theorem \ref{thm:VWSBO} (2) and Theorem \ref{thm:1532113}]
\label{thm:intro160150}
Suppose $(\sigma,V)\in \widehat{O(n)}$ and $(\tau,W)\in \widehat{O(n-1)}$.  
Assume $[V:W] \ne 0$.  
\begin{enumerate}
\item[{\rm{(1)}}] 
{\rm{(existence of symmetry breaking operators)}}\enspace
We have
\[
  \dim_{\mathbb{C}} \operatorname{Hom}_{G'}(I_{\delta}(V, \lambda)|_{G'}, J_{\varepsilon}(W, \nu))
  \ge 1
  \quad
  \text{for all 
  $\delta, \varepsilon \in \{\pm\}$, 
  and $\lambda, \nu \in {\mathbb{C}}$}.  
\]
\item[{\rm{(2)}}]
\index{B}{genericmultiplicityonetheorem@generic multiplicity-one theorem}
{\rm{(generic multiplicity-one)}}\enspace
\[
\dim_{\mathbb{C}} \operatorname{Hom}_{G'}(I_{\delta}(V, \lambda)|_{G'}, J_{\varepsilon}(W, \nu))
  =1
\]
for any $(\lambda,\nu, \delta, \varepsilon) \in ({\mathbb{C}}^2 \times \{\pm\}^2) \setminus \Psising$.  
\item[{\rm{(3)}}]
Let $\ell(\sigma)$ be the \lq\lq{norm}\rq\rq\
 of $\sigma$ defined by using its highest weight 
 (see \eqref{eqn:ONlength}).  
Then we have 
\[
  \dim_{\mathbb{C}} \operatorname{Hom}_{G'}
  (I_{\delta}(V, \lambda)|_{G'}, J_{\varepsilon}(W, \nu))>1
\]
 for any $(\lambda,\nu, \delta, \varepsilon) \in \Psising$
 such that $\nu \in {\mathbb{Z}}$ with $\nu \le -\ell(\sigma)$.  
\end{enumerate}
\end{theorem}

We prove Theorem  \ref{thm:intro160150}
 by constructing (generically) regular symmetry breaking operators
 $\Atbb \lambda \nu {\delta\varepsilon} {V,W}$:
 they are nonlocal operators 
 ({\it{e.g.}}, integral operators)
 for generic parameters, 
 whereas for some parameters they are local operators
 ({\it{i.e.}}, differential operators).  
See Theorem \ref{thm:152389}
 for the construction of the normalized operator $\Atbb \lambda \nu \pm {V,W}$;
Theorem \ref{thm:regexist}
 for \lq\lq{regularity}\rq\rq\
 (\cite[Def.~3.3]{sbon})
 of $\Atbb \lambda \nu \pm {V,W}$
 under a certain generic condition;
Theorem \ref{thm:170340}
 for a renormalization of $\Atbb \lambda \nu \pm {V,W}$
 when it vanishes;
Fact \ref{fact:153316}
 for the residue formula of $\Atbb \lambda \nu \pm {V,W}$
 when it reduces to a differential operator.

\vskip 1pc
In the case $[V:W]=0$, 
 symmetry breaking operators are \lq\lq{rare}\rq\rq\ 
 but 
there may exist {\it{sporadic}} symmetry breaking operators: 
\begin{theorem}
\label{thm:VW0SBO}
Assume $[V:W] = 0$.  
\begin{enumerate}
\item[{\rm{(1)}}]
{\rm{(vanishing for generic parameters, 
 Corollary \ref{cor:VWvanish})}}
If $(\lambda, \nu, \delta, \varepsilon) \not \in \Psising$,
then 
\[
   {\operatorname{Hom}}_{G'}(I_{\delta}(V,\lambda)|_{G'}, J_{\varepsilon}(W,\nu))=\{0\}.  
\] 
\item[{\rm{(2)}}]
{\rm{(localness theorem, Theorem \ref{thm:152347})}}
\index{B}{localnesstheorem@localness theorem}
Any nontrivial symmetry breaking operator
\[
   C^{\infty}(G/P, {\mathcal{V}}_{\lambda,\delta})
   \to 
   C^{\infty}(G'/P', {\mathcal{W}}_{\nu,\varepsilon})
\]
is a differential operator.  
\end{enumerate}
\end{theorem}

Combining Theorem \ref{thm:intro160150} (2) and Theorem \ref{thm:VW0SBO} (1)
 together with the existence condition
 of differential symmetry breaking operators 
 (see Theorem \ref{thm:vanishDiff}), 
 we determine the following multiplicity formul{\ae}
 for {\it{generic parameters}}:
\begin{theorem}
Suppose that $(\lambda,\nu, \delta,\varepsilon) \not \in \Psising$.   
Then there are no differential symmetry breaking operators and
\[
     {\operatorname{dim}}_{\mathbb{C}} {\operatorname{Hom}}_{G'}(I_\delta(V,\lambda)|_{G'},J_\varepsilon (W,\nu))
=
\begin{cases}
 1 \qquad&\text{if $[V:W]\ne0$, }
\\
 0 \qquad&\text{if $[V:W]=0$.  }
\end{cases}
\]
\end{theorem}

\bigskip
It deserves to be mentioned
 that the parameter set $({\mathbb{C}}^2 \times \{\pm\}^2)\setminus \Psising$ contains parameters $(\lambda,\nu)$
 for which the $G$-module $I_{\delta}(V,\lambda)$
 or the $G'$-module $J_{\varepsilon}(W,\nu)$
 is {\it{not}} irreducible.

\vskip 2pc
In the major part of this monograph,
 we focus our attention on the special case
\[
  (V,W)=(\Exterior^i({\mathbb{C}}^{n}), \Exterior^j({\mathbb{C}}^{n-1})).  
\]
The principal series representations
 of $G$ and the subgroup $G'$
 are written as $I_{\delta}(i,\lambda)$ for $I_{\delta}(\Exterior^i({\mathbb{C}}^{n}),\lambda)$
 and $J_{\varepsilon}(j,\nu)$ for $J_{\varepsilon}(\Exterior^j({\mathbb{C}}^{n-1}),\nu)$, 
 respectively.  
The representations
 $I_\delta(i,\lambda)$ of $G$
 and $J_{\varepsilon}(j,\nu)$ of $G'$ are of  interest in geometry as well as in automorphic forms and in the cohomology of arithmetic groups. 
In geometry, 
 given an arbitrary Riemannian manifold $X$, 
 one forms a natural family of representations
 of the conformal group $G$
 on the space ${\mathcal{E}}^i(X)$ of differential forms,
 to be denoted by ${\mathcal{E}}^i(X)_{\lambda',\delta'}$
 for $0 \le i \le \dim X$, $\lambda' \in {\mathbb{C}}$, 
 and $\delta' \in \{\pm\}$.  
Then the representations $I_\delta(i,\lambda)$
 are identified with such conformal representations
 in the case where $(G, X)=(O(n+1,1), S^n)$, 
 see {\it{e.g.}}, \cite[Chap.~2, Sect.~2]{KKP}
 for precise statement.  
In representation theory,
 all irreducible, unitarizable representations
 with nonzero $({\mathfrak{g}}, K)$-cohomology arise
 as subquotients of $I_{\delta}(i,\lambda)$ with $\lambda=i$
 for some $0 \le i \le n$
 and $\delta=(-1)^i$, 
 see Theorem \ref{thm:LNM20} (9).

Our main results of this article include 
 a complete solution
 to the general problem
 of constructing and classifying the elements
 of ${\operatorname{Hom}}_{G'}(\Pi|_{G'}, \pi)$
 (see \cite[Prob.~7.3 (3) and (4)]{xkvogan})
 in the following special setting:
\begin{alignat*}{2}
(G,G')&=(O(n+1,1),O(n,1)) \qquad\text{with $n \ge 3$, }
\\
(\Pi,\pi)&=(I_{\delta}(i,\lambda),J_{\varepsilon}(j,\nu)), 
\end{alignat*}
where $0 \le i \le n$, $0 \le j \le n-1$, 
 $\delta, \varepsilon \in \{\pm\}$, 
 and $\lambda, \nu \in {\mathbb{C}}$.  
Thus our main results include a complete solution 
 to the following question in conformal geometry:
\begin{problem}
\label{prob:conf}
\begin{enumerate}
\item[{\rm{(1)}}]
Find a necessary and sufficient condition on 6-tuples
 $(i,j,\lambda,\nu,\delta,\varepsilon)$ for the existence
 of conformally covariant, 
 symmetry breaking operators
\[
  A \colon {\mathcal{E}}^i(X)_{\lambda,\delta} 
           \to {\mathcal{E}}^i(X)_{\nu,\varepsilon}
\]
in the model space $(X,Y)=(S^n, S^{n-1})$.  
\item[{\rm{(2)}}]
Construct those operators explicitly
 in the (flat) coordinates.  
\item[{\rm{(3)}}]
Classify all such symmetry breaking operators.  
\end{enumerate}
\end{problem}
Partial results were known earlier:
when the operator $A$ is given by a {\it{differential}} operator,
 Juhl \cite{Juhl} solved Problem \ref{prob:conf}
 in the case $(i,j)=(0,0)$,
 see also \cite{KOSS}, 
 which has been recently extended
 in Kobayashi--Kubo--Pevzner \cite{KKP}
 for the general $(i,j)$.  
Problem \ref{prob:conf} was solved for all 
 (possibly, {\it{nonlocal}}) operators
 in our previous paper \cite{sbon}
 in the case $(i,j)=(0,0)$.  
The complete classification of (continuous) symmetry breaking operators 
 for the general $(i,j)$ is given in Theorem \ref{thm:1.1} (multiplicity)
 and Theorem \ref{thm:SBObasis} (construction of explicit generators),  
 and we have thus settled Problem \ref{prob:conf}
 in this monograph.  
For this introduction,
 we explain only the
 \lq\lq{multiplicity}\rq\rq\
 (Theorem \ref{thm:1.1}).  
For this,
 using the same notation
 as in \cite[Chap.~1]{sbon}, 
 we define the following two subsets on ${\mathbb{Z}}^2$:
\begin{align*}
L_{\operatorname{even}}:=&\left \{ (-i,-j):
0 \le j\leq i \mbox{ and } i\equiv j \mod 2 \right \},
\\
L_{\operatorname{odd}}:=&\left \{ (-i,-j)
: 0 \le j\leq i \mbox{ and } i \equiv j +1 \mod 2 \right \}.
\end{align*}

\begin{theorem}
[multiplicity, Theorem \ref{thm:1.1}]
Suppose $\Pi=I_{\delta}(i,\lambda)$
 and $\pi=J_{\varepsilon}(j,\nu)$
 for $0 \le i\le n$, $0 \le j \le n-1$,
 $\delta, \varepsilon \in \{\pm \}$, 
 and $\lambda,\nu \in {\mathbb{C}}$.  
Then 
we have the following.  
\begin{enumerate}
\item[{\rm{(1)}}]
\begin{alignat*}{2}
m(\Pi,\pi) \in &\{ 1,2 \} \qquad
&&\text{if $j=i-1$ or $i$}, 
\\
m(\Pi,\pi) \in &\{ 0,1 \} \qquad
&&\text{if $j=i-2$ or $i+1$}, 
\\
m(\Pi,\pi) =& 0 \qquad
&&\text{otherwise}.  
\end{alignat*}
\item[{\rm{(2)}}]
Suppose $j=i-1$ or $i$.  
Then $m(\Pi,\pi)=1$ generically,
 and $=2$ when the parameter belongs to the following exceptional 
 countable set.  
Without loss of generality,
 we take $\delta$ to be $+$.  
\begin{enumerate}
\item[{\rm{(a)}}]
Case $1 \le i \le n-1$.
\begin{alignat*}{2}
m(I_+(i,\lambda), J_+(i,\nu))
=
&2
\qquad
&&\text{if }
(\lambda, \nu) \in L_{\operatorname{even}}\setminus \{\nu=0\}
                   \cup \{(i,i)\}.   
\\
m(I_+(i,\lambda), J_-(i,\nu))
=
&2
&&\text{if }
(\lambda, \nu) \in L_{\operatorname{odd}}\setminus \{\nu=0\}.   
\\
m(I_+(i,\lambda), J_+(i-1,\nu))
=
&2
&&
\text{if }
(\lambda, \nu) \in L_{\operatorname{even}}\setminus \{\nu=0\} 
                   \cup \{(n-i,n-i)\}.   
\\
m(I_+(i,\lambda), J_-(i-1,\nu))
=
&2
&&
\text{if }
(\lambda, \nu) \in L_{\operatorname{odd}}\setminus \{\nu=0\}.   
\end{alignat*}
\item[{\rm{(b)}}]
Case $i=0$.  
\begin{align*}
m(I_+(0,\lambda), J_+(0,\nu))
=&
2
\qquad
\text{if }
(\lambda, \nu) \in L_{\operatorname{even}}.   
\\
m(I_+(0,\lambda), J_-(0,\nu))
=&
2
\qquad
\text{if }
(\lambda, \nu) \in L_{\operatorname{odd}}.   
\end{align*}
\item[{\rm{(c)}}]
Case $i=n$.  
\begin{align*}
m(I_+(n,\lambda), J_+(n-1,\nu))
=&
2
\qquad
\text{if }
(\lambda, \nu) \in L_{\operatorname{even}}.   
\\
m(I_+(n,\lambda), J_-(n-1,\nu))
=&
2
\qquad
\text{if }
(\lambda, \nu) \in L_{\operatorname{odd}}.   
\end{align*}
\end{enumerate}
\item[{\rm{(3)}}]
Suppose $j=i-2$ or $i+1$.  
Then $m(\Pi,\pi)=1$ 
 if one of the following conditions {\rm{(d)--(g)}}
 is satisfied, 
 and $m(\Pi,\pi)=0$ otherwise.  
\begin{enumerate}
\item[{\rm{(d)}}]
Case $j=i-2$, $2 \le i \le n-1$, $(\lambda, \nu)=(n-i,n-i+1)$, 
$\delta \varepsilon =-1$.  
\item[{\rm{(e)}}]
Case $(i,j)=(n,n-2)$, $-\lambda \in {\mathbb{N}}$, $\nu=1$, 
$\delta \varepsilon =(-1)^{\lambda+1}$.  
\item[{\rm{(f)}}]
Case $j=i+1$, $1 \le i \le n-2$, $(\lambda, \nu)=(i,i+1)$, 
$\delta \varepsilon =-1$.  
\item[{\rm{(g)}}]
Case $(i,j)=(0,1)$, $-\lambda \in {\mathbb{N}}$, $\nu=1$, 
$\delta \varepsilon =(-1)^{\lambda+1}$.  
\end{enumerate}
\end{enumerate}
\end{theorem}

\medskip
More than just an abstract formula of multiplicities, 
 we also obtain explicit generators
 of ${\operatorname{Hom}}_{G'}(I_{\delta}(i,\lambda)|_{G'}, J_{\varepsilon}(j,\nu))$
 for $j \in \{i-1,i\}$
 in Theorem \ref{thm:SBObasis}.  
The generators for $j \in \{i-2,i+1\}$ are always differential operators
 ({\it{localness theorem}}, 
 see Theorem \ref{thm:VW0SBO} (2)), 
 and they were constructed and classified in \cite{KKP}
 (see Fact \ref{fact:3.9}).

The principal series representations
 $I_{\delta}(i,\lambda)$ and $J_{\varepsilon}(j,\nu)$
 in the above theorem
 are not necessarily irreducible.  
For the study of symmetry breaking of the irreducible subquotients,
 we utilize the concrete generators of 
${\operatorname{Hom}}_{G'}(I_{\delta}(i,\lambda)|_{G'}, J_{\varepsilon}(j,\nu))$ and determine explicit formul{\ae}
 about 
\begin{enumerate}
\item[$\bullet$]
the $(K,K')$-spectrum
 of the normalized regular symmetry breaking operators
 $\Atbb \lambda \nu \pm {i,j}$
 on basic \lq\lq{$(K,K')$-types}\rq\rq\
 (Theorem \ref{thm:153315});
\item[$\bullet$]
the functional equations among symmetry breaking operators
 $\Atbb \lambda \nu \pm {i,j}$
(Theorems \ref{thm:TAA} and \ref{thm:ATA}).  
\end{enumerate}
Here, 
  the 
\index{B}{KspectrumKprime@$(K,K')$-spectrum}
$(K,K')$-spectrum is defined in Definition \ref{def:KKspec}.  
It resembles eigenvalues of symmetry breaking
 operators, 
 and serves as a clue to find the functional equations.

\bigskip
We now highlight symmetry breaking
 of irreducible representations
 that have the same infinitesimal character $\rho$
 with the trivial one-dimensional representation ${\bf{1}}$.  
Denote by 
\index{A}{IrrGrho@${\operatorname{Irr}}(G)_{\rho}$,
 set of irreducible admissible smooth representation of $G$
 with trivial infinitesimal character $\rho$\quad}
${\operatorname{Irr}}(G)_{\rho}$
 the (finite) set
 of equivalence classes 
 of irreducible admissible representations
 of $G$
 with trivial infinitesimal character $\rho\equiv \rho_G$.  
The principal series representations $I_\delta(i,i)$ of $G=O(n+1,1)$
 are reducible, 
 and any element in ${\operatorname{Irr}}(G)_{\rho}$
 is a subquotient
 of the representations $I_\delta(i,i)$
 for some $0 \le i \le n$
 and $\delta \in \{\pm\}$.  
To be more precise,
 we have the following.  
\begin{theorem}
[see Theorem \ref{thm:LNM20}]
Let $G=O(n+1,1)$ $(n \ge 1)$.  
\begin{enumerate}
\item[{\rm{(1)}}]
For $0 \le \ell \le n$ and $\delta \in \{ \pm \}$, 
 there are exact sequences of $G$-modules:
\begin{align*}
&0 \to \Pi_{\ell,\delta} \to I_{\delta}(\ell,\ell) \to \Pi_{\ell+1,-\delta}
\to 0, 
\\
&0 \to \Pi_{\ell+1,-\delta} \to I_{\delta}(\ell,n-\ell) \to \Pi_{\ell,\delta} \to 0.  
\end{align*}
These exact sequences split
 if and only if $n=2\ell$.

\item[{\rm{(2)}}]
Irreducible admissible representations of $G$ 
 with trivial infinitesimal character 
 can be classified as
\[
  {\operatorname{Irr}}(G)_{\rho}
  =\{
    \Pi_{\ell, \delta}
    :
    0 \le \ell \le n+1, \, \delta = \pm
\}.  
\]

\item[{\rm{(3)}}]
Every $\Pi_{\ell,\delta}$ $(0 \le \ell \le n+1, \delta = \pm)$
 is unitarizable.  
\end{enumerate}
\end{theorem}

\medskip
There are four one-dimensional representations of $G$, 
 and they are given by 
\begin{equation*}
\{
  \Pi_{0,+} \simeq {\bf{1}}, \quad \Pi_{0,-} \simeq \chi_{+-}, 
\quad
   \Pi_{n+1,+} \simeq \chi_{-+}, \quad 
   \Pi_{n+1,-} \simeq \chi_{--} (=\det)
\}.  
\end{equation*}
(See \eqref{eqn:chiab} for the definition of
\index{A}{1chipmpm@$\chi_{\pm\pm}$, one-dimensional representation of $O(n+1,1)$}
 $\chi_{\pm\pm}$.)
The other representations $\Pi_{\ell,\delta}$
 ($1 \le \ell \le n$, $\delta=\pm$)
 are infinite-dimensional representations.

For the subgroup $G'=O(n,1)$,
 we use the letters $\pi_{j,\varepsilon}$ to 
 denote the irreducible representations
 in ${\operatorname{Irr}}(G')_{\rho}$, 
 similar to $\Pi_{i,\delta}$
 in ${\operatorname{Irr}}(G)_{\rho}$.

With these notations,
 we determine
\[
  m(\Pi_{i,\delta},\pi_{j,\varepsilon}) = \dim_{\mathbb{C}} {\operatorname{Hom}}_{G'}
  (\Pi_{i,\delta}|_{G'},\pi_{j,\varepsilon})
\]
for all $\Pi_{i,\delta} \in {\operatorname{Irr}}(G)_{\rho}$
 and $\pi_{j,\varepsilon} \in {\operatorname{Irr}}(G')_{\rho}$
 as follows.  
\begin{theorem}
[vanishing, 
 see Theorem {\ref{thm:SBOvanish}}]
\label{thm:introSBOvanish}
\index{B}{vanishingtheorem@vanishing theorem}
Suppose $0 \le i \le n+1$, $0 \le j \le n$, 
 $\delta, \varepsilon \in \{\pm\}$.  
\begin{enumerate}
\item[{\rm{(1)}}]
If $j \not = i, i-1$ then
$
  {\operatorname{Hom}}_{G'}(\Pi_{i,\delta}|_{G'}, \pi_{j,\varepsilon})=\{0\}. 
$

\item[{\rm{(2)}}] 
If $\delta \varepsilon =-$, 
 then 
$ 
{\operatorname{Hom}}_{G'}(\Pi_{i,\delta}|_{G'}, \pi_{j,\varepsilon}) =\{0\}.  
$ 
\end{enumerate}
\end{theorem}

\medskip
\noindent

\begin{theorem} 
[multiplicity-one, see Theorem {\ref{thm:SBOone}}]
\label{thm:introSBOone}
\index{B}{multiplicityonetheorem@multiplicity-one theorem}
Suppose $0 \le i \le n+1$, $0 \le j \le n$
 and $\delta, \varepsilon \in \{\pm\}$.  
If $j=i-1$ or $i$
 and if $\delta \varepsilon =+$, 
 then 
\[
   \dim_{\mathbb{C}}
   {\operatorname{Hom}}_{G'}
   (\Pi_{i,\delta}|_{G'}, \pi_{j,\varepsilon}) =1.  
\]
\end{theorem}

\medskip

We can represent these results graphically as follows.  
We suppress the subscript, 
 and write $\Pi_i$ for $\Pi_{i,+}$, 
 and $\pi_j$ for $\pi_{j,+}$. 
 The first row are representations of $G$,
 the second row are representations of $G'$. 
The existence of nonzero symmetry breaking operators is represented by  arrows.

\begin{theorem}
[see Theorem \ref{thm:SBOfg}]
\label{thm:introSBOfg}
Symmetry breaking for  irreducible representations with infinitesimal character $\rho$  is represented graphically in the following form.  

{Symmetry breaking for $O(2m+1,1)\downarrow O(2m,1)$ }

\begin{center}
\begin{tabular}{c@{~}c@{~}c@{~}c@{~}c@{~}c@{~}c@{~}c@{~}c}
$\Pi_0$& &$\Pi_1$& &\dots & & $\Pi_{m-1} $& & $\Pi_{m}$ 
\\
$\downarrow$ &$\swarrow$& $\downarrow $& $\swarrow$ & &$\swarrow$ & $ \downarrow $&  $\swarrow $  &  $\downarrow$ 
\\
$\pi_0$& &$\pi_1$& &\dots & & $\pi_{m-1}$ & & $\pi_{m}$ 
\end{tabular}
\end{center}
\label{fig:introHasse1}

\bigskip

{Symmetry breaking for $O(2m+2,1) \downarrow O(2m+1,1)$ }
\begin{center}
\begin{tabular}{@{}c@{~}c@{~}c@{~}c@{~}c@{~}c@{~}c@{~}c@{~}c@{~}c@{~}c@{}}
$\Pi_0$& &$\Pi_1$& &\dots & & $\Pi_{m-1} $& & $\Pi_{m}$ & & $\Pi_{m+1}$
\\
$\downarrow$ & $\swarrow$ & $\downarrow $ & $\swarrow$ & & $\swarrow$ & $ \downarrow $& $\swarrow $ & $\downarrow$ & $\swarrow$ & \\
$\pi_0$& &$\pi_1$& &\dots & & $\pi_{m-1}$ & & $\pi_{m}$ 
\end{tabular}
\end{center}
\label{fig:introHasse2}
\end{theorem}

We believe that we are seeing in Theorem \ref{thm:SBOfg}
 only the \lq\lq{tip of the iceberg}\rq\rq, 
 and we present a conjecture that a similar statement holds
 in more generality, 
 see Conjecture \ref{conj:GPver1}.  
Suppose that $F$ and $F'$ are irreducible finite-dimensional 
 representations
 of $G$ and the subgroup $G'$, 
 respectively,  and 
 that
\[ \mbox{Hom}_{G'}(F|_{G'},F') \not = \{0\}.\] 
In Chapters \ref{sec:conjecture} and \ref{sec:aq}
 we describe  sequences of irreducible representations $\{\Pi_i\equiv \Pi_i(F) \}$
 and $\{\pi_j \equiv \pi_j(F') \}$ of $G$ and $G'$
 with the same infinitesimal characters with $F$ and $F'$,
 respectively.  
We refer to these  sequences as {\it{standard sequences}}
 that starting with $\Pi_0(F)=F$
 and $\pi_0(F')=F'$, 
 see Definition \ref{def:Hasse}. 
They generalize the standard sequence 
 with trivial infinitesimal character
 which we used in the formulation of Theorem \ref{thm:introSBOfg}.  
They are an analogue of a diagrammatic description
 of irreducible representations 
 with regular integral infinitesimal characters for the connected group $G_0=SO_0(n+1,1)$
 given in Collingwood
 \cite[p.~144, Fig.~6.3]{C}.  
In this generality,
 we conjecture
 that the results of symmetry breaking can be represented
 graphically exactly as in Theorem \ref{thm:introSBOfg} 
 for the representations
 with trivial infinitesimal character $\rho$.   
Again in the first row are representations of $G$, 
 and in the second row are representations of $G'$. 
Conjecture \ref{conj:GPver1} asserts 
 that symmetry breaking operators are represented by  arrows.

 \medskip
{Symmetry breaking for $O(2m+1,1)\downarrow O(2m,1)$ }

\begin{center}
\begin{tabular}{c@{~}c@{~}c@{~}c@{~}c@{~}c@{~}c@{~}c@{~}c@{~}c}
$\Pi_0(F)$ & &$\Pi_1(F)$ & &\dots & & $\Pi_{m-1}(F)$ & &$\Pi_{m}(F)$ 
\\
$\downarrow$ & $\swarrow$ & $\downarrow$ & $\swarrow$ &  & $\swarrow$ & $ \downarrow $ & $\swarrow $ &  $\downarrow$ 
\\
$\pi_0(F')$ & & $\pi_1(F')$ & &\dots & & $\pi_{m-1}(F')$ & & $\pi_{m}(F')$ 
\end{tabular}
\end{center}
\label{fig:introHa1}

 \medskip
{Symmetry breaking for $O(2m+2,1)\downarrow O(2m+1,1)$ }

\begin{center}
\begin{tabular}{@{}c@{~}c@{~}c@{~}c@{~}c@{~}c@{~}c@{~}c@{~}c@{~}c@{~}c@{~}c@{}}
$\Pi_0(F)$& &$\Pi_1(F)$ & & \dots & & $\Pi_{m-1}(F)$& & $\Pi_{m}(F)$ & & $\Pi_{m+1}(F)$
\\
$\downarrow$ &$\swarrow$& $\downarrow $& $\swarrow$ &  & $\swarrow$ & $ \downarrow $& $\swarrow $ & $\downarrow$ & $\swarrow$ & \\
$\pi_0(F')$& &$\pi_1(F')$& &\dots  & & $\pi_{m-1}(F')$ & & $\pi_{m}(F')$ 
\end{tabular}
\end{center}
\label{fig:introHa2}

\medskip \noindent
We present some supporting evidence for this conjecture
 in Chapter \ref{sec:conjecture}.

\bigskip

Applications of our formul{\ae} include
 some results about periods of representations.  
Suppose that $H$ is a subgroup of $G$. 
Following the  terminology used in automorphic forms and the relative trace formula,
 we say that a smooth representation $U $ of $G$ is
\index{B}{distinguishedH@distinguished, $H$-|textbf}
 {\it{$H$-distinguished}} if there is a nontrivial linear $H$-invariant linear functional 
\[F^{H}:U \rightarrow \bC. \]
If the $G$-module $U$ is $H$-distinguished,
 we say that $(F^{H}, H)$ is
 a 
\index{B}{period@period}
 {\it{period}} (or an $H$-period) of $U$.

Let $(G,H)=(O(n+1,1),O(m+1,1))$ with $m \le n$.  
For $0 \le i \le n+1$ and $0 \le j \le m+1$, 
 we denote by $\Pi_{i}$ and $\pi_j$
 the irreducible representations $\Pi_{i,+}$ of $G$
 and analogous ones of $H$ with trivial infinitesimal character $\rho $.

\begin{theorem} 
[see Theorems \ref{thm:period1} and \ref{thm:171517}]
\
\begin{enumerate}
\item[{\rm{(1)}}]
The irreducible representation $\Pi_i $ is $H$-distinguished if $i \le n-m$.  
\item[{\rm{(2)}}]
The outer tensor product representation
\[ \Pi_i \boxtimes  \pi_j\]
has a nontrivial $H$-period if $0 \le i-j \le n-m$.  
\end{enumerate}
\end{theorem}

\medskip
The period is given by the composition 
 of the normalized regular symmetry breaking operators
 (see Chapter \ref{sec:section7})
 with respect to the chain
 of subgroups:  
\[
  G =O(n+1,1) \supset O(n,1) \supset O(n-1,1) \supset \cdots \supset O(m+1,1)
  =H.  
\]
Using the above chain of subgroups we also define a vector $v$ in the 
\index{B}{minimalKtype@minimal $K$-type}
minimal $K$-type of $\Pi_i$.  
We prove 
\begin{theorem}
[see Theorem \ref{thm:period2}]
\label{thm:intro-period2}
Suppose that $G= O(n+1,1)$ and $\Pi_i$ $(0 \le i \le n)$ is the irreducible representation
 with trivial infinitesimal character $\rho$
 defined as above.
Then the value of the $O(n+1-i,1)$-period on $v \in \Pi_i$ is
\[
  \frac{\pi^{\frac 1 4 i(2n-i-1)}}{((n-i)!)^{i-1}}
\times
\begin{cases}
\frac 1 {(n-2i)!}
\qquad
&\text{if $2i < n+1$, }
\\
(-1)^{n+1} (2i-n-1)!
\qquad
&\text{if $2i \ge n+1$.  }
\end{cases}
\]
\end{theorem} 

\medskip
We also prove in Chapter \ref{sec:period} a generalization of a theorem of Sun \cite{S}.

\begin{theorem}
[see Theorem \ref{thm:gKOn}]
Let $(G,G')=(O(n+1,1),O(n,1))$, 
 $0 \le i \le n$, 
 and $\delta \in \{ \pm \}$.  
\begin{enumerate}
\item[{\rm{(1)}}]
The symmetry breaking operator
$
  T \colon
  \Pi_{i,\delta} \to \pi_{i,\delta}
$
 in Proposition \ref{prop:AiiAq}
 induces bilinear forms
\[
   B_T \colon 
   H^j({\mathfrak{g}}, K; \Pi_{i,\delta}) 
   \times 
   H^{n-j}({\mathfrak{g}}', K'; \pi_{n-i,(-1)^n \delta})
   \to {\mathbb{C}}
\]
 for all $j$.  
\item[{\rm{(2)}}]
The bilinear form $B_T$ is nonzero
 if and only if $j=i$ and $\delta=(-1)^i$.  

\end{enumerate}
\end{theorem}

\bigskip

Inspired by automorphic forms and number theory B.~Gross and D.~Prasad published  in  1992 conjectures about the multiplicities of irreducible tempered representations $(U,U')$ of $(SO(p,q), SO(p-1,q))$ \cite{GP}.  
Over time these conjectures have been modified and proved
 in some cases for automorphic forms
 and for $p$-adic orthogonal and unitary groups. See for example  Ast\'{e}risque volumes \cite{A1, A2}
 by W.~T.~Gan, B.~Gross, D.~Prasad, C.~M{\oe}glin and J.-L. Waldspurger 
 and the references therein as well as the work
 by R.~Beuzart-Plessis \cite{raphael} for the unitary groups.

We prove the multiplicity conjecture by B.~Gross and D.~Prasad for  tempered principal series representations of $(SO(n+1,1), SO(n,1))$  and also for 3 representations $\Pi, \pi, \varpi$
 of $SO(2m+2,1)$, $SO(2m+1,1)$ and $SO(2m,1)$ with infinitesimal character $\rho$.
More precisely we show:

\begin{theorem}
[Theorem \ref{thm:VWSBO}]  
Suppose that  
 $\Pi= I_\delta (V,\lambda)$, $\pi=J_\varepsilon(W,\nu )$ are (smooth) tempered principal series representations of $G=O(n+1,1)$ and $G'=O(n,1)$. 
Then
\[ {\operatorname{dim}}_{\mathbb{C}} {\operatorname{Hom}}_{G'}
   (\Pi|_{G'},\pi)  = 1  
\quad
\text{if and only if} 
\quad
  [V:W] \ne 0.
\] 
\end{theorem}
Restricting the principal series representations  to special orthogonal groups implies the conjecture of B.~Gross and D.~Prasad about multiplicities 
for tempered principal series representations
 (Theorem \ref{theorem:GPtemp}).

In 2000 B.~Gross and N.~Wallach \cite{GW} showed
 that the restriction of {\it{small}} discrete series representations
 of $G=SO(2p+1,2q)$ to $G'=SO(2p,2q)$ 
satisfies the Gross--Prasad conjectures
 \cite{GP}. 
In that case,
 both the groups $G$ and $G'$ admit discrete series representations.  
On the other hand,
 for the pair $(G,G')=(SO(n+1,1),SO(n,1))$, 
 only one of $G$ or $G'$ admits discrete series representations.  
Our results confirm the Gross--Prasad conjecture
 also for {\bf tempered} representations
 with trivial infinitesimal character $\rho$
 (Theorem \ref{thm:GPdisctemp}).

\bigskip \noindent
The article is roughly divided in three parts and an appendix:

In the first part, 
 Chapters \ref{sec:ps}--\ref{sec:SBOrho}, 
 we give an overview of the notation
 and the results about symmetry breaking operators.  
Notations and properties for principal series and irreducible representations of orthogonal groups are introduced in Chapter \ref{sec:ps}. 
Important concepts and properties of symmetry breaking operators are discussed 
 in Chapter \ref{sec:general}, 
 in particular,
 a classification scheme of all symmetry breaking operators
 is presented in Theorem \ref{thm:VWSBO}.  
This includes a number of theorems about the dimension
 of the space of symmetry operators for principal series representations
 which are stated and discussed also in Chapter \ref{sec:general}. 
The classification scheme is carried out
 in full details for symmetry breaking from principal series representations
 $I_\delta(i,\lambda)$ of $G$ to $J_\varepsilon(j,\nu)$
 of the subgroup $G'$, 
 and is used to obtain results on symmetry breaking of 
 irreducible representations with trivial infinitesimal character $\rho $
 in Chapter \ref{sec:SBOrho}.

The second part,
 Chapters \ref{sec:section7}--\ref{sec:holo},
 contains the proofs of the results discussed in Part one. 
This is the technical heart of this monograph.  
In Chapter \ref{sec:section7} the estimates and results about regular symmetry breaking operators in Theorems \ref{thm:intro160150} and \ref{thm:VW0SBO}
 are proved. 
Chapter \ref{sec:DSVO} is devoted to differential symmetry breaking operators. 
In the remaining chapters of this part we concentrate on the symmetry breaking  $I_\delta(i,\lambda) \rightarrow J_\varepsilon(j,\nu )$. 
We collect some technical results
 in Chapters \ref{sec:section9} and \ref{sec:psdetail}.  
The analytic continuation of the regular symmetry breaking operators, 
 their $(K,K')$-spectrum, 
 and the functional equation are discussed in Chapter \ref{sec:holo}.
Many of the results and techniques developed here are of independent interest, 
 and would be applied to other problems.

In the third part,
 Chapters \ref{sec:Gross-Prasad}--\ref{sec:conjecture}, 
 we use the results in Chapters \ref{sec:general} and \ref{sec:SBOrho}
 to prove some of the conjectures of Gross and Prasad about symmetry breaking for tempered representations of orthogonal groups in Chapter \ref{sec:Gross-Prasad}. 
We discuss periods of representations and a bilinear form
 on the ($\mathfrak{g},K$)-cohomology using symmetry breaking
 in Chapter \ref{sec:period}.  
It also includes a conjecture about symmetry breaking for a family of representations of irreducible representations
 with regular integral infinitesimal character
 in Chapter \ref{sec:conjecture},  
 which we plan to attack in a sequel to this monograph. 
A major portion of Part 3 can be read immediately after Part 1.

The appendix contains technical results used in the monograph. 
We provide three characterizations of irreducible representations
 of the group $G=O(n+1,1)$:
Langlands quotients (or subrepresentations), 
 cohomological parabolic induction, 
 and translation from ${\operatorname{Irr}}(G)_{\rho}$.  
The first two are discussed in Appendix I
 (Chapter \ref{sec:aq})
 and the third one is in Appendix III
 (Chapter \ref{sec:Translation}).  
For the second description,
 we recall the description of the Harish-Chandra modules
 of the irreducible representations of $O(n,1)$
 as the cohomological induction from a $\theta$-stable Levi subgroup
 and introduce $\theta$-stable coordinates for irreducible representations
 with regular integral infinitesimal character. This notation is used in the formulation of the conjecture in Chapter \ref{sec:conjecture}.  
We discuss the restriction of representations of the orthogonal group $O(n,1)$ to the special orthogonal group $SO(n,1)$
 in  Appendix II (Chapter \ref{sec:SOrest}). 
The results are used in Chapter \ref{sec:Gross-Prasad}
 about the Gross--Prasad conjecture.  
In Appendix III, 
 we discuss translation functor of $G=O(n+1,1)$
 which is not in the Harish-Chandra class
 when $n$ is even.

\bigskip
{\bf{Acknowlegements:}}\enspace
 Many of the results were obtained while the authors were supported by the Research in Pairs program   at the Mathematisches Forschungsinstitut MFO in Oberwolfach, Germany.
 
Research by the first author was partially supported by Grant-in-Aid for Scientific 
Research (A)
(25247006) and (18H03669), Japan Society for the Promotion of Science.  

Research by the second author was partially supported by  NSF grant DMS-1500644. 
Part of this research was conducted during a visit of the second author at  the Graduate School of Mathematics of The University of Tokyo, Komaba. 
She would like to thank it for its support and hospitality during her stay.

\vskip 3pc
{\bf{Notation:}}
\index{A}{N0natural@${\mathbb{N}}$}
\index{A}{R1real1@${\mathbb{R}}_+$}
\begin{align*}
& A \setminus B
\quad
&&\text{set theoretic complement of $B$ in $A$}
\\
&{\mathbb{N}}
&&\text{$\{\text{integers $\ge 0$}\}$}
\\
&{\mathbb{N}}_+
&&\text{$\{\text{positive integers}\}$}
\\
&{\mathbb{R}}_+
&&\text{$\{t \in {\mathbb{R}}:t>0\}$}
\\
& {\operatorname{Image}}\, (T)
&&\text{image of the operator $T$}
\\
& {\operatorname{Ker}}\, (T)
&&\text{kernel of the operator $T$}
\\
& E_{i j}
&&\text{the matrix unit}
\\
&[a]
&&\text{the largest integer that does not exceed $a$}
\\
&{\bf{1}}
&&\text{the trivial one-dimensional representation}
\\
& \pi^{\vee}
&&\text{the contragredient representation of $\pi$}
\\
& \pi_1 \boxtimes \pi_2
&&\text{the outer tensor product representation of a direct product group}
\\
& \pi_1 \otimes \pi_2
&&\text{the tensor product representation}
\\
& \rho (\equiv \rho_G)
&&\text{the infinitesimal character of the trivial representation ${\bf{1}}$}
\end{align*}

\newpage
\section{Review of principal series representations}
\label{sec:ps}
In this chapter
 we recall results about representations
 of the indefinite orthogonal group 
 $G = O(n+1,1)$.  

\subsection{Notation}
\label{subsec:Xi}

The object of our study is intertwining restriction operators
 ({\it{symmetry breaking operators}})
 between representations of $G=O(n+1,1)$
 and those of its subgroup $G'=O(n,1)$.  
Most of main results are stated 
 in a coordinate-free fashion, 
 whereas concrete description of symmetry breaking operators
 depends on coordinates.  
For the latter purpose,
 we choose subgroups of $G$ and $G'$
 in a compatible fashion.  
The notations here are basically taken from \cite{sbon}.  

\subsubsection{Subgroups of $G=O(n+1,1)$
 and $G'=O(n,1)$}
\label{subsec:subgpOn}
We define $G$ to be the indefinite orthogonal group 
 $O(n+1,1)$
 that preserves the quadratic form 
\begin{equation}
\label{eqn:Qn1}
Q_{n+1,1}(x)=x_0^2+ \cdots + x_{n}^2 - x_{n+1}^2
\end{equation}
of signature $(n+1,1)$.  
Let $G'$ be the stabilizer of the vector
 $e_n:= {}^{t\!}(0, \cdots,0,1,0)$.  
Then $G' \simeq O(n,1)$.

We take maximal compact subgroups of $G$ and $G'$, 
 respectively, 
 as 
\index{A}{Kmaxcpt@$K=O(n+1) \times O(1)$, maximal compact subgroup of $G$|textbf}
\index{A}{Kmaxcptsub@$K'= K \cap G'$|textbf}
\begin{alignat*}{2}
K 
&:= O(n+2) \cap G
&&
\simeq O(n+1) \times O(1),
\\
K'
&:= K \cap G'
=\left\{
\begin{pmatrix}
                A & & \\ & 1 & \\ & & \varepsilon
\end{pmatrix}
:
 A \in O(n),\,\, \varepsilon =\pm 1
\right\}
&&
\simeq O(n) \times O(1).  
\end{alignat*}

Let ${\mathfrak{g}} = {\mathfrak{o}}(n+1,1)$
 and ${\mathfrak{g}}' = {\mathfrak{o}}(n,1)$
 be the Lie algebras
 of 
\index{A}{GOindefinitelarge@$G=O(n+1,1)$, Lorentz group|textbf}
$G=O(n+1,1)$
 and 
\index{A}{GOindefinitesprime@$G'=O(n,1)$|textbf}
$G'=O(n,1)$, 
 respectively.  
We take a hyperbolic element 
\index{A}{H@$H$|textbf}
\begin{equation}
\label{eqn:Hyp}
 H:= E_{0,n+1} + E_{n+1,0} \in {\mathfrak{g}}', 
\end{equation}
 and set 
\index{A}{0caLiealg@${\mathfrak{a}}$, maximally split abelian subspace|textbf}
\[
  {\mathfrak {a}}:={\mathbb{R}}H.  
\]
Then ${\mathfrak{a}}$ is a maximally split abelian subspace
 of ${\mathfrak{g}}'$, 
 as well as that of ${\mathfrak{g}}$.
The eigenvalues of 
 $\operatorname{ad} (H) \in \operatorname{End}({\mathfrak{g}})$
 are $\pm 1$ and $0$, 
 and the eigenspaces give rise to the following two maximal nilpotent subalgebras
 of ${\mathfrak {g}}$:
\begin{equation}
\label{eqn:ngen}
{\mathfrak {n}}_+=\operatorname{Ker} (\operatorname{ad}(H)-1)
                 =\sum_{j=1}^n {\mathbb{R}} N_j^+, 
\quad
{\mathfrak {n}}_-=\operatorname{Ker} (\operatorname{ad}(H)+1)
                 =\sum_{j=1}^n {\mathbb{R}} N_j^-, 
\end{equation}
where $N_j^+$ and $N_j^-$
 ($1 \le j \le n$) are nilpotent elements
 of ${\mathfrak{g}}$ defined by 
\index{A}{Nj1@$N_j^+$|textbf}
\index{A}{Nj2@$N_j^-$|textbf}
\begin{align*}
N_j^+=&-E_{0,j} + E_{j,0} - E_{j,n+1} -E_{n+1,j}, 
\\
N_j^-=&-E_{0,j} + E_{j,0} + E_{j,n+1} +E_{n+1,j}.  
\end{align*}

\medskip
For $b={}^t(b_1, \cdots, b_n) \in {\mathbb{R}}^n$, 
 we define unipotent matrices in $G$ by
\index{A}{n1@$n_+\colon {\mathbb{R}}^n \to N_+$|textbf}
\index{A}{n2@$n_-\colon {\mathbb{R}}^n \to N_-$|textbf}
\begin{align}
n_+(b):=& \exp(\sum_{j=1}^n b_j N_j^+)=I_{n+2} 
       + 
\begin{pmatrix} 
     - \frac 1 2 Q(b) &-{}^{t\!}b & \frac 1 2 Q(b) 
\\
       b & 0 & -b
\\
     - \frac 1 2 Q(b) &-{}^{t\!}b & \frac 1 2 Q(b) 
\end{pmatrix},
\label{eqn:nplus}
\\
n_-(b):=&\exp(\sum_{j=1}^n b_j N_j^-)= I_{n+2} 
       + 
\begin{pmatrix} 
     - \frac 1 2 Q(b) &-{}^{t\!}b & -\frac 1 2 Q(b) 
\\
       b & 0 & b
\\
     \frac 1 2 Q(b) &{}^{t\!}b & \frac 1 2 Q(b) 
\end{pmatrix}, 
\label{eqn:nbar}
\end{align}
where we set
\begin{equation}
\label{eqn:Qn}
Q(b)\equiv |b|^2 =\sum_{l=1}^n b_l^2.  
\end{equation}
Then 
\index{A}{n1@$n_+\colon {\mathbb{R}}^n \to N_+$|textbf}
$n_+$
 and 
\index{A}{n2@$n_-\colon {\mathbb{R}}^n \to N_-$|textbf}
$n_-$ give coordinates
 of the nilpotent groups 
\index{A}{N1@$N_+ = \exp ({\mathfrak{n}}_+)$|textbf}
 $N_+ := \exp ({\mathfrak{n}}_+)$ and 
\index{A}{N2@$N_- = \exp ({\mathfrak{n}}_-)$|textbf}
 $N_- := \exp ({\mathfrak{n}}_-)$, 
 respectively.  
Then $N_+$ stabilizes ${}^{t\!}(1,0,\cdots,0,1)$, 
 whereas $N_-$ stabilizes ${}^{t\!}(1,0,\cdots,0,-1)$.

Since $H$ is contained in the Lie algebra ${\mathfrak {g}}'$, 
\[
  {\mathfrak{n}}_{\varepsilon}'
  :=
  {\mathfrak{n}}_{\varepsilon} \cap {\mathfrak{g}}'
  =
  \sum_{j=1}^{n-1} {\mathbb{R}} N_j^{\varepsilon}
\quad
  \text{for ${\varepsilon}=\pm$}
\]
are maximal nilpotent subalgebras 
 of ${\mathfrak{g}}'$.  
We set 
\index{A}{N1sub@$N_+' = N_+ \cap G'$|textbf}
 $N_+' := N_+ \cap G' = \exp ({\mathfrak{n}}_+')$
 and 
\index{A}{N2sub@$N_-' = N_- \cap G'$|textbf}
 $N_-' := N_- \cap G' = \exp ({\mathfrak{n}}_-')$.

We define a split abelian subgroup $A$
 and its centralizers $M$ and $M'$ in $K$ and $K'$, 
 respectively, 
 as follows:
\index{A}{M1@$M=O(n) \times O(1)$, the centralizer of ${\mathfrak {a}}$ in $K$|textbf}
\index{A}{M1prime@$M'=O(n-1) \times O(1)$|textbf}
\begin{alignat*}{2}
A&:= \exp ({\mathfrak {a}}), 
&&
\\
M&:=\left\{ \begin{pmatrix} \varepsilon & & \\ &  B & \\ & & \varepsilon \end{pmatrix} : B \in O(n), \varepsilon = \pm 1\right\}
&& \simeq O(n) \times {\mathbb{Z}}/2 {\mathbb{Z}}, 
\\
M'
&:=\left\{ \begin{pmatrix} \varepsilon & & & \\ &  B & & \\ & & 1 & \\ & & & \varepsilon\end{pmatrix} : B \in O(n-1), \varepsilon = \pm 1\right\}
&& \simeq O(n-1) \times {\mathbb{Z}}/2 {\mathbb{Z}}.  
\end{alignat*}
Then 
\index{A}{PLanglandsdecomp@$P=MAN_+$, Langlands decomposition of 
a minimal parabolic subgroup of $G$|textbf}
$P= MAN_+$ is a Langlands decomposition
 of a minimal parabolic subgroup $P$ of $G$.  
Likewise,
\index{A}{PLanglandsdecompsubgp@$P'=M'AN_+'$|textbf}
 $P'= M' AN_+'$ is that
 of a minimal parabolic subgroup $P'$ of $G'$.  
We note
 that $A$ is a common maximally split abelian subgroup
 in $P'$ and $P$
 because we have chosen $H \in {\mathfrak {g}}'$.  
The Langlands decompositions
 of the Lie algebras of $P$ and $P'$
 are given in a compatible way 
 as 
\[
{\mathfrak {p}}={\mathfrak {m}}+{\mathfrak {a}}+{\mathfrak {n}}_+, 
\quad
 {\mathfrak {p}}'={\mathfrak {m}}'+{\mathfrak {a}}+{\mathfrak {n}}_+'
 =({\mathfrak {m}} \cap {\mathfrak {g}}') + ({\mathfrak {a}} \cap {\mathfrak {g}}')+({\mathfrak {n}}_+ \cap {\mathfrak {g}}').  
\]
We set 
\index{A}{m2@$m_-={\operatorname{diag}}(-1,1,\cdots,1,-1)$|textbf}
\begin{equation}\label{eqn:m-}
m_- \ := \ 
  \left(
   \begin{array}{cc|c}
      -1  &        &   \\
        & I_{n}&      \\
           \hline
        &        &      -1
    \end{array}
    \right)\in M'.  
   \end{equation}
We note
 that $m_-$ does not belong to the identity component
 of $G'$.

\subsubsection{Isotropic cone $\Xi$}
\label{subsec:Xicone}
The isotropic cone
\index{A}{0nXi@$\Xi$, isotropic cone|textbf}
\[
     \Xi\equiv \Xi({\mathbb{R}}^{n+1,1})=
     \{(x_0, \cdots, x_{n+1}) \in {\mathbb{R}}^{n+2}:
        x_0^2+ \cdots + x_{n}^2 - x_{n+1}^2=0\}
     \setminus
     \{0\}    
\]
 is a homogeneous $G$-space 
 with the following fibration: 
\begin{alignat*}{6}
   &G/O(n) N_+ \,\,&& \simeq\,\,  && \Xi
   && \qquad\quad g O(n)N_+ && \mapsto && \ g p_+
\\
&  {\mathbb{R}}^{\times}\downarrow && && \downarrow {\mathbb{R}}^{\times}
   && \qquad\quad\qquad \rotatebox[origin=c]{-90}{$\mapsto$}  && && \  \rotatebox[origin=c]{-90}{$\mapsto$}
\\
     &\,G/P && \simeq \,\, && S^n,
     &&
     \qquad\quad\quad\,\,  gP &&\mapsto &&  [g p_+]
\end{alignat*}
where 
\index{A}{p1@$p_+={}^t(1,0,\cdots,0,1)$|textbf}
\begin{equation}
\label{eqn:p+}
p_+ :={} {}^{t\!} (1, 0, \cdots, 0, 1) \in \Xi. 
\end{equation}

The action of the subgroup $N_+$ on the isotropic cone
 $\Xi$ is given in the coordinates as 
\begin{equation}
\label{eqn:nbxi}
  n_+(b) \begin{pmatrix} \xi_0 \\ \xi \\ \xi_{n+1} \end{pmatrix}
 =
\begin{pmatrix}
 \xi_0 -(b,\xi)
\\
 \xi 
\\ 
 \xi_{n+1}-(b,\xi)
\end{pmatrix}
+
\frac{\xi_{n+1} -\xi_0}{2}
\begin{pmatrix} Q(b)\\ -2 b\\ Q(b)\end{pmatrix}, 
\end{equation}
where
 $b \in {\mathbb{R}}^n$, 
 $\xi \in {\mathbb{R}}^n$ and $\xi_0, \xi_{n+1}\in\mathbb{R}$.

The intersections
 of the isotropic cone $\Xi$
 with the hyperplanes
 $\xi_0 + \xi_{n+1}=2$
 or $\xi_{n+1}=1$
 can be identified with ${\mathbb{R}}^n$ or $S^n$, 
 respectively.  
We write down the embeddings
 $\iota_N:{\mathbb{R}}^n \hookrightarrow \Xi$
 and $\iota_K:S^n \hookrightarrow \Xi$
 in the coordinates as follows:
\begin{align}
\label{eqn:NXi}
\iota_N:&
\mathbb{R}^n \hookrightarrow \Xi,
\
{}^t(x,x_n) \mapsto 
{n_-} (x,x_n) p_+
   = \left(\begin{array}{l}
           1 - |x|^2 - x_n^2  \\
           2x                     \\
           2x_n                  \\
           1 + |x|^2 + x_n^2
     \end{array}\right), 
\\
\label{eqn:KXi}
\iota_K\colon&
S^n \to \Xi, 
\
\eta \mapsto (\eta,1).  
\end{align}

The composition
 of $\iota_N$ and the projection
\begin{equation*}
   \Xi \to \Xi/\mathbb{R}^\times \overset{\sim}{\to} S^n,
\quad
\xi \mapsto \frac{1}{\xi_{n+1}} (\xi_0,\dots,\xi_n)
\end{equation*}
yields the conformal compactification of ${\mathbb{R}}^n$:
\begin{equation}
\label{eqn:XiSn}
\mathbb{R}^n \hookrightarrow S^n,
\quad
r \omega \mapsto \eta = (s, \sqrt{1-s^2} \, \omega ) 
= \Bigl( \frac{1-r^2}{1+r^2}, \frac{2r}{1+r^2} \omega \Bigr).  
\end{equation}
Here $\omega \in S^{n-1}$
 and the inverse map is given by 
$
r = \sqrt{\frac{1-s}{1+s}}
$
 for $s \ne -1$.

\subsubsection{Characters $\chi_{\pm\pm}$ of the component group $G/G_0$}
\label{subsec:chcomp}

There are four connected components
 of the group $G=O(n+1,1)$.  
Let 
\index{A}{GOindefinitespecialidentity@$G_0=SO_0(n+1,1)$, the identity component
 of $O(n+1,1)$\quad|textbf}
$G_0$ denote the identity component of $G$.  
Then $G_0\simeq SO_0(n+1,1)$ and the quotient group $G/G_0$
\index{B}{componentgroup@component group $G/G_0$|textbf}
({\it{component group}})
 is isomorphic
 to ${\mathbb{Z}}/2 {\mathbb{Z}} \times {\mathbb{Z}}/2 {\mathbb{Z}}$.  
Accordingly, there are four one-dimensional representations of $G$, 
\index{A}{1chipmpm@$\chi_{\pm\pm}$, one-dimensional representation of $O(n+1,1)$|textbf}
\begin{equation}
\label{eqn:chiab}
\chi_{ab}\colon  G\to\{\pm1\}
\end{equation}
with $a,b\in\{\pm\}\equiv\{\pm1\}$ such that
\begin{equation*}
\chi_{ab}\left(\mathrm{diag}(-1,1,\cdots,1)\right)=a,\quad 
\chi_{ab}\left(\mathrm{diag}(1,\cdots,1,-1)\right)=b.
\end{equation*}
We note that $\chi_{--}$ is given 
 by the determinant, 
 $\mathrm{det}$, 
 of matrices in $O(n+1,1)$. 
Then the restriction of 
\index{A}{1chipmpm@$\chi_{\pm\pm}$, one-dimensional representation of $O(n+1,1)$|textbf}
$\chi_{--}$ to the subgroup $M\simeq O(n) \times O(1)$ is given by the outer tensor product representation:
\index{A}{1chipmm@$\chi_{--}=\det$}
\begin{equation}\label{eqn:chiM}
\chi_{--}\vert_{M}\simeq \det\boxtimes \,\mathbf{1},  
\end{equation}
where $\det$ in the right-hand side
 stands for the determinant
 for $n$ by $n$ matrices.

\subsubsection{The center ${\mathfrak{Z}}_G({\mathfrak{g}})$ 
 and the Harish-Chandra isomorphism}
\label{subsec:2.1.4}

For a Lie algebra ${\mathfrak{g}}$ over ${\mathbb{R}}$, 
 we denote by 
\index{A}{Ug@$U({\mathfrak{g}})$, enveloping algebra|textbf}
$U({\mathfrak{g}})$
 the universal enveloping algebra
 of the complexified Lie algebra 
$
   {\mathfrak{g}}_{\mathbb{C}}={\mathfrak{g}} \otimes_{\mathbb{R}} {\mathbb{C}}
$, 
 and by
\index{A}{Zg@${\mathfrak{Z}}({\mathfrak{g}})$|textbf}
 ${\mathfrak{Z}}({\mathfrak{g}})$ its center.  
For a real reductive Lie group $G$ with Lie algebra ${\mathfrak{g}}$, 
 we define a subalgebra of ${\mathfrak{Z}}({\mathfrak{g}})$
 of finite index 
 by 
\index{A}{ZGg@${\mathfrak{Z}}_G({\mathfrak{g}})$|textbf}
\[
 {\mathfrak{Z}}_G({\mathfrak{g}})
  :=
  U({\mathfrak{g}})^G
  =\{z \in U({\mathfrak{g}}) : 
    {\operatorname{Ad}}(g)z=z
\quad
\text{for all } g \in G\}.  
\]
Schur's lemma implies 
 that the algebra ${\mathfrak{Z}}_G({\mathfrak{g}})$ acts
 on any irreducible admissible smooth representation of $G$
 by scalars, 
 which we refer to as 
\index{B}{ZGginfinitesimalcharacter@${\mathfrak{Z}}_G({\mathfrak{g}})$-infinitesimal character|textbf}
 the ${\mathfrak{Z}}_G({\mathfrak{g}})$-{\it{infinitesimal character}}.  
If the reductive group $G$ is of Harish-Chandra class, 
 then the adjoint group ${\operatorname{Ad}}(G)$
 is contained in the inner automorphism group
 ${\operatorname{Int}}({\mathfrak{g}}_{{\mathbb{C}}})$, 
 and consequently, 
 ${\mathfrak{Z}}_G({\mathfrak{g}}) ={\mathfrak{Z}}({\mathfrak{g}})$.  
However, 
 special attention is required 
 when $G$ is not of Harish-Chandra class, 
 as we shall see below.

For the disconnected group $G=O(n+1,1)$, 
 ${\operatorname{Ad}}(G)$ is not contained in ${\operatorname{Int}}({\mathfrak{g}}_{{\mathbb{C}}})$
 and ${\mathfrak{Z}}_G({\mathfrak{g}})$ is of index two in 
 ${\mathfrak{Z}}({\mathfrak{g}})$
if $n$ is even, 
 whereas ${\operatorname{Ad}}(G) \subset {\operatorname{Int}}({\mathfrak{g}}_{{\mathbb{C}}})$ and ${\mathfrak{Z}}_G({\mathfrak{g}})={\mathfrak{Z}}({\mathfrak{g}})$
 if $n$ is odd.  
In both cases,
 via the standard coordinates
 of a Cartan subalgebra of ${\mathfrak{g}}_{{\mathbb{C}}} \simeq {\mathfrak{o}}
(n+2,{{\mathbb{C}}})$, 
 we have the following
\index{B}{HarishChandraisomorphism@Harish-Chandra isomorphism|textbf}
 Harish-Chandra isomorphisms
\begin{equation*}
\begin{array}{ccc}
   {\mathfrak{Z}}({\mathfrak{g}})
   & \simeq
   & S({{\mathbb{C}}}^{m+1})^{W_{\mathfrak{g}}}
\\
   \cup 
   & 
   & \cup
\\
  {\mathfrak{Z}}_G({\mathfrak{g}})
   & \simeq
   &\,\, S({{\mathbb{C}}}^{m+1})^{W_G}.  
\end{array}
\end{equation*}
Here  we identify a Cartan subalgebra ${\mathfrak{h}}_{{\mathbb{C}}}$
 of ${\mathfrak{g}}_{{\mathbb{C}}}\simeq {\mathfrak{o}}(n+2,{\mathbb{C}})$
 with ${{\mathbb{C}}}^{m+1}$
 where $m:=[\frac n 2]$, 
 and set 
\index{A}{Weylgroupg@$W_{\mathfrak{g}}$, Weyl group for ${\mathfrak{g}}_{\mathbb{C}}={\mathfrak{o}}(n+2,{\mathbb{C}})$|textbf}
\index{A}{WeylgroupG@$W_G$, Weyl group for $G=O(n+1,1)$|textbf}
\begin{align*}
W_{\mathfrak{g}}
:=& W(\Delta({\mathfrak{g}}_{\mathbb{C}},{\mathfrak{h}}_{\mathbb{C}}))
\simeq
\begin{cases}
{\mathfrak{S}}_{m+1} \ltimes ({\mathbb{Z}}/2 {\mathbb{Z}})^{m+1}
\quad
&\text{for $n=2m+1$}, 
\\
{\mathfrak{S}}_{m+1} \ltimes ({\mathbb{Z}}/2 {\mathbb{Z}})^{m}
\quad
&\text{for $n=2m$}, 
\end{cases}
\\
W_G:=& {\mathfrak{S}}_{m+1} \ltimes ({\mathbb{Z}}/2 {\mathbb{Z}})^{m+1}. 
\end{align*}

We shall describe
 the
\index{B}{infinitesimalcharacter@infinitesimal character|textbf}
 ${\mathfrak{Z}}_G({\mathfrak{g}})$-infinitesimal character
 by an element of ${{\mathbb{C}}}^N$ modulo $W_G$
 via the following isomorphism.  

\begin{equation}
\begin{array}{ccc}
    {\operatorname{Hom}}_{{\mathbb{C}}\operatorname{-alg}}({\mathfrak{Z}}({\mathfrak{g}}), {{\mathbb{C}}}) 
&\simeq 
&{\mathbb{C}}^N/W_{\mathfrak{g}}
\\
\rotatebox[origin=c]{-90}{$\twoheadrightarrow$}
&
&\rotatebox[origin=c]{-90}{$\twoheadrightarrow$}
\\
{\operatorname{Hom}}_{{\mathbb{C}}\operatorname{-alg}}({\mathfrak{Z}}_G({\mathfrak{g}}), {{\mathbb{C}}}) 
&\simeq 
&{\mathbb{C}}^N/W_G
\\
\label{eqn:HCpara}
\end{array}
\end{equation}

To define the notion of \lq\lq{regular}\rq\rq\
 or \lq\lq{singular}\rq\rq\
 about ${\mathfrak{Z}}_G({\mathfrak{g}})$-infinitesimal characters, 
 we use the action of the Weyl group $W_{\mathfrak{g}}$
 for the Lie algebra 
 ${\mathfrak{g}}_{\mathbb{C}}={\mathfrak{o}}(n+2,{\mathbb{C}})$
 rather than the Weyl group $W_G$
 for the disconnected group $G$ as below.  
 
\begin{definition}
\label{def:intreg}
Let $G=O(n+1,1)$ and $m:=[\frac n 2]$.  
Suppose $\chi\in {\operatorname{Hom}}_{{\mathbb{C}}\operatorname{-alg}}({\mathfrak{Z}}_G({\mathfrak{g}}), {{\mathbb{C}}})$
 is given by $\mu \in {\mathbb{C}}^{m+1} \mod W_G$ 
 via the Harish-Chandra isomorphism \eqref{eqn:HCpara}.  
We say $\chi$ is {\it{integral}}
 if 
\[
  \mu -\rho_G \in {\mathbb{Z}}^{m+1}, 
\]
see \eqref{eqn:rhoG} below
 for the definition of $\rho_G$, 
 or equivalently,
 if 
\begin{alignat*}{3}
\mu &\in {\mathbb{Z}}^{m+1}
\qquad
&&
\text{for $n=2m$}
\qquad
&&
\text{(even), }
\\
\mu &\in ({\mathbb{Z}}+\frac 1 2)^{m+1}
\qquad
&&
\text{for $n=2m+1$}
\qquad
&&
\text
{(odd).  }
\end{alignat*}
We note that this condition is stronger
 than the one which is usually referred
 to as \lq\lq{integral}\rq\rq:
\[
   \langle \mu, \alpha^{\vee} \rangle \in {\mathbb{Z}}
   \quad
  \text{for all $\alpha \in \Delta({\mathfrak{g}}_{\mathbb{C}}, {\mathfrak{h}}_{\mathbb{C}})$}
\]
where  $\alpha^{\vee}$ denotes the coroot of $\alpha$.  

For $\mu \in {\mathbb{C}}^{m+1}$, 
 we set 
\begin{align*}
  W_{\mu} \equiv (W_{\mathfrak{g}})_{\mu} 
          &:=\{w \in W_{\mathfrak{g}}: w \mu =\mu \}, 
\\
(W_G)_{\mu} 
          &:=\{w \in W_G: w \mu =\mu \}.  
\end{align*}
We say $\mu$ is 
\index{B}{Wgregular@$W_{\mathfrak{g}}$-regular|textbf}
{\it{$W_{\mathfrak{g}}$-regular}}
 (or simply, {\it{regular}})
 if $(W_{\mathfrak{g}})_{\mu}=\{e\}$, 
 and 
\index{B}{WGregular@$W_G$-regular|textbf}
{\it{$W_G$-regular}}
 if $(W_G)_{\mu}=\{e\}$.  
These definitions depend only
 on the $W_G$-orbit through $\mu$ 
 because $\# W_{\mu}= \# W_{\mu'}$ if $\mu' \in W_G \mu$.  
We say $\chi$ is \index{B}{regularintegralinfinitesimalcharacter@regular integral infinitesimal character|textbf}
 {\it{regular integral}}
 (respectively, 
\index{B}{singularintegralinfinitesimalcharacter@singular integral infinitesimal character|textbf}
{\it{singular integral}}) 
 infinitesimal character
 if $\chi$ is integral
 and $W_{\mu}=\{e\}$
 (respectively, $W_{\mu}\ne \{e\}$).  
In the coordinates of $\mu=(\mu_1, \cdots, \mu_{m+1})$, 
 $W_{\mu}=\{e\}$
 if and only if 
\begin{alignat*}{3}
\mu_i &\ne \pm \mu_j 
\quad
&&(1 \le \forall i < \forall j \le m+1)
&&\text{for $n$ even}, 
\\
\mu_i &\ne \pm \mu_j 
&&(1 \le \forall i < \forall j \le m+1), 
\,\,
\mu_k \ne 0
\quad
(1 \le \forall k \le m+1)
\quad
&&\text{for $n$ odd}.  
\end{alignat*}
\end{definition}

\begin{remark}
\label{rem:Fregint}
Suppose $G=O(n+1,1)$ with $n \ge 1$.  
Then the ${\mathfrak{Z}}_G({\mathfrak{g}})$-infinitesimal character
 of an irreducible finite-dimensional representation 
 of $G$ is regular integral, 
 and conversely,
 for any regular integral $\chi$, 
 there exists an irreducible finite-dimensional representation $F$
 of $G$ such that $\chi$ is the ${\mathfrak{Z}}_G({\mathfrak{g}})$-infinitesimal character of $F$.  
Here we remind from Definition \ref{def:intreg}
 above that by \lq\lq{regular}\rq\rq\
 we mean $W_{\mathfrak{g}}$-regular,
 and not $W_G$-regular.  
\end{remark}

The ${\mathfrak{Z}}_G({\mathfrak{g}})$-infinitesimal character
 of the trivial one-dimensional representation ${\bf{1}}$ of $G=O(n+1,1)$ 
 is given by 
\index{A}{1parho@$\rho_G$|textbf}
\begin{equation}
\label{eqn:rhoG}
     \rho 
\equiv 
     \rho_G =(\frac n 2, \frac n 2-1, \cdots, \frac n 2-[\frac n 2])
\in {\mathbb{C}}^{[\frac n 2]+1}/W_G.  
\end{equation}
The infinitesimal character $\rho_G$ will be also 
 referred to as the 
\index{B}{trivialinfinitesimalcharacter@trivial infinitesimal character|textbf}
 {\it{trivial infinitesimal character}}.  
\begin{definition}
\label{def:Irrrho}
We denote by 
\index{A}{IrrGrho@${\operatorname{Irr}}(G)_{\rho}$,
 set of irreducible admissible smooth representation of $G$
 with trivial infinitesimal character $\rho$\quad|textbf}
${\operatorname{Irr}}(G)_{\rho}$
 the set of equivalence classes
 of irreducible admissible smooth representations of $G$
 that have the trivial infinitesimal character $\rho$.  
\end{definition}
The finite set ${\operatorname{Irr}}(G)_{\rho}$ is classified 
 in Theorem \ref{thm:LNM20} (2)
 for $G=O(n+1,1)$
 and in Proposition \ref{prop:161648} (3)
 for the special orthogonal group $SO(n+1,1)$.

\subsection{Representations of the orthogonal group $O(N)$}
\label{subsec:repON}
We recall that the orthogonal group $O(N)$ has two connected components.  
In this section,
 we review a parametrization 
 of irreducible finite-dimensional representations
 of the {\it{disconnected}} group $O(N)$
 following Weyl 
 \cite[Chap.~V, Sect.~7]{Weyl97}.  
For later reference we include classical branching theorems 
 for the restriction 
 of representations
 for the pairs $O(N) \supset O(N-1)$
 and $O(N) \supset SO(N)$.  
The results will be applied to the four compact subgroups $K$, 
\index{A}{Kmaxcptsub@$K'= K \cap G'$}
$K'$, 
\index{A}{M1@$M=O(n) \times O(1)$, the centralizer of ${\mathfrak {a}}$ in $K$}
$M$ and 
\index{A}{M1prime@$M'=O(n-1) \times O(1)$}
$M'=M \cap K'$ of $G$
 introduced in Section \ref{subsec:subgpOn}, 
 which satisfy the following obvious inclusive relations:
\[
\begin{pmatrix}
 K & \supset & K'
\\
\cup & & \cup
\\
 M & \supset & M'
\end{pmatrix}
=
\begin{pmatrix}
 O(n+1) \times O(1) & \supset & O(n) \times O(1)
\\
\cup & & \cup
\\
 O(n) \times {\operatorname{diag}}(O(1)) & \supset & O(n-1) \times {\operatorname{diag}}(O(1))
\end{pmatrix}.  
\]

\subsubsection{Notation for irreducible representations of $O(N)$}
\label{subsec:ONWeyl}
For finite-dimensional irreducible representations
 of orthogonal groups, 
 we use the following notation. 
We set
\index{A}{0tLambda@$\Lambda^+(N)$, dominant weight|textbf}
\begin{equation}
\label{eqn:Lambda}
\Lambda^+(N):=\{ \lambda = (\lambda_1, \ldots, \lambda_N) \in {\mathbb{Z}}^N:
\lambda_1 \geq \lambda_2 \geq \cdots \geq \lambda_N \geq 0\}.
\end{equation}
We write 
\index{A}{FUNl@$\Kirredrep{U(N)}{\lambda}$, irreducible representation of $U(N)$ with highest weight $\lambda$|textbf}
 $\Kirredrep {U(N)}{\lambda}$
 for the irreducible finite-dimensional representation 
of $U(N)$ (or equivalently, the irreducible polynomial representation
 of $GL(N,{\mathbb{C}})$)
with highest weight $\lambda \in \Lambda^+(N)$. 
If $\lambda$ is of the form 
\begin{equation*}
(\underbrace{c_1,\cdots,c_1}_{m_1},
\underbrace{c_2,\cdots,c_2}_{m_2},\cdots,\underbrace{c_\ell,\cdots,c_\ell}_{m_\ell},0,\cdots,0),
\end{equation*}
then we also write $\lambda=\left(c_1^{m_1},c_2^{m_2},\cdots,c_\ell^{m_\ell}\right)$ as usual.

We define a subset of $\Lambda^+(N)$ by
\index{A}{0tExterior@$\Lambda^+(O(N))$|textbf}
\begin{equation*}
\Lambda^+(O(N)):= \{\lambda \in \Lambda^+(N):
\lambda_1' + \lambda_2' \leq N\},
\end{equation*}
where $\lambda_1':=\max\{i : \lambda_i \geq 1 \}$ and 
$\lambda_2': = \max\{i : \lambda_i \geq 2\}$ 
for $\lambda =(\lambda_1, \ldots, \lambda_N) \in \Lambda^+(N)$.
We note that $\lambda_1'$ equals the maximal column length
 in the corresponding Young diagram.

It is readily seen that $\Lambda^+(O(N))$ consists of elements of 
the following two types:
\index{B}{typeON1@type I, for $\Lambda^+(O(N))$|textbf}
\index{B}{typeON2@type II, for $\Lambda^+(O(N))$|textbf}
\begin{align}
\label{eqn:TypeI}
&\text{Type I: $(\lambda_1,\cdots,\lambda_k,\underbrace{0,\cdots,0}_{N-k})$,}
\\
\label{eqn:TypeII}
&\text{Type II: $(\lambda_1,\cdots,\lambda_k,\underbrace{1,\cdots,1}_{N-2k},\underbrace{0,\cdots,0}_k)$,}
\end{align}
 with $\lambda_1\geq \lambda_2\geq\cdots\geq \lambda_k > 0$ and $0\leq k\leq\left[\frac N2\right]$.

For any $\lambda \in \Lambda^+(O(N))$, 
 there exists a unique $O(N)$-irreducible summand, 
 to be denoted by 
\index{A}{FONl@$\Kirredrep{O(N)}{\lambda}$|textbf}
$\Kirredrep {O(N)}{\lambda}$, 
 of the $U(N)$-module $\Kirredrep {U(N)}{\lambda}$
 which contains the highest weight vector.  
Following Weyl (\cite[Chap.~V, Sect.~7]{Weyl97}),
we parametrize the set $\widehat{O(N)}$ of equivalence classes
of irreducible representations of $O(N)$ by
\begin{equation} 
\label{eqn:CWOn}
\Lambda^+(O(N)) \stackrel{\sim}{\longrightarrow} \widehat{O(N)},
\quad
\lambda \mapsto \Kirredrep {O(N)}{\lambda}.  
\end{equation}
By the Weyl unitary trick, 
 we may identify
 $\Kirredrep {O(N)}{\lambda}$ 
 with a holomorphic irreducible representation
 of $O(N,{\mathbb{C}})$, 
 to be denoted by 
\index{A}{FONCl@$\Kirredrep{O(N,{\mathbb{C}})}{\lambda}$}
 $\Kirredrep {O(N,{\mathbb{C}})}{\lambda}$, 
 on the same representation space.

\begin{definition}
\label{def:type}
We say $\Kirredrep{O(N)}{\lambda} \in \widehat{O(N)}$ is of type I
 (or type II), 
 if $\lambda \in \Lambda^+(O(N))$ is of 
\index{B}{type1@type I, representation of $O(N)$}
 type I (or type II), 
respectively.  
\end{definition}

We shall identify $\widehat{O(N)}$
 with $\Lambda^+(O(N))$ via \eqref{eqn:CWOn}, 
 and by abuse of notation, 
 we write $\sigma =(\sigma_1, \cdots, \sigma_N)\in \widehat{O(N)}$
 when $(\sigma_1, \cdots, \sigma_N)\in\Lambda^+(O(N))$.  

\begin{remark}
We shall also use the notation
\begin{align*}
&\Kirredrep{O(N)}
      {\sigma_1, \cdots, \sigma_k, \underbrace{0,\cdots,0}_{[\frac N2]-k}}_+
\text{ instead of }
 \Kirredrep{O(N)}
      {\sigma_1, \cdots, \sigma_k, \underbrace{0,\cdots,0}_{N-k}}, 
\\
&\Kirredrep{O(N)}
      {\sigma_1, \cdots, \sigma_k, \underbrace{0,\cdots,0}_{[\frac N2]-k}}_-
\text{ instead of }
 \Kirredrep{O(N)}
      {\sigma_1, \cdots, \sigma_k, \underbrace{1, \cdots,1}_{N-2k},\underbrace{0,\cdots,0}_{k}}, 
\end{align*}
 by putting the subscript $+$ or $-$
 for irreducible representations of type I
 or of type II,
respectively, 
 see Remark \ref{rem:FOn} in Appendix I.  
\end{remark}
We define a map by summing up 
the first $k$-entries 
 ($k \le [\frac N 2]$) of $\sigma$:
\index{A}{lengthrep@$\ell(\sigma)$|textbf}
\begin{equation}
\label{eqn:ONlength}
\ell \colon
\Lambda^+(O(N)) \to {\mathbb{N}}, 
\quad
\sigma
\mapsto 
\ell(\sigma):=\sum_{i=1}^{k}\sigma_i, 
\end{equation}
which induces a map
\[
  \ell \colon \widehat {O(N)} \to {\mathbb{N}}
\]
via the identification \eqref{eqn:CWOn}.  
By \eqref{eqn:type1to2}, 
 we have 
\begin{equation}
\label{eqn:lsigma}
\ell(\sigma)=\ell(\sigma \otimes \det).  
\end{equation}

\subsubsection{Branching laws for $O(N) \downarrow SO(N)$}
\label{subsec:OnSOn}
\begin{definition}
\label{def:OSO}
We say $\sigma \in \widehat{O(N)}$ is of 
\index{B}{typeX@type X, representation of ${O(N)}$\quad|textbf}
type X
 or 
\index{B}{typeY@type Y, representation of ${O(N)}$\quad|textbf}
type Y, 
 if the restriction $\sigma|_{SO(N)}$ to the special orthogonal group $SO(N)$ 
 is irreducible
 or reducible, 
 respectively.  
\end{definition}
With the convention as in Definition \ref{def:type},
 we recall a classical fact
 about the branching rule for the restriction
 $O(N) \downarrow SO(N)$. 
\begin{lemma}
[$O(N) \downarrow SO(N)$]
\label{lem:OSO}
Let $\sigma=(\sigma_1, \cdots, \sigma_N) \in \Lambda^+(O(N))$, 
 and $k$ $(\le [\frac N 2])$ be as in \eqref{eqn:TypeI} and \eqref{eqn:TypeII}.  
\begin{enumerate}
\item[{\rm{(1)}}]
{\rm{(type X)}}\enspace
The restriction of the irreducible $O(N)$-module 
 $\Kirredrep {O(N)}{\sigma}$
 to  $SO(N)$ is irreducible if and only if $N \ne2k$. 
In this case, 
the restricted $SO(N)$-module has highest weight
 $(\sigma_1, \cdots, \sigma_k, 0, \cdots, 0)$.  
\item[{\rm{(2)}}]
{\rm{(type Y)}}\enspace
If $N=2k$, 
 then the restriction $\Kirredrep {O(N)}{\lambda}|_{SO(N)}$ splits into two inequivalent 
 irreducible representations
 of $SO(N)$
 with highest weights
 $(\sigma_1, \cdots, \sigma_{k-1}, \sigma_k)$
 and $(\sigma_1, \cdots, \sigma_{k-1}, -\sigma_k)$.  
\end{enumerate}
\end{lemma}

\begin{example}
\label{ex:2.1}
The orthogonal group $O(N)$ acts irreducibly
 on the $\ell$-th exterior tensor space
\index{A}{0tExteriorCN@$\Exterior^{\ell}({\mathbb{C}}^N)$, exterior tensor|textbf}
 $\Exterior^{\ell}({\mathbb{C}}^N)$
 and on the space 
\index{A}{HSl@${\mathcal{H}}^s({\mathbb{C}}^N)$, spherical harmonics|textbf}
${\mathcal{H}}^s({\mathbb{C}}^N)$
of spherical harmonics of degree $s$.  
Via the parametrization \eqref{eqn:CWOn}, 
 these representations are described as follows:
\begin{alignat*}{2}
\Exterior^\ell ({\mathbb{C}}^N) &= \Kirredrep {O(N)}{1^{\ell}}
\qquad
&&(0\leq \ell \leq N), 
\\
\mathcal{H}^s({\mathbb{C}}^N)&=\Kirredrep {O(N)}{s,0,\cdots,0}
\qquad
&&(s \in {\mathbb{N}}).  
\end{alignat*}
The $O(N)$-module $\Exterior^{\ell}({\mathbb{C}}^N)$ is of type Y
 if and only if $N=2\ell$;
 the $O(N)$-module ${\mathcal{H}}^s({\mathbb{C}}^N)$ is of type Y
 if and only if $N=2$
 and $s\ne 0$.  
\end{example}

Irreducible $O(N)$-modules of types I and II are related by the following $O(N)$-isomorphism:
\begin{equation}
\label{eqn:type1to2}
   \Kirredrep {O(N)}{a_1,\cdots,a_k,1,\cdots,1,0,\cdots,0}
  =
   \det \otimes
   \Kirredrep {O(N)}{a_1,\cdots,a_k,0,\cdots,0}.
\end{equation}

Hence we obtain the following:
\begin{lemma}
\label{lem:typeY}
Let $\sigma \in \widehat{O(N)}$.  
Then $\sigma$ is of type Y
 if and only if $\sigma \otimes \det \simeq \sigma$.  
\end{lemma}

Then the following proposition is clear.  
\begin{proposition}
\label{prop:XYI}
Suppose $\sigma \in \widehat{O(n)}$.  
\begin{enumerate}
\item[{\rm{(1)}}]
If $\sigma$ is of type Y, 
 then $\sigma$ is of type I. 
\item[{\rm{(2)}}]
If $\sigma$ is of type II, 
 then $\sigma$ is of type X. 
\end{enumerate}
\end{proposition}

\subsubsection{Branching laws $O(N) \downarrow O(N-1)$}
\label{subsec:ONbranch}

Next we recall the classical branching laws for $O(N) \downarrow O(N-1)$.  
Let $\sigma=(\sigma_1, \cdots, \sigma_N) \in \Lambda^+(O(N))$
 and $\tau=(\tau_1, \cdots, \tau_{N-1}) \in \Lambda^+(O(N-1))$.  

\begin{definition}
\label{def:Young}
We denote by 
\index{A}{1pataorder@$\tau \prec \sigma$|textbf}
$\tau \prec \sigma$
 if  
\[
   \sigma_1 \ge \tau_1 \ge \sigma_2 \ge \tau_2 \ge \cdots \ge \tau_{N-1} \ge \sigma_N.  
\]
\end{definition}
Then the irreducible decomposition
 of representations of $O(N)$ with respect to the subgroup $O(N-1)$
 is given as follows:
\begin{fact}
[Branching rule for orthogonal groups]
\label{fact:ONbranch}
\index{B}{branching rule@branching rule, for $O(N) \downarrow O(N-1)$|textbf}
Let $(\sigma_1, \cdots,\sigma_N) \in \Lambda^+(O(N))$.  
Then the irreducible representation
 $\Kirredrep{O(N)}{\sigma_1, \cdots,\sigma_N}$ decomposes
 into a multiplicity-free sum
 of irreducible representations
 of $O(N-1)$ as follows:
\begin{equation}
\label{eqn:ONbranch}
   \Kirredrep{O(N)}{\sigma_1, \cdots, \sigma_N}|_{O(N-1)}
   =
  \bigoplus_{\tau \prec \sigma} 
   \Kirredrep{O(N-1)}{\tau_1, \cdots, \tau_{N-1}}.  
\end{equation}
\end{fact}

The commutant $O(1)$ of $O(N-1)$ in $O(N)$ acts
 on the irreducible summand
 $\Kirredrep {O(N-1)}{\tau_1, \cdots, \tau_{N-1}}$
 by 
$
   ({\operatorname{sgn}})^{\sum_{j=1}^N \sigma_j - \sum_{i=1}^{N-1} \tau_i}.
$

The following lemma is derived from Lemma \ref{lem:typeY} and Fact \ref{lem:branchII}.  
\begin{lemma}
\label{lem:branchII}
Let $\sigma \in \widehat {O(n)}$
 be of type I (see Definition \ref{def:type}).  
Then the following four conditions are equivalent:
\begin{enumerate}
\item[{\rm{(i)}}]
$\sigma \otimes \det \simeq \sigma$; 
\item[{\rm{(ii)}}]
$[\sigma|_{O(n-1)}:\tau]=[\sigma|_{O(n-1)}:\tau \otimes \det]$
 for all $\tau \in \widehat {O(n-1)}$; 
\item[{\rm{(iii)}}]
$n$ is even and $\sigma =\Kirredrep{O(n)}{s_1, \cdots, s_{\frac n 2}, 0, \cdots, 0}$
 with $s_{\frac n 2} \ne 0$;  
\item[{\rm{(iv)}}]
$\sigma|_{S O(n)}$ is reducible, 
 {\it{i.e.,}} $\sigma$ is of type Y
 (Definition \ref{def:OSO}).  
\end{enumerate}
\end{lemma}

\subsection{Principal series representations $I_{\delta}(V,\lambda)$
 of the orthogonal group $G=O(n+1,1)$}
\label{subsec:ps}

We discuss here (nonspherical) principal series representations
\index{A}{IdeltaV@$I_{\delta}(V, \lambda)$}
 $I_{\delta}(V,\lambda)$ of $G=O(n+1,1)$.  
We shall use the symbol 
\index{A}{JWepsilon@$J_{\varepsilon}(W, \nu)$}
 $J_{\varepsilon}(W,\nu)$
 for the principal series representations of the subgroup 
 $G'=O(n,1)$.

We recall the structure of  principal series representations for  rank one orthogonal groups.  
The main references are Borel--Wallach \cite{BW}
 and Collingwood \cite[Chap.~5, Sect.~2]{C}
 for the representations
 of the identity component $G_0=SO_0(n+1,1)$.  
We extend  here the results to the disconnected group.  
For representations of the disconnected group $G$, 
 see also \cite[Chap.~2]{sbon}
 for the spherical case
 ({\it{i.e.,}} $V={\bf{1}}$) and \cite[Chap.~2, Sect.~3]{KKP}
 for $V=\Exterior^i({\mathbb{C}}^n)$
 ($0 \le i \le n$).

\subsubsection{$C^{\infty}$-induced representations $I_\delta(V,\lambda)$}  
\label{subsec:smoothI}

We recall from Section \ref{subsec:subgpOn}
 that the Levi subgroup $M A$ of the minimal parabolic subgroup $P$
 of $G$ is a direct product group 
 $(O(n) \times O(1)) \times {\mathbb{R}}$.  
Then any irreducible representation
 of $M A$
 is the outer tensor product
 of irreducible representations
 of the three groups $O(n)$, $O(1)$, and ${\mathbb{R}}$.

One-dimensional representations $\delta$
 of
\index{A}{m2@$m_-={\operatorname{diag}}(-1,1,\cdots,1,-1)$}
 $O(1)=\{1,m_-\}$
 are labeled by $+$ or $-$, 
 where we write $\delta=+$ for the trivial representation ${\bf{1}}$,
 and $\delta=-$ for the nontrivial one
 given by $\delta(m_-)=-1$.  
Thus we identify $\widehat{O(1)}$ with the set $\{\pm\}$.

For $\lambda \in {\mathbb{C}}$, 
 we denote by 
\index{A}{charA@${\mathbb{C}}_{\lambda}$, character of $A$|textbf}
${\mathbb{C}}_{\lambda}$
 the one-dimensional representation
 of the split group $A$
 normalized 
 by the generator $H \in {\mathfrak{a}}$
 (see \eqref{eqn:Hyp}) as 
\index{A}{H@$H$}
\[
   A \mapsto {\mathbb{C}}^{\times}, 
  \qquad
  \exp(t H) \mapsto e^{\lambda t}.  
\]

Let $(\sigma, V)$ be an irreducible representation of $O(n)$, 
 $\delta \in \{\pm\}$, 
 and $\lambda \in{\mathbb{C}}$.  
We extend the outer tensor product representation
\index{A}{Vln@$V_{\lambda,\delta}=V \otimes \delta \otimes {\mathbb{C}}_{\lambda}$, representation of $P$|textbf}
\begin{equation}
\label{eqn:Vlmddelta}
 V_{\lambda,\delta}
:=V \boxtimes \delta \boxtimes{\mathbb{C}}_{\lambda}
\end{equation}
of the direct product group $MA \simeq O(n) \times O(1) \times {\mathbb{R}}
$
 to a representation
 of the parabolic subgroup $P=MAN_+$
 by letting the unipotent subgroup $N_+$ act trivially.  
The resulting irreducible $P$-module will be written 
as 
$
   V_{\lambda,\delta}
   =V \otimes \delta \otimes{\mathbb{C}}_{\lambda}
$
 by a little abuse of notation.  
We define  the induced representation of $G$ by 
\[
   I_\delta({V},\lambda)
  \equiv I({V}\otimes \delta, \lambda)
  := {\operatorname{Ind}}_P^G(V_{\lambda,\delta}).
\]
We refer to $\delta$ as the 
\index{B}{signature of the induced representation}
{\it{signature}} of the induced representation. 
If $\delta=+$ (the trivial character ${\bf{1}}$),
 we sometimes suppress the subscript.

If $(\sigma, V) \in \widehat{O(n)}$ is given 
 as $V=\Kirredrep{O(n)}{\sigma_1, \cdots,\sigma_n}$
with $(\sigma_1,\cdots,\sigma_n) \in \Lambda^+(O(n))$
 via \eqref{eqn:CWOn}, 
 then $I_{\delta}(V,\lambda)$ has 
a 
\index{B}{infinitesimalcharacter@infinitesimal character}
${\mathfrak {Z}}_G({\mathfrak{g}})$-infinitesimal character
\begin{equation}
\label{eqn:ZGinfI}
 (\sigma_1+\frac {n} 2-1, \sigma_2+\frac {n} 2-2, \cdots, 
  \sigma_k+\frac {n} 2-k, \frac {n} 2-k-1, \cdots, \frac {n} 2-[\frac n 2], 
  \lambda-\frac {n} 2)
\end{equation}
in the standard coordinates
 via the 
\index{B}{HarishChandraisomorphism@Harish-Chandra isomorphism}
 Harish-Chandra isomorphism, 
 see \eqref{eqn:HCpara}. 
We are using in this article
 unnormalized induction, 
{\it{i.e.}}, 
 the representation $I_{\delta}({V}, \frac n 2)$ is a unitarily induced principal series representation. 
Thus if $\lambda$ is purely imaginary, 
 the principal series representations
 $I_\delta(V,\lambda + \frac n 2)$ are 
\index{B}{temperedrep@tempered representation}
tempered.  
If $n$ is even, 
 then every irreducible tempered representation is isomorphic
 to a tempered principal series representation.  
If $n$ is odd, 
 then there is one family of discrete series representations
 parametrized by characters
 of the compact Cartan subgroup
 and every irreducible tempered representation
 is isomorphic to a tempered principal series representation
 or a discrete series representation.

\medskip

We denote by
\index{A}{Vlnd@${\mathcal{V}}_{\lambda,\delta}$, homogeneous vector bundle over $G/P$|textbf}
\begin{equation}
\label{eqn:Vlmdbdle}
  {\mathcal{V}}_{\lambda,\delta}
  :=
  G \times_P V_{\lambda,\delta}
\end{equation}
 the $G$-equivariant vector bundle over the real flag manifold $G/P$
 associated to the representation 
 $V_{\lambda,\delta}$ of $P$.  
We assume from now on that the principal series representations
 $I_\delta({V},\lambda)$
 are realized  on the Fr{\'e}chet space
 $C^{\infty}(G/P, {\mathcal{V}}_{\lambda,\delta})$
 of smooth sections
 for the vector bundle 
 ${\mathcal{V}}_{\lambda,\delta} \to G/P$.  
Thus 
\index{A}{IdeltaV@$I_{\delta}(V, \lambda)$|textbf}
$I_{\delta}(V,\lambda)$
 is the induced representation 
 $C^{\infty}{\operatorname{-Ind}}_P^G(V_{\lambda, \delta})$
 which is of moderate growth, 
see \cite[Chap.~3, Sect.~4]{sbon}.  
As usual, 
 we denote the representation space and the representation by the same letter. 
We trivialize the vector bundle ${\mathcal{V}}_{\lambda,\delta}$
 over $G/P$
 on the open Bruhat cell via the following map
\[
   \iota_N \colon {\mathbb{R}}^n {\underset {n_-}{\overset \sim \to}} N_-
   \overset \sim \to N_- \cdot o \subset G/P.  
\]
Then $I_{\delta}(V,\lambda)$ is realized in a subspace
 of $C^{\infty}({\mathbb{R}}^n) \otimes V$ by 
\begin{equation}
\label{eqn:iotaN}
   \iota_N^{\ast} \colon 
   I_{\delta}(V,\lambda)
   \hookrightarrow
   C^{\infty}({\mathbb{R}}^n) \otimes V,
   \quad
   F \mapsto f(b):=F(n_-(b)), 
\end{equation}
and this model is referred to as
 the {\it{noncompact picture}}, 
 or the 
\index{B}{Npicture@$N$-picture|textbf}
{\it{$N$-picture}}, 
 see Section \ref{subsec:psKN}.

\subsubsection{Tensoring with characters
 $\chi_{\pm\pm}$ of $G$}
The character group $(G/G_0)\hspace{-1mm}{\widehat{\hphantom{m}}}$
 of the component group
 $G/G_0 \simeq {\mathbb{Z}}/2{\mathbb{Z}} \times {\mathbb{Z}}/2{\mathbb{Z}}$
 acts on the set of admissible representations $\Pi$ of $G$, 
 by taking the tensor product
\begin{equation}
\label{eqn:Pichi}
   \Pi \mapsto \Pi \otimes \chi
\end{equation}
for $\chi \in (G/G_0)\hspace{-1mm}{\widehat{\hphantom{m}}}$.  
This action leaves the subsets 
\index{A}{Irredrep@${\operatorname{Irr}}(G)$}
 ${\operatorname{Irr}}(G)$ and ${\operatorname{Irr}}(G)_{\rho}$ 
(see Definition \ref{def:Irrrho})
invariant.  
We describe the action explicitly
 on principal series representations
 in Lemma \ref{lem:IVchi} below.  
The action on 
\index{A}{IrrGrho@${\operatorname{Irr}}(G)_{\rho}$,
 set of irreducible admissible smooth representation of $G$
 with trivial infinitesimal character $\rho$\quad}
${\operatorname{Irr}}(G)_{\rho}$
  will be given explicitly 
 in Theorem \ref{thm:LNM20} (5), 
 and on the space of symmetry breaking operators
 in Section \ref{subsec:actPont}.  

\begin{lemma}
\label{lem:IVchi}
Let $V \in \widehat {O(n)}$, 
 $\delta \in \{\pm\}$, 
 and $\lambda \in {\mathbb{C}}$.  
Let $\chi_{\pm\pm}$ be the one-dimensional representations
 of $G=O(n+1,1)$
 as defined in \eqref{eqn:chiab}.  
Then we have the following isomorphisms
 between representations of $G$:  
\begin{align*}
I_{\delta}(V,\lambda) \otimes \chi_{+-} \simeq & I_{-\delta}(V,\lambda), 
\\
I_{\delta}(V,\lambda) \otimes \chi_{-+} \simeq & I_{-\delta}(V \otimes \det,\lambda), 
\\
I_{\delta}(V,\lambda) \otimes \chi_{--} \simeq & I_{\delta}(V\otimes \det,\lambda).  
\end{align*}
\end{lemma}

\begin{proof}
For any $P$-module $U$
 and for any finite-dimensional $G$-module $F$, 
 there is an isomorphism of $G$-modules:
\[
  F \otimes {\operatorname{Ind}}_P^G(U)
  \simeq
  {\operatorname{Ind}}_P^G(F \otimes U).  
\]
Then Lemma \ref{lem:IVchi} follows from the restriction formula
 of the character $\chi$ of $G$ to the subgroup $M \simeq O(n) \times O(1)$ as below:
\[
   \chi_{+-}|_M \simeq {\bf{1}} \boxtimes {\operatorname{sgn}}, 
\quad   
   \chi_{-+}|_M \simeq \det \boxtimes {\operatorname{sgn}}, 
\quad
   \chi_{--}|_M \simeq \det \boxtimes {\bf{1}}.  
\]
\end{proof}

A special case of Lemma \ref{lem:IVchi} 
 for the exterior tensor representations $V=\Exterior^i({\mathbb{C}}^n)$
 will be stated in Lemma \ref{lem:LNM27}.

\subsubsection{$K$-structure of the principal series representation
 $I_{\delta}(V,\lambda)$}
\label{subsec:KstrI}
Let $(\sigma,V) \in \widehat {O(n)}$
 and $\delta \in \{\pm\}$ as before.  
By the Frobenius reciprocity law, 
 $K$-types of the principal series representation $I_\delta({V},\lambda)$
 are the irreducible representations
 of $K=O(n+1) \times  O(1)$
 whose restriction to $M \simeq O(n)\times O(1)$ contains the representation
 $V\boxtimes \delta$ of $M$.  
The classical branching theorem (Fact \ref{fact:ONbranch}) is used
 to determine $K$-types of the $G$-module $I_{\delta}(V,\lambda)$.  
We shall give an explicit $K$-type formula
 in the next section
 when $V$ is the exterior tensor representation $\Exterior^i({\mathbb{C}}^n)$
 of $O(n)$.  
For the general representation $(\sigma,V) \in \widehat {O(n)}$, 
 we do not use an explicit $K$-type formula of $I_{\delta}(V,\lambda)$, 
 but just mention an immediate corollary of Fact \ref{fact:ONbranch}:
\begin{proposition}
The $K$-types of principal series representations $I_{\delta}(V,\lambda)$
 of $O(n+1,1)$ have multiplicity one.  
\end{proposition}

\subsection{Principal series representations $I_{\delta}(i, \lambda)$}
\label{subsec:irred}
For $0 \le i \le n$, 
 $\delta \in \{\pm\}$, 
 and $\lambda \in {\mathbb{C}}$, 
 we denote the principal series representation 
 $I_\delta (\Exterior^i({\mathbb{C}}^n),\lambda)
 =
 C^{\infty}{\text{-}}{\operatorname{Ind}}_P^G(\Exterior^i({\mathbb{C}}^n)\otimes \delta \otimes {\mathbb{C}}_{\lambda})$
 of $G=O(n+1,1)$ simply by 
\index{A}{Ideltai@$I_{\delta}(i, \lambda)$|textbf}
$I_\delta(i,\lambda)$.  
Similarly,
 we write 
\index{A}{Jjepsilon@$J_{\varepsilon}(j, \nu)$|textbf}
$J_{\varepsilon}(j,\nu)$
 for the induced representation
 $C^{\infty}{\operatorname{-Ind}}_{P'}^{G'}(\Exterior^j({\mathbb{C}}^{n-1}) \otimes {\varepsilon} \otimes {\mathbb{C}}_{\nu})$
 of $G'=O(n,1)$ 
 for $0 \le j \le n-1$, 
 $\varepsilon \in \{\pm\}$, 
 and $\nu \in {\mathbb{C}}$. 
In the major part of this monograph, 
 we focus our attention on special families of principal series representations
 $I_\delta(i,\lambda)$ of $G$
 and $J_{\varepsilon}(j,\nu)$ of the subgroup $G'$.

In geometry, 
 $I_{\delta}(i,\lambda)$ is a family of representations
 of the conformal group $O(n+1,1)$ 
 of $S^n$ on the space ${\mathcal{E}}^i(S^n)$
 of differential forms
 ({\it{cf.}} \cite[Chap.~2, Sect.~2]{KKP})
 on one hand.  
In representation theory,
 any irreducible, unitarizable representations
 with nonzero $({\mathfrak{g}}, K)$-cohomologies arise
 as subquotients
 in $I_{\delta}(i,\lambda)$ with $\lambda=i$
 for some $0 \le i \le n$
 and $\delta=(-1)^i$, 
 see Theorem \ref{thm:LNM20} (9), 
 also Proposition \ref{prop:gKq} in Appendix I.

In this section we collect some basic properties
 of the principal series representations
\[
  I_{\delta}(i,\lambda)
\qquad
  \text{for $\delta \in \{\pm\}$, $0 \le i \le n$, $\lambda \in {\mathbb{C}}$, }
\]
which will be used throughout the article.  

\subsubsection{${\mathfrak{Z}}_G({\mathfrak{g}})$-infinitesimal character of $I_{\delta}(i,\lambda)$}

As we have seen in \eqref{eqn:ZGinfI} in the general case, 
the 
\index{B}{infinitesimalcharacter@infinitesimal character}
 ${\mathfrak{Z}}_G({\mathfrak{g}})$-infinitesimal character
 of the principal series representation $I_{\delta}(i,\lambda)$ is given by 
\begin{alignat*}{2}
&(\underbrace{\frac n 2, \frac n 2-1, \cdots, \frac n 2-i+1}_i; 
  \underbrace{\frac n 2-i-1, \cdots, \frac n 2-[\frac n 2]}_{[\frac n 2]-i}; \lambda-\frac n 2)
\quad
&&\text{if }
0 \le i \le \frac n 2, 
\\
&(\underbrace{\frac n 2, \frac n 2-1, \cdots, -\frac n 2+i+1}_{n-i}; 
  \underbrace{-\frac n 2+i-1, \cdots, \frac n 2-[\frac n 2]}_{i-[\frac {n+1} 2]}; \lambda-\frac n 2)
\quad
&&\text{if }
\frac n 2 \le i \le n.  
\end{alignat*}

In particular,
 the $G$-module $I_{\delta}(i,\lambda)$ has the trivial infinitesimal character $\rho_G$
 if and only if $\lambda=i$ or $n-i$.

\subsubsection{$K$-type formula of the principal series representations $I_\delta(i,\lambda)$}

By the Frobenius reciprocity,
 we can compute the $K$-type formula of $I_{\delta}(i,\lambda)$
 explicitly by using the classical branching law 
 (Fact \ref{fact:ONbranch})
 and Example \ref{ex:2.1} as follows:  
\begin{lemma}
[$K$-type formula of $I_{\delta}(i,\lambda)$]
\label{lem:KtypeIi}
Let $0 \le i \le n$ and $\delta \in\{\pm\}$.  
With the parametrization \eqref{eqn:CWOn}, 
 the 
\index{B}{Ktypeformula@$K$-type formula|textbf}
$K$-type formula
 of the principal series representation $I_{\delta}(i,\lambda)$
 of $G=O(n+1,1)$ is described 
 as below:
\begin{enumerate}
\item[{\rm{(1)}}] 
 for $i=0$, 
\[
\bigoplus_{a=0}^{\infty} \Kirredrep {O(n+1)}{a,0^n}
  \boxtimes
  (-1)^a \delta;
\]
\item[{\rm{(2)}}]
for $1 \le i \le n-1$, 
\[
  \bigoplus_{a=1}^{\infty} \Kirredrep {O(n+1)}{a,1^i,0^{n-i}}
  \boxtimes
  (-1)^a \delta
  \oplus
  \bigoplus_{a=1}^{\infty} \Kirredrep {O(n+1)}{a,1^{i-1},0^{n+1-i}}
  \boxtimes
  (-1)^{a+1} \delta;
\]
\item[{\rm{(3)}}] 
for $i=n$, 
\[
  \bigoplus_{a=1}^{\infty} (\det \otimes \Kirredrep {O(n+1)}{a,0^n})
  \boxtimes
  (-1)^{a+1} \delta.  
\]
\end{enumerate}
\end{lemma}
See Proposition \ref{prop:KtypeIV} for a more general $K$-type formula
 of the principal series representation $I_{\delta}(V,\lambda)$.  

\subsubsection{Basic $K$-types of $I_{\delta}(i,\lambda)$}
\label{subsec:KIilmd}
\index{B}{basicKtype@basic $K$-type|textbf}
Let $\delta \in \{ \pm \}$
 and $0 \le i\le n$.  
Following the notation \cite[Chap.~2, Sect.~3]{KKP}, 
 we define two irreducible representations
 of $K\simeq O(n+1) \times O(1)$ by:
\index{A}{0muflat@$\mub(i,\delta)$|textbf}
\index{A}{0musharp@$\mus(i,\delta)$|textbf}
\begin{align}
   \mub(i,\delta):=& \Exterior^i({\mathbb{C}}^{n+1}) \boxtimes \delta, 
\label{eqn:muflat}
\\
   \mus(i,\delta):=& \Exterior^{i+1}({\mathbb{C}}^{n+1}) \boxtimes (-\delta).  
\label{eqn:musharp}
\end{align}
This means:
\[
 \begin{cases}
\mub(i,+)&=\Exterior^i({\mathbb{C}}^{n+1}) \boxtimes {\bf{1}},
\\
\mub(i,-)
&=\Exterior^i({\mathbb{C}}^{n+1})\boxtimes {\operatorname{sgn}},
   \end{cases}
\quad
\begin{cases}
\mus(i,+)
&=\Exterior^{i+1}({\mathbb{C}}^{n+1}) \boxtimes {\operatorname{sgn}}, 
\\
\mus(i,-)
&=\Exterior^{i+1}({\mathbb{C}}^{n+1}) \boxtimes {\bf{1}}.  
\end{cases}
\]
The superscripts $\sharp$ and $\flat$ indicate
 that there are the following obvious $K$-isomorphisms
\begin{equation}
\label{eqn:flatsharp}
   \mus(i,\delta)=\mub(i+1,-\delta)
   \qquad
   (0 \le i \le n), 
\end{equation}
which will be useful
 in describing the standard sequence
 with trivial infinitesimal character $\rho_G$
 (Definition \ref{def:Pii} below), 
 see also Remark \ref{rem:flatsharp}.

By the $K$-type formula of the principal series representation 
 $I_{\delta}(i,\lambda)$ in Lemma \ref{lem:KtypeIi}, 
 the $K$-types $\mub(i,\delta)$ and $\mus(i,\delta)$
 occur in $I_{\delta}(i,\lambda)$
 with multiplicity one for any $\lambda \in {\mathbb{C}}$.  

\begin{definition}
\label{def:basicK}
We say $\mub (i,\delta)$ and $\mus (i,\delta)$
 are {\it{basic $K$-types}}
 of the principal series representations $I_{\delta}(i,\lambda)$
 of $G=O(n+1,1)$.  
\end{definition}

\subsubsection{Reducibility of $I_{\delta}(i,\lambda)$}
\label{subsec:irredIilmd}
The principal series representation $I_{\delta}(i,\lambda)$
 is generically irreducible.  
More precisely,
 we have the following.
\begin{proposition}
\label{prop:redIilmd}
Let $G=O(n+1,1)$, 
 $0 \le i \le n$, 
 $\delta \in \{\pm\}$, 
 and $\lambda \in {\mathbb{C}}$.  
\begin{enumerate}
\item[{\rm{(1)}}]
The principal series representation $I_{\delta}(i,\lambda)$
 is reducible 
 if and only if 
\begin{equation}
\label{eqn:redIilmd}
 \lambda \in \{i,n-i\} \cup (-{\mathbb{N}}_+) \cup (n+{\mathbb{N}}_+).  
\end{equation}
\item[{\rm{(2)}}]
Suppose $(n,\lambda) \ne (2i,i)$.  
If $\lambda$ satisfies \eqref{eqn:redIilmd}, 
 then the $G$-module $I_{\delta}(i,\lambda)$ 
 has a unique irreducible proper submodule
 (say, A)
 and has a unique irreducible subquotient 
 (say, B)
 and there is a nonsplitting exact sequence
 of $G$-modules:
\[
   0 \to A \to I_{\delta}(i,\lambda) \to B \to 0.  
\]
\item[{\rm{(3)}}]
Suppose $(n,\lambda) = (2i,i)$.  
Then the $I_{\delta}(i,\lambda)$ decomposes
 into the direct sum 
 of two irreducible representations of $G$
 which are not isomorphic to each other.  
\end{enumerate}
\end{proposition}
When $n \ne 2 i$, 
 the \lq\lq{only if}\rq\rq\ part of the first statement and the second one
 in Proposition \ref{prop:redIilmd} follow readily from 
 the corresponding results 
 (\cite{BW,C,Hirai62})
 for the connected group $SO_0(n+1,1)$ 
 and from Lemma \ref{lem:GPconn} below
 because $\Exterior^i({\mathbb{C}}^n)$ is irreducible
 as an $SO(n)$-module.  
We need some argument for $n=2i$
 where $\Exterior^i({\mathbb{C}}^{n})$ is reducible
 as an $SO(n)$-module, 
 see Examples \ref{ex:irrIilmd} and \ref{ex:Imm}
 in Appendix II
 for the proof of Proposition \ref{prop:redIilmd} (1) and (3), 
 respectively.
In Section \ref{subsec:KerKnSt}, 
 we discuss the description of proper submodules
 of reducible $I_{\delta}(i,\lambda)$ by using the Knapp--Stein operator
 \eqref{eqn:KSii}
 and its normalized one \eqref{eqn:Ttilde}.  
The \lq\lq{if}\rq\rq\ part of the first statement is proved there, 
 see Lemma \ref{lem:reducibleI}.

The composition series of $I_{\delta}(i,\lambda)$ 
with trivial infinitesimal character $\rho_G$
 ({\it{i.e.}}, for $\lambda=i$ or $n-i$)
 will be discussed
 in the next subsection 
 (see Theorem \ref{thm:LNM20}), 
 which will be extended in Theorem \ref{thm:1714107}
 to the case of
\index{B}{regularintegralinfinitesimalcharacter@regular integral infinitesimal character}
 regular integral infinitesimal characters.  

\subsubsection{Irreducible subquotients of $I_{\delta}(i,i)$}
\label{subsec:subIii}
Every irreducible representation of $G=O(n+1,1)$ with trivial infinitesimal character $\rho $ is equivalent to a subquotient of $I_{\delta}(i,i)$
 for some $0 \le i \le n$
 and $\delta\in \{\pm\}$, 
 or equivalently,
 of $I_+(i,i) \otimes \chi $ with $i \geq n/2$ and $\chi \in (G/G_0)\hspace{-1mm}{\widehat{\hphantom{m}}}$.
We recall now facts about the principal series representations
 $I_+(i,i)$, $I_-(i,i)$, $I_+(n-i,i)$ and $I_-(n-i,i)$
 of the orthogonal group $O(n+1,1)$
 and their composition factors.

We denote by 
$
\index{A}{Ideltaiflat@$I_{\delta}(i)^{\flat}$, submodule of $I_{\delta}(i,i)$|textbf}
I_{\delta}(i)^{\flat}
$ 
and 
$
\index{A}{Ideltaisharp@$I_{\delta}(i)^{\sharp}$, quotient of $I_{\delta}(i,i)$|textbf}
I_{\delta}(i)^{\sharp}
$
 the unique irreducible subquotients
 of $I_{\delta}(i, i)$
 containing the basic $K$-types
 $\mu^{\flat}(i,\delta)$ and $\mu^{\sharp}(i,\delta)$, 
 respectively.  
Then we have $G$-isomorphisms:
\begin{equation}
\label{eqn:LNS218}
  I_{\delta}(i)^{\sharp} \simeq I_{-\delta}(i+1)^{\flat}
\quad
  \text{for $0 \le i \le n-1$ and $\delta \in \{\pm \}$}, 
\end{equation}
see Theorem \ref{thm:LNM20} (1) below.  
For $0 \le \ell \le n+1$ and $\delta \in \{\pm \}$, 
 we set 
\index{A}{1Piidelta@$\Pi_{i,\delta}$, irreducible representations of $G$|textbf}
\begin{equation}
\label{eqn:Pild}
  \Pi_{\ell,\delta}
  :=
  \begin{cases}
  I_{\delta}(\ell)^{\flat} \quad&(0 \le \ell \le n), 
  \\
  I_{-\delta}(\ell-1)^{\sharp} \quad&(1 \le \ell \le n+1).  
  \end{cases}
\end{equation}
In view of \eqref{eqn:LNS218}, 
 the irreducible representation $\Pi_{\ell,\delta}$ of $G$ is well-defined.  

\begin{remark}
\label{rem:flatsharp}
The point here is 
 that each irreducible representation
 $\Pi_{\ell,\delta}$ ($1 \le \ell \le n$, $\delta=\pm$)
 can be realized in two different principal series representations:
\begin{align*}
   I_{\delta}(\ell,\ell) =& {\operatorname{Ind}}_P^G (\Exterior^\ell({\mathbb{C}}^n) \otimes \delta \otimes {\mathbb{C}}_\ell), 
\\
I_{-\delta}(\ell-1,\ell-1) 
=& {\operatorname{Ind}}_P^G (\Exterior^{\ell-1}({\mathbb{C}}^n) \otimes (-\delta) \otimes {\mathbb{C}}_{\ell-1}).  
\end{align*}
\end{remark}

\begin{theorem}
\label{thm:LNM20}
Let $G=O(n+1,1)$ $(n \ge 1)$.  
\begin{enumerate}
\item[{\rm{(1)}}]
For $0 \le \ell \le n$ and $\delta \in \{ \pm \}$, 
 we have exact sequences of $G$-modules:
\begin{align*}
&0 \to \Pi_{\ell,\delta} \to I_{\delta}(\ell,\ell) \to \Pi_{\ell+1,-\delta}
\to 0, 
\\
&0 \to \Pi_{\ell+1,-\delta} \to I_{\delta}(\ell,n-\ell) \to \Pi_{\ell,\delta} \to 0.  
\end{align*}
These exact sequences split
 if and only if $n=2\ell$.  

\item[{\rm{(2)}}]
Irreducible admissible smooth representations of $G$ 
 with trivial
\index{B}{infinitesimalcharacter@infinitesimal character}
 ${\mathfrak {Z}}_G({\mathfrak {g}})$-infinitesimal character
\index{A}{1parho@$\rho_G$}
 $\rho_G$
 can be classified as
\index{A}{IrrGrho@${\operatorname{Irr}}(G)_{\rho}$,
 set of irreducible admissible smooth representation of $G$
 with trivial infinitesimal character $\rho$\quad}
\[
  {\operatorname{Irr}}(G)_{\rho}
  =\{
    \Pi_{\ell, \delta}
    :
    0 \le \ell \le n+1, \, \delta = \pm
\}.  
\]

\item[{\rm{(3)}}]
For any $0 \le \ell \le n+1$
 and $\delta \in \{\pm\}$, 
 the 
\index{B}{minimalKtype@minimal $K$-type}
minimal $K$-type of the irreducible $G$-module $\Pi_{\ell,\delta}$
 is given 
by $\mub(\ell,\delta)=\Exterior^\ell({\mathbb{C}}^{n+1}) \boxtimes \delta$.  

\item[{\rm{(4)}}]
There are four one-dimensional representations of $G$, 
 and they are given by 
\begin{equation*}
\{
  \Pi_{0,+} \simeq {\bf{1}}, \quad \Pi_{0,-} \simeq \chi_{+-}, 
\quad
   \Pi_{n+1,+} \simeq \chi_{-+}, \quad 
   \Pi_{n+1,-} \simeq \chi_{--} (=\det)
\}.  
\end{equation*}
The other representations $\Pi_{\ell,\delta}$
 $(1 \le \ell \le n, \delta \in \{\pm\})$
 are infinite-dimensional.  
\item[{\rm{(5)}}]
There are isomorphisms as $G$-modules
 for any $0 \le \ell \le n+1$ and $\delta = \pm$:
\index{A}{1chipmpm@$\chi_{\pm\pm}$, one-dimensional representation of $O(n+1,1)$}
\begin{align*}
   \Pi_{\ell,\delta} \otimes \chi_{+-} \simeq \, & \Pi_{\ell,-\delta}, 
\\
   \Pi_{\ell,\delta} \otimes \chi_{-+} \simeq \, & \Pi_{n+1-\ell,\delta}, 
\\   
   \Pi_{\ell,\delta} \otimes \chi_{--} \simeq \, & \Pi_{n+1-\ell,-\delta}.  
\end{align*}
\item[{\rm{(6)}}]
Every $\Pi_{\ell,\delta}$ $(0 \le \ell \le n+1, \delta = \pm)$
 is unitarizable and self-dual.  
\item[{\rm{(7)}}]
For $n$ odd, 
 there are exactly two inequivalent 
\index{B}{discreteseries@discrete series representation}
discrete series representations
 of $G=O(n+1,1)$ with infinitesimal character $\rho_G$.  
Their smooth representations
 are given by 
\[
  \{
   \Pi_{\frac {n+1}{2},\delta}
   :
    \delta = \pm
  \}.  
\]

\par
All the other representations in the list (2)
 are nontempered representations of $G$.  
\item[{\rm{(8)}}]
For $n$ even, 
 there are exactly four inequivalent irreducible
\index{B}{temperedrep@tempered representation}
 tempered representations of $G=O(n+1,1)$
 with infinitesimal character $\rho_G$.  
Their smooth representations are given by 
\[
  \{
   \Pi_{\frac {n}{2},\delta}, \Pi_{\frac {n}{2}+1,\delta}
   :
    \delta = \pm
  \}.  
\]
\item[{\rm{(9)}}]
Irreducible and unitarizable $({\mathfrak {g}}, K)$-modules
 with nonzero $({\mathfrak {g}}, K)$-cohomologies are
 exactly given as the set
 of the underlying $({\mathfrak {g}}, K)$-modules of $\Pi_{\ell, \delta}$
 $(0 \le \ell \le n+1, \delta = \pm)$.  
\end{enumerate}
\end{theorem}

The exact sequences in Theorem \ref{thm:LNM20} (1) leads us
 to a labeling of the finite set ${\operatorname{Irr}}(G)_{\rho}$
 as follows: 
\begin{definition}
[standard sequence]
\label{def:Pii}
Let $G=O(n+1,1)$
 and $n=2m$ or $2m-1$.  
We refer to the sequence 
\begin{center}
\begin{tabular}{ccccccccccc}
& &$\Pi_{0,+}$&,&$\Pi_{1,+} $ & ,  &\dots & , & $\Pi_{m-1,+} $ &, & $\Pi_{m,+}$ 
\end{tabular}
\end{center}
as the 
\index{B}{standardsequence@standard sequence|textbf}
{\it standard sequence starting
 with the trivial one-dimensional representation 
\index{A}{01one@${\bf{1}}$, trivial one-dimensional representation}
$\Pi_{0,+}={\bf{1}}$}.  
Likewise,
 we refer to the sequence
\begin{center}
\begin{tabular}{ccccccccccc}
& &$\Pi_{0,-}$&,&$\Pi_{1,-} $ & ,  &\dots & , & $\Pi_{m-1,-} $ &, & $\Pi_{m,-}$ 
\end{tabular}
\end{center}
as the standard sequence starting with the one-dimensional representation 
 $\Pi_{0,-}=\chi_{+-}$.  
Sometimes we suppress the subscript $+$
 and write 
\index{A}{1Pii@$\Pi_{i}=\Pi_{i,+}$|textbf}
 $\Pi_{i}$ for $\Pi_{i,+}$ for simplicity.  
\end{definition}

More generally,
 we shall define the standard sequence 
 starting with other irreducible finite-dimensional representations of $G$
 in Chapter \ref{sec:conjecture}, 
 see Definition \ref{def:Hasse} and Example \ref{ex:171455}. 
An analogous sequence, 
 which we refer to as the Hasse sequence,
 will be defined also in Chapter \ref{sec:conjecture}, 
 see Definition-Theorem \ref{def:UHasse}.   

\vskip 1pc
We give some remarks on the proof of Theorem \ref{thm:LNM20}.  
Basic references are \cite{BW, C, KKP}.  
Theorem \ref{thm:LNM20} (1) generalizes the results proved
 in Borel--Wallach 
 \cite[pp.~128--129 in the new edition; p.~192 in the old edition]{BW}
 for the identity component group 
 $G_0=SO_0(n+1,1)$.  
(Unfortunately and confusingly the restriction
 of our representations $I_+(i,i)$
 to the connected component $G_0$ are denoted there by $I_i$
 when $n \ne 2i$.) 
See also Collingwood \cite[Chap.~5, Sect.~2]{C} for 
 the identity component group $G_0$;
\cite[p.~20]{KKP}
 for the disconnected group $G=O(n+1,1)$.

For the relationship
 between principal series representations of $G$
 and of its identity component group $G_0$, 
 we recall from \cite[Chap.~5]{sbon} the following.  
\begin{lemma}
\label{lem:GPconn}
For $G=O(n+1,1)$, 
 let $P_0:=P \cap G_0$.  
Then $P_0$ is connected, 
 and is a minimal parabolic subgroup of $G_0$.  
Then we have a natural bijection:
\[
G_0/P_0 \overset \sim \to G/P \ 
(\simeq S^n).  
\]
\end{lemma}
Then we can derive results
 for the disconnected group $G$ from those for the connected group $G_0$
 and {\it{vice versa}}
 by using the action of the Pontrjagin dual $(G/G_0)\hspace{-1mm}{\widehat{\hphantom{m}}}$
 of the component group $G/G_0$
 and the classical branching law
 $O(N) \downarrow SO(N)$
 (Section \ref{subsec:OnSOn}).  
In Appendix II (Chapter \ref{sec:SOrest}) 
 we discuss restrictions
 of representations of $O(n+1,1)$
 with respect to $SO(n+1,1)$
 in the same spirit.

In Proposition \ref{prop:161655} of Appendix I, 
 we will give a description
 of the underlying $({\mathfrak{g}},K)$-modules $(\Pi_{i,\pm})_K$
 of the $G$-irreducible subquotients $\Pi_{i,\pm}$
 in terms of the so-called 
\index{A}{Aqlmd@$A_{\mathfrak{q}}(\lambda)$}
$A_{\mathfrak{q}}(\lambda)$-modules,
 {\it{i.e.,}} cohomologically induced representations from 
 one-dimensional representations
 of a $\theta$-stable parabolic subalgebra ${\mathfrak{q}}$.

By using the description,
 Theorem \ref{thm:LNM20} (9) follows readily from 
 results of Vogan and Zuckerman \cite{VZ}, 
 see Proposition \ref{prop:gKq} in Appendix I.  
The unitarizability of the irreducible subquotients
 $\Pi_{i,\pm}$
 (Theorem \ref{thm:LNM20} (6)) traces
 back to T.~Hirai \cite{Hirai62}, 
 see also Howe and Tan \cite{HT93}.  
Alternatively,
 the unitarizability in Theorem \ref{thm:LNM20} (6)
 is deduced from the theory
 on $A_{\mathfrak{q}}(\lambda)$, 
 see \cite[Thm.~0.51]{KV}.

\begin{remark}
Analogous results for the special orthogonal group $SO(n+1,1)$
 will be given in Proposition \ref{prop:161648}
 in Appendix II, 
 where we denote the group by $\overline G$.  
\end{remark}

\newpage
\section{Symmetry breaking operators 
 for principal series representations---general theory}
\label{sec:general}

In this chapter we discuss important concepts
 and properties of symmetry breaking operators from 
 principal series representations
 $I_{\delta}(V,\lambda)$ 
 of the orthogonal group $G=O(n+1,1)$
 to $J_{\varepsilon}(W,\nu)$ of the subgroup $G'=O(n,1)$.  
In particular,
 we present a classification scheme
 (Theorem \ref{thm:VWSBO})
 of all symmetry breaking operators,
 which is built on the strategy
 of the classification 
 in the spherical case
 \cite{sbon}
 and also on a new phenomenon 
 for which we refer to as {\it{sporadic operators}}
 (Section \ref{subsec:SBOVW2}).  
The classification scheme is carried out
 in full details
 for symmetry breaking from principal series representations
 $I_{\delta}(V,\lambda)$ of $G$
 to $J_{\varepsilon}(W,\nu)$ of the subgroup $G'$, 
 which will play a crucial role 
 in understanding symmetry breaking 
 of all {\it{irreducible}} admissible representations
 of $G$ 
 having the trivial infinitesimal character
 (Chapters \ref{sec:SBOrho}, \ref{sec:Gross-Prasad}, and \ref{sec:period}).  
Various theorems stated in this chapter 
 will be proved in later chapters,
 in particular,
 in Chapter \ref{sec:section7}.

\subsection{Generalities}
\label{subsec:SBOgen}
We refer to nontrivial homomorphisms
 in 
\[
  {\operatorname{Hom}}_{G'}
  (I_{\delta}(V,\lambda)|_{G'}, J_{\varepsilon}(W,\nu))
\]
as intertwining restriction operators
 or 
\index{B}{symmetrybreakingoperators@symmetry breaking operators|textbf}
{\it{symmetry breaking operators}}.  
Here $\delta, \varepsilon \in {\mathbb{Z}}/2 {\mathbb{Z}}$
 in our setting
 where $(G,G')=(O(n+1,1), O(n,1))$.  
For a detailed introduction to symmetry breaking operators
 we refer to \cite{xkvogan} and 
 \cite[Chaps.~1 and 3]{sbon}.

\subsection{Summary of results}
\label{subsec:summarySBO}
We keep our setting
 where $(G,G')=(O(n+1,1), O(n,1))$.

For $(\sigma, V) \in \widehat {O(n)}$, 
 $\delta \in \{\pm\}$, 
 and $\lambda \in {\mathbb{C}}$, 
 we write $I_{\delta}(V, \lambda)$
 for the principal series representation
 of $G$ as in Section \ref{subsec:ps}.  
Similarly, 
 let $(\tau, W)$ be an irreducible representation of $O(n-1)$, 
 $\varepsilon \in \{\pm \}$, 
 and $\nu \in{\mathbb{C}}$.  
We extend the outer tensor product representation
\index{A}{Wnuepsi@$W_{\nu,\varepsilon}=W \otimes \varepsilon \otimes {\mathbb{C}}_{\nu}$, representation of $P'$|textbf}
\[
 W_{\nu,\varepsilon}
:=W \boxtimes {\varepsilon} \boxtimes{\mathbb{C}}_{\nu}
\]
of the direct product group $M' A \simeq O(n-1) \times O(1) \times {\mathbb{R}}
$ to $P' = M' A N_+'$
 by letting $N_+'$ act trivially.  
We also write $W_{\nu,\varepsilon}=W \otimes \varepsilon \otimes {\mathbb{C}}_{\nu}$
 when we regard it as a $P'$-module.  
We form a $G'$-equivariant vector bundle
\index{A}{Wnuepsical@${\mathcal{W}}_{\nu,\varepsilon}$, homogeneous vector bundle over $G'/P'$|textbf}
$
  {\mathcal{W}}_{\nu,\varepsilon}
  :=
  G' \times_{P'} W_{\nu,\varepsilon}
$
 over the real flag manifold $G'/P'$.  
The principal series representation 
\index{A}{JWepsilon@$J_{\varepsilon}(W, \nu)$|textbf}
 $J_{\varepsilon}(W, \nu)$
 of $G'=O(n,1)$
 is defined to be the induced representation
 $\operatorname{Ind}_{P'}^{G'}(W_{\nu,\varepsilon})$
 on the space 
 $C^{\infty}(G'/P',{\mathcal{W}}_{\nu,\varepsilon})$
 of smooth sections for the vector bundle.

For $(\sigma, V) \in \widehat {O(n)}$
 and $(\tau, W) \in \widehat {O(n-1)}$, 
 we set 
\index{A}{VWmult@$[V:W]$|textbf}
\begin{equation}
\label{eqn:multVW}
    [V:W]:= \dim_{\mathbb{C}} \operatorname{Hom}_{O(n-1)}(V|_{O(n-1)}, W).  
\end{equation}
If we want to emphasize the subgroup,
 we also write $[V|_{O(n-1)}:W]$ for $[V:W]$.  
We recall from Fact \ref{fact:ONbranch}
 on the classical branching rule for the restriction
 $O(N)\downarrow O(N-1)$
 that the multiplicity $[V:W]$ is either $0$ or $1$.

\subsubsection{Symmetry breaking operators
 when $[V:W] \ne 0$}
\label{subsec:SBOVW1}
Suppose $[V:W] \ne 0$.  
In this case 
 we prove the existence of nonzero symmetry breaking operators
 for all $\lambda, \nu \in {\mathbb{C}}$
 and for all signatures $\delta$, $\varepsilon \in \{ \pm \}$:
\begin{theorem}
[existence of symmetry breaking operators, 
 see Theorem {\ref{thm:existSBO}}]
\label{thm:160150}
Suppose $(\sigma,V)\in \widehat{O(n)}$ and $(\tau,W)\in \widehat{O(n-1)}$.  
Assume $[V:W] \ne 0$.  
Then we have
\index{A}{IdeltaV@$I_{\delta}(V, \lambda)$}
\[
  \dim_{\mathbb{C}} \operatorname{Hom}_{G'}(I_{\delta}(V, \lambda)|_{G'}, J_{\varepsilon}(W, \nu))
  \ge 1
  \quad
  \text{for all }
  \delta, \varepsilon \in {\mathbb{Z}}/2 {\mathbb{Z}}, 
  \,
  \lambda, \nu \in {\mathbb{C}}.  
\]
\end{theorem}
Theorem \ref{thm:160150} is proved in Section \ref{subsec:existSBO}
 by constructing symmetry breaking operators:
 generic ones are nonlocal 
 ({\it{e.g.}} integral operators)
 see Theorem \ref{thm:regexist} below, 
 whereas a few are local operators
 ({\it{i.e.}} differential operators, 
 see Theorem \ref{thm:1716110}).

\begin{definition}
\label{def:generic}
We say that the quadruple $(\lambda,\nu,\delta,\varepsilon)$
 is a 
\index{B}{genericparametercondition@generic parameter condition|textbf}
 {\it{generic parameter}} if
$(\lambda, \nu)\in {\mathbb{C}}^2$ and $\delta, \varepsilon 
 \in \{\pm\}$ satisfy
\begin{equation}
\label{eqn:nlgen}
\begin{cases}
\nu-\lambda \not \in 2{\mathbb{N}}
\quad
&
\text{when}
\quad
\delta \varepsilon =+;
\\
\nu-\lambda \not \in 2{\mathbb{N}}+1
\quad
&
\text{when}
\quad
\delta \varepsilon =-.  
\end{cases}
\end{equation}
\end{definition}

We recall from \eqref{eqn:singset}
 that the set of \lq\lq{special parameters}\rq\rq\ is given 
 as the complement of \lq\lq{generic parameters}\rq\rq, 
 namely,
\begin{alignat}{2}
  \Psising
  =
  \left\{(\lambda,\nu, \delta, \varepsilon) \in {\mathbb{C}}^2 \times \{\pm\}^2
   :\,\,\right.
&   \nu-\lambda \in 2 {\mathbb{N}} 
&&\text{when $\delta \varepsilon =+$}
\notag
\\
   \text{ or }\,\,
&  \nu-\lambda \in 2 {\mathbb{N}}+1
\quad 
&&
\text{when $\delta \varepsilon =-$}
\left.\right\}.  
\label{eqn:Psisp}
\end{alignat}

In the case $[V:W]\ne 0$, 
 we also prove the following 
\lq\lq{generic multiplicity-one theorem}\rq\rq, 
which extends \cite[Thm.\ 1.1]{sbon}
 in the scalar case 
($V=W={\mathbb{C}}$).

\begin{theorem}
[generic multiplicity-one theorem]
\label{thm:unique}
Suppose $(\sigma, V) \in \widehat{O(n)}$, 
$(\tau,W) \in \widehat{O(n-1)}$
 with $[V:W] \ne 0$.  
If $(\lambda, \nu, \delta, \varepsilon) \in {\mathbb{C}}^2 \times \{\pm\}^2$
 satisfies the generic parameter condition, 
 namely,
 $(\lambda, \nu, \delta, \varepsilon) \not \in \Psising$, 
 then 
\index{A}{IdeltaV@$I_{\delta}(V, \lambda)$}
\index{A}{JWepsilon@$J_{\varepsilon}(W, \nu)$}
\[
\dim_{\mathbb{C}} \operatorname{Hom}_{G'}(I_{\delta}(V, \lambda)|_{G'}, J_{\varepsilon}(W, \nu))
  =1.  
\]
\end{theorem}

Theorem \ref{thm:unique} gives a stronger estimate 
 than what the existing general theory guarantees:
\begin{enumerate}
\item[$\bullet$]
the dimension $\le 1$
 if both $I_{\delta}(V, \lambda)$ and $J_{\varepsilon}(W, \nu)$ are 
irreducible \cite{SunZhu}, 
\item[$\bullet$]
the dimension is uniformly bounded 
 with respect to $\sigma$, $\tau$, 
 $\delta$, $\varepsilon$, 
 $\lambda$, $\nu$
 \cite{xKOfm}.  
\end{enumerate}
We note
 that $I_{\delta}(V,\lambda)$ or $J_{\varepsilon}(W,\nu)$
 can be reducible
 even if $(\lambda, \nu, \delta, \varepsilon) \in \Psising$.  
Theorem \ref{thm:unique} will be proved
 in a strengthened form by giving an explicit generator 
 (see Theorem \ref{thm:genbasis}
 in Section \ref{subsec:existSBO}).  

\subsubsection{Differential symmetry breaking operators
 when $[V:W] \ne 0$}
\label{subsec:DSBOVW}

We realize the principal series representations
 $I_{\delta}(V,\lambda)$ and $J_{\varepsilon}(W,\nu)$
 in the Fr{\'e}chet spaces
 $C^{\infty}(G/P, {\mathcal{V}}_{\lambda,\delta})$
 and 
 $C^{\infty}(G'/P', {\mathcal{W}}_{\nu,\varepsilon})$.  

\begin{definition}
[differential symmetry breaking operator]
\label{def:DSBO}
A linear map 
\[
   D \colon C^{\infty}(G/P, {\mathcal{V}}_{\lambda,\delta})
     \to 
            C^{\infty}(G'/P', {\mathcal{W}}_{\nu,\varepsilon})
\]
 is called a 
\index{B}{differentialsymmetrybreakingoperators@differential symmetry breaking operator|textbf}
{\it{differential symmetry breaking operator}}
 if $D$ is a differential operator
 with respect to the inclusion $G'/P' \hookrightarrow G/P$
 and $D$ intertwines the action of the subgroup $G'$.  
See Definition \ref{def:diff}
 in Chapter \ref{sec:DSVO}
 for the notion of differential operators
 between two different manifolds.  
We denote by 
\index{A}{DiffG@${\operatorname{Diff}}_{G'}(I_{\delta}(V,\lambda)\vert_{G'}, J_{\varepsilon}(W,\nu))$|textbf}
\[
{\operatorname{Diff}}_{G'}
   (I_{\delta}(V,\lambda)|_{G'}, J_{\varepsilon}(W,\nu))
\]
 the subspace of 
 ${\operatorname{Hom}}_{G'}
   (I_{\delta}(V,\lambda)|_{G'}, J_{\varepsilon}(W,\nu))$
 consisting of differential symmetry breaking operators.  
\end{definition}

We retain the assumption 
 that $[V:W] \ne 0$.  
We give a necessary and sufficient condition
 for the existence of nonzero differential symmetry breaking operators:
\begin{theorem}
[existence of differential symmetry breaking operators]
\label{thm:1716110}
Suppose $(\sigma, V) \in \widehat{O(n)}$
 and $(\tau,W) \in \widehat{O(n-1)}$
 satisfy $[V:W] \ne 0$. 
Then the following two conditions
 on the parameters $\lambda, \nu \in {\mathbb{C}}$
 and $\delta, \varepsilon \in \{\pm\}$
 are equivalent:
\begin{enumerate}
\item[{\rm{(i)}}]
The quadruple $(\lambda, \nu, \delta, \varepsilon)$ does not satisfy
 the generic parameter condition \eqref{eqn:nlgen}, 
 namely,
 $(\lambda, \nu, \delta, \varepsilon)\in \Psising$.  
\item[{\rm{(ii)}}]
${\operatorname{Diff}}_{G'}
   (I_{\delta}(V,\lambda)|_{G'}, J_{\varepsilon}(W,\nu)) \ne \{0\}.  
$
\end{enumerate}
\end{theorem}

We shall prove Theorem \ref{thm:1716110}
 in Chapter \ref{sec:DSVO},
 see Theorem \ref{thm:existDSBO}.  

\vskip 1pc
\subsubsection{Sporadic symmetry breaking operators
 when $[V:W]=0$}
\label{subsec:SBOVW2}
This section treats the case $[V:W]=0$.  
In the holomorphic setting,
 we found in \cite{KP1} a phenomenon 
 that all symmetry breaking operators
 are given by differential operators
 ({\it{localness theorem}}).  
This phenomenon does not occur 
 in the real setting
 if both $V$ and $W$ are the trivial one-dimensional representations
 \cite{sbon}.  
However,
 we shall see that this phenomenon may occur
 in the real setting for vector bundles.  
Indeed, 
 the following theorem shows that 
 there may exist 
\index{B}{sporadicsymmetrybreakingoperator@sporadic symmetry breaking operator|textbf}
{\it{sporadic}} symmetry breaking operators
 which are differential operators
 in the case $[V:W]=0$:
\begin{theorem}
[localness theorem]
\label{thm:152347}
\index{B}{localnesstheorem@localness theorem|textbf}
Assume $[V:W] = 0$.  
Then 
\[
{\operatorname{Hom}}_{G'}(I_{\delta}(V,\lambda)|_{G'}, J_{\varepsilon}(W, \nu))
  =
{\operatorname{Diff}}_{G'}(I_{\delta}(V,\lambda)|_{G'}, J_{\varepsilon}(W, \nu))
\]
for all $(\lambda,\nu,\delta,\varepsilon) \in {\mathbb{C}}^2 \times \{\pm\}^2$, 
 that is, 
any symmetry breaking operator
 (if exists)
\[
   C^{\infty}(G/P, {\mathcal{V}}_{\lambda,\delta})
   \to 
   C^{\infty}(G'/P', {\mathcal{W}}_{\nu,\varepsilon})
\]
is a differential operator.  
\end{theorem}

Theorem \ref{thm:152347} is proved
 in Section \ref{subsec:SDone}.  
We call such operators {\it{sporadic}}
 because there is no regular symmetry breaking operator
 if $[V:W]=0$, 
 see Theorem \ref{thm:regexist} below.  
Another localness theorem is formulated 
 in Theorem \ref{thm:VWSBO} (2-b)
 (see also Proposition \ref{prop:1532102} in Chapter \ref{sec:DSVO})
 under the assumption
 that the parameter 
$(\lambda,\nu)
 \in {\mathbb{C}}^2$
 satisfies $\nu-\lambda \in {\mathbb{N}}$.  
\begin{example}
\label{ex:3.3}
Suppose $(V,W)=(\Exterior^i({\mathbb{C}}^n), \Exterior^j({\mathbb{C}}^{n-1}))$.  
Then $[V:W] \ne 0$
 if and only if $j=i-1$ or $i$.  
Hence Theorem \ref{thm:152347} tells that 
 there exists a nonlocal symmetry breaking operators
$
I_{\delta}(i,\lambda)  \to J_{\varepsilon}(j, \nu) 
$
 only if $j \in \{i-1, i\}$.  
(In fact,
 this is also a sufficient condition,
 see Theorem \ref{thm:152271}.)
On the other hand,
 there exist nontrivial differential symmetry breaking operators
 for some $(\lambda, \nu) \in {\mathbb{C}}^2$
 if and only if $j \in \{i-2, i-1, i, i+1\}$, 
 as is seen from the complete classification
 of differential symmetry breaking operators 
 (Fact \ref{fact:2.7}).  
Thus there exist sporadic (differential) symmetry breaking operators
 when $j=i-2$ or $i+1$.  
\end{example}

\begin{remark}
\label{rem:MVW}
The assumption $[V:W]\ne 0$
 in Theorems \ref{thm:unique} and \ref{thm:1716110}
 is {\it{not}} an intertwining property
 for 
\index{A}{M1prime@$M'=O(n-1) \times O(1)$}
$M'=M \cap G' \simeq O(n-1) \times O(1)$
 but for the subgroup $O(n-1)$
 which is of index two in $M'$.  
We note that 
for 
\index{A}{Vlnd@$V_{\delta}:= V \boxtimes {\delta}$|textbf}
$V_{\delta}:= V \boxtimes {\delta} \in \widehat M$
 and $W_{\varepsilon}:= W \boxtimes {\varepsilon} \in \widehat {M'}$, 
\index{A}{VWmult@$[V:W]$}
\[
\text{
$
\operatorname{Hom}_{M'}(V_{\delta}|_{M'}, W_{\varepsilon})\ne\{0\}
$
 if and only if $[V:W]\ne 0$ and $\delta=\varepsilon$.  
}
\]
Indeed the condition
$
\delta = \varepsilon
$ 
 is not included in the assumption of Theorem \ref{thm:unique}
 on the construction of regular symmetry breaking operators.  
The reason is clarified
 in Theorem \ref{thm:regexist} in the next subsection.  
\end{remark}

\subsubsection{Existence condition for regular symmetry breaking operators}
\label{subsec:regSBOexist}
A {\it{regular symmetry breaking operator}} is 
 an \lq\lq{opposite}\rq\rq\ notion
 to a differential symmetry breaking operator
 in the sense
 that the support of its distribution kernel contains 
an interior point
 in the real flag manifold,
 see \cite[Def.~3.3]{sbon}.  
(See also Definition \ref{def:regSBO} in our special setting.)  
In \cite[Cor.~3.6]{sbon}
 we give a necessary condition 
 for the existence of regular symmetry breaking operators
 in the general setting.
This condition is also sufficient in our setting:
\begin{theorem}
[existence of regular symmetry breaking operators]
\label{thm:regexist}
Suppose $V \in \widehat {O(n)}$ and $W \in \widehat{O(n-1)}$.  
Then the following three conditions on the pair $(V,W)$ are equivalent:
\begin{enumerate}
\item[{\rm{(i)}}]
$[V:W] \ne 0$.  
\item[{\rm{(ii)}}]
There exists a nonzero regular symmetry breaking operator from 
 the $G$-module $I_{\delta}(V,\lambda)$
 to the $G'$-module $J_{\varepsilon}(W,\nu)$
 for some $(\lambda, \nu,\delta, \varepsilon) \in {\mathbb{C}}^2 \times \{\pm\}^2$.  
\item[{\rm{(iii)}}]
For any $(\delta, \varepsilon) \in \{\pm\}^2$, 
 there is an open dense subset $U$ in ${\mathbb{C}}^2$
 such that a nonzero regular symmetry breaking operator exists from $I_{\delta}(V,\lambda)$
 to $J_{\varepsilon}(W,\nu)$
 for all $(\lambda, \nu) \in U$.  
\end{enumerate}
\end{theorem}
The proof will be given in Section \ref{subsec:regexist}.  
The open dense subset $U$ is explicitly given in Proposition \ref{prop:172425}.

\subsubsection{Integral operators, 
 analytic continuation, and normalization factors}
\label{subsec:VWSBO}

For an explicit construction of {\it{regular}} symmetry breaking operators,
 we use the reflection map 
\index{A}{1psin@$\psi_n$}
$\psi_n$ defined as follows: 

\index{A}{1psin@$\psi_n$|textbf}
\begin{equation}
\label{eqn:psim}
\psi_n \colon {\mathbb{R}}^n \setminus \{0\} \to O(n),
\quad
 x \mapsto I_n - \frac{2 x {}^{t\!} x}{|x|^2}.  
\end{equation}
Then $\psi_n(x)$ gives the reflection $\psi_n(x)$
 with respect to the hyperplane 
 $\{y \in {\mathbb{R}}^n:(x,y)=0\}$.  
Clearly,
 we have
\begin{equation}
\label{eqn:psix}
  \psi_n(x)=\psi_n(-x), 
\quad
\psi_n(x)^2=I_n, 
\quad
\text{ and }
\quad
\det \psi_n(x)=-1.  
\end{equation}

Suppose $(\sigma, V) \in \widehat{O(n)}$
 and $(\tau, W) \in \widehat{O(n-1)}$.  
For the construction of regular symmetry breaking operators,
 we need the condition $[V:W] \ne 0$, 
 see Theorem \ref{thm:regexist}.  
So let us assume 
$
[V:W] \ne 0.  
$
We fix a nonzero $O(n-1)$-homomorphism
\[
  \pr V W \colon V \to W, 
\]
which is unique up to scalar multiplication
 by Schur's lemma
 because $[V:W]=1$.  
We introduce a smooth map
\index{A}{RVW@$\Rij VW$|textbf}
\[
\Rij V W \colon {\mathbb{R}}^n \setminus \{ 0 \} 
            \to 
           \operatorname{Hom}_{\mathbb{C}}(V,W)
\]
by 
\begin{equation}
\label{eqn:RVW}
   \Rij VW := \pr VW  \circ \sigma \circ \psi_n.  
\end{equation}

In what follows,
 we use the coordinates $(x,x_n) \in {\mathbb{R}}^n = {\mathbb{R}}^{n-1} \oplus {\mathbb{R}}$
 where $x=(x_1, \cdots, x_{n-1})$, 
 and the $n$-th coordinate $x_n$
 will play a special role.  

We set
\index{A}{Actg1@$\Atcal {\lambda}{\nu}{+}{V,W}$|textbf}
\index{A}{Actg2@$\Atcal {\lambda}{\nu}{-}{V,W}$|textbf}
\begin{align}
\label{eqn:KVWpt}
\Atcal {\lambda}{\nu}{+}{V,W}
 :=&
\frac{1}{\Gamma(\frac{\lambda + \nu -n +1}{2})
         \Gamma(\frac{\lambda-\nu}{2})}
         (|x|^2+ x_n^2)^{-\nu}
|x_n|^{\lambda+\nu-n}
\Rij V W (x,x_n), 
\\
\label{eqn:KVWmt}
\Atcal {\lambda}{\nu}{-}{V,W}
 :=&
\frac{1}{\Gamma(\frac{\lambda + \nu -n +2}{2})
         \Gamma(\frac{\lambda-\nu+1}{2})}
         (|x|^2 +x_n^2)^{-\nu}
|x_n|^{\lambda+\nu-n}
 {\operatorname{sgn}} x_n
\Rij VW (x,x_n).  
\end{align}

\begin{theorem}
[regular symmetry breaking operators]
\label{thm:152389}
\index {B}{regularsymmetrybreakingoperator@regular symmetry breaking operator}
Suppose $[V:W] \ne 0$
 and $\gamma \in \{ \pm \}$.  
Then the distributions 
 $\Atcal \lambda \nu {\gamma} {V,W}$, 
 initially defined
 as $\operatorname{Hom}_{\mathbb{C}}(V,W)$-valued
 locally integrable functions
 on ${\mathbb{R}}^n$
 for $\operatorname{Re} \lambda \gg |\operatorname{Re} \nu|$, 
 extends to $P'$-invariant elements in 
$
   {\mathcal{D}}'(G/P, {\mathcal{V}}_{\lambda,\delta}^{\ast}) \otimes W_{\nu,\varepsilon}
$
 for all $(\lambda, \nu) \in {\mathbb{C}}^2$
 and $\delta, \varepsilon \in \{ \pm \}$
with $\delta \varepsilon = \gamma$.  
Then the distributions $\Atcal \lambda \nu {\gamma} {V,W}$
 induce a family of symmetry breaking operators
\index{A}{Ahtg0@$\Atbb \lambda \nu {\pm} {V,W}$|textbf}
\begin{alignat*}{2}
\Atbb \lambda \nu {\gamma} {V,W}
\colon \,&C^{\infty}(G/P, {\mathcal{V}}_{\lambda,\delta})
   \to 
   C^{\infty}(G'/P', {\mathcal{W}}_{\nu,\varepsilon}), 
\end{alignat*}
which depends holomorphically on $(\lambda, \nu)$
 in the entire ${\mathbb{C}}^2$.  
\end{theorem}

\begin{remark}
The denominator in \eqref{eqn:KVWpt}
 is different from the product
 of the denominators
 of the two distributions 
$\frac{(|x|^2 + x_n^2)^{-\nu}}{\Gamma(\frac{n-\nu}{2})}$
and 
$\frac{|x_n|^{\lambda+\nu-n}}{\Gamma(\frac{\lambda+\nu-n+1}{2})}$
 on ${\mathbb{R}}^n$
 that depend holomorphically 
 on $(\lambda,\nu)$ in the entire ${\mathbb{C}}^2$.  
In fact the product
\begin{equation}
\label{eqn:bad}
\frac{(|x|^2 + x_n^2)^{-\nu}}{\Gamma(\frac{n-\nu}{2})}
\times
\frac{|x_n|^{\lambda+\nu-n}}{\Gamma(\frac{\lambda+\nu-n+1}{2})}
\end{equation}
 does not always make sense
 as distributions on ${\mathbb{R}}^n$.  
For instance, 
 if $(\lambda,\nu)=(-1,n)$, 
 then the multiplication \eqref{eqn:bad} means 
 the multiplication 
 (up to nonzero scalar multiplication)
 of the Dirac delta functions
 $\delta(x_1, \cdots,x_n)$ 
 by $\delta(x_n)$, 
 which is not well-defined in the usual sense.  
\end{remark}

Theorem \ref{thm:152389} will be proved in Section \ref{subsec:holoAVW}.

We prove in Theorem \ref{thm:161243}
 that the normalization is optimal
 for $(V,W)=(\Exterior^i({\mathbb{C}}^n), \Exterior^j({\mathbb{C}}^{n-1}))$
 in the sense
 that the zeros of $\Atbb \lambda \nu \pm {V,W}$
 are of codimension $>1$
 in the parameter space
 of $(\lambda,\nu)$, 
 namely,
 discrete in ${\mathbb{C}}^2$
 in our setting.  
For the general $(V,W)$, 
 we shall give an upper and lower estimate
 of the null set of 
 the symmetry breaking operators 
 $\Atbb \lambda \nu {+} {V,W}$ 
 and $\Atbb \lambda \nu {-} {V,W}$
 in Theorem \ref{thm:1532113}.

\subsection{
Classification scheme of symmetry breaking operators: general case}
\label{subsec:170228}
\index{B}{classificationscheme@classication scheme, symmetry breaking operators}

In this section,
we give a general scheme 
 for the classification of all symmetry breaking operators
 $I_{\delta}(V,\lambda)|_{G'} \to J_{\varepsilon}(W,\nu)$
 between the two principal series representations
 of $G$ and the subgroup $G'$
 in full generality
 where $(\sigma,V)\in \widehat{O(n)}$
 and $(\tau,W)\in \widehat{O(n-1)}$.

We begin with conditions
on the parameter $(\lambda, \nu, \delta, \varepsilon)$
 for the existence of differential symmetry breaking operators.  
\begin{theorem}
[existence of differential symmetry breaking operators]
\index{B}{differentialsymmetrybreakingoperators@differential symmetry breaking operator}
\label{thm:vanDiff}
\qquad\qquad
\begin{enumerate}
\item[{\rm{(1)}}]
(Theorem {\ref{thm:vanishDiff}})
Suppose $\lambda, \nu \in {\mathbb{C}}$
 and $\delta, \varepsilon \in \{\pm\}$
 satisfy the generic parameter condition \eqref{eqn:nlgen}.  
Then,
\[
   {\operatorname{Diff}}_{G'}
   (I_{\delta}(V,\lambda)|_{G'}, J_{\varepsilon}(W,\nu))
   =\{0\}
\]
for any $(\sigma, V) \in \widehat{O(n)}$
 and $(\tau, W) \in \widehat{O(n-1)}$.  
\item[{\rm{(2)}}]
(Theorem {\ref{thm:existDSBO}})
Suppose $[V:W]\ne 0$.  
Then the converse statement holds,
 namely,
 if 
\index{A}{1psi@$\Psising$,
          special parameter in ${\mathbb{C}}^2 \times \{\pm\}^2$}
$(\lambda, \nu, \delta, \varepsilon) \in \Psising$
 (see \eqref{eqn:singset}), 
 then 
\[
   {\operatorname{Diff}}_{G'}
   (I_{\delta}(V,\lambda)|_{G'}, J_{\varepsilon}(W,\nu))
   \ne\{0\}.  
\]
\end{enumerate}
\end{theorem}
We give a proof for the first statement of Theorem \ref{thm:vanDiff}
 in Section \ref{subsec:vanDiff}, 
 and the second statement in Section \ref{subsec:extDiff}.  
Keeping Theorem \ref{thm:vanDiff} 
 on differential symmetry breaking operators in mind, 
 we state a general scheme for the classification
 of {\it{all}} symmetry breaking operators:

\begin{theorem}
[classification scheme
 of symmetry breaking operators]
\label{thm:VWSBO}
Let $n \ge 3$, 
 $(\sigma, V) \in \widehat{O(n)}$, 
 $(\tau,W) \in \widehat{O(n-1)}$, 
 $\lambda, \nu \in {\mathbb{C}}$
 and $\delta, \varepsilon \in \{\pm\}$.  
\begin{enumerate}
\item[{\rm{(1)}}]
Suppose 
\index{A}{VWmult@$[V:W]$}
$[V:W] =0$.  
Then 
\index{A}{DiffG@${\operatorname{Diff}}_{G'}(I_{\delta}(V,\lambda)\vert_{G'}, J_{\varepsilon}(W,\nu))$}
\[
  {\operatorname{Hom}}_{G'}(I_{\delta}(V,\lambda)|_{G'}, J_{\varepsilon}(W,\nu))  =
{\operatorname{Diff}}_{G'}(I_{\delta}(V,\lambda)|_{G'}, J_{\varepsilon}(W,\nu)).  
\]
\item[{\rm{(2)}}]
Suppose $[V:W]  \ne 0$.  
\begin{enumerate}
\item[{\rm{(2-a)}}]
{\rm{(generic case)}}\enspace
Suppose further 
 that $(\lambda, \nu, \delta, \varepsilon) \not\in \Psising$, 
 namely,
 it satisfies
 the generic parameter condition \eqref{eqn:nlgen}.  
Then 
\[
  {\operatorname{Hom}}_{G'}(I_{\delta}(V,\lambda)|_{G'}, J_{\varepsilon}(W,\nu))  =
{\mathbb{C}} \Atbb \lambda \nu {\delta\varepsilon}{V,W}.  
\]
In this case, 
 $\Atbb \lambda \nu {\delta\varepsilon}{V,W}$ is nonzero 
 and is not a differential operator.  
\item[{\rm{(2-b)}}]
{\rm{(special parameter case I, localness theorem)}}\enspace
\index{B}{localnesstheorem@localness theorem|textbf}
Suppose $\Atbb \lambda \nu {\delta\varepsilon}{V,W} \ne 0$
 and $(\lambda, \nu, \delta, \varepsilon)\in \Psising$
 ({\it{i.e.,}} does not satisfy the generic parameter condition \eqref{eqn:nlgen}).  
Then any symmetry breaking operator
 (in particular, 
 $\Atbb \lambda \nu {\delta\varepsilon}{V,W}$) is a differential operator and 
\[
  {\operatorname{Hom}}_{G'}(I_{\delta}(V,\lambda)|_{G'}, J_{\varepsilon}(W,\nu))  =
   {\operatorname{Diff}}_{G'}(I_{\delta}(V,\lambda)|_{G'}, J_{\varepsilon}(W,\nu))
\ni
\Atbb \lambda \nu {\delta\varepsilon}{V,W}.  
\]
\item[{\rm{(2-c)}}]
{\rm{(special parameter case II)}}\enspace
Suppose $\Atbb \lambda \nu {\delta\varepsilon}{V,W} = 0$.  
Then $(\lambda, \nu, \delta, \varepsilon)\in \Psising$, 
 and the 
 \index{B}{renormalized regular symmetry breaking operator@regular symmetry breaking operator, renormalized---}
renormalized operator 
\index{A}{Ahttgln0@$\Attbb {\lambda}{\nu}{\pm}{V,W}$}
$\Attbb \lambda \nu {\delta\varepsilon}{V,W}$
 (see Section \ref{subsec:AVW})
 gives a nonzero symmetry breaking operator
 which is not a differential operator.
We have
\[
  {\operatorname{Hom}}_{G'}(I_{\delta}(V,\lambda)|_{G'}, J_{\varepsilon}(W,\nu))  =
  {\mathbb{C}} \Attbb \lambda \nu {\delta\varepsilon}{V,W}
  \oplus
  {\operatorname{Diff}}_{G'}(I_{\delta}(V,\lambda)|_{G'}, 
 J_{\varepsilon}(W,\nu)).  
\]
In particular, 
\[
  \dim_{\mathbb{C}}
  {\operatorname{Hom}}_{G'}(I_{\delta}(V,\lambda)|_{G'}, J_{\varepsilon}(W,\nu))  \ge 2.  
\]
\end{enumerate}
\end{enumerate}
\end{theorem}

The first assertion of Theorem \ref{thm:VWSBO}
 is a restatement of Theorem \ref{thm:152347}.  
The case (2-a) is given in Theorem \ref{thm:genbasis}
 and the case (2-b) is in Proposition \ref{prop:1532102}. 
The first statement 
 for the case (2-c) is proved in Theorem \ref{thm:170340} (1).  
The direct sum decomposition is given in Corollary \ref{cor:Azero}.  
The last statement follows from the existence 
 of nonzero differential symmetry breaking operators
 for all special parameters
 (Theorem \ref{thm:vanDiff} (2)).

Theorems \ref{thm:1716110} and \ref{thm:VWSBO} lead us to a vanishing result
 of symmetry breaking operators
 as follows:  

\begin{corollary}
[vanishing of symmetry breaking operators]
\label{cor:VWvanish}
Let $(\sigma, V) \in \widehat{O(n)}$, 
 $(\tau, W) \in \widehat{O(n-1)}$, 
 $\lambda, \nu \in {\mathbb{C}}$
 and $\delta, \varepsilon \in \{\pm\}$.  
If $[V:W]=0$ and 
 $(\lambda, \nu)$ satisfies 
 the generic parameter condition \eqref{eqn:nlgen}, 
 then 
\[
   {\operatorname{Hom}}_{G'}
   (I_{\delta}(V,\lambda)|_{G'}, J_{\varepsilon}(W,\nu))
   =\{0\}.  
\]
\end{corollary}

\begin{proof}
By Theorem \ref{thm:VWSBO} (1), 
 we have
\[
   {\operatorname{Hom}}_{G'}
   (I_{\delta}(V,\lambda)|_{G'}, J_{\varepsilon}(W,\nu))
   ={\operatorname{Diff}}_{G'}
    (I_{\delta}(V,\lambda)|_{G'}, J_{\varepsilon}(W,\nu))
\]
because $[V:W]=0$.  
In turn, 
 the right-hand side reduces to zero 
 by Theorem \ref{thm:1716110}
 because of the generic parameter condition \eqref{eqn:nlgen}.  
\end{proof}

Theorem \ref{thm:VWSBO} gives a classification 
 of symmetry breaking operators
 up to the following two problems:
\begin{enumerate}
\item[$\bullet$]
the location of zeros of the normalized
 regular symmetry breaking operator
 $\Atbb \lambda \nu \gamma {V,W}$;
\item[$\bullet$]
the classification of {\it{differential}} symmetry breaking operators.  
\end{enumerate}

For $(V,W)=(\Exterior^i({\mathbb{C}}^n), \Exterior^j({\mathbb{C}}^{n-1}))$, 
 these two problems are solved explicitly
 in Theorem \ref{thm:161243} and Fact \ref{fact:2.7}, 
 respectively, 
 and thus we accomplish the complete classification
 of symmetry breaking operators.  
This will be stated
 in Theorem \ref{thm:1.1} (multiplicity formula)
 and in Theorem \ref{thm:SBObasis} (explicit generators).

\subsection{Summary: vanishing of regular 
 symmetry breaking operators $\Atbb \lambda \nu \pm {V,W}$}
\label{subsec:vanAVW}
As we have seen in the classification scheme
 (Theorem \ref{thm:VWSBO}) 
 for all symmetry breaking operators, 
 the parameter $(\lambda,\nu,\delta,\varepsilon)$ 
 for which the (generically) regular symmetry breaking operator
 $\Atbb \lambda \nu \pm {V,W}$
 vanishes plays a crucial role
 in the classification theory.  
For $(\lambda,\nu,\delta,\varepsilon) \in \Psising$, 
 we noted:
\begin{enumerate}
\item[$\bullet$]
when $\Atbb {\lambda_0} {\nu_0} \pm {V, W}=0$, 
 we can construct a nonzero symmetry breaking operator
 $\Attbb {\lambda_0} {\nu_0} \pm {V, W}$
 by \lq\lq{renormalization}\rq\rq\
 which is {\it{not}} a differential operator
 (Theorem \ref{thm:170340});
\item[$\bullet$]
when $\Atbb {\lambda_0} {\nu_0} \pm {V, W}\ne 0$, 
we prove a 
\index{B}{localnesstheorem@localness theorem}
 {\it{localness theorem}}
 asserting
 that all symmetry breaking operators are differential operators
 (Proposition \ref{prop:1532102}).  
\end{enumerate}

We obtain a condition 
 for the (non) vanishing of $\Atbb {\lambda} {\nu} \pm {V, W}$
 as follows.  
Using the same notation
 as in \cite[Chap.~1]{sbon}, 
 we define the following two subsets in ${\mathbb{Z}}^2$:
\index{A}{Leven@$L_{\operatorname{even}}$|textbf}
\index{A}{Lodd@$L_{\operatorname{odd}}$|textbf}
\begin{alignat}{2}
\label{eqn:Leven}
L_{\operatorname{even}}:=
&\left\{ \right.
 (-i,-j):
0 \le j\leq i \mbox{ and } i\equiv j 
&& \mod 2 \left. \right \},
\\
\label{eqn:Lodd}
L_{\operatorname{odd}}:=
&\left\{ \right.
  (-i,-j)
: 0 \le j\leq i \mbox{ and } i \equiv j +1 
&& \mod 2 \left. \right \}.
\end{alignat}

\begin{theorem}
\label{thm:1532113}
Let $(\sigma, V) \in \widehat {O(n)}$
 and $(\tau, W) \in \widehat {O(n-1)}$
 with $[V:W]\ne 0$.  
\begin{enumerate}
\item[{\rm{(1)}}]
There exists 
\index{A}{Nsigma@$N(\sigma)$}
$
N(\sigma) \in {\mathbb{N}}
$
such that
\begin{alignat*}{2}
\Atbb \lambda \nu {+} {V,W} =&0
\qquad
&&\text{if } (\lambda,\nu) \in L_{\operatorname{even}}
\text{ and } \nu \le -N(\sigma), 
\\
\Atbb \lambda \nu {-} {V,W} =&0
\qquad
&&\text{if } (\lambda,\nu) \in L_{\operatorname{odd}}
\text{ and } \nu \le -N(\sigma).  
\end{alignat*}
\item[{\rm{(2)}}]
If $\Atbb \lambda \nu {+} {V,W} =0$
 then 
 $\nu-\lambda \in 2{\mathbb{N}}$;
 if $\Atbb \lambda \nu {-} {V,W} =0$
then 
 $\nu - \lambda \in 2 {\mathbb{N}}+1$.    
\end{enumerate}
\end{theorem}

\begin{remark}
We shall show in Lemma \ref{lem:Nsigma}
 that $N(\sigma)$ can be taken
 to be $\ell(\sigma)$, 
 as defined in \eqref{eqn:lsigma}.  
\end{remark}

Theorem \ref{thm:1532113} (2) is a part
 of Theorem \ref{thm:VWSBO} (2), 
 and will be proved in Section \ref{subsec:1532113}.

Combining Theorems \ref{thm:VWSBO} and \ref{thm:1532113}, 
 we see that there exist infinitely many 
 $(\lambda,\nu)\in {\mathbb{C}}^2$
 such that the multiplicity 
$
   m(I_{\delta}(V,\lambda), J_{\varepsilon}(W,\nu))>1
$
 as follows:
\begin{corollary}
\label{cor:170821}
Let $(\sigma,V) \in \widehat {O(n)}$ and $(\tau,W) \in \widehat{O(n-1)}$
 satisfy $[V:W] \ne 0$.  
If 
\[
  (\lambda, \nu) \in 
  \begin{cases}
  L_{\operatorname{even}} \cap \{\nu \le - N(\sigma)\}
  \quad
  &\text{for \,$\delta \varepsilon =+$}, 
\\
  L_{\operatorname{odd}} \cap \{\nu \le - N(\sigma)\}
  \quad
  &\text{for \,$\delta \varepsilon =-$}, 
  \end{cases}
\]
then we have
\[
  \dim_{\mathbb{C}}
  {\operatorname{Hom}}_{G'}
  (I_{\delta}(V,\lambda)|_{G'}, J_{\varepsilon}(W,\nu))>1.  
\]
\end{corollary}
By Theorem \ref{thm:1532113},
 we get readily the following corollary, 
 to which we shall return in Chapter \ref{sec:conjecture}
 (see Example \ref{ex:317}).  
\begin{corollary}
\label{cor:IV1}
Suppose that 
$\Atbb \lambda \nu \delta{V,W}=0$.  
Then $\Atbb {n-\lambda}{n-1-\nu} \delta{V,W} \ne 0$.  
\end{corollary}

Theorem \ref{thm:1532113} means that
\index{A}{Leven@$L_{\operatorname{even}}$}
\index{A}{Lodd@$L_{\operatorname{odd}}$}
\begin{align*}
L_{\operatorname{even}}
\cap 
\{\nu \le -N(\sigma) \}
\subset
&
\{(\lambda,\nu) \in {\mathbb{C}}^2
:
\Atbb \lambda \nu {+} {V,W} =0\}
\subset
\{(\lambda,\nu) \in {\mathbb{C}}^2
:
\nu-\lambda \in 2 {\mathbb{N}}\}, 
\\
L_{\operatorname{odd}}
\cap 
\{\nu \le -N(\sigma) \}
\subset
&
\{(\lambda,\nu) \in {\mathbb{C}}^2
:
\Atbb \lambda \nu {-} {V,W} =0\}
\subset
\{(\lambda,\nu) \in {\mathbb{C}}^2
:
\nu-\lambda \in 2 {\mathbb{N}}+1\}.  
\end{align*}

We shall determine in Theorem \ref{thm:161243}
 the set $\{(\lambda,\nu) \in {\mathbb{C}}^2:
\Atbb \lambda \nu \gamma {V,W} =0\}$
 for $\gamma = \pm$
 in the special case where $(V,W)=(\Exterior^i({\mathbb{C}}^n), \Exterior^j({\mathbb{C}}^{n-1}))$.  
If $\sigma$ is the $i$-th exterior representation $\sigma^{(i)}$
 on $\Exterior^i({\mathbb{C}}^n)$, 
 then we can take $N(\sigma)$ to be 0
 if $i=0$ or $n$;
 to be 1 if $1 \le i \le n-1$.  
In this case, 
 the left inclusion is almost a bijection.  
On the other hand, 
 concerning the right inclusions,
 we refer to Theorem \ref{thm:VWSBO} (2-b), 
 which will be proved
 in Section \ref{subsec:170213},
 see Proposition \ref{prop:1532102}.

\subsection{The classification of symmetry breaking operators
 for differential forms}
\label{subsec:exhaust}

Let $(G,G')=(O(n+1,1),O(n,1))$ with $n \ge 3$
 as before.  
We consider the special case
\[
  (V,W)=(\Exterior^i({\mathbb{C}}^{n}), \Exterior^j({\mathbb{C}}^{n-1})).  
\]
Then the corresponding principal series representations
 $I_{\delta}(V,\lambda)$ of $G$
 and $J_{\varepsilon}(W,\nu)$ of the subgroup $G'$
 are denoted by $I_{\delta}(i,\lambda)$
 and $J_{\varepsilon}(j,\nu)$, 
 respectively.  
In this section 
 we summarize the complete classification
 of symmetry breaking operators from 
 the $G$-module $I_{\delta}(i,\lambda)$ 
 to the $G'$-module $J_{\varepsilon}(j,\nu)$.  
The main results are stated in Theorems \ref{thm:1.1} and \ref{thm:SBObasis}.  
Our results rely on the vanishing condition
 of the normalized regular symmetry breaking operators
 $\Atbb \lambda \nu \gamma {i,j}$
 (Theorem \ref{thm:161243}) 
 and the classification
 of differential symmetry breaking operators
 (Fact \ref{fact:2.7}).

\subsubsection{Vanishing condition
 for the regular symmetry breaking operators 
 $\Atbb \lambda \nu \gamma {i,j}$}
\label{subsec:Aijvanish}

We apply the general construction
 of the (normalized) symmetry breaking operators
 $\Atbb \lambda \nu \gamma {V,W}$
 in Theorem \ref{thm:152389}
 to the pair
 of representations
$
  (V,W)=(\Exterior^i({\mathbb{C}}^{n}), \Exterior^j({\mathbb{C}}^{n-1})).  
$  
Then we obtain (normalized) symmetry breaking operators, 
 to be denoted by 
 $\Atbb \lambda \nu \gamma {i,j}$, 
that depend holomorphically on $(\lambda,\nu)$
 in the entire complex plane ${\mathbb{C}}^2$
 if $j \in \{i-1,i\}$
 and $\gamma \in \{\pm\}$, 
 see Theorem \ref{thm:1522a}.

We determine the zero set of $\Atbb \lambda \nu {\gamma} {i,j}$
 explicitly as follows:
\begin{theorem}
[zeros of regular symmetry breaking operators
 $\Atbb \lambda \nu \pm {i,j}$]
\label{thm:161243}
~~~~~
\begin{enumerate}
\item[{\rm{(1)}}]
For $0 \le i \le n-1$, 
\begin{multline*}
\{(\lambda, \nu) \in {\mathbb{C}}^2
:
\Atbb \lambda \nu + {i,i} =0\}
\\
= 
\begin{cases}
L_{\operatorname{even}}
\quad
&
\text{if $i=0$}, 
\\
(L_{\operatorname{even}} \setminus \{\nu=0\}) \cup \{(i,i)\}
\quad
&
\text{if $1 \le i \le n-1$}.  
\end{cases}
\end{multline*}
\item[{\rm{(2)}}]
For $1 \le i \le n$, 
\begin{multline*}
\{(\lambda, \nu) \in {\mathbb{C}}^2
:
\Atbb \lambda \nu + {i,i-1} =0\}
\\
=
\begin{cases}
(L_{\operatorname{even}}\setminus \{\nu=0\}) \cup \{(n-i,n-i)\}
\quad
&
\text{if $1 \le i \le n-1$}, 
\\
L_{\operatorname{even}} 
\quad
&
\text{if $i=n$}.  
\end{cases}
\end{multline*}
\item[{\rm{(3)}}]
For $0 \le i \le n-1$, 
\begin{multline*}
\{(\lambda, \nu) \in {\mathbb{C}}^2
:
\Atbb \lambda \nu - {i,i} =0\}
\\
= 
\begin{cases}
L_{\operatorname{odd}}
\quad
&
\text{if $i=0$}, 
\\
L_{\operatorname{odd}}\setminus \{\nu=0\} 
\quad
&
\text{if $1 \le i \le n-1$}.  
\end{cases}
\end{multline*}
\item[{\rm{(4)}}]
For $1 \le i \le n$, 
\begin{multline*}
\{(\lambda, \nu) \in {\mathbb{C}}^2
:
\Atbb \lambda \nu - {i,i-1} =0\}
\\
=
\begin{cases}
L_{\operatorname{odd}}\setminus \{\nu=0\}
\quad
&
\text{if $1 \le i \le n-1$}, 
\\
L_{\operatorname{odd}} 
\quad
&
\text{if $i=n$}.  
\end{cases}
\end{multline*}
\end{enumerate}
\end{theorem}

Theorem \ref{thm:161243} will be proved in Section \ref{subsec:Aijzero}
 by using the residue formula
 of $\Atbb \lambda \nu{\pm}{i,j}$
 (\cite{xkresidue}).

A special case of Theorem \ref{thm:161243} includes the following.  
\begin{example}
\label{ex:3.19}
\begin{enumerate}
\item[{\rm{(1)}}]
For $0 \leq i \leq n$, 
$\Atbb i i + {i,i}=0$
 and 
$\Atbb {n-i} {n-i-1} + {i,i} \ne 0$.  
\item[{\rm{(2)}}]
For $0 \leq i \leq n-1$, 
$\Atbb i i + {n-i,n-i-1}  = 0$
 and $\Atbb {n-i}{n-i-1} + {n-i,n-i-1} \ne 0$. 
\end{enumerate}
\end{example}

\begin{remark}
In the case $i=0$, 
 $\Atbb \lambda \nu + {i,i}$
 is the scalar-valued
 symmetry breaking operator
 induced from the scalar-valued distribution 
 $\Atcal \lambda \nu + {}$, 
 as we recall from \eqref{eqn:KAlnn+}.  
Thus the case $i=0$ in (1)
 was proved in \cite[Thm.~8.1]{sbon}.  
\end{remark}

\subsubsection{Differential symmetry breaking operators}

We review from \cite{KKP} the notation
 of conformally equivariant {\it{differential}} operators
 ${\mathcal{E}}^i(S^n) \to {\mathcal{E}}^j(S^{n-1})$, 
 namely,
 {\it{differential}} symmetry breaking operators 
 $I_{\delta}(V,\lambda)|_{G'} \to J_{\varepsilon}(W,\nu)$
 with $(V,W)=(\Exterior^i({\mathbb{C}}^n), \Exterior^j({\mathbb{C}}^{n-1}))$.  
The complete classification of {\it{differential}} symmetry breaking operators
 was recently accomplished in \cite[Thm.~2.8]{KKP}
 based on the F-method
 \cite{xkHelg85}.  

\begin{fact}
[classification of differential symmetry breaking operators]
\label{fact:2.7}
Let $n \ge 3$.  
Suppose $0 \le i \le n$, 
 $0 \le j \le n-1$, 
 $\lambda, \nu \in {\mathbb{C}}$, 
 and $\delta, \varepsilon \in \{ \pm \}$.  
Then the following three conditions on 6-tuple 
 $(i,j, \lambda,\nu,\delta, \varepsilon)$ are equivalent.  
\begin{enumerate}
\item[{\rm{(i)}}]
\index{A}{DiffG@${\operatorname{Diff}}_{G'}(I_{\delta}(V,\lambda)\vert_{G'}, J_{\varepsilon}(W,\nu))$}
\index{A}{Ideltai@$I_{\delta}(i, \lambda)$}
\index{A}{Jjepsilon@$J_{\varepsilon}(j, \nu)$}
$
\operatorname{Diff}_{G'}(I_{\delta}(i,\lambda)|_{G'}, J_{\varepsilon}(j,\nu))\ne \{0\}$.  
\item[{\rm{(ii)}}]
$\dim_{\mathbb{C}} \operatorname{Diff}_{G'}(I_{\delta}(i,\lambda)|_{G'}, J_{\varepsilon}(j,\nu))=1$.  
\item[{\rm{(iii)}}]
$\nu-\lambda \in {\mathbb{N}}$, 
 $(-1)^{\nu-\lambda}=\delta \varepsilon$, 
 and one of the following conditions holds:
\begin{enumerate}
\item[{\rm{(a)}}]
$j=i-2$, 
$2 \le i \le n-1$, 
$(\lambda,\nu)=(n-i,n-i+1)$;
\item[{\rm{(a$'$)}}]
$(i,j)=(n,n-2)$, 
 $-\lambda \in {\mathbb{N}}$, 
 $\nu=1$;
\item[{\rm{(b)}}]
$j=i-1$, 
$1 \le i \le n$;
\item[{\rm{(c)}}]
$j=i$, 
$0 \le i \le n-1$;
\item[{\rm{(d)}}]
$j=i+1$, 
$1 \le i \le n-2$, 
$(\lambda,\nu)=(i,i+1)$;
\item[{\rm{(d$'$)}}]
$(i,j)=(0,1)$, 
 $-\lambda \in {\mathbb{N}}$, 
 $\nu=1$.  
\end{enumerate}
\end{enumerate} 
\end{fact}

The generators are explicitly constructed
 in \cite[(2.24)--(2.32)]{KKP}
(see \cite{Juhl, KOSS, sbon} for the $i=0$ case), 
 which we review quickly.  
Let $\widetilde C_{\ell}^{\alpha}(z)$ 
be the 
\index{B}{Gegenbauerpolynomial@Gegenbauer polynomial|textbf}
Gegenbauer polynomial of degree $\ell$, 
 normalized by
\index{A}{CGegenbauernorm@$\widetilde C_l^{\alpha}(z)$, normalized Gegenbauer polynomial|textbf}
\begin{equation}
\label{eqn:Gegen}
\widetilde C_{\ell}^{\alpha}(z)
:=
\frac{1}{\Gamma(\alpha + [\frac{\ell+1}{2}])}
\sum_{k=0}^{[\frac \ell 2]}
(-1)^k
\frac{\Gamma(\ell-k+\alpha)}{k! (\ell-2k)!}
(2 z)^{\ell-2k}
\end{equation}
 as in \cite[(14.3)]{KKP}.  
Then $\widetilde C_{\ell}^{\alpha}(z) \not \equiv 0$
 for all $\alpha \in {\mathbb{C}}$
 and $\ell \in {\mathbb{N}}$.

For $\ell \in {\mathbb{N}}$, 
 we inflate $\widetilde C_{\ell}^{\alpha} (z)$
 to a polynomial of two variables
 $x$ and $y$:
\begin{align}
  \widetilde C_{\ell}^{\alpha} (x,y)
  :=& x^{\frac \ell 2} \widetilde C_{\ell}^{\alpha} (\frac{y}{\sqrt x})
\notag
\\
  =& \sum_{k=0}^{[\frac \ell 2]} 
     \frac{(-1)^k \Gamma(\ell-k+\alpha)}
          {\Gamma(\alpha + [\frac{\ell +1}{2}])\Gamma(\ell-2k+1)k!}(2y)^{\ell -2k}x^k.  
\label{eqn:Gegentwo}
\end{align}
For instance, 
 $\widetilde C_{0}^{\alpha} (x,y)=1$, 
 $\widetilde C_{1}^{\alpha} (x,y)=2y$, 
 $\widetilde C_{2}^{\alpha} (x,y)=2(\alpha+1)y^2-x$, etc.  
Notice that $\widetilde C_{\ell}^{\alpha} (x^2,y)$ is a homogeneous
 polynomial of $x$ and $y$ of degree $\ell$.

For $\nu-\lambda \in {\mathbb{N}}$, 
we set
 a scalar-valued differential operator
 $\Ctbb \lambda \nu {}:C^{\infty}({\mathbb{R}}^n)
 \to 
 C^{\infty}({\mathbb{R}}^{n-1})$
 by 
\index{A}{Chlmdn@$\Ctbb \lambda \nu {}$, Juhl's operator|textbf}
\begin{equation}
\label{eqn:Cln}
   \Ctbb \lambda \nu {}
   :=
   \operatorname{Rest}_{x_n=0} \circ
   \widetilde C_{\nu-\lambda}^{\lambda -\frac{n-1}{2}}
   (-\Delta_{{\mathbb{R}}^{n-1}}, \frac{\partial}{\partial x_n}).  
\end{equation}

For $\mu \in {\mathbb{C}}$
 and $a \in {\mathbb{N}}$, 
 we set 
\index{A}{0cgammamualpha@$\gamma(\mu,a)$|textbf}
\begin{equation}
\label{eqn:gamma}
\gamma(\mu,a)
:=
\begin{cases}
1
&\text{if $a$ is odd, }
\\
\mu + \frac a 2
\qquad
&\text{if $a$ is even.  }
\end{cases}
\end{equation}
We are ready 
 to define matrix-valued differential operators
\[
   \Ctbb \lambda \nu {i,j}\colon
 {\mathcal{E}}^i({\mathbb{R}}^n)
 \to 
 {\mathcal{E}}^j({\mathbb{R}}^{n-1})
\]
which were introduced in \cite[(2.24) and (2.26)]{KKP}
 by the following formul{\ae}:
\index{A}{Ciiln@$\Cbb \lambda \nu {i,j}$, matrix-valued differential operator|textbf}
\begin{equation}
\label{eqn:Ciiln}
\Cbb \lambda \nu {i, i}
:=
\Ctbb{\lambda+1}{\nu-1}{} d_{\mathbb{R}^n}d_{\mathbb{R}^n}^{\ast}
-
\gamma(\lambda - \frac n 2, \nu -\lambda) 
\Ctbb{\lambda}{\nu-1}{} d_{\mathbb{R}^n} \iota_{\frac{\partial}{\partial x_n}}
+\frac 1 2(\nu-i)
\Ctbb{\lambda}{\nu}{}, 
\end{equation}

\index{A}{Ciiln@$\Cbb \lambda \nu {i,j}$, matrix-valued differential operator|textbf}
\begin{equation}
\label{eqn:Cijln}
\Cbb \lambda \nu {i,i-1}
:=
- \Ctbb{\lambda+1}{\nu-1}{} d_{\mathbb{R}^n} d_{\mathbb{R}^n}^{\ast}
 \iota_{\frac{\partial}{\partial x_n}}
-
\gamma(\lambda - \frac {n-1} 2, \nu -\lambda) 
\Ctbb{\lambda+1}{\nu}{} d_{\mathbb{R}^n}^{\ast}
+
\frac 1 2
(\lambda+i-n)
\Ctbb{\lambda}{\nu}{}\iota_{\frac{\partial}{\partial x_n}}.  
\end{equation}
Here 
$
   \iota_{Z}\colon {\mathcal{E}}^i({\mathbb{R}^n}) \to {\mathcal{E}}^{i-1}({\mathbb{R}^n})
$
 stands
 for the interior product 
 which is defined to be the contraction with a vector field $Z$.

We note that 
\begin{alignat*}{4}
&\Cbb \lambda \nu {0,0}
&&=\frac 1 2 \nu \Ctbb \lambda \nu {}, 
\qquad
&&\Cbb \nu \nu {i,i}
&&=\frac 1 2 (\nu-i){\operatorname{Rest}}_{x_n=0}, 
\\
&\Cbb \lambda \lambda {i,i-1}
&&=\frac 1 2 (\lambda+i-n){\operatorname{Rest}}_{x_n=0}
 \circ \iota_{\frac{\partial}{\partial x_n}}, 
\qquad
&&\Cbb \lambda \nu {n,n-1} 
&&= \frac 1 2 \nu \Ctbb \lambda \nu {}
    \circ \iota_{\frac{\partial}{\partial x_n}}.  
\end{alignat*}

The operators $\Cbb \lambda \nu {i,j}$ vanish
 for the following special values of $(\lambda, \nu)$:
\begin{alignat*}{3}
\Cbb \lambda \nu {i,i} = & 0
\qquad
&&\text{if and only if }
\lambda =  \nu = i\,\,
&&\text{ or }
\nu=i=0, 
\\
\Cbb \lambda \nu {i,i-1} = & 0
\qquad
&&\text{if and only if }
\lambda =  \nu =n-i\,\,
&&\text{ or }
\nu=n-i=0.  
\end{alignat*}

In order to provide {\it{nonzero}} operators, 
 following the notation as in \cite[(2.30)]{KKP}, 
 we renormalize $\Cbb \lambda \nu {i,j}$ as 
\index{A}{Ctiiln@$\Ctbb \lambda \nu {i,i}$|textbf}
\index{A}{Ctiim1ln@$\Ctbb \lambda \nu {i,i-1}$|textbf}
\begin{align}
\Ctbb \lambda \nu {i,i}
&:=
\begin{cases}
\operatorname{Rest}_{x_n=0}
\qquad
&\text{if $\lambda = \nu$}, 
\\
\Ctbb \lambda \nu {}
&\text{if $i=0$}, 
\\
\Cbb \lambda \nu {i,i}
&\text{otherwise}, 
\end{cases}
\label{eqn:Ciitilde}
\\
\Ctbb \lambda \nu {i,i-1}
&:=
\begin{cases}
\operatorname{Rest}_{x_n=0}\circ \iota_{\frac{\partial}{\partial x_n}}
\qquad
&\text{if $\lambda = \nu$}, 
\\
\Ctbb \lambda \nu {}\circ \iota_{\frac{\partial}{\partial x_n}}
&\text{if $i=n$}, 
\\
\Cbb \lambda \nu {i,i-1}
&\text{otherwise}.  
\end{cases}
\label{eqn:Cii-1tilde}
\end{align}

For $j=i-2$ or $i+1$, 
 we also set 
\index{A}{Ctii2ln@$\Ctbb \lambda {n-i+1} {i,i-2}$|textbf}
\index{A}{Ctii1ln@$\Ctbb \lambda {i+1} {i,i+1}$|textbf}
\begin{align*}
\Ctbb \lambda {n-i+1} {i,i-2}
:=&
\begin{cases}
- d_{\mathbb{R}^{n-1}}^{\ast} \circ \Ctbb \lambda 0 {n,n-1}
&\text{if $i=n$, $\lambda \in -{\mathbb{N}}$}, 
\\
\operatorname{Rest}_{x_n=0}
\circ 
\iota_{\frac{\partial}{\partial x_n}} 
d_{\mathbb{R}^{n}}^{\ast}
\qquad
&\text{if $2 \le i \le n-1$, $\lambda = n-i$}.   
\end{cases}
\\
\Ctbb \lambda {i+1} {i,i+1}
:=&
\begin{cases}
d_{{\mathbb{R}}^{n-1}}
\circ 
\Ctbb \lambda 0 {}
\qquad\hspace{6mm}
&\text{if $i=0$, $\lambda  \in -{\mathbb{N}}$}, 
\\
\operatorname{Rest}_{x_n=0}
\circ 
d_{\mathbb{R}^{n}}
&\text{if $1 \le i \le n-2$, $\lambda=i$}.  
\end{cases}
\end{align*}

With the notation as above,
 we can describe explicit generators
 of the space ${\operatorname{Diff}}_{G'}(I_{\delta}(i,\lambda)|_{G'}, J_{\varepsilon}(j,\varepsilon))$
 of differential symmetry breaking operators:
\begin{fact}
[basis, {\cite[Thm.~2.9]{KKP}}]
\label{fact:3.9}
\index{B}{differentialsymmetrybreakingoperators@differential symmetry breaking operator}
Suppose that 6-tuple $(i,j,\lambda,\nu, \delta, \varepsilon)$ is 
 one of the six cases
 in Fact \ref{fact:2.7} (iii).  
Then the differential symmetry breaking operators
 $I_{\delta}(i,\lambda) \to J_{\varepsilon}(j,\nu)$
 are proportional to 
\begin{alignat*}{2}
&j=i-2:\,\,&& \Ctbb {n-i}{n-i+1}{i,i-2} \, (2 \le i \le n-1);
              \quad
              \Ctbb {\lambda}{1}{n,n-2} \, (i=n), 
\\
&j=i-1:&&  \Ctbb \lambda \nu{i,i-1}, 
\\
&j=i:&& \Ctbb \lambda \nu {i,i}, 
\\
&j=i+1:&& \Ctbb {i}{i+1}{i,i+1} \, (1 \le i \le n-2);
          \quad
          \Ctbb {\lambda}{1}{0,1} \, (i=0).  
\end{alignat*}
\end{fact}
\begin{remark}
The scalar case ($i=j=0$)
 was classified
 in Juhl \cite{Juhl} for $n \ge 3$.  
See also \cite{KOSS}
 for a different approach 
 using the 
\index{B}{Fmethod@F-method}
F-method.  
The case $n=2$ (and $i=j=0$)
 is essentially equivalent to find 
 differential symmetry breaking operators from 
 the tensor product of two principal series representations
 to another principal series representation 
 for $SL(2,{\mathbb{R}})$.  
In this case,
 generic (but not all) operators are given
 by the 
\index{B}{RankinCohenbracket@Rankin--Cohen bracket}
Rankin--Cohen brackets, 
 and the complete classification
 was accomplished in \cite[Thms.~9.1 and 9.2]{KP2}.  
We note
 that the dimension of differential symmetry breaking operators
 may jump to two 
 at some singular parameters
 where $n=2$.  
\end{remark}

\subsubsection{Formula
 of the dimension of 
$
   {\operatorname{Hom}}_{G'}
   (I_{\delta}(i,\lambda)|_{G'},J_{\varepsilon}(j,\nu))
$}
For admissible smooth representations $\Pi$ of $G$
 and $\pi$ of the subgroup $G'$, 
 we set
\index{A}{mPipi@$m(\Pi, \pi)$, multiplicity|textbf}
\[
   m(\Pi, \pi) := \dim_{\mathbb{C}} {\operatorname{Hom}}_{G'}
                               (\Pi|_{G'}, \pi).  
\]
In this subsection
 we give a formula of the multiplicity $m(\Pi, \pi)$
 for $\Pi=I_{\delta}(i,\lambda)$
 and $\pi=J_{\varepsilon}(j,\nu)$.

\begin{theorem}
[multiplicity formula]
\label{thm:1.1}
Let $(G,G')=(O(n+1,1),O(n,1))$
 with $n \ge 3$.  
Suppose $\Pi=I_{\delta}(i,\lambda)$
 and $\pi=J_{\varepsilon}(j,\nu)$
 for $0 \le i\le n$, $0 \le j \le n-1$,
 $\delta$, $\varepsilon \in \{\pm \}$, 
 and $\lambda,\nu \in {\mathbb{C}}$.  
Then we have the following.  
\begin{enumerate}
\item[{\rm{(1)}}]
\begin{alignat*}{2}
m(\Pi,\pi) \in &\{ 1,2 \} \qquad
&&\text{if $j=i-1$ or $i$}, 
\\
m(\Pi,\pi) \in &\{ 0,1 \} \qquad
&&\text{if $j=i-2$ or $i+1$}, 
\\
m(\Pi,\pi) =& 0 \qquad
&&\text{otherwise}.  
\end{alignat*}
\item[{\rm{(2)}}]
Suppose $j=i-1$ or $i$.  
Then $m(\Pi,\pi)=1$ 
 except for the countable set described as below.  
\begin{enumerate}
\item[{\rm{(a)}}]
Case $1 \le i \le n-1$.
Then $m(I_{\delta}(i,\lambda), J_{\varepsilon}(i,\nu))=2$
 if and only if 

\begin{alignat*}{3}
&j=i, 
&&\delta \varepsilon=+, 
&&(\lambda, \nu) \in L_{\operatorname{even}}\setminus \{\nu=0\}
                   \cup \{(i,i)\},
\\
&j=i, 
&&\delta \varepsilon=-, 
&&(\lambda, \nu) \in L_{\operatorname{odd}}\setminus \{\nu=0\}, 
\\
&j=i-1, \quad
&&\delta \varepsilon=+, \quad
&&(\lambda, \nu) \in L_{\operatorname{even}}\setminus \{\nu=0\} 
                   \cup \{(n-i,n-i)\}, 
\intertext{or}
&j=i-1, 
&&\delta \varepsilon=-, 
&&(\lambda, \nu) \in L_{\operatorname{odd}}\setminus \{\nu=0\}.   
\end{alignat*}
\item[{\rm{(b)}}]
Case $i=0$.  
Then 
$
   m(I_{\delta}(0,\lambda), J_{\varepsilon}(0,\nu))=2
$
if 
$
\delta \varepsilon=+, (\lambda, \nu) \in L_{\operatorname{even}}
$
 or 
$
\delta \varepsilon=-, (\lambda, \nu) \in L_{\operatorname{odd}}.
$   
\item[{\rm{(c)}}]
Case $i=n$.  
Then 
$
   m(I_{\delta}(n,\lambda), J_{\varepsilon}(n-1,\nu))=2
$
 if 

$\delta \varepsilon=+, (\lambda, \nu) \in L_{\operatorname{even}}
$
 or 
$
\delta \varepsilon=-, 
(\lambda, \nu) \in L_{\operatorname{odd}}.   
$
\end{enumerate}
\item[{\rm{(3)}}]
Suppose $j=i-2$ or $i+1$.  
Then $m(\Pi,\pi)=1$ 
 if one of the following conditions {\rm{(d)--(g)}}
 is satisfied, 
 and $m(\Pi,\pi)=0$ otherwise.  
\begin{enumerate}
\item[{\rm{(d)}}]
Case $j=i-2$, $2 \le i \le n-1$, $(\lambda, \nu)=(n-i,n-i+1)$, 
$\delta \varepsilon =-1$.  
\item[{\rm{(e)}}]
Case $(i,j)=(n,n-2)$, $-\lambda \in {\mathbb{N}}$, $\nu=1$, 
$\delta \varepsilon =(-1)^{\lambda+1}$.  
\item[{\rm{(f)}}]
Case $j=i+1$, $1 \le i \le n-2$, $(\lambda, \nu)=(i,i+1)$, 
$\delta \varepsilon =-1$.  
\item[{\rm{(g)}}]
Case $(i,j)=(0,1)$, $-\lambda \in {\mathbb{N}}$, $\nu=1$, 
$\delta \varepsilon =(-1)^{\lambda+1}$.  
\end{enumerate}
\end{enumerate}
\end{theorem}

The proof of Theorem \ref{thm:1.1} will be given right after Theorem \ref{thm:SBObasis}, 
by using Fact \ref{fact:2.7}
 and Theorems \ref{thm:VWSBO} and \ref{thm:161243}, 
 whose proofs are deferred at later chapters.  

\subsubsection{Classification of symmetry breaking operators
 $I_{\delta}(i,\lambda) \to J_{\varepsilon}(j,\nu)$}
\label{subsec:SBOthm}
In this subsection,
 we give explicit generators
 of
\[
   {\operatorname{Hom}}_{G'}
   (I_{\delta}(i,\lambda)|_{G'},J_{\varepsilon}(j,\nu)), 
\]
of which the dimension is determined in Theorem \ref{thm:1.1}.  
For most of the cases,
 the regular symmetry breaking operators
 $\Atbb {\lambda} {\nu} \pm {i,j}$
 and the differential symmetry breaking operators
 $\Ctbb \lambda \nu {i,j}$ give the generators.  
However,  for the exceptional discrete set classified in Theorem \ref{thm:161243}, 
 we need more operators
 which are defined as follows:
 for $(\lambda_0, \nu_0)\in {\mathbb{C}}^2$
 such that $\Atbb {\lambda_0} {\nu_0} \pm {i,j}=0$, 
 we renormalize the regular symmetry breaking operators
 $\Atbb \lambda \nu \pm {i,j}$
 as follows 
 (see Section \ref{subsec:161702}).  
For $j=i$ of $i-1$, 
 we set 
\index{A}{Ahttsln1@$\Attbb \lambda \nu {+} {i,j}$|textbf}
\index{A}{Ahttsln2@$\Attbb \lambda \nu {-} {i,j}$|textbf}
\begin{align}
\label{eqn:Aij+re}
\Attbb {\lambda_0} {\nu_0} {+} {i,j}
:=&
\lim_{\lambda \to \lambda_0}
\Gamma(\frac{\lambda-\nu_0}{2})
\Atbb \lambda {\nu_0} {+} {i,j},
\\
\label{eqn:Aij+re}
\Attbb {\lambda_0} {\nu_0} {-} {i, j}
:=&
\lim_{\lambda \to \lambda_0}
\Gamma(\frac{\lambda-\nu_0+1}{2})
\Atbb \lambda {\nu_0} {-} {i,j}.  
\end{align}
Then $\Attbb \lambda \nu {\pm} {i,j}$ are well-defined
 and nonzero symmetry breaking operators
(Theorem \ref{thm:170340}).

For $j \in \{i-1,i\}$ and $\gamma \in \{\pm\}$, 
 the set
\[
  \{
(\lambda, \nu)\in {\mathbb{C}}^2
:
  \Atbb \lambda \nu \gamma {i,j}=0
\}
\]
 is classified in Theorem \ref{thm:161243}.  
Then we are ready to give an explicit basis
 of symmetry breaking operators:
\begin{theorem}
[generators]
\label{thm:SBObasis}
Suppose $j=i$ or $i-1$.  
\begin{enumerate}
\item[{\rm{(1)}}]
$m(I_{\delta}(i,\lambda),J_{\varepsilon}(j,\nu))=1$
 if and only if
 $\Atbb \lambda \nu {\delta\varepsilon} {i,j} \ne 0$.  
In this case
\[
   \operatorname{Hom}_{G'}
  (I_{\delta}(i,\lambda)|_{G'},J_{\varepsilon}(j,\nu))
  =
  {\mathbb{C}}\Atbb \lambda \nu {\delta\varepsilon} {i,j}.  
\]
\item[{\rm{(2)}}]
$m(I_{\delta}(i,\lambda),J_{\varepsilon}(j,\nu))=2$
 if and only if
 $\Atbb \lambda \nu {\delta\varepsilon} {i,j} = 0$.  
In this case
\[
   \operatorname{Hom}_{G'}
  (I_{\delta}(i,\lambda)|_{G'},J_{\varepsilon}(j,\nu))
  =
  {\mathbb{C}}\Attbb \lambda \nu {\delta\varepsilon} {i,j}
  \oplus
  {\mathbb{C}} \Ctbb \lambda \nu {i,j}.  
\]
\end{enumerate}
\end{theorem}

See Theorem \ref{thm:161243} for the necessary 
 and sufficient condition on $(i,j,\lambda,\nu,\gamma)$
 for $\Atbb \lambda \nu \gamma {i,j}$ to vanish.  

\begin{remark}
For $j=i+1$ or $i-2$, 
all symmetry breaking operators
 are differential operators
 by the localness theorem
\index{B}{localnesstheorem@localness theorem}
 (Theorem \ref{thm:152347}), 
 and the generators are given in Fact \ref{fact:3.9}. 
\end{remark}

\begin{proof}
[Proof of Theorems \ref{thm:1.1} and \ref{thm:SBObasis}]
We apply the general scheme
 of symmetry breaking operators
 (Theorem \ref{thm:VWSBO})
 to the special setting:
\[
   V=\Exterior^i({\mathbb{C}}^n)
\quad
  \text{and}
\quad
   W=\Exterior^j({\mathbb{C}}^{n-1}).  
\]
Then the theorems follow from the explicit description
 of the zero sets
 of the (normalized) regular symmetry breaking operators
 $\Atbb \lambda \nu \gamma {i,j}$
 (Theorem \ref{thm:161243})
 and the classification
 of {\it{differential}} symmetry breaking operators
 (Fact \ref{fact:2.7}).  
\end{proof}

\begin{remark}
The first statement ({\it{i.e.}}, $\delta \varepsilon=+$ case)
 of Theorem \ref{thm:1.1} (2) (b)
 was established in \cite[Thm.~1.1]{sbon},
 and the second statement ({\it{i.e.}}, $\delta \varepsilon=-$ case) of (b)
 can be proved similarly.  
In this article,
 we take another approach for the latter case:
we deduce results
 for all the matrix-valued cases
 (including the scalar-valued case with $\delta \varepsilon=-$) from 
 the scalar valued case with $\delta \varepsilon=+$.  
\end{remark}

\subsection{Consequences of main theorems in Sections \ref{subsec:170228} and \ref{subsec:exhaust}}
\label{subsec:3.6}
In this section
 we discuss symmetry breaking from principal series representations
 $\Pi=I_{\delta}(V,\lambda)$ of $G$
 to $\pi=J_{\varepsilon}(W,\nu)$ of the subgroup $G'$
 in the case
 where $\Pi$ and $\pi$ are {\it{unitarizable}}.  
Unitary principal series representations
 are treated in Section \ref{subsec:SBtemp}, 
 and complementary series representations
 are treated
 in Sections \ref{subsec:SBcomp} and \ref{subsec:singcomp}.  
We note that $\Pi$ and $\pi$ are irreducible
 in these cases.  
On the other hand,
 if $\lambda$ (resp. $\nu$) is integral, 
 then $\Pi$ (resp. $\pi$) may be reducible.  
We shall discuss symmetry breaking operators
 for the subquotients 
 in the next chapter
 in detail
 when they have the trivial infinitesimal character $\rho$.

\subsubsection{Tempered representations}
\label{subsec:SBtemp}
We recall the concept of tempered unitary representations
 of locally compact groups.  
\begin{definition}
[tempered unitary representation]
\label{def:tempered}
A unitary representation of a unimodular group $G$ is called 
\index{B}{temperedrep@tempered representation|textbf}
{\it{tempered}}
 if it is weakly contained
 in the regular representations 
 on $L^2(G)$.  
By a little abuse of notation, 
 we also say the smooth representation 
 $\Pi^{\infty}$
 is {\it{tempered}}.  
\end{definition}
Returning to our setting 
 where $(G,G')=(O(n+1,1), O(n,1))$, 
 we see that the principal series representations $I_\delta (V,\lambda)$ and $J_\varepsilon (W,\nu)$ are tempered
 if and only if $\lambda \in \sqrt{-1}{\mathbb{R}}+\frac n 2$ and $\nu \in \sqrt{-1}{\mathbb{R}}+\frac 1 2 (n-1)$, 
 respectively.  
We refer to them as {\it{tempered principal series representations}}.

We recall 
$
   [V:W]=\dim_{\mathbb{C}} {\operatorname{Hom}}_{O(n-1)}(V|_{O(n-1)}, W).
$
Then Theorem \ref{thm:VWSBO} implies the following:
\begin{theorem}
[tempered principal series representations]
\label{thm:tempVW}
Let $(\sigma, V) \in \widehat {O(n)}$, 
 $(\tau, W) \in \widehat {O(n-1)}$, 
 $\delta, \varepsilon \in \{\pm\}$, 
 and $\lambda \in \sqrt{-1}{\mathbb{R}}+\frac n 2$, 
 $\nu \in \sqrt{-1}{\mathbb{R}}+\frac 1 2(n-1)$
 so that $I_\delta(V,\lambda)$ and $J_\varepsilon(W,\nu)$ are tempered principal series representations. 
Then the following four conditions are equivalent:
\begin{enumerate}
\item
[{\rm{(i)}}]
\index{A}{VWmult@$[V:W]$}
$[V:W] \ne 0;$
\item
[{\rm{(i$'$)}}]
$[V:W] = 1;$
\item
[{\rm{(ii)}}]
$
 {\operatorname{Hom}}_{G'}
 (I_\delta(V,\lambda)|_{G'}, J_\varepsilon(W,\nu))\not = \{0\};
$
\item
[{\rm{(ii$'$)}}]
$
\dim_{\mathbb{C}} {\operatorname{Hom}}_{G'}(I_\delta(V,\lambda)|_{G'}, J_\varepsilon(W,\nu)) =1.  
$
\end{enumerate}
\end{theorem}

Applying Theorem \ref{thm:tempVW} to the exterior tensor representations
 $V=\Exterior^i({\mathbb{C}}^n)$ of $O(n)$
 and $W=\Exterior^j({\mathbb{C}}^{n-1})$ of $O(n-1)$, 
 we get:
\begin{corollary}
\label{cor:tempered}
Suppose $\lambda \in \sqrt{-1}{\mathbb{R}}+ \frac n 2$, 
 and $\nu \in \sqrt{-1}{\mathbb{R}}+\frac 1 2(n-1)$.  
Then  
\[
   \dim_{\mathbb{C}} {\operatorname{Hom}}_{G'}(I_\delta(i,\lambda)|_{G'}, J_\varepsilon (j,\nu))
=
\begin{cases}
1 \quad&\text{if $i=j$ or $j=i-1$}, 
\\
0      &\text{otherwise}.
\end{cases}
\]
\end{corollary}

\subsubsection{Complementary series representations}
\label{subsec:SBcomp}
We say that $I_\delta(V,\lambda)$ is a (smooth)
\index{B}{complementaryseries@complementary series representation|textbf}
 {\it{complementary series  representation}}
 if it has a Hilbert completion
 to a unitary complementary series representation.  
If the irreducible $O(n)$-module $(\sigma,V)$ is of
\index{B}{typeX@type X, representation of ${O(N)}$\quad}
 type X
 (see Definition \ref{def:OSO}), 
 {\it{i.e.,}} 
 the last digit of the highest weight of $V$ is not zero, 
 then the principal series representation 
 $I_\delta(V,\lambda)$ is irreducible at $\lambda = \frac n 2$, 
 and consequently, 
 there exist complementary series representations
 $I_\delta(V,\lambda)$
 for some interval $\lambda \in (\frac n 2-a, \frac n 2+a)$ with $a >0$.  
\begin{example}
\label{ex:Iicompl}
Suppose $(\sigma,V)$ is the $i$-th exterior tensor representation
 $\Exterior^i({\mathbb{C}}^n)$.  
We assume that this representation is of type X, 
 equivalently,
 $n \ne 2i$
 (see Example \ref{ex:2.1}).  
The first reduction point of the principal series representation
 of $I_{\delta}(i,\lambda)$ is given by $\lambda=i$ or $n-i$
 (see Proposition \ref{prop:redIilmd}).  
Therefore $I_\delta(i,\lambda) \equiv I_\delta(\Exterior^i({\mathbb{C}}^n),\lambda)$
 is a complementary series representation
 if 
\[
  {\operatorname{min}}(i,n-i) < \lambda < {\operatorname{max}}(i,n-i).  
\]
\end{example}

In the category of unitary representations, 
 the restriction of a tempered representation of $G$
 to a reductive subgroup $G'$
 decomposes into the direct integral 
 of irreducible unitary tempered representations
 of a reductive subgroup $G'$
 because it is weakly contained 
 in the regular representation.  
In particular, 
 complementary series representations of the subgroup $G'$ do not appear
 in the {\it{unitary}} branching law
 of the restriction of a unitary tempered principal series  representation 
 $I_\delta(V,\lambda)$, 
 whereas Theorem \ref{thm:VWSBO}
 in the category of admissible {\it{smooth}} representations shows that there are  nontrivial symmetry breaking operators 
\[ 
  \Atbb \lambda \nu {\delta\varepsilon} {V,W} : I_\delta(V,\lambda) \rightarrow J_\varepsilon(W,\nu)
\]
 to all complementary series representations $J_\delta(W,\nu)$
 of the subgroup $G'$ if $[V:W] \not = 0.$

Moreover,
 Theorem \ref{thm:VWSBO} (2) implies also 
 that there are nontrivial symmetry breaking operators from
 any (smooth) complementary series representation
 $I_{\delta}(V,\lambda)$
 of $G$ to all (smooth) tempered principal series representations
 $J_{\varepsilon}(W,\nu)$ of the subgroup $G'$
 as far as $[V:W] \ne 0$.

\subsubsection{Singular complementary series representations}
\label{subsec:singcomp}
We consider the complementary series representations $I_\delta(i,s)$ 
 for $i < s< \frac n 2$ 
 with an additional assumption
 that $s$ is an integer.  
These representations are irreducible
 and have {\it{singular}} integral infinitesimal characters. 
We may describe the underlying $({\mathfrak{g}},K)$-modules
 of these singular complementary series representations
 in terms of cohomological parabolic induction 
\index{A}{Aqlmd@$A_{\mathfrak{q}}(\lambda)$}
 $A_{\mathfrak{q}}(\lambda)$
 where the parameter $\lambda$ wanders outside
\index{B}{goodrange@good range}
 the good range
 relative to the $\theta$-stable parabolic subalgebra ${\mathfrak{q}}$
 (see \cite[Def.~0.49]{KV} for the definition).

For $0 \le r \le \frac{n+1}{2}$, 
 we denote by
\index{A}{qi@${\mathfrak{q}}_i$, $\theta$-stable parabolic subalgebra}
 ${\mathfrak{q}}_r$ the $\theta$-stable parabolic subalgebra
 of ${\mathfrak{g}}_{\mathbb{C}}={\mathfrak{o}}(n+2,{\mathbb{C}})$
 with Levi factor $SO(2)^r \times O(n-2r+1,1)$
 in $G=O(n+1,1)$
 (see Definition \ref{def:qi}).  
\begin{lemma}
\label{lem:compAq}
Let $0 \le i \le [\frac n 2]-1$.  
For $s \in \{i+1,i+2,\cdots,[\frac n 2]\}$, 
 we have an isomorphism as $({\mathfrak{g}},K)$-modules:
\[
  I_+(i,s)_K \simeq A_{\mathfrak{q}_{i+1}}(0,\cdots,0,s-i).  
\]
\end{lemma}
See Remark \ref{rem:goodrange} in Appendix I
 for the normalization 
 of the $({\mathfrak{g}},K)$-module $A_{\mathfrak{q}}(\lambda)$
 and Theorem \ref{thm:compint} for more details
 about Lemma \ref{lem:compAq}.  
See also \cite[Thm.~3]{KMemoirs92} for some more general cases.  
The restriction of these representations
 to the special orthogonal group $SO(n+1,1)$ stays irreducible
 (see Lemma \ref{lem:170305} in Appendix II). 
Bergeron and Clozel proved
 that there are  automorphic square integrable representations,
 whose component at infinity is isomorphic
 to a representation $I_\delta(i,s)|_{SO(n+1,1)}$
 (see \cite{BC, BLS}). 

A special case of Theorem \ref{thm:1.1} includes:
\begin{proposition}
\label{prop:3.32}
Suppose $s \in {\mathbb{N}}$ and $i < s \le [\frac n 2]$.  
Let $\delta$, $\varepsilon \in \{\pm\}$.  
\newline\noindent
{\rm{(1)}}\enspace
For $i < r \le [\frac{n-1}2]$, 
\[ 
   {\operatorname{Hom}}_{G'}(I_{\delta}(i,s)|_{G'},J_{\varepsilon}(i,r))
 = {\mathbb{C}}.  
\]
\newline\noindent
{\rm{(2)}}\enspace
For $0 \le i-1 <r \le [\frac{n-1}{2}]$, 
\[  
   {\operatorname{Hom}}_{G'}(I_{\delta}(i,s)|_{G'},J_{\varepsilon}(i-1,r)) = {\mathbb{C}}.  
\]
\end{proposition}

\begin{remark}
Proposition \ref{prop:3.32} may be viewed as symmetry breaking operators from 
 the Casselman--Wallach globalization
 of the irreducible $({\mathfrak{g}},K)$-module $A_{\mathfrak{q}}(\lambda)$ to that of the irreducible $({\mathfrak{g}}',K')$-module $A_{\mathfrak{q}'}(\nu)$
 in some special cases 
 where both $\lambda$ and $\nu$ are 
 outside the good range of parameters
 relative to the $\theta$-stable parabolic subalgebras.  
\end{remark}

In the next chapter,
 we treat the case
 with trivial infinitesimal character $\rho$,
 and thus the parameters stay in the good range
 relative to the $\theta$-stable parabolic subalgebras.  
In particular,
 we shall determine a necessary and sufficient condition
 for a pair $({\mathfrak{q}}, {\mathfrak{q}}')$
 of $\theta$-stable parabolic subalgebras ${\mathfrak{q}}$ of ${\mathfrak{g}}_{\mathbb{C}}$
 and ${\mathfrak{q}}'$ of its subalgebra ${\mathfrak{g}}_{\mathbb{C}}'$
 such that 
\[
  {\operatorname{Hom}}_{G'}(\Pi|_{G'}, \pi) \ne \{0\}, 
\]
 when the underlying $({\mathfrak{g}},K)$-module $\Pi_K$
 of $\Pi \in {\operatorname{Irr}}(G)$
 is isomorphic to 
\index{A}{Aqlmdac@$(A_{{\mathfrak{q}}})_{\pm\pm}$}
 $(A_{\mathfrak{q}})_{\pm\pm}$ 
 and the underlying $({\mathfrak{g}}',K')$-module
 of $\pi \in {\operatorname{Irr}}(G')$
 is $(A_{\mathfrak{q}'})_{\pm \pm}$, 
 see Theorems \ref{thm:SBOvanish} and \ref{thm:SBOone}
 for the multiplicity-formula, 
 and Proposition \ref{prop:161655} in Appendix I
 for the description of $\Pi_K$ in terms of $(A_{\mathfrak{q}})_{\pm\pm}$.  
In contrast to the case of Proposition \ref{prop:3.32}, 
 the irreducible $G$-module $\Pi$
 and $G'$-module $\pi$ do not coincide with principal series representations,
 but appear as their subquotients in this case, 
 see Theorem \ref{thm:LNM20} (1).

\subsection{Actions
 of $(G/G_0)
\hspace{1mm}
{\widehat{}}
\, 
\times (G'/G_0')
\hspace{1MM}
{\widehat{}}
\,$ on symmetry breaking operators}
\label{subsec:actPont}

In this section 
 we discuss the action of the character group of $G \times G'$
 on the set 
\[
  \{ {\operatorname{Hom}}_{G'}(\Pi|_{G'}, \pi) \}
\]
 of the spaces of symmetry breaking operators 
 where admissible smooth representations $\Pi$ of $G$
 and those $\pi$ of the subgroup $G'$ vary.  
Actual computations for the pair
 $(G,G')=(O(n+1,1),O(n,1))$ are carried out
 by using Lemma \ref{lem:IVchi}
 for principal series representations
 and Theorem \ref{thm:LNM20} (5) for their irreducible subquotients.  
\subsubsection{Generalities: The action of character group of 
 $G \times G'$ on $\{{\operatorname{Hom}}_{G'}(\Pi|_{G'}, \pi)\}$
 in the general case}
Let $G \supset G'$ be a pair of real reductive Lie groups.  
Then the character group of $G \times G'$
 acts on the set of vector spaces
 $\{{\operatorname{Hom}}_{G'}(\Pi|_{G'}, \pi)\}$
 where $\Pi$ runs over admissible smooth representations of $G$, 
 and $\pi$ runs over those of the subgroup $G'$.  
Here the action is given by 
\[
   {\operatorname{Hom}}_{G'}(\Pi|_{G'}, \pi)
   \mapsto
   {\operatorname{Hom}}_{G'}((\Pi \otimes \chi^{-1})|_{G'}, \pi \otimes \chi')
\]
for a character $\chi$ of $G$
 and $\chi'$ of the subgroup $G'$.

In what follows, 
 we regard a character of $G$ as a character of $G'$
 by restriction,
 and use the same letter
 to denote its restriction to the subgroup $G'$.  
Then for all characters $\chi$ and $\chi'$ of $G$, 
 we have the following isomorphisms:
\begin{alignat}{2} 
\mbox{Hom}_{G'} ((\Pi\otimes \chi)|_{G'}, \pi \otimes \chi' ) 
&\simeq\,\,
&&\mbox{Hom}_{G'}(\Pi|_{G'},\pi\otimes \chi^{-1}  \otimes \chi') 
\notag
\\
&\simeq 
&&\mbox{Hom}_{G'} ((\Pi\otimes (\chi')^{-1})|_{G'}, \pi \otimes \chi^{-1}) 
\notag
\\
&\simeq 
&&\mbox{Hom}_{G'} ((\Pi\otimes \chi \otimes (\chi')^{-1})|_{G'}, \pi).
\label{eqn:314}
\end{alignat}

The above isomorphisms define an equivalence relation
 on the set 
\[
   \{{\operatorname{Hom}}_{G'}(\Pi|_{G'}, \pi)\}
\]
 of the spaces of symmetry breaking operators
 where $\Pi$ and $\pi$ vary.

\subsubsection{Actions of the character group
 of the component group
 on 
 $\{{\operatorname{Hom}}_{G'}
  (I_{\delta}(i,\lambda)|_{G'}, J_{\varepsilon}(j,\nu))\}$}
\label{subsec:chiIilmd}
We apply the above idea to our setting
\[
  (G,G') = (O(n+1,1), O(n,1)).  
\]
Then the 
\index{B}{componentgroup@component group $G/G_0$}
component groups of $G$ and $G'$ are a finite abelian group given by
\begin{equation}
\label{eqn:chiabrest}
   G'/G_0' \simeq G/G_0 
         \simeq {\mathbb{Z}}/2{\mathbb{Z}} \times {\mathbb{Z}}/2{\mathbb{Z}}.
\end{equation}

We recall from \eqref{eqn:chiab}
 that the set of their one-dimensional representations
 is parametrized by 
\[
(G'/G_0')\hspace{-1mm}{\widehat{\hphantom{m}}} 
\simeq (G/G_0)\hspace{-1mm}{\widehat{\hphantom{m}}} 
           = \{\chi_{a b}: a, b \in \{\pm\}\}.  
\]
By abuse of notation,
 we shall use the same letters $\chi_{a b}$
 to denote the corresponding one-dimensional representations
 of $G$, $G'$, $G/G_0$, and $G'/G_0'$.

The action of the character group 
 ({\it{Pontrjagin dual}})
 $(G/G_0)\hspace{-1mm}{\widehat{\hphantom{m}}}$
 on the set of principal series representations
 can be computed
 by using Lemma \ref{lem:IVchi}.  
To describe the action of the Pontrjagin dual 
 $(G/G_0)\hspace{-1mm}{\widehat{\hphantom{m}}} \simeq
  (G'/G_0')\hspace{-1mm}{\widehat{\hphantom{m}}}$
 on the parameter set of the principal series representations $I_{\delta}(i,\lambda)$ of $G$
 and $J_{\varepsilon}(j,\nu)$ of the subgroup $G'$, 
 we define
\begin{alignat*}{3}
   &S:=\{0,1,\cdots,n\} \times {\mathbb{C}} \times {\mathbb{Z}}/2{\mathbb{Z}}, 
 \quad
   &&I(s):=I_{\delta}(i,\lambda)
\quad
  &&\text{for } s =(i,\lambda,\delta)\in S, 
\\
   &T:=\{0,1,\cdots,n-1\} \times {\mathbb{C}} \times {\mathbb{Z}}/2{\mathbb{Z}}, \quad
   &&J(t):=J_{\varepsilon}(j,\nu)
\quad
  &&\text{for } t =(j,\nu,\varepsilon)\in T.   
\end{alignat*}
We let the character group $(G/G_0)\hspace{-1mm}{\widehat{\hphantom{m}}}$ act on $S$
 by the following formula:
\begin{alignat*}{3}
\chi_{++} \cdot (i,\lambda,\delta):=&(i,\lambda,\delta),
\quad
&&\chi_{+-} \cdot (i,\lambda,\delta):=&&(i,\lambda,-\delta),
\\
\chi_{-+} \cdot (i,\lambda,\delta):=&(\tilde i,\lambda,-\delta),
\quad
&&\chi_{--} \cdot (i,\lambda,\delta):=&&(\tilde i,\lambda,\delta),
\end{alignat*}
where $\tilde i:=n-i$.  
The action of $(G'/G_0')\hspace{-1mm}{\widehat{\hphantom{m}}}$ on the set $T$
 is defined similarly,
 with obvious modification
\[
\tilde j:=n-1-j
\]
 when we discuss representations
 of the subgroup $G'=O(n,1)$.  
By Lemma \ref{lem:IVchi}
 and by the $O(n)$-isomorphism
 $\Exterior^i({\mathbb{C}}^n) \simeq \Exterior^{n-i}({\mathbb{C}}^n) \otimes \det$, 
 we obtain the following.  
\begin{lemma}
\label{lem:LNM27}
For all $\chi \in (G/G_0)\hspace{-1mm}{\widehat{\hphantom{m}}}\simeq (G'/G_0')\hspace{-1mm}{\widehat{\hphantom{m}}}$
 and for $s \in S$, $t \in T$, 
we have the following isomorphisms as $G$-modules
 and $G'$-modules, 
respectively:
\begin{align*}
   I(s) \otimes \chi &\simeq I(\chi \cdot s), 
\\
   J(t) \otimes \chi &\simeq J(\chi \cdot t).  
\end{align*}
\end{lemma}

Then the equivalence defined by the isomorphisms
 \eqref{eqn:314} implies
 that it suffices to consider symmetry breaking operators
 for $(\delta, \varepsilon) =(+,+)$
 and $(\delta, \varepsilon) =(+,-)$.
To be more precise,
 we obtain the following.  

\begin{proposition} 
Let $\lambda, \nu \in {\mathbb{C}}$.  
Then every symmetry breaking operator in 
\[
     \bigcup_{\delta,\varepsilon \in \{\pm\} }
     \bigcup_{0 \le i \le n} \bigcup_{0 \le j \le n-1}
     \operatorname{Hom}_{G'}(I_{\delta}(i,\lambda)|_{G'},J_{\varepsilon}(j,\nu))
\]
is equivalent to a symmetry breaking  operator in
\[ 
 \bigcup_{0 \le i \le [\frac n2]} \bigcup_{0 \le j \le n-1}
 \left(\operatorname{Hom}_{G'}(I_+(i,\lambda)|_{G'}, J_+(j,\nu)) 
 \cup
 \operatorname{Hom}_{G'}(I_+(i,\lambda)|_{G'}, J_-(j,\nu))
 \right).  
\]
\end{proposition}
\proof
We use a graph to prove this.
We set
\begin{eqnarray*} 
(\delta,\varepsilon) &:= & \mbox{Hom}_{G'}(I_{\delta} (i, \lambda)|_{G'}, J_{\varepsilon}(j,\nu)), 
\\
\begin{pmatrix} \delta \\ \varepsilon \end{pmatrix} 
&:=& \operatorname{Hom}_{G'}(I_{\delta} (n-i, \lambda)|_{G'}, J_{\varepsilon}(n-j-1,\nu)).  
\end{eqnarray*}
In the following graph the nodes are indexed by $(\delta,\varepsilon)$
 in first row and by
 $\begin{pmatrix} \delta \\ \varepsilon \end{pmatrix}$
 in the second row.  
The nodes are connected by a line if they are equivalent.  
By Lemma \ref{lem:LNM27}, 
 we obtain the graph by taking 
\index{A}{1chipm@$\chi_{+-}$}
$\chi=\chi'=\chi_{+-}$
 in \eqref{eqn:314}
 for horizontal equivalence,
 and 
\index{A}{1chipmm@$\chi_{--}=\det$}
$\chi=\chi'=\chi_{-+}$ in \eqref{eqn:314}
 for crossing equivalence
 (we omit here lines in the graph
 corresponding to 
 $\chi=\chi'=\chi_{--}$ in \eqref{eqn:314}
 for vertical equivalence):

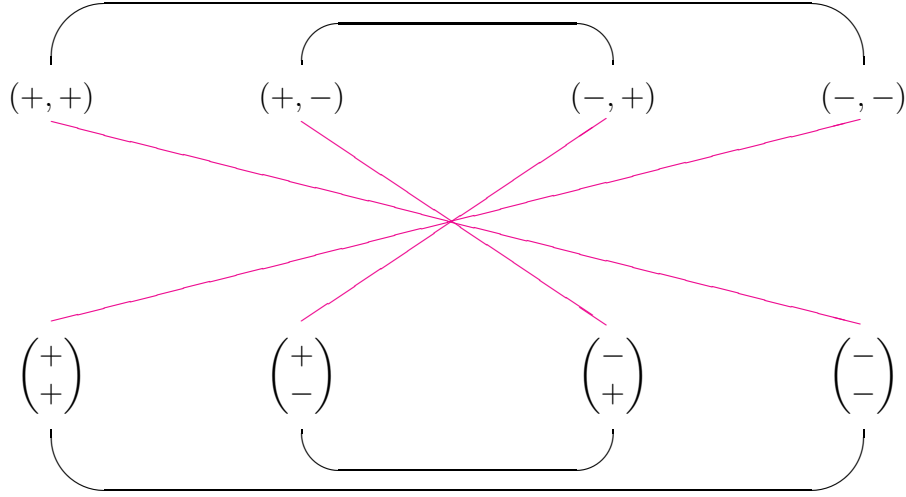
\begin{center}
\unitlength.9cm
\begin{picture}(14, 8.4)
\put(0,6){\makebox(2,1){$(+,+)$}}
\put(3.7,6){\makebox(2,1){$(+,-)$}}
\put(8.3,6){\makebox(2,1){$(-,+)$}}
\put(12,6){\makebox(2,1){$(-,-)$}}

\put(7,7){\oval(12,1.8)[t]}
\put(7,7){\oval(4.6,1.2)[t]}

\put(0,1.9){\makebox(2,1){$\begin{pmatrix}+\\+ \end{pmatrix}$}}
\put(3.7,1.9){\makebox(2,1){$\begin{pmatrix}+\\- \end{pmatrix}$}}
\put(8.3,1.9){\makebox(2,1){$\begin{pmatrix}-\\+ \end{pmatrix}$}}
\put(12,1.9){\makebox(2,1){$\begin{pmatrix}-\\- \end{pmatrix}$}}

\put(7,1.6){\oval(12,1.8)[b]}
\put(7,1.6){\oval(4.6,1.2)[b]}

\color{magenta}
\put(1,3.2){\line(4,1){11.95}}
\put(1,6.15){\line(4,-1){11.95}}
\put(4.7,3.2){\line(3,2){4.5}}
\put(4.7,6.15){\line(3,-2){4.5}}
\end{picture}
\color{black}

\end{center}

We observe that there are exactly two connected components of the graph, 
  and that 
$\mbox{Hom}_{G'}(I_+(i,\lambda)|_{G'}, J_+(j,\nu))$ and
$\mbox{Hom}_{G'}(I_+(i,\lambda)|_{G'}, J_-(j,\nu))$ are
 in a different connected component.  
Moreover,
 we may choose $i$ or $n-i$
 in the same equivalence classes, 
 and thus we may take $0 \le i \le \frac n 2$
 as a representative.  
\qed

\begin{example}
\label{ex:tempdual}
\begin{enumerate}
\item[{\rm{(1)}}]
Suppose $n=2m$ and $i=m$.  
Applying the isomorphism \eqref{eqn:314} to $(\Pi,\pi)=(I_{\delta}(m,\lambda), J_{\varepsilon}(m,\nu))$
 with $\chi=\chi'=\chi_{--}$, 
 we obtain a natural bijection:
\[
   {\operatorname{Hom}}_{G'}
    (I_{\delta}(m,\lambda)|_{G'}, J_{\varepsilon}(m,\nu))
    \simeq
    {\operatorname{Hom}}_{G'}
    (I_{\delta}(m,\lambda)|_{G'}, J_{\varepsilon}(m-1,\nu)).  
\]
We note that the $G$-module $I_{\delta}(m,\lambda)$
 at $\lambda=m$ splits into the direct sum
 of two irreducible smooth tempered representations
 (Theorem \ref{thm:LNM20} (1) and (8)).  
\item[{\rm{(2)}}]
Suppose $n=2m+1$ and $i=m$.  
Similarly to the first statement,
 we have a natural bijection:
\[
   {\operatorname{Hom}}_{G'}
    (I_{\delta}(m,\lambda)|_{G'}, J_{\varepsilon}(m,\nu))
    \simeq
    {\operatorname{Hom}}_{G'}
    (I_{\delta}(m+1,\lambda)|_{G'}, J_{\varepsilon}(m,\nu)).  
\]
In this case,
 the $G'$-module $J_{\varepsilon}(m,\nu)$ at $\nu=m$
 splits into the direct sum
 of two irreducible smooth tempered representations.  
\end{enumerate}
\end{example}

\subsubsection{Actions of characters of the component group 
 on ${\operatorname{Hom}}_{G'}(\Pi_{i,\delta}|_{G'}, \pi_{j,\varepsilon})$}
In the next chapter,
 we discuss 
\[
 {\operatorname{Hom}}_{G'}(\Pi|_{G'}, \pi)
\]
 for $\Pi \in {\operatorname{Irr}}(G)_{\rho}$
 and $\pi \in {\operatorname{Irr}}(G')_{\rho}$.  
In this case,
 \eqref{eqn:314} implies the following:
\begin{proposition}
[duality for symmetry breaking operators]
\label{prop:SBOdual}
There are natural isomorphisms
\begin{align*}
   {\operatorname{Hom}}_{G'}(\Pi_{i,\delta}|_{G'}, \pi_{j,\varepsilon})
   \,\,\simeq\,\,&
   {\operatorname{Hom}}_{G'}(\Pi_{n+1-i,\delta}|_{G'}, \pi_{n-j,\varepsilon})
\\
   \,\,\simeq\,\,&
   {\operatorname{Hom}}_{G'}(\Pi_{i,-\delta}|_{G'}, \pi_{j,-\varepsilon})
\\
   \,\,\simeq\,\,&
{\operatorname{Hom}}_{G'}(\Pi_{n+1-i,-\delta}|_{G'}, \pi_{n-j,-\varepsilon}).  
\end{align*}
\end{proposition}
\begin{proof}
By Theorem \ref{thm:LNM20} (5),  
 we have a natural $G$-isomorphism
\index{A}{1chipmpm@$\chi_{\pm\pm}$, one-dimensional representation of $O(n+1,1)$}
$\Pi_{i,\delta} \otimes \chi_{-+} \simeq \Pi_{n+1-i,\delta}$
 and a $G'$-isomorphism
 $\pi_{j,\varepsilon} \otimes \chi_{-+} \simeq \pi_{n-j,\varepsilon}$.  
Hence the first isomorphism is derived from \eqref{eqn:314}.  
By taking the tensor product with $\chi_{+-}$, 
 we get the last two isomorphisms
 again by Theorem \ref{thm:LNM20} (5).  
\end{proof}

\newpage
\section
{Symmetry breaking for irreducible representations
 with infinitesimal character $\rho$.}
\label{sec:SBOrho}
\index{B}{infinitesimalcharacter@infinitesimal character}

In this chapter, 
 we focus on symmetry breaking operators from 
 {\it{irreducible}} representations $\Pi$ of $G=O(n+1,1)$
 with 
\index{B}{infinitesimalcharacter@infinitesimal character}
${\mathfrak{Z}}_G({\mathfrak{g}})$-infinitesimal character
\index{A}{1parho@$\rho_G$}
 $\rho_G$ 
 to {\it{irreducible}} representations $\pi$ of the subgroup $G'=O(n,1)$
 with ${\mathfrak{Z}}_{G'}({\mathfrak{g}}')$-infinitesimal character
 $\rho_{G'}$.  
The main results are Theorems \ref{thm:SBOvanish} and \ref{thm:SBOone}, 
 where we determine the multiplicity 
 $\dim_{\mathbb{C}}{\operatorname{Hom}}_{G'}(\Pi|_{G'}, \pi)$
 for all pairs $(\Pi,\pi)$.  
A diagrammatic formulation 
 of the main results is given in Theorem \ref{thm:SBOfg}.

The proof uses basic properties
 of the normalized symmetry breaking operators
 for principal series representations
 of $G$ and $G'$, 
\begin{equation*}
\Atbb \lambda \nu {\delta \varepsilon}{i,j} \colon 
I_{\delta}(i,\lambda)  \rightarrow J_{\varepsilon}(j,\nu), 
\end{equation*}
in particular, 
 the $(K,K')$-spectrum 
 on basic $K$-types (Theorem \ref{thm:153315})
 and  their functional equations
 (Theorems \ref{thm:TAA} and \ref{thm:ATA}).

\subsection{Main Theorems} 
We recall from Theorem \ref{thm:LNM20}
 that irreducible admissible smooth representations of $G$ 
 with trivial
\index{B}{infinitesimalcharacter@infinitesimal character}
 ${\mathfrak {Z}}_G({\mathfrak {g}})$-infinitesimal character
 $\rho_G$
 are classified as
\index{A}{IrrGrho@${\operatorname{Irr}}(G)_{\rho}$,
 set of irreducible admissible smooth representation of $G$
 with trivial infinitesimal character $\rho$\quad}
\[
    {\operatorname{Irr}}(G)_{\rho}
   =
  \{
    \Pi_{i, \delta}
    :
    0 \le i \le n+1, \, \delta = \pm
\}.  
\]
Similarly,
 irreducible admissible smooth representations
 of the subgroup $G'=O(n,1)$
 with trivial ${\mathfrak {Z}}_{G'}({\mathfrak {g}}')$-infinitesimal character
 $\rho_{G'}$
 are classified as 
\[
  {\operatorname{Irr}}(G')_{\rho}
  =
  \{
  \pi_{j, \varepsilon}
    :
    0 \le j \le n, \, \varepsilon = \pm
\},   
\]
where we have used lowercase letters $\pi$
 for the subgroup $G'$ instead of $\Pi$.  
We also recall that the representation $\Pi_{i, \delta}$ of $G=O(n+1,1)$
 is 

\begin{enumerate}
\item[$\bullet$]
one-dimensional 
 if and only if $i=0$ or $n+1$;

\item[$\bullet$]
 the smooth representation of a discrete series representation 
 if $i=\frac{n+1}{2}$
 ($n$: odd);

\item[$\bullet$]
 that of a tempered representation
 if $i=\frac n 2$
 ($n$: even).  
\end{enumerate}

The following two theorems determine
 the dimension of
\[
  {\operatorname{Hom}}_{G'}(\Pi|_{G'},\pi)
\quad
\text{for $\Pi \in {\operatorname{Irr}}(G)_{\rho}$
      and $\pi \in {\operatorname{Irr}}(G')_{\rho}$.
}
\]

\begin{theorem}
[vanishing]
\label{thm:SBOvanish}
\index{B}{vanishingtheorem@vanishing theorem}
Suppose $0 \le i \le n+1$, $0 \le j \le n$, 
 $\delta, \varepsilon \in \{\pm\}$.  
\begin{enumerate}
\item[{\rm{(1)}}]
If $j \not = i, i-1$ then
$
  {\operatorname{Hom}}_{G'}(\Pi_{i,\delta}|_{G'}, \pi_{j,\varepsilon})=\{0\}. 
$

\item[{\rm{(2)}}] 
If $\delta \varepsilon =-$, 
 then 
$ 
{\operatorname{Hom}}_{G'}(\Pi_{i,\delta}|_{G'}, \pi_{j,\varepsilon}) =\{0\}.  
$ 
\end{enumerate}
\end{theorem}

\begin{theorem} 
[multiplicity-one]
\label{thm:SBOone}
\index{B}{multiplicityonetheorem@multiplicity-one theorem}
Suppose $0 \le i \le n+1$, $0 \le j \le n$
 and $\delta, \varepsilon \in \{\pm\}$.  
If $j=i-1$ or $i$
 and if $\delta \varepsilon =+$, 
 then 
\[
   \dim_{\mathbb{C}}
   {\operatorname{Hom}}_{G'}
   (\Pi_{i,\delta}|_{G'}, \pi_{j,\varepsilon}) =1.  
\]
\end{theorem}

The proof of Theorems \ref{thm:SBOvanish} and \ref{thm:SBOone}
 will be given in Chapter \ref{sec:pfSBrho}.  
The nonzero symmetry breaking operators from 
 $\Pi_{i,+}$ to $\pi_{j,+}$ ($j \in \{i-1,i\}$) will be applied 
 to construct
\index{B}{period@period}
 periods 
 in Chapter \ref{sec:period}
 (see Theorem \ref{thm:171517} for example).

\subsection{Graphic description of the multiplicity
 for irreducible representations
 with infinitesimal character $\rho$}

Using the action of the Pontrjagin dual 
 of the component group 
 $(G/G_0)\hspace{1mm}{\widehat{}}\, \times (G'/G_0')\hspace{1MM}{\widehat{}}\,\,$
 on ${\operatorname{Hom}}_{G'}(\Pi_{i,\delta}|_{G'}, \pi_{j,\varepsilon})$,
 see Proposition \ref{prop:SBOdual}, 
 we see that Theorems \ref{thm:SBOvanish} and \ref{thm:SBOone}
 are equivalent to their special case
 where $i \le \frac{n+1}{2}$ and $\delta=+$.  
Furthermore, 
 taking the vanishing result (Theorem \ref{thm:SBOvanish}) into account, 
 we focus on the case $j \le \frac n 2$ and $\varepsilon=+$.  
We then describe Theorems \ref{thm:SBOvanish} and \ref{thm:SBOone}
 graphically in this setting.  

We suppress the subscript, 
 and write $\Pi_i$ for $\Pi_{i,+}$, 
 and $\pi_j$ for $\pi_{j,+}$.  
Then $\Pi_i$ ($0 \le i \le \frac{n+1}{2}$)
 and $\pi_j$ ($0 \le j \le \frac{n}{2}$) are the standard sequence of representations with infinitesimal character $\rho$ of $G$,
 respectively $G'$ starting with the trivial one-dimensional representation
(Definition \ref{def:Pii}). 
In the diagrams below, 
 the first row are representations of $G$,
 the second row are representations of the subgroup $G'$. 
Arrows mean
 that there exist nonzero symmetry breaking operators.

\begin{theorem} 
\label{thm:SBOfg}
Symmetry breaking for the standard sequence
\index{B}{standardsequence@standard sequence}
 of irreducible representations
 starting at the trivial one-dimensional representations
 are represented graphically
 in Diagrams \ref{fig:Hasse1} and \ref{fig:Hasse2}. 

\medskip
\begin{figure}[htp]
\caption{Symmetry breaking for $O(2m+1,1)\downarrow O(2m,1)$}
\begin{center}
\begin{tabular}{c@{~}c@{~}c@{~}c@{~}c@{~}c@{~}c@{~}c@{~}c}
$\Pi_0$ & & $\Pi_1$ & &\dots & & $\Pi_{m-1}$ & & $\Pi_{m}$ 
\\
$\downarrow$ & $\swarrow$ & $\downarrow$ & $\swarrow$ & & $\swarrow$ & $ \downarrow $
&  $\swarrow $  &  $\downarrow$ 
\\
$\pi_0$ & & $\pi_1$ & &\dots & & $\pi_{m-1}$ & & $\pi_{m}$ 
\end{tabular}
\end{center}
\label{fig:Hasse1}
\end{figure}%

\medskip
\begin{figure}[htp]
\caption{Symmetry breaking for $O(2m+2,1) \downarrow O(2m+1,1)$ }
\begin{center}
\begin{tabular}{@{}c@{~}c@{~}c@{~}c@{~}c@{~}c@{~}c@{~}c@{~}c@{~}c@{~}c@{~}c@{}}
$\Pi_0$& &$\Pi_1$& &\dots & & $\Pi_{m-1} $& & $\Pi_{m}$ & &$\Pi_{m+1}$
\\
$\downarrow$ & $\swarrow$ & $\downarrow $ & $\swarrow$ & & $\swarrow$ & $ \downarrow $& $\swarrow $ & $\downarrow$ & $\swarrow$ &
\\
$\pi_0$& &$\pi_1$& &\dots & & $\pi_{m-1}$ & & $\pi_{m}$ &
\end{tabular}
\end{center}
\label{fig:Hasse2}
\end{figure}%
\end{theorem}

\newpage

{}~~~~~{}
\newpage
\section{Regular symmetry breaking operators}
\label{sec:section7}
\index {B}{regularsymmetrybreakingoperator@regular symmetry breaking operator}

Let 
\index{A}{IdeltaV@$I_{\delta}(V, \lambda)$}
$
I_{\delta}(V,\lambda)
$
 be a principal series representation
 of $G=O(n+1,1)$
 realized in the Fr{\'e}chet space 
 $C^{\infty}(G/P, {\mathcal{V}}_{\lambda,\delta})$, 
 and 
\index{A}{JWepsilon@$J_{\varepsilon}(W, \nu)$}
$
J_{\varepsilon}(W,\nu)
$
 that of $G'=O(n,1)$
 realized in $C^{\infty}(G'/P', {\mathcal{W}}_{\nu,\varepsilon})$
 as in Section \ref{subsec:smoothI}.  
In this chapter 
 we apply the general result in \cite[Chap.~3]{sbon} to construct a 
 \lq\lq{matrix-valued regular symmetry breaking operators}\rq\rq\
\index{A}{Ahtg0@$\Atbb \lambda \nu {\pm} {V,W}$}
$
   \Atbb \lambda \nu {\pm} {V,W}
   \colon
   I_{\delta}(V,\lambda) \to J_{\pm \delta}(W,\nu)
$
 that depend holomorphically on $(\lambda, \nu) \in {\mathbb{C}}^2$. 
We shall prove
 that the normalization \eqref{eqn:KVWpt} and \eqref{eqn:KVWmt}
 is optimal
 in the sense
 that the zeros of the operator
 $\Atbb \lambda \nu \pm {V,W}$
 are of codimension $>1$ in the parameter space of $(\lambda, \nu)$, 
 that is, 
 discrete in ${\mathbb{C}}^2$ in our setting.  
A key idea of the proof is a reduction to the scalar case.   

\subsection{Generalities}
We recall from the general theory \cite[Chap.~3]{sbon}
 on the distribution kernels 
 of symmetry breaking operators, 
 which will be the basic tool in this chapter.  
Furthermore,
 we discuss some subtle questions
 on the underlying topology 
 of representation spaces for symmetry breaking, 
 see Theorem \ref{thm:SBOBF}.  

\subsubsection{Distribution kernels of symmetry breaking operators}

Throughout this monograph,
 we shall regard distributions
 as the dual of compactly supported smooth densities
 rather than that of compactly supported smooth functions.
Thus we treat distributions as \lq\lq{generalized functions}\rq\rq, 
 and write their pairing 
 with test functions by using the integral symbol,  
 as if they were ordinary functions
 (with densities).

Let $G \supset G'$ be a pair of real reductive Lie groups, 
 and $P$, $P'$ their parabolic subgroups.  
We do not require an inclusive relation $P \supset P'$ in this subsection.  
Let $(\widetilde \sigma,V)$ be a finite-dimensional representation of $P$, 
 and $(\widetilde \tau,W)$ that of the subgroup $P'$.  
We form homogeneous vector bundles
 over flag manifolds
 by 
\begin{alignat*}{2}
{\mathcal{V}}:=& G \times_P V &&\to G/P, 
\\
{\mathcal{W}}:=& G' \times_{P'} W &&\to G'/P'.   
\end{alignat*}

We write ${\operatorname{Ind}}_P^G(\widetilde \sigma)$
 for the admissible smooth representation 
 of $G$ on the Fr{\'e}chet space $C^{\infty}(G/P,{\mathcal{V}})$, 
 and ${\operatorname{Ind}}_{P'}^{G'}(\widetilde \tau)$
 for that of the subgroup of $G'$ on $C^{\infty}(G'/P',{\mathcal{W}})$.

We denote by ${\mathcal{V}}^{\ast}$ 
 the dualizing bundle of ${\mathcal{V}}$, 
 which is a $G$-homogeneous vector bundle over $G/P$ 
 associated to the representation
\[
   V^{\ast}:=V^{\vee} \otimes |\det ({\operatorname{Ad}}_{{\mathfrak{g}}/{\mathfrak{p}}})|^{-1}
\]
 of the group $P$, 
where $V^{\vee}$ denotes the contragredient representation
 of $(\widetilde \sigma,V)$.  
Then the regular representation of $G$
 on the space ${\mathcal{D}}'(G/P, {\mathcal{V}}^{\ast})$
 of ${\mathcal{V}}^{\ast}$-valued distribution sections
 is the dual of the representation on $C^{\infty}(G/P, {\mathcal{V}})$.

The Schwartz kernel theorem guarantees 
 that any symmetry breaking operator can be expressed
 by using a distribution kernel.  
Conversely,
 distributions that give rise to symmetry breaking operators
 are characterized as follows.  
\begin{fact}
[{\cite[Prop.~3.2]{sbon}}]
\label{fact:SBOdistr}
There are natural linear bijections:
\[
{\operatorname{Hom}}_{G'}
 (C^{\infty}(G/P, {\mathcal{V}})|_{G'}, 
  C^{\infty}(G'/P', {\mathcal{W}}))
\simeq\,  {\mathcal{D}}'(G/P \times G'/P', {\mathcal{V}}^{\ast} \boxtimes {\mathcal{W}}
)^{\Delta(G')}.  
\]
Here ${\mathcal{V}}^{\ast} \boxtimes {\mathcal{W}}$ denotes
 the outer tensor product bundle 
 over the direct product manifold $G/P \times G'/P'$.  
\end{fact}
We note that the multiplication map
\[
  m \colon G \times G' \to G, 
\quad
  (x,y) \mapsto y^{-1} x
\]
 induces a linear bijection
\[
{\mathcal{D}}'(G/P \times G'/P',
  {\mathcal{V}}^{\ast} \boxtimes {\mathcal{W}})
^{\Delta(G')}
\underset{m^{\ast}}{\overset \sim \leftarrow}
({\mathcal{D}}'(G/P, {\mathcal{V}}^{\ast}) \otimes W)^{\Delta(P')}, 
\]
where
 the right-hand side stands for the space
 of $P'$-invariant vectors
 under the diagonal action
 on the tensor product of the $G$-module
 ${\mathcal{D}}'(G/P, 
  {\mathcal{V}}^{\ast})$
 and the $P'$-module $W$.

Thus Fact \ref{fact:SBOdistr} may be reformulated
 as the following linear bijection
\begin{equation}
\label{eqn:SBOdistr}
   {\operatorname{Hom}}_{G'}
   (C^{\infty}(G/P, {\mathcal{V}})|_{G'}, C^{\infty}(G'/P', {\mathcal{W}}))
   \simeq
   ({\mathcal{D}}'(G/P, {\mathcal{V}}^{\ast}) \otimes W)^{\Delta(P')}.  
\end{equation}

The point of Fact \ref{fact:SBOdistr} is 
 that the map
\[
   C^{\infty}(G/P, {\mathcal{V}}) 
  \to {\mathcal{D}}'(G'/P', {\mathcal{W}}), 
\qquad
   f \mapsto \int_X K(x,y) f(x)
\]
to the space ${\mathcal{D}}'(G'/P', {\mathcal{W}})$
 of {\it{distribution}} sections
 becomes automatically a continuous map
 to the space $C^{\infty}(G'/P', {\mathcal{W}})$
 of {\it{smooth}} sections 
 for any $K\in {\mathcal{D}}'(G/P \times G'/P', {\mathcal{V}}^{\ast} \boxtimes {\mathcal{W}})^{\Delta(G')}$.  
This observation leads us to the proof of the isomorphism 
 \eqref{eqn:SBOBF} in Theorem \ref{thm:SBOBF}.

\subsubsection{Invariant bilinear forms
 on admissible smooth representations
 and symmetry breaking operators}
\label{subsec:BFSBO}
We retain the setting
 of the previous subsection,
 in particular,
 we suppose
 that $G \supset G'$
 are a pair of real reductive Lie groups.

Let $(\Pi,U)$ and $(\pi,U')$ be admissible smooth representations
 of $G$ and $G'$, 
 respectively.  
We recall that the underlying topological vector space
 of any admissible smooth representation
 is a nuclear Fr{\'e}chet space.  
We define $\Pi \boxtimes \pi$ 
 to be the natural representation
 of the direct product group $G \times G'$
 on the space $U \widehat \otimes U'$.  
In this subsection,
 we study the space
 ${\operatorname{Hom}}_{G'}(\Pi \boxtimes \pi,{\mathbb{C}})$
 of continuous functionals
 that are invariant under the diagonal action
 of the subgroup $G'$.

For an admissible smooth representation $(\Pi,U)$ of $G$, 
 we denote by $\Pi^{\vee}$ the contragredient representation 
 of $\Pi$
 in the category of admissible smooth representations, 
 namely,
 the Casselman--Wallach minimal globalization of $(\Pi^{\vee})_K$
 (\cite[Chap.~11]{W}).  
The topological dual $U^{\vee}$ of $U$ is the space
 of distribution vectors,
 on which we can define a continuous representation of $G$.  
This is the maximal globalization of $(\Pi^{\vee})_K$
 in the sense of Casselman--Wallach, 
 which we refer to $(\Pi^{\vee})^{-\infty}$.  
Thus we have
\[
   (\Pi^{\vee})_K \subset \Pi^{\vee} \subset (\Pi^{\vee})^{-\infty}.  
\]
We shall use these symbols
 for a representation $\pi$ of the subgroup $G'$ below.  

\begin{example}
\label{ex:psdual}
Let $\widetilde \tau$ be a finite-dimensional representation
 of a parabolic subgroup $P'$ of $G'$, 
 and $\pi:={\operatorname{Ind}}_{P'}^{G'}(\widetilde \tau)$
 the representation on $C^{\infty}(G'/P', {\mathcal{W}})$.  
The dualizing bundle ${\mathcal{W}}^{\ast}$ is given
 as the $G'$-homogeneous vector bundle 
 over $G'/P'$
 associated to $\tau^{\ast}:= \widetilde \tau^{\vee} \otimes |\det ({\operatorname{Ad}}|_{{\mathfrak{g}}'/{\mathfrak{p}}'})|^{-1}$, 
 where $\widetilde \tau^{\vee}$ is the contragredient representation of $\widetilde \tau$.  
Then the smooth admissible representation $\pi^{\vee}$ of $G'$
 is given as a representation 
$
{\operatorname{Ind}}_{P'}^{G'}(\tau^{\ast})
$ 
 on $C^{\infty}(G'/P', {\mathcal{W}}^{\ast})$, 
 whereas $(\pi^{\vee})^{-\infty}$
 is given as a representation
 on ${\mathcal{D}}'(G'/P', {\mathcal{W}}^{\ast})$.  
\end{example}
Any symmetry breaking operator $T \colon \Pi|_{G'} \to \pi^{\vee}$
 induces a continuous bilinear form
\[
  \Pi \boxtimes \pi \to {\mathbb{C}}, 
\qquad
   u \otimes v \mapsto \langle T u, v \rangle,
\]
and we have a natural embedding
\begin{equation}
\label{eqn:SBOtoBF}
   {\operatorname{Hom}}_{G'}(\Pi|_{G'}, \pi^{\vee})
   \hookrightarrow
   {\operatorname{Hom}}_{G'}(\Pi \boxtimes \pi,{\mathbb{C}})
   \simeq
   {\operatorname{Hom}}_{G'}(\Pi|_{G'}, (\pi^{\vee})^{-\infty}).  
\end{equation}
Here the second isomorphism follows from the natural bijections 
 for nuclear Fr{\'e}chet spaces
 (\cite[Prop.~50.7]{Treves}):
\[
   {\operatorname{Hom}}_{{\mathbb{C}}}(U \otimes U',{\mathbb{C}})
   \simeq
   {\operatorname{Hom}}_{{\mathbb{C}}}(U, (U')^{\vee}), 
\]
  where ${\operatorname{Hom}}_{{\mathbb{C}}}$ denotes the space of continuous linear maps.

As an immediate consequence of Fact \ref{fact:SBOdistr}, 
 we have the following:
\begin{proposition}
\label{prop:BSBOps}
Suppose $\widetilde \sigma$ and $\widetilde \tau$ are finite-dimensional representations
 of parabolic subgroups $P$ and $P'$,
 respectively.  
Let $\Pi={\operatorname{Ind}}_P^G (\widetilde \sigma)$
 and $\pi={\operatorname{Ind}}_{P'}^{G'} (\widetilde \tau)$
 be admissible smooth representations
 of $G$ and $G'$, 
 respectively.  
Then the embedding in \eqref{eqn:SBOtoBF} is an isomorphism.  
\end{proposition}

\begin{proof}
We recall
 that ${\operatorname{Hom}}_{{\mathbb{C}}}(\cdot, {\mathbb{C}})$ denotes 
 the space of (continuous) functionals.  
Then ${\operatorname{Hom}}_{G'}(\Pi \boxtimes \pi,{\mathbb{C}})$
 is naturally isomorphic 
 to the spaces of $G'$-invariant elements
 of the following vector spaces
\[
   {\operatorname{Hom}}_{{\mathbb{C}}}
   (C^{\infty}(G/P \times G'/P', {\mathcal{V}} \boxtimes {\mathcal{W}}), {\mathbb{C}})
   \simeq
   {\mathcal{D}}'(G/P \times G'/P', {\mathcal{V}}^{\ast} \boxtimes {\mathcal{W}}^{\ast}), 
\]
 and so we have
\[
   {\operatorname{Hom}}_{G'}(\Pi \boxtimes \pi,{\mathbb{C}})
    \simeq
    {\mathcal{D}}'(G/P \times G'/P', {\mathcal{V}}^{\ast} \boxtimes {\mathcal{W}}^{\ast})^{\Delta(G')}.  
\]
Since $\tau^{\ast\ast} \simeq \tau$, 
 the right-hand side is canonically isomorphic to 
\[{\operatorname{Hom}}_{G'}
   (C^{\infty}(G/P, {\mathcal{V}})|_{G'}, C^{\infty}(G'/P',{\mathcal{W}}^{\ast}))
\simeq 
{\operatorname{Hom}}_{G'}(\Pi|_{G'}, \pi^{\vee})
\]
 by Fact \ref{fact:SBOdistr} and Example \ref{ex:psdual}.  
Hence Proposition \ref{prop:BSBOps} is proved.  
\end{proof}
More generally,
 we obtain the following.  
\begin{theorem}
\label{thm:SBOBF}
Let $G \supset G'$ be a pair of real reductive Lie groups.  
For any $\Pi \in {\operatorname{Irr}}(G)$
 and $\pi \in {\operatorname{Irr}}(G')$, 
 we have a canonical bijection:
\begin{equation}\label{eqn:SBOBF}
   {\operatorname{Hom}}_{G'}(\Pi|_{G'}, \pi^{\vee})
   \overset \sim \rightarrow
   {\operatorname{Hom}}_{G'}(\Pi \boxtimes \pi,{\mathbb{C}}).  
\end{equation}
\end{theorem}
By the second isomorphism \eqref{eqn:SBOtoBF}, 
Theorem \ref{thm:SBOBF} is deduced from the following proposition, 
 where we change the notation from 
 $\pi^{\vee}$ to $\pi$ for simplicity.

\begin{proposition}
\label{prop:SBOtarget}
Suppose $\Pi \in {\operatorname{Irr}}(G)$
 and $\pi \in {\operatorname{Irr}}(G')$, 
Let $\pi^{-\infty}$ be the representation of $G'$
 on distribution vectors.  
Then the natural embedding
\[
   {\operatorname{Hom}}_{G'}(\Pi|_{G'}, \pi)
   \hookrightarrow
   {\operatorname{Hom}}_{G'}(\Pi|_{G'}, \pi^{-\infty})
\]
 is a bijection.  
\end{proposition}
\begin{proof}
[Proof of Proposition \ref{prop:SBOtarget}]
We take $P$ and $P'$ to be minimal parabolic subgroups of $G$ and $G'$, 
 respectively.  
By Casselman's subrepresentation theorem
 (or equivalently, \lq\lq{quotient theorem}\rq\rq), 
 see \cite[Chap.~3, Sect.~8]{W} for instance,
 for any $\Pi \in {\operatorname{Irr}}(G)$, 
 there exists an irreducible finite-dimensional representation
 $(\widetilde \sigma, V)$ of $P$
 such that
 $\Pi_K$ is obtained as a quotient
 of ${\operatorname{Ind}}_P^{G}(\widetilde \sigma)_K$, 
 and therefore,
 there is a surjective continuous $G$-homomorphism
$p \colon C^{\infty}(G/P, {\mathcal{V}}) \to \Pi$
 by the automatic continuity theorem
 \cite[Chap.~11, Sect.~4]{W}.  
Likewise, 
 for any $\pi \in {\operatorname{Irr}}(G')$, 
 there exists an irreducible finite-dimensional representation
 $(\widetilde \tau, W)$ of $P'$
 such that
 $\pi_{K'}$ is a subrepresentation of 
 ${\operatorname{Ind}}_{P'}^{G'}(\widetilde \tau)_{K'}$,
 and therefore, 
 there is an injective continuous $G'$-homomorphism
$\iota \colon \pi^{-\infty} \hookrightarrow {\mathcal{D}}'(G'/P', {\mathcal{W}})$
 by the dual of the automatic continuity theorem.  
If $T \colon \Pi \to\pi^{-\infty}$ is a continuous $G'$-homomorphism, 
 then $T$ induces a continuous $G'$-homomorphism
\[
  \iota\circ T \circ p \colon C^{\infty}(G/P, {\mathcal{V}}) \to {\mathcal{D}}'(G'/P', {\mathcal{W}}).  
\] 
By Proposition \ref{prop:BSBOps}, 
 $\iota\circ T \circ p$ is actually a continuous $G'$-homomorphism, 
\[
  C^{\infty}(G/P, {\mathcal{V}}) \to C^{\infty}(G'/P', {\mathcal{W}}).  
\] 
Hence the image of $T$ is contained in the admissible smooth representation $\pi$.  
Since the topology of the admissible smooth representation $\pi$ coincides
 with the relative topology of $C^{\infty}(G'/P', {\mathcal{W}})$, 
 $T$ is actually a $G'$-homomorphism $\Pi|_{G'} \to \pi$.  
\end{proof}

\begin{remark}
\begin{enumerate}
\item[{\rm{(1)}}]
In \cite[Lem.~A.0.8]{AG},
 the authors proved 
 the injectivity of the map \eqref{eqn:SBOtoBF}.  
\item[{\rm{(2)}}]
Theorem \ref{eqn:SBOBF} simplifies
 part of the proof
 of \cite[Thm.~4.1]{xkShintani}
 on twelve equivalence conditions
 including the finiteness
 criterion for the dimension
 of continuous invariant bilinear forms.  
\end{enumerate}
\end{remark}

\subsection{Distribution kernels of symmetry breaking operators
 for $G=O(n+1,1)$}
We analyze the distribution kernels
 of symmetry breaking operators
 in coordinates.  
For this,
 we set up some structural results
 for $G=O(n+1,1)$.  

\subsubsection{Bruhat and Iwasawa decompositions for $G=O(n+1,1)$}
\label{subsec:KAN}

We recall from \eqref{eqn:psim}
 that the map $\psi_n \colon {\mathbb{R}}^n \setminus \{0\}\to O(n)$, 
 $x \mapsto \psi_n(x)$
 is defined as the reflection
 with respect to the hyperplane orthogonal to $x$.  
By using $\psi_n(x)$, 
 we give an explicit formula
 of the Bruhat decomposition
 $G=N_+ w MAN_+
 \cup MA N_+$
 and the Iwasawa decomposition
 $G=KAN_+$
 for an element of $N_-$
 for $G=O(n+1,1)$.  
Here we set
\index{A}{winversion@$w$, inversion element|textbf}
\begin{equation}
\label{eqn:w}
w:= \operatorname{diag} (1,\cdots,1,-1) \in N_K({\mathfrak{a}}).  
\end{equation}

Retain the notation as in Section \ref{subsec:subgpOn}.  
In particular, 
 we recall from \eqref{eqn:nplus} and \eqref{eqn:nbar}
 the definition 
of the diffeomorphisms
\index{A}{n1@$n_+\colon {\mathbb{R}}^n \to N_+$}
$
n_+ \colon {\mathbb{R}}^n \overset \sim \to N_+
$
 and 
\index{A}{n2@$n_-\colon {\mathbb{R}}^n \to N_-$}
$
n_- \colon {\mathbb{R}}^n \overset \sim \to N_-
$, 
 respectively.  
\begin{lemma}
[Bruhat decomposition]
\label{lem:0.1}
\index{B}{Bruhatdecomposition@Bruhat decomposition}
For $b \in {\mathbb{R}}^n - \{0\}$, 
\[
n_-(b) = n_+(a)
\begin{pmatrix}
-1 & & 
\\
   & \psi_n(b) &
\\
 & & 1
\end{pmatrix}
e^{tH} n, 
\]
where $a \in {\mathbb{R}}^n$
 and $t \in {\mathbb{R}}$ are given uniquely by 
 $a =-\frac{b}{|b|^2}$ and $e^t =|b|^2$, 
 respectively, 
 and $n \in N_+$.  
\end{lemma}

\begin{proof}
Suppose that $a \in {\mathbb{R}}^n$, 
 $\varepsilon = \pm 1$, 
 $B \in O(n)$, 
 $t \in {\mathbb{R}}$
 and $n \in N_+$ satisfies
\begin{equation}
\label{eqn:nBruhat}
n_-(b)= n_+(a) w 
\begin{pmatrix} \varepsilon & & 
\\
                            & B & 
\\
                            & & \varepsilon
\end{pmatrix} e^{t H} n. 
\end{equation}
Applying \eqref{eqn:nBruhat}
 to the vector $p_+ ={}^{t\!}(1,0,\cdots,0,1)\in \Xi$
 (see \eqref{eqn:p+}),  
 we have
\[
\begin{pmatrix}
1-|b|^2
\\
2b
\\
1+ |b|^2
\end{pmatrix}
=
\varepsilon 
e^t
\begin{pmatrix}
1-|a|^2
\\
2a
\\
-1 - |a|^2
\end{pmatrix}.  
\]
Hence $\varepsilon =-1$, 
 $e^t=\frac{1}{|a|^2}$, 
 and $a = -|a|^2 b$.  
Thus $|a|\,|b|=1$.  
In turn, 
 \eqref{eqn:nBruhat} amounts to 
\[
 n_+(a)^{-1} n_-(b)
=
\begin{pmatrix} -1 & & 
\\
                            & B & 
\\
                            & & 1
\end{pmatrix} 
e^{t H} n, 
\]
whence 
$
  B = I_n + 2 a {}^{t\!}b
    = I_n - \frac{2 b {}^{t\!}b}{|b|^2}
    =\psi_n(b).  
$
\end{proof}
For $b \in {\mathbb{R}}^n$, 
 we define $k(b) \in SO(n+1)$ by 
\index{A}{kb@$k(b)$|textbf}
\begin{equation}
\label{eqn:kb}
   k(b):=I_{n+1} + \frac{1}{1+|b|^2}
        \begin{pmatrix}
         -2 |b|^2 & -2\, {}^{t\!}b
         \\
         2 b   & -2\,b\,{}^{t\!}b
        \end{pmatrix}
       = \psi_{n+1}(1,b) \begin{pmatrix} -1 & \\ & I_n \end{pmatrix}.  
\end{equation}
\begin{lemma}
[Iwasawa decomposition]
\label{lem:0.2}
\index{B}{Iwasawadecomposition@Iwasawa decomposition}
For any $b \in {\mathbb{R}}^n$, 
 we have 
\index{A}{n2@$n_-\colon {\mathbb{R}}^n \to N_-$}
\begin{equation}
\label{eqn:nKAN}
   n_-(b)=k (b) e^{t H} n_+(a) 
\in KAN_+, 
\end{equation}
where $a \in {\mathbb{R}}^n$ 
 and $t \in {\mathbb{R}}$ are given by 
$
   a= \frac{-b}{1+|b|^2}
$
 and 
$
   e^t=1+|b|^2. 
$
\end{lemma}
\begin{proof}
We shall prove 
 that $k(b)$ in \eqref{eqn:nKAN} is given
 by the formula \eqref{eqn:kb}.  
Since $n_-(b)$ is contained 
 in the connected component of $G$, 
$k(b)=(k(b)_{ij})_{0 \le i,j \le n}$
 in \eqref{eqn:nKAN}
 belongs to the connected group $SO(n+1)$.  
We write $k(b)=(k(b)_0, k'(b))$
 where $k(b)_0 \in {\mathbb{R}}^{n+1}$
 and $
  k'(b):=(k(b)_{ij})_{\substack {0 \le i \le n \\ 1 \le j \le n}}
  \in M(n+1,n;{\mathbb{R}})
$.  
Applying \eqref{eqn:nKAN} to the vector
\index{A}{p1@$p_+={}^t(1,0,\cdots,0,1)$}
 $p_+ ={}^{t\!} (1,0,\cdots,0,1)$, 
we have
\[
\begin{pmatrix}
1-|b|^2 \\ 2 b\\ 1+|b|^2
\end{pmatrix}
=
e^t
\begin{pmatrix}
k(b)_0 \\ 1
\end{pmatrix}.  
\]
The last component shows
 $e^t=1+|b|^2$.  
In turn, 
 we get the first column vector $k(b)_0$ of $k(b)$.  
On the other hand, 
 we observe 
\[
  k(b)_{i j}
  = 
  (n_-(b) n_+(a)^{-1} e^{-t H})_{ij}
  =
  (n_-(b) n_+(a)^{-1} )_{ij}
\]
for $0 \le i \le n+1$ and $1 \le j \le n$.  
Hence we get
\[
\begin{pmatrix}
k'(b)
\\
 0 \cdots 0
\end{pmatrix}
=
\begin{pmatrix}
(1-|b|^2) {}^{t\!}a-{}^{t\!}b
\\
I_n +2 b {}^{t\!}a
\\
(1+|b|^2) {}^{t\!}a + {}^{t\!}b
\end{pmatrix}, 
\]
which implies
\[
  a=-\frac{b}{1+|b|^2}
\quad
\text{ and }
\quad
  k'(b)
  =\begin{pmatrix}
   -\frac{1-|b|^2}{1+|b|^2} {}^{t\!}b-{}^{t\!}b
  \\
   I_n -\frac{2 b {}^{t\!}b}{1+|b|^2}
   \end{pmatrix}
  =\begin{pmatrix}
   \frac{-2}{1+|b|^2}{}^{t\!}b
   \\
   I_n-\frac{2 b\, {}^{t\!}b}{1+|b|^2}
   \end{pmatrix}.
\]
In particular,
  we have shown 
 that $k(b)$ in \eqref{eqn:nKAN} is given by the formula \eqref{eqn:kb}.  
\end{proof}

\subsubsection{Distribution kernels for symmetry breaking operators}
We apply Fact \ref{fact:SBOdistr} to the pair
 $(G,G')=(O(n+1,1),O(n,1))$
 and a pair of the minimal parabolic subgroups $P$ and $P'$.  
With the notation of Fact \ref{fact:SBOdistr}, 
 we shall take
\begin{alignat*}{2}
\widetilde \sigma =\, & V \otimes \delta \otimes {\mathbb{C}}_{\lambda}
\quad
&&\text{on $V_{\lambda,\delta}$}
\\
\widetilde \tau =\, & W \otimes \varepsilon \otimes {\mathbb{C}}_{\nu}
\quad
&&\text{on $W_{\varepsilon,\nu}$}
\end{alignat*}
as (irreducible) representations of $P$ and $P'$, 
 respectively,
 for $(\sigma, V) \in \widehat{O(n)}$, 
 $\delta \in \{\pm\}$, 
 and $\lambda \in {\mathbb{C}}$
 and $(\tau, W) \in \widehat{O(n-1)}$, 
$\varepsilon \in \{\pm\}$, 
 and $\nu \in {\mathbb{C}}$.  
We recall from \eqref{eqn:Vlmdbdle}
 that 
$
{\mathcal{V}}_{\lambda, \delta}
=G \times_P V_{\lambda,\delta}
$
 is a homogeneous vector bundle over the real flag variety $G/P$.  
The dualizing bundle 
\index{A}{Vlndast@${\mathcal{V}}_{\lambda,\delta}^{\ast}$, dualizing bundle|textbf}
\index{A}{Vlnd@${\mathcal{V}}_{\lambda,\delta}$, homogeneous vector bundle over $G/P$|textbf}
$
{\mathcal{V}}_{\lambda, \delta}^{\ast}
$
 of ${\mathcal{V}}_{\lambda, \delta}$, 
 is given by a $G$-homogeneous vector bundle over $G/P$
 associated to the representation
 of $P/N_+ \simeq MA \simeq O(n) \times {\mathbb{Z}}/2{\mathbb{Z}} \times {\mathbb{R}}$:
\[
 V_{\lambda, \delta}^{\ast}
:=(V_{\lambda,\delta})^{\vee} \otimes {\mathbb{C}}_{2\rho}
\simeq V^{\vee} \boxtimes {\delta} \boxtimes{\mathbb{C}}_{n-\lambda}, 
\]
where 
\index{A}{VWV@$V^{\vee}$, contragredient representation of $V$|textbf}
$V^{\vee}$ denotes the contragredient representation
 of $(\sigma,V)$.  
Then the regular representation
 of $G$ on the space
\index{A}{distribution@${\mathcal{D}}'$, distribution}
$
   {\mathcal{D}}'(G/P, {\mathcal{V}}_{\lambda, \delta}^{\ast})
$
 of ${\mathcal{V}}_{\lambda, \delta}^{\ast}$-valued distribution sections
 is the dual of the representation 
 $I_{\delta}(V,\lambda)$ of $G$
 on $C^{\infty}(G/P, {\mathcal{V}}_{\lambda, \delta})$
 as we discussed in Example \ref{ex:psdual}.

In this special setting, 
 Fact \ref{fact:SBOdistr} amounts to the following.  
\begin{fact}
\label{fact:kernel}
There is a natural bijective map:
\begin{equation}
\label{eqn:ker}
\operatorname{Hom}_{G'}
  (
   I_{\delta}(V,\lambda)|_{G'}, 
   J_{\varepsilon}(W,\nu)
  )
\overset \sim \to 
   ({\mathcal{D}}'
    (G/P, {\mathcal{V}}_{\lambda, \delta}^{\ast})
    \otimes
    W_{\nu,\varepsilon}
    )
    ^{\Delta (P')}, 
\quad
   T \mapsto K_T.  
\end{equation}
\end{fact}

\vskip 0.8pc
In \cite[Def.~3.3]{sbon}, 
 we defined regular symmetry breaking operators
 in the general setting.  
In our special setting,
 there is only one open $P'$-orbit 
 in the real flag manifold $G/P$,
 and thus the definition is reduced to the following.  
\begin{definition}
[regular symmetry breaking operator]
\label{def:regSBO}
A symmetry breaking operator
 $T \colon I_{\delta}(V,\lambda) \to J_{\varepsilon}(W,\nu)$
 is 
\index{B}{regularsymmetrybreakingoperator@regular symmetry breaking operator|textbf}
{\it{regular}} 
 if the support of the distribution kernel $K_T$ is $G/P$.  
\end{definition}

\subsubsection{Distribution sections for dualizing bundle
 ${\mathcal{V}}_{\lambda,\delta}^{\ast}$ over $G/P$}
\label{subsec:GPsec}
This section provides a concrete description
 of the right-hand side of \eqref{eqn:ker}
 in the coordinates on the open 
\index{B}{Bruhatcell@Bruhat cell}
Bruhat cell.

We begin with a description of 
 the $G$- and ${\mathfrak{g}}$-action 
 on ${\mathcal{D}}'
    (G/P, {\mathcal{V}}_{\lambda, \delta}^{\ast})$
 in the coordinates.  
We identify ${\mathcal{D}}'(G/P, {\mathcal{V}}_{\lambda, \delta}^{\ast})$
 with a subspace of $V^{\vee}$-valued distribution on $G$
 via the following map: 
\[
   {\mathcal{D}}'(G/P, {\mathcal{V}}_{\lambda, \delta}^{\ast})
   \simeq 
  ({\mathcal{D}}'(G) \otimes V_{\lambda, \delta}^{\ast})^{\Delta(P)}
   \subset 
   {\mathcal{D}}'(G) \otimes V^{\vee}.  
\]
We recall
 that the Bruhat decomposition of $G$ is given by 
 $G=N_+ w P \cup P$ 
 where 
\index{A}{winversion@$w$, inversion element}
$w= \operatorname{diag}(1,\cdots,1,-1) \in G$, 
 see \eqref{eqn:w}.  
Since the real flag manifold $G/P$ is covered by the two open subsets
 $N_+ w P/P$ and $N_- P/P$, 
 distribution sections on $G/P$ are determined uniquely
 by the restriction  
 to these two open sets:
\begin{equation}
\label{eqn:Dpm}
{\mathcal{D}}'(G/P, {\mathcal{V}}_{\lambda, \delta}^{\ast})
\hookrightarrow
{\mathcal{D}}'(N_+ w P/P, {\mathcal{V}}_{\lambda, \delta}^{\ast}|_{N_+ w P/P})
\oplus
{\mathcal{D}}'(N_- P/P,  {\mathcal{V}}_{\lambda, \delta}^{\ast}|_{N_- P/P}).  
\end{equation}
By a little abuse of notation,
 we use the letters $n_+$ and $n_-$
 to denote the induced diffeomorphisms
$
   {\mathbb{R}}^n \overset \sim \to 
   N_+ w P/P
$
 and 
$
   {\mathbb{R}}^n \overset \sim \to 
   N_- P/P
$, 
respectively.  
Via the following trivialization of the two restricted bundles:
\begin{alignat*}{9}
&  {\mathbb{R}}^n  \times V^{\vee}
\;\;
&& \overset \sim \to 
\;\;
&& {\mathcal{V}}_{\lambda,\delta}^{\ast}|_{N_+ w P/P}
\;\;
&& \subset 
\;\;
&& {\mathcal{V}}_{\lambda,\delta}^{\ast}
\;\;
&& \supset
\;\;
&& {\mathcal{V}}_{\lambda,\delta}^{\ast}|_{N_- P/P}
\;\;
&& \overset \sim \leftarrow
\;\;
&& {\mathbb{R}}^n \times V^{\vee}
\\
& \downarrow
&&
&& \downarrow
&& 
&& \downarrow
&&
&& \downarrow
&&
&&\downarrow
\\
&  {\mathbb{R}}^n  
\;\;
&& \underset{n_+} {\overset \sim \to} 
\;\;
&& N_+ w P / P
\;\;
&& \subset 
\;\;
&& G/P
\;\;
&& \supset
\;\;
&& N_- P / P
\;\;
&& \underset {n_-}{\overset \sim \leftarrow}
\;\;
&& {\mathbb{R}}^n, 
\end{alignat*}
the injection \eqref{eqn:Dpm} is restated as the following map:

\begin{equation}
\label{eqn:fFpm}
   {\mathcal{D}}'(G/P, {\mathcal{V}}_{\lambda, \delta}^{\ast})
   \hookrightarrow
   \left({\mathcal{D}}'({\mathbb{R}}^n) \otimes V^{\vee}\right)
   \oplus
   \left({\mathcal{D}}'({\mathbb{R}}^n) \otimes V^{\vee}\right), 
\qquad
   f \mapsto (F_{\infty}, F)
\end{equation}
where 
\index{A}{n1@$n_+\colon {\mathbb{R}}^n \to N_+$}
\index{A}{n2@$n_-\colon {\mathbb{R}}^n \to N_-$}
\begin{equation*}
F_{\infty}(a):= f(n_+(a) w), 
\qquad
F(b):= f(n_-(b)).  
\end{equation*}

\begin{lemma}
\label{lem:152341}
Let 
\index{A}{1psin@$\psi_n$}
$
\psi_n \colon {\mathbb{R}}^n\setminus \{0\} \to O(n)
$ be the map
 taking the reflection defined in \eqref{eqn:psim}.  

\begin{enumerate}
\item[{\rm{(1)}}]
The image of the injective map \eqref{eqn:fFpm}
 is characterized by the following identity
 in ${\mathcal{D}}'({\mathbb{R}}^n\setminus \{0\}) \otimes V^{\vee}$:
\begin{equation}
\label{eqn:F12}
   F(b)
   =
   \delta \sigma^{\vee}(\psi_n(b)^{-1})
   |b|^{2\lambda-2n}
   F_{\infty}(-\frac{b}{|b|^2})
\qquad
\text{on }
{\mathbb{R}}^n \setminus \{0\}.  
\end{equation}

\item[{\rm{(2)}}]
{\rm{(first projection)}}
$f \in {\mathcal{D}}'(G/P, {\mathcal{V}}_{\lambda, \delta}^{\ast})$
 is supported at the singleton 
\index{A}{p1@$p_+={}^t(1,0,\cdots,0,1)$}
$
\{[p_+]\} = \{ e P/P\}
$
 if and only if $F_{\infty}=0$.  
\item[{\rm{(3)}}]
{\rm{(second projection)}}
The second projection $f \mapsto F$ is injective.  
\end{enumerate}
\end{lemma}

\begin{proof}
(1) \enspace
The image of the map \eqref{eqn:Dpm} is characterized
 by the compatibility condition
 on the intersection $(N_+ w P \cap N_- P)/P$, 
 namely, 
 the pair $(F_{\infty}, F)$
 in \eqref{eqn:fFpm} should satisfy:
\[
  F (b)=\sigma_{\lambda,\delta}^{\ast}(p)^{-1} F_{\infty}(a)
\]
for all $(a,b,p) \in {\mathbb{R}}^n \times {\mathbb{R}}^n \times P$
 such that $n_+(a) w p =n_-(b)$.  
In this case, 
 $b \ne 0$
 because $N_+ w P \not \ni e$.  
By Lemma \ref{lem:0.1}, 
we have
\[
a = - \frac{b}{|b|^2}, 
\quad
p= \begin{pmatrix} -1 & & \\
                              & \psi_n(b) & \\
                              &  &  -1 \end{pmatrix} 
        e^{tH}, 
\]
where $e^t = |b|^2$.  
Then 
\begin{align*}
F(b) =& f(n_-(b))
\\
       =& \sigma_{\lambda, \delta}^{\ast}(p^{-1})
        f(n_+(a) w)
\\
       =& \delta|b|^{2\lambda-2n} 
          \sigma^{\vee} (\psi_n(b)) F_{\infty}(a).  
\end{align*}
\begin{enumerate}
\item[(2)]
Clear from $G \setminus N_+ w P=P$.  
\item[(3)]
Since $P' N_-P=G$
\cite[Cor.~5.5]{sbon}, 
the third statement follows from \cite[Thm.~3.16]{sbon}.  
\end{enumerate}
\end{proof}

The regular representation of $G$
 on ${\mathcal{D}}'(G/P, {\mathcal{V}}_{\lambda, \delta}^{\ast})$ induces
 an action
 on the pairs $(F_{\infty}, F)$
 of $V^{\vee}$-valued distributions
 through Lemma \ref{lem:152341} (1).  
We need an explicit formula
 of the action of the parabolic subgroup $P=MAN_+$
 or its Lie algebra
$
  {\mathfrak{p}}
 ={\mathfrak{m}}
 +{\mathfrak{a}}
 +{\mathfrak{n}}_+
$, 
 which is given in the following two elementary lemmas.

We begin with the first projection $f \mapsto F_{\infty}$ in \eqref{eqn:fFpm}.  
Since the action of $P$ on $G/P$ leaves the open subset 
 $N_+ w P /P = P w P /P$ invariant, 
 we can define the geometric action
 of the group $P$ on ${\mathcal{D}}'(N_+ w P /P, {\mathcal{V}}_{\lambda,\delta}^{\ast})$
 as follows.  
We recall 
$
M=O(n) \times \{1, m_-\}
$
 (see \eqref{eqn:m-}).  
We collect
 some basic formul{\ae} for the coordinates 
 $n_{\varepsilon}\colon {\mathbb{R}}^n \overset \sim \to N_{\varepsilon}$:
 for $\varepsilon = \pm$
 (by abuse of notation, 
 we also write as $\varepsilon = \pm 1$), 
\index{A}{H@$H$|textbf}
\begin{align}
{n_{\varepsilon}} (Bb)
   =& \begin{pmatrix} 1 \\ & B \\ && 1 \end{pmatrix}
\,\,
n_{\varepsilon} (b)
      \begin{pmatrix} 1 \\ & B^{-1} \\ && 1 \end{pmatrix}
\quad
\text{for $B \in O(n)$, }
\label{eqn:nAb}
\\
 {n_{\varepsilon}} (-b)
   =& m_- {n_{\varepsilon}} (b) m_-^{-1},
\label{eqn:minvn}
\\
 {n_{\varepsilon}} (e^{{\varepsilon}t} b)
   =& e^{tH} {n_{\varepsilon}} (b) e^{-tH}.  
\label{eqn:invm}
\end{align}
\begin{lemma}
\label{lem:152343}
We let $P=M A N_+$ act on 
$
   {\mathcal{D}}'({\mathbb{R}}^n) \otimes V^{\vee}
$
by
\index{A}{m2@$m_-={\operatorname{diag}}(-1,1,\cdots,1,-1)$}
\begin{alignat}{3}
&\left(\pi \begin{pmatrix} 1 & & \\ & B & \\ & & 1 \end{pmatrix} F_{\infty}\right)(a)
&&=\sigma^{\vee}(B) F_{\infty}(B^{-1}a)
\quad
&&\text{ for $B \in O(n)$, }
\label{eqn:inftyM}
\\
&(\pi (m_-) F_{\infty})(a)
&&=\delta F_{\infty}(-a), 
&&
\label{eqn:inftym}
\\
&(\pi (e^{tH}) F_{\infty})(a)
&&= e^{(\lambda-n)t} F_{\infty}(e^{-t}a)
\quad
&&\text{ for all $t \in {\mathbb{R}}$, }
\label{eqn:inftyA}
\\
&(\pi (n_+(c)) F_{\infty})(a)
&&= F_{\infty}(a-c)
\quad
&&\text{ for all $c \in {\mathbb{R}}^{n}$. }
\label{eqn:inftyN}
\end{alignat}
Then the first projection $f \mapsto F_{\infty}$ in \eqref{eqn:fFpm}
 is a $P$-homomorphism.  
\end{lemma}
\begin{proof}
We give a proof for \eqref{eqn:inftyA} on the action
 of the split abelian group $A$.  
Let $t \in {\mathbb{R}}$. 
By \eqref{eqn:invm}
 and 
\index{A}{winversion@$w$, inversion element}
$
e^{-tH} w =w e^{tH}
$, 
 we have 
\[
f(e^{-tH} n_+(a) w)
=f(n_+(e^{-t}a)e^{-tH}w)
=e^{(\lambda-n)t} f(n_+(e^{-t}a)w)
=e^{(\lambda-n)t} F_{\infty}(e^{-t} a), 
\]
whence we get the desired formula.  
The proof for the actions of $M$ and $N_+$ is similar.  
\end{proof}

Next,
 we consider the second projection $f \mapsto F$
 in \eqref{eqn:fFpm}.  
In this case,
 the group $N_+$ does not preserve the open subset
 $N_-P/P$ in $G/P$, 
 and therefore we shall use the action 
 of the Lie algebra ${\mathfrak {n}}_+$ instead
 (see \eqref{eqn:On} below).  
We denote by 
\index{A}{Euler@$E$, Euler homogeneity operator\quad|textbf}
$
E
$ 
 the Euler homogeneity operator
 $\sum_{\ell=1}^n x_\ell \frac{\partial}{\partial x_\ell}$.  
\begin{lemma}
\label{lem:20161104}
We let the group $MA$ and the Lie algebra ${\mathfrak{n}}_+$
 act on
 ${\mathcal{D}}'({\mathbb{R}}^n) \otimes V^{\vee}$
 by 
\index{A}{Nj1@$N_j^+$}
\begin{alignat}{2}
& \left(\pi \begin{pmatrix} 1 & & \\ & B & \\ & & 1 \end{pmatrix}  F\right)(b)
&&=\sigma^{\vee}(B) F(B^{-1} b)
\quad
\text{for $B \in O(n)$}, 
\label{eqn:0M}
\\
& (\pi(m_-) F)(b)
&&=
\delta F(-b), 
\label{eqn:0m}
\\
& (\pi (e^{tH})F)(b)
&&=
e^{(n-\lambda) t}F(e^t b)
\quad
\text{for all $t \in {\mathbb{R}}$}, 
\label{eqn:0A}
\\
& d\pi (N_j^+)F(b)
&&=
\left((\lambda-n)b_j - b_j E 
 + \frac 1 2 |b|^2 \frac{\partial}{\partial b_j}\right) F
\quad
\text{for $1 \le j \le n$}.  
\label{eqn:On}
\end{alignat}
Here $b=(b_1, \cdots, b_n)$.  
Then the second projection $f \mapsto F$
 in \eqref{eqn:fFpm}
 is an $(M A, {\mathfrak {n}}_+)$-homomorphism.  
\end{lemma}

\begin{proof}
See \cite[Prop.~6.4]{sbon}
for \eqref{eqn:On}.  
The other formul\ae\ are easy, 
 and we omit the proof.  
\end{proof}

\subsubsection{Pair of distribution kernels
 for symmetry breaking operators}
\label{subsec:Tpair}

We extend Lemma \ref{lem:152341}
 to give a local expression
 of the distribution kernels
 of symmetry breaking operators
 via the isomorphism \eqref{eqn:ker}.  
Suppose $(\tau,W) \in \widehat{O(n-1)}$, 
 $\nu \in {\mathbb{C}}$, 
and $\varepsilon \in \{ \pm \}$.  
We define 
\begin{equation}
\label{eqn:DVWinfty}
   ({\mathcal{D}}'({\mathbb{R}}^n)
    \otimes 
    \operatorname{Hom}_{{\mathbb{C}}}(V,W))^{\Delta(P')}
   \equiv
({\mathcal{D}}'({\mathbb{R}}^n)
    \otimes 
    \operatorname{Hom}_{{\mathbb{C}}}(V_{\lambda, \delta},W_{\nu,\varepsilon})
   )^{\Delta(P')}
\end{equation}
 to be the space of $\operatorname{Hom}_{{\mathbb{C}}}(V,W)$-valued
 distributions ${\mathcal{T}}_{\infty}$ on ${\mathbb{R}}^n$
satisfying the following four conditions:
\begin{alignat}{2}
\label{eqn:152345a}
& \tau(B) \circ {\mathcal{T}}_{\infty}(B^{-1} y,y_n) \circ \sigma^{-1}(B)
  ={\mathcal{T}}_{\infty}(y,y_n)
\qquad
&&\text{for all } B \in O(n-1), 
\\
\label{eqn:152345b}
& {\mathcal{T}}_{\infty}(-y,-y_n)= \delta \varepsilon {\mathcal{T}}_{\infty}(y,y_n), 
&&
\\
\label{eqn:152345c}
& {\mathcal{T}}_{\infty}(e^t y,e^t y_n)= e^{(\lambda+\nu-n)t}{\mathcal{T}}_{\infty}(y,y_n)
\quad
&&\text{for all }t \in {\mathbb{R}}, 
\\
\label{eqn:152345d}
& {\mathcal{T}}_{\infty}(y-z,y_n)={\mathcal{T}}_{\infty}(y,y_n)
\quad
&&\text{for all }z \in {\mathbb{R}}^{n-1}.  
\end{alignat}
For the open Bruhat cell $N_- P \subset G$, 
 we consider the following.  
\begin{definition}
\label{def:solVW}
We define 
\index{A}{SolRVW@${\mathcal{S}}ol
({\mathbb{R}}^n;V_{\lambda, \delta}, W_{\nu, \varepsilon})$|textbf}
${\mathcal{S}}ol
({\mathbb{R}}^n;V_{\lambda, \delta}, W_{\nu, \varepsilon})
\subset
{\mathcal{D}}'({\mathbb{R}}^n)
    \otimes 
    \operatorname{Hom}_{{\mathbb{C}}}(V,W)
$
to be the space of $\operatorname{Hom}_{\mathbb{C}}(V,W)$-valued
 distributions 
${\mathcal{T}}$ on ${\mathbb{R}}^n$
 satisfying 
the following invariance 
under the action of the Lie algebras
 ${\mathfrak{a}}$, ${\mathfrak{n}}_+'$, 
 and the group $M' \simeq O(n-1) \times O(1)$:
\begin{alignat}{2}
&(E-(\lambda-\nu-n)){\mathcal{T}}=0, 
&&
\label{eqn:Fainv}
\\
&\left((\lambda-n)x_j - x_j E + \frac 1 2 (|x|^2 + x_n^2)\frac{\partial}{\partial x_j}\right){\mathcal{T}}=0
\quad
&&(1 \le j \le n-1), 
\label{eqn:Fninv}
\\
&\tau(m) \circ {\mathcal{T}}(m^{-1}b)\circ \sigma(m^{-1})={\mathcal{T}}(b)
\quad
&&\text{for all } m \in O(n-1), 
\label{eqn:FMinv}
\\
&{\mathcal{T}}(-b)= \delta \varepsilon {\mathcal{T}}(b).  
&&
\label{eqn:Fparity}
\end{alignat}
\end{definition}

Applying Lemma \ref{lem:152341}
 to the right-hand side of \eqref{eqn:ker}, 
 we have the following:
\begin{proposition}
\label{prop:Tpair}
Let $(\sigma, V) \in \widehat M$, 
$(\tau,W) \in \widehat {M'}$, 
 $\delta$, $\varepsilon \in \{ \pm \}$, 
 and $\lambda, \nu \in {\mathbb{C}}$.  
\begin{enumerate}
\item[{\rm{(1)}}]
There is a one-to-one correspondence between a symmetry breaking operator
\[
  {\mathbb{T}} 
  \in 
  \operatorname{Hom}_{G'}(I_{\delta}(V,\lambda)|_{G'},J_{\varepsilon}(W,\nu))
\]
 and a pair $({\mathcal{T}}_{\infty}, {\mathcal{T}})$
 of $\operatorname{Hom}_{{\mathbb{C}}}(V,W)$-valued distributions
  on ${\mathbb{R}}^n$
 subject to the following three conditions:
\begin{align}
{\mathcal{T}}_{\infty}
\in &
({\mathcal{D}}'({\mathbb{R}}^n) \otimes \operatorname{Hom}_{{\mathbb{C}}}(V_{\lambda,\delta},W_{\nu,\varepsilon}))^{\Delta(P')}, 
\label{eqn:Tinfty}
\\
{\mathcal{T}}
\in &
{\mathcal{S}}ol
({\mathbb{R}}^n;V_{\lambda, \delta}, W_{\nu, \varepsilon}), 
\label{eqn:Tzero}
\\
{\mathcal{T}}(b)
= &
\delta Q(b)^{\lambda-n}
{\mathcal{T}}_{\infty}
\left(-\frac{b}{|b|^2}\right)
\circ
\sigma(\psi_n(b))
\quad
\text{on }
{\mathbb{R}}^n \setminus \{0\}.  
\label{eqn:Tpatching}
\end{align}

\item[{\rm{(2)}}]
${\mathcal{T}}$ determines ${\mathbb{T}}$ uniquely.  
\item[{\rm{(3)}}]
Suppose that ${\mathbb{T}} \leftrightarrow 
({\mathcal{T}}_{\infty}, {\mathcal{T}})
\index{A}{TinftyT@$({\mathcal{T}}_{\infty}, {\mathcal{T}})$}
$
 is the correspondence in (1).  
Then the following three conditions are equivalent:
\begin{enumerate}
\item[{\rm{(i)}}]
${\mathcal{T}}_{\infty}=0$.  
\item[{\rm{(ii)}}]
$\operatorname{Supp} {\mathcal{T}} \subset \{0\}$.  
\item[{\rm{(iii)}}]
${\mathbb{T}}$ is a differential operator
 (see Definition \ref{def:diff}).  
\end{enumerate}
\end{enumerate}
\end{proposition}

\begin{proof}
The first statement follows from Fact \ref{fact:kernel}, 
 Lemmas \ref{lem:152341} (1), 
 \ref{lem:152343}
 and \ref{lem:20161104}.  
The second statement is immediate from 
 Lemma \ref{lem:152341} (3).  
The third one is proved in \cite{KP1}, 
 see Section \ref{subsec:diff}
 for more details about differential operators
 between two manifolds.  
\end{proof}

\begin{remark}
\label{rem:7.5}
The advantage of using ${\mathcal{T}}$ is
 that the second projection
\[
 \operatorname{Hom}_{G'}
 (I_{\delta}(V,\lambda)|_{G'}, J_{\varepsilon}(W, \nu))
\overset \sim \to 
{\mathcal{S}}ol
({\mathbb{R}}^n;\sigma_{\lambda, \delta}, \tau_{\nu, \varepsilon}),
\quad
 {\mathbb{T}} \mapsto {\mathcal{T}}
\]
 is bijective, 
 and therefore, 
it is sufficient
 to use ${\mathcal{T}}$ 
 in order to describe a symmetry breaking operator ${\mathbb{T}}$.  
This was the approach 
 that we took in \cite{sbon}.  
In this monograph,
 we shall use both ${\mathcal{T}}_{\infty}$ and ${\mathcal{T}}$.  
The advantage of using ${\mathcal{T}}_{\infty}$ is
 that the group $P'$ leaves $N_+ w P/P$ invariant,
 and consequently,
 we can easily determine ${\mathcal{T}}_{\infty}$
 (see Proposition \ref{prop:20150828-1231} below), 
 although the first projection
\[
 \operatorname{Hom}_{G'}
 (I_{\delta}(V,\lambda)|_{G'}, J_{\varepsilon}(W, \nu))
\to 
 ({\mathcal{D}}'({\mathbb{R}}^n) \otimes \operatorname{Hom}_{{\mathbb{C}}}(V,W))^{\Delta(P')}, 
\quad
 {\mathbb{T}} \mapsto {\mathcal{T}}_{\infty}
\]
 is neither injective nor surjective.  
We shall return to this point
 in Section \ref{subsec:holoAVW}.  
\end{remark}

\subsection{Distribution kernels near infinity}
\label{subsec:infty}
Let $({\mathcal{T}}_{\infty}, {\mathcal{T}})$ be
 as in Proposition \ref{prop:Tpair}.  
This section determines ${\mathcal{T}}_{\infty}$ 
 up to scalar multiplication.  
The main result is Proposition \ref{prop:20150828-1231}, 
 which also determines uniquely the restriction of ${\mathcal{T}}$
 to ${\mathbb{R}}^n \setminus \{0\}$
 up to scalar multiplication.  

\begin{example}
\label{ex:152341}
For $\sigma = {\bf{1}}$, 
$\tau={\bf{1}}$, $\delta=+1$, and 
\[
   {\mathcal{T}}_{\infty}(y,y_n)=|y_n|^{\lambda+\nu-n}, 
\]
 we have from \eqref{eqn:Tpatching}
\[
   {\mathcal{T}}(x,x_n) = (|x|^2+x_n^2)^{-\nu} |x_n|^{\lambda+\nu-n}.  
\]
\end{example}

We begin with the following classical result
 on homogeneous distributions
 of one variable:
\begin{lemma}
\label{lem:Riesz}
\begin{enumerate}
\item[{\rm{(1)}}]
Both $\frac{1}{\Gamma(\frac \mu 2)}|t|^{\mu-1}$
 and $\frac{1}{\Gamma(\frac {\mu+1} 2)}|t|^{\mu-1} \operatorname{sgn} t$
 are nonzero distributions
 on ${\mathbb{R}}$
 that depend holomorphically
 on $\mu$ in the entire complex plane ${\mathbb{C}}$.  
\item[{\rm{(2)}}]
Suppose $k \in {\mathbb{N}}$.  
Then 
\begin{alignat*}{2}
\frac{|t|^{\mu-1}}{\Gamma(\frac \mu 2)}
=&
\frac{(-1)^k}{2^k (2k-1)!!}
\delta^{(2k)}(t)
\qquad
&&\text{if }
\mu=-2k, 
\\
\frac{|t|^{\mu-1}\operatorname{sgn} t}{\Gamma(\frac {\mu+1} 2)}
=&
\frac{(-1)^k(k-1)!}{(2k-1)!}
\delta^{(2k-1)}(t)
\qquad
&&\text{if }
\mu=-2k-1.  
\end{alignat*}
\item[{\rm{(3)}}]
Suppose $\mu \in {\mathbb{C}}$ and $\gamma = \pm 1$.  
Then any distribution $g(t)$ on ${\mathbb{R}}$
 satisfying the homogeneity condition 
\[
  \text{$g(at)=a^{\mu-1}g (t)$
 for all 
 $a >0$, 
 and 
 $g(-t)= \gamma g(t)$
}
\]
 is a scalar multiple
 of $\frac{1}{\Gamma(\frac \mu 2)}|t|^{\mu-1}$ $(\gamma =1)$,
 or of $\frac{1}{\Gamma(\frac {\mu+1} 2)}|t|^{\mu-1} \operatorname{sgn} t$
 $(\gamma =-1)$.  
\end{enumerate}
\end{lemma}

For $(\sigma, V) \in \widehat{O(n)}$ and $(\tau, W) \in \widehat{O(n-1)}$,
 we recall that
$ [V:W] $
is the dimension of
$ 
{\operatorname{Hom}}_{O(n-1)}(V|_{O(n-1)},W)
$.  
Suppose $[V:W] \ne 0$, 
 or equivalently, 
 $[V:W] =1$.  
We fix a generator 
\[
 \pr V W \in \operatorname{Hom}_{O(n-1)}(V|_{O(n-1)},W)
\]
 which is unique
 up to nonzero scalar multiplication
 by Schur's lemma.  
In light of the $\Gamma$-factors in Lemma \ref{lem:Riesz}, 
we introduce $\operatorname{Hom}_{{\mathbb{C}}}(V,W)$-valued distributions 
\index{A}{Actt0inf@$(\Attcal \lambda \nu {\pm} {V,W})_{\infty}$|textbf}
$
(\Attcal \lambda \nu \pm {V,W})_{\infty}
$
 on ${\mathbb{R}}^n$ 
 that depend holomorphically on $(\lambda,\nu) \in {\mathbb{C}}^2$
 by 
\begin{align}
(\Attcal \lambda \nu + {V,W})_{\infty}(x,x_n)
:=&
\frac{1}{\Gamma(\frac{\lambda+ \nu-n+1}{2})}|x_n|^{\lambda+\nu-n} \pr V W, 
\label{eqn:AVW+u}
\\
(\Attcal \lambda \nu - {V,W})_{\infty}(x,x_n)
:=&
\frac{1}{\Gamma(\frac{\lambda+ \nu-n+2}{2})}|x_n|^{\lambda+\nu-n}
\operatorname{sgn} x_n \pr V W.  
\label{eqn:AVW-u}
\end{align}
We regard $\pr V W =0$
 if $[V:W] = 0$.  

\begin{remark}
The notation $(\Attcal \lambda \nu \gamma {V,W})_{\infty}$
 with double tildes
 is used here
 because it will be compatible
 with the 
\index{B}{renormalized regular symmetry breaking operator@regular symmetry breaking operator, renormalized---}
{\it{renormalization}}
 $\Attbb \lambda \nu \gamma {V,W}$
 of the normalized symmetry breaking operator
 $\Atbb \lambda \nu \gamma {V,W}$
 which we will introduce 
 in the next sections.  
\end{remark}

Let $\gamma=\delta \varepsilon$.  
If there exists 
$
   {\mathcal {T}}_{\gamma} 
   \in 
   {\mathcal{S}}ol({\mathbb{R}}^n;V_{\lambda,\delta}, W_{\nu,\varepsilon})
$
 such that the pair 
 $((\Attcal \lambda \nu \varepsilon {V,W})_{\infty}, {\mathcal{T}}_{\gamma})$
 satisfies the compatibility condition \eqref{eqn:Tpatching},  
 then the restriction ${\mathcal{T}}_{\gamma}|_{{\mathbb{R}}^n \setminus \{0\}}$
 must be of the form 
\index{A}{Actt0@$\Attcal \lambda \nu {\pm} {V,W}$|textbf}
$
(\Attcal \lambda \nu \gamma {V,W})' \in 
 {\mathcal{D}}'({\mathbb{R}}^n \setminus \{0\}) \otimes {\operatorname{Hom}}_{\mathbb{C}}(V,W)
$
where we set 
\begin{align}
(\Attcal {\lambda}{\nu}{+}{V,W})'
 :=&
\frac{1}{\Gamma(\frac{\lambda + \nu-n+1}{2})}
         (|x|^2+ x_n^2)^{-\nu}
|x_n|^{\lambda+\nu-n}
\Rij V W (x,x_n), 
\label{eqn:AVW+}
\\
(\Attcal {\lambda}{\nu}{-}{V,W})'
 :=&
\frac{1}{\Gamma(\frac{\lambda+\nu-n+2}{2})}
         (|x|^2 +x_n^2)^{-\nu}
|x_n|^{\lambda+\nu-n}
 {\operatorname{sgn}} x_n
\Rij VW (x,x_n), 
\label{eqn:AVW-}
\end{align}
with $\Rij V W = \pr V W \circ \sigma \circ \psi_n$
 (see \eqref{eqn:RVW}).  
We have used the notation
$
   (\Attcal {\lambda}{\nu}{\gamma}{V,W})'
$
 instead of 
$
   \Attcal {\lambda}{\nu}{\gamma}{V,W}
$
 because it is defined 
 only on ${\mathbb{R}}^n \setminus \{0\}$
 and may not extend to ${\mathbb{R}}^n$
 (see Proposition \ref{prop:20170213} below).

Then we have:
\begin{proposition}
\label{prop:20150828-1231}
{\rm{(1)}}\enspace
For any $(\sigma,V) \in \widehat{O(n)}$, 
 $(\tau,W) \in \widehat{O(n-1)}$,
 $\delta, \varepsilon \in \{ \pm \}$, 
 and $\lambda, \nu \in {\mathbb{C}}$, 
 we have
\[
({\mathcal{D}}'({\mathbb{R}}^n)
   \otimes 
   \operatorname{Hom}_{\mathbb{C}} (V_{\lambda,\delta}, W_{\nu, \varepsilon}))
^{\Delta(P')}
=
 {\mathbb{C}}(\Attcal \lambda \nu {\delta\varepsilon}{V,W})_{\infty}.  
\]
\newline
{\rm{(2)}}\enspace
If $[V:W] \ne 0$
 then 
 $(\Attcal \lambda \nu \pm {V,W})' \ne 0$
 for all $\lambda, \nu \in {\mathbb{C}}$.  
\newline
{\rm{(3)}}\enspace
If ${\mathcal{T}} \in {\mathcal{S}}ol({\mathbb{R}}^n;V_{\lambda,\delta}, W_{\nu, \varepsilon})$, 
 then ${\mathcal{T}}|_{{\mathbb{R}}^n \setminus \{0\}}$ is a scalar multiple
 of $(\Attcal \lambda \nu {\delta \varepsilon}{V,W})'$.  
\end{proposition}

\begin{proof}
Suppose 
$
   F 
   \in 
   ({\mathcal{D}}'({\mathbb{R}}^n) \otimes {\operatorname{Hom}}_{\mathbb{C}}(V_{\lambda,\delta}, W_{\nu,\varepsilon}))^{\Delta(P')}$.  
\newline
(1)\enspace
Let $p_n \colon {\mathbb{R}}^n \to {\mathbb{R}}$ be
 the $n$-th projection, 
 and 
$
   p_n^{\ast} \colon {\mathcal{D}}'({\mathbb{R}}) \to {\mathcal{D}}'({\mathbb{R}}^n)
$
 the pull-back of distributions.  
By the $N_+'$-invariance
 \eqref{eqn:152345d}, 
 $F$ depends only on the last coordinate, 
 namely, 
 $F$ is of the form $p_n^{\ast}f$
 for some $f \in {\mathcal{D}}'({\mathbb{R}}) \otimes \operatorname{Hom}_{\mathbb{C}}(V,W)$.  
In turn,
 the $O(n-1)$-invariance \eqref{eqn:152345a} implies
\[
   f \in {\mathcal{D}}'({\mathbb{R}}) \otimes \operatorname{Hom}_{O(n-1)}(V|_{O(n-1)},W).  
\]
In particular, 
 $F=0$
 if $[V:W]=0$.  

{}From now,
 we assume $[V:W] \ne 0$.  
Then $f$ is of the form $h(y_n) \pr V W$
 for some $h(t) \in {\mathcal{D}}'({\mathbb{R}})$.  
By \eqref{eqn:152345b} and \eqref{eqn:152345c}, 
 $h$ is a homogeneous distribution
 of degree $\lambda+\nu-n$
 and of parity $\delta \varepsilon$.  
Then $h(t)$ is determined by Lemma \ref{lem:Riesz}, 
 and we get the desired result.  
\newline
(2)\enspace
The assertion follows from the nonvanishing statement
 for the distribution
 of one-variable
 (see Lemma \ref{lem:Riesz} (1)).  
\newline
(3)\enspace
The third statement follows from 
 the first assertion and Proposition \ref{prop:Tpair}.  
\end{proof}

\subsection{Vanishing condition of differential symmetry breaking operators:
Proof of Theorem \ref{thm:vanDiff} (1)}
\label{subsec:vanDiff}

In this section, we prove a necessary condition 
 for the existence of nonzero differential symmetry breaking operators
 as stated in Theorem \ref{thm:vanDiff} (1):
\begin{theorem}
[vanishing of differential symmetry breaking operators]
\label{thm:vanishDiff}
~~~\newline
Suppose that $V$ and $W$ are finite-dimensional representations
 of $O(n)$ and $O(n-1)$, 
respectively,
 $\delta$, $\varepsilon \in \{\pm\}$, 
 and $(\lambda,\nu) \in {\mathbb{C}}^2$.  
If $(\lambda, \nu, \delta, \varepsilon)$
 satisfies the 
\index{B}{genericparametercondition@generic parameter condition}
generic parameter condition
 \eqref{eqn:nlgen}, 
 namely, 
 $\nu - \lambda \not \in 2{\mathbb{N}}$
 for $\delta \varepsilon=+$, 
 or $\nu - \lambda \not \in 2{\mathbb{N}}+1$
 for $\delta \varepsilon=-$, 
then 
\[
  {\operatorname{Diff}}_{G'}
  (I_{\delta}(V,\lambda)|_{G'},J_{\varepsilon}(W,\nu))=\{0\}.  
\]
\end{theorem}

\begin{remark}
In the above theorem, 
 we do not impose any assumption on $V$ and $W$.  
In Chapter \ref{sec:DSVO}, 
 we give a converse implication
 under the assumption $[V:W] \ne 0$, 
 see Theorem \ref{thm:existDSBO}.  
\end{remark}

For the proof of Theorem \ref{thm:vanishDiff}, 
 we use the following properties of distributions supported at the origin:
\begin{lemma}
\label{lem:distpos}
Let $F$ be any $\operatorname{Hom}_{{\mathbb{C}}}
 (V,W)$-valued distribution on ${\mathbb{R}}^n$
 supported at the origin
 and satisfying the Euler homogeneity differential equation \eqref{eqn:Fainv}.  
\begin{enumerate}
\item[{\rm{(1)}}]
Assume $\nu - \lambda \not \in {\mathbb{N}}$.  
Then $F$ must be zero. 
\item[{\rm{(2)}}]
Assume $\nu - \lambda \in {\mathbb{N}}$.  
Then $F(-x)=(-1)^{\nu -\lambda}F(x)$.   
\end{enumerate}
\end{lemma}
\begin{proof}
Let $\delta(x) \equiv \delta(x_1, \cdots, x_n)$ be the Dirac delta function 
on ${\mathbb{R}}^n$.  
For a multi-index $\alpha=(\alpha_1, \cdots, \alpha_n) \in {\mathbb{N}}^n$, 
 we define another distribution by 
\[
\delta^{(\alpha)}(x_1, \cdots, x_n)
:=
\frac{\partial^{|\alpha|}}{\partial x_1^{\alpha_1} \cdots \partial x_n^{\alpha_n}}  
\delta(x_1, \cdots, x_n)
\]
where $|\alpha|=\alpha_1+ \cdots +\alpha_n$.   
By the structural theory 
 of distributions,
 $F$ must be of the following form 
\[
   F=\sum_{\alpha \in {\mathbb{N}}^n} a_{\alpha} 
   \delta^{(\alpha)}(x_1,\cdots,x_n)
\qquad
\text{(finite sum)}
\]
with some $a_{\alpha} \in \operatorname{Hom}_{{\mathbb{C}}}
 (V,W)$
 for $\alpha \in {\mathbb{N}}^n$.  
Since $\delta^{(\alpha)}(x_1,\cdots,x_n)$ is 
 a homogeneous distribution
 of degree $-n-|\alpha|$, 
 the Euler homogeneity operator 
\index{A}{Euler@$E$, Euler homogeneity operator\quad}
$E$ acts
 as the scalar multiplication
 by $-(n+|\alpha|)$, 
 and thus 
\[
  E F
  = 
  -\sum_{\alpha \in {\mathbb{N}}^n} 
  (n+|\alpha|)
   a_{\alpha} 
   \delta^{(\alpha)}(x_1,\cdots,x_n).  
\]
Since $\{ \delta^{(\alpha)}(x_1,\cdots,x_n) \}_{\alpha \in {\mathbb{N}}^n}$
 are linearly independent distributions,
 the differential equation \eqref{eqn:Fainv}, 
 namely, 
 $E F = (\lambda - \nu - n) F$
 implies that
\[
  a_{\alpha}=0
\qquad
\text{whenever }
  -n-|\alpha| \ne \lambda - \nu - n.  
\]
Thus we conclude:
\begin{enumerate}
\item[(1)]
If $\nu - \lambda \not \in {\mathbb{N}}$, 
 we get $a_{\alpha}=0$ 
 for all $\alpha \in {\mathbb{N}}^n$, 
 whence $F=0$.  
\item[(2)]
If $\nu - \lambda \in {\mathbb{N}}$, 
 then $a_{\alpha}$ can survive 
 only when $|\alpha|= \nu - \lambda$.  
Then $F(-x)=(-1)^{|\alpha|}F(x)=(-1)^{\nu - \lambda}F(x)$
 because $\delta(x)=\delta(-x)$.  
\end{enumerate}
Therefore Lemma \ref{lem:distpos} is proved.  
\end{proof}

\begin{proof}
[Proof of Theorem \ref{thm:vanishDiff}]
Immediate from the characterization 
 of differential symmetry breaking operators
 (Proposition \ref{prop:Tpair} (3)) and from Lemma \ref{lem:distpos}.  
\end{proof}

\subsection{Upper estimate of the multiplicities}
\label{subsec:SDone}
We recall from the general theory 
 \cite{xKOfm}
 that there exists a constant $C >0$
 such that
\begin{equation}
\label{eqn:IJdimC}
  \dim_{\mathbb{C}}{\operatorname{Hom}}_{G'}
  (I_{\delta}(V,\lambda)|_{G'}, J_{\varepsilon}(W,\nu))\le C
\end{equation}
for any $(\sigma,V) \in \widehat {O(n)}$, 
 $(\tau,W) \in \widehat {O(n-1)}$, 
 $\delta, \varepsilon \in \{\pm\}$, 
 and $(\lambda,\nu) \in {\mathbb{C}}^2$. 
Moreover,
 we also know
 that the left-hand side of \eqref{eqn:IJdimC} is either 0 or 1
 if both the $G$-module $I_{\delta}(V,\lambda)$
 and the $G'$-module $J_{\varepsilon}(W,\nu)$ are irreducible
 \cite{SunZhu}.  
In this section,
 we give a more precise upper estimate
 of the dimension of (continuous) symmetry breaking operators
 by that of {\it{differential}} symmetry breaking operators.  
Owing to the 
\index{B}{dualitytheorem@duality theorem, between differential symmetry breaking operators and Verma modules|textbf}
\lq\lq{duality theorem}\rq\rq\
 (see \cite[Thm.~2.9]{KP1}, 
 see also Fact \ref{fact:dualityDSV} in the next chapter), 
 the latter object can be studied algebraically
 as a branching problem for generalized Verma modules, 
 and is completely classified in \cite{KKP}
 when $(V,W)=(\Exterior^i({\mathbb{C}}^n), \Exterior^j({\mathbb{C}}^{n-1}))$.  
The proof for the upper estimate leads us 
 to complete the proof of a 
\index{B}{localnesstheorem@localness theorem}
 localness theorem
 (Theorem \ref{thm:152347}), 
 namely,
 a sufficient condition
 for all symmetry breaking operators
 to be differential operators. 

\begin{theorem}
[upper estimate of dimension]
\label{thm:SDone}
For any $V \in \widehat{O(n)}$,  
 $W \in \widehat{O(n-1)}$, 
 $\delta, \varepsilon \in \{ \pm \}$, 
 and $(\lambda, \nu)\in {\mathbb{C}}^2$, 
 we have
\begin{equation*}
   \dim_{\mathbb{C}}
   \operatorname{Hom}_{G'}
  (
   I_{\delta}(V,\lambda)|_{G'}, J_{\varepsilon}(W,\nu))
   \le 
   1+ 
   \dim_{\mathbb{C}}
  \operatorname{Diff}_{G'}
  (
   I_{\delta}(V,\lambda)|_{G'}, 
   J_{\varepsilon}(W,\nu)).  
\end{equation*}
\end{theorem}
\begin{proof}
Let $({\mathcal{T}}_{\infty},{\mathcal{T}})$ be the pair
 of distribution kernels 
 of a symmetry breaking operator ${\mathbb{T}}$ as in Proposition \ref{prop:Tpair}.  
Then the first projection
 ${\mathbb{T}} \mapsto {\mathcal{T}}_{\infty}$
 induces an exact sequence:
\begin{equation*}
0 \to
 \operatorname{Diff}_{G'}
  (
   I_{\delta}(V, \lambda)|_{G'}, 
   J_{\varepsilon}(W,\nu)
  )
  \to 
 \operatorname{Hom}_{G'}
  (
   I_{\delta}(V, \lambda)|_{G'}, 
   J_{\varepsilon}(W,\nu)
  )
\to 
{\mathbb{C}}(\Attcal \lambda \nu {\delta \varepsilon} {V,W})_{\infty}, 
\end{equation*}
 by Proposition \ref{prop:Tpair} (3)
 and Proposition \ref{prop:20150828-1231}.  
Thus Theorem \ref{thm:SDone} is proved.  
\end{proof}

We are ready to prove a localness theorem 
  stated in Theorem \ref{thm:152347}.  
\begin{proof}
[Proof of Theorem \ref{thm:152347}]
If $[V:W]=0$
 then 
$
   ({\mathcal{D}}'({\mathbb{R}}^n) 
    \otimes 
    \operatorname{Hom}_{\mathbb{C}}(V_{\lambda,\delta}, W_{\nu, \varepsilon}))^{P'}=\{0\}
$
 by Proposition \ref{prop:20150828-1231}
 because $\pr V W =0$.  
Hence we get Theorem \ref{thm:152347}
 by the exact sequence
 in the above proof.  
\end{proof}

We also prove a part of Theorem \ref{thm:unique}, 
 a generic uniqueness result.  
\begin{corollary}
\label{cor:160150upper}
Suppose $(\sigma, V) \in \widehat{O(n)}$, 
 $(\tau, W) \in \widehat{O(n-1)}$, 
 $\delta, \varepsilon \in \{\pm \}$, 
 and $(\lambda,\nu)\in {\mathbb{C}}^2$.  
If $(\lambda, \nu, \delta, \varepsilon)$
 satisfies the generic parameter condition
 \eqref{eqn:nlgen}, 
 namely, 
 if $\nu - \lambda \not \in 2{\mathbb{N}}$
 for $\delta \varepsilon=+$, 
 or $\nu - \lambda \not \in 2{\mathbb{N}}+1$
 for $\delta \varepsilon=-$, 
then 
\[
  \dim_{\mathbb{C}} 
  \operatorname{Hom}_{G'}
  (I_{\delta}(V,\lambda)|_{G'},J_{\varepsilon}(W,\nu))
  \le 1.  
\]
\end{corollary}

\begin{proof}
[Proof of Corollary \ref{cor:160150upper}]
Owing to Theorem \ref{thm:SDone}, 
 we obtain Corollary \ref{cor:160150upper} by 
 Theorem \ref{thm:vanishDiff}.  
\end{proof}

We shall see 
 that the inequality in Corollary \ref{cor:160150upper}
 is actually the equality 
 by showing the lower estimate
 of the multiplicities in Theorem \ref{thm:existSBO}
 below.  
\subsection{Proof of Theorem \ref{thm:152389}: Analytic continuation of 
symmetry breaking operators
 $\Atbb \lambda \nu {\pm} {V, W}$}
\label{subsec:holoAVW}
The goal of this section
 is to complete the proof of Theorem \ref{thm:152389}
 about the analytic continuation of $\Atbb \lambda \nu {\pm}{V,W}$.  
For $(\sigma,V) \in \widehat {O(n)}$ and 
 $(\tau,W) \in \widehat{O(n-1)}$ such that $[V:W]\ne 0$
 and for $\delta, \varepsilon \in \{\pm\}$, 
 we set $\gamma =\delta \varepsilon$
 and construct a family
 of matrix-valued symmetry breaking operators,
 to be denoted by
\[
  \Atbb \lambda \nu {\gamma} {V, W} \colon
  I_{\delta}(V,\lambda) \to J_{\varepsilon}(W,\nu), 
\]
 which are initially defined
 for $\operatorname{Re} \lambda \gg |\operatorname{Re} \nu|$
 in Lemma \ref{lem:Astep1}.  
We show that they have a holomorphic continuation 
 to the entire plane $(\lambda, \nu) \in {\mathbb{C}}^2$, 
 and thus complete the proof of Theorem \ref{thm:152389}.

\vskip 0.8pc
Here is a strategy.  
\par\noindent
{\bf{Step 0.}}
(distribution kernel near infinity)\enspace

We define ${\operatorname{Hom}}_{\mathbb{C}}(V,W)$-valued distributions
 $(\Atcal \lambda \nu {\gamma} {V, W})_{\infty}$
 on ${\mathbb{R}}^n$ as a multiplication 
 of $(\Attcal \lambda \nu {\gamma} {V, W})_{\infty}$
 (see \eqref{eqn:AVW+u} and \eqref{eqn:AVW-u})
 by appropriate holomorphic functions
 of $\lambda$ and $\nu$
 (Section \ref{subsec:Ainfty}).  
The distributions $(\Atcal \lambda \nu \gamma{V,W})_{\infty}$ depend holomorphically
 on $(\lambda,\nu)$ in the entire plane ${\mathbb{C}}^2$
 (but may vanish at special $(\lambda,\nu)$).  

\vskip 0.8pc
\par\noindent
{\bf{Step 1.}}\enspace
(very regular case)
\enspace
For $\operatorname{Re} \lambda \gg |\operatorname{Re} \nu|$, 
 we define ${\operatorname{Hom}}_{\mathbb{C}}(V,W)$-valued,
 locally integrable functions
 $\Atcal \lambda \nu {\pm} {V,W}$ on ${\mathbb{R}}^n$
 such that the restriction $\Atcal \lambda \nu {\pm} {V,W}|_{\mathbb{R}^n \setminus \{0\}}$
 satisfies the compatibility condition \eqref{eqn:Tpatching}.  
We then prove that the pair
$
((\Atcal \lambda \nu {\gamma} {V,W})_{\infty}, 
  \Atcal \lambda \nu {\gamma} {V,W})
$
 belongs to 
$
   ({\mathcal{D}}'(G/P, {\mathcal{V}}_{\lambda,\delta}^{\ast}) \otimes W_{\nu,\varepsilon})^{\Delta(P')}
$
 for $\delta \varepsilon = \gamma$
 if ${\operatorname{Re}} \lambda \gg |{\operatorname{Re}} \nu|$.

\vskip 0.8pc
\par\noindent
{\bf{Step 2.}}\enspace
(meromorphic continuation and possible poles of $\Atcal \lambda \nu \pm {V,W}$)\enspace
We find polynomials $p_{\gamma}^{V,W}(\lambda,\nu)$
 such that $p_{\gamma}^{V,W}(\lambda,\nu) \Atcal \lambda \nu \gamma {V,W}$
 is a family of distributions on ${\mathbb{R}}^n$
 that depend {\it{holomorphically}} on $(\lambda, \nu) \in {\mathbb{C}}^2$
 (see Proposition \ref{prop:20151209}).

\vskip 0.8pc
\par\noindent
{\bf{Step 3.}}\enspace
(holomorphic continuation of $\Atcal \lambda \nu \pm {V,W}$)\enspace
We prove that there are actually no poles
 of the distributions $\Atcal \lambda \nu \gamma {V,W}$
 by inspecting the residue formula
 of the {\it{scalar-valued}} symmetry breaking operators
 and the zeros
 of the polynomials $p_{\gamma}^{V,W}(\lambda,\nu)$.  
Thus $\Atcal \lambda \nu \gamma {V,W}$ are
 distributions on ${\mathbb{R}}^n$
 that depend holomorphically on $(\lambda,\nu) \in {\mathbb{C}}^2$.

Thus the pair $((\Atcal \lambda \nu {\gamma} {V,W})_{\infty}, \Atcal \lambda \nu {\gamma} {V,W})$
 gives an element
 of ${\mathcal{D}}'(G/P, {\mathcal{V}}_{\lambda,\delta}^{\ast}) \otimes W_{\nu,\varepsilon}$
 for $\delta \varepsilon =\gamma$
 which is invariant under the diagonal action of $P'$, 
 yielding a regular symmetry breaking operator $\Atbb \lambda \nu {\gamma} {V,W}$
 that depends holomorphically on $(\lambda,\nu) \in {\mathbb{C}}^2$
 by Proposition \ref{prop:Tpair}.

\vskip 1pc
The key idea for Steps 1 and 2
 is a reduction 
 to {\it{scalar-valued}} symmetry breaking operators
 which will be discussed
 in Section \ref{subsec:shiftA}
 (Lemma \ref{lem:152392}).

\subsubsection{Normalized distributions
 $(\Atcal \lambda \nu {\gamma} {V,W})_{\infty}$
 at infinity}
\label{subsec:Ainfty}
This is for Step 0.  
We note that the map ${\mathbb{T}} \mapsto {\mathcal{T}}_{\infty}$
 in Proposition \ref{prop:Tpair}
 is neither injective nor surjective
 in general.  
In particular,
 the nonzero distribution $(\Attcal \lambda \nu \pm {V,W})_{\infty}$
 on ${\mathbb{R}}^n$
 (see \eqref{eqn:AVW+u} and \eqref{eqn:AVW-u})
 does not always extend
 to the compactification $G/P$ as an element 
in $({\mathcal{D}}'(G/P, {\mathcal{V}}_{\lambda,\delta}^{\ast}) \otimes W_{\nu,\varepsilon})^{\Delta(P')}$, 
 see Proposition \ref{prop:20170213}.  
However,
 we shall see in Section \ref{subsec:pfAholo}
 that the following {\it{renormalization}} extends to a distribution 
on the compact manifold $G/P$
 for any $\lambda$, $\nu \in {\mathbb{C}}$.  
\begin{alignat*}{2}
  (\Atcal \lambda \nu + {V,W})_{\infty}
  :=& \frac{1}{\Gamma(\frac{\lambda-\nu}{2})}(\Attcal \lambda \nu + {V,W})_{\infty}
  &&= \frac{1}{\Gamma(\frac{\lambda-\nu}{2}) \Gamma(\frac{\lambda+\nu-n+1}{2})}
       |x_n|^{\lambda + \nu -n} \pr V W, 
\\
  (\Atcal \lambda \nu - {V,W})_{\infty}
  :=& \frac{1}{\Gamma(\frac{\lambda-\nu+1}{2})}(\Attcal \lambda \nu - {V,W})_{\infty}
  &&= \frac{1}{\Gamma(\frac{\lambda-\nu+1}{2}) \Gamma(\frac{\lambda+\nu-n+2}{2})}
   |x_n|^{\lambda + \nu -n} \operatorname{sgn} x_n\pr V W.  
\end{alignat*}
By definition, 
 $(\Atcal \lambda \nu \pm {V,W})_{\infty}$ are distributions
 on ${\mathbb{R}}^n$ 
 that depend holomorphically on $(\lambda, \nu)$ in the entire ${\mathbb{C}}^2$.  
Inspecting the poles of $\Gamma(\frac{\lambda-\nu}{2})$ and 
$\Gamma(\frac{\lambda-\nu+1}{2})$, 
 we immediately have the following:

\begin{lemma}
\label{lem:ABzero}
Suppose $[V:W] \ne 0$.  
Then, 
$(\Atcal \lambda \nu + {V,W})_{\infty}=0$
 if and only if $\nu - \lambda \in 2 {\mathbb{N}}$;
 $(\Atcal \lambda \nu - {V,W})_{\infty}=0$
 if and only if $\nu - \lambda \in 2 {\mathbb{N}}+1$.  
\end{lemma}

\subsubsection{Preliminary results in the scalar-valued case}
\label{subsec:shiftA}

As we have seen 
 in Section \ref{subsec:Ainfty}, 
 the analytic continuation 
 of the distribution $(\Atcal \lambda \nu {\gamma} {V,W})_{\infty}$
 at infinity is easy.  
In order to deal with the nontrivial case,
 {\it{i.e.,}}
 the distribution kernel $\Atcal \lambda \nu {\gamma} {V,W}$
 near the origin, 
 we begin with some basic properties
 of the {\it{scalar-valued}} symmetry breaking operators.  
We recall from \cite[(7.8)]{sbon} 
 that the (scalar-valued) distribution kernels
 $\Atcal \lambda \nu \pm {} \in {\mathcal{D}}'({\mathbb{R}}^n)$
 are initially defined
 as locally integrable functions on ${\mathbb{R}}^n$
 by 
\index{A}{Act1@$\Atcal \lambda \nu + {}$}
\index{A}{Act2@$\Atcal \lambda \nu - {}$}
\begin{align}
\Atcal \lambda \nu + {} (x, x_n)
=& \frac{1}{\Gamma (\frac{\lambda+\nu-n+1}{2}) \Gamma (\frac{\lambda-\nu}{2})}
    (|x|^2+ x_n^2)^{-\nu} |x_n|^{\lambda + \nu -n},
\label{eqn:KAlnn+}
\\
\Atcal{\lambda}{\nu}{-}{} (x, x_n)
=& \frac{1}
   {\Gamma (\frac{\lambda+\nu-n+2}{2}) \Gamma (\frac{\lambda-\nu+1}{2})}
   (|x|^2+ x_n^2)^{-\nu} |x_n|^{\lambda+\nu-n}\operatorname{sgn} x_n,
\label{eqn:KAlnn-}
\end{align}
respectively
 for $\operatorname{Re}\lambda \gg |\operatorname{Re} \nu|$.  
(In \cite{sbon}, 
we used the notation $\widetilde K_{\lambda,\nu}^{\mathbb{A}}$
 for the scalar-valued distribution kernel $\Atcal{\lambda}{\nu}{+}{}$.)
More precisely,
 we have:
\begin{fact}
[{\cite[Chap.~7]{sbon}}]
\label{fact:sbonA}
$\Atcal{\lambda}{\nu}{\pm}{}$ are locally integrable 
 on ${\mathbb{R}}^n$
 if $\operatorname{Re}\left(\lambda - \nu\right)>0$
 and $\operatorname{Re}\left(\lambda + \nu\right)>n-1$, 
 and extend as distributions
 on ${\mathbb{R}}^n$
 that depend holomorphically on $\lambda$, $\nu$
 in the entire $(\lambda,\nu) \in {\mathbb{C}}^2$.  
\end{fact}
The distributions $\Atcal {\lambda}{\nu}{+}{}$ were
 thoroughly studied in \cite[Chap. 7]{sbon}, 
 and analogous results for $\Atcal \lambda \nu -{}$
 can be proved exactly in the same way.

We introduce polynomials 
 $p_{\pm,N}(\lambda,\nu)$ of the two-variables $\lambda$ and $\nu$
 by 
\index{A}{p1N@$p_{+,N}(\lambda,\nu)$|textbf}
\index{A}{p2N@$p_{-,N}(\lambda,\nu)$|textbf}
\begin{alignat}{2}
p_{+,N}(\lambda,\nu):=& \prod_{j=1}^{N} (\lambda-\nu-2j)
\quad
&&\text{for $N \in {\mathbb{N}}_+$}, 
\label{eqn:p+N}
\\
p_{-,N}(\lambda,\nu):=& (\lambda+\nu-n)\prod_{j=0}^{N} (\lambda-\nu-1-2j)
\quad
&&\text{for $N \in {\mathbb{N}}$}.  
\label{eqn:p-N}
\end{alignat}
We use a trick to raise the regularity 
 of the distribution $\Atcal \lambda \nu +{}(x,x_n)$
 at the origin 
 by shifting the parameter.  
The resulting distributions are under control
 by the polynomials $p_{\pm,N}(\lambda,\nu)$ as follows:
\begin{lemma}
\label{lem:pAshift}
We have the following identities
 as distributions on ${\mathbb{R}}^n$:
\begin{align*}
p_{+,N}(\lambda, \nu) \Atcal \lambda \nu + {} (x,x_n)
=&
2^{N} (|x|^2+x_n^2)^N \Atcal {\lambda-N}{\nu+N}{+}{} (x,x_n), 
\\
p_{-,N}(\lambda, \nu) \Atcal \lambda \nu - {} (x,x_n)
=&
2^{N+2} (|x|^2+x_n^2)^N x_n \Atcal {\lambda-N-1}{\nu+N}{+}{} (x,x_n). 
\end{align*}
\end{lemma}

\begin{proof}
For $\operatorname{Re}\lambda \gg  |\operatorname{Re}\nu|$, 
 we have from the definition \eqref{eqn:KAlnn+}, 
\begin{align*}
  (|x|^2+ x_n^2)^N \Atcal {\lambda-N}{\nu+N}{+}{} (x,x_n)
=&
\frac{\Gamma(\frac{\lambda-\nu}{2})}{\Gamma(\frac{\lambda-\nu}{2}-N)}
\Atcal{\lambda}{\nu}{+}{} (x,x_n)
\\
=&
\frac{1}{2^N} p_{+,N}(\lambda,\nu) \Atcal \lambda \nu +{}(x,x_n).  
\end{align*}
Since both sides depend holomorphically on $(\lambda, \nu)\in {\mathbb{C}}^2$, 
 we get the first assertion.  
The proof of the second assertion goes similarly.  
\end{proof}

\begin{lemma}
\label{lem:1.16}
If $(\lambda,\nu) \in {\mathbb{C}}^2$ satisfies $p_{+,N}(\lambda, \nu)=0$, 
then 
\[
  h(x,x_n) \Atcal {\lambda-N}{\nu+N}{+}{}=0
\quad
 \text{in } {\mathcal{D}}'
   ({\mathbb{R}}^n), 
\]
for all homogeneous polynomials $h(x,x_n)$
 of degree $2N$.  
\end{lemma}

\begin{proof}
It follows from $p_{+, N}(\lambda, \nu)=0$
 that 
\[
   (\nu+N)-(\lambda-N) \in \{0,2,4,\cdots,2N-2\}.  
\]
By the residue formula
 of the scalar-valued symmetry breaking operator
 $\Atcal {\lambda'} {\nu'} {+} {}$ 
 (see \cite[Thm.~12.2 (2)]{sbon}), 
 we have 
\[
  \Atcal {\lambda-N}{\nu+N}{+}{} = q \, \Ctcal{\lambda-N}{\nu+N}{}
\]
 for some constant 
\index{A}{qAC@$q_C^A$}
$q\equiv q_C^A(\lambda-N,\nu+N)$
 depending on $\lambda-N$ and $\nu+N$.  
Since $\Ctcal{\lambda-N}{\nu+N}{}$ is a distribution
 of the form $D \delta(x_1, \cdots, x_n)$
 where 
\index{A}{CGegenbauernorm@$\widetilde C_l^{\alpha}(z)$, normalized Gegenbauer polynomial}
$
   D=\widetilde C_{2N-2j}^{\lambda - N - \frac{n-1}{2}}(-\Delta_{\mathbb{R}^{n-1}}, \frac{\partial}{\partial x_n})$ is a differential operator
 of homogeneous degree $2N-2j (< 2N)$, 
 see \eqref{eqn:Gegentwo}, 
 an iterated use of the Leibniz rule shows
\[
   \text{$h(x,x_n)D \delta(x_1, \cdots, x_n)=0$
 in ${\mathcal{D}}'
   ({\mathbb{R}}^n)$}
\]
 for any homogeneous polynomial $h(x,x_n)$
 of degree $2N$.  
\end{proof}

\begin{lemma}
\label{lem:152293-copy}
If $(\lambda,\nu) \in {\mathbb{C}}^2$ satisfies $p_{-,N}(\lambda,\nu)=0$, 
then 
\[
   x_n h(x,x_n) \Atcal {\lambda-N-1}{\nu+N}{+}{} (x,x_n) =0
\quad\text{ in }
{\mathcal{D}}'({\mathbb{R}}^n)
\]
for all homogeneous polynomial $h(x,x_n)$ of degree $2N$.  
\end{lemma}
\begin{proof}
It follows from $p_{-, N}(\lambda, \nu)=0$
 that $(\nu+N)-(\lambda-N-1) \in \{0,2,\cdots,2N\}$
 or $(\lambda-N-1)+(\nu+N) =n-1$.  
By using again the residue formula
 of the scalar-valued symmetry breaking operator
 $\Atcal {\lambda'}{\nu'}{+}{}$
 in \cite[Thm.~12.2]{sbon}, 
 we see that the distribution kernel $\Atcal {\lambda-N-1}{\nu+N}{+}{}(x,x_n)$ is
 a scalar multiple
 of the following distributions:
\begin{alignat*}{2}
& \delta(x_n)
\qquad
&&\text{if } \lambda + \nu=n, 
\\
& D \delta(x_1, \cdots,x_n)
\qquad
&&\text{if } \lambda - \nu =2j+1 \qquad (0 \le j \le N),   
\end{alignat*}
where 
$
    D 
    = 
    {\widetilde C}_{2N-2j}^{\lambda-N-1-\frac{n-1}2}(-\Delta_{\mathbb{R}^{n-1}}, \frac \partial{\partial x_n})
$
 is a differential operator
 of homogeneous degree $2N-2j$
($< 2N+1$).  
Then the multiplication by a homogeneous polynomial $x_n h(x,x_n)$
 of degree $2N+1$
 annihilates these distributions.  
Hence the lemma follows.  
\end{proof}

\subsubsection{Step 1. Very regular case}
\label{subsec:Astep1}

We recall from \eqref{eqn:RVW}
 that 
\index{A}{RVW@$\Rij VW$}
\index{A}{1psin@$\psi_n$}
$
\Rij VW = \pr VW \circ \sigma \circ \psi_n$
$\in C^{\infty}({\mathbb{R}}^n \setminus \{0\}) \otimes
 {\operatorname{Hom}}_{\mathbb{C}}(V,W)$.  
For $\operatorname{Re}\lambda \gg |\operatorname{Re}\nu|$, 
 we define
$\Atcal \lambda \nu \pm {V,W}
 \in C({\mathbb{R}}^n \setminus \{0\}) \otimes \operatorname{Hom}_{\mathbb{C}}(V,W)$
 by
\begin{alignat*}{2}
\Atcal \lambda \nu + {V,W}
:=&
\frac{1}{\Gamma(\frac{\lambda+\nu-n+1}{2})\Gamma(\frac{\lambda-\nu}{2})}
(|x|^2+x_n^2)^{-\nu} |x_n|^{\lambda+\nu-n}
\Rij VW (x,x_n), 
\\
\Atcal \lambda \nu - {V,W}
:=&
\frac{1}{\Gamma(\frac{\lambda+\nu-n+2}{2})\Gamma(\frac{\lambda-\nu+1}{2})}
(|x|^2+x_n^2)^{-\nu} |x_n|^{\lambda+\nu-n}
\operatorname{sgn} x_n
\Rij VW (x,x_n), 
\end{alignat*}
 (see \eqref{eqn:KVWpt} and \eqref{eqn:KVWmt}), 
 respectively.  
The goal of this section
 is to prove the following lemma
 in the matrix-valued case  
 for $\operatorname{Re} \lambda \gg |\operatorname{Re} \nu|$.  
\begin{lemma}
\label{lem:Astep1}
Let $(\sigma,V) \in \widehat {O(n)}$, 
 $(\tau,W) \in \widehat {O(n-1)}$
 and $\delta$, $\varepsilon \in \{ \pm \}$.  
Suppose $\operatorname{Re}\left(\lambda -\nu\right) >0$
 and $\operatorname{Re}\left(\lambda + \nu\right) > n-1$.  

\begin{enumerate}
\item[{\rm{(1)}}]
$\Atcal \lambda \nu \pm {V,W}$ are
 $\operatorname{Hom}_{\mathbb{C}}(V,W)$-valued
 locally integrable functions on ${\mathbb{R}}^n$.  

\item[{\rm{(2)}}]
The pair $((\Atcal \lambda \nu {\delta \varepsilon} {V,W})_{\infty}, \Atcal \lambda \nu {\delta \varepsilon} {V,W})$ defines an element
 of $({\mathcal{D}}'(G/P, {\mathcal{V}}_{\lambda,\delta}^{\ast}) \otimes W_{\nu,\varepsilon})^{\Delta(P')}$, 
 and thus yield a symmetry breaking operator
 $\Atbb \lambda \nu {\delta \varepsilon} {V,W}\colon I_{\delta}(V,\lambda) \to J_{\varepsilon}(W,\nu)$.  

\end{enumerate}
\end{lemma}

\begin{proof}
We fix inner products on $V$ and $W$
 that are invariant by $O(n)$ and $O(n-1)$, 
 respectively.  
\index{A}{opnor@$\mid\mid \cdot \mid\mid_{\operatorname{op}}$, operator norm}
Let $\| \cdot \|_{\operatorname{op}}$
 denote the operator norm
 for linear maps 
 between (finite-dimensional) Hilbert spaces.  
In view of the definition $\Rij V W = \pr V W \circ \sigma \circ \psi_n$
 (see \eqref{eqn:RVW}), 
 we have
\[
   \| \Rij V W(x,x_n)  \|_{\operatorname{op}}
  \le
  \| \sigma \circ \psi_n (x,x_n) \|_{\operatorname{op}}
  = 1
\quad
  \text{for all }
  (x,x_n) \in {\mathbb{R}}^n \setminus \{0\}.  
\]
Hence the first statement is reduced to the scalar case
 as stated in Fact \ref{fact:sbonA}.

The compatibility condition \eqref{eqn:Tpatching}
 can be verified readily from the definition
 of $(\Atcal \lambda \nu \pm {V,W})_{\infty}$
 and $\Atcal \lambda \nu \pm {V,W}$.  
Hence the pair
$
   ((\Atcal \lambda \nu {\delta\varepsilon} {V,W})_{\infty}, 
   \Atcal \lambda \nu {\delta\varepsilon} {V,W})
$
 defines an element of 
 ${\mathcal{D}}'(G/P, {\mathcal{V}}_{\lambda,\delta}^{\ast}) \otimes W_{\nu, \varepsilon}$
 by Lemma \ref{lem:152341}.  
The invariance under the diagonal action of $P'$
 follows from Proposition \ref{prop:20150828-1231}
 for $(\Atcal \lambda \nu {\delta\varepsilon}{V,W})_{\infty}$
 and from a direct computation 
 for $\Atcal \lambda \nu {\delta\varepsilon}{V,W}$
 when $\operatorname{Re} \lambda \gg |\operatorname{Re}\nu|$
 because 
both $(\Atcal \lambda \nu {\delta\varepsilon}{V,W})_{\infty}$
 and $\Atcal \lambda \nu {\delta\varepsilon} {V,W} \in L_{\operatorname{loc}}^1({\mathbb{R}}^n)$.  
\end{proof}

\subsubsection{Step 2. Reduction to the scalar-valued case}
We shall prove:
\begin{proposition}
\label{prop:20151209}
Let $(\sigma,V) \in\widehat{O(n)}$
 and $(\tau,W) \in\widehat{O(n-1)}$.  
Then the distributions  
$\Atcal \lambda \nu {\pm} {V,W}$, 
 initially defined
 as an element of $L_{\operatorname{loc}}^1({\mathbb{R}}^n) \otimes \operatorname{Hom}_{\mathbb{C}}(V,W)$
 for $\operatorname{Re} \lambda \gg |\operatorname{Re} \nu|$
 in Lemma \ref{lem:Astep1}, 
 extend meromorphically 
 in the entire plane $(\lambda,\nu) \in {\mathbb{C}}^2$.  
\end{proposition}

In order to prove Proposition \ref{prop:20151209}, 
 we need to control the singularity
 of $\sigma \circ \psi_n \in C^{\infty}({\mathbb{R}}^n \setminus \{0\}) \otimes
 {\operatorname{End}}_{\mathbb{C}}(V)$
 at the origin.  
We formulate a necessary lemma:
\begin{lemma}
\label{lem:1.12}
For any irreducible representation $(\sigma, V)$ of $O(n)$, 
 there exists $N \in {\mathbb{N}}$ 
such that
\[
  g(x,x_n):= (|x|^2+x_n^2)^{N} \sigma(\psi_n(x,x_n))
\]
 is an $\operatorname{End}(V)$-valued homogeneous polynomial of degree $2N$.  
\end{lemma}

\begin{definition}
\label{def:Nsigma}
For $\sigma \in \widehat {O(n)}$, 
 we denote by 
\index{A}{Nsigma@$N(\sigma)$|textbf}
$N(\sigma)$
  the smallest integer $N$
 satisfying the conclusion of Lemma \ref{lem:1.12}.  
\end{definition}

We prove Lemma \ref{lem:1.12}
 by showing the following estimate of the integer $N(\sigma)$.  
Let $\ell(\sigma)$ be as defined in \eqref{eqn:ONlength}.  

\begin{lemma}
\label{lem:Nsigma}
\index{A}{lengthrep@$\ell(\sigma)$}
$N(\sigma) \le \ell(\sigma)$
 for all $\sigma \in \widehat {O(n)}$.  
\end{lemma}

\begin{proof}
[Proof of Lemma \ref{lem:Nsigma}]
Suppose 
\index{A}{0tExterior@$\Lambda^+(O(N))$}
$(\sigma_1, \cdots, \sigma_n) \in \Lambda^+(O(n))$, 
 and let $(\sigma,V)$ be the irreducible finite-dimensional 
 representation $\Kirredrep {O(n)}{\sigma_1, \cdots, \sigma_n}$
 of $O(n)$
 via the Cartan--Weyl isomorphism
 \eqref{eqn:CWOn}.  
It is convenient to set $\sigma_{n+1}=0$.  
Since the exterior representations
 of $GL(n,{\mathbb{C}})$ on $\Exterior^j({\mathbb{C}}^n)$
 have 
 highest weights $(1,\cdots,1,0,\cdots,0)$, 
 and since
\[
   \sum_{j=1}^{n}(\sigma_j-\sigma_{j+1})
   (\underbrace{1,\cdots,1}_j,0,\cdots,0)
  =(\sigma_1,\cdots,\sigma_n), 
\]
we can realize the irreducible representation of $GL(n,{\mathbb{C}})$
 with highest weight $(\sigma_1,\cdots,\sigma_n)$
 as a subrepresentation
 of the tensor product representation
\[
  \bigotimes_{j=1}^{n} (\Exterior^j({\mathbb{C}}^n))^{\sigma_j-\sigma_{j+1}}.  
\]
This is a polynomial representation of homogeneous degree
\[
   \sum_{j=1}^n j (\sigma_j-\sigma_{j+1})=\sum_{j=1}^n \sigma_j.  
\]
We set $N:=\sum_{j=1}^n \sigma_j$.  
Then the matrix coefficients of this $GL(n,{\mathbb{C}})$-module
 are given by homogeneous polynomials
 of degree $N$
 of $z_{ij}$
 ($1 \le i,j \le n$)
 where $z_{ij}$ are the coordinates
 of $GL(n,{\mathbb{C}})$.  
Since the representation $(\sigma, V)$ of $O(n)$
 arises as a subrepresentation
 of this $GL(n,{\mathbb{C}})$-module, 
 the formula \eqref{eqn:psim}
 of $\psi_n$ shows that the matrix coefficients of $\sigma(\psi_n(x,x_n))$
 is a polynomial of $x$ and $x_n$ 
 after multiplying $(|x|^2 + x_n^2)^N$.

We note
 that $\det \psi_n(x,x_n)=-1$
 for all $(x,x_n) \in {\mathbb{R}}^n \setminus \{0\}$
 by \eqref{eqn:psix}.  
Therefore, 
 we may assume
 that $(\sigma,V)$ is of 
\index{B}{type1@type I, representation of $O(N)$}
type I
 by \eqref{eqn:type1to2}, 
 namely, 
 $\sigma_{k+1}=\cdots=\sigma_n$
 for some $k$
 with $2k \le n$.  
In this case $N=l(\sigma)$
 by the definition \eqref{eqn:ONlength}.  
By \eqref{eqn:lsigma}, 
 we have shown the lemma.  
\end{proof}

The estimate in Lemma \ref{lem:Nsigma} is not optimal.

\begin{example}
\label{ex:Nsigma}
\begin{enumerate}
\item[{\rm{1)}}]
$N(\sigma)=0$ if $(\sigma, V)$ is a one-dimensional representation.  
\item[{\rm{2)}}]
$N(\sigma)=1$ if $\sigma$ is the exterior representation
 on $V= \Exterior^i({\mathbb{C}}^n)$
 ($1 \le i \le n-1$).  
See \eqref{eqn:exrep}
 and Lemma \ref{lem:psidet} (2)
 for the proof.  
\end{enumerate}
\end{example}

Let $N \equiv N(\sigma) \in {\mathbb{N}}$ 
 and $g \in {\operatorname{Pol}}[x_1, \cdots, x_n] \otimes {\operatorname{End}}_{\mathbb{C}}(V)$
 be as in Lemma \ref{lem:1.12}, 
 and $\pr V W \colon V \to W$ be a nonzero $O(n-1)$-homomorphism.  
We define 
$
   g^{V,W} 
   \in{\operatorname{Pol}}[x_1, \cdots, x_n] 
   \otimes 
   {\operatorname{Hom}}_{\mathbb{C}}(V,W)
$
 by
\index{A}{gVW@$g^{V,W}$|textbf}
\begin{equation}
   g^{V,W}
  :=
  \pr V W \circ g.  
\end{equation}
With notation of 
\index{A}{RVW@$\Rij VW$}
$
 \Rij VW
$ as in \eqref{eqn:RVW}, 
 we have 
\begin{align}
\label{eqn:gVW}
g^{V,W}(x,x_n)=&(|x|^2+x_n^2)^N \Rij VW (x,x_n)
\\
\notag
              =&(|x|^2+x_n^2)^N \pr VW \circ \sigma(\psi_n(x,x_n)).  
\end{align}
Then $g^{V,W}$ is a $\operatorname{Hom}_{\mathbb{C}}(V,W)$-valued
 polynomial of homogeneous degree $2N$.

The following lemma will imply 
 that the singularity at the origin
 of the matrix-valued distributions
 $\Atcal \lambda \nu {\pm}{V,W}$
 is under control by the scalar-valued case:
\begin{lemma}
\label{lem:152392}
Suppose $\operatorname{Re} \lambda \gg |\operatorname{Re} \nu|$.  
Let $p_{\pm}(\lambda,\nu)$ be the polynomials of $\lambda$ and $\nu$
 defined in \eqref{eqn:p+N} and \eqref{eqn:p-N}.  
Then, 
\index{A}{p1N@$p_{+,N}(\lambda,\nu)$}
\index{A}{p2N@$p_{-,N}(\lambda,\nu)$}
\begin{align}
\label{eqn:152392a}
p_{+,N}(\lambda,\nu) \Atcal \lambda \nu {+} {V,W}(x,x_n)
=&
2^N \Atcal {\lambda-N} {\nu+N}{+}{}(x,x_n)g^{V,W}(x,x_n), 
\\
\label{eqn:152392b}
p_{-,N}(\lambda,\nu) \Atcal \lambda \nu {-} {V,W}(x,x_n)
=&
2^{N+2} x_n \Atcal {\lambda-N-1} {\nu+N}{+}{}(x,x_n) g^{V,W}(x,x_n). 
\end{align}  
\end{lemma}
\begin{proof}
For $\operatorname{Re} \lambda \gg |\operatorname{Re} \nu|$, 
 both $\Atcal \lambda \nu {\pm} {V,W}$
 and $\Atcal \lambda \nu {\pm} {}$
 are locally integrable in ${\mathbb{R}}^n$.  
By definition, 
 we have
\begin{equation*}
(|x|^2+ x_n^2)^N \Atcal {\lambda}{\nu}{\gamma}{V,W}
 =
\Atcal {\lambda}{\nu}{\gamma}{} (x,x_n) g^{V,W} (x,x_n)
\end{equation*}
for $\gamma=\pm$.  
By Lemma \ref{lem:pAshift}, 
 we have
\[
  (|x|^2+x_n^2)^N
  (p_{+,N}(\lambda,\nu) \Atcal \lambda \nu +{V,W}(x,x_n)
  -2^N \Atcal {\lambda-N} {\nu+N} +{}(x,x_n) g^{V,W}(x,x_n))
  =0.  
\]
Hence we get the equality \eqref{eqn:152392a}
 as ${\operatorname{Hom}}_{\mathbb{C}}(V,W)$-valued 
 locally integrable functions in ${\mathbb{R}}^n$.  
Similarly,
we obtain
\[
  (|x|^2+ x_n^2)^N \Atcal {\lambda-N-1}{\nu+N}{+}{} (x,x_n) x_n
=
\frac{1}{2^{N+2}}
p_{-,N}(\lambda, \nu)
\Atcal{\lambda}{\nu}{-}{} (x,x_n).  
\]
Thus the second statement follows.  
\end{proof}

We are ready to prove the main result
 of this section.  
\begin{proof}
[Proof of Proposition \ref{prop:20151209}]
Since $g^{V,W}(x,x_n)$ is a polynomial of $(x,x_n)=(x_1, \cdots, x_n)$, 
 the multiplication of any distributions
 on ${\mathbb{R}}^n$ 
 by $g^{V,W}$
 is well defined.  
Therefore, 
 the right-hand sides of \eqref{eqn:152392a} and \eqref{eqn:152392b}
 make sense as distributions on ${\mathbb{R}}^n$
 that depend holomorphically in $(\lambda, \nu) \in {\mathbb{C}}^2$.  

Taking their quotients
 by the polynomials $p_{\pm,N}(\lambda, \nu)$, 
 we set
\index{A}{Actg1@$\Atcal {\lambda}{\nu}{+}{V,W}$|textbf}
\index{A}{Actg2@$\Atcal {\lambda}{\nu}{-}{V,W}$|textbf}
\begin{align}
\label{eqn:pAVW+}
\Atcal \lambda \nu {+} {V,W}(x,x_n)
:=&
\frac{2^N}{p_{+,N}(\lambda,\nu)}
\Atcal {\lambda-N} {\nu+N}{+}{}(x,x_n)
g^{V,W}(x,x_n), 
\\
\label{eqn:pAVW-}
\Atcal \lambda \nu {-} {V,W}(x,x_n)
:=&
\frac{2^{N+2}}{p_{-,N}(\lambda,\nu)}
\Atcal {\lambda-N-1} {\nu+N}{+}{}(x,x_n)
x_n
g^{V,W}(x,x_n).  
\end{align}  
Then $\Atcal \lambda \nu {\pm} {V,W}$ are $\operatorname{Hom}_{\mathbb{C}}(V,W)$-valued distributions
 on ${\mathbb{R}}^n$
 which depend meromorphically on $(\lambda,\nu) \in {\mathbb{C}}^2$
 because $\Atcal {\lambda'} {\nu'}{+}{} (x,x_n)$ is a family of scalar-valued distributions 
 on ${\mathbb{R}}^n$ 
 that depend holomorphically on $(\lambda',\nu') \in {\mathbb{C}}^2$
(Fact \ref{fact:sbonA})
 and $g^{V,W}(x,x_n)$ is a polynomial.  
By Lemma \ref{lem:152392}, 
 they coincide locally integrable functions
 on ${\mathbb{R}}^n$
 that are defined in \eqref{eqn:KVWpt} and \eqref{eqn:KVWmt}, 
 respectively, 
 when $\operatorname{Re} \lambda \gg |\operatorname{Re} \nu|$.  
Thus Proposition \ref{prop:20151209} is proved.  
\end{proof}

\subsubsection{Step 3. Proof of holomorphic continuation}

In this section,
 we show
 that there are no poles of $\Atcal \lambda \nu {\pm} {V,W}$.  
\begin{lemma}
\label{lem:Aholo}
$\Atcal \lambda \nu {\pm} {V,W}$ are distributions on ${\mathbb{R}}^n$
 that depend holomorphically on $(\lambda, \nu)\in {\mathbb{C}}^2$.  
\end{lemma}
\begin{proof}
By \eqref{eqn:pAVW+} and \eqref{eqn:pAVW-}, 
 the only possible places that the distribution
 $\Atcal \lambda \nu {\gamma} {V,W}$ may have poles 
are the zeros of the denominators,
 namely,
\begin{alignat*}{2}
p_{+,N}(\lambda,\nu)
&=\prod_{j=1}^N (\lambda - \nu -2j) \qquad 
&& \gamma =+, 
\\
p_{-,N}(\lambda,\nu)
&=
(\lambda + \nu -n) \prod_{j=0}^N (\lambda - \nu -1-2j) 
\qquad
&& \gamma =-, 
\end{alignat*} 
however,
 we have proved
 that they are not actually poles
 by Lemmas \ref{lem:1.16} and \ref{lem:152293-copy}, 
 respectively.  
Hence $\Atcal \lambda \nu {\gamma} {V,W}$ are distributions
 that depend holomorphically on $(\lambda,\nu) \in {\mathbb{C}}^2$.  
\end{proof}

\subsubsection{Proof of Theorem \ref{thm:152389}}
\label{subsec:pfAholo}

We are ready to prove
 that the matrix-valued symmetry breaking operator
 $\Atbb \lambda \nu {\pm} {V,W}$ has a holomorphic continuation
 in the entire plane $(\lambda,\nu) \in {\mathbb{C}}^2$.  
 
\begin{proof}
[Proof of Theorem \ref{thm:152389}]
Suppose $(\sigma,V) \in \widehat {O(n)}$.  
Let 
\index{A}{Nsigma@$N(\sigma)$}
$
N \equiv N(\sigma) \in {\mathbb{N}}
$ be as in Lemma \ref{lem:1.12}.  
We recall from \eqref{eqn:gVW}
 that the $\operatorname{Hom}_{\mathbb{C}} (V,W)$-valued function 
\index{A}{gVW@$g^{V,W}$}
\begin{equation*}
 g^{V,W}(x,x_n)= (|x|^2+x_n^2)^{N} \operatorname{pr}_{V\to W} \circ\sigma(\psi_n(x,x_n))
\end{equation*}
 is actually a ${\operatorname{Hom}}_{\mathbb{C}}(V,W)$-valued polynomial
 of homogeneous degree $2N$.

We know 
 that the pair $((\Atcal \lambda \nu {\pm} {V,W})_{\infty},\Atcal \lambda \nu {\pm} {V,W})$
 satisfies the following properties:
\begin{enumerate}
\item[{\rm{(1)}}]
$(\Atcal \lambda \nu {\pm} {V,W})_{\infty}$ is
 a $\operatorname{Hom}_{\mathbb{C}}(V,W)$-valued distribution 
 on ${\mathbb{R}}^n$ satisfying \eqref{eqn:Tinfty}
 that depend holomorphically
 in $(\lambda,\nu) \in {\mathbb{C}}^2$.  
\item[{\rm{(2)}}]
$\Atcal \lambda \nu {\pm} {V,W}$
 is a $\operatorname{Hom}_{\mathbb{C}}(V,W)$-valued distribution 
 on ${\mathbb{R}}^n$
 that depend holomorphically
 on $(\lambda, \nu) \in {\mathbb{C}}^2$.  
\item[{\rm{(3)}}]
For $\delta$, $\varepsilon \in \{\pm\}$, 
 $\Atcal \lambda \nu {\delta\varepsilon} {V,W}
\in {\mathcal{S}}ol({\mathbb{R}}^n; V_{\lambda, \delta}, 
 W_{\nu, \varepsilon})$.  
Moreover,
 the conditions \eqref{eqn:Tzero} and \eqref{eqn:Tpatching} are satisfied
 when $\operatorname{Re}\lambda \gg |\operatorname{Re}\nu|$.  
\end{enumerate}

All the equations concerning 
$
{\mathcal{S}}ol({\mathbb{R}}^n; V_{\lambda, \delta}, 
 W_{\nu, \varepsilon})
$
 depend holomorphically
 on $(\lambda, \nu)$ in the entire ${\mathbb{C}}^2$.  
On the other hand,
 for $\gamma \in \{\pm \}$, 
 the properties (1) and (2) tell
 that the pair 
$
   ( (\Atcal \lambda \nu {\gamma} {V,W})_{\infty}, 
       \Atcal \lambda \nu {\gamma} {V,W}
   )
$
 depends holomorphically on $(\lambda, \nu)$ in the entire ${\mathbb{C}}^2$.  
Hence the property (3) holds 
 in the entire $(\lambda, \nu)\in {\mathbb{C}}^2$
 by analytic continuation.  
In turn, 
 Proposition \ref{prop:Tpair} implies
 that the pair
$
   ((\Atcal \lambda \nu {\gamma} {V,W})_{\infty}, 
\Atcal \lambda \nu {\gamma} {V,W})
$
 gives an element of $({\mathcal{D}}'(G/P, {\mathcal{V}}_{\lambda,\delta}^{\ast}) \otimes W_{\nu,\varepsilon})^{\Delta(P')}$
 for all $(\lambda, \nu)\in {\mathbb{C}}^2$, 
 and we have completed the proof of Theorem \ref{thm:152389}.  
\end{proof}

\subsection{Existence condition for regular symmetry breaking operators :
Proof of Theorem \ref{thm:regexist}}
\label{subsec:regexist}
In Theorem \ref{thm:152389}, 
 we have assumed 
 $[V:W] \ne 0$ for the construction 
 of symmetry breaking operators.  
In this section 
 we complete the proof of Theorem \ref{thm:regexist}, 
 which asserts
 that the condition $[V:W] \ne 0$ is necessary
 and sufficient for the existence
 of regular symmetry breaking operators.

Suppose $[V:W]\ne 0$.  
Let $\Atbb \lambda \nu {\delta\varepsilon}{V,W} \colon 
 I_{\delta}(V,\lambda) \to J_{\varepsilon}(W,\nu)$
 be the normalized symmetry breaking operator
 which is obtained by the analytic continuation
 of the integral operator
 in Section \ref{subsec:holoAVW}.  
We study the support
 of its distribution kernel 
 $\Atcal \lambda \nu {\delta\varepsilon}{V,W}$.  
We define subsets $U_{+}^{{\operatorname{reg}}}$ and $U_{-}^{{\operatorname{reg}}}$in ${\mathbb{C}}^2$ by
\begin{align}
\label{eqn:Ureg+}
U_{+}^{{\operatorname{reg}}}:=&
\{(\lambda,\nu) \in {\mathbb{C}}^2
:
 n-\lambda-\nu-1 \not\in 2 {\mathbb{N}}, \nu-\lambda \not\in 2 {\mathbb{N}}\}, 
\\
\label{eqn:Ureg-}
U_{-}^{{\operatorname{reg}}}:=&
\{(\lambda,\nu) \in {\mathbb{C}}^2
:
 n-\lambda-\nu-2 \not\in 2 {\mathbb{N}}, \nu-\lambda-1 \not\in 2 {\mathbb{N}}\}.  
\end{align}
\begin{proposition}
\label{prop:172425}
Suppose $V \in \widehat{O(n)}$ and $W \in \widehat{O(n-1)}$
 satisfy $[V:W]\ne 0$.  
Let $\delta, \varepsilon \in \{\pm\}$.  
Then $\Atbb \lambda \nu{\delta\varepsilon}{V,W}$
 is a nonzero 
\index{B}{regularsymmetrybreakingoperator@regular symmetry breaking operator}
{\it{regular symmetry breaking operator}}
 in the sense of Definition \ref{def:regSBO}
 for all $(\lambda,\nu) \in U_{\delta \varepsilon}^{{\operatorname{reg}}}$.  
\end{proposition}
\begin{proof}[Proof of Proposition \ref{prop:172425}]
As in Proposition \ref{prop:Tpair},
 the distribution kernel of the operator $\Atbb \lambda \nu {\delta\varepsilon}{V,W}$ can be expressed by 
 a pair $((\Atcal \lambda \nu{\delta\varepsilon}{V,W})_{\infty}, \Atcal \lambda \nu{\delta\varepsilon}{V,W})$
 of ${\operatorname{Hom}}_{\mathbb{C}}(V,W)$-valued distributions
 on ${\mathbb{R}}^n$
 corresponding to the open covering
 $G/P=N_+ w P/P \cup N_-P /P$.  
Then it suffices to show 
 ${\operatorname{Supp}}(\Atcal \lambda \nu{\delta\varepsilon}{V,W})_{\infty}
 ={\mathbb{R}}^n$
 for $(\lambda,\nu)\in U_{\delta\varepsilon}^{{\operatorname{reg}}}$.  
If $(\lambda,\nu) \in U_{\delta\varepsilon}^{{\operatorname{reg}}}$, 
then $(\lambda,\nu,\delta,\varepsilon)\not\in \Psising$, 
 and therefore $(\Atcal \lambda \nu {\delta\varepsilon}{V,W})_{\infty} \ne 0$ by 
 Lemma \ref{lem:ABzero}.  
Moreover,
 if $n-\lambda-\nu-1 \not\in 2{\mathbb{N}}$
 for $\delta \varepsilon =+$
 (or if $n-\lambda-\nu-2 \not\in 2{\mathbb{N}}$
 for $\delta \varepsilon =-$), 
 then we deduce ${\operatorname{Supp}}(\Atcal \lambda \nu {\delta\varepsilon}{V,W})_{\infty}={\mathbb{R}}^n$ from Lemma \ref{lem:Riesz}
 about the support of the Riesz distribution.  
Hence Proposition \ref{prop:172425} is proved.   
\end{proof}
\begin{definition}
[normalized regular symmetry breaking operator]
\label{def:AregSBO}
We shall say $\Atbb \lambda \nu {\delta\varepsilon}{V,W} \colon 
I_{\delta}(V,\lambda) \to J_{\varepsilon}(W,\nu)$ is a holomorphic family
 of the normalized
\index{B}{genericallyregularsymmetrybreakingoperator@generically regular symmetry breaking operator|textbf}
 {\it{(generically) regular symmetry breaking operators}}.  
For simplicity,
 we also call it a holomorphic family
 of the normalized 
\index{B}{regularsymmetrybreakingoperator@regular symmetry breaking operator|textbf}
 regular symmetry breaking operators
 by a little abuse of terminology.  
We are ready to complete the proof of Theorem \ref{thm:regexist}.  
\end{definition}
\begin{proof}[Proof of Theorem \ref{thm:regexist}]
The implication (i) $\Rightarrow$ (iii) follows from the explicit construction
 of the (normalized) regular symmetry breaking operators
 $\Atbb \lambda \nu {\pm} {V,W}$
 in Theorem \ref{thm:152389}, 
 and from Proposition \ref{prop:172425}.

(iii) $\Rightarrow$ (ii) Clear.

Let us prove the implication (ii) $\Rightarrow$ (i).  
We use the notation as in Section \ref{subsec:Xi}
 which is adopted from \cite[Chap.~5] {sbon}.  
Then there exists a unique open orbit of $P'$ on $G/P$, 
 and the isotropy subgroup
 at 
\index{A}{qplusvector@$q_+ ={}^{t}(0,\cdots,0,1,1)$}
$[q_+] =[{}^{t}\!(0,\cdots,0,1,1)]\in \Xi/{\mathbb{R}}^{\times} \simeq G/P$
 is given by 
\[
   \{\begin{pmatrix} 1 & & & \\ & B & & \\ & & 1 & \\ & & & 1\end{pmatrix}:
 B \in O(n-1)\}
 \simeq O(n-1).  
\]
Then the implication (ii) $\Rightarrow$ (i) follows from
 the necessary condition
 for the existence of regular symmetry breaking operators proved in 
 \cite[Prop.~3.5]{sbon}.

Thus Theorem \ref{thm:regexist} is proved.  
\end{proof}

\subsection{Zeros of $\Atbb \lambda \nu {\pm} {V,W}$ : Proof of Theorem \ref{thm:1532113}}
\label{subsec:1532113}
This section discusses the zeros
 of the analytic continuation of the symmetry breaking operator
$
   \Atbb \lambda \nu {\gamma} {V,W}
   \colon
   I_{\delta}(V,\lambda) \to J_{\varepsilon}(W,\nu)  
$
 with $\delta \varepsilon =\gamma$.  
\begin{proof}
[Proof of Theorem \ref{thm:1532113}]
(1)\enspace
Let 
\index{A}{Nsigma@$N(\sigma)$}
$
N:=N(\sigma)
$
 as in Definition \ref{def:Nsigma}.  
We first observe that 
\begin{align*}
(\lambda-N, \nu+N) \in L_{\operatorname{even}}\quad
&
\text{ if $(\lambda, \nu) \in L_{\operatorname{even}}$
 and $\nu \le -N$,}
\\
(\lambda-N-1, \nu+N) \in L_{\operatorname{even}}\quad
&
\text{ if $(\lambda, \nu) \in L_{\operatorname{odd}}$
 and $\nu \le -N$.}
\end{align*}
Then the scalar-valued distributions
 $\Atcal {\lambda-N}{\nu+N}{+}{}$
 and $\Atcal {\lambda-N-1}{\nu+N}{+}{}$
 vanish, 
 respectively
 by \cite[Thm.~8.1]{sbon}.  
By Lemma \ref{lem:152392}, 
 the $\operatorname{Hom}_{\mathbb{C}}(V,W)$-valued distributions
 $p_{+,N}(\lambda,\nu)\Atcal \lambda \nu {+} {V,W}$
 and $p_{-,N}(\lambda,\nu)\Atcal \lambda \nu {-} {V,W}$ vanish, 
respectively,
 because the multiplication of distributions
 by the polynomial $g^{V,W}(x,x_n)$
 is well-defined.  
Since $p_{+,N}(\lambda,\nu) \ne 0$
 for 
\index{A}{Leven@$L_{\operatorname{even}}$}
$(\lambda,\nu) \in L_{\operatorname{even}}$
 and $p_{-,N}(\lambda,\nu) \ne 0$
 for 
\index{A}{Lodd@$L_{\operatorname{odd}}$}
$(\lambda,\nu) \in L_{\operatorname{odd}}$, 
 the first assertion follows from 
 Proposition \ref{prop:Tpair} (2).  
\par\noindent
(2)\enspace
If the symmetry breaking operator $\Atbb \lambda \nu {\gamma} {V,W}$ vanishes, 
 then its distribution kernel is zero, 
 and in particular,
 $(\Atcal \lambda \nu {\gamma} {V,W})_{\infty}=0$
 (see Proposition \ref{prop:Tpair}).  
This implies $\nu-\lambda \in 2{\mathbb{N}}$
 for $\gamma=+$, 
 and $\nu-\lambda \in 2 {\mathbb{N}}+1$
 for $\gamma=-$, 
 owing to Lemma \ref{lem:ABzero}. 
Hence Theorem \ref{thm:1532113} is proved.   
\end{proof}

\subsection{Generic multiplicity-one theorem:
 Proof of Theorem \ref{thm:unique}}
\label{subsec:generic}
\index{B}{genericmultiplicityonetheorem@generic multiplicity-one theorem}
We recall from \eqref{eqn:Psisp}
 the definition of \lq\lq{generic parameter}\rq\rq\
 \eqref{eqn:nlgen}
 that $(\lambda, \nu, \delta, \varepsilon) \not \in
 \Psising$
if and only if
\begin{equation*}
\text{
  $\nu - \lambda \not \in 2{\mathbb{N}}$
 for $\delta \varepsilon = +$;
\quad
  $\nu - \lambda \not \in 2{\mathbb{N}}+1$
 for $\delta \varepsilon = -$.
}
\end{equation*}
We are ready to classify symmetry breaking operators
 for generic parameters.  
The main result
 of this section 
 is Theorem \ref{thm:genbasis}, from which Theorem \ref{thm:unique}
 follows.  
\begin{theorem}
[generic multiplicity-one theorem]
\label{thm:genbasis}
Suppose $(\sigma,V) \in \widehat{O(n)}$, 
 $(\tau,W) \in \widehat{O(n-1)}$
 with 
\index{A}{VWmult@$[V:W]$}
$
[V:W]\ne 0.  
$
Assume $(\lambda, \nu)\in {\mathbb{C}}^2$
 and $\delta, \varepsilon \in \{\pm \}$
 satisfy the 
\index{B}{genericparametercondition@generic parameter condition}
 generic parameter condition, 
 namely, 
 $(\lambda, \nu, \delta, \varepsilon) \not \in
 \Psising$.  
Then the normalized operator $\Atbb \lambda \nu {\delta \varepsilon} {V,W}$ is nonzero
 and is not a differential operator.  
Furthermore we have
\[
    {\operatorname{Hom}}_{G'}
    (I_{\delta}(V,\lambda)|_{G'}, J_{\varepsilon}(W,\nu))
    =
   {\mathbb{C}} \Atbb \lambda \nu {\delta \varepsilon} {V,W}.  
\]
\end{theorem}

\begin{proof}
By Theorem \ref{thm:152389}, 
 $\Atbb \lambda \nu {\pm}{V,W}$ is a symmetry breaking operator
 for all $\lambda, \nu \in {\mathbb{C}}$.  
The generic assumption on $(\lambda, \nu, \delta, \varepsilon)$
 implies $\Atbb \lambda \nu {\delta \varepsilon} {V,W} \ne 0$
 by Theorem \ref{thm:1532113} (2).  
On the other hand,
 by Theorem \ref{thm:SDone} and Corollary \ref{cor:160150upper},
 we see that $\Atbb \lambda \nu {\delta \varepsilon} {V,W}$ is not 
 a differential operator
 and $\dim_{{\mathbb{C}}}{\operatorname{Hom}}_{G'}
    (I_{\delta}(V,\lambda)|_{G'}, J_{\varepsilon}(W,\nu)) \le 1$.  
Thus we have proved Theorem \ref{thm:genbasis}.  
\end{proof}

The generic multiplicity-one theorem 
 given in Theorem \ref{thm:unique} is the second statement
 of Theorem \ref{thm:genbasis}.

\subsection{Lower estimate of the multiplicities}
\label{subsec:existSBO}

In this section
 we do not assume the generic parameter condition
 (Definition \ref{def:generic}), 
 and allow the case $(\lambda, \nu, \delta,\varepsilon) \in \Psising$.  
In this generality, 
 we give a lower estimate 
 of the dimension of the space of symmetry breaking operators.  
\begin{theorem}
\label{thm:existSBO}
Let $(\sigma, V) \in \widehat {O(n)}$
 and $(\tau,W) \in \widehat {O(n-1)}$
 satisfying 
$
[V:W]\ne 0.  
$
For any $\delta, \varepsilon \in \{ \pm \}$ 
 and $(\lambda, \nu)\in {\mathbb{C}}^2$, 
 we have
\[
   \dim_{\mathbb{C}}
   {\operatorname{Hom}}_{G'}(I_{\delta}(V,\lambda)|_{G'},J_{\varepsilon}(W,\nu))\ge 1.  
\]
\end{theorem}

We use a general technique from \cite[Lem.~11.10]{sbon}
 to prove that the multiplicity function
 is upper semicontinuous.

As before, 
 we denote by 
$
   ((\Atcal \lambda \nu {\gamma} {V,W})_{\infty}, 
    \Atcal \lambda \nu {\gamma} {V,W})
$ 
 the pair of ${\operatorname{Hom}}_{\mathbb{C}}(V,W)$-valued distributions
 on ${\mathbb{R}}^n$
that represents the symmetry breaking operator $\Atbb \lambda \nu {\gamma} {V,W}$
 via Proposition \ref{prop:Tpair}.

We fix $(\lambda_0, \nu_0)\in {\mathbb{C}}^2$,
 and define ${\operatorname{Hom}}_{\mathbb{C}}(V,W)$-valued
 distributions on ${\mathbb{R}}^n$
 for $k, \ell \in {\mathbb{N}}$ as follows:
\begin{align*}
F_{k \ell}:=& \left.\frac{\partial^{k+\ell}}{\partial \lambda^k \partial \nu^\ell}\right|_{\substack{\lambda=\lambda_0 \\ \nu=\nu_0}}
          \Atcal \lambda \nu {\gamma}{V,W}, 
\\
(F_{k \ell})_{\infty}:=& \left.\frac{\partial^{k+\ell}}{\partial \lambda^k \partial \nu^\ell}\right|_{\substack{\lambda=\lambda_0 \\ \nu=\nu_0}}
          (\Atcal \lambda \nu {\gamma}{V,W})_{\infty}.   
\end{align*}

\begin{lemma}
\label{lem:upsemi}
Let $\gamma \in \{ \pm \}$ and $m$ a positive integer
 such that 
\[
   ((F_{k \ell})_{\infty}, F_{k \ell})=(0,0)
\quad\text{for all $(k, \ell)\in {\mathbb{N}}^2$ with $k+\ell<m$.}
\]  
Then for any $(k,\ell)$ with $k+\ell=m$, 
 the pair $((F_{k\ell})_{\infty}, F_{k\ell})$ defines a symmetry breaking operator
 $I_{\delta}(V,\lambda) \to J_{\varepsilon}(W,\nu)$
 for $(\delta, \varepsilon)$ with $\delta \varepsilon=\gamma$.  
\end{lemma}

\begin{proof}
Since both the equations \eqref{eqn:Tinfty}--\eqref{eqn:Tpatching}
 and the pairs
$
   ((\Atcal \lambda \nu {\gamma}{V,W})_{\infty}, 
    \Atcal \lambda \nu {\gamma}{V,W})
$
 satisfying \eqref{eqn:Tinfty}--\eqref{eqn:Tpatching}
 depend holomorphically on $(\lambda,\nu)$
 in the entire ${\mathbb{C}}^2$, 
 we can apply \cite[Lem.~11.10]{sbon} 
 to conclude that the pair $((F_{k \ell})_{\infty},F_{k \ell})$ satisfies
 \eqref{eqn:Tinfty}--\eqref{eqn:Tpatching}
 at $(\lambda, \nu)=(\lambda_0, \nu_0)$ 
 for any $(k, \ell)\in {\mathbb{N}}^2$ with $k+\ell=m$.  
Then $((F_{k \ell})_{\infty},F_{k \ell})$ gives an element
 in 
$
   {\operatorname{Hom}}_{G'}
   (I_{\delta}(V,\lambda_0)|_{G'},J_{\varepsilon}(W,\nu_0))
$
by Proposition \ref{prop:Tpair}.  
\end{proof}

\begin{definition}
\label{def:AVWder}
Suppose we are in the setting of Lemma \ref{lem:upsemi}.  
For $(k,\ell)$ with $k+\ell=m$
 and $\delta, \varepsilon \in \{ \pm \}$ with $\delta \varepsilon = \gamma$, 
 we denote by 
\index{A}{Ahtgln2p@$\frac{\partial^{k+l}}{\partial \lambda^k \partial \nu^l}
\mid_{\substack {\lambda=\lambda_0 \\ \nu = \nu_0}} \Atbb{\lambda}{\nu}{\gamma}{V,W}$|textbf}
\[
 \left. \frac{\partial^{k+\ell}}{\partial \lambda^k \partial \nu^\ell}
 \right|_{\substack {\lambda=\lambda_0 \\ \nu = \nu_0}} 
 \Atbb \lambda \nu \gamma {V, W}
 \in 
 {\operatorname{Hom}}_{G'}(I_{\delta}(V,\lambda_0)|_{G'}, J_{\varepsilon}(W,\nu_0)), 
\]
the symmetry breaking operator
 associated to the pair $((F_{k \ell})_{\infty}, F_{k \ell})$.  
\end{definition}

\begin{proof}
[Proof of Theorem \ref{thm:existSBO}]
Set $\gamma:=\delta \varepsilon$.  
Then the pair 
$
   ((\Atcal \lambda \nu {\gamma} {V,W})_{\infty}, 
    \Atcal \lambda \nu {\gamma} {V,W})
$
 of ${\operatorname{Hom}}_{\mathbb{C}}(V,W)$-valued distributions
 depends holomorphically 
 on $(\lambda, \nu)$ in the entire ${\mathbb{C}}^2$
 and satisfies \eqref{eqn:Tinfty}--\eqref{eqn:Tpatching}
 for all $(\lambda, \nu)\in {\mathbb{C}}^2$.  
Moreover, 
the pair 
$
   ((\Atcal \lambda \nu {\gamma} {V,W})_{\infty}, 
    \Atcal \lambda \nu {\gamma} {V,W})
$
 is nonzero
 as far as $\nu-\lambda \not \in {\mathbb{N}}$
 by Lemma \ref{lem:ABzero}.  
This implies that, 
 given $(\lambda_0,\nu_0) \in {\mathbb{C}}^2$, 
 there exists $(k,\ell)\in {\mathbb{N}}^2$
 for which $((F_{k \ell})_{\infty},F_{k \ell})$ is nonzero. 
Take $(k, \ell) \in {\mathbb{N}}^2$
 such that $k+\ell$ attains the minimum
 among all $(k,\ell)$
 for which the pair $((F_{k \ell})_{\infty},F_{k \ell})$ is nonzero.  
By Lemma \ref{lem:upsemi}, 
$
 \left.\frac{\partial^{k+\ell}}{\partial \lambda^k \partial \nu^\ell}\right|_{\substack {\lambda=\lambda_0 \\ \nu = \nu_0}} 
 \Atbb {\lambda} {\nu} \gamma {V, W}
$
 is a symmetry breaking operator.  
\end{proof}

\subsection{Renormalization of symmetry breaking operators
 $\Atbb \lambda \nu \gamma {V,W}$}
\label{subsec:exAVW}

In this section we construct a nonzero symmetry breaking operator
 $\Attbb {\lambda_0} {\nu_0} \gamma {V, W}$
 by \lq\lq{renormalization}\rq\rq\
 when $\Atbb {\lambda_0} {\nu_0} \gamma {V, W}=0$.  
We shall also prove
 that the renormalized operator is {\it{not}} a differential operator.  
The main results are stated
 in Theorem \ref{thm:170340}.  

\subsubsection{Expansion of $\Atbb {\lambda} {\nu} \gamma {V, W}$
 along $\nu=\text{constant}$}
\label{subsec:expand}
We fix $\gamma \in \{\pm\}$ and $(\lambda_0, \nu_0)\in {\mathbb{C}}^2$
 such that
\[
  \nu_0 - \lambda_0
  =
  \begin{cases}
  2\ell \qquad &\text{for $\gamma=+$, }
\\
  2\ell+1 \qquad &\text{for $\gamma=-$, }
  \end{cases}
\]
with $\ell \in {\mathbb{N}}$.  
For every $(\sigma, V) \in \widehat{O(n)}$ and $(\tau, W) \in \widehat{O(n-1)}$, the distribution kernel $\Atcal \lambda \nu \gamma {V,W}$
 of the symmetry breaking operator
 $\Atbb \lambda \nu \gamma {V,W}$ is
 a ${\operatorname{Hom}}_{\mathbb{C}}(V,W)$-valued
 distribution on ${\mathbb{R}}^n$
 that depend holomorphically on $(\lambda,\nu) \in {\mathbb{C}}^2$
 by Theorem \ref{thm:152389}.  
We fix $\nu=\nu_0$ and expand $\Atcal \lambda {\nu_0} \gamma {V,W}$
 with respect to $\lambda$ near $\lambda=\lambda_0$
 as 
\begin{equation}
\label{eqn:Anear0}
   \Atcal \lambda {\nu_0} \gamma {V,W}
   =
   F_0 + (\lambda-\lambda_0) F_1 + (\lambda-\lambda_0)^2 F_2 + \cdots
\end{equation}
 with ${\operatorname{Hom}}_{\mathbb{C}}(V,W)$-valued
 distributions $F_0$, $F_1$, $F_2, \cdots$
 on ${\mathbb{R}}^n$.  
By definition, 
\[
\text{
$\Atcal {\lambda_0} {\nu_0} \gamma {V,W} \ne 0$
 if and only if 
$F_0 \ne 0$.
}
\]

For the next term $F_1$, 
 we have the following two equivalent expressions:
\begin{equation}
\label{eqn:F1}
F_1
=
\lim_{\lambda \to \lambda_0}\frac 1 {\lambda-\lambda_0}
(\Atcal \lambda {\nu_0} \gamma {V,W} - \Atcal {\lambda_0} {\nu_0} {\gamma} {V,W}), 
\end{equation}
 and 
\begin{equation}
\label{eqn:F2}
F_1
=
\left. \frac{\partial}{\partial \lambda}\right|_{\lambda=\lambda_0}
 \Atcal \lambda {\nu_0} \gamma {V,W}.  
\end{equation}

\subsubsection{
Renormalized regular symmetry breaking operator
 $\Attbb \lambda \nu \gamma {V,W}$}
\label{subsec:AVW}
 \index{B}{renormalized regular symmetry breaking operator@regular symmetry breaking operator, renormalized---|textbf}

We consider the following renormalized operators
\begin{alignat}{2}
\label{eqn:rAVW+}
  \Attbb \lambda \nu + {V,W}:=&\Gamma(\frac{\lambda-\nu}{2}) \Atbb \lambda \nu + {V,W}
  \qquad
  &&\text{for $\nu-\lambda \not \in 2 {\mathbb{N}}$, }
\\
\label{eqn:rAVW-}
  \Attbb \lambda \nu - {V,W}:=&\Gamma(\frac{\lambda-\nu+1}{2}) \Atbb \lambda \nu - {V,W}
  \qquad
  &&\text{for $\nu-\lambda \not \in 2 {\mathbb{N}}+1$.  }
\end{alignat}
Since $\Atbb \lambda \nu {\gamma} {V,W}$ depend holomorphically 
 on $(\lambda, \nu)$ in ${\mathbb{C}}^2$, 
\index{A}{Ahttgln0@$\Attbb {\lambda}{\nu}{\pm}{V,W}$|textbf}
 $\Attbb \lambda \nu \gamma {V,W}$ are obviously well-defined 
 as symmetry breaking operators
 $I_{\delta}(V,\lambda)\to J_{\varepsilon}(W,\nu)$
 if $\gamma = \delta \varepsilon$, 
 because the gamma factors do not have poles
 in the domain of definitions \eqref{eqn:rAVW+} and \eqref{eqn:rAVW-}.

On the other hand, 
 Theorem \ref{thm:1532113} (2) implies 
 that the gamma factors in \eqref{eqn:rAVW+} or \eqref{eqn:rAVW-}
 have poles 
 if $\Atbb {\lambda_0} {\nu_0} \gamma {V,W}=0$.  
Nevertheless we shall see in Theorem \ref{thm:170340} below
 that the renormalization $\Attbb {\lambda_0} {\nu_0} \gamma {V,W}$
 still makes sense
 if $\Atbb {\lambda_0} {\nu_0} \gamma {V,W}=0$.

\begin{theorem}
\label{thm:170340}
Suppose $[V:W]\ne 0$
 and let $({\lambda_0}, {\nu_0}) \in {\mathbb{C}}^2$
 such that $\Atbb {\lambda_0} {\nu_0} {\gamma} {V,W}=0$.  
\begin{enumerate}
\item[{\rm{(1)}}]
There exists $\ell \in {\mathbb{N}}$
 such that
\[
 \nu_0 - \lambda_0 =
 \begin{cases}
 2\ell \quad &\text{when $\gamma =+$, }
\\
 2\ell+1 \quad &\text{when $\gamma =-$. }
 \end{cases}
\]
\item[{\rm{(2)}}]
We set 
\begin{equation}
\label{eqn:Ader}
  \Attbb {\lambda_0} {\nu_0} {\gamma} {V,W}
  :=
  \frac{2(-1)^\ell}{\ell!}
  \left. \frac{\partial}{\partial \lambda}\right|_{\lambda=\lambda_0}
  \Atbb {\lambda} {\nu_0} {\gamma} {V,W}.  
\end{equation}
Then $\Attbb {\lambda_0} {\nu_0} {\gamma} {V,W}$
 gives a nonzero symmetry breaking operator from
 $I_{\delta}(V,\lambda)$ to $J_{\varepsilon}(W,\nu_0)$
 with $\delta \varepsilon = \gamma$.  
\item[{\rm{(3)}}]
We fix $\nu=\nu_0$.  
Then $\Attbb {\lambda} {\nu_0} {\gamma} {V,W}$ defined
 by \eqref{eqn:rAVW+} and \eqref{eqn:rAVW-}
 for $\lambda \ne \lambda_0$, 
 and by \eqref{eqn:Ader} for $\lambda = \lambda_0$, 
 is a family of symmetry breaking operators from 
 $I_{\delta}(V,\lambda)$ to $J_{\varepsilon}(W,\nu_0)$
 with $\delta \varepsilon =\gamma$ 
 that depend holomorphically
 on $\lambda$ 
 in the entire complex plane ${\mathbb{C}}$.  
In particular,
 we have 
\begin{equation}
\label{eqn:rSBOlim}
\Attbb {\lambda_0} {\nu_0} {\gamma} {V,W}
=
  \lim_{\lambda \to \lambda_0} \Attbb {\lambda} {\nu_0} {\gamma} {V,W}.  
\end{equation}
\item[{\rm{(4)}}]
$\Attbb {\lambda_0} {\nu_0} {\gamma} {V,W}$ is not a differential operator.  
\end{enumerate}
\end{theorem}

\begin{proof}
\begin{enumerate}
\item[(1)]
The assertion is already given in Theorem \ref{thm:1532113} (2).  
\item[(2)]
The assertion follows from Lemma \ref{lem:upsemi}.  
\item[(3)]
By the first statement,
 we see $(\lambda, \nu_0,\delta, \varepsilon)$ 
 with $\delta \varepsilon =\gamma$ satisfies the generic parameter condition 
\eqref{eqn:nlgen}
 if and only if $\lambda \ne \lambda_0$ and that
\[
  \Attbb {\lambda} {\nu_0} {\gamma} {V,W}
  =
  \Gamma(\frac {\lambda-\lambda_0}2 -\ell)
  \Atbb {\lambda} {\nu_0} {\gamma} {V,W}
\quad
  \text{if $\lambda \ne \lambda_0$.  }
\]

We expand the distribution 
$
   \Atcal \lambda {\nu_0} \gamma {V,W}
$
 as in \eqref{eqn:Anear0} near $\lambda=\lambda_0$.  
By the assumption 
 that $\Atbb {\lambda_0} {\nu_0} {\gamma} {V,W}=0$,
 it follows from the two expressions \eqref{eqn:F1} and \eqref{eqn:F2}
 of the second term $F_1$
 that  
\begin{align*}
  F_1 &= \lim_{\lambda \to \lambda_0} \frac 1 {\lambda-\lambda_0}
        \Atcal \lambda {\nu_0} \gamma{V,W}
      =\lim_{\lambda \to \lambda_0} \frac 1 {(\lambda-\lambda_0)
                                          \Gamma(\frac {\lambda-\lambda_0} 2-\ell)}
        \Attcal \lambda {\nu_0} \gamma{V,W}, 
\\
 F_1 &= \frac {(-1)^\ell \ell!}{2}
        \Attcal {\lambda_0} {\nu_0} \gamma{V,W}.  
\end{align*}
In light that $\lim_{\mu \to 0} \mu \Gamma (\frac {\mu}2- \ell)
 = \frac {2(-1)^\ell}{\ell !}$, 
 we obtain
\[
   \lim_{\lambda \to \lambda_0}\Attbb {\lambda} {\nu_0} {\gamma} {V,W}
   =
   \Attbb {\lambda_0} {\nu_0} {\gamma} {V,W}.  
\]
Since $\Attbb {\lambda} {\nu_0} {\gamma} {V,W}$ depends holomorphically
 on $\lambda$ in ${\mathbb{C}}\setminus \{\lambda_0\}$, 
 and since it is continuous at $\lambda=\lambda_0$, 
 $\Attbb {\lambda} {\nu_0} {\gamma} {V,W}$ is holomorphic in $\lambda$ 
 in the entire complex plane ${\mathbb{C}}$.  
\item[(4)]
Let $( (\Attcal {\lambda_0} {\nu_0} \gamma{V,W})_{\infty}, \Attcal {\lambda_0} {\nu_0} \gamma{V,W})$ be the pair
 of the distribution kernels 
 for $\Attbb {\lambda_0} {\nu_0} {\gamma} {V,W}$
 via Proposition \ref{prop:Tpair} (1).  
Then as in the above proof, 
 we have
\[
   (\Attcal {\lambda_0} {\nu_0} {\gamma} {V,W})_{\infty}
   =
    \lim_{\lambda \to \lambda_0}
    (\Attcal {\lambda} {\nu_0} {\gamma} {V,W})_{\infty}.  
\]
By Proposition \ref{prop:20150828-1231} (2), 
 the right-hand side is not zero.  
Hence $\Attbb {\lambda_0} {\nu_0} {\gamma} {V,W}$ is not 
 a differential operator
 by Proposition \ref{prop:Tpair} (3).  
\end{enumerate}
\end{proof}

We are ready to complete the proof
 of Theorem \ref{thm:VWSBO} (2-C).  
\begin{corollary}
\label{cor:Azero}
Let $\gamma \in \{\pm\}$.  
Suppose $\Atbb {\lambda} {\nu} {\gamma} {V,W}=0$.  
Then the following holds.  
\begin{equation}
\label{eqn:HAD}
  {\operatorname{Hom}}_{G'}(I_{\delta}(V,\lambda)|_{G'}, J_{\varepsilon}(W,\nu))  =
  {\mathbb{C}} \Attbb \lambda \nu {\delta\varepsilon}{V,W}
  \oplus
  {\operatorname{Diff}}_{G'}(I_{\delta}(V,\lambda)|_{G'}, 
 J_{\varepsilon}(W,\nu)).  
\end{equation}
\end{corollary}

\begin{proof}
[Proof of Corollary \ref{cor:Azero}]
By Theorem \ref{thm:170340}, 
 the renormalized operator $\Attbb \lambda \nu{\delta\varepsilon}{V,W}$
 is well-defined and nonzero.  
Moreover,
 the right-hand side of \eqref{eqn:HAD}
 is a direct sum, 
 and is contained in the left-hand side.

Conversely,
 take any ${\mathbb{T}} \in  {\operatorname{Hom}}_{G'}(I_{\delta}(V,\lambda)|_{G'}, J_{\varepsilon}(W,\nu))$, 
 and write $({\mathcal{T}}_{\infty}, {\mathcal{T}})$
 for the corresponding pair
 of distribution kernels for ${\mathbb{T}}$ 
 via Proposition \ref{prop:Tpair}.  
Let $\gamma := \delta\varepsilon$.  
Then Proposition \ref{prop:20150828-1231} tells 
 that ${\mathcal{T}}_{\infty}$ must be proportional
 to $(\Attcal {\lambda} {\nu} {\gamma} {V,W})_{\infty}$, 
 namely, 
 ${\mathcal{T}}_{\infty}=C(\Attcal {\lambda} {\nu} {\gamma} {V,W})_{\infty}$
 for some $C \in {\mathbb{C}}$.  
This implies 
 that the distribution kernel 
 ${\mathcal{T}} - C \Attcal {\lambda} {\nu} {\gamma} {V,W}$
 of the symmetry breaking operator
 ${\mathbb{T}}-C \Attbb {\lambda} {\nu} {\gamma} {V,W}$
 is supported
 at the origin,
 and consequently
 ${\mathbb{T}}-C \Attbb {\lambda} {\nu} {\gamma} {V,W}$
 is a differential operator
 by Proposition \ref{prop:Tpair}.  
\end{proof}

\newpage
\section{Differential symmetry breaking operators}
\label{sec:DSVO}

In this chapter,
 we analyze the space
\[
 {\operatorname{Diff}}_{G'}
  (I_{\delta}(V,\lambda)|_{G'},J_{\varepsilon}(W,\nu))
\]
of 
\index{B}{differentialsymmetrybreakingoperators@differential symmetry breaking operator}
 differential symmetry breaking operators
between principal series representations of $G=O(n+1,1)$ and $G'=O(n,1)$
 for arbitrary $V \in \widehat{O(n)}$ and $W \in \widehat{O(n-1)}$
 with 
\index{A}{VWmult@$[V:W]$}
$[V:W]\ne 0$.

The goal of this chapter is to prove Theorem \ref{thm:existDSBO} below.  
We recall from \eqref{eqn:singset} that 
 the set of \lq\lq{special parameters}\rq\rq\ is denoted by 
\index{A}{1psi@$\Psising$,
          special parameter in ${\mathbb{C}}^2 \times \{\pm\}^2$}
\[
 \Psising
=\{(\lambda,\nu,\delta,\varepsilon) \in {\mathbb{C}}^2 \times \{\pm\}^2:
\text{$\nu-\lambda \in 2 {\mathbb{N}}$
 ($\delta\varepsilon = +$)
 or 
$\nu-\lambda \in 2 {\mathbb{N}}+1$
 ($\delta\varepsilon = -$)}
\}.  
\]
\begin{theorem}
\label{thm:existDSBO}
Let $(G, G')=(O(n+1,1), O(n,1))$.  
Suppose $(\sigma,V) \in \widehat{O(n)}$ and $(\tau,W) \in \widehat{O(n-1)}$
 satisfy $[V:W]\ne 0$.  
\begin{enumerate}
\item[{\rm{(1)}}]
The following two conditions on $\lambda, \nu \in {\mathbb{C}}$
 and $\delta, \varepsilon \in \{\pm\}$
 are equivalent:
\begin{enumerate}
\item[{\rm{(i)}}]
$(\lambda,\nu,\delta,\varepsilon) \in \Psising$.  
\item[{\rm{(ii)}}]
${\operatorname{Diff}}_{G'}
  (I_{\delta}(V,\lambda)|_{G'},J_{\varepsilon}(W,\nu)) \ne \{0\}.
$
\end{enumerate}
\item[{\rm{(2)}}]
If $2 \lambda \not \in {\mathbb{Z}}$ 
then {\rm{(i)}} (or equivalently, (ii) ) implies
\begin{enumerate}
\item[{\rm{(ii)$'$}}]
$
\dim_{\mathbb{C}} {\operatorname{Diff}}_{G'}
  (I_{\delta}(V,\lambda)|_{G'},J_{\varepsilon}(W,\nu)) =1.  
$
\end{enumerate}
\end{enumerate}
\end{theorem}

The implication (ii) $\Rightarrow$ (i) in Theorem \ref{thm:existDSBO}
 holds without the assumption $[V:W] \ne 0$
 as we have seen in Theorem \ref{thm:vanishDiff}.  
Thus the remaining part is to show the opposite implication (i) $\Rightarrow$ (ii)
 and the second statement, 
 which will be carried out
 in Sections \ref{subsec:extDiff} and \ref{subsec:mfDiff}, 
 respectively.

\begin{remark}
In the setting
 where $(V,W)=(\Exterior^i({\mathbb{C}}^n), \Exterior^j({\mathbb{C}}^{n-1}))$, 
 an explicit construction
 and the complete classification
 of the space
$
  {\operatorname{Diff}}_{G'}
  (I_{\delta}(V,\lambda)|_{G'},J_{\varepsilon}(W,\nu))
$
 were carried out in \cite{KKP}
 without the assumption $[V:W] \ne 0$, 
 see Fact \ref{fact:3.9}.   
\end{remark}

\subsection{Differential operators between two manifolds}
\label{subsec:diff}

To give a rigorous definition 
 of 
\index{B}{differentialsymmetrybreakingoperators@differential symmetry breaking operator|textbf}
{\it{differential}} symmetry breaking operators,
 we need the notion of differential operators
 between two manifolds, 
 which we now recall.

For any smooth vector bundle ${\mathcal{V}}$
 over a smooth manifold $X$, 
 there exists the unique (up to isomorphism)
 vector bundle $J^k {\mathcal{V}}$
 over $X$
 (called the $k$-th 
\index{B}{jetprolongation@jet prolongation}
{\it{jet prolongation}} of ${\mathcal{V}}$)
 together with the canonical differential operator
\[
   J^k \colon C^{\infty}(X,{\mathcal{V}}) \to C^{\infty}(X,J^k{\mathcal{V}})
\]
of order $k$.  
We recall that a linear operator
$
   D \colon C^{\infty}(X,{\mathcal{V}}) \to C^{\infty}(X,{\mathcal{V}}')
$
between two smooth vector bundles ${\mathcal{V}}$ and ${\mathcal{V}}'$
 over $X$ is called
 a differential operator
 of order at most $k$, 
 if there is a bundle morphism
$
   Q \colon J^k{\mathcal{V}} \to {\mathcal{V}}'
$
 such that $D=Q_{\ast} \circ J^k$, 
 where $Q_{\ast} \colon C^{\infty}(X,J^k{\mathcal{V}}) \to C^{\infty}(X,{\mathcal{V}}')$ is the induced homomorphism.  
We need a generalization
of this classical definition to the case
 of linear operators
 acting between vector bundles
 over two {\it{different}} smooth manifolds.  
\begin{definition}
[differential operators between two manifolds
 \cite{KOSS, KP1}]
\label{def:diff}
~~~\newline
Suppose that $p \colon Y \to X$
 is a smooth map between two smooth manifolds $Y$ and $X$.  
Let ${\mathcal{V}} \to X$ and ${\mathcal{W}} \to Y$ 
 be  two smooth vector bundles.  
A linear map 
$
   D \colon  C^{\infty}(X,{\mathcal{V}}) \to C^{\infty}(Y,{\mathcal{W}})
$
 is said to be a {\it{differential operator}}
 of order at most $k$
 if there exists a bundle map 
 $Q \colon p^{\ast}(J^k{\mathcal{V}}) \to {\mathcal{W}}$ such that 
\[
   D= Q_{\ast} \circ p^{\ast} \circ J^k.  
\]
\end{definition}
Alternatively, 
one can give the following equivalent definitions
 of differential operators acting between vector bundles
 over two manifolds $Y$ and $X$ with morphism $p$:
\begin{enumerate}
\item[$\bullet$]
based on local properties that generalize Peetre's theorem
 \cite{Pe}
 in the $X=Y$ case
 (\cite[Def.~2.1]{KP1});
\item[$\bullet$]
based on the Schwartz kernel theorem
 (\cite[Lem.~2.3]{KP1});
\item[$\bullet$]
by local expression in coordinates
(\cite[Ex.~2.4]{KP1}).  
\end{enumerate}
Here is a local expression 
 in the case where $p$ is an immersion:
\begin{example}
[{\cite[Ex.~2.4 (2)]{KP1}}]
\label{ex:diffXY}
Suppose that $p \colon Y \hookrightarrow X$
 is an immersion.  
\newline
Choose an atlas
 of local coordinates
 $\{(y_i,z_j)\}$ on $X$
 such that $Y$ is given locally
 by $z_j=0$
 for all $j$.  
Then every differential operator
 $D \colon C^{\infty} (X, {\mathcal{V}}) \to C^{\infty} (Y, {\mathcal{W}})$
 is locally
 of the form
\[
   D=\sum_{\alpha,\beta} g_{\alpha \beta} (y)
     \left. 
     \frac{\partial^{|\alpha|+|\beta|}}{\partial y ^{\alpha}\partial z^{\beta}}
     \right|_{z=0}
\qquad
\text{(finite sum), }
\]
where $g_{\alpha \beta} (y)$ are ${\operatorname{Hom}}(V,W)$-valued smooth functions on $Y$.  
\end{example}

Let $X$ and $Y$ be two smooth manifolds
 acted by $G$ and its subgroup $G'$, 
respectively,
 with a $G'$-equivariant smooth map 
 $p \colon Y \to X$.  
When ${\mathcal{V}} \to X$ is a $G$-equivariant vector bundle
 and ${\mathcal{W}} \to Y$ is a $G'$-equivariant one,
 we denote by 
\[
   {\operatorname{Diff}}_{G'}
   (C^{\infty}(X,{\mathcal{V}})|_{G'}, C^{\infty}(Y,{\mathcal{W}}))   
\]
 the space of differential symmetry breaking operators, 
 namely,
 differential operators
 in the sense of Definition \ref{def:diff}
 that are also $G'$-homomorphisms.  

\subsection{Duality for differential symmetry breaking operators}
\label{subsec:dualVerma}

We review briefly the duality theorem
 between differential symmetry breaking operators
 and morphisms for branching of generalized Verma modules.  
See \cite[Sect.~2]{KP1} for details.  

Let $G$ be a (real) Lie group.  
We denote by 
\index{A}{envelopingalgebra@$U({\mathfrak{g}})$, enveloping algebra|textbf}
$U({\mathfrak{g}})$
 the universal enveloping algebra of the complexified Lie algebra
 ${\mathfrak{g}}_{\mathbb{C}}={\operatorname{Lie}}(G) \otimes_{{\mathbb{R}}}
 {\mathbb{C}}$.  
Analogous notations will be applied to other Lie groups.

Let $H$ be a (possibly disconnected) closed subgroup of $G$.  
Given a finite-dimensional representation $F$ of $H$, 
we set 
\index{A}{indhgV@${\operatorname{ind}}_{\mathfrak{h}}^{\mathfrak{g}}(V)
=U({\mathfrak{g}}) \otimes_{U({\mathfrak{h}})}V$|textbf}
\begin{equation}
\label{eqn:VermaH}
{\operatorname{ind}}_{\mathfrak{h}}^{\mathfrak{g}}(F)
:=U({\mathfrak{g}}) \otimes_{U({\mathfrak{h}})}F.  
\end{equation}

The diagonal $H$-action 
 on the tensor product $U({\mathfrak{g}}) \otimes_{\mathbb{C}} F$
 induces an action of $H$
 on $U({\mathfrak{g}}) \otimes_{U({\mathfrak{h}})} F$, 
 and thus ${\operatorname{ind}}_{\mathfrak{h}}^{\mathfrak{g}}(F)$
 is endowed with a $({\mathfrak{g}},H)$-module structure.

When $X$ and $Y$ are homogeneous spaces $G/H$ and $G'/H'$, 
respectively,
 with $G' \subset G$ and $H' \subset H \cap G'$, 
 we have a natural $G'$-equivariant smooth map
 $G'/H' \to G/H$ induced from 
the inclusion map $G' \hookrightarrow G$.  
In this case, 
 the following duality theorem
 (\cite[Thm.~2.9]{KP1}, 
 see also \cite[Thm.~2.4]{KOSS})
 is a generalization
 of the classical duality 
 in the case
 where $G=G'$ are complex reductive Lie groups
 and $H=H'$ are Borel subgroups:

\begin{fact}
[duality theorem]
\label{fact:dualityDSV}
Let $F$ and $F'$ be finite-dimensional representations
 of $H$ and $H'$, 
respectively,
 and we define equivariant vector bundles ${\mathcal{V}}= G \times_H F$
 and ${\mathcal{W}}= G' \times_{H'} F'$
 over $X$ and $Y$, 
 respectively.  
Then there is a canonical linear isomorphism:
\begin{equation}
\label{eqn:dualDSV}
   {\operatorname{Hom}}_{{\mathfrak{g}}',H'}
   ({\operatorname{ind}}_{\mathfrak{h}'}^{\mathfrak{g}'}
   (F'^{\vee}), 
    {\operatorname{ind}}_{\mathfrak{h}}^{\mathfrak{g}}
    (F^{\vee})|_{{\mathfrak{g}}',H'}
    )
\simeq
   {\operatorname{Diff}}_{G'}
   (C^{\infty}(X,{\mathcal{V}})|_{G'}, C^{\infty}(Y,{\mathcal{W}})).  
\end{equation}
\end{fact}
Applying Fact \ref{fact:dualityDSV}
 to our special setting, 
 we obtain the following:
\begin{proposition}
\label{prop:dualityDSV}
Let $(G, G')=(O(n+1,1), O(n,1))$, 
 $V \in \widehat{O(n)}$, $W \in \widehat{O(n-1)}$,
 $\lambda, \nu\in {\mathbb{C}}$, 
 and $\delta, \varepsilon \in \{\pm\}$.  
Let 
\index{A}{Vln@$V_{\lambda,\delta}=V \otimes \delta \otimes {\mathbb{C}}_{\lambda}$, representation of $P$}
$V_{\lambda,\delta}= V \otimes \delta \otimes {\mathbb{C}}_{\lambda}$
 be the irreducible representation of $P$
 with trivial $N_+$-action as before, 
 and $V_{\lambda,\delta}^{\vee}$ the contragredient representation.  
Similarly,
\index{A}{Wnuepsi@$W_{\nu,\varepsilon}=W \otimes \varepsilon \otimes {\mathbb{C}}_{\nu}$, representation of $P'$}
 $W_{\nu,\varepsilon}^{\vee}$ be the contragredient $P'$-module
 of $W_{\nu,\varepsilon}=W \otimes \varepsilon \otimes {\mathbb{C}}_{\nu}$.  
Then there is a canonical linear isomorphism:
\begin{equation}
\label{eqn:dualVW}
  {\operatorname{Hom}}_{{\mathfrak{g}}',P'}
  ({\operatorname{ind}}_{\mathfrak{p}'}^{\mathfrak{g}'}
   (W_{\nu,\varepsilon}^{\vee}), 
    {\operatorname{ind}}_{\mathfrak{p}}^{\mathfrak{g}}
    (V_{\lambda,\delta}^{\vee})|_{{\mathfrak{g}}',P'}
    )
\simeq
  {\operatorname{Diff}}_{G'}
   (I_{\delta}(V,\lambda)|_{G'},J_{\varepsilon}(W,\nu)).  
\end{equation}
\end{proposition}

\subsection{Parabolic subgroup compatible 
 with a reductive subgroup}
\label{subsec:PGcompati}
In this section we treat the general setting 
 where $G$ is a real reductive Lie group
 and $G'$ is a reductive subgroup,
 and study basic properties
 of differential symmetry breaking operators
 between principal series representation $\Pi$
 of $G$
 and $\pi$ of the subgroup $G'$.  
We shall prove in Theorem \ref{thm:imgDSBO} below
 that the image of any nonzero differential symmetry breaking operator
 is infinite-dimensional 
 if $\Pi$ is induced from a parabolic subgroup $P$
 which is compatible with the subgroup $G'$
 (see Definition \ref{def:compatiP}).

Let us give a basic setup.  
Suppose that $G$ is a real reductive Lie group
 with Lie algebra ${\mathfrak{g}}$.  
Take a hyperbolic element $H$ of ${\mathfrak{g}}$, 
 and we define the direct sum decomposition, 
 referred sometimes to as the Gelfand--Naimark decomposition
 (cf. \cite{GN}):
\[
   {\mathfrak{g}}={\mathfrak{n}}_- + {\mathfrak{l}} + {\mathfrak{n}}_+
\]
 where ${\mathfrak{n}}_-$, ${\mathfrak{l}}$, and ${\mathfrak{n}}_+$
 are the sum of eigenspaces
 of ${\operatorname{ad}}(H)$ with negative, zero and positive eigenvalues, 
respectively.  
We define a parabolic subgroup $P \equiv P(H)$ of $G$ by
\[
  P=L N_+
\qquad
 \text{(Levi decomposition)}, 
\]
 where $L = \{g \in G: {\operatorname{Ad}}(g)H=H\}$
 and $N_+=\exp ({\mathfrak{n}}_+)$.  
The following \lq\lq{compatibility}\rq\rq\ gives a sufficient condition
 for the \lq\lq{discrete decomposability}\rq\rq\
 of the generalized Verma module ${\operatorname{ind}}_{\mathfrak{p}}^{\mathfrak{g}}(V^{\vee})$
 when restricted to the subalgebra ${\mathfrak{g}}'$, 
 which concerns with the left-hand side
 of the duality \eqref{eqn:dualDSV}
 (see \cite[Thm.~4.1]{xktransgp12}):
\begin{definition}
[{\cite{xktransgp12}}]
\label{def:compatiP}
Suppose $G'$ is a reductive subgroup of $G$ with Lie algebra ${\mathfrak{g}}'$. A parabolic subgroup $P$ of $G$ is said to be $G'$-{\it{compatible}}
 if there exists a hyperbolic element $H$ in ${\mathfrak{g}}'$
 such that $P=P(H)$.  
\end{definition}

If $P$ is $G'$-compatible, 
 then $P' :=P \cap G'$
 is a parabolic subgroup of the reductive subgroup $G'$
 with Levi decomposition $P' = L' N_+'$
 where $L' :=L \cap G'$ and $N_+' :=N_+ \cap G'$.  

\begin{theorem}
\label{thm:imgDSBO}
Let $G$ be a real reductive Lie group,
 $P$ a parabolic subgroup 
 which is compatible with a reductive subgroup $G'$, 
 and $P' := P \cap G'$.  
Suppose that ${\mathcal{V}}$ is a $G$-equivariant vector bundle
 of finite rank over the real flag manifold $G/P$, 
and that ${\mathcal{W}}$ is a $G'$-equivariant one over $G'/P'$.  
Then for any nonzero differential operator
 $D \colon C^{\infty}(G/P,{\mathcal{V}}) \to C^{\infty}(G'/P',{\mathcal{W}})$, 
 we have 
\[
  \dim_{\mathbb{C}} {\operatorname{Image}} D= \infty.  
\]
\end{theorem}
As we shall see
 in the proof below,
Theorem \ref{thm:imgDSBO}
 follows from the definition
of differential operators (Definition \ref{def:diff})
 without the assumption
 that $D$ intertwines the $G'$-action.  

\begin{proof}
[Proof of Theorem \ref{thm:imgDSBO}]
We set $Y=G'/P'$ and $X=G/P$.  
Then $Y \subset X$ because $P' = P \cap G'$.  
There exist countably many disjoint open subsets $\{U_j\}$ of $X$
 such that $Y \cap U_j \ne \emptyset$.  
It suffices to show
 that for every $j$ there exists $\varphi_j \in C^{\infty}(X,{\mathcal{V}})$
 such that ${\operatorname{Supp}} (\varphi_j) \subset U_j$
 and $D \varphi_j \ne 0$
 because ${\operatorname{Supp}} (D \varphi_j) \subset U_j \cap Y$
 and because $\{U_j \cap Y\}$ is a set of disjoint open sets of $Y$.  
We fix $j$, 
 and write $U$ simply for $U_j$.  
By shrinking $U$ if necessary,
 we trivialize the bundles ${\mathcal{V}}|_U$ and ${\mathcal{W}}|_{U \cap Y}$.  Then we see from Example \ref{ex:diffXY}
 that $D$ can be written locally as the matrix-valued operators:
\[
   D=\sum_{\alpha,\beta} g_{\alpha \beta} (y)
     \left. 
     \frac{\partial^{|\alpha|+|\beta|}}{\partial y ^{\alpha}\partial z^{\beta}}
     \right|_{z=0}.  
\]
Take a multi-index $\beta$
 such that $g_{\alpha \beta} (0) \ne 0$
 on $U$ for some $\alpha$.  
We fix $\alpha$ such that $|\alpha|=\alpha_1+ \cdots +\alpha_{\dim Y}$
attains its maximum 
 among all multi-indeces $\alpha$ with $g_{\alpha \beta} (y) \not \equiv 0$.  
Take $v$ in the typical fiber ${\mathcal{V}}$ at $(y, z)=(0,0)$
 such that $g_{\alpha \beta}(0) v \ne 0$.  
By using a cut function,
 we can construct easily $\varphi \in C^{\infty}(X, {\mathcal{V}})$ 
 such that ${\operatorname{Supp}}(\varphi) \subset U$
 and that $\varphi(y,z) \equiv y^{\alpha} z^{\beta} v$
 in a neighbourhood of $(y,z)=(0,0)$.  
Then we have 
\[
  D \varphi \ne 0.  
\]
Thus Theorem \ref{thm:imgDSBO} is proved.  
\end{proof}

\subsection{Character identity for branching in the parabolic BGG category}
\label{subsec:GroOp}

We retain the general setting 
 as in Section \ref{subsec:PGcompati}, 
 and discuss the duality theorem 
 in Section \ref{subsec:dualVerma}.  
To study the left-hand side of \eqref{eqn:dualVW}, 
 we use the results \cite{xktransgp12, KOSS} on the restriction
 of parabolic Verma modules
 ${\operatorname{ind}}_{\mathfrak{p}}^{\mathfrak{g}}(F)$
  with respect to a reductive subalgebra ${\mathfrak{g}}'$
 under the assumption
 that ${\mathfrak{p}}$ is compatible with ${\mathfrak{g}}'$.  
For later purpose, 
 we need to formulate the results
 in \cite{xktransgp12, KOSS}
 in a slightly more general form as below, 
 because a parabolic subgroup $P$ of a real reductive Lie group
 is not always connected.

Suppose that $P =L N_+$ is a parabolic subgroup of $G$
 which is compatible with a reductive subgroup $G'$.  
We set ${\mathfrak{n}}_-':= {\mathfrak{n}}_- \cap {\mathfrak{g}}'$.  
Then the $L'$-module structure on the nilradical ${\mathfrak{n}}_-$
 descends to the quotient ${\mathfrak{n}}_-/{\mathfrak{n}}_-'$, 
 and extends to the (complex) symmetric tensor algebra
 $S(({\mathfrak{n}}_-/{\mathfrak{n}}_-') \otimes_{\mathbb{R}} {\mathbb{C}})$.

For an irreducible $L$-module $F$
 and an irreducible $L'$-module $F'$, 
 we set 
\begin{equation}
\label{eqn:nFH}
  n(F,F') :=
  \dim_{\mathbb{C}}{\operatorname{Hom}}_{L'}
  (F',F|_L \otimes S(({\mathfrak{n}}_-/{\mathfrak{n}}_-') \otimes_{\mathbb{R}} {\mathbb{C}})).  
\end{equation}

Then we have the following branching rule
 in the Grothendieck group 
 of the parabolic BGG category
 of $({\mathfrak{g}}',P')$-modules
 ({{\cite[Prop.~5.2]{xktransgp12}},{\cite[Thm.~3.5]{KOSS}}}):

\begin{fact}
[character identity for branching to a reductive subalgebra]
\label{fact:Grobranch}
Suppose that $P=L N_+$ is a $G'$-compatible parabolic subgroup of $G$
 (Definition \ref{def:compatiP}).  
Let $F$ be an irreducible finite-dimensional $L$-module.  
\begin{enumerate}
\item[{\rm{(1)}}]
$n(F,F')<\infty$ for all irreducible finite-dimensional $L'$-modules $F'$.  
\item[{\rm{(2)}}]
We inflate $F$ to a $P$-module 
 by letting $N_+$ act trivially,
 and form a $({\mathfrak{g}},P)$-module
 ${\operatorname{ind}}_{\mathfrak{p}}^{\mathfrak{g}}(F)
  = U({\mathfrak{g}}) \otimes_{U(\mathfrak{p})} F$.  
Then we have the following identity
 in the Grothendieck group
 of the parabolic BGG category
 of $({\mathfrak{g}}',P')$-modules:
\[
{\operatorname{ind}}_{\mathfrak{p}}^{\mathfrak{g}}(F)|_{{\mathfrak{g}}',P'}
\simeq
\bigoplus_{F'}
n(F,F') 
{\operatorname{ind}}_{\mathfrak{p}'}^{\mathfrak{g}'}(F').  
\]
In the right-hand side, 
 $F'$ runs over all irreducible finite-dimensional $P'$-modules, 
 or equivalently, 
 all irreducible finite-dimensional $L'$-modules
 with trivial $N_+'$-actions.  
\end{enumerate}
\end{fact}

\begin{proof}
The argument is parallel to the one
 in \cite[Thm.~3.5]{KOSS}
 for $({\mathfrak{g}}',{\mathfrak{p}}')$-modules,
 which is proved by using \cite[Prop.~5.2]{xktransgp12}.  
\end{proof}
\subsection{Branching laws for generalized Verma modules}
\label{subsec:branchVerma}

In this section we refine the character identity
 (identity in the Grothendieck group)
 in Section \ref{subsec:GroOp}
 to obtain actual branching laws.  
The idea works in the general setting
 (cf.~ \cite[Sect.~3]{KOSS}), 
however,
 we confine ourselves with the pair $(G,G')=(O(n+1,1),O(n,1))$
 for actual computations below.  
In particular,
 under the assumption $2 \lambda \not \in {\mathbb{Z}}$, 
 we give an explicit irreducible decomposition 
 of the $({\mathfrak{g}},P)$-module
$
   {\operatorname{ind}}_{\mathfrak{p}}^{\mathfrak{g}}(V_{\lambda,\delta}^{\vee})
$
 when we regard it 
 as a $({\mathfrak{g}}',P')$-module:

\begin{theorem}
[branching law for generalized Verma modules]
\label{thm:Vermadeco}
Let $V \in \widehat {O(n)}$, 
 $\lambda \in {\mathbb{C}}$, 
 and $\delta \in \{\pm\}$.  
Assume $2 \lambda \not \in {\mathbb{Z}}$.  
Then the $({\mathfrak{g}},P)$-module 
 ${\operatorname{ind}}_{\mathfrak{p}}^{\mathfrak{g}}
  (V_{\lambda,\delta}^{\vee})$ decomposes
 into the multiplicity-free direct sum of irreducible $({\mathfrak{g}}',P')$-modules
 as follows:
\begin{equation}
\label{eqn:Vermadeco}
 {\operatorname{ind}}_{\mathfrak{p}}^{\mathfrak{g}}
  (V_{\lambda,\delta}^{\vee})|_{{\mathfrak{g}}',P'}
  \simeq
  \bigoplus_{a=0}^{\infty}
  \bigoplus_{[V:W]\ne 0}
  {\operatorname{ind}}_{\mathfrak{p}'}^{\mathfrak{g}'}
   ((W_{\lambda+a,(-1)^a\delta})^{\vee}).  
\end{equation}
Here $W$ runs over all irreducible $O(n-1)$-modules
 such that $[V:W] \ne 0$.  
\end{theorem}

\begin{proof}
[Proof of Theorem \ref{thm:Vermadeco}]
The hyperbolic element $H$ defined in \eqref{eqn:Hyp}
 is contained in ${\mathfrak{g}}' = {\mathfrak{o}}(n,1)$, 
 and therefore, 
 the parabolic subgroup $P$
 is compatible with the reductive subgroup $G'=O(n,1)$
 in the sense of Definition \ref{def:compatiP}.  
We then apply Fact \ref{fact:Grobranch}
 to 
\[
  (F,{\mathfrak{n}}_-,{\mathfrak{n}}_-')
  =
  (V_{\lambda,\delta}^{\vee}, \sum_{j=1}^{n}{\mathbb{R}}N_j^-, \sum_{j=1}^{n-1}{\mathbb{R}}N_j^-).  
\]
Since ${\mathfrak{n}}_-/{\mathfrak{n}}_-' \simeq {\mathbb{R}}N_n^-$, 
 the $a$-th symmetric tensor space amounts to 
\[
  S^a(({\mathfrak{n}}_-/{\mathfrak{n}}_-') \otimes_{\mathbb{R}} {\mathbb{C}})
   \simeq
  {\bf{1}} \boxtimes (-1)^a \boxtimes {\mathbb{C}}_{-a}
\]
 as a module of $L' \simeq O(n-1) \times O(1) \times {\mathbb{R}}$.  
Therefore we have an $L'$-isomorphism:
\[
F|_{L'} \otimes S^a(({\mathfrak{n}}_-/{\mathfrak{n}}_-') \otimes_{\mathbb{R}} {\mathbb{C}})
\simeq
\bigoplus_{\substack {W \in \widehat{O(n-1)} \\ [V:W] \ne 0}} W^{\vee} \boxtimes (-1)^a\delta \boxtimes {\mathbb{C}}_{-\lambda-a}, 
\]
where we observe $[V^{\vee}:W^{\vee}] \ne 0$
 if and only if $[V:W] \ne 0$.  
Thus the identity \eqref{eqn:Vermadeco}
 in the level of the Grothendieck group
 of $({\mathfrak{g}}',P')$-modules
 is deduced from Fact \ref{fact:Grobranch}.

In order to prove the identity \eqref{eqn:Vermadeco}
 as $({\mathfrak{g}}',P')$-modules, 
 we use the following two lemmas.  
\end{proof}

\begin{lemma}
\label{lem:indinf}
Assume $2 \lambda \not \in {\mathbb{Z}}$.  
Then any ${\mathfrak{Z}}({\mathfrak{g}}')$-infinitesimal characters
 of the summands in \eqref{eqn:Vermadeco}
 are all distinct.  
\end{lemma}

\begin{lemma}
\label{lem:indirred}
Assume $2 \lambda \not \in {\mathbb{Z}}$.  
Then any summand 
$
   {\operatorname{ind}}_{\mathfrak{p}'}^{\mathfrak{g}'}
   ((W_{\lambda+a,(-1)^a \delta})^{\vee})
$
 in \eqref{eqn:Vermadeco}
 is irreducible as a $({\mathfrak{g}}',P')$-module.  
\end{lemma}

\begin{proof}
[Proof of Lemma \ref{lem:indinf}]
Via the Cartan--Weyl bijection \eqref{eqn:CWOn}
 for the disconnected group $O(N)$
 ($N=n,n-1$), 
 we write $V=\Kirredrep{O(n)}{\mu}$
 and $W=\Kirredrep{O(n-1)}{\mu'}$
 for $\mu=(\mu_1, \cdots,\mu_n) \in \Lambda^{+}(O(n))$
 and $\mu'=(\mu_1', \cdots,\mu_{n-1}') \in \Lambda^{+}(O(n-1))$.  
By the classical 
\index{B}{branching rule@branching rule, for $O(N) \downarrow O(N-1)$}
branching law
 for the restriction $O(n) \downarrow O(n-1)$
 (Fact \ref{fact:ONbranch}), 
 $[V:W] \ne 0$
 if and only if
\begin{equation}
\label{eqn:mubranch}
  \mu_1 \ge \mu_1' \ge \mu_2 \ge \cdots \ge \mu_{n-1}' \ge \mu_{n}.   
\end{equation}
Since any irreducible $O(N)$-module
 is self-dual, 
 we have $W^{\vee} \simeq \Kirredrep{O(n-1)}{\mu'}$.  
Therefore, 
 the 
\index{B}{infinitesimalcharacter@infinitesimal character}
${\mathfrak{Z}}({\mathfrak{g}}')$-infinitesimal character
 of the ${\mathfrak{g}}'$-module 
${\operatorname{ind}}_{\mathfrak{p}'}^{\mathfrak{g}'}
 (W^{\vee} \otimes (-1)^a\delta \otimes {\mathbb{C}}_{-\lambda-a})
$ 
 is given by
\[
  (-\lambda-a+\frac{n-1}{2}, 
   \mu_1'+\frac{n-3}{2},
   \mu_2'+ \frac{n-5}{2}, 
   \cdots, 
   \mu_{[\frac{n-1}{2}]}'+\frac{n-1}{2}- [\frac{n-1}{2}])
\]
 modulo the Weyl group
 ${\mathfrak{S}}_m \ltimes ({\mathbb{Z}}/2 {\mathbb{Z}})^m$
 for the disconnected group $G'=O(n,1)$
 where $m=[\frac{n+1}{2}]$.  
Hence,
 if $2 \lambda \not \in {\mathbb{Z}}$, 
they are all distinct
 when $a$ runs over ${\mathbb{N}}$
 and $\mu'$ runs over $\Lambda^+(O(n-1))$
 subject to \eqref{eqn:mubranch}.  
Thus Lemma \ref{lem:indinf} is proved.  
\end{proof}
\begin{proof}
[Proof of Lemma \ref{lem:indirred}]
By the criterion of Conze-Berline and Duflo
 \cite{BeDu}, 
 the ${\mathfrak{g}}'$-module ${\operatorname{ind}}_{\mathfrak{p}'}^{\mathfrak{g}'}
 (\tau_{\nu} \otimes {\mathbb{C}}_{-\lambda-a})
$ 
 is irreducible
 if $\tau_{\nu}$ is an irreducible ${\mathfrak{so}}(n-1)$-module
 with highest weight $(\nu_1, \cdots,\nu_{[\frac{n-1}{2}]})$
 satisfying
\[
\langle
  -\lambda-a+\frac{n-1}{2}, 
   \nu_1+\frac{n-3}{2},
   \nu_2+ \frac{n-5}{2}, 
   \cdots, 
   \nu_{[\frac{n-1}{2}]}+\frac{n-1}{2}- [\frac{n-1}{2}], 
   \beta^{\vee}
\rangle
  \not \in{\mathbb{N}}_+, 
\]
 where $\beta^{\vee}$ is the coroot of $\beta$, 
 and $\beta$ runs over the set 
\[
\Delta^+({\mathfrak{g}}_{\mathbb{C}}) 
 \setminus \Delta^+({\mathfrak{l}}_{\mathbb{C}})
 =\{e_1 \pm e_j : 2 \le j \le [\frac{n+1}{2}]\}
 (\cup \{e_1\}, 
 \text{when $n$ is even}).  
\]
This condition is fulfilled
 if $2 \lambda \not \in {\mathbb{Z}}$
 because $\nu_1$, $\cdots$, $\nu_{[\frac{n-1}{2}]}$
 $\in \frac 1 2 {\mathbb{Z}}$
 and $a \in {\mathbb{N}}$.  
Hence ${\operatorname{ind}}_{\mathfrak{p}'}^{\mathfrak{g}'}
 ((W_{\lambda+a,(-1)^a\delta})^{\vee})
$
 is an irreducible ${\mathfrak{g}}'$-module
 if $W^{\vee}(\simeq W) \in O(n-1)$
 is of type X
 (Definition \ref{def:OSO}), 
 namely,
 if $W^{\vee}$ is irreducible 
 as an ${\mathfrak{so}}(n-1)$-module.  
On the other hand, 
if $W^{\vee} \in O(n-1)$ is of type Y, 
 then ${\operatorname{ind}}_{\mathfrak{p}'}^{\mathfrak{g}'}
 ((W_{\lambda+a,(-1)^a\delta})^{\vee})
$ splits into the direct sum
 of two irreducible ${\mathfrak{g}}'$-module
 according to the decomposition 
 of $W^{\vee}$
 into irreducible ${\mathfrak{so}}(n-1)$-modules.  
Since these two ${\mathfrak{g}}'$-submodules 
 are not stable
 by the $L'$-action,
 we conclude
 that ${\operatorname{ind}}_{\mathfrak{p}'}^{\mathfrak{g}'}
 ((W_{\lambda+a,(-1)^a\delta})^{\vee})
$
 is irreducible as a $({\mathfrak{g}'},L')$-module,
 in particular,
 as a $({\mathfrak{g}'},P')$-module.  
Thus Lemma \ref{lem:indirred} is proved.  
\end{proof}
\subsection{Multiplicity-one theorem for differential symmetry breaking operators: Proof of Theorem \ref{thm:existDSBO} (2)}
\label{subsec:mfDiff}
Combining Proposition \ref{prop:dualityDSV} (duality theorem)
 with the branching law
 for generalized Verma modules
 (Theorem \ref{thm:Vermadeco}), 
 we obtain a 
\index{B}{genericmultiplicityonetheoremforDSBO@
generic multiplicity-one theorem,
 for differential symmetry breaking operator}
generic multiplicity-one theorem
 for differential symmetry breaking operators
 as follows:
\begin{corollary}
\label{cor:DSVOgeneric}
Suppose $V \in \widehat {O(n)}$ and $W \in \widehat{O(n-1)}$
 satisfy $[V:W] \ne 0$.  
Suppose that $(\lambda, \nu, \delta, \varepsilon) \in \Psising$
 (see \eqref{eqn:singset}).  
Assume further $2 \lambda \not \in {\mathbb{Z}}$.  
Then 
\[
   \dim_{\mathbb{C}} 
   {\operatorname{Diff}}_{G'}
   (I_{\delta}(V,\lambda)|_{G'},J_{\varepsilon}(W,\nu))
   =1. 
\]
\end{corollary}
This gives a proof of the second statement
 of Theorem \ref{thm:existDSBO}.  

\subsection{Existence of differential symmetry breaking operators:
Extension to special parameters}
\label{subsec:extDiff}
What remains to prove is the implication (i) $\Rightarrow$ (ii)
 in Theorem \ref{thm:existDSBO}
 for special parameters,
 namely,
 for $2 \lambda \in {\mathbb{Z}}$.  
We shall use the general idea given in \cite[Lem.~11.10]{sbon}
 and deduce the implication (i) $\Rightarrow$ (ii)
 for the special parameters from Corollary \ref{cor:DSVOgeneric}
 for the regular parameters,  
 and thus complete the proof of Theorem \ref{thm:existDSBO} (1).

Let ${\operatorname{Diff}}^{\operatorname{const}}({\mathfrak{n}}_{-})$
 denote the ring of holomorphic differential operators
 on ${\mathfrak{n}}_{-}$
 with constant coefficients
 and $\langle \,\, , \,\, \rangle$ denote the natural pairing 
 ${\mathfrak{n}}_-=\sum_{j=1}^{n}{\mathbb{R}} N_j^-$
 and ${\mathfrak{n}}_+=\sum_{j=1}^{n}{\mathbb{R}} N_j^+$.  
Then the symbol map
\[
  {\operatorname{Symb}} \colon 
  {\operatorname{Diff}}^{\operatorname{const}}({\mathfrak{n}}_{-})
  \to 
  {\operatorname{Pol}}({\mathfrak{n}}_{+}), 
  \quad
  D_z \mapsto Q(\zeta)
\]
 given by the characterization
\[
   D_z e^{\langle z, \zeta \rangle}
   =
   Q(\zeta) e^{\langle z, \zeta \rangle}
\]
 is a ring isomorphism 
 between ${\operatorname{Diff}}^{\operatorname{const}}({\mathfrak{n}}_{-})$
 and the polynomial ring ${\operatorname{Pol}}({\mathfrak{n}}_{+})$.

The 
\index{B}{Fmethod@F-method}
F-method (\cite[Thm.~4.1]{KP1})
 characterizes the \lq\lq{Fourier transform}\rq\rq\
 of differential symmetry breaking operators
 by certain systems
 of differential equations.  
It tells that 
 any element in ${\operatorname{Diff}}_{G'}
   (I_{\delta}(V,\lambda)|_{G'},J_{\varepsilon}(W,\lambda+a))$
 is given as a ${\operatorname{Hom}}_{\mathbb{C}}(V,W)$-valued
 differential operator $D$
 on the Bruhat cell $N_- \simeq {\mathbb{R}}^n$
 as 
\[
  D={\operatorname{Rest}}_{x_n=0} \circ ({\operatorname{Symb}}^{-1} \otimes {\operatorname{id}})(\psi), 
\]
 where $\psi(\zeta_1, \cdots, \zeta_n)$
 is a ${\operatorname{Hom}}_{\mathbb{C}}(V,W)$-valued homogeneous polynomial
 of degree $a$
 satisfying a system of linear (differential) equations
 (cf.~\cite[(4.3) and (4.4)]{KP1})
 that depend holomorphically on $\lambda \in {\mathbb{C}}$.

If we write the solution $\psi(\zeta)$ as 
\[
  \psi(\zeta)=\sum_{\beta_1 + \cdots + \beta_n=a} 
              \varphi(\beta) \zeta_1^{\beta_1} \cdots \zeta_n^{\beta_n}, 
\]
then the system of differential equations
 for $\psi(\zeta)$ in the F-method
 amounts to a system 
 of linear (homogeneous) equations 
 for the coefficients $\{\varphi(\beta):|\beta|=a\}$.  
We regard $\varphi = (\varphi(\beta)) \in {\mathbb{C}}^k$
 where $k := \#\{\beta \in {\mathbb{N}}^n: |\beta|=\alpha\}$, 
 and use the following elementary lemma
 on the global basis of solutions:
\begin{lemma}
\label{lem:holoeq}
Let $Q_{\lambda}\varphi =0$ be a system
 of linear homogeneous equations
 of $\varphi \in {\mathbb{C}}^k$
 such that $Q_{\lambda}$ depends holomorphically on $\lambda \in {\mathbb{C}}$.  
Assume that there exists a nonempty open subset $U$ of ${\mathbb{C}}$
 such that the space of solutions to $Q_{\lambda} \varphi =0$
 is one-dimensional 
 for every $\lambda$ 
 in $U$.  
Then there exists $\varphi_{\lambda} \in {\mathbb{C}}^k$
 that depend holomorphically on $\lambda$
 in the entire ${\mathbb{C}}$
 such that $Q_{\lambda}\varphi_{\lambda} =0$
 for all $\lambda \in {\mathbb{C}}$.   
\end{lemma}

\begin{proof}
We may regard the equation $Q_{\lambda} \varphi =0$
 as a matrix equation 
 where $Q_{\lambda}$ is an $l$ by $k$ matrix
 ($l \ge k$) 
 whose entries are holomorphic functions of $\lambda \in {\mathbb{C}}$.  
By assumption,
 we have
\[
 {\operatorname{rank}} \, Q_{\lambda}=k-1
\qquad
 \text{for all }\,\, \lambda \in U.  
\]
We can choose a nonempty open subset $U'$ of $U$
 and $k$ row vectors in $Q_{\lambda}$
 such that the corresponding square submatrix $P_{\lambda}$
 is of rank $k-1$, 
 provided $\lambda$ belongs to $U'$.  
Then at least one of row vectors
 in the cofactor of $P_{\lambda}$ is nonzero, 
 which we choose and denote by $\varphi_{\lambda}$.  
Clearly, 
 $\varphi_{\lambda}$ depends holomorphically
 on the entire ${\lambda}\in {\mathbb{C}}$, 
 and $Q_{\lambda} \varphi_{\lambda}=0$
 for all $\lambda \in U'$.

Since both $Q_{\lambda}$ and $\varphi_{\lambda}$ depend holomorphically
 on $\lambda$
 in the entire ${\mathbb{C}}$, 
 the equation $Q_{\lambda}\varphi_{\lambda}=0$ holds
 for all $\lambda \in {\mathbb{C}}$.  
\end{proof}

We note that the solution $\varphi_{\lambda}$
 in Lemma \ref{lem:holoeq} may vanish 
 for some $\lambda \in {\mathbb{C}}$.  
However,
 the following nonvanishing result holds for all $\lambda \in {\mathbb{C}}$.  
\begin{proposition}
\label{prop:holononvan}
Suppose we are in the setting of Lemma \ref{lem:holoeq}.  
Then 
\begin{equation}
\label{eqn:holononvan}
\dim_{\mathbb{C}}
\{\varphi\in {\mathbb{C}}^k:Q_{\lambda}\varphi=0\} \ge 1
\qquad
\text{for all $\lambda \in {\mathbb{C}}$.}
\end{equation}
\end{proposition}
\begin{proof}
Let $\varphi_{\lambda}$ be as in Lemma \ref{lem:holoeq}.  
Then it suffices to show
 \eqref{eqn:holononvan} for $\lambda$
 belonging to the discrete set $\{\lambda \in {\mathbb{C}}:\varphi_{\lambda}=0\}$.  
Take any $\lambda_0$ such that
 $\varphi_{\lambda_0}=0$.  
Let $k$ be the smallest positive integer 
 such that 
\[
 \psi_{\lambda_0}
 :=\left.\frac{\partial^k}{\partial\lambda^k}\right|_{\lambda=\lambda_0}\varphi_{\lambda} \ne 0
\quad
\text{and}
\quad
 \left.\frac{\partial^{j}}{\partial\lambda^{j}}\right|_{\lambda=\lambda_0}\varphi_{\lambda} = 0
\quad
\text{for $0 \le j \le k-1$}.  
\]
By the Leibniz rule,
 $\left.\frac{\partial^k}{\partial\lambda^k}\right|_{\lambda=\lambda_0}
  (Q_{\lambda} \varphi_{\lambda}) = 0
$
 yields $Q_{\lambda_0} \psi_{\lambda_0}=0$, 
 because $\left. \frac{\partial^{j}}{\partial\lambda^{j}}\right|_{\lambda=\lambda_0}
               \varphi_{\lambda} = 0$
 for all $0 \le j \le k-1$.  
Therefore $\psi_{\lambda_0}$
 is a nonzero solution to $Q_{\lambda_0} \varphi=0$, 
showing \eqref{eqn:holononvan}
 for ${\lambda}=\lambda_0$.  
Hence Proposition \ref{prop:holononvan} is proved.  
\end{proof}

As in the proof of Theorem \ref{thm:existSBO}, 
 the implication (i) $\Rightarrow$ (ii) in Theorem \ref{thm:existDSBO}
 follows from Corollary \ref{cor:DSVOgeneric}
 (generic parameters) and the extension result to special parameters
 (Proposition \ref{prop:holononvan}).  
Thus we have completed a proof
 of Theorem \ref{thm:existDSBO}, 
 and in particular,
 of Theorem \ref{thm:vanDiff} (2).

\subsection{Proof of Theorem \ref{thm:VWSBO} (2-b)}
\label{subsec:170213}

In this section, 
 we give a proof of Theorem \ref{thm:VWSBO} (2-b), 
 namely,
 we prove the following proposition.  

\begin{proposition}
[localness theorem]
\label{prop:1532102}
\index{B}{localnesstheorem@localness theorem}
Suppose $[V:W]\ne 0$.  
Suppose that 
\index{A}{1psi@$\Psising$,
          special parameter in ${\mathbb{C}}^2 \times \{\pm\}^2$}
 $(\lambda,\nu,\delta,\varepsilon) \in \Psising$, 
 namely, 
 $(\lambda, \nu)\in {\mathbb{C}}^2$
 and $\delta, \varepsilon \in \{ \pm \}$
 satisfy
\[
  \text{
   $\nu-\lambda \in 2 {\mathbb{N}}$
   when $\delta \varepsilon =+;$
   $\nu-\lambda \in 2 {\mathbb{N}}+1$
   when $\delta \varepsilon =-$.  
  }
\]
Assume further that $\Atbb \lambda \nu {\delta\varepsilon}{V,W} \ne 0$.  
The we have 
\[
 {\operatorname{Hom}}_{G'}
 (I_{\delta}(V, \lambda)|_{G'}, 
  J_{\varepsilon}(W, \nu)
 )
 =
 {\operatorname{Diff}}_{G'}
 (I_{\delta}(V, \lambda)|_{G'}, 
  J_{\varepsilon}(W, \nu)
 ).  
\]
\end{proposition}

We need two lemmas from \cite{sbon}.  
\begin{lemma}
[{\cite[Lem.~11.10]{sbon}}]
\label{lem:sbon1110}
Suppose $D_\mu$ is a differential operator with holomorphic parameter $\mu$,
and $F_\mu$ is a distribution on $\mathbb{R}^n$
that depends holomorphically on $\mu$
having the following expansions:
\begin{align*}
&D_\mu = D_0 + \mu D_1 + \mu^2 D_2 + \dotsb ,
\\
&F_\mu = F_0 + \mu  F_1 + \mu^2 F_2 + \dotsb, 
\end{align*}
where $D_j$ are differential operators
 and $F_i$ are distributions on ${\mathbb{R}}^n$.  
Assume that there exists $\varepsilon>0$
 such that $D_\mu F_\mu = 0$ for any
complex number $\mu$ with $0 < |\mu| < \varepsilon$.  
Then the distributions $F_0$ and $F_{1}$ satisfy
the following differential equations:
\[
D_0 F_0 = 0  \quad\text{and}\quad 
D_0 F_1 + D_1 F_0 = 0.
\]
\end{lemma}
\begin{lemma}
[{\cite[Lem.~11.11]{sbon}}]
\label{lem:sbon1111}
Suppose $h \in \mathcal{D}' (\mathbb{R}^n)$
is supported at the origin.  
Let $E$ be the 
\index{B}{Eulerhomogeneityoperator@Euler homogeneity operator}
Euler homogeneity operator
 $\sum_{\ell=1}^n x_{\ell} \frac {\partial}{\partial x_{\ell}}$ as before.  
If
$(E + A)^2  h = 0$
for some $A \in \mathbb{Z}$ then
$(E + A) h = 0$.  
\end{lemma}

The argument below is partly similar
 to the one in Section \ref{subsec:AVW}, 
 however, 
 we note that the renormalization 
$
   \Attbb {\lambda_0} {\nu_0} {\gamma}{V,W}
$
 in Theorem \ref{thm:170340} is not defined 
 under our assumption
 that $\Atbb {\lambda_0} {\nu_0} {\gamma}{V,W} \ne 0$
 and $(\lambda_0, \nu_0, \delta, \varepsilon) \in \Psising$.  
Instead, 
we shall use the distribution
 $(\Attcal {\lambda} {\nu} {\gamma}{V,W})'$
 on ${\mathbb{R}}^n \setminus \{0\}$, 
 of which we recall \eqref{eqn:AVW+} and \eqref{eqn:AVW-}
 for the definition.   

\begin{proof}
[Proof of Proposition \ref{prop:1532102}]
Take any symmetry breaking operator
\[
 {\mathbb{T}} \in {\operatorname{Hom}}_{G'}(I_{\delta}(V,\lambda_0)|_{G'}, 
 J_{\varepsilon}(W,\nu_0)).  
\]
We write $({\mathcal{T}}_{\infty}, {\mathcal{T}})$
 for the pair of distribution kernels of ${\mathbb{T}}$
 as in Proposition \ref{prop:Tpair}.  
We set $\gamma:= \delta \varepsilon$.

It follows from Proposition \ref{prop:20150828-1231} (3)
 that 
$
   {\mathcal{T}}|_{{\mathbb{R}}^n\setminus \{0\}}
   =
   c' (\Attcal {\lambda_0} {\nu_0} \gamma {V,W})'$
 for some $c' \in {\mathbb{C}}$.  

Suppose $\Atbb {\lambda_0} {\nu_0} {\gamma}{V,W} \ne 0$
 and $\nu_0 - \lambda_0 \in 2 {\mathbb{N}}$
 ($\gamma=+$)
 or $\in 2 {\mathbb{N}}+1$
 ($\gamma=-$).  
As in \eqref{eqn:Anear0}, 
 we expand $\Atcal {\lambda} {\nu_0} {\gamma}{V,W}$
 near $\lambda = \lambda_0$:
\[
  \Atcal {\lambda} {\nu_0} {\gamma}{V,W}
  =
  F_0 + (\lambda-\lambda_0) F_1 + (\lambda-\lambda_0)^2 F_2
  + \cdots,
\]
where $F_j \in {\mathcal{D}}'({\mathbb{R}}^n) \otimes {\operatorname{Hom}}_{\mathbb{C}}(V,W)$.  
We note that $F_0 \ne 0$
 because $\Atbb {\lambda_0} {\nu_0} {\gamma}{V,W} \ne 0$.  
We define a nonzero constant $c$
 by 
\begin{equation}
\label{eqn:1702110}
  c:= \lim_{\mu \to 0}
      \mu \Gamma(\frac \mu 2-l)
    = \frac{2(-1)^l}{l!}.  
\end{equation}
In view of the relation
\[
   \Atcal \lambda \nu + {V,W}|_{{\mathbb{R}}^n \setminus \{0\}}
   =
   \frac{1}{\Gamma(\frac{\lambda-\nu}{2})} (\Attcal \lambda \nu + {V,W})', 
\quad
  \Atcal \lambda \nu - {V,W}|_{{\mathbb{R}}^n \setminus \{0\}}
   =
   \frac{1}{\Gamma(\frac{\lambda-\nu+1}{2})} (\Attcal \lambda \nu - {V,W})',
\]
 we get
\[
  c F_1|_{{\mathbb{R}}^n \setminus \{0\}}
  =
  (\Attcal {\lambda_0} {\nu_0} \gamma {V,W})', 
\]
as in the proof of Theorem \ref{thm:170340} (3).  
We set 
\[
   D_0 := E-\lambda_0 + \nu_0 + n
       = \sum_{j=1}^n x_j \frac{\partial}{\partial x_j}
         -\lambda_0 + \nu_0 + n.  
\]
Applying Lemma \ref{lem:sbon1110}
 to the differential equation \eqref{eqn:Fainv}:
\[
  (E-\lambda + \nu_0 + n) \Atcal \lambda {\nu_0} \gamma {V,W}
  =
  (D_0-(\lambda-\lambda_0))\Atcal \lambda {\nu_0} \gamma {V,W}
  =0, 
\]
we get 
\begin{equation}
\label{eqn:DF01}
D_0 F_0 =0, 
\quad
D_0 F_1 - F_0 =0. 
\end{equation}

We set
\[
   h:={\mathcal{T}}- c c' F_1
  \in 
  {\mathcal{D}}'({\mathbb{R}}^n)
    \otimes 
    \operatorname{Hom}_{{\mathbb{C}}}(V,W). 
\]
Then ${\operatorname{Supp}}\, h \subset \{0\}$.  
Moreover, 
$D_0^2 h =0$
 by $D_0 {\mathcal{T}}=0$ and \eqref{eqn:DF01}.

Applying Lemma \ref{lem:sbon1111}, 
 we get $D_0 h=0$.  
It turn, 
 $c c' F_0=0$
 again by $D_0 {\mathcal{T}}=0$
 and \eqref{eqn:DF01}.  
Therefore,
 if $\Atcal {\lambda_0} {\nu_0} \gamma {V,W} \ne 0$,
 or equivalently,
 if $\Atbb {\lambda_0} {\nu_0} \gamma {V,W} \ne 0$, 
 then we conclude $c'=0$
 because $F_0 \ne 0$.  
Thus ${\mathcal{T}}$ is supported at the origin, 
 and therefore ${\mathbb{T}}$ is a differential operator
 (see Proposition \ref{prop:Tpair} (3)).  

Hence Proposition \ref{prop:1532102} is proved.  
\end{proof}

The above proof implies 
 that the distribution $(\Attcal \lambda \nu \gamma {V,W})' \in {\mathcal{D}}'({\mathbb{R}}^n \setminus \{0\}) \otimes {\operatorname{Hom}}_{\mathbb{C}}(V,W)$
 in \eqref{eqn:AVW+} and \eqref{eqn:AVW-} does not always extend
 to an element of ${\mathcal{S}}{\it{ol}}({\mathbb{R}}^n; V_{\lambda, \delta}, 
                            W_{\nu, \varepsilon})$
 ($\gamma = \delta \varepsilon$):

\begin{proposition}
\label{prop:20170213}
Let $\gamma \in \{\pm\}$. 
Suppose $(\lambda, \nu)\in {\mathbb{C}}^2$ satisfies
\[
  \text{
  $\nu-\lambda \in 2{\mathbb{N}}$
  when $\gamma =+;$
  \,\,
  $\nu-\lambda \in 2{\mathbb{N}}+1$
  when $\gamma =-$.  
 }
\]
If $\Atbb \lambda \nu \gamma {V,W} \ne 0$, 
 then for $\delta, \varepsilon \in \{ \pm \}$
 with $\delta \varepsilon =\gamma$, 
 the restriction map
\[
   {\mathcal{S}}{\it{ol}}({\mathbb{R}}^n; V_{\lambda, \delta}, 
                            W_{\nu, \varepsilon})
   \to 
   {\mathcal{D}}'({\mathbb{R}}^n \setminus \{0\})
    \otimes 
    \operatorname{Hom}_{{\mathbb{C}}}(V,W)
\]
 is identically zero.  
\end{proposition}

\newpage
\section{Minor summation formul{\ae}
 related to exterior tensor $\Exterior^i({\mathbb{C}}^n)$}
\label{sec:section9}

This chapter collects some combinatorial formul{\ae}, 
 which will be used in later chapters to compute the $(K,K')$-spectrum
 for symmetry breaking operators
 between differential forms 
 on spheres $S^n$ and $S^{n-1}$, 
 namely,
 between principal series representations
 $I_{\delta}(V,\lambda)$ of $G$
 and $J_{\varepsilon}(W,\nu)$ of its subgroup $G'$
 in the setting where $(V,W)=(\Exterior^i({\mathbb{C}}^n), \Exterior^j({\mathbb{C}}^{n-1}))$.  
\subsection{Some notation on index sets}
\label{subsec:index}

Let $n$ be a positive integer.  
We shall use the following convention 
 of index sets:
\index{A}{Ini@${\mathfrak{I}}_{n,i}$, index set|textbf}
\begin{equation}
\label{eqn:Indexni}
{\mathfrak{I}}_{n,i}
:=\{I \subset \{1,\cdots,n\}
  : \# I =i\}.   
\end{equation}

\begin{convention}
\label{conv:index}
We use calligraphic uppercase letters 
 ${\mathcal {I}}$, ${\mathcal{J}}$
 instead of Roman uppercase letters $I$, $J$
 if the index set may contain 0.  
That is, 
if we write ${\mathcal{I}} \in {\mathfrak{I}}_{n+1,i}$, 
then 
\[
\text{${\mathcal{I}} \subset \{0,1,\cdots,n\}$
\,\, with \,\,$\# {\mathcal{I}} =i$.  
}
\]
\end{convention}
In later applications for symmetry breaking 
 with respect to $(G,G')=(O(n+1,1),O(n,1))$, 
 the notation ${\mathfrak{I}}_{n+1,i}$
 for subsets of $\{0,1,\cdots,n\}$ will be used
 when we describe the basis
 of the basic $K$-types and $K'$-types, 
 whereas the notation ${\mathfrak{I}}_{n,i}$, 
 ${\mathfrak{I}}_{n-1,i}$
 will be used 
 when we discuss representations of $M$ and $M'$, 
respectively.

\subsubsection{Exterior tensors $\Exterior^i({\mathbb{C}}^n)$}
Let $\{e_1, \cdots, e_n\}$ be the standard basis of ${\mathbb{C}}^n$.  
For $I = \{k_1,k_2,\cdots, k_i\}\in {\mathfrak{I}}_{n,i}$
 with $k_1 < k_2 < \cdots < k_i$,
 we set
\[
   e_I := e_{k_1} \wedge \cdots \wedge e_{k_i} \in {\Exterior}^i({\mathbb{C}}^n).  
\]
Then $\{e_I:I\in {\mathfrak{I}}_{n,i}\}$ forms
 a basis of the exterior tensor space $\Exterior^i({\mathbb{C}}^n)$.  
We define linear maps
\index{A}{prij@$\pr i j$, projection|textbf}
\[
\pr i j \colon {\Exterior}^i({\mathbb{C}}^n) \to {\Exterior}^j({\mathbb{C}}^{n-1}), 
\quad
(j=i-1,i)
\]
by 
\begin{align}
\label{eqn:Tii1}
\pr i i (e_I)=&
\begin{cases}
e_I \hspace{25mm} &\text{if}\ n \notin I, 
\\
0 \qquad &\text{if}\ n \in I, 
\end{cases}
\\
\label{eqn:Tii2}
\pr i {i-1} (e_I)=&
\begin{cases}
0  &\text{if}\ n \notin I, 
\\
(-1)^{i-1} e_{I \setminus \{n\}} \quad &\text{if}\ n \in I.  
\end{cases}
\end{align}

Then we have the direct sum decomposition 
\begin{equation}
\label{eqn:extij}
  {\Exterior}^i({\mathbb{C}}^n) \simeq {\Exterior}^i({\mathbb{C}}^{n-1})
                                     \oplus
                                     {\Exterior}^{i-1}({\mathbb{C}}^{n-1}).  
\end{equation}
\subsubsection{Signatures for index sets}
Let $N \in {\mathbb{N}}_+$.  
In later sections, 
 $N$ will be $n-1$, $n$ or $n+1$.  

For a subset $I \subset \{1,\cdots, N\}$, 
 we define a signature $\varepsilon_I(k)$ by 
\index{A}{0epsilonI@$\varepsilon_I$|textbf}
\[
 \varepsilon_I(k)
 :=
\begin{cases}
 1
 \qquad
  &\text{if}\ k \in I, 
\\
 -1
\qquad
  &\text{if}\ k \not\in I, 
\end{cases}
\]
 and a quadratic polynomial 
\index{A}{QIb@$Q_I(b)$, quadratic polynomial|textbf}
$
Q_I(y)
$
 by 
\begin{equation}
\label{eqn:QI}
Q_I(y):= \sum_{\ell \in I} {y_{\ell}}^2
\qquad
\text{for }\,\, 
 y =(y_1, \cdots, y_N)\in {\mathbb{R}}^N.  
\end{equation}
We note that
\[
  2 Q_I(y) - |y|^2
  = 
  \sum_{k =1}^N \varepsilon_I(k) {y_k}^2.  
\]
For $I , J \subset {\mathfrak{I}}_{N,i}$, 
 we set
\[
|I-J|:= \# I - \#(I \cap J)=\# J - \#(I \cap J).  
\]
By definition,
 $|I-J|=0$ if and only if $I=J$;
 $|I-J|=1$
 if and only if there exist $K \in{\mathfrak{I}}_{N,i-1}$
 and $p$, $q \not \in K$
 with $p \ne q$ 
such that $I = K \cup \{p\}$
 and $J = K \cup \{q\}$.

\begin{definition}
\label{def:sign}
For $I \subset \{1,2,\cdots,n\}$
 and $p,q \in {\mathbb{N}}$, 
 we set 
\begin{align*}
\index{A}{sgnIp@$\operatorname{sgn}(I;p)$|textbf}
\operatorname{sgn}(I;p)
:=&(-1)^{\# \{r \in I:r < p\}}, 
\\
\index{A}{sgnIpq@$\operatorname{sgn}(I;p,q)$|textbf}
\operatorname{sgn}(I;p,q)
:=&(-1)^{\# \{r \in I:\operatorname{min}(p,q)<r < \operatorname{max}(p,q)\}}.  
\end{align*}
\end{definition}

The following lemma is readily seen from the definition.  
\begin{lemma}
\label{lem:sgn}
For $I \subset \{1,2,\cdots,n\}$
 and $p,q \in {\mathbb{N}}$, 
 we have 
\begin{equation*}
\operatorname{sgn}(I;p)
\operatorname{sgn}(I;q)
=
\begin{cases}
\operatorname{sgn}(I;p,q)
\quad
&
\text{ if }\operatorname{min}(p,q) \notin I, 
\\
-\operatorname{sgn}(I;p,q)
\quad
&
\text{ if }\operatorname{min}(p,q) \in I.  
\end{cases}
\end{equation*}
\end{lemma}

For $y=(y_1,\cdots,y_N) \in {\mathbb{R}}^N$, 
 we define quadratic polynomials 
\index{A}{SIJ@$S_{I J}$|textbf}
$
S_{I J}(y)
$ 
by 
\begin{equation}
\label{eqn:SIJ}
S_{I J}(y)
:=
\begin{cases}
\sum_{k=1}^{N} \varepsilon_I(k) y_k^2
\quad
&\text{ if } I =J, 
\\
2 \operatorname{sgn}(K;p,q) y_p y_q
\quad
&\text{ if } I = K \cup \{p\}, J = K \cup \{q\}, 
\\
0
\quad
&\text{ if } |I- J| \ge 2, 
\\
\end{cases}
\end{equation}
where we write $I=K \cup \{p\}$
 and $J=K \cup \{q\}$
 $(p \ne q)$
 when $|I - J|=1$.  

It is convenient to set

\begin{equation}
\label{eqn:Sempty}
S_{\emptyset \emptyset}(y)=-\sum_{k=1}^N y_k^2.  
\end{equation}

\subsection{Minor determinant 
 for $\psi:{\mathbb{R}}^N \setminus \{0\}\to O(N)$
}
\label{subsec:psiN}

We introduce the following map:
\index{A}{1psinl@$\psi_n(\cdot ;\lambda)$}
\begin{equation}
\label{eqn:psilmd}
\psi_N \colon \mathbb{R}^N \times {\mathbb{C}}
        \to M(N,{\mathbb{C}}),
\quad
(y; \lambda) \mapsto I_N - \lambda\, y \, {}^t \! y.  
\end{equation}
Here we have used a similar notation to the map $\psi_N(y)$ defined in \eqref{eqn:psim}.  
In fact,
 the map \eqref{eqn:psilmd} may be thought of as an extension
 of the previous one,
 since its special value 
 at $\lambda = \dfrac{2}{|y|^2}$ recovers \eqref{eqn:psim} by 
\begin{equation}
   \psi_N(y)=\psi_N(y;\frac{2}{|y|^2})
\qquad
\text{for $y \in {\mathbb{R}}^N \setminus \{0\}$}.  
\end{equation}

For $I,J \subset \{1,2,\cdots,N\}$
 with $\# I= \# J$, 
 the minor determinant
 of $A=(A_{i j})_{1 \le i, j \le N} \in M(N,{\mathbb{R}})$
 is denoted by 
\[
\det A_{I J}:=\det(A_{i j})_{\substack{i \in I \\ j \in J}}.  
\]
Then the exterior representation
\[
  \sigma:O(N) \to GL_{\mathbb{C}}(\Exterior^k ({\mathbb{C}}^N))
\]
is given by
\begin{equation}
\label{eqn:exrep}
\sigma(A) e_J 
=
\sum_{J' \in {\mathfrak{I}}_{N, k}}
(\det A)_{J' J}e_{J'}.  
\end{equation}
It follows from \eqref{eqn:exrep}
 that for $A, B \in O(N)$
 we have
\begin{equation}
\label{eqn:klminor}
\det(A B)_{J J'}
=\sum_{J'' \in {\mathfrak{I}}_{N, j}}
(\det A)_{J J''}
(\det B)_{J'' J'}
\end{equation}

\begin{lemma}
\label{lem:psidet}
Suppose $I$, $J \subset \{1,\cdots,N\}$
 with $\# I = \# J$.  
\par\noindent
{\rm{(1)}}\enspace
For $(y;\lambda) \in {\mathbb{R}}^N \times {\mathbb{C}}$,
\[
  \det \psi_N(y; \lambda)_{IJ}
  =
  \begin{cases}
  1 -\lambda Q_I(y)
  \quad
  &\text{if } I=J,
  \\
  -\lambda \operatorname{sgn}(K;p,q) y_p y_q
  &\text{if } I=K \cup \{p\}, \, J = K \cup \{q\},
  \\
  0
  &\text{if } |I-J| \ge 2.  
  \end{cases}
\]
\par\noindent
{\rm{(2)}}\enspace
For $y \in {\mathbb{R}}^N \setminus \{0\}$,
\index{A}{1psin@$\psi_n$}
\index{A}{0epsilonI@$\varepsilon_I$}
\index{A}{SIJ@$S_{I J}$}
\begin{align*}
  \det \psi_N(y)_{I J}
=&
  -\frac {1}{|y|^2} S_{I J}(y)
\\
=&
 \frac{1}{|y|^2} \times
 \begin{cases}
  -\sum_{l =1}^{N} \varepsilon_I(l) y_l^2 
 \quad
  &\text{if } I=J,
  \\
  -2 \operatorname{sgn}(K;p,q) y_p y_q
  &\text{if } I=K \cup \{p\}, J=K \cup \{q\},
  \\
  0
  &\text{if } |I-J| \ge 2.  
  \end{cases}
\end{align*}
\end{lemma}

\begin{proof}
(1) \enspace
Suppose $I=J$.  
Since the symmetric matrix $y \ {}^{t\!}y$ is of rank 1,
 its characteristic polynomial has zeros of order $N-1$:
\[
    \det (\mu I_N - y \ {}^{t\!}y)
   =\mu^N-\mu^{N-1}(\operatorname{Trace} y \ {}^{t\!}y)
   =\mu^N-\mu^{N-1}\sum_{j=1}^{N} y_j^2, 
\]
 and therefore
\[
    \det (I_N - \lambda y \ {}^{t\!}y)= 1- \lambda \sum_{j=1}^{N}y_j^2.  
\]
Applying this to the principal minor
 of size $\# I$,
 we get the first formula.  

\par\indent
Next suppose $|I-J|=1$.  
We may write as $I=K \cup \{p\}$, $J=K \cup \{q\}$.  
Then the $q$-th column vector of the minor matrix
 $(I_N - \lambda y \ {}^{t\!}y)_{IJ}$
 is of the form $-\lambda y_q(y_i)_{i\in I}$.  
Adding this vector multiplied by the scalar
 $(- y_j/y_q)$ to the $j$-th column vector
 for $j \in J \setminus \{q\}$,
 we get 
\begin{align*}
  \det \psi_N(y; \lambda)_{I J}
  =&
  \operatorname{sgn}(K;q) 
  \det (-\lambda y_q(y_i)_{i \in I}, (\delta_{i j})_{\substack{i \in I \\ j \in K}})
\\
=&
-\lambda \operatorname{sgn}(K;p)
 \operatorname{sgn}(K;q)
 y_p y_q.  
\end{align*}
Hence the second formula follows from Lemma \ref{lem:sgn}.  
The third one is proved similarly.  
\par\noindent
(2)\enspace
Substitute $\lambda=\frac{2}{|y|^2}$.  
\end{proof}
As a special case of Lemma \ref{lem:psidet} (2) with $N=n+1$, 
 we have the following:
\begin{lemma}
\label{lem:psid}
For ${\mathcal{I}}, {\mathcal{J}} \in {\mathfrak{I}}_{n+1,i}$
 and $b \in {\mathbb{R}}^n$, 
 we have
\[
 \det \psi_{n+1}(1,b)_{{\mathcal{I}} {\mathcal{J}}}
 =
\frac{-1}{1+|b|^2} S_{{\mathcal{I}} {\mathcal{J}}}(1,b).  
\]
Here $(1,b):=(1,b_1, \cdots, b_n) \in {\mathbb{R}}^{n+1}$.  
\end{lemma}

\subsection{Minor summation formul{\ae}}
\label{subsec:minorsum}

We collect minor summation formul{\ae} 
 that we shall need
 in computing the $(K,K')$-spectrum of symmetry breaking operators
 for \lq\lq{basic $K$-types}\rq\rq.

We recall from \eqref{eqn:QI}
 that $Q_I(b)=\sum_{k \in I} b_k^2$.  
\begin{lemma}
\label{lem:msum}
Suppose $I \in {\mathfrak{I}}_{n,i}$.  
For $b \in {\mathbb{R}}^n$ and $s,t \in {\mathbb{C}}$, 
 we have:
\begin{alignat}{2}
\label{eqn:msum}
&{\text{\rm{(1)}}}\,\,
&&\sum_{J \in {\mathfrak{I}}_{n,i}}
\det \psi_n(b;s)_{I J}
\det \psi_n(b;t)_{I J}
=
1-(s+t)Q_I(b) + s t |b|^2 Q_I(b).  
\\
&{\text{\rm{(2)}}}
&&\sum_{J \in {\mathfrak{I}}_{n,i}}
\det \psi_{n+1}(1,b;s)_{I\cup \{0\},  J\cup \{0\}}
\det \psi_n(b;t)_{I J}
\notag
\\
&
&&\qquad\qquad\qquad\qquad
=
1- s -(s+t-st)Q_I(b) + s t |b|^2 Q_I(b).  
\label{eqn:msum2}
\\
\label{eqn:msum3}
&{\text{\rm{(3)}}}
&&
\sum_{J \in {\mathfrak{I}}_{n,i}}
  \det \psi_{n+1}(1,b;s)_{I J}
\,
  \det \psi_n(b;t)_{I J}
  =
   1-(s+t) Q_I(b) + s t |b|^2 Q_I(b).  
\end{alignat}
\end{lemma}
\begin{proof}
(1)\enspace
By Lemma \ref{lem:psidet}, 
the left-hand side is equal
 to 
\begin{align*}
&(1-s Q_I(b))(1-t Q_I(b))
+\sum_{\substack{k \in I\\ l \notin I}} s t b_k^2 b_l^2
\\
=&
1-(s+t)Q_I(b) + s t Q_I(b)^2 + st Q_I(b)(|b|^2-Q_I(b)), 
\end{align*}
whence the equation \eqref{eqn:msum}.  
\par\noindent
(2)\enspace
By Lemma \ref{lem:psidet}, 
the left-hand side is equal
 to 
\[
(1-s (1+Q_I(b)))(1-t Q_I(b))
+s t \sum_{\substack{k \in I\\ l \notin I}} b_k^2 b_l^2, 
\]
whence the equation \eqref{eqn:msum2}.  
\par\noindent
(3)\enspace
By Lemma \ref{lem:psidet}, 
 the left-hand side is equal to 
\begin{align*}
   &(1-s Q_I(b))(1-t Q_I(b))
   + s t \sum_{\substack{k \in I\\ l \notin I}} b_k^2 b_l^2
\\
   =&
   1-(s+t) Q_I(b) + s t Q_I(b)^2 + s t Q_I(b) (|b|^2-Q_I(b)), 
\end{align*}
whence the equation \eqref{eqn:msum3}.  
\end{proof}

The following proposition 
 will be used in obtaining the closed formul{\ae}
 of the $(K,K')$-spectrum
 of the Knapp--Stein intertwining operators (Proposition \ref{prop:TminK})
 and the ones of the regular symmetry breaking operators
 (Theorem \ref{thm:153315}).  
\begin{proposition}
\label{prop:msum}
For $I \in {\mathfrak{I}}_{n,i}$, 
we have:
\begin{alignat}{2}
&\text{\rm{(1)}}
&&\sum_{J \in {\mathfrak{I}}_{n,i}}
\det \psi_n(b;\frac{2}{1+|b|^2})_{I J}
\det \psi_n(b)_{I J}
=
1- \frac{2 Q_I(b)}{(1+|b|^2) |b|^2}.  
\notag
\\
&\text{\rm{(2)}}
&&\sum_{J \in {\mathfrak{I}}_{n,i}}
\det \psi_{n+1}(1,b)_{I\cup \{0\},  J\cup \{0\}}
\det \psi_n(b)_{I J}
=
\frac{-1+|b|^2}{1+|b|^2} + \frac{2 Q_I(b)}{(1+|b|^2) |b|^2}.  
\notag
\\
&\text{\rm{(3)}}
&&\sum_{J \in {\mathfrak{I}}_{n,i}}
\left(\det \psi_n(b;\frac{2}{1+|b|^2})_{I J}
+
\det \psi_{n+1}(1,b)_{I\cup \{0\},  J\cup \{0\}}
\right)
\det \psi_n(b)_{I J}
=
\frac{2 |b|^2}{1+|b|^2}.  
\notag
\\
\label{eqn:kbsum}
&\text{\rm{(4)}}
&&\sum_{J \in {\mathfrak{I}}_{n,i}} \det \psi_{n+1}(1,b)_{I J}
\det \psi_n(b)_{I J}
=
1-\frac{2Q_I(b)}{(1+|b|^2)|b|^2}.  
\end{alignat}
\end{proposition}

\begin{proof}
The assertions (1), (2), and (4) are special cases
 of Lemma \ref{lem:msum} (1), (2), and (3), 
 respectively,
 with $s = \dfrac{2}{1+|b|^2}$ and $t=\dfrac {2}{|b|^2}$.  
The third one follows from the first two.  
\end{proof}

\begin{lemma}
\label{lem:S0n}
For $I \in {\mathfrak{I}}_{n-1,i-1}$, 
\[
  \sum_{J \in {\mathfrak{I}}_{n,i}}
  \det \psi_{n+1}(1,b)_{I \cup \{0\},  J}
  \det \psi_{n}(b)_{I \cup \{n\},  J}
  =
  \frac{2(-1)^{i+1} b_n}{1+|b|^2}.  
\]
\end{lemma}

\begin{proof}
Since $0 \notin J$, 
the summand vanishes
 except for the following two cases:
\par\noindent
Case 1)\enspace $J=I \cup \{ n \}$.  
\par\noindent
Case 2)\enspace $J=I \cup \{ p \}$
  for some $p \in \{1,2,\cdots,n-1\} \setminus I$.

By Lemma \ref{lem:psidet}, 
 we get 
\begin{align*}
&  (1+|b|^2)|b|^2
  \sum_{J \in {\mathfrak{I}}_{n,i}}
  \det \psi_{n+1}(1,b)_{I \cup \{0\},  J}
\ \det \psi_{n}(b)_{I \cup \{n\},  J}  
\\
=& (-2 \operatorname{sgn}(I;0,n) b_n)
   (|b|^2 -2Q_{I}(b)-2b_n^2)
\\
  &+
  \sum_{p \in \{1,2,\cdots,n-1\}\setminus I}
  (-2 \operatorname{sgn}(I;0,p) b_p)
  (-2 \operatorname{sgn}(I;p,n) b_p b_n)
\\
=&2(-1)^{i+1} b_n(2 Q_{I}(b)+2b_n^2-|b|^2)
  +
  4(-1)^{i+1}b_n (|b|^2 - Q_{I}(b) -b_n^2)
\\
=& 2(-1)^{i+1} |b|^2 b_n.  
\end{align*}
Hence Lemma \ref{lem:S0n} is proved.  
\end{proof}
\begin{lemma}
\label{lem:Sn0}
For $I \in {\mathfrak{I}}_{n-1,i}$, 
\[
  \sum_{J \in {\mathfrak{I}}_{n,i}}
  \det \psi_{n+1}(1,b)_{I \cup \{n\}, J \cup \{0\}} \det \psi_n(b)_{I J}
=
 \frac{2(-1)^{i+1} b_n}{1+|b|^2}.  
\]
\end{lemma}
\begin{proof}
Since $0 \notin I$, 
$| (I \cup \{n\}) - (J \cup \{0\}) |\le 1$
 holds 
 in the following two cases:
\par\noindent
Case 1. \enspace
$I=J$.
\par\noindent
Case 2. \enspace
$I=K \cup \{p\}$
 and $J=K \cup \{n\}$
 for some $K \in {\mathfrak{I}}_{n-1,i-1}$.  
\par
In Case 1, 
\[
   \det \psi_{n+1}(1,b)_{I \cup \{n\}, J \cup \{0\}} 
   \det \psi_n(b)_{I I}
=
 \frac{2(-1)^{i+1} b_n}{1+|b|^2}
 \times 
(1-\frac{2 Q_I(b)}{|b|^2}).  
\]
\par
In Case 2, 
\begin{align*}
  &\det \psi_{n+1}(1,b)_{K \cup \{ p, n\}, K \cup \{ 0, n \}}
   \det \psi_{n}(b)_{K \cup \{ p \}, K \cup \{ n \}}
\\
&= \frac{-2 {\operatorname{sgn}}(K \cup \{n\};0,p)b_p}{1+|b|^2}
   \times 
   \frac{-2 {\operatorname{sgn}}(K;p,n)b_p b_n}{|b|^2}
\\
&= (-1)^{i-1} \frac{4}{(1+|b|^2) |b|^2} b_p^2 b_n.  
\end{align*}
Adding the term in Case 1 and taking the summation of the terms
 over $p \in I$ in Case 2, 
 we get the lemma.  
\end{proof}

\newpage
\section{The Knapp--Stein intertwining operators revisited:
 Renormalization and $K$-spectrum}
\label{sec:psdetail}

In this chapter,
 we discuss the classical Knapp--Stein operators,
 which may be viewed
 as a baby case of symmetry breaking operators
 ({\it{i.e.}},  $G=G'$ case).  
We determine the $(K,K)$-spectrum
 ($K$-spectrum, for short)
 of the matrix-valued Knapp--Stein operators
 $\Ttbb \lambda {n-\lambda} V \colon I_{\delta}(V,\lambda) \to I_{\delta}(V,n-\lambda)$, 
 see \eqref{eqn:KnappStein}
 in the case where $V=\Exterior^i({\mathbb{C}}^n)$.  
We also study the renormalization
of the operator $\Ttbb \lambda {n-\lambda} V$
 when it vanishes,
 see Section \ref{subsec:renormT}.  

\subsection{Basic $K$-types in the compact picture}
\label{subsec:minKK}
Let $(\mu,U)$ be an irreducible representation
 of a compact Lie group $K$, 
 and $(\sigma,V)$ that of a subgroup $M$.  
The classical Frobenius reciprocity tells
 that $\mu$ occurs in the induced representation 
 ${\operatorname{Ind}}_M^K \sigma$
 if and only if ${\operatorname{Hom}}_M(\mu|_M, \sigma) \ne \{0\}$.  
In this section we provide a concrete realization
 of $(\mu,U)$
 in the space $C^{\infty}(K/M, {\mathcal{V}})$
 of global sections 
 for the $K$-equivariant vector bundle 
 ${\mathcal{V}}=K \times_M V$
 which we will use later.  
\begin{lemma}
\label{lem:KM}
\begin{enumerate}
\item[{\rm{(1)}}]
Let $(\mu,U)$ be a finite-dimensional representation
 of a compact Lie group $K$.  
The left regular representation
 on $C^{\infty}(K,U)$ is defined by $f(\cdot) \mapsto f({\ell}^{-1} \cdot)$
 for $f \in C^{\infty}(K,U)$ and ${\ell} \in K$, 
 where we regard $U$ just as a vector space.  
By assigning to $u \in U$, 
 the function $f_u \colon K \to U$ is defined by 
 $f_u(k):=\mu(k)^{-1}u$.  
Then the $K$-module $U$ can be embedded 
 as a submodule of the left regular representation $C^{\infty}(K,U)$
 by
\[
  U \to C^{\infty}(K,U),
  \qquad
  u \mapsto f_u.  
\]
\item[{\rm{(2)}}]
Let $V$ be a vector space over ${\mathbb{C}}$, 
 and $\pr U V :U \to V$ a linear map.  
Then we have a $K$-homomorphism
\[
  U \to C^{\infty}(K,V),
  \qquad
  u \mapsto \pr U V \circ f_u.  
\]
\item[{\rm{(3)}}]
Suppose that $\sigma:M \to GL_{\mathbb{C}}(V)$ is
 a representation of a subgroup $M$ of $K$ 
 and that $\pr UV$ is an $M$-homomorphism.  
Then we have a well-defined $K$-homomorphism
\[
  U \to C^{\infty}(K/M,{\mathcal{V}}),
  \qquad
  u \mapsto \pr U V \circ f_u,  
\]
where we identify the space
 of smooth sections
 for ${\mathcal{V}}:=K \times_M V$
 over $K/M$ 
with the space of $M$-invariant elements
\[
  C^{\infty}(K,V)^M
  :=
  \{F \in C^{\infty}(K,V)
    :
    F(\cdot\, m)=\sigma(m)^{-1}F(\cdot)
   \quad
   \text{for all }
   m \in M\}.  
\]
\end{enumerate}
\end{lemma}
\begin{proof}
The detailed formulation of each statement gives a proof by itself.  
\end{proof}
Applying Lemma \ref{lem:KM}
 to differential forms
 on the sphere, 
 we obtain:
\begin{example}
\label{ex:KMi}
Let $K:=O(n+1)$, 
 and $\sigma$ be the $i$-th exterior tensor representation
 of the subgroup $M:=O(n)$ on $V:= \Exterior^i({\mathbb{C}}^{n})$.  
Then the vector bundle ${\mathcal{V}}=K \times_M V$ is identified
 with the $i$-th exterior tensor
 of the cotangent bundle
 of the $n$-sphere $S^n \simeq K/M$, 
 and we may identify $C^{\infty}(K,V)^M \simeq  C^{\infty}(K/M,{\mathcal{V}})$
 with the space ${\mathcal{E}}^i(S^n)$
 of differential $i$-forms on $S^n$.  
Suppose that $\mu$ is the $k$-th exterior tensor representation
 of $K=O(n+1)$ on $U:= \Exterior^k({\mathbb{C}}^{n+1})$.  
For $k=i$ or $i+1$, 
 the projection 
$
  \pr {k}i \colon \Exterior^{k}({\mathbb{C}}^{n+1}) \to \Exterior^i({\mathbb{C}}^n)
$, 
 see \eqref{eqn:Tii1} and \eqref{eqn:Tii2}, 
 is an $M$-homomorphism, 
 and therefore,
 Lemma \ref{eqn:fI} gives a concrete realization
 of the $K$-module $U=\Exterior^k({\mathbb{C}}^{n+1})$ in ${\mathcal{E}}^i(S^n)\simeq C^{\infty}(K,V)^M$ as below.  
Let $\{e_0, e_1, \cdots, e_n\}$ be the standard basis of ${\mathbb{C}}^{n+1}$, 
 and $\{e_{\mathcal{I}}: {\mathcal{I}} \in {\mathfrak {I}}_{n+1,k}\}$
 the standard basis of $\Exterior^k({\mathbb{C}}^{n+1})$.

We treat the cases $k=i$ and $i+1$, 
 separately.  
In what follows,
 we use Convention \ref{conv:index}
 for the index set ${\mathfrak {I}}_{n+1,k}$.
See also Section \ref{subsec:psiN} 
 for minor determinant $(\det A)_{I J}$
 of $A \in M(N,{\mathbb{R}})$.  
\par\noindent
{\bf{Case 1.}}\enspace
Suppose $k=i$.  
Then 
\index{A}{01oneI@${\bf{1}}^{\mathcal{I}}$|textbf}
$
{\bf{1}}^{\mathcal{I}}:=\pr i i \circ f_{e_{\mathcal{I}}}
$ is a map given by 
\begin{equation}
\label{eqn:fI}
O(n+1) \to \Exterior^i({\mathbb{C}}^n), 
\qquad
  k \mapsto 
  {\bf{1}}^{\mathcal{I}}(k)=\sum_{J \in {\mathfrak {I}_{n,i}}}
  (\det k)_{{\mathcal{I}}J} e_J.  
\end{equation}
Thus ${\bf{1}}^{\mathcal{I}}$ is regarded as an element
 of 
$
   C^{\infty}(O(n+1),\Exterior^i({\mathbb{C}}^n))^{O(n)}\simeq {\mathcal{E}}^i(S^n).  
$
\end{example}

\par\noindent
{\bf{Case 2.}}\enspace
Suppose $k=i+1$.  
Then 
\index{A}{hI@$h^{\mathcal{I}}$|textbf}
$
h^{\mathcal{I}} := (-1)^i \pr {i+1} i \circ f_{e_{\mathcal{I}}}
$
 is a map given by 
\begin{equation}
\label{eqn:hJ}
O(n+1) \to \Exterior^i({\mathbb{C}}^n),
\quad
k \mapsto 
  h^{\mathcal{I}}(k) =\sum_{J \in {\mathfrak{I}}_{n,i}}
          (\det k)_{{\mathcal{I}},J \cup \{0\}}e_J, 
\end{equation}
which is again regarded as an element
 of ${\mathcal{E}}^i(S^n)$.  
We remark
 that the projection
\[
   \pr {i+1} j \colon
   \Exterior^{i+1}({\mathbb{C}}^{n+1}) \to \Exterior^i({\mathbb{C}}^n)
\]
is given by \lq\lq{removing}\rq\rq\
 $e_0$, 
 whereas the projection in \eqref{eqn:Tii2} was by \lq\lq{removing}\rq\rq\
 $e_n$.

By Lemma \ref{lem:KM}, 
 we obtain injective $O(n+1)$-homomorphisms
\begin{alignat*}{3}
&\Exterior^{i}({\mathbb{C}}^{n+1})
&&\to 
{\mathcal{E}}^i(S^n), 
\qquad
&&e_{\mathcal{I}} \mapsto {\bf{1}}^{\mathcal{I}}, 
\\
&\Exterior^{i+1}({\mathbb{C}}^{n+1})
&&\to 
{\mathcal{E}}^i(S^n),  
\qquad
&&e_{\mathcal{I}} \mapsto h^{\mathcal{I}}.  
\end{alignat*}

\subsection{$K$-picture and $N$-picture of principal series representations}
\label{subsec:psKN}

Let $(\sigma, V) \in \widehat {O(n)}$, 
 $\delta \in \{\pm 1\}$, 
 and $\lambda \in {\mathbb{C}}$.  
We recall from Section \ref{subsec:smoothI}
 that the principal series representation
\[
 I_{\delta}(V,\lambda)
 =
 {\operatorname{Ind}}_P^G(V \otimes \delta \otimes {\mathbb{C}}_{\lambda})
\]
 of $G=O(n+1,1)$
 is realized on the Fr{\'e}chet space
 $C^{\infty}(G/P, {\mathcal{V}}_{\lambda,\delta})$
 of smooth sections for the homogeneous vector bundle
 $G \times_P V_{\lambda, \delta}$ 
over the real flag manifold $G/P$, 
 see \eqref{eqn:Vlmdbdle}. 

\subsubsection{Explicit $K$-finite vectors in the $N$-picture} 
In this subsection we review
 the 
\index{B}{Kpicture@$K$-picture}
$K$-picture and
\index{B}{Npicture@$N$-picture}
$N$-picture
 of the principal series representation $I_{\delta}(V,\lambda)$,
 and provide a concrete formula
 connecting the two pictures.

As we saw in \eqref{eqn:iotaN}, 
 the noncompact picture ($N$-picture) of $I_{\delta}(V,\lambda)$
 is given by
\[
  \iota_N^{\ast} \colon 
  I_{\delta}(V,\lambda) \hookrightarrow C^{\infty}({\mathbb{R}}^n) \otimes V, 
\quad
  F \mapsto f(b):=F(n_-(b)), 
\]
 as the pull-back of sections
 via the coordinate map of the open Bruhat cell
$
  \iota_N \colon {\mathbb{R}}^n \hookrightarrow G/P, 
$
$b \mapsto n_-(b) \cdot o$, 
 where $n_- \colon {\mathbb{R}}^n \overset \sim \to N_-$
 is defined in \eqref{eqn:nbar}.

Next,
 let $V_{\delta}$ denote the outer tensor product representation
 $V \boxtimes \delta$ of $M=O(n) \times O(1)$.  
Then the diffeomorphism
\index{A}{0iotaK@$\iota_K$|textbf}
 $\iota_{K} \colon K/M \overset \sim \to G/P$ induces an isomorphism
 $\iota_K^{\ast}({\mathcal{V}}_{\lambda,\delta}) \simeq K \times_M V_{\delta}$
 as $K$-equivariant vector bundles
 over $K/M$, 
 and hence $K$-isomorphisms
 between the space of sections:
\[
  \iota_K^{\ast} \colon
  I_{\delta}(V,\lambda)
  \overset \sim \to
  C^{\infty}(K/M, K \times_M V_{\delta})
  \simeq 
  (C^{\infty}(K) \otimes V_{\delta})^{M}, 
\]
which is referred to as the {\it{$K$-picture}}
 of $I_{\delta}(V,\lambda)$.

The transform from the $K$-picture to the $N$-picture
 is given by 
\begin{equation}
\label{eqn:KtoN}
  \iota_{\lambda}^{\ast}:=\iota_{N}^{\ast} \circ (\iota_{K}^{\ast})^{-1}
  \colon 
  (C^{\infty}(K) \otimes V_{\delta})^M 
  \hookrightarrow 
  C^{\infty}({\mathbb{R}}^n) \otimes V.  
\end{equation}
Then the three realizations
 of the principal series representation
 $I_{\delta}(V,\lambda)$ of $G$ are summarized as below.  
\index{A}{0iotalmdd@$\iota_{\lambda}^{\ast}$}
$$
\xymatrix@R-=1.3pc@C+=0.9cm{
& C^{\infty}(G/P, {\mathcal{V}}_{\lambda,\delta})
\ar[lddd]
\ar[rddd]^{\;\; \iota_N^{\ast}}
&{} 
\\
\qquad\qquad\qquad\qquad\qquad{}_{\iota_K^{\ast}}
&
&
\\
&  
&{}
\\
\text{($K$-picture)}
\quad
(C^{\infty}(K) \otimes V_{\delta})^M 
\ar[rr]_{\hspace{25pt}\iota_{\lambda}^{\ast}}
&{}
&C^{\infty}({\mathbb{R}}^n) \otimes V
\quad
\text{($N$-picture)}
}
$$

To compute $\iota_{\lambda}^{\ast}$, 
 we recall from Lemma \ref{lem:0.2}
 that the map
\[
     k \colon {\mathbb{R}}^n \to SO(n+1) \subset K =O(n+1) \times O(1), 
\]
 see \eqref{eqn:kb}, 
 induces the following commutative diagram:

\vskip 1pc
\[
\begin{xy}
(0,0) *{{\mathbb{R}}^n}="A",
(23,-1)*{\underset{n_-}{\overset \sim \longrightarrow} N_-
        \hookrightarrow G/P {\overset \sim \longleftarrow}}="D",
(45,0)*{K/M}="B",
(45,20)*{ K }="C",
\ar "C";"B"
\ar^k "A";"C" 
\end{xy}
\]

\begin{lemma}
\label{lem:fKN}
Suppose $F \in (C^{\infty}(K) \otimes V_{\delta})^M$.  
Then we have 
\begin{equation}
\label{eqn:fKN}
(\iota_{\lambda}^{\ast} F)(b)
=
(1+|b|^2)^{-\lambda}F(k(b))
\quad
\text{for all $b \in {\mathbb{R}}^n$}.  
\end{equation}
Here $k(b) \in SO(n+1)$ is viewed as an element of $K$
 on the right-hand side.  
\end{lemma}
\begin{proof}
We define $t \in {\mathbb{R}}$
 by $e^t=1+|b|^2$.  
It follows from Lemma \ref{lem:0.2}
 that 
\index{A}{n2@$n_-\colon {\mathbb{R}}^n \to N_-$}
\index{A}{kb@$k(b)$}
\begin{align*}
(\iota_{\lambda}^{\ast} F)(b)=({\iota_K^{\ast}}^{-1} F)(n_-(b))
=&({\iota_K^{\ast}}^{-1} F)
\left(
\begin{pmatrix}
k(b) & 0
\\
0& 1
\end{pmatrix}
e^{t H}
n_+(\frac{-b}{1+|b|^2})
\right)
\\
=&
(1+|b|^2)^{-\lambda}
F\left(
\begin{pmatrix}
k(b) & 0
\\
0& 1
\end{pmatrix}
\right). 
\end{align*}
Hence the lemma is verified.  
\end{proof}

\subsubsection{Basic $K$-types in the $N$-picture}
\label{subsec:minKtype}
We recall from the $K$-type formula
 (Lemma \ref{lem:KtypeIi})
 that the principal series representation
\index{A}{Ideltai@$I_{\delta}(i, \lambda)$}
$
I_{\delta}(i,\lambda)
$
 of $G=O(n+1,1)$ contains two 
\index{B}{basicKtype@basic $K$-type}
\lq\lq{basic $K$-types}\rq\rq\
 $\mub(i, \delta)=\Exterior^i({\mathbb{C}}^{n+1}) \boxtimes \delta$
 and $\mus(i,\delta)=\Exterior^{i+1}({\mathbb{C}}^{n+1}) \boxtimes (-\delta)$
 for $0 \le i \le n$.

In this section,
 we write down explicit $K$-finite vectors
 belonging to $\mub(i,\delta)$ and $\mus(i,\delta)$
 in the noncompact picture.

Let ${\bf{1}}^{\mathcal{I}}$ and $h^{\mathcal{I}}$ be the elements
 in ${\mathcal{E}}^i(S^n) \simeq C^{\infty}(O(n+1),V)^{O(n)}$
 constructed in Example \ref{ex:KMi}, 
 where we take $V$ to be $\Exterior^i({\mathbb{C}}^n)$.  
We note that the pair
\[
   (K,M) = (O(n+1) \times O(1), O(n) \times O(1))
\]
is not exactly the same with the pair
 $(O(n+1), O(n))$ 
 in Example \ref{ex:KMi}, 
however, 
 the diffeomorphism
 $O(n+1)/O(n) \overset \sim \to K/M$
 induces the following isomorphisms
\[
   {\mathcal{E}}^i(S^n)
  \simeq
   C^{\infty}(O(n+1) \otimes V)^{O(n)} \overset \sim \leftarrow
  C^{\infty}(K \otimes V_{\delta})^M.  
\]
Thus we may regard that 
 $\{{\bf{1}}^{\mathcal{I}}: {\mathcal{I}} \in {\mathfrak{I}}_{n+1,i}\}$
 is a basis of $\mub (i,\delta)$
 and 
 $\{h ^{\mathcal{I}}: {\mathcal{I}} \in {\mathfrak{I}}_{n+1,i+1}\}$
 is a basis of $\mus (i,\delta)$.  
Applying the map
 $\iota_{\lambda}^{\ast} \colon C^{\infty}(K \otimes V_{\delta})^M 
 \hookrightarrow 
 C^{\infty}({\mathbb{R}}^n) \otimes V$ 
 (see \eqref{eqn:fKN}), 
 we set
\begin{alignat*}{2}
{\bf{1}}_{\lambda}^{\mathcal{I}}:=& \iota_{\lambda}^{\ast} {\bf{1}}^{\mathcal{I}}
\qquad
&&\text{for ${\mathcal{I}} \in {\mathfrak{I}}_{n+1,i}$, }
\\
h_{\lambda}^{\mathcal{I}}:=& \iota_{\lambda}^{\ast} h^{\mathcal{I}}
\qquad
&&\text{for ${\mathcal{I}} \in {\mathfrak{I}}_{n+1,i+1}$.  }
\end{alignat*}
By Lemma \ref{lem:KM}
 and Example \ref{ex:KMi}, 
 we have shown the following.

\begin{proposition}
[basic $K$-type $\mub$ and $\mus$]
\label{prop:minKN}
We define linear maps by 
\begin{alignat}{3}
\label{eqn:minlmd}
  \Exterior^i({\mathbb{C}}^{n+1})
  &\to 
  C^{\infty}({\mathbb{R}}^{n}, \Exterior^i({\mathbb{C}}^{n})), 
  \quad
  &&e_{\mathcal{I}} \mapsto {\bf{1}}_{\lambda}^{\mathcal{I}}
  \quad
  &&\text{for ${\mathcal{I}} \in {\mathfrak{I}}_{n+1,i}$}, 
\\
\Exterior^{i+1}({\mathbb{C}}^{n+1})
&\to 
C^{\infty}
({\mathbb{R}}^n,  
 \Exterior^{i}({\mathbb{C}}^{n}) 
), 
\quad
&&e_{\mathcal{I}} \mapsto 
h_{\lambda}^{\mathcal{I}}
\quad
  &&\text{for ${\mathcal{I}} \in {\mathfrak{I}}_{n+1,i+1}$}.  
\notag
\end{alignat}
Then, 
 for $\delta=\pm$, 
 the images give 
 the unique $K$-types 
$
   \mub(i,\delta)=\Exterior^i({\mathbb{C}}^{n+1}) \boxtimes \delta
$
 and 
$\mus(i,\delta)=\Exterior^{i+1}({\mathbb{C}}^{n+1}) \boxtimes (-\delta)$
 respectively, 
 of the principal series representation 
$I_{\delta}(i,\lambda)
=\operatorname{Ind}_P^G
  (\Exterior^i
  ({\mathbb{C}}^{n}) \otimes \delta \otimes {\mathbb{C}}_{\lambda})$
 of $G$ in the $N$-picture.  
\end{proposition}
An explicit formula 
 for ${\bf{1}}_{\lambda}^{\mathcal{I}}$
 and $h_{\lambda}^{\mathcal{I}}$
 is given as follows.  
\begin{lemma}
\label{lem:NKtype}
Let $S_{{\mathcal{I}} {\mathcal{J}}}(b)$
 be the quadratic polynomial 
 of $b=(b_1, \cdots, b_n)$ defined in \eqref{eqn:SIJ}.  
\begin{enumerate}
\item[{\rm{(1)}}]
Let $0 \le i \le n$.  
For ${\mathcal{I}} \in {\mathfrak{I}}_{n+1,i}$
 and $\lambda \in {\mathbb{C}}$, 
 we have 
\begin{align}
  {\bf{1}}_{\lambda}^{\mathcal{I}}(b)
  =& (1+|b|^2)^{-\lambda}
      \sum_{J \in {\mathfrak{I}}_{n,i}}
      \det \psi_{n+1}(1,b)_{{\mathcal{I}} J} e_J
\notag
\\
  =&-(1+|b|^2)^{-\lambda-1}
  \sum_{J \in {\mathfrak{I}}_{n,i}}
  S_{{\mathcal{I}}J}(1,b) e_J.  
\label{eqn:minI}
\end{align}

If $i=0$, 
 we regard ${\mathcal{I}}=\emptyset$
 and ${\mathbf{1}}_{\lambda}^{\emptyset} =(1+|b|^2)^{-\lambda}$
 (see \eqref{eqn:Sempty}).  
\item[{\rm{(2)}}]
Let $0 \le i \le n$.  
For ${\mathcal{I}} \in {\mathfrak{I}}_{n+1,i+1}$
 and $\lambda \in {\mathbb{C}}$, 
 we have 
\begin{align}
h_{\lambda}^{\mathcal{I}}(b)
=&
-(1+|b|^2)^{-\lambda}
      \sum_{J \in {\mathfrak{I}}_{n,i}}
      \det \psi_{n+1}(1,b)_{{\mathcal{I}}, J\cup\{0\}} e_J
\notag
\\
=&
(1+|b|^2)^{-\lambda-1}
\sum_{J \in {\mathfrak{I}}_{n,i}}
S_{{\mathcal{I}},J \cup \{0\}}(1,b)e_J.  
\label{eqn:hIlmd}
\end{align}
\end{enumerate}
\end{lemma}

We note that Lemma \ref{lem:NKtype} implies 
\begin{align}
\label{eqn:I10}
{\mathbf{1}}_{\lambda}^{\mathcal{I}}(0)
&=
 \begin{cases}
 e_{\mathcal{I}} \quad &\text{0 $\not\in {\mathcal{I}}$}, 
\\
 0 \quad &\text{0 $\in {\mathcal{I}}$}, 
 \end{cases} 
\\
\label{eqn:hI0}
h_{\lambda}^{\mathcal{I}}(0)
&=
\begin{cases}
e_{{\mathcal{I}} \setminus \{0\}}\qquad &0 \in {\mathcal{I}}, 
\\
0
\qquad &0 \notin {\mathcal{I}}.  
\end{cases}
\end{align}

\begin{proof}
[Proof of Lemma \ref{lem:NKtype}]
Suppose $b \in {\mathbb{R}}^n$, 
 and let $k(b) \in SO(n+1)$ be as defined 
 in \eqref{eqn:kb}.  
By \eqref{eqn:fI} and \eqref{eqn:hJ}, 
 respectively,
 the formula \eqref{eqn:fKN} of $\iota_{\lambda}^{\ast}$ tells
 that 
\index{A}{01oneIlmd@${\bf{1}}_{\lambda}^{\mathcal{I}}$|textbf}
\index{A}{hIlmd@$h_{\lambda}^{{\mathcal{I}}}$|textbf}

\begin{alignat*}{2}
{\bf{1}}_{\lambda}^{\mathcal{I}}(b)
&=(\iota_{\lambda}^{\ast} {\bf{1}}^{\mathcal{I}})(b)
&&=(1+|b|^2)^{-\lambda}{\bf{1}}^{\mathcal{I}}(k(b))
\\
&
&&=(1+|b|^2)^{-\lambda} \sum_{J \in {\mathfrak{I}}_{n,i}}
(\det k(b))_{{\mathcal{I}} J}e_J, 
\\
h_{\lambda}^{\mathcal{I}}(b)
&=(\iota_{\lambda}^{\ast} h^{\mathcal{I}})(b)
&&= (1+|b|^2)^{-\lambda} h^{\mathcal{I}}(k(b))
\\
&
&&= (1+|b|^2)^{-\lambda} \sum_{J \in {\mathfrak{I}}_{n,i}}
   (\det k(b))_{{\mathcal{I}}, J \cup \{0\}} e_J.  
\end{alignat*}

It follows from Lemma \ref{lem:psidet} (2)
 that, for ${\mathcal{I}},{\mathcal{J}} \subset \{0,1,\cdots,n\}$
 with 
 $\# {\mathcal{I}} = \# {\mathcal{J}}=i$, 
 the minor determinant
 of 
\index{A}{kb@$k(b)$}
$
k(b)
$ is given by 
\begin{equation}
\label{eqn:kbdet}
(\det k(b))_{{\mathcal{I}} {\mathcal{J}}}
=
-\varepsilon_{\mathcal{J}}(0) 
(\det \psi_{n+1}(1,b))_{{\mathcal{I}} {\mathcal{J}}} 
=
 \varepsilon_{\mathcal{J}} (0) \frac {S_{{\mathcal{I}} {\mathcal{J}}} (1,b)}{1 + |b|^2}, 
\end{equation}
where we set $\varepsilon_{\mathcal{J}} (0)=-1$
 for $0 \notin {\mathcal{J}}$
 and $\varepsilon_{\mathcal{J}} (0)=1$ for $0 \in {\mathcal{J}}$.

Now the second formul{\ae} in Lemma \ref{lem:NKtype} 
 are also shown.  
\end{proof}

\subsection{Knapp--Stein intertwining operator}
\label{sec:KS}

In this section
 we summarize some basic results
 on the matrix-valued Knapp--Stein intertwining operators, 
 see \cite{KS, KS2}.  
In the general framework of symmetry breaking operators
 for the restriction $G \downarrow G'$, 
 this classical case may be thought of as a special case
 where $G=G'$, 
 and the proof is much easier than the general case 
 $G \supsetneqq G'$.  
Nevertheless, 
 we sketch a proof 
 of results which we need
 in other chapters.

\subsubsection{Knapp--Stein intertwining operator}
\label{subsec:KerKS}

For $(\sigma,V)\in \widehat{O(n)}$, 
 $\delta, \varepsilon \in \{ \pm \}$
 and $\lambda,\nu \in {\mathbb{C}}$, 
 we consider intertwining operators between two principal series representations 
 $I_{\delta}(V,\lambda)$
 and $I_{\varepsilon}(V,\nu)$
 of $G=O(n+1,1)$.  
They are determined by distribution kernels,
 and Fact \ref{fact:kernel}
 (see \cite[Prop.~3.2]{sbon}) with $G=G'$ and $V=W$ 
 gives a linear isomorphism
\begin{equation}
\label{eqn:GGkernel}
  {\operatorname{Hom}}_G(I_{\delta}(V,\lambda), I_{\varepsilon}(V,\nu))
  \simeq
  ({\mathcal{D}}'(G/P, {\mathcal{V}}_{\lambda, \delta}^{\ast})
  \otimes 
  V_{\nu, \varepsilon})^{\Delta(P)}, 
\end{equation}
 where $P$ acts diagonally
 on the $(G \times P)$-module
 ${\mathcal{D}}'(G/P, {\mathcal{V}}_{\lambda,\delta}^{\ast})
  \otimes V_{\nu,\varepsilon}$.  
As in Proposition \ref{prop:Tpair} (2), 
 the restriction to the open Bruhat cell 
 determines invariant distributions 
 in the right-hand side, 
 and thus we have an injective homomorphism
\[
  ({\mathcal{D}}'(G/P, {\mathcal{V}}_{\lambda,\delta}^{\ast})
   \otimes 
   V_{\nu,\varepsilon})^{\Delta(P)}
  \hookrightarrow
  {\mathcal{D}}'({\mathbb{R}}^n)\otimes {\operatorname{End}}_{\mathbb{C}}(V), 
  \quad
  f \mapsto F(x):=f(n_-(x)), 
\]
where we have used the canonical isomorphism
 $V^{\vee} \otimes V \simeq {\operatorname{End}}_{\mathbb{C}}(V)$.  
Different from the case $G \supsetneqq G'$
 for symmetry breaking operators, 
 there are strong constraints
 on the parameter
 for the existence of nonzero elements
 in \eqref{eqn:GGkernel}.  
In fact,
 it follows readily from the $P$-invariance
 that $F|_{{\mathbb{R}}^n \setminus \{0\}}$ is nonzero
 only if $\nu=n-\lambda$, 
 and in this case it is proportional
 to $|x|^{2 \lambda-2n} \sigma(\psi_n(x))$, 
 where we recall from \eqref{eqn:psim} the definition
 of $\psi_n \colon {\mathbb{R}}^n \setminus \{0\} \to O(n)$.  
We normalize as 
\begin{equation}
\label{eqn:KT}
\Ttcal{\lambda}{n-\lambda}V(x)
:=
 \frac{1}{\Gamma(\lambda-\frac n 2)}
|x|^{2\lambda-2n} \sigma(\psi_n(x)).  
\end{equation}

\begin{remark}
The normalization of the Knapp--Stein operator
 is not unique,
 and different choices are useful
 for different purposes.  
See for example Knapp--Stein \cite{KS} or Langlands \cite{LLNM}.  
\end{remark}
With the normalization \eqref{eqn:KT},
 we now review the Knapp--Stein intertwining operators
 in this setting as follows.

\begin{lemma}
[normalized Knapp--Stein operator]
\label{lem:nKS}
The distribution \eqref{eqn:KT} belongs
 to $L_{\operatorname{loc}}^1({\mathbb{R}}^n) \otimes {\operatorname{End}}_{\mathbb{C}}(V)$
 if $\operatorname{Re} \lambda \gg 0$, 
 and extends to an element of 
\index{A}{Vlndast@${\mathcal{V}}_{\lambda,\delta}^{\ast}$, dualizing bundle}
 $({\mathcal{D}}'(G/P, {\mathcal{V}}_{\lambda,\delta}^{\ast}) \otimes V_{n-\lambda, \delta})^{\Delta(P)}$.  
Furthermore,
 it has an analytic continuation
 to the entire $\lambda \in {\mathbb{C}}$.  
\end{lemma}
By definition, 
 the (normalized) 
\index{B}{KnappSteinoperator@Knapp--Stein operator|textbf}
Knapp--Stein intertwining operator
\index{A}{IdeltaV@$I_{\delta}(V, \lambda)$}
\index{A}{TVln@$\Ttbb \lambda{n-\lambda}{V}$, normalized Knapp--Stein intertwining operator|textbf}
\begin{equation}
\label{eqn:KnappStein}
  \Ttbb{\lambda}{n-\lambda} V
  :
  I_{\delta}(V,\lambda) \to I_{\delta}(V,n-\lambda) 
\end{equation}
is defined 
 in the $N$-picture
 of the principal series representation 
 by the formula
\[
(\Ttbb{\lambda}{n-\lambda} V f)(x)
 = \int_{{\mathbb{R}}^n} \Ttcal{\lambda}{n-\lambda} V (x-y) f(y) dy.  
\]

When $(\sigma,V)$ is the $i$-th exterior representation
 on $\Exterior^i({\mathbb{C}}^n)$,
 we write simply 
 $\Ttbb \lambda{n-\lambda}{i}$ and $\Ttcal \lambda{n-\lambda}{i}$
 for 
 the operator $\Ttbb \lambda{n-\lambda}{V}$
 and the distribution $\Ttcal \lambda{n-\lambda}{V}$, 
respectively.

The Knapp--Stein operator \eqref{eqn:KnappStein} gives
 a continuous $G$-homomorphism 
 $I_{\delta}(i,\lambda) \to I_{\varepsilon}(j,\nu)$
  when $j=i$ 
 (and $\delta =\varepsilon$, $\nu=n-\lambda$).  
On the other hand,
 there exist $G$-intertwining operators
 $I_{\delta}(i,\lambda) \to I_{\varepsilon}(j,\nu)$
 also when $i \ne j$
 for special parameters.  
Like sporadic symmetry breaking operators
 ({\it{cf}}. Theorem \ref{thm:152347}),
 they are given by differential operators
 as follows.  
\begin{fact}
\label{fact:Di}
Suppose that $0 \le i \le n-1$.  
\begin{enumerate}
\item[{\rm{(1)}}]
We can identify $I_{(-1)^i}(i,i)$
 with the space ${\mathcal{E}}^i(S^n)$
 of differential $i$-forms
 endowed with the natural action of the conformal group 
 $G=O(n+1,1)$. 
\item[{\rm{(2)}}]
The exterior derivative
 $d \colon {\mathcal{E}}^i(S^n) \to {\mathcal{E}}^{i+1}(S^n)$
 induces a $G$-intertwining operator
\[
   D_i \colon I_{(-1)^i}(i,i) \rightarrow I_{(-1)^{i+1}}(i+1,i+1).  
\]
The kernel of $D_i$ is $\Pi_{i,(-1)^{i}}$, 
 and the image is $\Pi_{i+1,(-1)^{i+1}}$. 
\end{enumerate}
\end{fact}

This follows from \cite[Thm.~12.2]{KKP}.  
We note that the existence
 of such an intertwining operator is assured
 {\it{a priori}} by the composition series
 of the principal series representation
 (Theorem \ref{thm:LNM20}), 
 see also \cite{C}.

\subsubsection{$K$-spectrum
 of the Knapp--Stein intertwining operator}
\label{subsec:KSspec}
\index{B}{Kspectrum@$K$-spectrum}

This section gives an explicit formula for the eigenvalues
 of the (normalized) Knapp--Stein intertwining operator
\index{A}{Tiln@$\Ttbb \lambda{n-\lambda}{i}$|textbf}
\begin{equation}
\label{eqn:KSii}
\Ttbb \lambda {n-\lambda}i :I_{\delta}(i,\lambda) \to I_{\delta}(i,n-\lambda)
\end{equation}
 on the basic $K$-types
 $\mub(i,\delta)$
 and $\mus(i,\delta)$
 (see \eqref{eqn:muflat} and \eqref{eqn:musharp}, 
respectively).  
For $0 \le i \le n$ and $\lambda \in {\mathbb{C}}$, 
 we set 
\begin{equation}
\label{eqn:specT}
  c^{\natural}(i,\lambda)
  =
  \frac{\pi^{\frac n2}}{\Gamma(\lambda+1)}
  \times
  \begin{cases}
  \lambda-i
  \qquad
  &\text{if $\natural = \flat$}, 
\\
  n-i-\lambda
  \qquad
  &\text{if $\natural = \sharp$}.  
  \end{cases}
\end{equation}

\begin{proposition}
\label{prop:TminK}
Suppose $0 \le i \le n$, 
 $\lambda \in {\mathbb{C}}$
 and $\delta \in \{ \pm \}$.  
Then the (normalized) Knapp--Stein intertwining operator
\index{A}{Tiln@$\Ttbb \lambda{n-\lambda}{i}$}
\index{A}{Ideltai@$I_{\delta}(i, \lambda)$}
\[
\Ttbb{\lambda}{n-\lambda}i
:
I_{\delta}(i,\lambda) \to I_{\delta}(i,n-\lambda)
\]
acts on the 
\index{B}{basicKtype@basic $K$-type}
basic $K$-types
\index{A}{0muflat@$\mub(i,\delta)$}
 $\mub (i,\delta)=\Exterior^i({\mathbb{C}}^{n+1}) \boxtimes \delta$
 and 
\index{A}{0musharp@$\mus(i,\delta)$}
 $\mus(i,\delta)=\Exterior^{i+1}({\mathbb{C}}^{n+1}) \boxtimes (-\delta)$
 as the scalar multiplication:
\[
  \Ttbb \lambda {n-\lambda}i  \circ \iota_{\lambda}^{\ast}
  =
  c^{\natural}(i,\lambda)\iota_{n-\lambda}^{\ast}
\quad
  \text{on $\mu^{\natural}(i,\delta)$
        for $\natural = \flat$ or $\sharp$}.  
\]
In other words, 
we have 
\index{A}{01oneIlmd@${\bf{1}}_{\lambda}^{\mathcal{I}}$}
\index{A}{hIlmd@$h_{\lambda}^{{\mathcal{I}}}$}
\begin{alignat*}{2}
\Ttbb \lambda {n-\lambda} i({\bf{1}}_{\lambda}^{\mathcal{I}})
=&
\frac{(\lambda-i)\pi^{\frac n 2}}{\Gamma(\lambda+1)}
{\bf{1}}_{n-\lambda}^{\mathcal{I}}  
\quad
&&\text{for all ${\mathcal{I}} \in {\mathfrak {I}}_{n+1,i}$}, 
\\
\Ttbb \lambda {n-\lambda}i
(h_{\lambda}^{\mathcal{I}})
=&
\frac{(n-i-\lambda) \pi^{\frac n 2}}{\Gamma(\lambda+1)}
h_{n-\lambda}^{\mathcal{I}}
\qquad
&&\text{for all ${\mathcal{I}} \in {\mathfrak{I}}_{n+1,i+1}$}.  
\end{alignat*}
\end{proposition}

\begin{remark}
Proposition \ref{prop:TminK}
 in the $i=0$ case for $\mub (i,\delta)$
 was proved in \cite[Prop.~4.6]{sbon}.  
\end{remark}
We will give a proof of Proposition \ref{prop:TminK}
 in Section \ref{subsec:Tspec}.  

We recall from Theorem \ref{thm:LNM20}
 that the composition series of $I_{\delta}(i,i)$ and $I_{\delta}(i,n-i)$
 are described by the following exact sequences of $G$-modules:
\begin{align*}
&0 \to \Pi_{i,\delta} \to I_{\delta}(i,i) \to \Pi_{i+1,-\delta} \to 0, 
\\
&0 \to \Pi_{i+1,-\delta} \to I_{\delta}(i,n-i) \to \Pi_{i,\delta} \to 0, 
\end{align*}
which do not split if $i \ne \frac n 2$.  
Thus Proposition \ref{prop:TminK} implies:
\begin{proposition}
\label{prop:Timage}
Suppose $G=O(n+1,1)$ and $i \ne \frac n 2$.  
Then the kernels and the images
 of the $G$-homomorphisms 
 $\Ttbb \lambda{n-\lambda}i \colon I_{\delta}(i,\lambda) \to I_{\delta}(i,n-\lambda)$
 for $\lambda=i$, 
 $n-i$ are given by  
\begin{alignat*}{2}
{\operatorname{Ker}}(\Ttbb i{n-i}i) 
&\simeq \Pi_{i,\delta}
&& \simeq {\operatorname{Image}}(\Ttbb {n-i}i i)
\\
{\operatorname{Image}}(\Ttbb i{n-i}i) 
&\simeq \Pi_{i+1,-\delta}
&&\simeq {\operatorname{Ker}}(\Ttbb {n-i}i i).  
\end{alignat*}
\end{proposition}

\subsubsection{Vanishing of the Knapp--Stein operator}

There are a few exceptional parameters
 $(i,\lambda)$ 
 for which $\Ttbb \lambda {n-\lambda} i$ vanishes:

\begin{proposition}
\label{prop:Tvanish}
Suppose $G=O(n+1,1)$, 
 $0 \le i \le n$, 
 and $\lambda \in {\mathbb{C}}$.  
Then the normalized Knapp--Stein intertwining operator 
$\Ttbb \lambda {n-\lambda}i$
 is zero
 if and only if 
 $\lambda=i=\frac n 2$.  
\end{proposition}

\begin{proof}
See \cite{xkresidue}.  
\end{proof}

A renormalization of the Knapp--Stein intertwining operator
 $\Ttbb \lambda {n-\lambda} i$
 for $n=2i$ will be discussed in Section \ref{subsec:renormT}.  

\subsubsection{Integration formula for the $(K,K)$-spectrum}
\label{subsec:Tspec}
In this subsection,
 we give a proof of Proposition \ref{prop:TminK}.  
Let $\natural = \flat$ or $\sharp$.  
Since the multiplicity
 of the $K$-type
 $\mu^{\natural}(i,\delta)$
 in the principal series representation $I_{\delta}(i,\lambda)$
 is one,
 there exists a constant $c^{\natural}(i,\lambda)$
 depending on $i$ and $\lambda$
 such that
\begin{equation}
\label{eqn:TcT}
\Ttbb \lambda{n-\lambda}i \circ \iota_{\lambda}^{\ast}
=
c^{\natural}(i,\lambda)
\iota_{n-\lambda}^{\ast}
\quad
\text{on $\mu^{\natural}(i,\delta)$}.  
\end{equation}
We shall show
 that the constants $c^{\natural}(i,\lambda)$
 in the equation \eqref{eqn:TcT}
 are given by the formul{\ae} \eqref{eqn:specT}.  
The first step is to give an integral formula 
 for the constants 
\index{A}{cyflat@$c^{\flat}(i,\lambda)$}
\index{A}{cyfsharp@$c^{\sharp}(i,\lambda)$}
$
c^{\natural}(i,\lambda)
$
for $\natural = \flat$ and $\sharp$:
\begin{lemma}
\label{lem:161783}
Suppose $0 \le i \le n$
 and $\lambda \in {\mathbb{C}}$
 with $\operatorname{Re} \lambda \gg 0$. 
Then we have  
\begin{align*}
  c^{\flat}(i,\lambda)
  &=
  \frac{1}{\Gamma(\lambda-\frac n 2)}
  \int_{{\mathbb{R}}^n} |b|^{2\lambda-2n}
                        (1+|b|^2)^{-\lambda} 
  \left(1-\frac{2}{|b|^2(1+|b|^2)}\sum_{k=1}^{i} b_k^2\right) d b,
\\
  c^{\flat}(i,\lambda) - c^{\sharp}(i,\lambda)
  &=
  \frac{2}{\Gamma(\lambda-\frac n 2)}
  \int_{{\mathbb{R}}^n} |b|^{2\lambda-2n+2}
                        (1+|b|^2)^{-\lambda-1} d b.  
\end{align*}
\end{lemma}

\begin{proof}
[Proof of Lemma \ref{lem:161783}]
We first consider \eqref{eqn:TcT} 
 for $\natural = \flat$.  
Then we have 
\[
  \Ttbb \lambda {n-\lambda} i ({\bf{1}}_{\lambda}^{\mathcal{I}})
= c^{\flat}(i,\lambda) {\bf{1}}_{n-\lambda}^{\mathcal{I}}
\qquad
\text{for all ${\mathcal{I}} \in {\mathfrak{I}}_{n+1,i}$.  }
\]
Take ${\mathcal{I}} \in {\mathfrak{I}}_{n+1,i}$
 such that $0 \not\in {\mathcal{I}}$.  
Then \eqref{eqn:I10} tells that 
\begin{equation}
\label{eqn:T1c}
(\Ttbb \lambda{n-\lambda}i {\bf{1}}_{\lambda}^{\mathcal{I}})
(0)
=
c^{\flat}(i,\lambda)e_{\mathcal{I}}.  
\end{equation}

Let us compute the left-hand side.  
In view of the distribution kernel \eqref{eqn:KT}
 of the normalized Knapp--Stein operator $\Ttbb \lambda {n-\lambda}i$, 
 for $\operatorname{Re}\lambda \gg 0$
 we have
\begin{equation*}
(\Ttbb \lambda {n-\lambda} i {\bf{1}}_{\lambda}^{\mathcal{I}})(0)
= \frac{1}{\Gamma(\lambda - \frac n 2)} \int_{\mathbb{R}^n}
   |-b|^{2 \lambda-2n} \sigma(\psi_n(-b)){\mathbf{1}}_{\lambda}^{\mathcal{I}}(b)db.  
\end{equation*}
By \eqref{eqn:exrep} and the formula \eqref{eqn:minI}
 of ${\bf{1}}_{\lambda}^{\mathcal{I}}(b)$, 
 the integrand amounts to 
\begin{equation*}
\sum_{J,J' \in {\mathfrak{I}}_{n,i}}
   |b|^{2 \lambda-2n} (1+|b|^2)^{-\lambda}
   (\det \psi_{n+1}(1,b))_{{\mathcal{I}} J}
   (\det \psi_n(b))_{J' J} e_{J'}.  
\end{equation*}
Comparing the coefficients
 of $e_{\mathcal{I}}$ in the both sides of \eqref{eqn:T1c}, 
 we get 
\[
   c^{\flat}(i,\lambda)
   =
   \frac{1}{\Gamma(\lambda-\frac n 2)}
   \int_{\mathbb{R}^n}
   |b|^{2 \lambda-2n} (1+|b|^2)^{-\lambda}
   g_{\mathcal{I}}(b)d b, 
\]
where we set 
\begin{equation}
   g_{\mathcal{I}}(b):=
   \sum_{J \in {\mathfrak{I}}_{n,i}}
   (\det \psi_{n+1} (1,b))_{{\mathcal{I}} J} 
   (\det \psi_n (b))_{{\mathcal{I}} J}
   =1-\frac{2 Q_{\mathcal{I}}(b)}{(1+|b|^2) |b|^2}.  
\label{eqn:Timinor}
\end{equation}
The second equality was proved as the minor summation formula
 in Proposition \ref{prop:msum} (4), 
 where we recall $Q_{\mathcal{I}}(b)=\sum_{l \in {\mathcal{I}}}b_l^2$.  
Therefore, 
 by taking ${\mathcal{I}}=\{1,2, \cdots, n\}$, 
 we get the first assertion
 of Lemma \ref{lem:161783}.

Next,
 we consider \eqref{eqn:TcT} for $\natural =\sharp$.  
Then we have 
\[
  \Ttbb {\lambda}{n-\lambda}i (h_{\lambda}^{\mathcal{I}})
  =
  c^{\sharp}(i,\lambda) h_{n-\lambda}^{\mathcal{I}}
\quad
\text{for all ${\mathcal{I}} \in {\mathfrak {I}}_{n+1,i+1}$.}
\]  
Take $I \in {\mathfrak {I}}_{n,i}$, 
 and set ${\mathcal{I}}:=I \cup \{0\} \in {\mathfrak {I}}_{n+1,i+1}$.  
By \eqref{eqn:hI0}, 
 we have
\begin{equation}
\label{eqn:Thc}
\Ttbb {\lambda}{n-\lambda} i (h_{\lambda}^{\mathcal{I}})(0)
  =
  c^{\sharp}(i,\lambda) e_I.  
\end{equation}
By \eqref{eqn:KT}, 
 we have
\begin{equation*}
\Ttbb {\lambda}{n-\lambda}i (h_{\lambda}^{\mathcal{I}})(0)
  = \frac{1}{\Gamma(\lambda-\frac n 2)}
     \int_{\mathbb{R}^n} |-b|^{2 \lambda-2n} \sigma(\psi_n(-b)) h_{\lambda}^{\mathcal{I}} (b) d b.  
\end{equation*}
Comparing the coefficients of $e_{I}$ in the both sides of the equation \eqref{eqn:Thc}, 
 we get from \eqref{eqn:hIlmd} and \eqref{eqn:exrep}
\[
  c^{\sharp}(i,\lambda)
 =- \frac{1}{\Gamma(\lambda-\frac n 2)}
    \int_{\mathbb{R}^n} |b|^{2\lambda-2n} (1+|b|^2)^{-\lambda}
    g_I'(b)
        d b, 
\]
 where we set 
\[
g_{I}'(b):=\sum_{J \in {\mathfrak{I}}_{n,i}}
(\det \psi_{n+1}(1,b))_{I \cup \{0\}, J\cup \{0\}}
       (\det \psi_{n}(b))_{I J}.  
\]

We note 
 that 
\[
  \det \psi_{n+1}(1,b)_{I J}
  =
  \det \psi_{n}(b;\frac{2}{1+|b|^2})_{I J}
\]
if $I, J \in {\mathfrak{I}}_{n,i}$ is regarded as elements
 of ${\mathfrak{I}}_{n+1,i}$
 in the left-hand side.  
Then we have
\[
   g_I(b)+g_I'(b)=\frac{2|b|^2}{1+|b|^2}
\]
{}from Proposition \ref{prop:msum} (3), 
 and thus we get 
\[
 c^{\flat}(i,\lambda)-c^{\sharp}(i,\lambda)
 =\frac{2}{\Gamma(\lambda-\frac n 2)} 
   \int_{\mathbb{R}^n} |b|^{2\lambda+2-2n} (1+|b|^2)^{-\lambda-1} d b.  
\]
Now Lemma \ref{lem:161783} is proved.  
\end{proof}
The second step is to compute the integrals in Lemma \ref{lem:161783}.  
\begin{lemma}
\label{lem:cTcompute}
For $\operatorname{Re}\lambda \gg 0$, 
$c^{\flat}(i,\lambda)$ and $c^{\sharp}(i,\lambda)$ take 
the form \eqref{eqn:specT}.  
\end{lemma}
\begin{proof}
Let $B(\lambda,\nu)$ denote the Beta function.  
By the change of variables
 $r^2= \frac{x}{1-x}$, 
 we have
\begin{equation}
\label{eqn:Rabnew}
  \int_0^{\infty}
  r^a (1+r^2)^b d r 
  = 
  \frac 1 2
  \int_0^1 x^{\lambda-1} (1-x)^{\nu-1} d x
  = \frac 1 2 B(\lambda,\nu), 
\end{equation}
where $a = 2 \lambda-1$ and $b=-\lambda-\nu$.  
Then Lemma \ref{lem:161783}
 in the polar coordinates tells that 
\[
c^{\flat}(i,\lambda)
=
\frac{\operatorname{vol}(S^{n-1})}{2\Gamma(\lambda-\frac n2)}
(B(\lambda-\frac n2,\frac n 2) 
- \frac{2i}{n}B(\lambda- \frac n 2,\frac n 2+1))
\]
by \eqref{eqn:Rabnew} and by the following observation:
\[
   \int_{S^{n-1}} |\omega_i|^2 d \omega
   =
  \frac 1 n \operatorname{vol}(S^{n-1})
\quad
  (1 \le i \le n).  
\]
Since
$
\operatorname{vol}(S^{n-1})
=\frac{2 \pi^{\frac n 2}}{\Gamma(\frac n 2)}, 
$
we get the first statement.

By the second formula of Lemma \ref{lem:161783}, 
 we have 
\begin{align*}
 c^{\flat}(i,\lambda)-c^{\sharp}(i,\lambda)
 =&\frac{1}{\Gamma(\lambda-\frac n 2)}
   {\operatorname{vol}}(S^{n-1}) B(\lambda-\frac n2, \frac n2 +1)
\\
=&\frac{(2\lambda-n) \pi^{\frac n 2}}
       {\Gamma(\lambda+1)}.  
\end{align*}
Thus the closed formula \eqref{eqn:specT}
 for $c^{\sharp}(i,\lambda)$ is also proved.  
\end{proof}

\begin{proof}
[Proof of Proposition \ref{prop:TminK}]
The assertion follows from Lemmas \ref{lem:161783} and \ref{lem:cTcompute}
 for $\operatorname{Re}\lambda\gg 0$.  
For general $\lambda \in {\mathbb{C}}$, 
Proposition holds by the analytic continuation.  
\end{proof}

\subsection{Renormalization of the Knapp--Stein intertwining operator}
\label{subsec:renormT}
Because of the vanishing of the normalized Knapp--Stein intertwining operators
 in the middle degree 
 when $n$ is even (Proposition \ref{prop:Tvanish}), 
 intertwining operators from $I_{\delta}(\frac n 2, \lambda)$ to $I_{\delta}(\frac n 2, n-\lambda)$ require special attention.  
In this case, 
we set $n=2m$ and renormalize the Knapp--Stein intertwining operator
 of $G=O(2m+1,1)$ at the middle degree by
\index{A}{TVttiln@$\Tttbb \lambda{n-\lambda}{\frac n 2}$, renormalized Knapp--Stein intertwining operator|textbf}
\begin{equation}
\label{eqn:Ttilde}
   \Tttbb \lambda {2m-\lambda}{m}
   :=
   \frac{1}{\lambda-m} \Ttbb \lambda {2m-\lambda}{m}.  
\end{equation}
Then $\Tttbb \lambda {2m-\lambda} {m} \colon I_{\delta}(m, \lambda) \to I_{\delta}(m, 2m-\lambda)$ depends
holomorphically in the entire $\lambda \in {\mathbb{C}}$,
 and is vanishing nowhere.

If $\lambda=m$,
 then $\Tttbb \lambda {2m-\lambda} m$ acts
 as an endomorphism
 of $I_{\delta}(m,m)$.  
On the other hand,
 we know from Theorem \ref{thm:LNM20} (1)
 that the principal series representation $I_{\delta}(m, m)$
 decomposes into the direct sum
 of two irreducible 
\index{B}{temperedrep@tempered representation}
 tempered representations
 of $G$ as follows:
\[
   I_{\delta}(m, m)
   \simeq 
   I_{\delta}(m)^{\flat} \oplus I_{\delta}(m)^{\sharp}
   \equiv
   \Pi_{m, \delta} \oplus \Pi_{m+1, -\delta}.  
\]
\begin{lemma}
\label{lem:161745}
Let $n=2m$ and $G=O(2m+1,1)$.  
Then the 
\index{B}{KnappSteinoperatorrenorm@ Knapp--Stein operator, renormalized---|textbf}
 renormalized Knapp--Stein operator
 $\Tttbb {m}{m}{m}$
 acts on $I_{\delta}(m, m)\simeq \Pi_{m,\delta} \oplus \Pi_{m+1,-\delta}$
 as
\[
   \frac{\pi^{m}}{m!}
   ({\operatorname{id}}_{\Pi_{m, \delta}} 
    \oplus 
    (-{\operatorname{id}})_{\Pi_{m+1, -\delta}} ).  
\]
\end{lemma}
\begin{proof}
Since the irreducible $G$-module $\Pi_{m, \delta}$ is not isomorphic to 
the irreducible $G$-module $\Pi_{m+1, -\delta}$, 
 the renormalized Knapp--Stein intertwining operator $\Tttbb {m}{m}{m}$ acts on each irreducible summand
 by scalar multiplication.  
Therefore, 
 it is sufficient to find the scalars
 on specific $K$-types occurring 
 in each summand.  
By Proposition \ref{prop:TminK}, 
 the renormalized Knapp--Stein intertwining operator
 $\Tttbb {\lambda}{2m-\lambda}{m}$ acts on vectors
 that belong to the $K$-types
 $\mub (m, \delta)(\subset \Pi_{m,\delta})$
 and $\mus(m, \delta)(\subset \Pi_{m+1,-\delta})$
 by the scalars
\[
   \text{
   $\frac{1}{\lambda-m} \frac{(\lambda-m)\pi^{m}}{\Gamma(\lambda+1)}$
\quad
   and 
\quad
   $\frac{1}{\lambda-m} \frac{(2m-m- \lambda)\pi^m}{\Gamma(\lambda+1)}$, 
}
\]
respectively.  
Taking the limit as $\lambda$ tends to $m$, 
 we get the lemma.  
\end{proof}

\subsection{Kernel of the Knapp--Stein operator}
\label{subsec:KerKnSt}
In this section,
 we discuss the proper submodules
 of the principal series representation
 $I_{\delta}(i,\lambda)$ of $G=O(n+1,1)$
 at reducible points
 (see \eqref{eqn:redIilmd}).  

We consider the composition of the Knapp--Stein operators,
 $\Ttbb {n-\lambda} \lambda i \circ \Ttbb \lambda{n-\lambda} i
\in {\operatorname{End}}_{G} (I_{\delta}(i,\lambda))$.  
By Proposition \ref{prop:TminK}, 
 its $K$-spectrum on the basic $K$-type
 $\mub (i,\delta)$ is given as 
\[
   \Ttbb {n-\lambda} \lambda i \circ \Ttbb \lambda{n-\lambda} i
   ({\bf{1}}_{\lambda}^{\mathcal{I}})
   =
   \frac{(\lambda-i)(n-\lambda-i) \pi^n}
        {\Gamma(\lambda+1)\Gamma(n-\lambda+1)}
        ({\bf{1}}_{\lambda}^{\mathcal{I}})
\quad
\text{for all ${\mathcal{I}} \in {\mathfrak{I}}_{n+1,i}$}.  
\]
Since the principal series representation $I_{\delta}(i,\lambda)$
 is generically irreducible,
 we conclude
\begin{equation}
\label{eqn:TTc}
\Ttbb {n-\lambda} \lambda i \circ \Ttbb \lambda{n-\lambda} i
  = 
  \frac{(\lambda-i)(n-\lambda-i) \pi^n}
        {\Gamma(\lambda+1)\Gamma(n-\lambda+1)}
  {\operatorname{id}}
\quad
\text{on $I_{\delta}(i,\lambda)$}
\end{equation}
 for generic $\lambda$ by Schur's lemma, 
 and then for all $\lambda \in {\mathbb{C}}$
 by analytic continuation.

\begin{lemma}
\label{lem:reducibleI}
Let $G=O(n+1,1)$, 
 $0 \le i \le n$, 
 and $\delta \in \{\pm\}$.  
Assume
\[
   \lambda \in \{i,n-i\} \cup (-{\mathbb{N}}_+) \cup (n + {\mathbb{N}}_+).  
\]
Then $I_{\delta}(i,\lambda)$ is reducible.  
\end{lemma}
\begin{proof}
If $(n,\lambda) = (2i,i)$, 
 we already know 
 that $I_{\delta}(i,\lambda)$ is reducible,
 see Lemma \ref{lem:161745}.  

Assume now $(n,\lambda) \ne (2i,i)$.  
Then Proposition \ref{prop:Tvanish} tells 
 that neither $\Ttbb {n-\lambda} \lambda i$
 nor $\Ttbb \lambda{n-\lambda} i$ vanishes.  
On the other hand, 
 by \eqref{eqn:TTc}, 
 the assumption on $\lambda$ implies 
\[
   \Ttbb {n-\lambda} \lambda i \circ \Ttbb \lambda{n-\lambda} i =0, 
\]
which shows that at least one of the $G$-modules $I_{\delta}(i,\lambda)$
 or $I_{\delta}(i,n-\lambda)$ is reducible.  
By Lemma \ref{lem:LNM27}, 
 we conclude 
 that both $I_{\delta}(i,\lambda)$
 and its contragredient representation $I_{\delta}(i,n-\lambda)$
 are reducible.  
\end{proof}
Lemma \ref{lem:reducibleI} gives an alternative proof 
 for the \lq\lq{if part}\rq\rq\
 of Proposition \ref{prop:redIilmd} (1).

\begin{proposition}
\label{prop:KerKnSt}
Let $G=O(n+1,1)$, 
 $0 \le i \le n$, 
 $\delta \in \{\pm\}$, 
and $\lambda \in {\mathbb{C}}$.  
Assume further that $I_{\delta}(i,\lambda)$ is reducible,
 namely,
\[
  \lambda \in \{i,n-i\} \cup (-{\mathbb{N}}_+) \cup (n+{\mathbb{N}}_+).  
\]
\begin{enumerate}
\item[{\rm{(1)}}]
Suppose $(n,\lambda) \ne (2i,i)$.  
Then the unique proper submodule
 of $I_{\delta}(i,\lambda)$ is given 
 as the kernel of the Knapp--Stein operator
 $\Ttbb \lambda {n-\lambda} i \colon I_{\delta}(i,\lambda) \to I_{\delta}(i,n-\lambda)$.  
\item[{\rm{(2)}}]
Suppose $(n,\lambda) = (2i,i)$.  
Then $\Ttbb \lambda {n-i} i =0$, 
 and there are two proper submodules
 of $I_{\delta}(i,\lambda)$, 
 which are given 
 as the kernel of
 $\Tttbb i i i \pm \frac{\pi^i}{i!} {\operatorname{id}}
 \in {\operatorname{End}}_{G} (I_{\delta}(i,i))$
 where $\Tttbb i i i$ is the renormalized Knapp--Stein operator.  
\end{enumerate}
\end{proposition}
\begin{proof}
\begin{enumerate}
\item[(1)]
There is a unique irreducible submodule of $I_{\delta}(i,\lambda)$
 for the parameter $\lambda$
 under consideration.  
Hence ${\operatorname{Ker}} (\Ttbb \lambda{n-\lambda}i)$ 
 is the unique irreducible submodule
 by the proof 
 of Lemma \ref{lem:reducibleI}.  
\item[(2)]
This is already proved in Lemma \ref{lem:161745}.  
\end{enumerate}
\end{proof}

\newpage
\section{Regular symmetry breaking operators
 $\Atbb \lambda \nu {\delta\varepsilon} {i,j}$ from $I_{\delta}(i,\lambda)$ to $J_{\varepsilon}(j,\nu)$}
\label{sec:holo}

In this chapter we apply the general results
 developed in Chapter \ref{sec:section7}
 on the analytic continuation
 of integral symmetry breaking operators
 $\Atbb \lambda \nu {\delta\varepsilon} {V,W} 
  \colon 
  I_{\delta}(V,\lambda) \to J_{\varepsilon}(W,\nu)$
 to the special setting where 
\begin{equation}
\label{eqn:VWexterior}
  (V,W)=(\Exterior^i({\mathbb{C}}^n), \Exterior^j({\mathbb{C}}^{n-1})), 
\end{equation}
 and construct a holomorphic family
 of (normalized) regular symmetry breaking operators
\[
  \Atbb \lambda \nu {\delta\varepsilon} {i,j}
  \colon 
   I_{\delta}(i,\lambda) \to J_{\varepsilon}(j,\nu), 
\]
 which exist if and only if $j=i-1$ or $i$
 (Theorems \ref{thm:152271} and \ref{thm:1522a}).  
Then the goal of this chapter 
 is to determine
\begin{enumerate}
\item[$\bullet$]
the parameter $(\lambda,\nu)$ 
for which $\Atbb \lambda \nu {\pm} {i,j}$ vanishes
 (Section \ref{subsec:Aijzero});
\item[$\bullet$]
the $(K,K')$-spectrum of $\Atbb \lambda \nu {\pm} {i,j}$
 (Sections \ref{subsec:Kspec}--\ref{subsec:Apm1});
\item[$\bullet$]
functional equations of $\Atbb \lambda \nu {\pm} {i,j}$
 (Sections \ref{subsec:FI}--\ref{subsec:161702}).  
\end{enumerate}
Thus we will complete the proof of Theorem \ref{thm:161243}
 that determines the zeros
 of the normalized operators
 $\Atbb \lambda \nu {\delta \varepsilon} {i,j}$.  
This is the last missing piece 
 in the classification scheme
 (Theorem \ref{thm:VWSBO}), 
 and thus we complete the proof of the classification
 of the space ${\operatorname{Hom}}_{G'}(I_{\delta}(i,\lambda)|_{G'}, J_{\varepsilon}(j,\nu))$ 
 of {\it{all}} symmetry breaking operators
 as stated in Theorems \ref{thm:1.1} and \ref{thm:SBObasis}.

The $(K,K')$-spectrum resembles eigenvalues
 of a symmetry breaking operator
 (Definition \ref{def:KKspec}), 
 for which we find an integral expression 
 and determine the explicit formula
 for basic $K$- and $K'$-types
 (Theorem \ref{thm:153315}).

The matrix-valued functional equations 
 among various intertwining operators are determined explicitly 
 in Theorems \ref{thm:TAA} and \ref{thm:ATA}
 by using the formula of the $(K,K')$-spectrum, 
 which in turn will play a crucial role in analyzing the behavior
 of the symmetry breaking operators
 at reducible places
 (Chapter \ref{sec:pfSBrho}).

Degenerate cases where the normalized operators $\Atbb \lambda \nu {\pm}{i,j}$ 
 vanish will be discussed in Sections \ref{subsec:161702} and \ref{subsec:Rest}.

As an application of the matrix-valued functional 
equations (Theorems \ref{thm:TAA} and \ref{thm:ATA})
 and the residue formul{\ae}
 of $\Atbb \lambda \nu {\pm}{i,j}$
 (Fact \ref{fact:153316}), 
 we determine when the {\it{differential}}
 symmetry breaking operators
 $\Ctbb \lambda \nu{i,j}$
 ($j=i,i-1$)
 are surjective in Section \ref{subsec:ImageC}.  

\subsection{Regular symmetry breaking operators $\Atbb \lambda \nu {\pm}{i,j}$}
\label{subsec:Aij}
In this section,
 we give the existence condition
 and an explicit construction
 of (generically) regular symmetry breaking operators from 
 $G$-modules $I_{\delta}(V,\lambda)$
 to $G'$-modules $J_{\varepsilon}(W,\nu)$
 in the setting \eqref{eqn:VWexterior}
 by applying the general results
 of Chapters \ref{sec:general} and \ref{sec:section7}, 
 in particular,  
 Theorems \ref{thm:regexist} and \ref{thm:152389}
 and their proofs.  

\subsubsection{Existence condition for regular symmetry breaking operators}
\label{subsec:reg}
We recall from Definition \ref{def:regSBO}
 the notion of {\it{regular}} symmetry breaking operators.~We also recall from \eqref{eqn:Ureg+} and \eqref{eqn:Ureg-}
 the definition of the open dense subsets
 $U_{\pm}^{\operatorname{reg}}$ in ${\mathbb{C}}^2$.  
Then the existence condition of regular symmetry breaking operators
 in the setting \eqref{eqn:VWexterior} is stated as follows.

\begin{theorem}
\label{thm:152271}
Suppose $0 \le i \le n$
 and $0 \le j \le n-1$.  
Then the following three conditions on the pair $(i,j)$ are equivalent:
\begin{enumerate}
\item[{\rm{(i)}}]
 there exists a nonzero regular symmetry breaking operator from the $G$-module
 $I_{\delta}(i,\lambda)$ to the $G'$-module $J_{\varepsilon}(j,\nu)$
 for some $(\lambda,\nu,\delta,\varepsilon) \in {\mathbb{C}}^2 \times \{ \pm \}^2$;
\item[{\rm{(ii)}}]
 for any $(\delta,\varepsilon) \in \{ \pm \}^2$, 
 there exists a nonzero regular symmetry breaking operator from 
 $I_{\delta}(i,\lambda)$ to $J_{\varepsilon}(j,\nu)$
 for all $(\lambda,\nu) \in U_{\delta \varepsilon}^{{\operatorname{reg}}}$;
\item[{\rm{(iii)}}]
 $j=i$ or $i-1$.  
\end{enumerate}
\end{theorem}

\begin{proof}
As we have seen in the decomposition \eqref{eqn:extij},  
 $[V:W] \ne 0$ in the setting \eqref{eqn:VWexterior}
 if and only if $j=i-1$ or $i$.  
Then Theorem \ref{thm:152271} follows from Theorem \ref{thm:regexist}
 and Proposition \ref{prop:172425}.  
\end{proof}

\subsubsection{Construction of $\Atbb \lambda \nu {\pm}{i,j}$
 for $j \in \{i-1,i\}$}

In this section 
 we apply Theorem \ref{thm:152389}
 about the construction
 of the (generically) regular symmetry breaking operators
 $\Atbb \lambda \nu {\pm}{V,W}$
 in the setting \eqref{eqn:VWexterior}
 with $j=i-1$ or $i$.  
In particular,
 we give concrete formul{\ae}
 of the matrix-valued distribution kernels
 $\Atcal \lambda \nu {\pm}{i,j}$
 for the operators.

Let $j=i-1$ or $i$.  
We recall from \eqref{eqn:Tii1} and \eqref{eqn:Tii2}
 that the projection $\pr i j \colon \Exterior^i({\mathbb{C}}^n) \to \Exterior^j({\mathbb{C}}^{n-1})$
 defines an element of 
\[
   {\operatorname{Hom}}_{O(n-1)}(V,W)
   =
   {\operatorname{Hom}}_{O(n-1)}({\Exterior}^i({\mathbb{C}}^n),{\Exterior}^j({\mathbb{C}}^{n-1})).  
\]
Denote by $\sigma \equiv \sigma^{(i)}$ 
 the $i$-th exterior representation of $O(n)$ on ${\Exterior}^i({\mathbb{C}}^n)$.  
Then the matrix-valued function $\Rij V W$ 
(see \eqref{eqn:RVW})
 amounts to the following map  
\[
\Rij i j : {\mathbb{R}}^n \setminus \{ 0 \} 
            \to 
           \operatorname{Hom}_{\mathbb{C}}({\Exterior}^i({\mathbb{C}}^n),
                            {\Exterior}^j({\mathbb{C}}^{n-1}) )
\]
given by 
\index{A}{Rij@$\Rij ij$|textbf}
\index{A}{prij@$\pr i j$, projection}
\index{A}{1psin@$\psi_n$}
\begin{equation}
\label{eqn:Rij}
   \Rij ij := \pr i j  \circ \sigma \circ \psi_n
\end{equation}
where we recall from \eqref{eqn:psim}
 that $\psi_n \colon {\mathbb{R}}^n \setminus\{0\} \to O(n)$
 is the map of taking \lq\lq{reflection}\rq\rq.  

\medskip
Applying the general formul{\ae} \eqref{eqn:KVWpt} and \eqref{eqn:KVWmt}
 of the distribution kernels 
 $\Atcal \lambda \nu {\pm}{V,W}$
 in the setting \eqref{eqn:VWexterior}, 
 we obtain 
$
    {\operatorname{Hom}}_{\mathbb{C}}
    (\Exterior^i({\mathbb{C}}^n),\Exterior^j({\mathbb{C}}^{n-1}))
$-valued locally integrable functions
 on ${\mathbb{R}}^n$
 for $\operatorname{Re}\lambda \gg |\operatorname{Re} \nu|$
 as follows.  
\begin{align}
\Atcal {\lambda}{\nu}{+}{i,j}
 :=&
\frac{1}{\Gamma(\frac{\lambda + \nu -n +1}{2})
         \Gamma(\frac{\lambda-\nu}{2})}
(|x|^2+ x_n^2)^{-\nu}
|x_n|^{\lambda+\nu-n}
\Rij ij (x,x_n), 
\label{eqn:Kijtilde}
\\
\Atcal {\lambda}{\nu}{-}{i,j}
 :=&
\frac{1}{\Gamma(\frac{\lambda + \nu -n +2}{2})
         \Gamma(\frac{\lambda-\nu+1}{2})}
(|x|^2 +x_n^2)^{-\nu}
|x_n|^{\lambda+\nu-n}
 {\operatorname{sgn}} x_n
\Rij i j(x,x_n).  
\label{eqn:Kijpmtilde}
\end{align}
Then, 
 as a special case
 of Theorem \ref{thm:152389}, 
 we obtain:
\begin{theorem}
[holomorphic continuation of integral operators]
\label{thm:1522a}
Let $(V,W)$ be as in \eqref{eqn:VWexterior}
 with $j=i,i-1$, 
 and $\delta, \varepsilon \in \{\pm \}$.  
Then the distributions $\Atcal \lambda \nu {\delta\varepsilon} {i,j}$, 
 initially defined as ${\operatorname{Hom}}_{\mathbb{C}}(V,W)$-valued 
 locally integrable functions
 on ${\mathbb{R}}^n$ 
 for $\operatorname{Re}\lambda \gg |\operatorname{Re} \nu|$, 
 extends to 
$
   ({\mathcal{D}}'(G/P, {\mathcal{V}}_{\lambda,\delta}^{\ast}) \otimes W_{\nu,\varepsilon})^{\Delta(P')}
$
 that depends holomorphically on $(\lambda, \nu)$
 in ${\mathbb{C}}^2$.  
Then the matrix-valued distribution kernels
\index{A}{Acts0@$\Atcal \lambda \nu {\pm} {i,j}$}
$
\Atcal \lambda \nu {\delta\varepsilon} {i,j}
$
 induce a family of symmetry breaking operators
\index{A}{Ahtsln0@$\Atbb \lambda \nu {\pm} {i,j}$|textbf}
\begin{equation}
\label{eqn:Aijdef}
  \Atbb \lambda \nu {\delta\varepsilon} {i,j}
  \colon
  I_{\delta}(i,\lambda) \to J_{\varepsilon}(j,\nu)
\end{equation}
 for all $(\lambda, \nu) \in {\mathbb{C}}^2$.  
\end{theorem}

Then $\Atbb \lambda \nu {\pm} {i,j}$ is the normalized (generically) regular 
 symmetry breaking operator (Definition \ref{def:AregSBO})
 in the sense that there exists an open dense subset
 $U_{\gamma}$ in ${\mathbb{C}}^2$
 for $\gamma \in \{\pm\}$
 such that the support of the distribution kernel 
 of $\Atbb \lambda \nu \gamma {i,j}$ equals 
 the whole flag manifold $G/P$
 as far as $(\lambda,\nu) \in U_{\gamma}$, 
 see Proposition \ref{prop:172425}.  
By a little abuse of terminology,
 we say that $\{\Atbb \lambda \nu \pm {i,j}\}$
 is a family of 
\index {B}{regularsymmetrybreakingoperator@regular symmetry breaking operator}
 {\it{normalized regular symmetry breaking operators}}.

\subsection{Zeros of $\Atbb \lambda \nu {\pm} {i,j}$ : Proof of Theorem \ref{thm:161243}}
\label{subsec:Aijzero}
In this section 
 we determine the exact place of the zeros
 of the normalized regular symmetry breaking operators
 $\Atbb \lambda \nu {\delta\varepsilon} {i,j}$, 
 and thus give a proof 
 of Theorem \ref{thm:161243}.  
In particular,
 we see that the Gamma factors
 in the normalization \eqref{eqn:Kijtilde} and \eqref{eqn:Kijpmtilde}
 are optimal 
 in the sense
 that the zeros of $\Atbb \lambda \nu {\gamma} {i,j}$
 are of codimension two in ${\mathbb{C}}^2$, 
 namely,
 form a discrete subset of ${\mathbb{C}}^2$. 
The proof of Theorem \ref{thm:161243}
 consists of the following steps.  
\par\noindent
{\bf{Step 0.}}\enspace
(existence condition)\enspace
Regular symmetry breaking operators from 
 $I_{\delta}(i,\lambda)$ to $J_{\varepsilon}(j,\nu)$ exist
 if and only if $j \in \{i-1,i\}$
 (Theorem \ref{thm:152271}).  
\par\noindent
{\bf{Step 1.}}\enspace
(generically nonzero)\enspace
If $\Atbb \lambda \nu {\delta\varepsilon} {i,j}=0$, 
 then $(\lambda, \nu, \delta, \varepsilon)$ belongs to the set 
 $\Psising$ 
 of special parameters
 (Theorem \ref{thm:genbasis}).  
\par\noindent
{\bf{Step 2.}}\enspace
(residue formula)\enspace
If $(\lambda, \nu, \delta, \varepsilon)\in \Psising$, 
 then $\Atbb \lambda \nu {\delta\varepsilon} {i,j}$ is proportional
 to the {\it{differential}} symmetry breaking operator
 $\Cbb \lambda \nu {i,j}$ 
 with explicit proportional constant
 (Fact \ref{fact:153316}).  

\subsubsection{Residue formula of the regular symmetry breaking operator
 $\Atbb \lambda \nu {\pm}{i,j}$}

Generalizing the residue formula
 of the scalar-valued regular symmetry breaking operators
 $\Atbb \lambda \nu {+}{0,0}$
 for spherical principal series representations
 given in \cite{xkEastwood} (see also \cite[Thm.~12.2]{sbon}), 
 we determined the residue
 of the matrix-valued regular symmetry breaking operators
 $\Atbb \lambda \nu {\pm}{i,j}$
 in \cite{xkresidue}, 
 as follows:
\begin{fact}
[residue formula {\cite[Thm.~1.3]{xkresidue}}]
\label{fact:153316}
\index{B}{residueformula@residue formula|textbf}
Let $\Cbb \lambda \nu {i,j}$ be the differential symmetry breaking operators
 defined in \eqref{eqn:Ciiln}
 and \eqref{eqn:Cijln} for $j=i-1$ or $i$.  
\begin{enumerate}
\item[{\rm{(1)}}]
Suppose $\nu-\lambda=2\ell$ with $\ell \in {\mathbb{N}}$.  
Then, 
\index{A}{Ciiln@$\Cbb \lambda \nu {i,j}$, matrix-valued differential operator}
\begin{equation}
\label{eqn:Aijres1}
\Atbb \lambda \nu {+} {i,j} 
=
\frac{(-1)^{i-j+\ell} \pi^{\frac{n-1}{2}} \ell !}{2^{2\ell-1}\Gamma(\nu+1)}
\Cbb{\lambda}{\nu}{i,j}.  
\end{equation}
\item[{\rm{(2)}}]
Suppose $\nu-\lambda=2\ell+1$ with $\ell \in {\mathbb{N}}$.  
Then, 
\[
\Atbb \lambda \nu {-} {i,j} 
=
\frac{(-1)^{i-j+\ell+1} \pi^{\frac{n-1}{2}} \ell !}{2^{2\ell+2}\Gamma(\nu+1)}
\Cbb{\lambda}{\nu}{i,j}.  
\]
\end{enumerate}
\end{fact}

We may unify the two formul\ae\
 in Fact \ref{fact:153316}
 into one formula:
for $\nu -\lambda \in {\mathbb{N}}$ and $j \in \{i,i-1\}$, 
\begin{equation}
\label{eqn:qAC}
\Atbb \lambda \nu {(-1)^{\nu-\lambda}} {i,j}
=
\frac{2(-1)^{i-j} \pi^{\frac {n-1}2}}{q(\nu-\lambda)\Gamma(\nu+1)}
\Cbb \lambda \nu {i,j}, 
\end{equation}
where we set,
 for $m \in {\mathbb{N}}$, 
\index{A}{qm@$q(m)$|textbf}
\begin{equation}
\label{eqn:qm}
q(m):=
\begin{cases}
\dfrac{(-1)^{\ell} 2^{2\ell}}{\ell!}
\qquad
&\text{if $m=2\ell$}, 
\\
&
\\
\dfrac{(-1)^{\ell+1} 2^{2\ell+3}}{\ell!}
\qquad
&\text{if $m=2\ell+1$}.  
\end{cases}
\end{equation}

\subsubsection{Zeros of $\Atbb \lambda \nu {\pm} {i,j}$}

The zeros of the operators $\Atbb \lambda \nu {\delta\varepsilon} {i,i}$
 for the special parameter
 in 
\index{A}{1psi@$\Psising$,
          special parameter in ${\mathbb{C}}^2 \times \{\pm\}^2$}
 $\Psising$
 (see \eqref{eqn:singset} for the definition)
 were determined in \cite{xkresidue}
 as a corollary of the residue formula
 (Fact \ref{fact:153316}), 
 which we recall now.  
\begin{corollary}
[zeros of $\Atbb \lambda \nu {\pm} {i,i}$
 for $\Psising$, 
{\cite[Thm.~8.1]{xkresidue}}]
\label{cor:Avanish}
\begin{enumerate}
\item[{\rm{(1)}}]
Suppose $\nu - \lambda \in 2 {\mathbb{N}}$.  
\newline
$\Atbb \lambda \nu {+} {i,i}=0$
if and only if
\begin{equation*}
(\lambda, \nu)\in 
\begin{cases}
L_{\operatorname{even}}
&\text{for $i=0$}, 
\\
(L_{\operatorname{even}} \setminus \{\nu =0\})\cup \{(i,i)\}
\quad
&\text{for $1 \le i \le n-1$}.   
\end{cases}
\end{equation*}
$\Atbb \lambda \nu {+} {i,i-1}=0$
if and only if
\begin{equation*}
(\lambda, \nu)\in 
\begin{cases}
(L_{\operatorname{even}} \setminus \{\nu =0\})\cup \{(n-i,n-i)\}
\quad
&\text{for $1 \le i \le n-1$}, 
\\
L_{\operatorname{even}}
&\text{for $i=n$}.   
\end{cases}
\end{equation*}
\item[{\rm{(2)}}]
Suppose $\nu - \lambda \in 2 {\mathbb{N}}+1$.  
\newline
$\Atbb \lambda \nu {-} {i,i}=0$
if and only if
\begin{equation*}
(\lambda, \nu)\in 
\begin{cases}
L_{\operatorname{odd}}
&\text{for $i=0$}, 
\\
L_{\operatorname{odd}} \setminus \{\nu =0\}
\quad
&\text{for $1 \le i \le n-1$}.   
\end{cases}
\end{equation*}
$\Atbb \lambda \nu {-} {i,i-1}=0$
if and only if
\begin{equation*}
(\lambda, \nu)\in 
\begin{cases}
L_{\operatorname{odd}} \setminus \{\nu =0\}
\quad
&\text{for $1 \le i \le n-1$}, 
\\
L_{\operatorname{odd}}
&\text{for $i=n$}.   
\end{cases}
\end{equation*}
\end{enumerate}
\end{corollary}
We are ready to complete the proof
 of Theorem \ref{thm:161243}
 on the zeros of the analytic continuation
 $\Atbb \lambda \nu \gamma {i,j}$
 of regular symmetry breaking operators.  
\begin{proof}
[Proof of Theorem \ref{thm:161243}]
We apply Theorem \ref{thm:genbasis}
 to the exterior representations
 \eqref{eqn:VWexterior},
 and see that  
$\Atbb \lambda \nu \gamma {i,j}=0$
 only if 
\begin{equation}
\label{eqn:AC}
\nu-\lambda \in 2 {\mathbb{N}}
\quad
(\gamma=+)
\quad
\text{or}
\quad
\nu-\lambda \in 2 {\mathbb{N}}+1
\quad
(\gamma=-).  
\end{equation}
Then Theorem \ref{thm:161243} follows from
 Corollary \ref{cor:Avanish}.  
\end{proof}

\subsection{$(K,K')$-spectrum for symmetry breaking operators}
\label{subsec:Kspec}

The second goal of this chapter is to formulate the concept
 of the $(K,K')$-spectrum
 for symmetry breaking operators (Definition \ref{def:KKspec}), 
 and give an explicit formula
 of the $(K,K')$-spectrum
\index{A}{SAln3ij@$S(\Atbb \lambda \nu {\varepsilon} {i,j})$, $(K,K')$-spectrum}\begin{equation}
\label{eqn:SAabcd}
S(\Atbb \lambda \nu \varepsilon {i,j})
=
\begin{pmatrix}
a_{\varepsilon}^{i,j}(\lambda,\nu)
&
b_{\varepsilon}^{i,j}(\lambda,\nu)
\\
c_{\varepsilon}^{i,j}(\lambda,\nu)
&
d_{\varepsilon}^{i,j}(\lambda,\nu)
\end{pmatrix}, 
\end{equation}
 (see \eqref{eqn:Smat}), 
 for the regular symmetry breaking operator
$
    \Atbb \lambda \nu \varepsilon{i,j}
    :
    I_{\delta}(i,\lambda)
    \to 
    J_{\delta \varepsilon}(j,\nu)
$
 with respect to basic $K$-types $\mu^{\natural}(i,\delta)$
 and $K'$-types $\mu^{\natural}(j,\delta \varepsilon)'$
 (see \eqref{eqn:muflat} and \eqref{eqn:musharp})
 for $\natural=\flat$ or $\sharp$.  
We will discuss the $(K,K')$-spectrum
 in Sections \ref{subsec:Kspec}--\ref{subsec:Apm1}.  
The main results are Theorem \ref{thm:153315}
 which will be proved in  Proposition \ref{prop:CApm1} (vanishing results)
 and Theorems \ref{thm:minKscalar}
 and \ref{thm:CApm1}.

One of the algebraic clues
 that we introduced in the study 
 of symmetry breaking operators $A$ in \cite{sbon}
 was an explicit formula
 of the \lq\lq{eigenvalues}\rq\rq\
 of $A$
 on spherical vectors.  
In the setting of this article,
 there is no spherical vector
 in the principal series representation
 $I_{\delta}(i,\lambda)$ if $i>0$
 or $J_{\varepsilon}(j,\nu)$ if $j>0$.  
In this section,
 we extend the idea of \cite{sbon}
 to the 
\index{B}{KspectrumKprime@$(K,K')$-spectrum|textbf}
$(K,K')$-{\it{spectrum}}
 for symmetry breaking operators
 with focus on 
\index{B}{basicKtype@basic $K$-type}
basic $K$-types.

\subsubsection{Generalities:
 $(K,K')$-spectrum of symmetry breaking operators}
We begin with a general setup.  
Let $(G,G')$ be a pair of real reductive Lie groups.  
Suppose $\Pi$ is a continuous representation of $G$, 
 and $\pi$ is that of the subgroup $G'$.  
We define a subset of $\widehat K \times \widehat{K'}$
 by 
\[
   {\mathcal{D}}(\Pi,\pi):=
   \{(\mu, \mu') \in \widehat K \times \widehat{K'}:
     [\Pi|_K:\mu], 
 [\pi|_{K'}:\mu'], 
 [\mu|_{K'}:\mu']
 \in \{0,1\}
\}.  
\]

Here is a sufficient condition for ${\mathcal{D}}(\Pi,\pi)$ 
 to be nonempty:
\begin{proposition}
\label{prop:170715}
Let $P=L N$
 and $P'=L' N'$ be parabolic subgroups of $G$ and its subgroup $G'$, 
 respectively.  
Suppose that $\Pi={\operatorname{Ind}}_P^G(\sigma \otimes {\mathbb{C}}_{\lambda})$
 and $\pi={\operatorname{Ind}}_{P'}^{G'}(\tau \otimes {\mathbb{C}}_{\nu})$
 are the induced representations from irreducible
 finite-dimensional representations
 $\sigma \otimes {\mathbb{C}}_{\lambda}$ of $L \simeq P/N$
 and $\tau \otimes {\mathbb{C}}_{\nu}$ of $L' \simeq P'/N'$, 
 respectively.  
\begin{enumerate}
\item[{\rm{(1)}}]
{\rm{(spherical principal series)}}\enspace
If $\sigma$ and $\tau$ are the trivial one-dimensional representations, 
then ${\mathcal{D}}(\Pi,\pi) \ni ({\bf{1}}_K, {\bf{1}}_{K'})$.  
\item[{\rm{(2)}}]
If $(K, L \cap K)$, $(K', L' \cap K')$ and $(K, K')$ are strong Gel'fand pairs, 
 in particular,
 if they are symmetric pairs, 
 then ${\mathcal{D}}(\Pi,\pi) = \widehat K \times \widehat {K'}$.  
\end{enumerate}
\end{proposition}
\begin{proof}
\begin{enumerate}
\item[(1)]
Clear from the Frobenius reciprocity.  
\item[(2)]
Immediate from the multiplicity-free property for strong Gel'fand pairs.  
\end{enumerate}
\end{proof}
The following is an example of Proposition \ref{prop:170715} (2).  
\begin{example}
\label{ex:5.2}
Let $(G,G')=(O(n+1,1),O(n,1))$, 
 and we consider $\Pi= I_{\delta}(V,\lambda)$, 
$\pi=J_{\varepsilon}(W, \nu)$
 for any $(\sigma, V) \in \widehat {O(n)}$
 and any $(\tau, W) \in \widehat {O(n-1)}$.  
Then ${\mathcal{D}}(\Pi,\pi)= \widehat K \times \widehat{K'}$.  
\end{example}

Now we introduce a $(K,K')$-spectrum
 for symmetry breaking operators as follows.  

\begin{definition}
[$(K,K')$-spectrum]
\label{def:KKspec}
Let $(\mu,\mu') \in {\mathcal{D}}(\Pi,\pi)$.  
If 
$[\Pi|_K:\mu]
=[\pi|_{K'}:\mu']
=[\mu|_{K'}:\mu']
=1$, 
 then we fix a nonzero $K$-homomorphism
 $\varphi:\mu \hookrightarrow \Pi$
 and nonzero $K'$-homomorphisms
 $\varphi':\mu' \hookrightarrow \pi$
 and $\iota:\mu' \hookrightarrow \mu$
 that are unique up to scalar multiplication.  
Suppose $A \in {\operatorname{Hom}}_{G'}(\Pi|_{G'}, \pi)$.  
Then by Schur's lemma,
 there exists a constant 
 $S_{\mu,\mu'}(A) \in {\mathbb{C}}$
 such that 
\begin{equation}
\label{eqn:Smumu}
A \circ \varphi \circ \iota
=
S_{\mu,\mu'}(A)\circ \varphi'
\qquad
\text{on $\mu'$}.  
\end{equation}
If one of $[\Pi|_K:\mu]$, 
 $[\pi|_{K'}:\mu']$,
 or $[\mu|_{K'}:\mu']$ is $0$, 
 then we just set
\[
   S_{\mu,\mu'}(A)=0
   \quad
   \text{for any }
   A \in   \operatorname{Hom}_{G'}(\Pi|_{G'}, \pi).  
\]
Thus we have defined a map
\begin{equation}
\label{eqn:KKspec}
 S \colon 
 \operatorname{Hom}_{G'}(\Pi|_{G'}, \pi)
 \times 
 {\mathcal{D}}(\Pi, \pi)
 \to {\mathbb{C}}, 
\qquad
(A, (\mu,\mu')) \mapsto S_{\mu,\mu'}(A).  
\end{equation}
We say $S_{\mu,\mu'}(A)$ is the {\it{$(K,K')$-spectrum}}
 of the symmetry breaking operator $A$
 for $(\mu,\mu') \in \widehat K \times \widehat {K'}$.  
We note that it is independent
 of the choice of the normalizations of $\varphi$, $\varphi'$, and $\iota$
 whether $S_{\mu,\mu'}(A)$ vanishes or not.  
\end{definition}

\subsection{Explicit formula of $(K,K')$-spectrum on basic $K$-types
 for regular symmetry breaking operators
 $\Atbb \lambda \nu {\pm}{i,j}$}
We return to our setting
 where $(G,G')=(O(n+1,1),O(n,1))$, 
 and thus
\[
K=O(n+1) \times O(1)
\supset
K'=O(n) \times O(1).  
\]
We consider a pair of representations
  $\Pi = I_{\delta}(i,\lambda)$
 of $G=O(n+1,1)$
 and $\pi=J_{\varepsilon}(j,\nu)$ of the subgroup $G'=O(n,1)$. 
In this case 
 ${\mathcal{D}}(\Pi,\pi)=\widehat K \times \widehat{K'}$
 as we saw in Example \ref{ex:5.2}, 
 however, 
 the following finite subset 
\[
   {\mathcal{D}}^{\flat, \sharp}
   \equiv
   {\mathcal{D}}^{\flat, \sharp}(\Pi,\pi)
   :=
   \{\mub (i,\delta), \mus (i,\delta)\}
   \times
   \{\mub (j,\varepsilon)', \mus (j,\varepsilon)'\}
   \subset \widehat K \times \widehat {K'}
\]
 will be sufficient for the later analysis of symmetry breaking operators.  
Here we recall from \eqref{eqn:muflat} and \eqref{eqn:musharp}
 that 
\index{A}{0muflat@$\mub(i,\delta)$}
$\mub (i,\delta)$
 and 
\index{A}{0musharp@$\mus(i,\delta)$}
$\mus (i,\delta)$ are 
 \lq\lq{basic $K$-types}\rq\rq\ of the principal series representation
 $I_{\delta}(i,\lambda)$ of $G$
 and that $\mub (j,\varepsilon)'$ and $\mus(j,\varepsilon)'$
 are those for $J_{\varepsilon}(j,\nu)$
 of the subgroup $G'$.

Then the $(K,K')$-spectrum restricted
 to the subset ${\mathcal{D}}^{\flat,\sharp}$
 is described as a $2 \times 2$ matrix: 
\begin{equation}
\label{eqn:Smat}
S:\operatorname{Hom}_{G'}(I_{\delta}(i,\lambda)|_{G'},J_{\varepsilon}(j,\nu))
  \to 
  M(2,{\mathbb{C}}),
  \qquad
  A 
  \mapsto
  \begin{pmatrix} a & b \\ c & d \end{pmatrix}
\end{equation}
 by taking $a$, $b$, $c$, $d$
 to be $S_{\mu,\mu'}(A)$ as follows:

\begin{equation*}
\begin{tabular}{c|rl|rl}
{$S_{\mu,\mu'}(A)$}
&
&{$\mu$}
&
&{$\mu'$}
\\
\hline
  $a$ 
& $\mub(i,\delta)$
& $=\Exterior^i({\mathbb{C}}^{n+1}) \boxtimes \delta$ 
& $\mub(j,\varepsilon)'$
& $=\Exterior^j({\mathbb{C}}^{n}) \boxtimes \varepsilon$
\\
$b$ 
& $\mub(i,\delta)$
& $=\Exterior^i({\mathbb{C}}^{n+1}) \boxtimes \delta$ 
& $\mus(j,\varepsilon)'$
& $=\Exterior^{j+1}({\mathbb{C}}^{n}) \boxtimes (-\varepsilon)$
\\
$c$ 
& $\mus(i,\delta)$
& $=\Exterior^{i+1}({\mathbb{C}}^{n+1}) \boxtimes (-\delta)$
& $\mub(j,\varepsilon)'$
&$=\Exterior^{j}({\mathbb{C}}^{n}) \boxtimes \varepsilon$
\\
$d$ 
& $\mus(i,\delta)$
& $=\Exterior^{i+1}({\mathbb{C}}^{n+1}) \boxtimes (-\delta)$
\qquad
& $\mus(j,\varepsilon)'$
&=$\Exterior^{j+1}({\mathbb{C}}^{n}) \boxtimes (-\varepsilon)$
\end{tabular}
\end{equation*}
To be more precise,
 we need a normalization
of the map $\varphi$, $\varphi'$
 and $\iota$
 in Definition \ref{def:KKspec} in this setting.  
For this,
 we realize the $K$-types
 $\mub(i,\delta)=\Exterior^i({\mathbb{C}}^{n+1}) \boxtimes \delta$
 and $\mus(i,\delta)=\Exterior^{i+1}({\mathbb{C}}^{n+1}) \boxtimes (-\delta)$
 in $I_{\delta}(i,\lambda)$
 as in Proposition \ref{prop:minKN}.  
Similarly,
 $\mub(j,\varepsilon)'=\Exterior^j({\mathbb{C}}^{n}) \boxtimes \varepsilon$
 and 
 $\mus(j,\varepsilon)'=\Exterior^{j+1}({\mathbb{C}}^{n}) \boxtimes (-\varepsilon)$ 
 are realized in $J_{\varepsilon}(j,\nu)$.  
When $\mu'$ and $\mu$ are representations
 on the exterior tensor spaces $\Exterior^l({\mathbb{C}}^{n})$
 and $\Exterior^k({\mathbb{C}}^{n+1})$
 ($l=k$ or $k-1$) respectively,
 we normalize an $O(n)$-homomorphism 
\[
   \iota_{l \to k}\colon \Exterior^l({\mathbb{C}}^n) 
         \hookrightarrow \Exterior^k({\mathbb{C}}^{n+1})
\]
 such that $\pr k l \circ \iota_{l \to k}= {\operatorname{id}}$, 
 where the projection 
\index{A}{prij@$\pr i j$, projection}
$\pr k l : \Exterior^k({\mathbb{C}}^{n+1}) \to \Exterior^l({\mathbb{C}}^{n})$
 is defined in \eqref{eqn:Tii1} and \eqref{eqn:Tii2}.  
With these normalizations,
 the map \eqref{eqn:Smat} is defined.  
We obtain the following closed formula
 of the $(K,K')$-spectrum for the normalized regular 
 symmetry breaking operators 
\index{A}{Ahtsln0@$\Atbb \lambda \nu {\pm} {i,j}$}
$\Atbb \lambda \nu {\pm}{i,j} \colon I_{\delta}(i,\lambda) \to J_{\pm \delta}(j,\nu)$.  
\begin{theorem}
[$(K,K')$-spectrum for $\Atbb \lambda \nu {\pm} {i,j}$]
\label{thm:153315}
Suppose $(\lambda, \nu)\in {\mathbb{C}}^2$.  
Then the $(K,K')$-spectrum of the analytic continuation
 $\Atbb \lambda \nu {\pm}{i,j}$
 of regular symmetry breaking operators
 takes the following form on basic $K$-types:
\index{A}{SAln3ij@$S(\Atbb \lambda \nu {\varepsilon} {i,j})$, $(K,K')$-spectrum}

\begin{alignat*}{2}
S(\Atbb \lambda \nu {+} {i,i})
=&
\frac{\pi^{\frac{n-1}{2}}}{\Gamma(\lambda+1)}
\begin{pmatrix} \lambda-i & 0 \\ 0 & \nu-i \end{pmatrix}
\quad
&&
\text{for $0 \le i \le n-1$};
\\
S(\Atbb \lambda \nu {-} {i,i})
=&
\frac{\pi^{\frac{n-1}{2}}}{\Gamma(\lambda+1)}
\begin{pmatrix}0 & 0 \\ 2(-1)^{i+1} & 0 \end{pmatrix}
\quad
&&
\text{for $0 \le i \le n-1$};
\\
S(\Atbb \lambda \nu {+} {i,i-1})
=&
\frac{\pi^{\frac{n-1}{2}}}{\Gamma(\lambda+1)}
\begin{pmatrix} n-\nu-i & 0 \\ 0 & \lambda-n+i \end{pmatrix}
\quad
&&
\text{for $1 \le i \le n$};
\\
S(\Atbb \lambda \nu {-} {i,i-1})
=&
\frac{\pi^{\frac{n-1}{2}}}{\Gamma(\lambda+1)}
\begin{pmatrix} 0 & -2 \\ 0 & 0 \end{pmatrix}
\quad
&&
\text{for $1 \le i \le n$}.
\end{alignat*}
\end{theorem}
The vanishing result (an easy part) of Theorem \ref{thm:153315}
 will be shown in Proposition \ref{prop:CApm1}, 
 and the remaining nontrivial part will be proved in Theorems \ref{thm:minKscalar} and \ref{thm:CApm1}.

\subsection{Proof of vanishing results on $(K,K')$-spectrum}
In this section,
 we formulate and prove vanishing results
 for $(K,K')$-spectrum
 that hold for general symmetry breaking operators.  
\begin{proposition}
\label{prop:CApm1}
Suppose $j \in \{i-1,i\}$, 
$\delta, \varepsilon \in \{\pm\}$, 
 and $\lambda, \nu \in {\mathbb{C}}$.  
Let $A:I_{\delta}(i,\lambda) \to J_{\varepsilon}(j,\nu)$
 be an arbitrary symmetry breaking operator.  
Then the $(K,K')$-spectrum $S(A)$
 for basic $K$-types takes the following form:
\begin{equation*}
\begin{tabular}{c|cccc}
$j$
&$i$
&$i$
&$i-1$
&$i-1$
\\
$\delta\varepsilon$
&$+$
&$-$
&$+$
&$-$
\\
\hline
$S(A)$ 
& $\begin{pmatrix} \ast & 0 \\ 0 & \ast \end{pmatrix}$
& $\begin{pmatrix} 0    & 0 \\ \ast & 0 \end{pmatrix}$
& $\begin{pmatrix} \ast & 0 \\ 0 & \ast \end{pmatrix}$
& $\begin{pmatrix} 0 & \ast \\ 0 & 0 \end{pmatrix}$
\end{tabular}
\end{equation*}
\end{proposition}

\begin{proof}
Without loss of generality,
 we may assume $\delta=+$.  
The $K$-modules $\mub(i,+)$ and $\mus(i,+)$
 (see \eqref{eqn:muflat} and \eqref{eqn:musharp})
 decompose into the sum of irreducible representations
 of the subgroup $K'$:
\begin{alignat*}{4}
 \mub(i,+) &= \Exterior^{i}({\mathbb{C}}^{n+1}) \boxtimes {\bf{1}}
&&\simeq
&& \Exterior^{i}({\mathbb{C}}^{n}) \boxtimes {\bf{1}}
&&    \oplus 
 \Exterior^{i-1}({\mathbb{C}}^{n}) \boxtimes {\bf{1}}, 
\\
 \mus(i,+) &= \Exterior^{i+1}({\mathbb{C}}^{n+1}) \boxtimes {\operatorname{sgn}}
&&\simeq
&& \Exterior^{i+1}({\mathbb{C}}^{n}) \boxtimes {\operatorname{sgn}}
&&    \oplus
    \Exterior^{i}({\mathbb{C}}^{n}) \boxtimes {\operatorname{sgn}}.  
\end{alignat*}

Using the notion $\mu^{\natural}(j,\pm)'$
 with $\natural=\flat$ or $\sharp$
 for $K'$-types,
 we may rewrite these decompositions as
\begin{align}
\mub(i,+)|_{K'}
\simeq & \mub(i,+)' \oplus \mub(i-1,+)'
\label{eqn:smbra1}
\\
\simeq & \mus(i-1,-)' \oplus \mus(i-2,-)', 
\notag
\\
\mus(i,+)|_{K'}
\simeq & \mus(i,+)' \oplus \mus(i-1,+)'
\label{eqn:smbra2}
\\
\simeq & \mub(i+1,-)' \oplus \mub(i,-)'.  
\notag
\end{align}
The second isomorphisms follow from \eqref{eqn:flatsharp}.

For simplicity,
 we discuss the symmetry breaking operator
 $A\colon I_{\delta}(i,\lambda) \to J_{\varepsilon}(j,\nu)$
 in the case $j=i$, $\delta=+$, and $\varepsilon=-$.  
Then the branching rule \eqref{eqn:smbra1} tells
 that neither the $K'$-type $\mub (i,-)'$ nor $\mus(i,-)'$
 occurs in the $K$-type $\mub (i,+)$ of $I_+(i,\lambda)$.  
Likewise, 
 \eqref{eqn:smbra2} tells
 that the $K'$-type $\mus(i,-)'$
 does not occur in the $K$-type $\mus (i,+)$.  
Hence the matrix $S(A)$ in \eqref{eqn:Smat} must be of the form
 $\begin{pmatrix} 0 & 0 \\ \ast & 0\end{pmatrix}$.

The vanishing statements in the other cases are proved similarly.  
\end{proof}

\subsection{Proof of Theorem \ref{thm:153315}
 on $(K,K')$-spectrum
\\
 for the normalized symmetry breaking operator 
 $\Atbb{\lambda}{\nu}{+}{i,j}\colon
I_{\delta}(i,\lambda) \to J_{\delta}(j,\nu)$}
\label{subsec:Apm1-5}

In this section,
 we determine the $(K,K')$-spectrum $a_{\varepsilon}^{i,j}(\lambda,\nu)$
 and $d_{\varepsilon}^{i,j}(\lambda,\nu)$
 for $j=i, i-1$ 
 in \eqref{eqn:SAabcd}
 when $\varepsilon=+$.  
The case $\varepsilon=-$ will be discussed separately
 in Section \ref{subsec:Apm1}.  
By definition \eqref{eqn:Smumu}, 
 the constants $a_{\varepsilon}^{i,j}(\lambda,\nu)$
 and $d_{\varepsilon}^{i,j}(\lambda,\nu)$ are characterized 
 by the following equations:

\index{A}{Ahtsln1@$\Atbb \lambda \nu {+} {i,j}$}
\index{A}{0iotalmdd@$\iota_{\lambda}^{\ast}$}
\begin{alignat}{2}
\label{eqn:Aaij}
\Atbb \lambda \nu + {i,j}
\circ
\iota_{\lambda}^{\ast}
\circ
\iota_{j \to i}
=&
a_+^{i,j}(\lambda,\nu)
\iota_{\nu}^{\ast}
\quad
&&\text{on $\Exterior^{j}({\mathbb{C}}^{n})$,}
\\
\label{eqn:Adij}
\Atbb \lambda \nu + {i,j}
\circ
\iota_{\lambda}^{\ast}
\circ
\iota_{j+1 \to i+1}
=&
d_+^{i,j}(\lambda,\nu)
\iota_{\nu}^{\ast}
\quad
&&\text{on $\Exterior^{j+1}({\mathbb{C}}^{n})$}, 
\end{alignat}
where $\Atbb \lambda \nu + {i,j} \colon I_{\delta}(i,\lambda) \to J_{\delta}(j,\nu)$
 is the normalized symmetry breaking operator,
 $\iota_{\lambda}^{\ast}$ is the transform from the $K$-picture
 to the $N$-picture
 (see \eqref{eqn:KtoN}), 
 and 
\index{A}{0iotato@$\iota_{j \to i}$|textbf}
$
\iota_{j \to i} \colon \Exterior^j({\mathbb{C}}^{n}) \to \Exterior^i({\mathbb{C}}^{n+1})
$ is the normalized injective $O(n)$-homomorphism
 such that
\index{A}{prij@$\pr i j$, projection}
 $\pr i j \circ \iota_{j \to i} = {\operatorname{id}}$.  
The main results of this section
 are part of Theorem \ref{thm:153315}, 
 which is given as follows:

\begin{theorem}
\label{thm:minKscalar}
Suppose $\lambda, \nu\in {\mathbb{C}}$.  
\begin{align}
a_+^{i,i}(\lambda,\nu)
=&
\frac{\pi^{\frac{n-1}{2}} (\lambda-i)}
     {\Gamma(\lambda+1)}.  
\notag
\\
a_+^{i,i-1}(\lambda,\nu)
=&
\frac{\pi^{\frac{n-1}{2}} (n-\nu-i)}
     {\Gamma(\lambda+1)}.  
\notag
\\
\label{eqn:pmIi-5}
   d_+^{i,i}(\lambda,\nu)
   =&
   \frac
   {\pi^{\frac {n-1}{2}} (\nu-i)}
   {\Gamma(\lambda+1)}.  
\\
\label{eqn:152361}
d_+^{i,i-1}(\lambda,\nu)=&\frac{\pi^{\frac{n-1}{2}} (\lambda-n+i)}{\Gamma(\lambda+1)}.  
\end{align}
\end{theorem}

\begin{remark}
Theorem \ref{thm:minKscalar} generalizes 
 \cite[Thm~1.10]{sbon}
 in the spherical case
 ($i=j=0$ and $\delta = \varepsilon=+$).  
\end{remark}

The proof of Theorem \ref{thm:minKscalar} is divided
 into the following two steps:
\begin{enumerate}
\item[$\bullet$]
integral expression of 
\index{A}{a1ij@$a_+^{i,j}(\lambda, \nu)$}
$
a_+^{i,j}(\lambda,\nu)
$
 and 
\index{A}{d1ij@$d_+^{i,j}(\lambda, \nu)$}
$d_+^{i,j}(\lambda,\nu)
$
 (Section \ref{subsec:KKspec1});
\item[$\bullet$]
computation of the integral
 (Section \ref{subsec:KKspec2}).  
\end{enumerate}

\subsubsection{Integral expression of $(K,K')$-spectrum}
\label{subsec:KKspec1}

As the first step of the proof, 
 we give an integral expression
 of the $(K,K')$-spectrum $a_+^{i,j}(\lambda,\nu)$ and $d_+^{i,j}(\lambda,\nu)$.  
For $I \in {\mathfrak{I}}_{n,i}$, 
we recall from \eqref{eqn:QI}
 that the quadratic form $Q_I(b)$ is defined
 to be $\sum_{k \in I} {b_k}^2$, 
 and set 
\index{A}{QIb@$Q_I(b)$, quadratic polynomial}
\begin{align}
\alpha_I(b):=&1-\frac{2 Q_I(b)}{(1+|b|^2)|b|^2}, 
\label{eqn:alphaI}
\\
\delta_I(b):=&1-\frac{2 |b|^2}{1+|b|^2}-\frac{2 Q_I(b)}{(1+|b|^2)|b|^2}.  
\label{eqn:deltaI}
\intertext{Consider the following integrals:}
A_I(\lambda,\nu):=&\int_{{\mathbb{R}}^n}\Atcal \lambda \nu + {}(b)(1+|b|^2)^{-\lambda}\alpha_I(b) d b,
\notag
\\
D_I(\lambda,\nu)
:=&\int_{{\mathbb{R}}^n}\Atcal \lambda \nu + {}(b)(1+|b|^2)^{-\lambda}\delta_I(b) d b.
\notag  
\end{align}

Then the $(K,K')$-spectrum $a_+^{i,j}(\lambda,\nu)$ and $d_+^{i,j}(\lambda,\nu)$
 in \eqref{eqn:Aaij} and \eqref{eqn:Adij}, 
 respectively, 
 is given by the integrals
 $A_I(\lambda,\nu)$ and $D_I(\lambda,\nu)$
 as follows:
\begin{proposition}
[integral expression of $(K,K')$-spectrum]
\label{prop:1617113}
\begin{alignat*}{2}
a_+^{i,i}(\lambda, \nu)
=&
A_I(\lambda, \nu)
\qquad
&&\text{for any $I \in {\mathfrak {I}}_{n,i}$ with $n \not\in I$,}
\\
a_+^{i,i-1}(\lambda, \nu)
=&
A_I(\lambda, \nu)
\qquad
&&\text{for any $I \in {\mathfrak {I}}_{n,i}$ with $n \in I$, }
\\
d_+^{i,i}(\lambda, \nu)
=&
D_I(\lambda, \nu)
\qquad
&&\text{for any $I \in {\mathfrak {I}}_{n,i}$ with $n \not\in I$,}
\\
d_+^{i,i-1}(\lambda, \nu)
=&
-D_I(\lambda, \nu)
\qquad
&&\text{for any $I \in {\mathfrak {I}}_{n,i}$ with $n \in I$.}
\end{alignat*}
\end{proposition}

In order to prove Proposition \ref{prop:1617113}, 
 we use the $N$-picture of the principal series representations
 $I_{\delta}(i,\lambda)$ and $J_{\varepsilon}(j,\nu)$.  
By Proposition \ref{prop:minKN}
 for the vectors ${\bf{1}}_{\lambda}^{\mathcal{I}}$
 and $h_{\lambda}^{\mathcal{I}}$ belonging to the basic $K$-types,
 the equation \eqref{eqn:Aaij} means
 that for ${\mathcal{I}} \in {\mathfrak{I}}_{n+1,i}$
\index{A}{01oneIlmd@${\bf{1}}_{\lambda}^{\mathcal{I}}$}
\index{A}{hIlmd@$h_{\lambda}^{{\mathcal{I}}}$}
\begin{alignat*}{2}
\Atbb \lambda \nu {+} {i,i} {\bf{1}}_{\lambda}^{\mathcal{I}}
&=
a_+^{i,i}(\lambda,\nu) {{\bf{1}}_{\nu}'}^{{\mathcal{I}}}
\hphantom{(-1)^{i-}mii}
\quad
&&(n \notin {\mathcal{I}}), 
\\
\Atbb \lambda \nu {+} {i,i-1}
{\bf{1}}_{\lambda}^{\mathcal{I}}
&=
(-1)^{i-1} a_+^{i,i-1}(\lambda,\nu) {{\bf{1}}_{\nu}'}^{\mathcal{I}\setminus \{n\}}
\quad
&&(n\in {\mathcal{I}}).  
\end{alignat*}
The signature in the second formula
 arises from the definition \eqref{eqn:Tii2}
 of the projection $\pr i {i-1}$.

To compute the constants $a_+^{i,j}(\lambda,\nu)$, 
 we take $I \in {\mathfrak {I}}_{n, i}$
 and set ${\mathcal{I}}:=I$, 
 regarded as an element of ${\mathfrak {I}}_{n+1, i}$, 
 where we recall Convention \ref{conv:index} of index sets.  
Since $0 \not \in {\mathcal{I}}$, 
 it follows from \eqref{eqn:I10} that 
\begin{alignat*}{2}
\Atbb \lambda \nu {+} {i,i}
 ({\bf{1}}_{\lambda}^{I})(0)
 =&
 a_+^{i,i}(\lambda,\nu) e_{I}
\qquad
&&\text{if }\,\,
 n \not \in I, 
\\
\Atbb \lambda \nu {+} {i,i-1}
 ({\bf{1}}_{\lambda}^{I})(0)
 =& 
 (-1)^{i-1} a_+^{i,i-1}(\lambda,\nu) e_{I \setminus \{n\}}
 \qquad
&&\text{if }\,\,
 n \in I.   
\end{alignat*}
Likewise, 
 the equation \eqref{eqn:Adij} means 
that
for ${\mathcal{I}} \in {\mathfrak{I}}_{n+1,i+1}$
\begin{alignat*}{2}
\Atbb \lambda \nu {+} {i,i} h_{\lambda}^{\mathcal{I}}
=&
d_+^{i,i}(\lambda,\nu) {h'}_{\nu}^{{\mathcal{I}}}
\qquad\qquad\qquad
&&(n \notin {\mathcal{I}}), 
\\
\Atbb \lambda \nu {+} {i,i-1}
h_{\lambda}^{\mathcal{I}}
=&(-1)^i d_+^{i,i-1}(\lambda,\nu) {h_{\nu}'}^{\mathcal{I}\setminus\{n\}}
\quad
&&(n\in {\mathcal{I}}).  
\end{alignat*}
In this case, 
 we take $I \in {\mathfrak{I}}_{n,i}$
 and set ${\mathcal{I}}:=I \cup \{0\} \in {\mathfrak {I}}_{n+1, i+1}$.  
Then \eqref{eqn:hI0} implies 
\begin{alignat}{2}
\label{eqn:CApm-5}
(\Atbb \lambda \nu + {i,i} h_{\lambda}^{I \cup \{0\}})(0)
=&
d_+^{i,i}(\lambda,\nu) e_I
\quad
&&\text{if $n \not \in I$}, 
\\
(\Atbb \lambda \nu + {i,i-1} h_{\lambda}^{I \cup \{0\}})(0)
=&
(-1)^i d_+^{i,i-1}(\lambda,\nu) e_{I\setminus \{n\}}
\quad
&&\text{if $n \in I$.  }
\notag
\end{alignat}

Let us compute 
$
 \Atbb \lambda \nu {+} {i,j}
 ({\bf{1}}_{\lambda}^{I})(0)
$
 and 
$
 \Atbb \lambda \nu {+} {i,j}
 (h_{\lambda}^{I \cup \{0\}})(0)
$
 for $j=i$ and $i-1$.  
If $\operatorname{Re} \lambda \gg |\operatorname{Re} \nu|$, 
 then the matrix-valued distribution kernel 
 $\Atcal \lambda \nu + {i,j}$
 (see \eqref{eqn:Kijtilde})
 of the regular symmetry breaking operator
 $\Atbb \lambda \nu +{i,j}$ is decomposed as 
\[
   \Atcal \lambda \nu + {i,j}= \Atcal \lambda \nu + {} \Rij ij, 
\]
where 
\index{A}{Act1@$\Atcal \lambda \nu + {}$}
$\Atcal \lambda \nu + {}$ is the scalar-valued,
 locally integrable function
 defined in \eqref{eqn:KAlnn+}
 and the matrix-valued function 
\index{A}{Rij@$\Rij ij$}
$
\Rij ij \in C^{\infty}({\mathbb{R}}^n \setminus \{0\})
 \otimes 
 {\operatorname{Hom}}_{{\mathbb{C}}}(\Exterior^i({\mathbb{C}}^n), \Exterior^j({\mathbb{C}}^{n-1}))$
 is defined in \eqref{eqn:Rij}.  
Hence,
 we have
\begin{align*}
(\Atbb \lambda \nu +{i,j} \psi)(0)
=&\int_{\mathbb{R}^n} \Atcal \lambda \nu +{}(-b) \Rij ij (-b)\psi(b) d b
\\
=&\int_{\mathbb{R}^n} \Atcal \lambda \nu +{}(b) \Rij ij (b)\psi(b) d b
\end{align*}
in the $N$-picture for any $\psi \in \iota_{\lambda}^{\ast}({\mathcal{E}}^i(S^n))\subset C^{\infty}(\mathbb{R}^n) \otimes \Exterior^i(\mathbb{C}^n)$.  
Thus Proposition \ref{prop:1617113} is a consequence
 of the following two lemmas
 on the computation of $\Rij i j(b) \psi(b) \in \Exterior^j({\mathbb{C}}^{n-1})$
 for $\psi={\bf{1}}_{\lambda}^{I}$
 or $h_{\lambda}^{I \cup \{0\}}$ and 
 for $j=i$ or $i-1$.  
\begin{lemma}
\label{lem:eIRone}
Suppose $I \in {\mathfrak{I}}_{n,i}$.  
\begin{enumerate}
\item[{\rm{(1)}}]
If $n \not \in I$, 
then the coefficient of $e_I$
 in 
$
  \Rij i i (b) {\bf{1}}_{\lambda}^{I}(b)
$
is given by 
\[
  (1+|b|^2)^{-\lambda} \alpha_I(b)
  =
 (1+|b|^2)^{-\lambda}(1-\frac{2Q_I(b)}{(1+|b|^2)|b|^2}), 
\]
where we recall 
 $Q_I(b)=\sum_{l \in I}b_l^2$ from \eqref{eqn:QI}.  
\item[{\rm{(2)}}]
If $n \in I$, 
then the coefficient of $e_{I\setminus \{ n \}}$
 in 
$
  \Rij i {i-1} (b) {\bf{1}}_{\lambda}^{I}(b)
$
is given by 
\[
  (-1)^{i-1}(1+|b|^2)^{-\lambda} \alpha_I(b)
  = 
  (1+|b|^2)^{-\lambda} (1-\frac{2Q_I(b)}{(1+|b|^2)|b|^2}).  
\]
\end{enumerate}
\end{lemma}

\begin{lemma}
\label{lem:eIpm-5}
Suppose $I \in {\mathfrak{I}}_{n,i}$.  
\begin{enumerate}
\item[{\rm{(1)}}]
If $n \not \in I$, 
then the coefficient of $e_I$
 in 
$
  \Rij i i (b) h_{\lambda}^{I \cup \{0\}}(b)
$
is given by 
\[
  (1+|b|^2)^{-\lambda} \delta_I(b)
  =
 (1+|b|^2)^{-\lambda-1}(1-|b|^2-\frac{2Q_I(b)}{|b|^2}).  
\]
\item[{\rm{(2)}}]
If $n \in I$, 
then the coefficient of $e_{I\setminus \{ n \}}$
 in 
$
  \Rij i {i-1} (b) h_{\lambda}^{I \cup \{0\}}(b)
$
is given by 
\[
  (-1)^{i-1}(1+|b|^2)^{-\lambda} \delta_I(b)
  = 
  (-1)^{i-1} (1+|b|^2)^{-\lambda-1} (1-|b|^2-\frac{2Q_I(b)}{|b|^2}).  
\]
\end{enumerate}
\end{lemma}

\begin{proof}
[Proof of Lemma \ref{lem:eIRone}]
Let $\sigma$ be the $i$-th exterior representation
 on $\Exterior^i({\mathbb{C}}^n)$.  
We recall from \eqref{eqn:Rij}
 $\Rij ij=\pr ij \circ \sigma \circ \psi_n$.  
We identify $I \in {\mathfrak{I}}_{n,i}$ with ${\mathcal{I}} \in {\mathfrak{I}}_{n+1,i}$
 such that $n \not \in {\mathcal{I}}$ as usual, 
 and apply the formula \eqref{eqn:minI}
 of ${\bf{1}}_{\lambda}^I$.  
Then we have
\[
  \sigma(\psi_n(b)) {\mathbf{1}}_{\lambda}^{I}(b)
=
(1+|b|^2)^{-\lambda}
\sigma(\psi_n(b))
\sum_{J \in {\mathfrak{I}}_{n,i}} (\det \psi_{n+1}(1,b))_{I J} e_{J}.  
\]
By the formula \eqref{eqn:exrep}
 of the matrix coefficients of the exterior tensor representation, 
 the coefficient of $e_I$ 
 in $\sigma(\psi_n(b)) {\mathbf{1}}_{\lambda}^{I}(b)$
 amounts to 
\[
   (1+|b|^2)^{-\lambda}
   \sum_{J \in {\mathfrak {I}}_{n,i}}
   (\det \psi_{n+1}(1,b))_{I J}
   (\det \psi_n(b))_{I J},   
\]
which is equal to 
\[
  (1+|b|^2)^{-\lambda}(1-\frac{2Q_I(b)}{(1+|b|^2)|b|^2})
   =
   (1+|b|^2)^{-\lambda} \alpha_I(b)
\]
by the minor summation formula \eqref{eqn:kbsum}
 in Proposition \ref{prop:msum}.  
Hence the lemma follows from $\pr ii (e_I) = e_I$ $(n \notin I)$
 and $\pr i{i-1} (e_I) = (-1)^{i-1}e_{I\setminus\{n\}}$ $(n \in I)$
 (see \eqref{eqn:Tii1} and \eqref{eqn:Tii2}).  
\end{proof}

\begin{proof}
[Proof of Lemma \ref{lem:eIpm-5}]
The proof goes in parallel to that of Lemma \ref{lem:eIRone}.  
For the sake of completeness, 
 we give a proof.

By \eqref{eqn:hIlmd} and \eqref{eqn:exrep}, 
we have
\begin{align*}
&\sigma(\psi_n(b)) h_{\lambda}^{I \cup \{0\}}(b)
\\
=& -(1+|b|^2)^{-\lambda}
  \sigma(\psi_n(b))
  \sum_{J \in {\mathfrak{I}}_{n,i}}
  (\det \psi_{n+1}(1,b))_{I \cup \{0\}, J \cup \{0\}} e_J
\\
=& -(1+|b|^2)^{-\lambda}
  \sum_{J \in {\mathfrak{I}}_{n,i}}
  \sum_{J' \in {\mathfrak{I}}_{n,i}}
  (\det \psi_{n+1}(1,b))_{I \cup \{0\}, J \cup \{0\}}
  \det \psi_n(b)_{J' J}e_{J'}.  
\end{align*}

Hence the coefficient of $e_I$
 in $\sigma(\psi_n(b)) h_{\lambda}^{I \cup \{0\}}(b)$
 is equal to 
\begin{equation*}
  - (1+|b|^2)^{-\lambda}
  \sum_{J \in {\mathfrak{I}}_{n,i}}
  (\det \psi_{n+1}(1,b))_{I \cup \{0\}, J \cup \{0\}}
   \det \psi_n(b)_{I J}, 
\end{equation*}
which amounts to 
\[
  (1+|b|^2)^{-\lambda-1}(1-|b|^2-\frac{2Q_I(b)}{|b|^2})
  =
 (1+|b|^2)^{-\lambda} \delta_I(b)
\]
by the minor summation formula
 in Proposition \ref{prop:msum} (2).  
Thus we have shown the lemma.  
\end{proof}

Therefore we have completed the proof of Proposition \ref{prop:1617113}.

\subsubsection{Integral formula of the $(K,K')$-spectrum}
\label{subsec:KKspec2}
As the second step,
 we compute the integrals
 $A_I(\lambda,\nu)$ and $D_I(\lambda,\nu)$
 in Section \ref{subsec:KKspec1}.  
We begin with the following integral formul{\ae}:
Denote by $d \omega$ the standard measure
 on the unit sphere $S^{n-1}=\{\omega=(\omega_1, \cdots, \omega_n) \in {\mathbb{R}}^n:
\sum_{j=1}^{n} {\omega_j}^2=1\}$.

For $a,b \in {\mathbb{C}}$
 with $\operatorname{Re} a, \operatorname{Re} b>-1$, 
 we set
\begin{equation}
\label{eqn:intsphere}
S(a,b) 
\equiv
S_n(a,b)
:=
\int_{S^{n-1}} |\omega_n|^a |\omega_{n-1}|^b d \omega.  
\end{equation}
Then we have 
\begin{equation}
   S_n(a,0) = \int_{S^{n-1}} |\omega_n|^a d \omega
=
\frac{2 \pi^{\frac {n-1}{2}} \Gamma(\frac{a+1}{2})}{\Gamma(\frac{a+n}{2})}, 
\label{eqn:1.11}
\end{equation}
see \cite[Lemma 7.6]{sbon}, 
 for instance.  
More generally,
 we have the following.  
\begin{lemma}
\label{lem:intsphere}
Suppose ${\operatorname{Re}}\, a >-1$ and ${\operatorname{Re}}\, b >-1$.  
Then we have 
\begin{equation}
\label{eqn:Sab}
S(a,b)=
\frac
{2 \pi^{\frac{n-2}{2}}
\Gamma(\frac{a+1}{2})\Gamma(\frac{b+1}{2})}{\Gamma(\frac{a+b+n}{2})}.  
\end{equation}
\end{lemma}
It is convenient to write down the following recurrence relations
 that are derived readily from \eqref{eqn:Sab}:
\begin{align}
\label{eqn:Sab2}
S(a,2)=&
\frac{1}{a+n}S(a,0), 
\\
\label{eqn:Sab3}
S(a+2,0)=&
\frac{a+1}{a+n}S(a,0).  
\end{align}

\begin{proof}
[Proof of Lemma \ref{lem:intsphere}]
For any $f \in C(S^{n-1})$, 
 the polar coordinates give the following expression
 of the integral:
\begin{equation}
\int_{S^{n-1}} f (\omega) d \omega
=\int_{-1}^{1} \int_{S^{n-2}} 
 f(\sqrt{1-t^2} \eta, t) (1-t^2)^{\frac {n-3}{2}} d \eta d t.
\label{eqn:polarS}
\end{equation}
Then we have
\begin{align*}
S(a,b)
=& \int_{-1}^1 \int_{S^{n-2}}
  |\sqrt{1-t^2} \eta_{n-1} |^b |t|^a (1-t^2)^{\frac{n-3}{2}} d\eta d t 
\\
=& \int_{S^{n-2}}
  |\eta_{n-1}|^b d \eta \int_{-1}^1 |t|^a (1-t^2)^{\frac{n+b-3}{2}} d t. 
\end{align*}
The first term equals $S_{n-1}(b,0)$, 
 see \eqref{eqn:1.11}.  
The second term is given by the Beta function:
\begin{equation}
\int_{0}^{1} t^{2A-1}(1-t^2)^{B-1} d t 
=\frac{\Gamma(A)\Gamma(B)}{2\Gamma(A+B)}.  
\label{eqn:Beta}
\end{equation}
Here we get the lemma.  
\end{proof}

\begin{lemma}
\label{lem:1616116}
Let $\Atcal \lambda \nu +{}$ be the 
 (scalar-valued) locally integrable function on ${\mathbb{R}}^n$
 defined in \eqref{eqn:KAlnn+}
 for $\operatorname{Re}\left(\lambda-\nu\right)>0$
 and $\operatorname{Re}\left(\lambda+\nu\right)>n-1$.  
\begin{enumerate}
\item[{\rm{(1)}}]
We have 
\index{A}{Act1@$\Atcal \lambda \nu + {}$}
\[
   \int_{{\mathbb{R}}^n} \Atcal \lambda \nu + {}(b)
                         (1+|b|^2)^{-\lambda} d b
  =
  \frac{\pi^{\frac{n-1}{2}}}{\Gamma(\lambda)}.  
\]
\item[{\rm{(2)}}]
Let $l \in \{1,2,\cdots,n\}$.  
Then we have 
\begin{multline*}
\int_{{\mathbb{R}}^n} \Atcal \lambda \nu + {}(b)
                         (1+|b|^2)^{-\lambda} 
                       \frac{2 {b_{\ell}}^2}{(1+|b|^2)|b|^2} d b
\\
  =
\frac{\pi^{\frac{n-1}{2}}}{\Gamma(\lambda+1)}
\times
\begin{cases}
  1
\quad
&\text{if $1\le \ell \le n-1$}, 
\\
\lambda + \nu -n+1
\quad
&\text{if $\ell=n$}.  
\end{cases}
\end{multline*}
\end{enumerate}
\end{lemma}

\begin{proof}
(1)\enspace
This formula was given in \cite[Prop.~7.4]{sbon}, 
 but we give a proof here 
 in order to illustrate our notation
 for later purpose.  
By \eqref{eqn:Rabnew}, 
 the left-hand side amounts to 
\begin{align*}
&\frac{1}{\Gamma(\frac{\lambda+\nu-n+1}{2})\Gamma(\frac{\lambda-\nu}{2})} 
\int_0^{\infty} r^{\lambda-\nu-1} (1+r^2)^{-\lambda} d r
          \int_{S^{n-1}}|\omega_n|^{\lambda+\nu-n} d \omega
\\
=&
\frac{1}{2\Gamma(\frac{\lambda+\nu-n+1}{2})\Gamma(\frac{\lambda-\nu}{2})}
B(\frac{\lambda-\nu}{2},\frac{\lambda+\nu}2)S(\lambda+\nu-n,0), 
\end{align*}
which equals $\frac{\pi^{\frac {n-1}2}}{\Gamma(\lambda)}$
 by \eqref{eqn:Sab}.  
\par\noindent
(2)\enspace
By a similar computation as above,
 the ratio of the two integrals is given as 
\begin{equation*}
\frac{\text{the left-hand side of (2)}}
      {\text{the left-hand side of (1)}}
 =
 \frac{2 \int_0^{\infty} r^{\lambda-\nu-1} (1+r^2)^{-\lambda-1} d r
       \int_{S^{n-1}}|\omega_n|^{\lambda+\nu-n} |\omega_{\ell}|^2 d \omega}
      {\int_0^{\infty} r^{\lambda-\nu-1} (1+r^2)^{-\lambda} d r
       \int_{S^{n-1}}|\omega_n|^{\lambda+\nu-n} d \omega}.  
\end{equation*}
The right-hand side depends on 
 whether $\ell=n$ or not.  
It amounts to 
\begin{align*}
&\frac{2B(\frac{\lambda-\nu}{2}, \frac{\lambda+\nu}2+1)}
     {B(\frac{\lambda-\nu}{2}, \frac{\lambda+\nu}{2})}
\cdot
\frac{1}{S(\lambda+\nu-n,0)}
\times
\begin{cases}
  S(\lambda+\nu-n,2)
\\
  S(\lambda+\nu-n+2,0)
\end{cases}
\\
=&
\frac{\lambda+\nu}
     {\lambda}
\cdot
\frac{1}{\lambda+\nu}
\times
\begin{cases}
  1
\quad
&\text{if $1\le \ell \le n-1$}, 
\\
\lambda + \nu -n+1
\quad
&\text{if $\ell=n$}
\end{cases}
\end{align*}
by the recurrence relations
 \eqref{eqn:Sab2} and \eqref{eqn:Sab3}.  
\end{proof}

\begin{lemma}
\label{lem:1619112}
$A_I(\lambda,\nu)-D_I(\lambda,\nu)=\frac{\pi^{\frac {n-1}2}(\lambda-\nu)}{\Gamma(\lambda+1)}$.  
\end{lemma}
\begin{proof}
By the definitions \eqref{eqn:alphaI} and \eqref{eqn:deltaI}, 
 we have $\alpha_I(b) - \delta_I(b)=\dfrac{2|b|^2}{1+|b|^2}$.  
Thus we have
\begin{align*}
   A_I(\lambda,\nu)-D_I(\lambda,\nu)
   =&
   2 \int_{{\mathbb{R}}^{n}}
   \Atcal \lambda \nu + {}(b)
   (1+|b|^2)^{-\lambda-1}|b|^2 d b
\\
   =& \frac{B(\frac{\lambda-\nu}2,\frac{\lambda+\nu}2+1)
             S(\lambda+\nu-n,0)}
           {\Gamma(\frac{\lambda-\nu}{2})\Gamma(\frac{\lambda+\nu-n+1}{2})}, 
\end{align*}
as in the proof of Lemma \ref{lem:1616116} (1).  
Thus the lemma follows from \eqref{eqn:Sab}.  
\end{proof}

\begin{proof}
[Proof of Theorem \ref{thm:minKscalar}]
It follows from Lemma \ref{lem:1616116} that 
\[
  A_I(\lambda,\nu)
  =
  \frac{\pi^{\frac{n-1}{2}}}{\Gamma(\lambda+1)}
  \times 
\begin{cases}
\lambda-i
&
\text{if $n \not\in I$, }
\\
\lambda-(i-1)-(\lambda+\nu-n+1)
\quad
&
\text{if $n \in I$, }
\end{cases}
\]
whence the first two formul{\ae} of Theorem \ref{thm:minKscalar}
 are proved 
 by Proposition \ref{prop:1617113}.  

By Lemma \ref{lem:1619112}, 
 we have 
\begin{align*}
  D_I(\lambda,\nu)
  &=
  A_I(\lambda,\nu)
  -
  \frac{\pi^{\frac{n-1}{2}}(\lambda-\nu)}{\Gamma(\lambda+1)}
\\
  &=
  \frac{\pi^{\frac{n-1}{2}}}{\Gamma(\lambda+1)}
  \times
\begin{cases}
(\lambda-i)-(\lambda-\nu)
&
\text{if $n \not\in I$, }
\\
(n-\nu-i)-(\lambda-\nu)
\quad
&
\text{if $n \in I$, }
\end{cases}
\end{align*}
whence the last two formul{\ae}
 of Theorem \ref{thm:minKscalar}
 by Proposition \ref{prop:1617113}.  
\end{proof}

\begin{remark}
\label{rem:KKdual}
Alternatively,
 one could derive the last two formul{\ae}
 of Theorem \ref{thm:minKscalar} from the first two
 by using the duality theorem for symmetry breaking operators
 given in Proposition \ref{prop:SBOdual}.  
\end{remark}

\subsection{Proof of Theorem \ref{thm:153315}
 on the $(K,K')$-spectrum for $\Atbb{\lambda}{\nu}{-}{i,j}:
I_{\delta}(i,\lambda) \to J_{-\delta}(j,\nu)$}
\label{subsec:Apm1}

In this section,
 we determine the $(K,K')$-spectrum
 $b_-^{i,i-1}(\lambda,\nu)$ and  $c_-^{i,i}(\lambda,\nu)$
 in \eqref{eqn:SAabcd}
 for the normalized regular symmetry breaking operators
\index{A}{Ahtsln2@$\Atbb \lambda \nu {-} {i,j}$}
$
\Atbb \lambda {\nu}{-} {i,j}
\colon
I_{\delta}(i,\lambda) \to J_{-\delta}(j,\nu)
$
 with $j \in \{i-1, i\}$.  
By definition,
 these constants 
\index{A}{b2ii-1@$b_-^{i,i-1}(\lambda,\nu)$|textbf}
$
   b_-^{i,i-1}(\lambda,\nu)
$
 and $c_-^{i,i}(\lambda,\nu)$
 are characterized by the following equations:
\begin{alignat}{2}
\label{eqn:Abij}
\Atbb \lambda \nu - {i,i-1}
\circ
\iota_{\lambda}^{\ast}
=&
b_-^{i,i-1}(\lambda,\nu)
\iota_{\nu}^{\ast}
\circ
\pr ii
\quad
&&\text{on $\Exterior^{i}({\mathbb{C}}^{n+1})$, }
\\
\label{eqn:Acij}
\Atbb \lambda \nu - {i,i}
\circ
\iota_{\lambda}^{\ast}
=&
(-1)^{i}
c_-^{i,i}(\lambda,\nu)
\circ
\iota_{\nu}^{\ast}
\circ
\pr {i+1}i
\quad
&&\text{on $\Exterior^{i+1}({\mathbb{C}}^{n+1})$}.  
\end{alignat}
The main results of this section are given as follows:
\begin{theorem}
\label{thm:CApm1}
Suppose $\lambda, \nu \in {\mathbb{C}}$.  
Then we have
\begin{align}
\label{eqn:Cpmii1}
   b_-^{i,i-1}(\lambda,\nu)
   =&
   -
   \frac
   {2\pi^{\frac {n-1}{2}}}
   {\Gamma(\lambda+1)}, 
\\
\label{eqn:pmIi}
   c_-^{i,i}(\lambda,\nu)
   =&
   \frac
   {2 (-1)^{i+1} \pi^{\frac {n-1}{2}}}
   {\Gamma(\lambda+1)}.  
\end{align}
\end{theorem}

This is the remaining part
 of Theorem \ref{thm:153315}, 
 and the proof of Theorem \ref{thm:153315} will be complete
 when Theorem \ref{thm:CApm1} is shown.  
The proof of Theorem \ref{thm:CApm1} is parallel to 
 that of Theorem \ref{thm:minKscalar}, 
 and thus will be discussed briefly.  
We begin with an integral expression
 of the constants 
 $b_-^{i,i-1}(\lambda,\nu)$
 and 
 $c_-^{i,i}(\lambda,\nu)$ 
 as follows.  
\begin{proposition}
[integral expression of $(K,K')$-spectrum]
\label{prop:1617115}
\begin{align*}
b_-^{i,i-1}(\lambda,\nu)=& -2\int_{{\mathbb{R}}^n} \Atcal \lambda \nu - {}(b)
                           (1+|b|^2)^{-\lambda-1}b_n d b, 
\\
c_-^{i,i}(\lambda,\nu)=& 2(-1)^{i+1}\int_{{\mathbb{R}}^n} \Atcal \lambda \nu - {}(b)
                           (1+|b|^2)^{-\lambda-1}b_n d b.  
\end{align*}
\end{proposition}
Admitting Proposition \ref{prop:1617115} for the time being, 
 we complete the proof 
 of Theorem \ref{thm:CApm1}.  
\begin{proof}
[Proof of Theorem \ref{thm:CApm1}]
Theorem \ref{thm:CApm1} is an immediate consequence
 of Proposition \ref{prop:1617115}
 and the following lemma.  
\end{proof}

\begin{lemma}
\label{lem:Abnint}
\index{A}{Act2@$\Atcal \lambda \nu - {}$}
\[
\int_{{\mathbb{R}}^n} \Atcal \lambda \nu -{}(b) (1+|b|^2)^{-\lambda-1}b_n d b
=
\frac{\pi^{\frac{n-1}{2}}}{\Gamma(\lambda+1)}.  
\]
\end{lemma}
\begin{proof}
We use the identity
\[ 
   b_n \Atcal \lambda \nu -{}(b)=\Atcal {\lambda+1} \nu +{}(b).  
\]
Then the lemma follows from Lemma \ref{lem:1616116}.  
\end{proof}
The rest of this section is devoted to the proof
 of Proposition \ref{prop:1617115}.  
In the $N$-picture, 
 the equation \eqref{eqn:Abij} amounts to 
\[
\Atbb \lambda \nu {-} {i,i-1} ({\bf{1}}_{\lambda}^{{\mathcal{I}}})
=
\begin{cases}
b_-^{i,i-1}(\lambda,\nu) {h'}_{\nu}^{{\mathcal{I}}}
\qquad
&\text{if}\ n \notin {\mathcal{I}}, 
\\
0
\qquad
&\text{if}\ n \in {\mathcal{I}}, 
\end{cases}
\]
for all ${\mathcal{I}} \in {\mathfrak{I}}_{n+1,i}$, 
 whereas \eqref{eqn:Acij} amounts to 
\[
\Atbb \lambda \nu {-} {i,i} h_{\lambda}^{\mathcal{I}}
=
\begin{cases}
c_-^{i,i}(\lambda,\nu) {{\bf{1}}'}_{\nu}^{{\mathcal{I}} \setminus \{n\}}
\qquad
&\text{if}\ n \in {\mathcal{I}}, 
\\
0
\qquad
&\text{if}\ n \notin {\mathcal{I}}, 
\end{cases}
\]
for all ${\mathcal{I}} \in {\mathfrak{I}}_{n+1,i+1}$.  
In particular, 
 we have
\begin{alignat}{2}
\label{eqn:CApm2}
\Abb \lambda \nu {-} {i,i-1} 
({\bf{1}}_{\lambda}^{I \cup \{0\}})(0)
=&
b_-^{i,i-1}(\lambda,\nu) e_{I}
\quad
&&\text{for any $I \in {\mathfrak{I}}_{n-1,i-1}$}, 
\\
\label{eqn:CApm}
(\Atbb \lambda \nu -{i,i} h_{\lambda}^{I \cup \{n\}})(0)
=&
c_-^{i,i}(\lambda,\nu) e_I
\quad
&&\text{for any $I \in {\mathfrak{I}}_{n-1,i}$}
\end{alignat}
by \eqref{eqn:I10} and \eqref{eqn:hI0}
 because $0 \notin I$.

The distribution kernel $\Atcal \lambda \nu -{i,j}$
 of the regular symmetry breaking operator
 $\Atbb \lambda \nu -{i,j}$ is decomposed as
\[
   \Atcal \lambda \nu -{i,j} = \Atcal \lambda \nu -{} \Rij ij, 
\]
 where
\index{A}{Act2@$\Atcal \lambda \nu - {}$}
 $\Atcal \lambda \nu -{}$ is the scalar-valued, 
locally integrable function defined
 in \eqref{eqn:KAlnn-}
 and the matrix-valued function 
\index{A}{Rij@$\Rij ij$}
 $\Rij ij$ is defined in \eqref{eqn:Rij}.  
Then we have
\begin{align*}
(\Atbb \lambda \nu -{i,j} \psi)(0)
=&\int_{\mathbb{R}^n} \Atcal \lambda \nu -{}(-b) \Rij ij (-b)\psi(b) d b
\\
=&-\int_{\mathbb{R}^n} \Atcal \lambda \nu -{}(b) \Rij ij (b)\psi(b) d b
\end{align*}
in the $N$-picture for any $\psi \in \iota_{\lambda}^{\ast}({\mathcal{E}}^i(S^n))$.  
Hence Proposition \ref{prop:1617115} is a consequence
 of \eqref{eqn:CApm2}, \eqref{eqn:CApm}, 
 and of the following two lemmas.  

\begin{lemma}
\label{lem:TeIpm}
Suppose $I \in {\mathfrak{I}}_{n-1,i}$.  
Then the coefficient of $e_{I}$
 in 
$
   \Rij i {i-1}(b) {\bf{1}}_{\lambda}^{I \cup \{0\}}(b)
$
is equal to
\[
   2(1+|b|^2)^{-\lambda-1}b_n.  
\]
\end{lemma}

\begin{proof}
Using  the formula \eqref{eqn:minI}
 of ${\bf{1}}_{\lambda}^{\mathcal{I}}(b)$, 
 we have for $I \in {\mathfrak{I}}_{n,i}$
\begin{align*}
& (1+|b|^2)^{\lambda} \Rij i {i-1}(b) {\bf{1}}_{\lambda}^{I\cup\{0\}}(b)
\\
=& \Rij i{i-1}(b) 
   \sum_{J \in {\mathfrak{I}}_{n,i}}
   \det \psi_{n+1}(1,b)_{I\cup\{0\}, J} e_J
\\
=& \pr i {i-1} 
   \sum_{J \in {\mathfrak{I}}_{n,i}}
   \sum_{J' \in {\mathfrak{I}}_{n,i}}
   \det \psi_{n+1}(1,b)_{I\cup\{0\}, J} 
   \det \psi_{n}(b)_{J' J}
   e_{J'}
\\
=&(-1)^{i-1}\sum_{J \in {\mathfrak{I}}_{n,i}}
   \sum_{J' \in {\mathfrak{I}}_{n,i} \setminus {\mathfrak{I}}_{n-1,i}}
   \det \psi_{n+1}(1,b)_{I\cup\{0\}, J} 
   \det \psi_{n}(b)_{J' J}
   e_{J' \setminus \{n\}}.  
\end{align*}
Here, 
 for $J' \in {\mathfrak{I}}_{n,i}$, 
 we mean by $J' \not \in {\mathfrak{I}}_{n-1,i}$ the condition
 that $n \not \in J'$.  
Hence the coefficient of $e_I$
in 
$
   \Rij i {i-1}(b) {\bf{1}}_{\lambda}^{I \cup \{0\}}(b)
$
amounts to 
\[
(-1)^{i-1}\sum_{J \in {\mathfrak{I}}_{n,i}}
   \det \psi_{n+1}(1,b)_{I\cup\{0\}, J} 
   \det \psi_{n}(b)_{I\cup \{n\}, J}.  
\]
Now the lemma follows from Lemma \ref{lem:S0n}.  
\end{proof}

\begin{lemma}
\label{lem:eIpm}
Suppose $I \in {\mathfrak{I}}_{n-1,i}$.  
The coefficient of $e_I$
 in 
$
  \Rij ii (b) h_{\lambda}^{I \cup \{n\}}(b)
$
is given by 
\[
  2(-1)^i (1+|b|^2)^{-\lambda-1}b_n.  
\]
\end{lemma}

\begin{proof}
By \eqref{eqn:hIlmd} and \eqref{eqn:exrep}, 
we have
\begin{align*}
&\sigma(\psi_n(b)) h_{\lambda}^{I \cup \{n\}}(b)
\\
=&-(1+|b|^2)^{-\lambda}
  \sigma(\psi_n(b))
  \sum_{J \in {\mathfrak{I}}_{n,i}}
  \det \psi_{n+1}(1,b)_{I \cup \{n\}, J \cup \{0\}}e_J
\\
=&-(1+|b|^2)^{-\lambda}
  \sum_{J \in {\mathfrak{I}}_{n,i}}
  \sum_{J' \in {\mathfrak{I}}_{n,i}}
  \det \psi_{n+1}(1,b)_{I \cup \{n\}, J \cup \{0\}} 
  \det \psi_n(b)_{J' J}e_{J'}.  
\end{align*}
Applying the projection 
 $\pr ii:\Exterior^{i}({\mathbb{C}}^{n}) \to \Exterior^{i}({\mathbb{C}}^{n-1})$
 (see \eqref{eqn:Tii1}), 
 we find that the coefficient of $e_I$
 in $\Rij ii (b) h_{\lambda}^{I \cup \{n\}}(b)$
 is equal to 
\[
  -(1+|b|^2)^{-\lambda}
  \sum_{J \in {\mathfrak{I}}_{n,i}}
  \det \psi_{n+1}(1,b)_{I \cup \{n\}, J \cup \{0\}}
  \det \psi_n(b)_{I J}.  
\]
Hence the lemma follows from the minor summation formula
 in Lemma \ref{lem:Sn0}.  
\end{proof}

\subsection{Matrix-valued functional equations}
\label{subsec:FI}
The third goal of this chapter 
 is to obtain explicit matrix-valued functional equations
 for the regular symmetry breaking operators
 $\Atbb \lambda \nu \pm {i,j}$.  
We retain the setting 
 where $(G,G')=(O(n+1,1),O(n,1))$.  
By the 
\index{B}{genericmultiplicityonetheorem@generic multiplicity-one theorem}
generic multiplicity-one theorem
 (Theorem \ref{thm:unique}),
 two symmetry breaking operators from 
 the $G$-module $I_{\delta}(V,\lambda)$
 to the $G'$-module $J_{\varepsilon}(W,\nu)$
 must be proportional to each other
 if $[V:W]\ne0$ and $(\lambda,\nu,\delta,\varepsilon)$ does not belong 
 to the set $\Psising$
 of special parameters.  
In Sections \ref{subsec:FI} and \ref{subsec:161702},
 we consider the case 
\[
    (V,W)=
    (\Exterior^i({\mathbb{C}}^n), \Exterior^j({\mathbb{C}}^{n-1})),
\quad
    j \in \{i-1,i\}, 
\]
 and compare the (normalized) regular symmetry breaking operator
 $\Atbb \lambda \nu \gamma {i,j}$
 with its composition of the Knapp--Stein intertwining operator
 for $G$ or for the subgroup $G'$
 as in the following diagrams:
\begin{eqnarray*}
\xymatrix{
I_{\delta}(i,\lambda)\ar[rr]^{\Atbb \lambda \nu \gamma {i,j}}
\ar[drr]_{\Atbb \lambda{n-1-\nu}\gamma{i,j}}
&& J_{\varepsilon}(j,\nu)\ar[d]^{\Ttbb \nu {n-1-\nu}j}
\\
&&J_{\varepsilon}(j,n-1-\nu),
}
&&
\xymatrix{
I_{\delta}(i,\lambda) \ar[d]_{\Ttbb \lambda {n-\lambda}i}
\ar[drr]^{\Atbb \lambda \nu \gamma {i,j}}
&&
\\
I_{\delta}(i,n-\lambda)\ar[rr]_{\Atbb {n-\lambda}\nu\gamma{i,j}}
&&J_{\varepsilon}(j,\nu)
}
\end{eqnarray*}
where $\gamma = \delta \varepsilon$.  
We obtain closed formul{\ae}
 of the proportional constants 
 for the two operators in each diagram
 in Theorems \ref{thm:TAA} and \ref{thm:ATA}.  
The zeros of the proportional constants provide us 
 crucial information
 on the kernels and the images of the symmetry breaking operators
 $\Atbb \lambda \nu {\delta \varepsilon} {i,j} \colon
 I_{\delta}(i,\lambda) \to J_{\varepsilon}(j,\nu)$
 at reducible places
 of the principal series representations, 
 which will be investigated 
 in Chapter \ref{sec:pfSBrho}.

\subsubsection{Main results : Functional equations of $\Atbb \lambda \nu {\varepsilon} {i,j}$}
\label{subsec:fiA}
Suppose $j \in \{i-1, i\}$.  
Let $\Atbb \lambda \nu {\delta\varepsilon} {i,j} \colon I_{\delta}(i,\lambda)
 \to J_{\varepsilon}(j,\nu)$ be the normalized symmetry breaking operators
 as defined in \eqref{eqn:Aijdef}, 
 and $\Ttbb \nu {n-1-\nu} {j} \colon J_{\varepsilon}(j,\nu) \to J_{\varepsilon}(j,n-1-\nu)$
 be the normalized 
\index{B}{KnappSteinoperator@Knapp--Stein operator}
Knapp--Stein operators
 as defined in \eqref{eqn:KSii}
 for principal series representations
 of the subgroup $G'$.  
Then we obtain:
\begin{theorem}
[functional equation]
\label{thm:TAA}
Suppose $(\lambda,\nu) \in {\mathbb{C}}^2$
 and $\gamma \in \{ \pm \}$.  
Then 
\begin{alignat*}{2}
\Ttbb {\nu}{n-1-\nu}{i}
\circ
\Atbb \lambda \nu {\gamma} {i,i}
=&
\frac{\pi^{\frac{n-1}{2}} (\nu-i)}{\Gamma(\nu +1)}
\Atbb \lambda {n-1-\nu} {\gamma} {i, i}
\qquad
&&\text{for $0 \le i \le n-1$, }
\\
\Ttbb{\nu}{n-1-\nu}{i-1}
\circ
\Atbb \lambda \nu {\gamma} {i,i-1}
=&
\frac{\pi^{\frac{n-1}{2}} (n-\nu-i)}{\Gamma(\nu +1)}
\Atbb \lambda {n-1-\nu} {\gamma}{i,i-1}
\qquad
&&\text{for $1 \le i \le n$.  }  
\end{alignat*}
\end{theorem}
In the next theorem,
 we use the same letter $\Ttbb \lambda {n-\lambda} i$
 to denote the normalized Knapp--Stein intertwining operators
 $\Ttbb \lambda {n-\lambda} i \colon I_{\delta}(i,\lambda) \to I_{\delta}(i,n-\lambda)$
 for the group $G$.  
Then we obtain:
\begin{theorem}
[functional equation]
\label{thm:ATA}
\index{B}{functionalequation@functional equation}
Suppose $(\lambda,\nu) \in {\mathbb{C}}^2$
 and $\gamma \in \{ \pm \}$.  
Then 
\index{A}{Tiln@$\Ttbb \lambda{n-\lambda}{i}$|textbf}
\begin{alignat*}{2}
\Atbb {n-\lambda} \nu {\gamma} {i,i}
\circ
\Ttbb {\lambda}{n-\lambda}{i}
=&
\frac{\pi^{\frac{n}{2}} (n-\lambda-i)}{\Gamma(n-\lambda +1)}
\Atbb \lambda {\nu} {\gamma} {i,i}
\qquad
&&\text{for $0 \le i \le n-1$, }
\\
\Atbb {n-\lambda} \nu {\gamma} {i,i-1}
\circ
\Ttbb{\lambda}{n-\lambda}{i}
=&
\frac{\pi^{\frac{n}{2}} (\lambda-i)}{\Gamma(n-\lambda +1)}
\Atbb \lambda {\nu} {\gamma} {i,i-1}
\qquad
&&\text{for $1 \le i \le n$.}
\end{alignat*}
\end{theorem}

\begin{remark}
Theorems \ref{thm:TAA} and \ref{thm:ATA} generalize 
 the functional equations 
 which we proved in the scalar case
 \cite[Thm.~8.5]{sbon}.  
Matrix-valued functional identities
\index{B}{factorizationidentities@factorization identity}
 (factorization identities)
 for {\it{differential}} symmetry breaking operators
 were recently proved explicitly
 in \cite[Chap. 13]{KKP}.  
Alternatively,
 we could deduce a large part of the identities 
\cite[Chap.~13]{KKP} from Theorems \ref{thm:TAA}
 and \ref{thm:ATA}
 by using the residue formula 
 of the normalized symmetry breaking operators
 $\Atbb \lambda \nu {\pm} {i,j}$
 given in Fact \ref{fact:153316}, see \cite{xkresidue}.  
\end{remark}

\subsubsection{Proof of functional equations}
\label{subsec:fi}

In this section
 we give a proof 
 of the functional equations
 that are stated
 in Theorems \ref{thm:TAA} and \ref{thm:ATA}.

We apply Proposition \ref{prop:TminK}
 on the $K'$-spectrum 
 of the Knapp--Stein intertwining operator
 to the subgroup $G'=O(n,1)$.  
Then the $K'$-spectrum of the (normalized) Knapp--Stein intertwining operator
$
\Ttbb \nu{n-1-\nu}j : J_{\varepsilon}(j,\nu) \to J_{\varepsilon}(j,n-1-\nu)
$
 of $G'$ is given by 
\[
   \Ttbb \nu {n-1-\nu}i \circ \iota_{\nu}^{\ast}
  =
  c^{\natural}(j,\nu)' \iota_{n-\nu}^{\ast}
\qquad
\text{on $\mu^{\natural}(j,\varepsilon)'$}
\]
for $\natural =\flat$ or $\sharp$, 
 where 
\[
  c^{\flat}(j,\nu)'=\frac{(\nu-j) \pi^{\frac {n-1} 2}}{\Gamma(\nu+1)}, 
\qquad
c^{\sharp}(j,\nu)'=\frac{(n-1-j-\nu) \pi^{\frac {n-1} 2}}{\Gamma(\nu+1)}.   
\]

\begin{proof}
[Proof of Theorem \ref{thm:TAA}]
For $j=i$ or $j-1$
 and for $(\lambda, \nu)\in {\mathbb{C}}^2$
 with $\nu-\lambda \not\in {\mathbb{N}}$,
 we recall from Theorem \ref{thm:SBObasis}
 and Corollary \ref{cor:160150upper}
 that 
\[
  \operatorname{Hom}_{G'}(I_+(i,\lambda)|_{G'}, J_{\varepsilon}(j,n-1-\nu))
  =
  {\mathbb{C}} \Atbb \lambda {n-1-\nu}{\varepsilon} {i,j}.  
\]
Hence, 
 there exists a constant
$p_A ^{TA}(i,j,\varepsilon;\lambda, \nu) \in {\mathbb{C}}$
 such that
\begin{equation}
\label{eqn:TApA}
  \Ttbb \nu{n-1-\nu} j
  \circ
  \Atbb \lambda \nu \varepsilon {i,j}
  =
  p_A ^{TA}(i,j,\varepsilon;\lambda, \nu)
  \Atbb \lambda {n-1-\nu}{\varepsilon} {i,j}
\end{equation}
 if $n-1-\nu-\lambda \not\in {\mathbb{N}}$.  
We compute $p_A ^{TA}(i,j,\varepsilon;\lambda, \nu)$
 by using the $(K,K')$-spectrum $S_{\mu,\mu'}$
 (see Section \ref{subsec:Kspec})
 for \eqref{eqn:TApA}
 with an appropriate choice of basic $K$-types $\mu \in \widehat K$ and $\mu' \in \widehat K'$.  
We recall from Theorem \ref{thm:153315}
 an explicit formula
 of the $(K,K')$-spectrum
\[
S(\Atbb \lambda \nu \varepsilon {i,j})
=
\begin{pmatrix}
a_{\varepsilon}^{i,j}(\lambda,\nu)
&
b_{\varepsilon}^{i,j}(\lambda,\nu)
\\
c_{\varepsilon}^{i,j}(\lambda,\nu)
&
d_{\varepsilon}^{i,j}(\lambda,\nu)
\end{pmatrix}
\]
for the regular symmetry breaking operator
$\Atbb \lambda \nu \varepsilon {i,j}:
I_+(i,\lambda) \to J_{\varepsilon}(j,\nu)$
 with respect to basic $K$-types.  

\par\noindent
{\bf{Case 1.}}\enspace
$j=i$ and $\varepsilon=+$.  
Take $(\mu,\mu')=(\mub(i,+), \mub(i,+)')$. 
Then the computation of $S_{\mu,\mu'}$
 on the both sides of \eqref{eqn:TApA}
 leads us to the following identity:
\[
   p_A ^{TA}(i,i,+;\lambda, \nu)
   =
   c^{\flat}(i,\nu)'
   \cdot
   \frac{a_+^{i,i}(\lambda, \nu)}{a_+^{i,i}(\lambda, n-1-\nu)}
   =
   \frac{\pi^{\frac{n-1}{2}}(\nu-j)}{\Gamma(\nu+1)}
   \cdot 1.  
\]
\par\noindent
{\bf{Case 2.}}\enspace
$j=i$ and $\varepsilon=-$.  
Take $(\mu,\mu')=(\mus(i,+), \mub(i,-)')$. 
By the same argument as above,
 we have
\[
   p_A ^{TA}(i,i,-;\lambda, \nu)
   =
   c^{\flat}(i,\nu)'
   \cdot
   \frac{c_-^{i,i}(\lambda, \nu)}{c_-^{i,i}(\lambda, n-1-\nu)}
   =
   \frac{\pi^{\frac{n-1}{2}}(\nu-j)}{\Gamma(\nu+1)}
   \cdot 1.  
\]
\par\noindent
{\bf{Case 3.}}\enspace
$j=i-1$ and $\varepsilon=+$.  
Take $(\mu,\mu')=(\mus(i,+), \mus(i-1,+)')$. 
\[
   p_A ^{TA}(i,i-1,+;\lambda, \nu)
   =
   c^{\sharp}(i-1,\nu)'
   \cdot
   \frac{d_+^{i,i-1}(\lambda, \nu)}{d_+^{i,i-1}(\lambda, n-1-\nu)}
   =
   \frac{\pi^{\frac{n-1}{2}}(n-\nu-i)}{\Gamma(\nu+1)}
   \cdot 1.  
\]
\par\noindent
{\bf{Case 4.}}\enspace
$j=i-1$ and $\varepsilon=-$.  
Take $(\mu,\mu')=(\mub(i,+), \mus(i-1,-)')$. 
\[
   p_A ^{TA}(i,i-1,-;\lambda, \nu)
   =
   c^{\sharp}(i-1,\nu)'
   \cdot
   \frac{b_-^{i,i-1}(\lambda, \nu)}{b_-^{i,i-1}(\lambda, n-1-\nu)}
   =
   \frac{\pi^{\frac{n-1}{2}}(n-\nu-i)}{\Gamma(\nu+1)}
   \cdot 1.  
\]
Since both sides of \eqref{eqn:TApA} depend
 holomorphically in the entire $(\lambda,\nu)\in {\mathbb{C}}^2$, 
 the identity \eqref{eqn:TApA} holds
 for all $(\lambda,\nu)\in {\mathbb{C}}^2$.  
Hence Theorem \ref{thm:TAA} is proved.  
\end{proof}

\begin{proof}
[Proof of Theorem \ref{thm:ATA}]
The proof of Theorem \ref{thm:ATA} goes similarly.  
Since $\Atbb {n-\lambda}{\nu} \varepsilon{i,j} \circ
\Ttbb \lambda{n-\lambda}i$
$\in \operatorname{Hom}_{G'}(I_+(i,\lambda)|_{G'}, J_{\varepsilon}(j,\nu))$, 
 there exists a constant 
\[
   p_A^{AT}(i,j,\varepsilon;\lambda,\nu)\in {\mathbb{C}}
\]
 such that
\begin{equation}
\label{eqn:ATpA}
  \Atbb {n-\lambda} \nu \varepsilon {i,j}
  \circ
  \Ttbb \lambda {n-\lambda} i
  =
  p_A ^{AT}(i,j,\varepsilon;\lambda, \nu)
  \Atbb \lambda \nu \varepsilon {i,j}
\end{equation}
by the generic multiplicity-one theorem 
 (Theorem \ref{thm:genbasis})
for $j \in \{i-1,i\}$, 
 $\varepsilon \in \{\pm\}$, 
 and $(\lambda,\nu) \in {\mathbb{C}}^2$
 with $\nu-\lambda \not \in {\mathbb{N}}$.  
\par\noindent
{\bf{Case 1.}}\enspace
$j=i$ and $\varepsilon=+$.  
Take $(\mu,\mu')=(\mus(i,+), \mus(i,+)')$. 
Applying both sides of \eqref{eqn:ATpA}
 to the basic $K'$-type $\mu'=\mus(i,+)'$
 via the inclusion $\mu' \hookrightarrow \mu = \mus(i,+)$, 
 we get the following identities from 
 Proposition \ref{prop:TminK}
 and Theorem \ref{thm:minKscalar}:

\index{A}{cyfsharp@$c^{\sharp}(i,\lambda)$}
\index{A}{d1ij@$d_+^{i,j}(\lambda, \nu)$}
\[
   p_A ^{AT}(i,i,+;\lambda, \nu)
   =
   c^{\sharp}(i,\lambda)
   \cdot
   \frac{d_+^{i,i}(n-\lambda, \nu)}{d_+^{i,i}(\lambda, \nu)}
   =
   \frac{\pi^{\frac{n}{2}}(n-\lambda-i)}{\Gamma(\lambda+1)}
   \cdot 
   \frac{\Gamma(\lambda+1)}{\Gamma(n-\lambda+1)}.  
\]
The other three cases are proved similarly as below.  
\par\noindent
{\bf{Case 2.}}\enspace
$j=i$ and $\varepsilon=-$.  
Take $(\mu,\mu')=(\mus(i,+), \mub(i,-)')$. 
\index{A}{cz2ii@$c_-^{i,i}(\lambda, \nu)$}
\[
   p_A ^{AT}(i,i,-;\lambda, \nu)
   =
   c^{\sharp}(i,\lambda)
   \cdot
   \frac{c_-^{i,i}(n-\lambda, \nu)}{c_-^{i,i}(\lambda, \nu)}
   =
   \frac{\pi^{\frac{n}{2}}(n-\lambda-i)}{\Gamma(\lambda+1)}
   \cdot 
   \frac{\Gamma(\lambda+1)}{\Gamma(n-\lambda+1)}.  
\]
\par\noindent
{\bf{Case 3.}}\enspace
$j=i-1$ and $\varepsilon=+$.  
Take $(\mu,\mu')=(\mub(i,+), \mub(i-1,+)')$. 
\index{A}{cyflat@$c^{\flat}(i,\lambda)$}
\[
   p_A ^{AT}(i,i-1,+;\lambda, \nu)
   =
   c^{\flat}(i,\lambda)
   \cdot
   \frac{a_+^{i,i-1}(n-\lambda, \nu)}{a_+^{i,i-1}(\lambda, \nu)}
   =
   \frac{\pi^{\frac{n}{2}}(\lambda-i)}{\Gamma(\lambda+1)}
   \cdot 
   \frac{\Gamma(\lambda+1)}{\Gamma(n-\lambda+1)}.  
\]
\par\noindent
{\bf{Case 4.}}\enspace
$j=i-1$ and $\varepsilon=-$.  
Take $(\mu,\mu')=(\mub(i,+), \mus(i-1,-)')$. 
\[
   p_A ^{AT}(i,i-1,-;\lambda, \nu)
   =
   c^{\flat}(i,\lambda)
   \cdot
   \frac{b_-^{i,i-1}(n-\lambda, \nu)}{b_-^{i,i-1}(\lambda, \nu)}
   =
   \frac{\pi^{\frac{n}{2}}(\lambda-i)}{\Gamma(\lambda+1)}
   \cdot 
   \frac{\Gamma(\lambda+1)}{\Gamma(n-\lambda+1)}.  
\]
Thus Theorem \ref{thm:ATA} is proved.  
\end{proof}

\subsection{Renormalized symmetry breaking operator $\Attbb \lambda \nu + {i,j}$}
\label{subsec:161702}
In Theorem \ref{thm:170340}, 
 we constructed a {\it{renormalized}} symmetry breaking operator
 $\Attbb \lambda \nu \pm {V,W}$
 when the normalized regular symmetry breaking operator $\Atbb \lambda \nu \pm {V,W}$
vanishes.~We apply it 
 to the special case
 $(V,W)=(\Exterior^i({\mathbb{C}}^n), \Exterior^j({\mathbb{C}}^{n-1}))$, 
 and obtain for those $(\lambda,\nu)$
 for which $\Atbb \lambda \nu \gamma {i,j}=0$
 the renormalized symmetry breaking operator
 $\Attbb \lambda \nu \gamma {i,j}$
 as the analytic continuation of the following:
\begin{equation}
\label{eqn:Aijrenorm}
\Attbb \lambda \nu \gamma {i,j}
=
\begin{cases}
\Gamma(\frac{\lambda-\nu}{2})\Atbb \lambda \nu + {i,j}
\quad
&\text{if $\gamma=+$}, 
\\
\Gamma(\frac{\lambda-\nu+1}{2})\Atbb \lambda \nu - {i,j}
\quad
&\text{if $\gamma=-$}.
\end{cases}
\end{equation}  
We recall 
 that for $j \in \{i-1,i\}$
 and $\gamma \in \{\pm\}$, 
 we have determined in Theorem \ref{thm:161243}
 precisely the zero set
\[
   \{
    (\lambda, \nu)\in {\mathbb{C}}^2
    :   
    \Atbb \lambda \nu \gamma {i,j}=0
   \}.  
\]
In this section,
 we discuss functional equations
 and $(K,K')$-spectrum
 of the 
\index{B}{renormalized regular symmetry breaking operator@regular symmetry breaking operator, renormalized---}
renormalized operators $\Attbb \lambda \nu \pm{i,j}$
 only in the few cases
 that are necessary for later arguments.

\subsubsection{Functional equations
 for the renormalized operator $\Attbb \lambda i + {i,i}$}
In this subsection,
 we treat the case $j=i$.  
For $\nu=i(=j)$, 
 $\Atbb \lambda i \gamma {i,i}=0$
 if and only if $\lambda=i \in \{0,1,\cdots,n-1\}$
 and $\gamma=+$ by Theorem \ref{thm:161243}.  
Then the renormalized operator 
$
   \Attbb \lambda i + {i,i} \colon I_{\delta}(i,\lambda) \to J_{\delta}(i,i)
$
 is the analytic continuation 
 of the following:
\begin{equation}
\label{eqn:Aii2tilde}
   \Attbb \lambda i + {i,i}=\Gamma(\frac{\lambda-i}2)\Atbb \lambda i + {i,i}.  
\end{equation}
Then $\Attbb \lambda i + {i,i}\colon I_{\delta}(i,\lambda) \to J_{\delta}(i,i)$
 is a $G'$-homomorphism
 that depends holomorphically on $\lambda$
 in the entire complex plane ${\mathbb{C}}$
 by Theorem \ref{thm:170340} (3).

We determine 
\index{B}{functionalequation@functional equation}
 functional equations
 and 
\index{B}{KspectrumKprime@$(K,K')$-spectrum}
$(K,K')$-spectrum 
 $S(\Attbb \lambda i + {i,i})$
 (see \eqref{eqn:Smat})
 on basic $K$- and $K'$-types
 for the renormalized operator
\index{A}{Ahttsln1@$\Attbb \lambda \nu {+} {i,j}$}
 $\Attbb \lambda i + {i,i}$ as follows.

\begin{lemma}
[functional equations and the $(K,K')$-spectrum
 for $\Attbb {n-\lambda} i + {i,i}$]
\label{lem:Aiiren}
~~~\newline
Suppose $0 \le i \le n-1$
 and $\lambda \in {\mathbb{C}}$.  
Then we have 
\begin{align}
\label{eqn:TA2tilde}
\Ttbb i {n-1-i} i \circ \Attbb {\lambda} i + {i,i} 
=&
0, 
\\
\label{eqn:ATAi2}
\Attbb {n-\lambda} i + {i,i} 
\,\,\,
\circ \Ttbb \lambda {n-\lambda} {i}
=&
\frac{2 \pi^{\frac n 2}\Gamma(\frac {n-\lambda-i}2+1)}
     {\Gamma(\frac{\lambda-i}2) \Gamma(n-\lambda+1)} \Attbb \lambda i + {i,i}, 
\\
\label{eqn:SAiitilde}
S(\Attbb \lambda i + {i,i}) \,=& \frac{\pi^{\frac {n-1} 2}}{\Gamma(\lambda+1)}
\begin{pmatrix} 2 & 0 \\ 0 & 0 \end{pmatrix}.  
\end{align}
\end{lemma}

\begin{proof}
Applying Theorem \ref{thm:TAA}
 with $\nu =i$ ($0 \le i \le n$), 
 we have
\[
   \Ttbb i {n-1-i} i \circ \Atbb {\lambda} i + {i,i} 
   =0 
\quad \text{for all $\lambda \in {\mathbb{C}}$.  }
\]
Taking the limit as $\lambda$ tends to $i$
 in the following equation:
\[
   \Ttbb i {n-1-i} i \circ 
   \Gamma(\frac{\lambda-i}{2})
   \Atbb {\lambda} i + {i,i} =0, 
\]
we get the desired formula \eqref{eqn:TA2tilde}
 by the definition \eqref{eqn:Aii2tilde} of the renormalization 
 $\Attbb {\lambda} i + {i,i}$.

Similarly,
 the formul{\ae} \eqref{eqn:ATAi2} and \eqref{eqn:SAiitilde}
 for the renormalized operator
 $\Attbb {\lambda} i + {i,i}$
 follow from the limit of the corresponding results
 for $\Atbb {\lambda} i + {i,i}$ given in Theorems \ref{thm:ATA}
 and \ref{thm:153315}, 
 respectively.  
\end{proof}

\subsubsection{Functional equations at middle degree for $n$ even}

For $n$ even (say, $n=2m$), 
 at the \lq\lq{middle degree}\rq\rq\ $i=\frac n 2( =m$), 
 we observe
 that the Knapp--Stein operator
$\Ttbb \lambda {2m-\lambda}{m} \colon I_+(m,\lambda) \to I_+(m,2m-\lambda)$
 vanishes 
 if $\lambda =m$
 (see Proposition \ref{prop:Tvanish}), 
 and so the functional equation \eqref{eqn:ATAi2} is trivial.  
Instead we use the
\index{B}{KnappSteinoperatorrenorm@ Knapp--Stein operator, renormalized---}
 renormalized  Knapp--Stein operator
\index{A}{TVttiln@$\Tttbb \lambda{n-\lambda}{\frac n 2}$, renormalized Knapp--Stein intertwining operator}
 $\Tttbb \lambda {2m-\lambda} {m}$
 defined in  \eqref{eqn:Ttilde}
 for another functional equation, 
 see \eqref{eqn:AT2tilde} below.  
We recall from Lemma \ref{lem:161745}
 that $\Tttbb \lambda{2m-\lambda}m$ is an endomorphism
 of $I_{\delta}(m,m)$ 
 when $\lambda=m$, 
 but is not proportional
 to the identity operator
 when $\lambda=m$.  

\begin{lemma}
[functional equation for $\Attbb m m + {i,i}$]
\label{lem:AmmT}
\index{B}{functionalequation@functional equation}
Let $(G,G')=(O(n+1,1),O(n,1))$ with $n=2m$.  
Then we have
\begin{equation}
\label{eqn:AT2tilde}
\Attbb {m}{m} + {m, m}
\circ
\Tttbb m m m
=
\frac{\pi^m}{m!} 
\Attbb {m}{m} + {m, m}.  
\end{equation}
\end{lemma}

\begin{proof}
By Theorem \ref{thm:170340} and \eqref{eqn:Ttilde}, 
\begin{align}
\Attbb {m}{m} + {m, m}
\circ
\Tttbb {m}{m} {m}
=&
(\lim_{\lambda \to m}
  \Gamma(\frac{(2m-\lambda)-m}{2})
  \Atbb {2m-\lambda}{m} + {m, m})
\circ
(\lim_{\lambda \to m}
  \frac{1}{\lambda - m}
   \Ttbb {\lambda}{2m-\lambda} {m})
\notag
\\
=&
\lim_{\lambda \to m}
  \frac{\Gamma(\frac{m-\lambda}{2})}
       {\lambda-m}
  \Atbb {2m-\lambda}{m} + {m, m}
  \circ
  \Ttbb {\lambda}{2m-\lambda} {m}.  
\notag
\end{align}
In turn, the functional equation
 in Theorem \ref{thm:ATA} shows that 
 the right-hand side amounts to
\begin{align}
\lim_{\lambda \to m}
  \frac{\Gamma(\frac{m-\lambda}{2})}
       {\lambda-m}
  \frac{\pi^{m} (2m-\lambda-m)}{\Gamma(2m-\lambda+1)}
  \Atbb {\lambda}{m} + {m, m}
=&
 \frac{-\pi^{m}}{\Gamma(m+1)}
 (\lim_{\lambda \to m} 
  \frac{\Gamma(\frac{m - \lambda}{2})}{\Gamma(\frac{\lambda-m}{2})})\Attbb {m}{m} + {m, m}
\notag
\\
=&
 \frac{\pi^{m}}{\Gamma(m+1)}
 \Attbb {m}{m} + {m, m}.  
\end{align}
Hence the formula \eqref{eqn:AT2tilde} is proved.  
\end{proof}

In contrast to Lemma \ref{lem:AmmT}
 where we needed to treat the renormalized operator $\Attbb {m}{m} + {m, m}$
 because $\Atbb {m}{m} + {m, m}=0$, 
 the normalized operator $\Atbb {m}{m} - {m, m}$ does not vanish
 (Theorem \ref{thm:161243} (3)).  
In this case,
 the functional equations for $\Atbb {m}{m} - {m, m}$ are given as follows:

\begin{lemma}
[functional equation for $\Atbb m m - {m, m}$]
\label{lem:AmmT-}
We retain the setting that $(G,G')=(O(n+1,1),O(n,1))$
 with $n=2m$.  
Then we have 
\begin{align}
\label{eqn:ATAm-}
\Atbb m m - {m,m}\circ \Tttbb m m m = & - \frac {\pi^m}{m!}\Atbb m m - {m,m}, 
\\
\label{eqn:TATm-}
\Ttbb m {m-1} m  \circ \Atbb m m - {m,m} = & 0.  
\end{align}
\end{lemma}

\begin{proof}
By the definition of $\Tttbb \lambda{2m-\lambda}{m,m}$
 in \eqref{eqn:Ttilde}
 and the functional equation in Theorem \ref{thm:ATA}, 
 we have
\begin{align*}
\Atbb m m - {m,m}\circ \Tttbb m m m 
= & \lim_{\lambda \to m} \Atbb {2m-\lambda} m - {m,m} 
    \circ \frac {1}{\lambda-m} \Ttbb \lambda {2m-\lambda} m
\\
= & \lim_{\lambda \to m} \frac {\pi^m (m-\lambda)}{(\lambda-m)\Gamma(2m-\lambda+1)}
  \Atbb {\lambda} m - {m,m} 
\\ 
=&\frac {-\pi^m}{m!} \Atbb m m - {m,m}.  
\end{align*}
Hence the first statement is verified.  
The second statement is a special case
 of Theorem \ref{thm:TAA}.  
\end{proof}

\subsubsection{Functional equations for the renormalized operator
 $\Attbb \lambda {n-i} + {i,i-1}$}
In this subsection, 
 we treat the case $j=i-1$.  
For $j=i-1$ and $\nu=n-i$, 
 $\Atbb \lambda {n-i} \gamma {i,j} =0$
 if and only if $\gamma=+$ and $\lambda =n-i$
 by Theorem \ref{thm:161243}.  
In this case, 
 the renormalized symmetry breaking operator
 $\Attbb \lambda {n-i} + {i,i-1} \colon I_{\delta}(i,\lambda) \to J_{\delta}(i-1, n-i)$ is obtained 
 as the analytic continuation of the following:
\index{A}{Ahttsln1@$\Attbb \lambda \nu {+} {i,j}$}
\[
   \Attbb \lambda {n-i} + {i,i-1}=\Gamma(\frac{\lambda-n+i}2)\Atbb \lambda {n-i} + {i,i-1}, 
\]
see Theorem \ref{thm:170340} (3).

We determine functional equations
 and $(K,K')$-spectrum $S(\Attbb \lambda {n-i}+{i,i-1})$
 (see \eqref{eqn:Smat})
 on basic $K$- and $K'$-types
 for the renormalized operator $\Attbb \lambda {n-i}+{i,i-1}$
 as follows.  
 
\begin{lemma}
[functional equations and the $(K,K')$-spectrum
 for $\Attbb \lambda {n-i}+{i,i-1}$]
\label{lem:Aii-ren}
~~~\newline
Suppose $1 \le i \le n$ and $\lambda \in {\mathbb{C}}$.  
Then we have 
\begin{align}
\label{eqn:1802119}
\Attbb {n-\lambda} {n-i} + {i,i-1} \circ \Ttbb \lambda {n-\lambda} {i}
=&
-\frac{2\pi^{\frac n 2}\Gamma(\frac{i-\lambda}2+1)}{\Gamma(n-\lambda+1)\Gamma(\frac{\lambda-n+i}{2})}
 \Attbb \lambda {n-i} + {i,i-1}, 
\\
\label{eqn:180294}
\Ttbb {n-i}{i-1} {i-1}\circ
\Attbb {\lambda} {n-i} + {i,i-1} 
=&
0, 
\\
S(\Attbb \lambda {n-i} + {i,i-1}) =& \frac{\pi^{\frac {n-1} 2}}{\Gamma(\lambda+1)}
\begin{pmatrix} 0 & 0 \\ 0 & 2 \end{pmatrix}.  
\notag
\end{align}
\end{lemma}

\begin{proof}
The functional equations follow from Theorems \ref{thm:TAA} and \ref{thm:ATA}.  The formula for the $(K,K')$-spectrum
 is derived from Theorem \ref{thm:153315}.  
\end{proof}

\subsubsection{Functional equations at middle degree for $n$ odd}
For $n$ odd 
 (say, $n=2m+1$), 
 the Knapp--Stein operator $\Ttbb \nu {n-\nu-1} j \colon J_{\varepsilon}(j,\nu) \to J_{\varepsilon}(j,n-1-\nu)$ 
 for the subgroup $G'=O(n,1)$ vanishes
 at the middle degree $j=\frac 1 2 (n-1) (=m)$ 
 if $\nu=m$ by Proposition \ref{prop:Tvanish}.  
We note that the exact sequence in Theorem \ref{thm:LNM20} (1)
 for $G'=O(2m+1,1)$ splits, 
 and we have a direct sum decomposition 
\[
  J_{\varepsilon}(m,m) \simeq \pi_{m,\varepsilon} \oplus \pi_{m+1,-\varepsilon}
\]
 of two irreducible tempered representations of $G'$.  
In this case,
 the functional equations \eqref{eqn:TA2tilde}
 in Lemma \ref{lem:Aiiren}
 and \eqref{eqn:180294} in Lemma \ref{lem:Aii-ren} are trivial,
 and we replace them
 by the following functional equations
 for the {\it{renormalized}} Knapp--Stein operator
\index{A}{TVttiln@$\Tttbb \lambda{n-\lambda}{\frac n 2}$, renormalized Knapp--Stein intertwining operator}
 $\Tttbb m m m$.  

\begin{lemma}
\label{lem:TAT3tilde}
For $(G,G')=(O(n+1,1),O(n,1))$ with $n=2m+1$
 and for $\lambda \in {\mathbb{C}}$, 
 we have
\begin{align}
\label{eqn:TAA3tilde}
   \Tttbb {m}{m} {m} \circ \Attbb {\lambda}{m} + {m, m}
   =&
   \frac{\pi^m}{m!} \Attbb {\lambda}{m} + {m, m}, 
\\
\label{eqn:180293}
   \Tttbb {m}{m} {m} \circ \Attbb {\lambda}{m} + {m+1, m}
   =&
   -\frac{\pi^m}{m!} \Attbb {\lambda}{m} + {m+1, m}.  
\end{align}
\end{lemma}

Lemma \ref{lem:TAT3tilde} tells that 
\begin{alignat*}{3}
&{\operatorname{Image}}
(\Attbb {\lambda}{m} + {m, m}\colon I_{\delta}(m,\lambda) \to J_{\delta}(m,m)) 
&&\subset \pi_{m,\delta}, 
\\
&{\operatorname{Image}}
(\Attbb {\lambda}{m} + {m+1, m}\colon I_{\delta}(m+1,\lambda) 
\to J_{\delta}(m,m)) 
&&\subset \pi_{m+1,-\delta}, 
\end{alignat*}
for all $\lambda \in {\mathbb{C}}$
by Lemma \ref{lem:161745}.

\begin{proof}
The functional equations in Theorem \ref{thm:TAA} tell
 that 
\begin{align*}
(\frac{1}{\nu-m} \Ttbb \nu {2m-\nu} m) 
\circ 
\Gamma(\frac{\lambda-m}{2}) \Atbb \lambda \nu + {m,m} 
=& 
\frac{\pi^m}{\Gamma(\nu+1)} \Gamma(\frac{\lambda-m}{2}) 
\Atbb \lambda {2m-\nu} + {m,m}, 
\\
(\frac{1}{\nu-m} \Ttbb \nu {2m-\nu} m) 
\circ 
\Gamma(\frac{\lambda-m}{2}) \Atbb \lambda \nu + {m+1,m} 
=& 
-\frac{\pi^m}{\Gamma(\nu+1)} \Gamma(\frac{\lambda-m}{2}) 
\Atbb \lambda {2m-\nu} + {m+1,m}.  
\end{align*}
Taking the limit as $\nu$ tends to $m$, 
 we get Lemma \ref{lem:TAT3tilde}.  
\end{proof}

\subsection{Restriction map $I_{\delta}(i,\lambda) \to J_{\delta}(i,\lambda)$}
\label{subsec:Rest}
The restriction of (smooth) differential forms
 to a submanifold defines an obvious continuous map
 between Fr{\'e}chet spaces,
 which intertwines the conformal representation
 (see \cite[Lem.~8.9]{KKP}).  
We end this section with the most elementary symmetry breaking operator,
 namely,
 the restriction map
 for the pair $(G'/P', G/P) \subset (S^{n-1},S^n)$.  
\begin{lemma}
\label{lem:152268}
The restriction map from $G/P$ to the submanifold $G'/P'$
 induces obvious symmetry breaking operators
\[
{\operatorname{Rest}}_{\lambda,\lambda,+}^{i,i}
:
I_{\delta}(i,\lambda) \to J_{\delta}(i,\lambda).  
\]
Then the $(K,K')$-spectrum
 for basic $K$- and $K'$-types (see \eqref{eqn:Smat})
 is given by
\begin{equation}
\label{eqn:SRest}
  S({\operatorname{Rest}}_{\lambda,\lambda,+}^{i,i})
  =
  \begin{pmatrix} 1 & 0 \\ 0 & 1\end{pmatrix}.
\end{equation}
\end{lemma}
\begin{proof}
We recall from Proposition \ref{prop:minKN}
 that 
\[
   \{ {\bf{1}}_{\lambda}^{\mathcal{I}}
    :
   {\mathcal{I}} \in {\mathfrak {I}}_{n+1,i}
   \}
\]
forms a basis of the basic $K$-type $\mub(i,\delta)$
 of the principal series representation $I_{\delta}(i,\lambda)$
 in the $N$-picture.  
Let $\mathcal{I} \in {\mathfrak {I}}_{n+1,i}$
 and $(x,x_n) \in {\mathbb{R}}^{n-1} \oplus {\mathbb{R}}={\mathbb{R}}^n$.  
By \eqref{eqn:minI}, 
 we have
\[
{\bf{1}}_{\lambda}^{\mathcal{I}}(x,x_n)
=-(1+|x|^2+x_n^2)^{-\lambda-1}
\sum_{J \in {\mathfrak{I}}_{n,i}} S_{{\mathcal{I}}J}(1,x,x_n)e_J.  
\]
Then an elementary computation by using \eqref{eqn:SIJ} shows
\begin{equation*}
{\bf{1}}_{\lambda}^{\mathcal{I}}(x,0)
=
\begin{cases}
{{\bf{1}}_{\lambda}'}^{\mathcal{I}}(x)
\qquad
&\text{if }n \notin {\mathcal{I}} \in {\mathfrak {I}}_{n+1,i}, 
\\
{{\bf{1}}_{\lambda}'}^{\mathcal{I}\setminus \{n\}}(x) \wedge e_n
\qquad
&\text{if }n \in {\mathcal{I}} \in {\mathfrak {I}}_{n+1,i}.  
\end{cases}
\end{equation*}
The case for the basic $K$-type
 $\mus(i,\delta)$ is similar, 
 where we recall from Proposition \ref{prop:minKN}
 that $\{h_{\lambda}^{\mathcal{I}}: {\mathcal{I}} \in {\mathfrak {I}}_{n+1,i+1}\}$ forms its basis, 
 for which we can compute the restriction $x_n=0$.  
Thus Lemma \ref{lem:152268} is shown.  
\end{proof}

\subsection{Image of the differential symmetry breaking operator
 $\Ctbb \lambda \nu {i,j}$}
\label{subsec:ImageC}

In Theorem \ref{thm:imgDSBO}, 
 we have proved 
 that the image of any nonzero differential symmetry breaking operator from
 principal series representation
 is infinite-dimensional.  
As an application of the functional equations
 of the (generically) regular symmetry breaking operators
 $\Atbb \lambda \nu \pm {i,j}$
 (Theorems \ref{thm:TAA} and \ref{thm:ATA})
 and of the residue formul{\ae}
 of $\Atbb \lambda \nu \pm {i,j}$
 (Fact \ref{fact:153316},  see \cite{xkresidue}), 
 we end this chapter
 with a necessary and sufficient condition
 for the renormalized differential symmetry breaking operator
 $\Ctbb \lambda \nu {i,j}$ to be surjective
 when $j=i$, $i-1$, 
see Theorems \ref{thm:180380} and \ref{thm:180386}.

\subsubsection{Surjectivity condition of $\Ctbb \lambda \nu {i,j}$}
Suppose $j \in \{i,i-1\}$.  
We recall from \eqref{eqn:Ciitilde} and \eqref{eqn:Cii-1tilde}
 that the renormalized differential symmetry 
 breaking operator 
$\Ctbb \lambda \nu {i,j}\colon I_{\delta}(i,\lambda)
 \to J_{\varepsilon}(j,\nu)$
 is defined 
 for $(\lambda, \nu) \in {\mathbb{C}}^2$
 with $\nu-\lambda \in {\mathbb{N}}$
 and $\delta \varepsilon =(-1)^{\nu-\lambda}$.  
Moreover, 
 $\Ctbb \lambda \nu {i,j}$ is nonzero
 for any $(i,j, \lambda,\nu)$
 with $j \in \{i,i-1\}$
 and $\nu-\lambda \in {\mathbb{N}}$.

In what follows, 
 we shall sometimes encounter the condition that 
\index{A}{Leven@$L_{\operatorname{even}}$|textbf}
\index{A}{Lodd@$L_{\operatorname{odd}}$|textbf}
$(\lambda,n-1-\nu) \in L_{\operatorname{even}} \cup L_{\operatorname{odd}}$, 
 which is equivalent to 
\begin{equation}
\label{eqn:he}
(\lambda,\nu) \in {\mathbb{Z}}^2
\quad
\text{and}
\quad
\lambda+\nu \le n-1 \le \nu.  
\end{equation}
\begin{theorem}
\label{thm:180380}
Suppose $0 \le i \le n-1$, 
 $\nu-\lambda \in {\mathbb{N}}$, 
 and $\delta$, $\varepsilon \in \{\pm\}$
 with $(-1)^{\nu-\lambda} = \delta \varepsilon$.  
Then the following two conditions (i) and (ii) on $(i, \lambda,\nu)$
 are equivalent:
\begin{enumerate}
\item[{\rm{(i)}}]
\index{A}{Ctiiln@$\Ctbb \lambda \nu {i,i}$}
$
   \Ctbb \lambda \nu {i,i} \colon I_{\delta}(i,\lambda) \to J_{\varepsilon}(i,\nu)
$
 is not surjective.  
\item[{\rm{(ii)}}]
One of the following conditions holds:
\begin{enumerate}
\item[{\rm{(ii-a)}}]
$1 \le i \le n-1$, 
$\nu=i$, 
 and ${\mathbb{Z}} \ni \lambda < i;$
\item[{\rm{(ii-b)}}]
$n$ is odd, 
 $i=0$, 
 and \eqref{eqn:he}$;$
\item[{\rm{(ii-c)}}]
$n$ is odd,
 $1 \le i \le n-1$, 
 \eqref{eqn:he}, 
 and $\nu \ne n-1;$
\item[{\rm{(ii-d)}}]
$n$ is odd,
 $1 \le i < \frac 1 2 (n-1)$, 
 $(\lambda, \nu)=(i,n-1-i)$.  
\end{enumerate}
\end{enumerate}
\end{theorem}

\begin{theorem}
\label{thm:180386}
Suppose $1 \le i \le n$, 
 $\nu-\lambda \in {\mathbb{N}}$, 
 and $\delta$, $\varepsilon \in \{\pm\}$
 with $(-1)^{\nu-\lambda} = \delta \varepsilon$.  
Then the following two conditions (i) and (ii)
 on $(i, \lambda,\nu)$ are equivalent:
\begin{enumerate}
\item[{\rm{(i)}}]
\index{A}{Ctiim1ln@$\Ctbb \lambda \nu {i,i-1}$}
$
   \Ctbb \lambda \nu {i,i-1} \colon I_{\delta}(i,\lambda) \to J_{\varepsilon}(i-1,\nu)
$
 is not surjective.  
\item[{\rm{(ii)}}]
One of the following conditions holds:
\begin{enumerate}
\item[{\rm{(ii-a)}}]
$1 \le i \le n-1$, $\nu=n-i$, and ${\mathbb{Z}} \ni \lambda<n-i;$
\item[{\rm{(ii-b)}}]
$n$ is odd,
 $1 \le i \le n-1$, 
 \eqref{eqn:he}, 
 and $\nu \ne n-1;$
\item[{\rm{(ii-c)}}]
$n$ is odd,
 $i=n$, 
 and \eqref{eqn:he}$;$
\item[{\rm{(ii-d)}}]
$n$ is odd,
 $\frac 1 2 (n+1) < i \le n-1$, 
 and $(\lambda,\nu)=(n-i,i-1)$.  
\end{enumerate}
\end{enumerate}
\end{theorem}

For the proof of Theorems \ref{thm:180380} and \ref{thm:180386}, 
 we first derive
 the functional equations for $\Ctbb \lambda \nu {i,j}$
 in Theorem \ref{thm:180389} from those
 for the regular symmetry breaking operators
 $\Atbb \lambda \nu \pm{i,j}$
 in Chapter \ref{sec:holo}
 and from the matrix-valued residue formul{\ae} \cite{xkresidue}.  
The results cover most of the cases
 where the Knapp--Stein intertwining operators
 $\Ttbb \nu {n-1-\nu}j$ do not vanish.  
A special attention is required
 when $\Ttbb \nu {n-1-\nu}j=0$.  
In this case,
 the principal series representation $J_{\varepsilon}(j,\nu)$ splits into the direct sum
 of two irreducible representations
 of the subgroup $G'=O(n,1)$, 
 and we shall treat this case separately 
 in Section \ref{subsec:10.6.3}.  
The proof of Theorems \ref{thm:180380} and \ref{thm:180386} will be completed 
 in Section \ref{subsec:10.6.4}.  

\subsubsection{Functional equation for $\Ctbb \lambda \nu {i,j}$}
Suppose $0 \le i \le n$, 
 $0 \le j \le n-1$,
 and $j=i$ or $i-1$.  
We set $p^{i,j}(\lambda,\nu)$ by 
\begin{align}
\label{eqn:pijlmdnu}
p^{i,i}(\lambda,\nu)
&:=
\begin{cases}
1 
&\text{if $i=0$ or $\lambda = \nu$}, 
\\
\frac 1 2(\nu-i)
\qquad\quad
&\text{otherwise}, 
\end{cases}
\\
\notag
p^{i,i-1}(\lambda,\nu)
&:=
\begin{cases}
1 
\quad
&\text{if $i=n$ or $\lambda = \nu$}, 
\\
\frac 1 2(\nu+i-n)
\quad
&\text{otherwise}.  
\end{cases}
\end{align}

\begin{theorem}
[functional equation for $\Ctbb \lambda \nu {i,j}$]
\label{thm:180389}
Suppose $0 \le i \le n$, 
 $0 \le j \le n-1$, 
 and $j \in \{i,i-1\}$.  
For $(\lambda, \nu) \in {\mathbb{C}}^2$
 with $\nu-\lambda \in {\mathbb{N}}$, 
 we have
\begin{equation}
\label{eqn:180389}
  \Ttbb \nu {n-1-\nu} {j} \circ \Ctbb \lambda \nu {i,j}
=
  q(\nu-\lambda) p^{i,j}(\lambda,\nu) \Atbb \lambda {n-1-\nu} {(-1)^{\nu-\lambda}}{i,j}, 
\end{equation}
 where 
\index{A}{qm@$q(m)$}
 $q(m)$ is a nonzero number defined in \eqref{eqn:qm}.  
\end{theorem}

\begin{proof}
We set
\begin{equation}
\label{eqn:qijnu}
p_{i,j}(\nu)
:=
\begin{cases}
\frac 12(\nu-i) 
\quad
&\text{if $j=i$}, 
\\
\frac 1 2(n-\nu-i)
\quad
&\text{if $j=i-1$}.  
\end{cases}
\end{equation}
By the functional equation 
 for the regular symmetry breaking operator
$\Atbb \lambda \nu \pm {i,j}$ 
given in Theorem \ref{thm:TAA}, 
 we have for $\gamma \in \{\pm\}$
\[
  \Ttbb \nu{n-1-\nu} j \circ \Atbb \lambda \nu \gamma {i,j}
= \frac{2 \pi^{\frac{n-1}{2}} p_{i,j}(\nu)}{\Gamma(\nu+1)}
 \Atbb \lambda {n-1-\nu} \gamma {i,j}.  
\]
Suppose $\nu-\lambda \in {\mathbb{N}}$.  
Applying the residue formula \eqref{eqn:qAC}
 of $\Atbb \lambda \nu {\pm}{i,j}$ to the left-hand side,
 we get 
\begin{equation}
\label{eqn:TCqpAij}
  \Ttbb \nu{n-1-\nu} j \circ \Cbb \lambda \nu {i,j}
= (-1)^{i-j} q(\nu-\lambda) p_{i,j}(\nu)
  \Atbb \lambda {n-1-\nu} {(-1)^{\nu-\lambda}}{i,j}.  
\end{equation}
On the other hand, 
 by using $p^{i,j}(\lambda,\nu)$ and $p_{i,j}(\nu)$, 
 the relation between the unnormalized operators
 $\Cbb \lambda \nu{i,j}$ and 
 the renormalized operators $\Ctbb \lambda \nu {i,j}$
 defined in \eqref{eqn:Ciitilde} and \eqref{eqn:Cii-1tilde}
 are given as the following unified formula:
\[
  p^{i,j}(\lambda,\nu) \Cbb \lambda \nu {i,j} 
 =
  (-1)^{i-j} p_{i,j}(\nu) \Ctbb \lambda \nu {i,j}
  \quad
  \text{for $j \in \{i,i-1\}$.}
\]
Multiplying both sides of the equation \eqref{eqn:TCqpAij}
 by $p^{i,j}(\lambda,\nu)$, 
 we get the desired formula.  
\end{proof}

\begin{proposition}
\label{prop:10.34}
Suppose $\nu-\lambda \in {\mathbb{N}}$, 
 and $\delta$, $\varepsilon \in \{\pm\}$
 with $(-1)^{\nu-\lambda} = \delta \varepsilon$.  
\begin{enumerate}
\item[{\rm{(1)}}]
Suppose $0 \le i \le n-1$.  
Then the following two conditions on $(i, \lambda,\nu)$ are equivalent:
\begin{enumerate}
\item[{\rm{(i)}}]
The image of 
$\hphantom{i}
   \Ctbb \lambda \nu {i,i} \colon I_{\delta}(i,\lambda) \to J_{\varepsilon}(i,\nu)
$
 is contained in ${\operatorname{Ker}}\,(\Ttbb \nu {n-1-\nu}i)$.  
\item[{\rm{(ii)}}]
One of the following conditions holds:
\begin{enumerate}
\item[{\rm{(ii-a)}}]
$1 \le i \le n-1$, $\nu=i$, and ${\mathbb{Z}} \ni \lambda < i;$
\item[{\rm{(ii-b)}}]
$n$ is odd, $i=0$, and \eqref{eqn:he}$;$
\item[{\rm{(ii-c)}}]
$n$ is odd, $1 \le i \le n-1$, \eqref{eqn:he}, 
 and $\nu \ne n-1;$
\item[{\rm{(ii-d)}}]
$n$ is odd, $1 \le i \le \frac 1 2(n-1)$, 
 and $(\lambda,\nu) = (i,n-1-i)$.  
\end{enumerate}
\end{enumerate}
\item[{\rm{(2)}}]
Suppose $1 \le i \le n$.  
Then the following two conditions on $(i, \lambda,\nu)$ are equivalent:
\begin{enumerate}
\item[{\rm{(iii)}}]
The image of 
$
   \Ctbb \lambda \nu {i,i-1} \colon I_{\delta}(i,\lambda) \to J_{\varepsilon}(i-1,\nu)
$
 is contained in ${\operatorname{Ker}}\,(\Ttbb \nu {n-1-\nu}{i-1})$.  
\item[{\rm{(iv)}}]
One of the following holds:
\begin{enumerate}
\item[{\rm{(iv-a)}}]
$1 \le i \le n-1$, $\nu=n-i$, and ${\mathbb{Z}} \ni \lambda < n-i;$
\item[{\rm{(iv-b)}}]
$n$ is odd, $1 \le i \le n-1$, \eqref{eqn:he}, 
 and $\nu \ne n-1;$
\item[{\rm{(iv-c)}}]
$n$ is odd, $i=n$, and \eqref{eqn:he}$;$
\item[{\rm{(iv-d)}}]
$n$ is odd, $\frac 12 (n+1) \le i \le n-1$, 
 and $(\lambda,\nu) = (n-i,i-1)$.  
\end{enumerate}
\end{enumerate}
\end{enumerate}
\end{proposition}
The difference of this proposition from Theorems \ref{thm:180380} and \ref{thm:180386}
 is that the cases $i=\frac 1 2 (n-1)$ in (1)
 and $i=\frac 1 2 (n+1)$ in (2)
 are included in Proposition \ref{prop:10.34}.  
In these cases,
 the Knapp--Stein intertwining operator $\Ttbb \nu {n-1-\nu} j$
 vanishes 
 where $j=i$ in (1)
 and $=i-1$ in (2), 
 and the conditions (i) and (iii) do not provide 
 any information of $\operatorname{Image}(\Ctbb \lambda \nu{i,j})$.  
In these special cases,
 we shall study $\operatorname{Image}(\Ctbb \lambda \nu{i,j})$
 separately
 in Section \ref{subsec:10.6.3}
 by using the renormalized Knapp--Stein operators
 $\Tttbb \nu{n-1-\nu} j$.  

\begin{proof}
By the functional equation \eqref{eqn:180389} in Theorem \ref{thm:180389}, 
 we see that 
\[
 {\operatorname{Image}}\,(\Ctbb \lambda \nu {i,j})
 \subset {\operatorname{Ker}}\,(\Ttbb \nu {n-1-\nu}j)
\]
 if and only if $p^{i,j}(\lambda,\nu)=0$
 or $\Atbb \lambda {n-1-\nu}{(-1)^{\nu-\lambda}}{i,j}=0$.  
Suppose $0 \le i \le n$, $0 \le j \le n-1$, 
 and $j \in \{i,i-1\}$.  
By definition \eqref{eqn:pijlmdnu}, 
\begin{align*}
 p^{i,i}(\lambda,\nu) = 0 &\Leftrightarrow \lambda \ne \nu =i
\quad
\text{and}
\quad
1 \le i \le n-1,
\\
 p^{i,i-1}(\lambda,\nu) = 0 &\Leftrightarrow 
\lambda \ne \nu =n-i
\quad
\text{and}
\quad
1 \le i \le n-1.  
\end{align*}
On the other hand, 
 we claim 
\begin{align*}
\Atbb \lambda{n-1-\nu}{(-1)^{\nu-\lambda}}{i,i}=0
&\Leftrightarrow
\text{(ii-b), (ii-c), or (ii-d) holds};
\\
\Atbb \lambda{n-1-\nu}{(-1)^{\nu-\lambda}}{i,i-1}=0
&\Leftrightarrow
\text{(iv-b), (iv-c), or (iv-d) holds.}
\end{align*}
Let us verify the first equivalence 
 for $1 \le i \le n-1$.  
The vanishing condition 
of $\Atbb \lambda {\nu}{\pm}{i,j}$
 given in Theorem \ref{thm:161243} (1) and (3) shows that
 $\Atbb \lambda{n-1-\nu}{(-1)^{\nu-\lambda}}{i,i}=0$
 if and only if one of the following three conditions holds:
\begin{enumerate}
\item[$\bullet$]
$i=0$,
 $(\lambda,n-1-\nu) \in {\mathbb{Z}}^2$,
 $(\nu+1-n)-\lambda \equiv \nu-\lambda \mod 2$, 
 and $0 \le \nu+1-n \le -\lambda$;
\item[$\bullet$]
$i\ne 0$, $(\lambda,n-1-\nu) \in {\mathbb{Z}}^2$,
 $(\nu+1-n)-\lambda \equiv \nu-\lambda \mod 2$, 
 and $0 < \nu+1-n \le -\lambda$;
\item[$\bullet$]
$i\ne 0$, $\nu-\lambda \in 2{\mathbb{Z}}$, 
 and $(\lambda,n-1-\nu) =(i,i)$.  
\end{enumerate}
These conditions amount to (ii-b), (ii-c), and (ii-d)
 in Proposition \ref{prop:10.34} (1), 
 respectively.  
The second equivalence is shown similarly.  
Hence Proposition \ref{prop:10.34} is proved.  
\end{proof}
\begin{remark}
For $\lambda=\nu$, 
 the above conditions are fulfilled 
 if and only if $(\lambda,\nu)=(\frac 1 2 (n-1),\frac 1 2 (n-1))$
 and $i=\frac 1 2 (n-1)$
 in Proposition \ref{prop:10.34} (1) 
 or $i=\frac 1 2 (n+1)$ in Proposition \ref{prop:10.34} (2).  
This is exactly when $\Ttbb \nu {n-1-\nu}j$
 ($j=i$, $i-1$) vanishes.  
\end{remark}

\subsubsection{The case when $\Ttbb \nu {n-1-\nu}j=0$}
\label{subsec:10.6.3}

By Proposition \ref{prop:Tvanish}, 
 the Knapp--Stein operators
 $\Ttbb \nu {n-1-\nu} {j}$
 for the subgroup $G'=O(n,1)$
 vanishes if and only if $n$ is odd and 
\[
\nu=j=\frac{n-1}2.
\]
We note in this case that 
 $\nu-i=0$ for $i=j$
 and $\nu+i-n=0$ for $i=j+1$, 
 and therefore the definition \eqref{eqn:pijlmdnu} tells
\begin{equation*}
p^{i,j}(\lambda,\nu)
=
\begin{cases}
1 
\quad
&\text{if $\lambda = \nu$}, 
\\
0
\quad
&\text{if $\lambda \ne \nu$}.  
\end{cases}
\end{equation*}

When $n=2m+1$
 and $j= \frac 1 2 (n-1)$ ($=m$), 
 we use the renormalized Knapp--Stein operator, 
 see \eqref{eqn:Ttilde}, 
 given by 
\[
  \Tttbb \nu {2m-\nu}m = \frac 1 {\nu-m} \Ttbb \nu {2m-\nu}m.  
\]

\begin{proposition}
\label{prop:1803105}
Suppose $(G,G')=(O(n+1,1),O(n,1))$ with $n=2m+1$.  
Let $i=m$ or $m+1$.  
Then the composition $\Tttbb \nu {2m-\nu}m \circ \Ctbb \lambda \nu {i,m} \colon I_{\delta}(i,\lambda) \to J_{\varepsilon}(m,2m-\nu)$
 for $\nu-\lambda \in {\mathbb{N}}$
 and $\delta \varepsilon \in (-1)^{\nu-\lambda}$ is given as follows.   
\begin{enumerate}
\item[{\rm{(1)}}]
For $\nu-\lambda \in {\mathbb{N}}_+$, 
\[
  \Tttbb \nu {2m-\nu}m \circ \Ctbb \lambda \nu {i,m}
  = 
  \frac 1 2 q(\nu-\lambda) \Atbb \lambda {2m-\nu} {(-1)^{\nu-\lambda}} {i,m}.  
\]
In particular,
 if $m-\lambda \in {\mathbb{N}}_+$, 
 then 
\[
  \Tttbb m m m \circ \Ctbb \lambda m {i,m}
  = 
  (-1)^{i-m} \frac {\pi^m} {m!} \Ctbb \lambda m {i,m}.  
\]

\item[{\rm{(2)}}]
For $\nu=\lambda=m$, 
\[
     \Tttbb m m m \circ \Ctbb m m {i,m}
  = 
     \Attbb m m + {i,m} 
     + (-1)^{i-m+1} \frac{\pi^m}{m!} \Ctbb m m {i,m}.  
\]
\end{enumerate}
\end{proposition}
\begin{proof}
(1)\enspace
The functional equation \eqref{eqn:180389}
 with $n=2m+1$ and $j=m$ shows
\[
   \Ttbb \nu {2m-\nu} m \circ \Ctbb \lambda \nu {i,m}
  = 
  q(\nu-\lambda) p^{i,m}(\lambda,\nu) \Atbb \lambda {2m-\nu}{(-1)^{\nu-\lambda}} {i,m}.  
\]
By \eqref{eqn:pijlmdnu}, 
 we have for $i \in \{m,m+1\}$ and $\lambda \ne \nu$, 
\[
  p^{i,m}(\lambda,\nu) = \frac 1 2 (\nu-m).  
\]
Hence the first equation is verified.  
For the second statement, 
 we substitute $\nu=m$.  
Then the second equation follows from the residue formula \eqref{eqn:qAC}
 and from the fact 
 $\Ctbb \lambda m {i,m}=\Cbb \lambda m {i,m}$
 when $\lambda \ne m$.  
\par\noindent
(2)\enspace
The case $i=m+1$ will be shown in Lemma \ref{lem:TCAC}.  
The case $i=m$ is similar by using
\[
     \lim_{\nu \to m} \frac 1 {\nu-m} \Atbb \nu{2m-\nu} +{m, m}
     = 
     \Attbb m m + {m,m}
     - \frac{\pi^m}{m!} \Ctbb m m {m,m} 
\]
as it will be explained in \eqref{eqn:AClimit} of Chapter \ref{sec:pfSBrho}.  
\end{proof}

\subsubsection{Proof of Theorems \ref{thm:180380} and \ref{thm:180386}}
\label{subsec:10.6.4}
Suppose $0 \le j \le n-1$.  
Then the principal series representation $J_{\delta}(j,\nu)$
 is reducible as a module of $G'=O(n,1)$
 if and only if
\begin{equation}
\label{eqn:Jreducible}
\nu \in \{j,n-1-j\} \cup (-{\mathbb{N}}_+) \cup (n + {\mathbb{N}}), 
\end{equation}
see Proposition \ref{prop:redIilmd} (1).  
Suppose $\nu$ satisfies \eqref{eqn:Jreducible}.  
Then the proper submodules
 of $J_{\delta}(j,\nu)$ are described as follows:
\par\noindent
{Case 1.}\enspace
$(n,\nu) \ne (2j+1,j)$, 
 equivalently,
 $\Ttbb \nu {n-1-\nu}j \ne 0$.  
\par
In this case, 
 the unique proper submodule
 of $J_{\delta}(j,\nu)$ is given as the kernel 
 of the Knapp--Stein operator
 $\Ttbb \nu {n-1-\nu}j \colon J_{\delta}(j,\nu) \to J_{\delta}(j,n-1-\nu)$.  
\par\noindent
{Case 2.}\enspace
$(n,\nu) = (2j+1,j)$, 
 equivalently,
 $\Ttbb \nu {n-1-\nu}j = 0$.  
\par
In this case, 
 there are two proper submodules
 of $J_{\delta}(j,\nu)$, 
 which are given as the kernel 
 of 
 $\Tttbb j j j \pm \frac{\pi^j}{j!}{\operatorname{id}}
  \in {\operatorname{End}}_{G'}(J_{\delta}(j,j))$, 
 see Lemma \ref{lem:161745}.  
\begin{proof}
[Proof of Theorems \ref{thm:180380} and \ref{thm:180386}]
Assume $(n,\nu) \ne (2j+1,j)$.  
This excludes the case 
 where ${\mathbb{Z}} \ni \lambda \le j$ from the conditions (ii) 
 ($i=j$) and (iv) ($i=j+1$)
 in Proposition \ref{prop:10.34}.  
In this case 
 Theorems \ref{thm:180380} and \ref{thm:180386}
 are immediate consequences of Proposition \ref{prop:10.34}.

Assume now $(n,\nu,j)=(2m+1,m,m)$ for some $m \in {\mathbb{N}}_+$.  
Then $\Ctbb \lambda m {i,m}$ is not surjective
 if $\lambda < m$, 
 and is surjective if $\lambda=m$
 by Proposition \ref{prop:1803105} (1) and (2), 
 respectively.  
Thus Theorems \ref{thm:180380} and \ref{thm:180386} are proved.  
\end{proof}

\newpage
\section{Symmetry breaking operators
 for irreducible representations
 with infinitesimal character $\rho$ :
 Proof of Theorems \ref{thm:SBOvanish}
 and \ref{thm:SBOone}}
\label{sec:pfSBrho}
In the first half of this chapter, 
 we give a proof of Theorems \ref{thm:SBOvanish}
 and \ref{thm:SBOone}
 that determine the dimension of the space of symmetry breaking operators from
 {\it{irreducible}} representations $\Pi$ of $G=O(n+1,1)$
 to {\it{irreducible}} representations $\pi$ of the subgroup $G'=O(n,1)$
 when both $\Pi$ and $\pi$ have the trivial infinitesimal characters $\rho$,
 or equivalently
 by Theorem \ref{thm:LNM20} (2), 
 when 
\index{A}{IrrGrho@${\operatorname{Irr}}(G)_{\rho}$,
 set of irreducible admissible smooth representation of $G$
 with trivial infinitesimal character $\rho$\quad}
\begin{alignat*}{2}
\Pi 
&\in {\operatorname{Irr}}(G)_{\rho} 
&&=\{\Pi_{i,\delta} : 0 \le i \le n+1, \,\delta \in \{\pm\}\}, 
\\
\pi 
&\in {\operatorname{Irr}}(G')_{\rho} 
&&=\{\pi_{j,\varepsilon}: 0 \le j \le n, \,\varepsilon \in \{\pm\}\}.  
\end{alignat*} 
The proofs of Theorems \ref{thm:SBOvanish} and \ref{thm:SBOone} are completed
 in Section \ref{subsec:Pfvanish} and Section \ref{subsec:IJdual}, 
 respectively.  
In the latter half of this chapter,
 we give a concrete construction of such symmetry breaking operators from 
 $\Pi$ to $\pi$.  
We pursue such constructions 
 more than what we need for the proof
 for Theorems \ref{thm:SBOvanish} and \ref{thm:SBOone}:
 some of the results will be used 
 in calculating \lq\lq{periods}\rq\rq\
 in Chapter \ref{sec:period}.  
Our proof uses the symmetry breaking operators
 for {\it{principal series representations}}
 and their basic properties
 that we have developed in the previous chapters.

\subsection{Proof of the vanishing result (Theorem \ref{thm:SBOvanish})}
\label{subsec:Pfvanish}
This section gives a proof of the vanishing theorem
 of symmetry breaking operators (Theorem \ref{thm:SBOvanish}).  
In the same circle of the ideas,
 we also give a proof of multiplicity-free results
 (Proposition \ref{prop:mfPipj}).  
In order to study symmetry breaking for irreducible representations
 $\Pi_{i,\delta}$ of $G$, 
 we embed 
$
   {\operatorname{Hom}}_{G'}(\Pi_{i,\delta}|_{G'}, \pi_{j,\varepsilon})
$
 into the space of symmetry breaking operators
 between principal series representations
 as follows:
\begin{lemma}
\label{lem:170519}
Let $\delta, \varepsilon \in \{\pm\}$.  
Then we have natural embeddings:
\begin{enumerate}
\item[{\rm{(1)}}]
for $0 \le i \le n$ and $0 \le j \le n-1$, 
\begin{equation}
{\operatorname{Hom}}_{G'}(\Pi_{i,\delta}|_{G'}, \pi_{j,\varepsilon})
\subset
{\operatorname{Hom}}_{G'}(I_{\delta}(i,n-i)|_{G'}, J_{\varepsilon}(j,j));
\label{eqn:PpIJ1}
\end{equation}
\item[{\rm{(2)}}]
for $1 \le i \le n+1$ and $0 \le j \le n-1$, 
\begin{equation}
{\operatorname{Hom}}_{G'}(\Pi_{i,\delta}|_{G'}, \pi_{j,\varepsilon})
\subset
{\operatorname{Hom}}_{G'}(I_{-\delta}(i-1,i-1)|_{G'}, J_{\varepsilon}(j,j));
\label{eqn:PpIJ2}
\end{equation}
\item[{\rm{(3)}}]
for $0 \le i \le n$ and $1 \le j \le n$,
\begin{equation}
{\operatorname{Hom}}_{G'}(\Pi_{i,\delta}|_{G'}, \pi_{j,\varepsilon})
\subset
{\operatorname{Hom}}_{G'}(I_{\delta}(i,n-i)|_{G'}, J_{-\varepsilon}(j-1,n-j)); 
\label{eqn:PpIJ3}
\end{equation}
\item[{\rm{(4)}}]
for $1 \le i \le n+1$ and $1 \le j \le n$, 
\begin{equation}
{\operatorname{Hom}}_{G'}(\Pi_{i,\delta}|_{G'}, \pi_{j,\varepsilon})
\subset
{\operatorname{Hom}}_{G'}(I_{-\delta}(i-1,i-1)|_{G'}, J_{-\varepsilon}(j-1,n-j)).  
\label{eqn:PpIJ4}
\end{equation}
\end{enumerate}
\end{lemma}
\begin{proof}
We recall from Theorem \ref{thm:LNM20} (1)
 that there are surjective $G$-homomorphisms
\[
  I \twoheadrightarrow \Pi_{i,\delta}
\quad
  \text{for }
  I= I_{\delta}(i,n-i) 
  \text{ or }
  I_{-\delta}(i-1,i-1)
\]
 and injective $G'$-homomorphisms
\[
  \pi_{j,\varepsilon} \hookrightarrow J
 \quad \text{for $J=J_{\varepsilon}(j,j)$ or $J_{-\varepsilon}(j-1,n-j)$}.  
\]
Then the composition 
 $I \twoheadrightarrow \Pi_{i,\delta} \to \pi_{j,\varepsilon} \hookrightarrow J_{\varepsilon}(j,j)$
 yields the embeddings \eqref{eqn:PpIJ1}--\eqref{eqn:PpIJ4}.  
\end{proof}

\begin{proposition}
\label{prop:Ppvanish}
If $j \notin \{i-1,i\}$, 
 then 
$
{\operatorname{Hom}}_{G'}(\Pi_{i,\delta}|_{G'}, \pi_{j,\varepsilon})
=
\{0\}.  
$
\end{proposition}
\begin{proof}
Assume ${\operatorname{Hom}}_{G'}(\Pi_{i,\delta}|_{G'}, \pi_{j,\varepsilon})
\ne
\{0\}$.

Suppose first $1 \le i \le n$.  
By Theorem \ref{thm:1.1} (1), 
 we get $j \in \{i-3, i-2, i-1,i\}$ from \eqref{eqn:PpIJ2}, 
 and $j \in \{i-1, i, i+1,i+2\}$ from \eqref{eqn:PpIJ3}.  
Hence we conclude $j \in \{i-1,i\}$.

Suppose next $i=0$ or $n+1$.  
By using Theorem \ref{thm:1.1} (1) again, 
 we get $j \in \{0,1\}$ from \eqref{eqn:PpIJ1} for $i=0$, 
 and $j \in \{n-1,n\}$ from \eqref{eqn:PpIJ2} for $i=n+1$.  
Since $\dim_{\mathbb{C}} \Pi_{0,\delta}= \dim_{\mathbb{C}} \Pi_{n+1,\delta}=1$
 whereas both $\pi_{1,\varepsilon}$ and $\pi_{n-1,\varepsilon}$
 are infinite-dimensional irreducible representations of $G'$
 (Theorem \ref{thm:LNM20} (4)), 
 we have an obvious vanishing result:
\[
{\operatorname{Hom}}_{G'}(\Pi_{i,\delta}|_{G'}, \pi_{j,\varepsilon})
=\{0\}
\quad
\text{ if $(i,j)=(0,1)$ or $(n+1,n-1)$}.  
\]

Hence we conclude $j \in \{i-1,i\}$
 for $i=0$ or $n+1$, too.  
\end{proof}

\begin{proposition}
\label{prop:170526}
If $\delta \varepsilon =-$, 
 then 
$
   {\operatorname{Hom}}_{G'}(\Pi_{i,\delta}|_{G'}, \pi_{j,\varepsilon})
   =
   \{0\}.
$
\end{proposition}
\newpage
\begin{proof}
We have already proved the assertion in the case 
 $j \not \in \{i-1,i\}$
 in Proposition \ref{prop:Ppvanish}.  
Therefore it suffices to prove 
 the assertion in the case $j=i-1$ and $i$.  
We begin with the case $j=i-1$.

Suppose $2 \le i \le n$.  
Then by Theorem \ref{thm:1.1} (3), 
\[
   {\operatorname{Hom}}_{G'}(I_{\delta}(i,n-i)|_{G'}, J_{-\varepsilon}(i-2,n-i+1))
   =
   \{0\}
\]
because $\delta (-\varepsilon) =+$.  
This implies 
$
   {\operatorname{Hom}}_{G'}
   (\Pi_{i,\delta}|_{G'}, \pi_{i-1,\varepsilon})
   =
   \{0\}
$ from \eqref{eqn:PpIJ3}.

For the case
 $(i,j)=(1,0)$, 
 we know from \cite[Thm.~2.5 (1-a)]{sbon}
 that
\[
   {\operatorname{Hom}}_{G'}
   (\Pi_{1,-}|_{G'}, \pi_{0,+})
   =
   \{0\}.  
\]
($F(0)=\pi_{0,+}$ and $T(0)=\Pi_{1,-}$ with the notation therein.)
It then follows from  Proposition \ref{prop:SBOdual}
 that 
$
   {\operatorname{Hom}}_{G'}
   (\Pi_{1,+}|_{G'}, \pi_{0,-})=\{0\}
$.

For the case $(i,j)=(n+1,n)$, 
 we use the fact 
 that both $\Pi_{i,\delta}$ and $\pi_{j,\varepsilon}$ are 
 one-dimensional.  
In fact,
 we have isomorphisms
 $\Pi_{n+1,\delta} \simeq \chi_{-,\delta}$
 and $\pi_{n,\varepsilon} \simeq \chi_{-,\varepsilon}|_{G'}$
 by Theorem \ref{thm:LNM20} (4).  
Thus the vanishing assertion is straightforward for $j=i-1$
 ($1 \le i \le n+1$).

The case $j=i$ is derived from the case $j=i-1$
 by duality
 (see Proposition \ref{prop:SBOdual}). 
\end{proof}

By Propositions \ref{prop:Ppvanish} and \ref{prop:170526}, 
 we have completed the proof of Theorem \ref{thm:SBOvanish}.

\subsection{Construction of symmetry breaking operators from 
 $\Pi_{i,\delta}$ to $\pi_{i,\delta}$: 
Proof of Theorem \ref{thm:SBOone}}
\label{subsec:StageCii}
In this section we prove
 the existence and the uniqueness 
 (up to scalar multiplication) of symmetry breaking operators from 
 the irreducible $G$-module $\Pi_{i,\delta}$ 
 to the irreducible $G'$-module $\pi_{j,\varepsilon}$
 when $j \in \{i-1,i\}$ and $\delta \varepsilon= +$, 
 and thus complete the proof of Theorem \ref{thm:SBOone}.  
Moreover,
 we investigate their $(K,K')$-spectrum
 for 
\index{B}{minimalKtype@minimal $K$-type}
minimal $K$- and $K'$-types,
 and also give an explicit construction of such operators.  

\subsubsection{Generators of symmetry breaking operators
 between principal series representations
 having the trivial infinitesimal character $\rho$}
\label{subsec:StageBii}

We have determined explicit generators 
 of symmetry breaking operators
 $I_{\delta}(i,\lambda) \to J_{\varepsilon}(j,\nu)$
 in Theorem \ref{thm:SBObasis}.  
In this subsection,
 we extract some special cases
 which will be used for the proof of Theorem \ref{thm:SBOone}.

The following lemma is used for the proof
 of the multiplicity-free theorem
 (Proposition \ref{prop:mfPipj} below), 
 and also for an explicit construction of nonzero symmetry breaking operators
 $\Pi \to \pi$ with $\Pi \in {\operatorname{Irr}}(G)_{\rho}$
 and $\pi \in {\operatorname{Irr}}(G')_{\rho}$
 (Proposition \ref{prop:172459}).

\begin{lemma}
\label{lem:170529}
\begin{enumerate}
\item[{\rm{(1)}}]
Suppose $0 \le i \le n-1$
 and $\delta \varepsilon=+$.  
Then 
\begin{equation*}
{\operatorname{Hom}}_{G'}
     (I_{\delta}(i,n-i)|_{G'}, J_{\delta}(i,i))
\simeq
\begin{cases}
{\mathbb{C}}\Atbb {n-i}i+ {i,i}
\quad
&\text{if $2i \ne n$},  
\\
{\mathbb{C}}\Attbb {n-i}i+ {i,i}
\oplus
{\mathbb{C}}\Ctbb {n-i}i {i,i}
\quad
&\text{if $2i = n$}.  
\end{cases}
\end{equation*}
\item[{\rm{(2)}}]
Suppose $1 \le i \le n-1$
 and $\delta \varepsilon=+$.  
Then 
\begin{equation*}
{\operatorname{Hom}}_{G'}
     (I_{-\delta}(i-1,i-1)|_{G'}, J_{\varepsilon}(i,i))
=
{\mathbb{C}}\Ctbb {i-1} i {i-1,i}.  
\end{equation*}
\item[{\rm{(3)}}]
Suppose $0 \le i \le n-1$
 and $\delta \in \{\pm\}$.  
Then we have
\begin{equation*}
{\operatorname{Hom}}_{G'}
     (I_{\delta}(i,i)|_{G'}, J_{\delta}(i,i))
=
{\mathbb{C}}\Attbb {i}i+ {i,i}
\oplus
{\mathbb{C}} \Ctbb i i{i,i}.    
\end{equation*}
\item[{\rm{(4)}}]
Suppose $1 \le i \le n$
 and $\delta \varepsilon=+$.  
Then 
\begin{multline*}
{\operatorname{Hom}}_{G'}
     (I_{-\delta}(i-1,i-1)|_{G'}, J_{-\varepsilon}(i-1,n-i))
\\
\simeq
\begin{cases}
{\mathbb{C}}\Atbb {i-1}{n-i}+ {i-1,i-1}
\quad
&\text{if $n \ne 2i-1$},  
\\
{\mathbb{C}}\Attbb {i-1}{n-i}+ {i-1,i-1}
\oplus
{\mathbb{C}} \Ctbb {i-1}{n-i}{i-1,i-1}
\quad
&\text{if $n=2i-1$}.  
\end{cases}
\end{multline*}
\end{enumerate}
\end{lemma}
\begin{proof}
We determined the dimension of the left-hand side
 by Theorem \ref{thm:1.1} (2) and (3).  
Then the lemma follows from Theorem \ref{thm:SBObasis}
 for (1), (3), (4);
 and from Fact \ref{fact:3.9} for (2).  
\end{proof}

\begin{remark}
\label{rem:170529}
In the $N$-picture
 where the open Bruhat cells
 for the pair of the real flag manifolds
 $G'/P' \subset G/P$
 are represented by ${\mathbb{R}}^{n-1} \subset {\mathbb{R}}^{n}$, 
 we have $\Ctbb {n-i}{i}{i,i} ={\operatorname{Rest}}_{x_n=0}$
 in Lemma \ref{lem:170529} (1), 
 $\Ctbb {i-1}{i}{i-1,i} ={\operatorname{Rest}}_{x_n=0}
 \circ d_{\mathbb{R}^n}$ in (2), 
$\Ctbb {i}{i}{i,i} ={\operatorname{Rest}}_{x_n=0}$
 in (3), 
 and 
 $\Ctbb {i-1}{i-1}{i-1,i-1} =
{\operatorname{Rest}}_{x_n=0}$
 in (4).  
\end{remark}

The following lemma is used for an alternative construction
 (see Proposition \ref{prop:172459} below)
 of symmetry breaking operators $\Pi_{i,\delta} \to \pi_{i,\delta}$. 
\begin{lemma}
\label{lem:172459}
Suppose $1 \le i \le n$
 and $\delta \in \{\pm\}$.  
Then we have
\begin{equation*}
{\operatorname{Hom}}_{G'}
     (I_{\delta}(i,i)|_{G'}, J_{-\delta}(i-1,n-i))
=
{\mathbb{C}}\Atbb {i}{n-i}- {i,i-1}.  
\end{equation*}
\end{lemma}

\begin{proof}
By Theorem \ref{thm:1.1} (2), 
 $\Atbb {i}{n-i}- {i,i-1} \ne 0$, 
 and therefore the lemma follows from Theorem \ref{thm:SBObasis}.  
\end{proof}

\subsubsection{Multiplicity-free property of symmetry breaking}
In this subsection,
 we prove the following multiplicity-free property:
\begin{proposition}
\label{prop:mfPipj}
For any $0 \le i \le n+1$, 
 $0 \le j \le n$, 
 and $\delta, \varepsilon \in \{\pm\}$, 
 we have
\begin{equation}
\label{eqn:mfPipj}
   \dim_{\mathbb{C}}{\operatorname{Hom}}_{G'}
   (\Pi_{i,\delta}|_{G'},\pi_{j,\varepsilon}) \le 1.  
\end{equation}
\end{proposition}
Proposition \ref{prop:mfPipj} is a very special case 
 of the multiplicity-free theorem 
 which was proved in Sun--Zhu \cite{SunZhu}, 
 however,
 we give a different proof based on Lemmas \ref{lem:170519}
 and \ref{lem:170529}
 because the following short proof illustrates the idea of this chapter.  
\begin{proof}
[Proof of Proposition \ref{prop:mfPipj}]
Owing to the vanishing results (Theorem \ref{thm:SBOvanish}), 
 it suffices to show \eqref{eqn:mfPipj}
when $j \in \{i-1,i\}$ and $\delta \varepsilon=+$.  
Moreover,
 the case $j=i-1$ can be reduced to the case $j=i$
 by the duality
 between the spaces of symmetry breaking operators
 (Proposition \ref{prop:SBOdual}).  
Henceforth, 
 we assume $j=i \in \{0,1,\ldots,n\}$
 and $\delta \varepsilon =+$.  
Then, 
 owing to the embedding results given in Lemma \ref{lem:170519}, 
the multiplicity-free property
 \eqref{eqn:mfPipj} holds for $1 \le i \le n-1$
 by Lemma \ref{lem:170529} (2), 
 and for $i =0$ and $n$ by Lemma \ref{lem:170529} (1) and (4).  
Thus Proposition \ref{prop:mfPipj} is proved.  
\end{proof}

\subsubsection{Multiplicity-one property: Proof of Theorem \ref{thm:SBOone}}
\label{subsec:pfSBOone}

In proving Theorem \ref{thm:SBOone}, 
 we use the following proposition,
 whose proof is deferred at the next subsection.  
\begin{proposition}
\label{prop:existii}
 ${\operatorname{Hom}}_{G'}(\Pi_{i,\delta}|_{G'}, \pi_{i,\delta}) \ne \{0\}$
 for all $0 \le i \le n$ and $\delta \in \{\pm\}$.  
\end{proposition}
\begin{remark}
\label{rem:existii}
Obviously Proposition \ref{prop:existii} holds for $i=0$
 because $\Pi_{0,\delta}|_{G'} \simeq \pi_{0,\delta}$ as $G'$-modules
 for $\delta \in \{\pm\}$.  
Indeed, 
 the $G$-modules $\Pi_{0,+}$ and $\Pi_{0,-}$ are the one-dimensional representations
 ${\bf{1}}$ and respectively
 $\chi_{+-}$ 
 (Theorem \ref{thm:LNM20} (4)),
 and likewise for the $G'$-modules $\pi_{0,\pm}$.  
\end{remark}

Before giving a proof of Proposition \ref{prop:existii}, 
 we show
 that Proposition \ref{prop:existii} implies Theorem \ref{thm:SBOone}.  
\begin{proof}
[Proof of Theorem \ref{thm:SBOone}]
By the duality among the spaces
 of symmetry breaking operators
 (Proposition \ref{prop:SBOdual}), 
 we may and do assume $j=i$ and $\delta=\varepsilon=+$
 because $\widetilde j:=n-j$
 and $\widetilde i:=n+1-i$
 satisfy $\widetilde j = \widetilde i-1$
 if and only if $j=i$.  
Then Theorem \ref{thm:SBOone} follows from Propositions \ref{prop:mfPipj}
 (uniqueness) and \ref{prop:existii} (existence).  
\end{proof}
For later purpose,
 we need a refinement of Proposition \ref{prop:existii}
 by providing information of $(K,K')$-spectrum
 in Proposition \ref{prop:AiiAq} below.  
For this,
 we fix some terminology:
\begin{definition}
[minimal $K$-type]
\label{def:minK}
We set $m:=[\frac{n+1}2]$.  
Suppose $\mu \in \widehat K$.  
To describe an irreducible finite-dimensional representation $\mu$
 of $K=O(n+1) \times O(1)$, 
 we use the notation in Section \ref{subsec:fdimrep}
 in Appendix I
 rather than the previous one
 in Section \ref{subsec:ONWeyl}, 
 and write
\[
  \mu = \Kirredrep{O(n+1)}{\sigma_1, \cdots,\sigma_m}_{\varepsilon} \boxtimes \delta
\]
for $\sigma =(\sigma_1, \cdots,\sigma_m) \in \Lambda^+(m)$
 and $\varepsilon,\delta \in \{\pm\}$.  
We define $\|\mu\|>0$ by 
\[
\|\mu\|^2 = \sum_{j=1}^m (\sigma_j + n + 1-2j)^2
\quad
(=\|\sigma + 2 \rho_c\|^2), 
\]
where $2 \rho_c =(n-1,n-3,\cdots,n+1-2m)$ is the sum
 of positive root
 for ${\mathfrak{k}}_{\mathbb{C}}={\mathfrak{o}}(n+1,{\mathbb{C}})$
 in the standard coordinates.  
For a nonzero admissible representation $\Pi$ of $G$, 
 the set of {\it{minimal $K$-types}} of $\Pi$ is 
\[
   \{\mu \in \widehat K
   :
   \text{$\mu$ occurs in $\Pi$, 
 and $\|\mu\|$ is minimal with this property}\}, 
\]
see \cite[Chap.~2]{KV}
 or \cite[Def.~5.4.18]{Vogan81}.  
\end{definition}

We then observe:
\begin{remark}
[minimal $K$-type]
\label{rem:minK}
\index{B}{minimalKtype@minimal $K$-type|textbf}
The basic $K$-type 
 (see Definition \ref{def:basicK}) of the principal series representation
 $I_{\delta}(i,\lambda)$
 is the unique minimal $K$-type
 of the irreducible $G$-module $\Pi_{i,\delta}$, 
 as stated in Theorem \ref{thm:LNM20} (3).  
\end{remark}
\begin{proposition}
\label{prop:AiiAq}
Let $(G,G')=(O(n+1,1),O(n,1))$, 
 $0 \le i \le n$
 and $\delta \in \{\pm\}$.  
Then there exists a nonzero symmetry breaking operator
\begin{equation}
\label{eqn:Pipi}
   A_{i,i} \colon \Pi_{i,\delta} \to \pi_{i,\delta}
\end{equation}
 such that its $(K,K')$-spectrum
 for the minimal $K'$- and $K$-types $\mub (i,\delta)' (\hookrightarrow \mub(i,\delta))$
 is nonzero.  
\end{proposition}

Proposition \ref{prop:AiiAq} is an {\it{existence}} theorem,
 however,
 we shall prove it by {\it{constructing}}
 nonzero symmetry breaking operators 
 $\Pi_{i,\delta} \to \pi_{i,\delta}$, 
 see Proposition \ref{prop:172459}
 in the next subsection.  
Alternative constructions are also given in Sections \ref{subsec:irrIJ} and \ref{subsec:3irrIJ}, 
 and thus we construct symmetry breaking operators
 $\Pi_{i,\delta} \to \pi_{i,\delta}$
 in the following three ways:
\begin{center}\begin{tabular}{ll}
$\bullet$\ $\Atbb i{n-i}- {i,i-1} \colon I_{\delta}(i,i) \to J_{-\delta}(i-1,n-i)$, &
(Proposition \ref{prop:172459}), \\[0.5cm]
$\bullet$\ $\Attbb {i}{i}+ {i,i} \colon I_{\delta}(i,i) \to J_{\delta}(i,i)$, &
(Proposition \ref{prop:SBOi4}), \\[0.5cm]
$\bullet$\ $\Atbb {n-i}{i}+ {i,i} \colon I_{\delta}(i,n-i) \to J_{\delta}(i,i)$, &
(Proposition \ref{prop:170530}).  
\end{tabular}\end{center}

\subsubsection{First construction $\Pi_{i,\delta} \to \pi_{i,\delta}$
 ($1 \le i \le n$)}
\label{subsec:IJdual}
In this subsection,
 we construct a nonzero symmetry breaking operator
\[
   \Pi_{i,\delta} \to \pi_{i,\delta}
   \qquad
  \text{for $1 \le i \le n$, $\delta \in \{\pm\}$, }
\]
by using Lemma \ref{lem:172459}.  
\begin{proposition}
\label{prop:172459}
Suppose $1 \le i \le n$
 and $\delta \in \{\pm\}$.  
Then the normalized symmetry breaking operator
\begin{equation*}
    \Atbb {i}{n-i}- {i,i-1}
    \colon 
   I_{\delta}(i,i) \to J_{-\delta}(i-1,n-i)
\]
satisfies the following:
\begin{enumerate}
\item[{\rm{(1)}}]
${\operatorname{Image}}(\Atbb {i}{n-i}- {i,i-1})_{K'}=(\pi_{i,\delta})_{K'}$
 as $({\mathfrak{g}}',K')$-modules;
\item[{\rm{(2)}}]
$\Atbb {i}{n-i}- {i,i-1}|_{\Pi_{i,\delta}} \ne 0$.  
\end{enumerate}
In particular,
 it induces a symmetry breaking operator
 $\Pi_{i,\delta} \to \pi_{i,\delta}$
 as in the diagram below.  
Moreover,
 the $(K,K')$-spectrum
 of the resulting operator 
 for the minimal $K'$- and $K$-types $\mub (i,\delta)'$
 $(\hookrightarrow \mub(i,\delta))$ is nonzero.  
\begin{eqnarray*}
\xymatrix@R=2pt{
I_{\delta}(i,i)\ar[rr]^{\Atbb {i}{n-i}- {i,i-1}\quad}
&& J_{-\delta}(i-1,n-i)
\\
\bigcup
&&\bigcup
\\
\Pi_{i,\delta}
\ar@{-->}[rr]
&&\pi_{i,\delta}
}
\end{eqnarray*}
\end{proposition}

\begin{convention}
\label{conv:Image}
Hereafter,
 by abuse of notation,
 we shall write simply as 
${\operatorname{Image}}(\Atbb {i}{n-i}- {i,i-1})=\pi_{i,\delta}$
 if their underlying $({\mathfrak{g}}',K')$-modules coincide
 ({\it{cf}}. Proposition \ref{prop:172459} (1)).  
\end{convention}
\newpage
\begin{proof}
[Proof of Proposition \ref{prop:172459}]
(1)\enspace
First we observe
\[
   {\operatorname{Image}}(\Atbb {i}{n-i}- {i,i-1})
  \subset 
   {\operatorname{Ker}}(\Ttbb {n-i}{i-1} {i-1})
\]
because Theorem \ref{thm:TAA} with $\nu=n-i$ tells the functional equation
$
   \Ttbb {n-i}{i-1}{i-1} \circ \Atbb i {n-i}-{i,i-1}=0.  
$

When $n \ne 2i-1$, 
 we conclude ${\operatorname{Image}}(\Atbb i {n-i}-{i,i-1})=\pi_{i,\delta}$ 
 by Proposition \ref{prop:Timage}
 because $\pi_{i,\delta}$ is irreducible as a $G'$-module.  
When $n=2i-1$, 
 the Knapp--Stein operator $\Ttbb {n-i}{i-1}{i-1}$
 vanishes 
 (Proposition \ref{prop:Tvanish}).  
Instead we use the following renormalized Knapp--Stein operator
 (see \eqref{eqn:Ttilde}):
\[
   \Tttbb \nu{n-1-\nu}{i-1}
   =
  \frac {1}{\nu-i+1}\Ttbb \nu{n-1-\nu}{i-1}.  
\]
Then the functional equation given in Theorem \ref{thm:TAA} implies
\[
   \left(
   \Tttbb {i-1}{i-1}{i-1} + \frac{\pi^{i-1}}{(i-1)!}\, {\operatorname{id}}
   \right)
   \circ \Atbb \lambda{i-1}-{i,i-1}=0.  
\]
By Lemma \ref{lem:161745} applied to the subgroup $G'=O(n,1)$ $(=O(2i-1,1))$, 
 we conclude 
 ${\operatorname{Image}}(\Atbb \lambda {n-i}- {i,i-1})
={\operatorname{Image}}(\Atbb \lambda {i-1}- {i,i-1})=\pi_{i,\delta}$
 in the case $n=2i-1$, too.  
\newline\noindent
(2)\enspace
The second statement follows from the fact
 that the $(K,K')$-spectrum of $\Atbb \lambda \nu - {i,i-1}$
 (Theorem \ref{thm:153315})
 for the basic $K$-types
 $(\mu,\mu')=(\mub(i,\delta),\mus(i-1,-\delta)')$
 does not vanish.  
The last assertion is derived from the following observation 
 (see \eqref{eqn:flatsharp}):
there are isomorphisms
 of representations of $K'=O(n) \times O(1)$, 
\[
  \mus (i-1,-\delta)'\simeq \mub(i,\delta)'.  
\]
Hence Proposition \ref{prop:172459} is proved.  
\end{proof}

\begin{proof}
[Proof of Proposition \ref{prop:AiiAq}]
Clear from Proposition \ref{prop:172459}
 and Remark \ref{rem:existii}.  
\end{proof}

Thus,
 the proof of Theorem \ref{thm:SBOone} has been completed.  

\vskip 1pc
For the rest of this chapter,
 we give alternative constructions
 of symmetry breaking operators for later purposes.  

\subsubsection{Second construction $\Pi_{i,\delta} \to \pi_{i,\delta}$
 ($0 \le i \le n-1$)}
\label{subsec:irrIJ}

In this subsection,
 we provide another construction
 of a nonzero symmetry breaking operator
\[
    \Pi_{i,\delta} \to \pi_{i,\delta}
   \quad
   \text{for $0 \le i \le n-1$, $\delta \in \{\pm\}$}, 
\]
by using Lemma \ref{lem:170529} (3).  
\begin{proposition}
\label{prop:SBOi4}
Suppose $0 \le i \le n-1$ and $\delta \in \{\pm\}$.  
Then the renormalized operator
\[
  \Attbb i i + {i,i}\colon I_{\delta}(i,i) \to J_{\delta}(i,i)
\]
 satisfies the following:
\begin{enumerate}
\item[{\rm{(1)}}]
${\operatorname{Image}}(\Attbb i i + {i,i})=\pi_{i,\delta}$;
\item[{\rm{(2)}}]
$\Attbb i i + {i,i}|_{\Pi_{i,\delta}}\ne 0.$
\end{enumerate}
In particular,
 it induces a symmetry breaking operator
 $\Pi_{i,\delta} \to \pi_{i,\delta}$ 
 as in the diagram below.  
Moreover,
 the $(K,K')$-spectrum
 of the resulting operator 
 for the minimal $K'$- and $K$-types $\mub (i,\delta)'$
 $(\hookrightarrow \mub(i,\delta))$ is nonzero.  

\begin{eqnarray*}
\xymatrix@R=2pt{
I_{\delta}(i,i)\ar[rr]^{\Attbb {i}{i}+ {i,i}\quad}
&& J_{\delta}(i,i)
\\
\bigcup
&&\bigcup
\\
\Pi_{i,\delta}
\ar@{-->}[rr]
&&\pi_{i,\delta}
}
\end{eqnarray*}
\end{proposition}
\begin{proof}
[Proof of Proposition \ref{prop:SBOi4}]
(1)\enspace
By the functional equation \eqref{eqn:TA2tilde}, 
we have
\[
  {\operatorname{Image}}(\Attbb \lambda i + {i,i})
 \subset {\operatorname{Ker}} (\Ttbb i {n-1-i}i).  
\]
When $n \ne 2i+1$, 
 we conclude 
 ${\operatorname{Image}}(\Attbb \lambda i + {i,i})=\pi_{i,\delta}$
 by Proposition \ref{prop:Timage}.

When $n=2i+1$, 
 the Knapp--Stein operator $\Ttbb i {n-1-i}i\equiv \Ttbb i i i$ vanishes 
 (Proposition \ref{prop:Tvanish}).  
Instead we use the functional equation \eqref{eqn:TAA3tilde}
 for the renormalized operators $\Tttbb i i i$ and $\Attbb \lambda i +{i,i}$, 
 which tells
 that
\[
   {\operatorname{Image}}(\Attbb \lambda i + {i,i})
   \subset {\operatorname{Ker}} 
   \left(\Tttbb i i i - \frac{\pi^i}{\Gamma(i+1)}\, {\operatorname{id}}\right).  
\]
By Lemma \ref{lem:161745}, 
 we conclude 
$
 {\operatorname{Image}}(\Attbb \lambda i + {i,i})=\pi_{i,\delta} 
$
 because $\Attbb \lambda i + {i,i}$ is nonzero
 and $\pi_{i,\delta}$ is irreducible as a $G'$-module.  
\newline\noindent
(2)\enspace
The assertion follows readily from the $(K,K')$-spectrum
 of the renormalized operator $\Attbb \lambda i + {i,i}$
 (see \eqref{eqn:SAiitilde})
 for the basic $K$- and $K'$-types
 $(\mu,\mu')=(\mub(i,\delta),\mub(i,\delta)')$.  
\end{proof}

\subsubsection{Third construction $\Pi_{i,\delta} \to \pi_{i,\delta}$}
\label{subsec:3irrIJ}

We give yet another construction
 of a nonzero symmetry breaking operator
 $\Pi_{i,\delta} \to \pi_{i,\delta}$
 in the case $n \ne 2i$.  
In the case $n=2i$, 
 the normalized operator 
 $\Atbb {n-i} i + {i,i}$ vanishes.  
We shall discuss this case separately
 in Section \ref{subsec:IJmm+}, 
see Proposition \ref{prop:Ammimage}.  

\begin{proposition}
\label{prop:170530}
If $2i \ne n$, 
 then 
$
   \Atbb {n-i}i+ {i,i}
   \in 
   {\operatorname{Hom}}_{G'}
   (I_\delta(i,n-i)|_{G'}, J_\delta(i,i))
$ 
satisfies
\[
   \Atbb {n-i}i+ {i,i}|_{\Pi_{i+1,-\delta}} \equiv 0
\quad
   \text{and}
\quad
   {\operatorname{Image}}(\Atbb {n-i}i+ {i,i})
   = \pi_{i,\delta}.  
\]
Thus it induces a symmetry breaking operator 
$\Pi_{i,\delta} \to \pi_{i,\delta}$
 as in the diagram below.  
Moreover,
 the $(K,K')$-spectrum
 of the resulting operator 
 for the minimal $K'$- and $K$-types $\mub (i,\delta)'$
 $(\hookrightarrow \mub(i,\delta))$ is nonzero.  
\begin{eqnarray*}
\xymatrix@C=20pt@R=2pt{
I_{\delta}(i,n-i)\ar[rr]^{\Atbb {n-i}{i}+ {i,i}\,\,}
\ar[dd]
&& J_{\delta}(i,i)
\\
&&\bigcup
\\
\Pi_{i,\delta} \simeq
I_{\delta}(i,n-i)/\Pi_{i+1,-\delta}
\ar@{-->}[rr]
&&\pi_{i,\delta}
}
\end{eqnarray*}
\end{proposition}

\begin{proof}
Since $\Atbb i i + {i,i} =0$
 by Theorem \ref{thm:161243} (1), 
 the composition 
$
     \Atbb {n-i}i+ {i,i} \circ \Ttbb i{n-i}{i}
$
 vanishes by the functional equation (Theorem \ref{thm:ATA}).  
Thus $\Atbb {n-i}i+ {i,i}$ is identically zero
 on ${\operatorname{Image}} (\Ttbb i{n-i}{i}) \simeq \Pi_{i+1,-\delta}$
 (see Proposition \ref{prop:Timage}).

For the second assertion,
 we use another functional equation 
 (Theorem \ref{thm:TAA})
 to get $\Ttbb i{n-1-i}{i} \circ \Atbb {n-i}i+ {i,i}=0$.  
Hence
\[
     {\operatorname{Image}}(\Atbb {n-i}i+ {i,i})
     \subset
     {\operatorname{Ker}}(\Ttbb i{n-1-i}i)
     \simeq 
     \pi_{i,\delta}
\]
 by Proposition \ref{prop:Timage}.  
Since $\Atbb {n-i}i+ {i,i} \ne 0$
(see Theorem \ref{thm:161243} (1))
 and since $\pi_{i,\delta}$ is irreducible,
 the underlying $({\mathfrak{g}}',K')$-modules
 of ${\operatorname{Image}}(\Atbb {n-i}i+ {i,i})$
 and $\pi_{i,\delta}$ coincide.  
\end{proof}

\subsection{Splitting of $I_{\delta}(m,m)$ and its symmetry breaking for $(G,G')=(O(2m+1,1),O(2m,1))$}
\label{subsec:2m}

Suppose $n$ is even, 
 say $n=2m$.  
A distinguished feature in this setting
 is that the principal series representation 
 $I_{\delta}(m,\lambda)$ of $G =O(2m+1,1)$
 splits into the direct sum
 of two irreducible $G$-modules
 when $\lambda=m$: for $\delta\in \{\pm\}$, 
\begin{equation}
\label{eqn:I+n2}
I_{\delta}(m,m) \simeq \Pi_{m,\delta} \oplus \Pi_{m+1,-\delta}, 
\end{equation}
both of which are smooth irreducible 
\index{B}{temperedrep@tempered representation}
tempered representations of $G$,
 see Theorem \ref{thm:LNM20} (1) and (8).  
Accordingly, 
 the space of symmetry breaking operators has a direct sum decomposition:
\begin{multline}
\label{eqn:mmdeco}
{\operatorname{Hom}}_{G'}(I_{\delta}(m,m)|_{G'}, J_{\varepsilon}(m, m))
\\
\simeq 
{\operatorname{Hom}}_{G'}(\Pi_{m,\delta}|_{G'}, J_{\varepsilon}(m, m))
\oplus 
{\operatorname{Hom}}_{G'}(\Pi_{m+1,-\delta}|_{G'}, J_{\varepsilon}(m, m)), 
\end{multline}
for each $\varepsilon \in \{\pm\}$.  
The left-hand side of \eqref{eqn:mmdeco} 
 has been understood by the classification
 of symmetry breaking operators
 given in Theorem \ref{thm:SBObasis}
 (see \eqref{eqn:IJmm} as below). 
On the other hand,
 the target space $J_{\varepsilon}(m,m)$ is not irreducible
 as a $G'$-module.  
We recall from Theorem \ref{thm:LNM20} (1)
 that the principal series representation
 $J_\varepsilon(m,\nu)$
 of $G'=O(2m,1)$ at $\nu=m$ has a nonsplitting exact sequence
 of $G'$-modules:
\begin{equation}
\label{eqn:Jmm}
0 \to \pi_{m,\varepsilon} \to J_{\varepsilon}(m,m) \to \pi_{m+1,-\varepsilon} \to 0.  
\end{equation}
With this in mind,
 we shall take a closer look 
 at the right-hand side of \eqref{eqn:mmdeco}
 and determine each summand as follows:

\begin{equation}
\label{tbl:PiJmm}
\begin{tabular}{c|c|c}
& {$\delta\varepsilon=+$}
&{$\delta\varepsilon=-$}
\\
\hline
  ${\operatorname{Hom}}_{G'}(\Pi_{m,\delta}|_{G'}, J_{\varepsilon}(m, m))$ 
& ${\mathbb{C}}$
& $\{0\}$
\\
${\operatorname{Hom}}_{G'}(\Pi_{m+1,-\delta}|_{G'}, J_{\varepsilon}(m, m))$ 
& ${\mathbb{C}}$
& ${\mathbb{C}}$ 
\end{tabular}
\end{equation}
See Section \ref{subsec:IJmm+}
 for the left column of \eqref{tbl:PiJmm} in detail,
 and for Section \ref{subsec:IJmm-} for the right column.

\subsubsection
{${\operatorname{Hom}}_{G'}(I_{\delta}(m,m)|_{G'}, J_{\varepsilon}(m,m))$
 with $\delta \varepsilon=+$}
\label{subsec:IJmm+}
We begin with the case $\delta\varepsilon =+$.  
Without loss of generality,
 we may and do assume $\delta=\varepsilon=+$.

Then Lemma \ref{lem:170529} (3) 
 with Remark \ref{rem:170529} tells that 
\begin{equation}
\label{eqn:IJmm}
{\operatorname{Hom}}_{G'}(I_{+}(m,m)|_{G'}, 
                          J_{+}(m,m))
=
{\mathbb{C}} \Attbb {m} {m} + {m,m}
\oplus
{\mathbb{C}} {\operatorname{Rest}}.  
\end{equation}

The first generator $\Attbb {m} {m} + {m,m}$ is defined
 as the renormalization (Theorem \ref{thm:170340})
\begin{equation}
\label{eqn:Amm2tilde}
\Attbb {m} {m} + {m,m}
=
\lim_{\lambda \to m}
\Gamma\left(\frac {\lambda-m}{2}\right)
\Atbb {\lambda} {m} + {m,m}
\end{equation}
of the normalized regular symmetry breaking operator
 $\Atbb {\lambda} {m} + {m,m}$
 which vanishes at $\lambda=m$
 (Theorem \ref{thm:161243}).  
The second generator,
 ${\operatorname{Rest}}\equiv {\operatorname{Rest}}_{x_n=0}$, 
 is the obvious symmetry breaking operator
({\it{cf.}}~Lemma \ref{lem:152268}), 
 given by ${\operatorname{Rest}}_{x_n=0}$
 in the $N$-picture.  
By using the second generator,
 we obtain the following.

\begin{proposition}
\label{prop:IJmmbase}
Let $(G,G')=(O(2m+1,1),O(2m,1))$.  
Then we have 
\begin{align*}
{\operatorname{Hom}}_{G'}(\Pi_{m,+}|_{G'},J_+(m,m))
&=
{\mathbb{C}} {\operatorname{Rest}}|_{\Pi_{m,+}}, 
\\
{\operatorname{Hom}}_{G'}(\Pi_{m+1,-}|_{G'},J_+(m,m))
&=
{\mathbb{C}} {\operatorname{Rest}}|_{\Pi_{m+1,-}}.  
\end{align*}
\end{proposition}

\begin{proof}
By the direct sum decompositions \eqref{eqn:IJmm} and \eqref{eqn:I+n2}, 
 we have
\begin{align*}
  2 =& \dim_{\mathbb{C}} \operatorname{Hom}_{G'}(I_+(m, m)|_{G'}, J_+(m, m))
\\
    =&\dim_{\mathbb{C}} \operatorname{Hom}_{G'}(\Pi_{m, +}|_{G'}, J_+(m, m))
 + \dim_{\mathbb{C}} \operatorname{Hom}_{G'}(\Pi_{m+1, -}|_{G'}, J_+(m, m)).   
\end{align*}
On the other hand,
 we know from Lemma \ref{lem:152268}
 that ${\operatorname{Rest}}|_{\Pi_{m,+}}\ne 0$
 and ${\operatorname{Rest}}|_{\Pi_{m+1,-}}\ne 0$.  
Hence we have proved the proposition.  
\end{proof}
We have not used the other generator
 $\Attbb m m + {m,m}$ in \eqref{eqn:IJmm}
 for the previous proposition.  
For the sake of completeness,
 we investigate its restriction
 to each of the irreducible components in \eqref{eqn:I+n2}.

\begin{proposition}
\label{prop:Ammrest}
Retain the notation as in \eqref{eqn:IJmm}.  
\begin{alignat*}{2}
&\Attbb {m}{m} + {m, m}
|_{\Pi_{m+1,-}}
&&
\equiv
0.  
\\
&\Attbb {m}{m} + {m, m}
|_{\Pi_{m,+}}
&&=
\frac{2 \pi^{m-\frac{1}{2}}}{m!} {\operatorname{Rest}}|_{\Pi_{m,+}}.  
\end{alignat*}
\end{proposition}
We also determine
 the image of the nonzero symmetry breaking operators
 $\Attbb m m + {m,m}$ and ${\operatorname{Rest}}$
 on each irreducible summand
in \eqref{eqn:I+n2}.  
\begin{proposition}
\label{prop:Ammimage}
With Convention \ref{conv:Image}, 
 we have 
\begin{align*}
{\operatorname{Image}}(\Attbb {m}{m} + {m, m}
|_{\Pi_{m,+}})
=&
{\operatorname{Image}}
 ({\operatorname{Rest}}|_{\Pi_{m,+}})
=\pi_{m,+},   
\\
{\operatorname{Image}}
 ({\operatorname{Rest}}|_{\Pi_{m+1,-}})
=&
J_+(m, m).  
\end{align*}
\end{proposition}
For the proof of Propositions \ref{prop:Ammrest} and \ref{prop:Ammimage}, 
 we use Lemma \ref{lem:AmmT} about functional equations
 with appropriate renormalizations.  
We set 
\begin{equation}
\label{eqn:cm}
   c(m):=\frac{\pi^{m}}{m!}. 
\end{equation}
\begin{proof}
[Proof of Proposition \ref{prop:Ammrest}]
It follows from the functional equation \eqref{eqn:AT2tilde}
 for the 
\index{B}{KnappSteinoperatorrenorm@ Knapp--Stein operator, renormalized---}
renormalized Knapp--Stein operator
 $\Tttbb m m m$
 that 
\[
  \Attbb {m}{m} + {m, m}
  \circ
  (c(m)\, {\operatorname{id}} - \Tttbb {m}{m}{m})=0.  
\]

On the other hand, 
 Lemma \ref{lem:161745} implies 
 that the renormalized Knapp--Stein operator satisfies 
\[
  c(m)\, {\operatorname{id}} - \Tttbb {m}{m}{m}
  =
  0\,\, {\operatorname{id}}_{\Pi_{m, +}} 
  \oplus
  2\,\, {\operatorname{id}}_{\Pi_{m+1, -}}, 
\]
which implies ${\operatorname{Image}}(c(m)\, {\operatorname{id}} - \Tttbb {m}{m}{m})=\Pi_{m+1, -}$.  
Therefore, 
 $\Attbb {m}{m} + {m, m}$ is 
 identically zero
 on the irreducible $G$-submodule $\Pi_{m+1,-}$.  

To see the second statement,
 we use Proposition \ref{prop:IJmmbase}, 
 which shows that $\Attbb {m}{m}+{m, m}|_{\Pi_{m, +}}$ must be 
 proportional to ${\operatorname{Rest}}|_{\Pi_{m,+}}$.  
Comparing the $(K,K')$-spectrum
 of the two operators $\Attbb {m}{m}+{m, m}$
 and ${\operatorname{Rest}}$
 with respect to basic $K'$- and $K$-types
 $\mub(m, +)' \hookrightarrow \mub(m, +)$ 
 (see the formula \eqref{eqn:SAiitilde} for $\Attbb m m + {m,m}$ and Lemma \ref{lem:152268}
 for ${\operatorname{Rest}}$), 
 we get the second statement.  
\end{proof}

\begin{proof}
[Proof of Proposition \ref{prop:Ammimage}]
By the functional equation \eqref{eqn:TA2tilde}, 
\[
  {\operatorname{Image}}
  (\Attbb m m + {m,m}|_{\Pi_{m,+}})
  \subset
  {\operatorname{Ker}}
  (\Ttbb m {m-1} {m})
  =\pi_{m,+}.  
\]
Since $\Attbb m m + {m,m}|_{\Pi_{m,+}}$ is nonzero, 
 and since $\pi_{m,+}$ is an irreducible $G'$-module, 
 we get the first statement.  
For the second one, 
 we compare the $(K,K')$-spectrum
 of $\Attbb m m + {m,m}$
 (see \eqref{eqn:SAiitilde})
 and that of ${\operatorname{Rest}}$
 (see \eqref{eqn:SRest}) in Lemma \ref{lem:152268}.  
\end{proof}

\subsubsection
{${\operatorname{Hom}}_{G'}(I_\delta(m,m)|_{G'},J_\varepsilon(m,m))$
 with $\delta \varepsilon=-$}
\label{subsec:IJmm-}
The case $\delta \varepsilon=-$ is much simpler 
 because the space of symmetry breaking operators
 is one-dimensional:
\[
   {\operatorname{Hom}}_{G'}(I_\delta(m,m)|_{G'},J_\varepsilon (m,m))
   =
   {\mathbb{C}} \Atbb m m - {m,m}, 
\]
 see Theorem \ref{thm:SBObasis}.  
Without loss of generality,
 we may and do assume 
 $(\delta, \varepsilon)=(+,-)$.  
The restriction of the generator $\Atbb m m - {m,m}$
 to each irreducible component in \eqref{eqn:I+n2}
 is given as follows.  
\begin{proposition}
\label{prop:161756}
Let $(G,G')=(O(2m+1,1),O(2m,1))$.  
Then we have
\begin{align*}
\Atbb m m - {m,m}|_{\Pi_{m,+}}\equiv & 0.  
\\
{\operatorname{Image}}(\Atbb m m - {m,m}|_{\Pi_{m+1,-}}) = & \pi_{m,-}.  
\end{align*}
\end{proposition}

The proof of Proposition \ref{prop:161756} relies 
 on the functional equations given in Lemma \ref{lem:AmmT-}.  

\begin{proof}
The functional equation \eqref{eqn:ATAm-} implies
\[
   \Atbb m m - {m,m} \circ (\Tttbb m m m + c(m){\operatorname{id}})
   =0.   
\]
By Lemma \ref{lem:161745}, 
 ${\operatorname{Image}}(\Tttbb m m m + c(m){\operatorname{id}})
 =\Pi_{m,+}$.  
Hence the first statement is proved.

The second statement follows from the functional equation \eqref{eqn:TATm-}
 and ${\operatorname{Ker}}(\Ttbb m {m-1} m: J_-(m,m) \to J_-(m,m-1))
 =\pi_{m,-}$
 (see Proposition \ref{prop:Timage}). 
\end{proof}

\subsection{Splitting of $J_{\varepsilon}(m,m)$
 and symmetry breaking operators for 
 $(G,G')=(O(2m+2,1),O(2m+1,1))$}
\label{subsec:2m+1}

Suppose $n$ is odd, 
say $n=2m+1$.  
In contrast to the $n$ even case treated
 in Section \ref{subsec:2m}, 
 a distinguished feature in this setting
 is that the principal series representation
 $J_{\varepsilon}(m,\nu)$ of the subgroup $G' =O(2m+1,1)$
 splits into the direct sum
 of two irreducible 
 tempered representations
 when $\nu=m$: for $\varepsilon\in \{\pm\}$, 
\begin{equation}
\label{eqn:JmmsplitI+n2}
J_{\varepsilon}(m,m) \simeq \pi_{m,\varepsilon} \oplus \pi_{m+1,-\varepsilon}, 
\end{equation}
 see Theorem \ref{thm:LNM20} (1) and (8).  
Accordingly, 
 the space of symmetry breaking operators has a direct sum decomposition:
\begin{multline}
\label{eqn:m+1mdeco}
{\operatorname{Hom}}_{G'}(I_{\delta}(i,\lambda)|_{G'}, J_{\varepsilon}(m, m))
\\
\simeq 
{\operatorname{Hom}}_{G'}(I_{\delta}(i,\lambda)|_{G'}, \pi_{m,\varepsilon})
\oplus 
{\operatorname{Hom}}_{G'}(I_{\delta}(i,\lambda)|_{G'},
 \pi_{m+1,-\varepsilon})
\end{multline}
for any $\lambda \in {\mathbb{C}}$.  
The left-hand side of \eqref{eqn:m+1mdeco} is understood
 via explicit generators
 given in Theorem \ref{thm:SBObasis} (classification).  
In this section,
 we examine the following two cases:
\index{A}{Ciiln@$\Cbb \lambda \nu {i,j}$, matrix-valued differential operator}
\begin{align}
\label{eqn:IJnoddmm}
{\operatorname{Hom}}_{G'}(I_{\delta}(m+1,m)|_{G'}, J_{\delta}(m, m))
= &
{\mathbb{C}}
\Attbb m m + {m+1,m}
\oplus 
{\mathbb{C}}
\Ctbb m m {m+1,m},
\\
\label{eqn:IJnoddmm-}
{\operatorname{Hom}}_{G'}(I_{-\varepsilon}(m,m)|_{G'}, J_{\varepsilon}(m, m))
= &
{\mathbb{C}} \Atbb m m -{m,m}, 
\end{align}
in connection with the decomposition 
 in the right-hand side of \eqref{eqn:IJnoddmm}.

We retain the notation \eqref{eqn:cm}
 in the previous section,
 that is, 
\[
  c(m) = \frac{\pi^m}{m!}.  
\]
Then the irreducible $G'$-modules 
 $\pi_{m,\varepsilon}$ and $\pi_{m+1,-\varepsilon}$
 in \eqref{eqn:JmmsplitI+n2}
 are the eigenspaces
 of the renormalized Knapp--Stein operator
 $\Tttbb m m m$ for the subgroup $G'$
 with eigenvalues $c(m)$ and $-c(m)$, 
respectively, 
 by Lemma \ref{lem:161745}.

The case \eqref{eqn:IJnoddmm} will be discussed
 in Section \ref{subsec:IJnoddmm}
 and the case \eqref{eqn:IJnoddmm-}
 in Section \ref{subsec:IJnoddmm-}.  
In particular,
 we shall see in Section \ref{subsec:SBOrho-}, 
 that both 
 ${\mathbb{A}}':= \frac 1 2 \Attbb m m + {m+1,m}+
 c(m) \Ctbb m m {m+1,m}$
 in \eqref{eqn:IJnoddmm}
 and $\frac 1 2 (-1)^{m+1}\Atbb m m - {m,m}$
 in \eqref{eqn:IJnoddmm-}
 yield the same symmetry breaking operator
\[
  A_{m+1,m} \colon \Pi_{m+1,\delta} \to \pi_{m,\delta}, 
\]
 which will be utilized in the construction 
 of nonzero periods in Chapter \ref{sec:period}, 
 see Theorem \ref{thm:period2}.  

\subsubsection
{${\operatorname{Hom}}_{G'}(I_{\delta}(m+1,m)|_{G'}, J_{\delta}(m,m))$
 for $n=2m+1$}
\label{subsec:IJnoddmm}
We recall from Theorem \ref{thm:161243} (2)
 that the regular symmetry breaking operator
 $\Atbb \lambda \nu +{i,j}$ vanishes
 when $(n,i,j, \lambda, \nu)=(2m+1,m+1,m,m,m)$, 
and therefore, 
 the left-hand side of \eqref{eqn:m+1mdeco} 
 at $\lambda=m$ is two-dimensional 
 by Theorem \ref{thm:1.1} (2).  
More precisely, 
 the classification
 of symmetry breaking operators
 given in Theorem \ref{thm:SBObasis} shows \eqref{eqn:IJnoddmm}.

On the other hand, 
 we recall from Theorem \ref{thm:LNM20} (1)
 that the principal series representation $I_{\delta}(m+1,m)$ has
 a nonsplitting exact sequence of $G$-modules:
\[
 0 \to \Pi_{m+2,-\delta} \to I_{\delta}(m+1,m) \to \Pi_{m+1,\delta} \to 0.  
\]
The irreducible $G$-submodule $\Pi_{m+2,-\delta}$ is the image
 of the Knapp--Stein operator $\Ttbb {m+1}m {m+1}$ for the group $G$.  
With this in mind,
 we shall take a closer look at the right-hand side of 
 \eqref{eqn:m+1mdeco}.

We introduce the following element in \eqref{eqn:IJnoddmm}:
\begin{equation}
\label{eqn:Aprimm}
  {\mathbb{A}}'
:=
  \frac 1 2 \Attbb m m + {m+1,m}
  +
  c(m) \Ctbb m m {m+1,m}.  
\end{equation}
The main result of this subsection is the following.

\begin{proposition}
\label{prop:1802103}
Let $(G,G')=(O(2m+2,1),O(2m+1,1))$.  
Then ${\mathbb{A}}' \colon I_{\delta}(m+1,m) \to J_{\delta}(m,m)$ is 
 a symmetry breaking operator satisfying 
\begin{align*}
 {\mathbb{A}}'  \circ \Ttbb {m+1} m {m+1} =&0, 
\\
 \Tttbb m m {m,m} \circ {\mathbb{A}}' 
  =& c(m) {\mathbb{A}}',  
\\
 S({\mathbb{A}}' 
  )=& c(m) \begin{pmatrix} 1 & 0 \\ 0 & 0 \end{pmatrix}.  
\end{align*}
\end{proposition}

Proposition \ref{prop:1802103} follows from the corresponding results
 for the renormalized operator $\Attbb m m + {m+1,m}$
 (Lemma \ref{lem:reAm+1m} below)
 and for the differential operator $\Ctbb m m {m+1,m}$
 (Lemmas \ref{lem:KKspecC}, \ref{lem:TCAC}, and \ref{lem:CTAtilde}).  
We begin with the functional equations
 and the $(K,K')$-spectrum
 of the first generator 
 $\Attbb m m + {m+1,m}$ in \eqref{eqn:IJnoddmm}.  
\begin{lemma}
\label{lem:reAm+1m}
Retain the setting 
 where $(G,G')=(O(2m+2,1),O(2m+1,1))$.  
Then the renormalized regular symmetry breaking operator
 $\Attbb m m + {m+1,m}$ satisfies the following:
\begin{align*}
 \Attbb m m + {m+1,m} \circ \Ttbb {m+1} m {m+1} 
  =& -2 c(m) \Attbb {m+1}m+{m+1,m}, 
\\
 \Tttbb {m} m {m}  \circ \Attbb m m + {m+1,m}
 =& -c(m) \Attbb m m + {m+1,m}, 
\\
 S(\Attbb m m + {m+1,m})=& c(m) \begin{pmatrix} 0 & 0 \\ 0 & 2 \end{pmatrix}.  
\end{align*}
\end{lemma}

\begin{proof}
See Lemma \ref{lem:Aii-ren} for the first and third equalities, 
 and Lemma \ref{lem:TAT3tilde} for the second.  
\end{proof}

For the differential symmetry breaking operator
 $\Ctbb m m {m+1,m}$ in \eqref{eqn:IJnoddmm},
 we recall from \cite[Thm.~1.3]{xkresidue}
 the 
\index{B}{residueformula@residue formula}
 residue formula of the regular symmetry breaking operators
 $\Atbb \lambda \nu {\delta \varepsilon}{i,j}
 \colon I_{\delta}(i,\lambda) \to J_{\varepsilon}(j,\nu)$
 when 
\index{A}{1psi@$\Psising$,
          special parameter in ${\mathbb{C}}^2 \times \{\pm\}^2$}
$(\lambda, \nu,\delta, \varepsilon) \in \Psising$
 for $j=i-1$ and $i$, 
 see Fact \ref{fact:153316}.  
Applying \eqref{eqn:Aijres1} to 
 $(n, i,j,\lambda,\nu)=(2m+1,m+1,m,\lambda,\lambda)$, 
 we obtain 
\begin{equation}
\label{eqn:resAC2}
\Atbb \lambda \lambda + {m+1,m}
=
\frac{(m-\lambda) \pi^{m}}{\Gamma(\lambda+1)}\Ctbb \lambda \lambda {m+1,m}, 
\end{equation}
where we recall from \eqref{eqn:Cijln}
 that the differential symmetry breaking operator
 $\Cbb \lambda \nu {i,i-1}$ vanishes
 for the parameter
 that we are dealing with, 
 namely,
 when $\lambda=\nu=n-i$.  
So we use the renormalized operator 
 $\Ctbb \lambda \nu {i,i-1}$ instead.  
We note that 
$
   \Ctbb \lambda \lambda {i,i-1}
  ={\operatorname{Rest}}_{x_n=0}\circ \iota_{\frac\partial{\partial x_n}}
$.

\begin{lemma}
\label{lem:KKspecC}
The $(K,K')$-spectrum of $\Ctbb m m {m+1,m}$ is given by
\[
   S(\Ctbb m m {m+1,m})
   = \begin{pmatrix} 1 & 0 \\ 0 & -1 \end{pmatrix}.  
\]
\end{lemma}
\begin{proof}
By the residue formula \eqref{eqn:resAC2}, 
 we have
\[
  \lim_{\lambda \to m}
  \frac 1 {\lambda-m} \Atbb \lambda \lambda + {m+1,m}
  = 
  -c(m)\Ctbb m m {m+1,m}.  
\]
Now the lemma follows from the $(K,K')$-spectrum
 of the regular symmetry breaking operator
 $\Atbb \lambda \nu {\pm} {i,j}$
 given in Theorem \ref{thm:153315}.  
\end{proof}

The symmetry breaking operator $\Atbb \lambda \nu + {m+1,m}$
 vanishes at $(\lambda,\nu)=(m,m)$.  
We recall from Lemma \ref{lem:upsemi}
 and Definition \ref{def:AVWder}
 that 
\[
  (\Atbb m m + {m+1,m})_{k,l}
  :=
   \left.
   \frac{\partial^{k+l}}{\partial \lambda^{k}\partial \nu^{l}}
   \right|
   _{\substack{\lambda=m \\ \nu=m}}
   \Atbb \lambda \nu + {m+1,m}
   \in {\operatorname{Hom}}_{G'}
   (I_{\delta}(m+1,m)|_{G'}, J_{\delta}(m,m))
\]
 for $(k,l)=(1,0)$ and $(0,1)$.

The base change of the vector space
 ${\operatorname{Hom}}_{G'}
   (I_{\delta}(m+1,m)|_{G'}, J_{\delta}(m,m))$, 
 see \eqref{eqn:IJnoddmm}, 
 is given as follows.  

\begin{lemma}
\label{lem:twobases}
\begin{enumerate}
\item[{\rm{(1)}}]
$2 (\Atbb m m + {m+1,m})_{1,0}= \Attbb m m + {m+1,m}$, 
\item[{\rm{(2)}}]
$(\Atbb m m + {m+1,m})_{1,0} + (\Atbb m m + {m+1,m})_{0,1}
 = -c(m) \Ctbb m m {m+1,m}$.  
\end{enumerate}
\end{lemma}
\begin{proof}
The first assertion is immediate from the definition 
 of the renormalized operator 
 $\Attbb m m + {m+1,m}$, 
 see \eqref{eqn:Ader}.  
The second assertion follows from the residue formula 
 \eqref{eqn:resAC2}.  
\end{proof}
It follows from Lemma \ref{lem:twobases}
 that 
\begin{equation}
\label{eqn:AClimit}
  \lim_{\lambda \to m}
  \frac 1 {\lambda-m} \Atbb \lambda {2m-\lambda} + {m+1,m}
  = 
  \Attbb m m + {m+1,m}
  +c(m) \Ctbb m m {m+1,m}.  
\end{equation}

Now we give functional equations
 of the differential symmetry breaking operators
 $\Ctbb m m {m+1,m}$
 and the (renormalized) Knapp--Stein operators for $G'$ and $G$
 as follows.

\begin{lemma}
\label{lem:TCAC}
$\Tttbb m m m \circ \Ctbb m m {m+1,m}
=
\Attbb m m + {m+1,m} + c(m) \Ctbb m m {m+1,m}$.  
\end{lemma}
\begin{proof}
By the functional equation in Theorem \ref{thm:TAA}, 
 we have
\[
 \Ttbb \lambda {2m-\lambda} {m}
 \circ 
 \Atbb \lambda \lambda + {m+1,m}
  = 
  \frac{(m-\lambda)\pi^{m}}{\Gamma(\lambda+1)}
  \Atbb \lambda {2m-\lambda} +  {m+1,m}.  
\]

Hence we get from the residue formula \eqref{eqn:resAC2}
\[
  \Ttbb \lambda {2m-\lambda} {m}
 \circ 
  \Ctbb \lambda \lambda {m+1,m}
  = 
  \Atbb \lambda {2m-\lambda} +  {m+1,m}.  
\]
Now Lemma \ref{lem:TCAC} follows from \eqref{eqn:AClimit}.  
\end{proof}

\begin{lemma}
\label{lem:CTAtilde}
$
\Ctbb m m {m+1,m} \circ \Ttbb {m+1} m {m+1} 
=
\Attbb {m+1} m + {m+1,m}.  
$
\end{lemma}
\begin{proof}
By the functional equation in Theorem \ref{thm:ATA}, 
 we have
\[
 \Atbb \lambda \lambda + {m+1,m}
 \circ 
 \Ttbb {2m+1-\lambda}\lambda {m+1}
 = 
  \frac{\pi^{m+\frac 1 2}(m-\lambda)}{\Gamma(\lambda+1)}
  \Atbb {2m+1-\lambda}\lambda +  {m+1,m}.  
\]
By the residue formula \eqref{eqn:resAC2}
 and by analytic continuation, 
 we get 
\[
  \Ctbb \lambda \lambda {m+1,m}
  \circ 
  \Ttbb {2m+1-\lambda} \lambda {m+1}
  = 
  \pi^{\frac 12}\Atbb {2m+1-\lambda}\lambda+{m+1,m}.  
\]
Since $\Attbb {m+1} m + {m+1,m} = \pi^{\frac 1 2} \Atbb {m+1} m + {m+1,m}$
 by the definition \eqref{eqn:Aijrenorm}
 of the renormalized operator $\Attbb {\lambda} \nu \pm {i,j}$, 
 the lemma is proved.  
\end{proof}

\subsubsection
{${\operatorname{Hom}}_{G'}(I_{-\varepsilon}(m,m)|_{G'}, J_{\varepsilon}(m,m))$
 for $n=2m+1$}
\label{subsec:IJnoddmm-}

In this subsection,
 we examine 
\[
  {\operatorname{Hom}}_{G'}
  (I_{-\varepsilon}(m,m)|_{G'}, J_{\varepsilon}(m,m))
 ={\mathbb{C}} \Atbb m m - {m,m}, 
\]
as stated in \eqref{eqn:IJnoddmm-}, 
 which is derived from Theorems \ref{thm:1.1} and \ref{thm:SBObasis}.  
We recall from Theorem \ref{thm:LNM20} (1)
 that there is a nonsplitting exact sequence of $G$-modules:
\[
   0 \to \Pi_{m,-\delta} \to I_{-\delta}(m,m) \to \Pi_{m+1,\delta} \to 0.  
\]
Concerning the regular symmetry breaking operator $\Atbb m m -{m,m}$, 
 we have the following.  
\begin{lemma}
\label{lem:IHP080328}
Let $(G,G')=(O(2m+2,1),O(2m+1,1))$.  
Then we have
\begin{align*}
   \Tttbb m m {m} \circ \Atbb m m - {m,m}
 =& c(m) \Atbb m m - {m,m}, 
\\
   \Atbb m m - {m,m} \circ \Ttbb {m+1} m {m} =&0, 
\\
   S(\Atbb m m - {m,m})=& 2(-1)^{m+1} c(m) \begin{pmatrix} 0 & 0 \\ 1 & 0\end{pmatrix}.  
\end{align*}
\end{lemma}

\begin{proof}
The proof of first formula parallels to that of Lemma \ref{lem:TAT3tilde}, 
 and the second formula is a special case of Theorem \ref{thm:ATA}.  
The third formula follows from Theorem \ref{thm:153315}.  
\end{proof}

\subsection{Symmetry breaking operators from $\Pi_{i,\delta}$ to $\pi_{i-1,\delta}$}
\label{subsec:SBOrho-}

In Sections \ref{subsec:StageCii} and \ref{subsec:2m}, 
 we constructed nontrivial symmetry breaking operators from 
 the irreducible representation $\Pi_{i,\delta}$ of $G=O(n+1,1)$
 to the irreducible one $\pi_{i,\delta}$ of $G'=O(n,1)$.  
This is sufficient for the proof of Theorem \ref{thm:SBOone}
 by the duality theorem
 (Proposition \ref{prop:SBOdual})
 between symmetry breaking operators
 for the indices:
\[
   \text{$(i,j)$ and $(\widetilde i,\widetilde j):=(n+1-i, n-j)$}.  
\]
Nevertheless,
 we give in this section an explicit construction
 of the normalized symmetry breaking operators
 $\Pi_{i,\delta} \to \pi_{j,\varepsilon}$
 also for $j=i-1$, 
 and determine their $(K,K')$-spectrum 
 of symmetry breaking operators from $\Pi_{i,\delta}$ to $\pi_{i-1,\delta}$.  
The results will be used in the computation
 of {\it{periods}} of admissible smooth representations 
 in Chapter \ref{sec:period}.

We begin with some basic properties
 of the regular symmetry breaking operator
\[
   \Atbb \lambda \nu \delta {i,i-1} \colon
 I_\delta(i,\lambda) \to J_\delta (i-1,\nu)
\]
for $(\lambda,\nu)=(n-i,i-1)$.  

\begin{proposition}
\label{prop:Aii1ker}
Suppose $1 \le i \le n$ and $\delta \in \{\pm\}$.  
\begin{enumerate}
\item[{\rm{(1)}}]
$\Pi_{i+1,-\delta} \subset{\operatorname{Ker}}~(\Atbb {n-i}{i-1} +{i,i-1})$.  
\item[{\rm{(2)}}]
${\operatorname{Image}}~(\Atbb {n-i}{i-1}+{i,i-1}) \simeq \pi_{i-1,\delta}$
 if $n \ne 2i-1 ;$
$\Atbb {n-i}{i-1}+{i,i-1} =0$
 if $n = 2i-1$.  
\end{enumerate}
\end{proposition}
\begin{proof}
(1)\enspace
Applying the functional equation
 given in Theorem \ref{thm:ATA} with $\lambda=i$, 
 we see that the symmetry breaking operator $\Atbb {n-i}{\nu} \gamma{i,i-1}$ vanishes 
 on the image of the Knapp--Stein intertwining operator
 $\Ttbb i {n-i} i \colon I_{\delta}(i,i) \to I_{\delta}(i,n-i)$, 
 namely, 
 on the irreducible submodule $\Pi_{i+1,-\delta}$
 (see Theorem \ref{thm:LNM20} (1)).  
\newline\noindent
{\rm{(2)}}\enspace
By Theorem \ref{thm:161243}, 
 $\Atbb {n-i}{i-1}+{i,i-1}=0$
 if and only if $n=2i-1$.

Suppose from now 
 that $n \ne 2i-1$.  
Applying the functional equation given in Theorem \ref{thm:TAA} with $(\lambda, \nu, \gamma)=(n-i,i-1,+)$, 
 we see that the composition 
 $\Ttbb {\nu} {n-1-\nu} {i-1} \circ \Atbb {\lambda}{\nu} +{i,i-1}$
 is a scalar multiple of the symmetry breaking operator $\Atbb {n-i}{n-i} +{i,i-1}$, 
 which vanishes by Theorem \ref{thm:161243}.  
In turn, 
 applying Proposition \ref{prop:Timage} to $G'=O(n,1)$, 
 we get
\[
    {\operatorname{Ker}}(\Ttbb \nu{n-1-\nu}{i-1} \colon
                         J_{\delta}(i-1,\nu) \to J_{\delta}(i-1,n-1-\nu))
    \simeq \pi_{i-1,\delta}
\]
 because $i-1 \ne \frac 1 2(n-1)$. 
Hence the second statement is also proved.  
\end{proof}

Since $\Atbb {n-i} {i-1}+ {i,i-1}=0$
 for $n=2i-1$, 
 we treat this case separately as follows.  
Suppose $n=2m+1$.  
We recall that there are a nonsplitting exact sequence of $G$-modules
\[
  0 \to \Pi_{m,-\delta} \to I_{-\delta}(m,m) \to \Pi_{m+1,\delta} \to 0
\]
 and a direct sum decomposition of irreducible $G'$-modules
\[
   J_{\delta}(m,m) \simeq \pi_{m,\delta} \oplus \pi_{m+1,-\delta}.  
\]
We use the following regular symmetry breaking operator
\[
  \Atbb m m - {m,m} \colon I_{-\delta}(m,m) \to J_{\delta}(m,m).  
\]
\begin{proposition}
\label{prop:IHP180328}
Suppose $(G,G')=(O(2m+2,1),O(2m+1,1))$
 and $\delta \in \{\pm\}$.  
\begin{enumerate}
\item[{\rm{(1)}}]
${\operatorname{Ker}}(\Atbb m m - {m,m}) \supset \Pi_{m,-\delta}$.  
\item[{\rm{(2)}}]
${\operatorname{Image}}(\Atbb m m - {m,m})=\pi_{m,\delta}$.  
\end{enumerate}
\end{proposition}
\begin{proof}
The assertions follow from Lemma \ref{lem:IHP080328}.  
\end{proof}

It follows from Proposition \ref{prop:Aii1ker}
 that if $n \ne 2i-1$ then the normalized symmetry breaking operator $\Atbb {n-i}{n-i} +{i,i-1}$
 yields a surjective $G'$-homomorphism
\begin{equation}
\label{eqn:Pipi-}
  A_{{i},{i-1}} \colon \Pi_{i,\delta} \to \pi_{i-1,\delta}
\end{equation}
by the following diagram.  
\begin{eqnarray*}
\xymatrix@C=16pt{
&&I_{\delta}(i,n-i)\ar[rr]^{\Atbb {n-i}{n-i}+{i,i-1}}
\ar@{}[dr]^\circlearrowleft
\ar@{->>}[d]
&&\pi_{i-1,\delta} \ar@{^{(}->}[rr]
&&J_{\delta}(i-1,i-1)
\\
\Pi_{i,\delta}
&&\ar[ll]_{\sim\qquad\quad} I_{\delta}(i,n-i)/\Pi_{i+1,-\delta}\ar@{-->}[urr]
&&
&&
}
\end{eqnarray*}
If $n=2i-1$, 
 we set $(n,i)=(2m+1,m+1)$.  
Then, 
similarly to the case $n \ne 2i-1$, 
 Proposition \ref{prop:1802103} shows
 that the symmetry breaking operator 
 ${\mathbb{A}}' \colon I_{\delta}(m+1,m) \to J_{\delta}(m,m)$ 
 defined in \eqref{eqn:Aprimm}
 yields a surjective $G'$-homomorphism
\begin{equation}
\label{eqn:Pmpm-}
   A_{m+1,m} \colon \Pi_{m+1,\delta} \to \pi_{m,\delta}
\end{equation}
 by the following diagram.  
\begin{eqnarray*}
\xymatrix@C=16pt{
&&I_{\delta}(m+1,m)\ar[rr]^{{\mathbb{A}}'} 
\ar@{->}[d]
\ar@{}[dr]^\circlearrowleft
&&\pi_{m,\delta} 
\subset J_{\delta}(m,m)
\\
\Pi_{m+1,\delta}
&&\ar[ll]_{\sim\quad\quad} I_{\delta}(m+1,m)/\Pi_{m+2,-\delta}\ar@{-->}[urr]
&&
}
\end{eqnarray*}

In order to define the $(K,K')$-spectrum,
 we need to fix an inclusive map from 
 the $K'$-type into the $K$-type,
 see Definition \ref{def:KKspec}.  
In our setting, 
 we use the natural embedding of the minimal $K$- and $K'$-types
\begin{equation}
\label{eqn:KKemb}
  \mub (i,\delta) \hookleftarrow \mub (i-1,\delta)'
\end{equation}
 of the irreducible representations $\Pi_{i,\delta}$ and $\pi_{i-1,\delta}$
 of $G$ and the subgroup $G'$, 
 respectively,
 as in Section \ref{subsec:Apm1-5}.  
Then we get the following formula
 for the 
\index{B}{KspectrumKprime@$(K,K')$-spectrum}
$(K,K')$-spectrum.  

\begin{proposition}
\label{prop:Aidown}
Let $(G, G')=(O(n+1,1), O(n,1))$ and $1 \le i \le n+1$. 
Then the symmetry breaking operator
\[
     A_{i,i-1} \colon \Pi_{i,\delta} \to \pi_{i-1,\delta}
\]
 acts on $\mub(i-1,\delta)'$ $(\hookrightarrow \mub(i,\delta))$
 as the following scalar:
\begin{equation*}
\begin{cases}
 \frac{\pi^{\frac{n-1}{2}}(n-2i+1)}{(n-i)!}
\qquad
 &\text{if $n \ne 2i-1$}, 
\\ 
 \frac{\pi^{\frac{n-1}{2}}}{(n-i)!}
  & \text{if $n=2i-1$}.
\end{cases}
\end{equation*}
\end{proposition}
\begin{proof}
For $n \ne 2i-1$, 
 the assertion follows directly from the $(1,1)$-component
 of the matrix $S(\Atbb{\lambda}{\nu}{+}{i,i-1})$
 in Theorem \ref{thm:153315}
 with $(\lambda,\nu)=(n-i,i-1)$.  
For $n=2i-1$, 
 the (1,1)-component of $S({\mathbb{A}}')$ in Proposition \ref{prop:1802103}
 with $(n,i)=(2m+1,m+1)$ shows the desired formula.  
\end{proof}

\begin{remark}
\label{rem:20180406}
When $n=2i-1$, 
 we set $(n,i)=(2m+1,m+1)$ as above.  
In this case we may use $\Atbb m m - {m,m}$ in Lemma \ref{lem:IHP080328}
 for an alternative construction 
 of $A_{m+1,m} \in {\operatorname{Hom}}_{G'}(\Pi_{m+1,\delta}|_{G'}, \pi_{m,\delta})$.  
To see this, 
 we recall from Section \ref{subsec:IJnoddmm-}
 the following natural inclusion 
\[
{\operatorname{Hom}}_{G'}(\Pi_{m+1,\delta}|_{G'}, \pi_{m,\delta})
\subset
 {\operatorname{Hom}}_{G'}(I_{-\delta}(m,m)|_{G'}, J_{\delta}(m,m))
 ={\mathbb{C}} \Atbb m m - {m,m}, 
\]
and therefore any element
 in 
$
   {\operatorname{Hom}}_{G'}(\Pi_{m+1,\delta}|_{G'}, \pi_{m,\delta})
$
 is proportional to the one
 which is induced from $\Atbb m m - {m,m}$.  
On the other hand, 
 Proposition \ref{prop:IHP180328} tells 
 that the symmetry breaking operator 
 $\Atbb m m - {m,m}$ yields
 a surjective $G'$-homomorphism
 $\Pi_{m+1,\delta} \to \pi_{m,\delta}$
 by the following diagram.  
\begin{eqnarray*}
\xymatrix@C=22pt{
&&I_{-\delta}(m,m)\ar[rr]^{\quad\Atbb m m -{m,m}}
\ar@{->>}[d]
\ar@{}[dr]^\circlearrowleft
&&\pi_{m,\delta} \ar@{^{(}->}[rr]
&&J_{\delta}(m,m)
\\
\Pi_{m+1,\delta}
&&\ar[ll]_{\sim\quad\quad} I_{-\delta}(m,m)/\Pi_{m,-\delta}\ar@{-->}[urr]
&&
&&
}
\end{eqnarray*}
By Lemma \ref{lem:IHP080328}, 
 $\frac 1 2 (-1)^{m+1} \Atbb m m - {m,m}$ has 
 the $(K,K')$-spectrum 
 for the basic $K$- and $K'$-types
\[
  S\left(\frac 1 2 (-1)^{m+1} \Atbb m m - {m,m}\right) = c(m) \begin{pmatrix} 0 & 0 \\ 1 & 0\end{pmatrix}.  
\]
In view of the (2,1)-component, 
 the resulting symmetry breaking operator from $\Pi_{m+1,\delta}$ to $\pi_{m,\delta}$
 has the $(K,K')$-spectrum $c(m)$
 for the embedding 
 of the $K$- and $K'$-types
 $\mub (m+1,\delta) \hookleftarrow \mub (m,\delta)'$.  
This is the same with the $(K,K')$-spectrum of $A_{m+1,m}$
 which is induced from 
 ${\mathbb{A}}' \in {\operatorname{Hom}}_{G'}(I_{\delta}(m+1,m)|_{G'}, J_{\delta}(m,m))$.  
Hence $\frac 1 2 (-1)^{m+1} \Atbb m m - {m,m}$ induces
 the same symmetry breaking operator with $A_{m+1,m}$.  
\end{remark}

\newpage
\section{Application I: Some conjectures by B.~Gross and D.~Prasad:  Restrictions of tempered representations of $SO(n+1,1)$ to $SO(n,1)$}
\label{sec:Gross-Prasad}

Inspired by automorphic forms and $L$-functions,
 B.~Gross and D.~Prasad published in 1992 conjectured
 about the restriction of irreducible
\index{B}{temperedrep@tempered representation}
 tempered representations
 of special orthogonal groups $SO(p+1,q)$
 to a special orthogonal subgroup $SO(p,q)$, 
 see \cite{GP}.  
B.~Sun and C.-B.~Zhu \cite{SunZhu} proved
 that in this case the multiplicities are at most one,
 and B.~Gross and D.~Prasad conjectured that given a Vogan packet of tempered representations of $SO_{n+2} \times SO_{n+1}$ 
 there exist exactly one group $SO(p+1,q) \times SO(p,q)$
 with $p+q=n+1$
 and one (tempered) representation $\overline{U}_1 \boxtimes \overline{U}_2$
 of this group
 with $m( \overline{U}_1 \boxtimes \overline{U}_2,\mathbb C)=1$.
They also stated a conjectured algorithm
 to determine the  group
 and the  representation $\overline{U}_1 \boxtimes \overline{U}_2$
in the  Vogan packet with   $m( \overline{U}_1 \boxtimes \overline{U}_2,\mathbb C)=1$.

In this chapter
 we prove that the algorithm  of B.~Gross and D.~Prasad predicts the multiplicity correctly  for representations in  Vogan packets of  tempered principal series representations 
 of $SO(n+1,1) \times SO(n,1)$ 
 as well as for the 3 irreducible representations
 $\overline{\Pi}, \overline{\pi}, \overline{\varpi}$ of 
 $SO(2m+2,1)$, $SO(2m+1,1)$, $SO(2m,1)$
 with trivial infinitesimal character $\rho$.

\medskip
The Gross--Prasad conjectures are stated only for representations of special orthogonal groups in  \cite{GP}. 
Thus we are considering in this chapter symmetry breaking for tempered representations of 
\index{A}{GOindefinitespecial@$\overline{G}=SO(n+1,1)$}
$\overline{G}\times \overline{G'}=SO(n+1,1) \times SO(n,1)$ and not as in the previous chapters for $G \times G'=O(n+1,1) \times O(n,1)$. 
We refer to  Appendix II
 (Chapter \ref{sec:SOrest}) for notation and 
 for results about the restriction of representations from orthogonal groups to special orthogonal groups.
 
\subsection{Vogan packets of tempered induced representations}
\label{subsec:Vpacket}

We use a bar
 over representations to distinguish
 between representations of the special orthogonal group
 and those of the orthogonal group.

\medskip
 Every tempered principal series representation of $SO(n+1,1)$ is of the form
\[
  \overline{I}_{\delta}(\overline{V},\lambda) 
  \equiv
  \mbox{Ind}_{\overline{P}}^{\overline{G}}(\overline{V} \boxtimes \delta ,\lambda ) 
\quad \text{for }\,\,(\overline{\sigma},\overline{V}) \in \widehat{SO(n)}, 
\,\, \delta \in \{\pm\}, 
\,\, \lambda \in \frac n 2 + \sqrt{-1}{\mathbb{R}}, 
\]
which is the smooth representation
 of a unitarily induced principal series representation from a finite-dimensional representation of the minimal parabolic subgroup 
\index{A}{PLanglandsdecompbar@$\overline P = \overline M A N_+$}
$\overline P$
 of $\overline G=SO(n+1,1)$.

For $n$ even,
 we assume that the central element $-I_{n+2}$
 of the special orthogonal group $\overline G=SO(n+1,1)$ acts nontrivially
 on the principal series representation  
 ${\overline{I}}_{\delta}(\overline{V},\lambda)$, 
 and thus $\overline I_{\delta}(\overline V, \lambda)$ is a genuine representation
 of $\overline G$,
{\it{i.e.}},
 that $-I_{n+2}$ is not in the kernel of $\overline{V} \boxtimes \delta.$ 
For $n$ odd,
 $\overline G=SO(n+1,1)$ does not have
 a nontrivial center,
 and we do not need an assumption
 on the pair $(\overline{V},\delta)$.

We observe if $n$ is odd,
 the Langlands parameter of the representations
 of $SO(n,0)$ factors through the identity component of its $L$-group,
 and it defines a representation of $SO(n-2p,2p)$ and not of $O(n-2p,2p)$, 
 see \cite{Arthur}.

The Langlands parameter  of the induced representations $ \overline{I}_{\delta}(\overline{V},\lambda)$ factors through the Levi subgroup of a maximal parabolic subgroup of the Langlands dual group 
\index{A}{GOindefinitesprimelang@${}^L G$, Langlands dual group|textbf}
$^L G$ \cite{LLNM}. 
This parabolic subgroup corresponds to a maximal parabolic subgroup
 of $SO(n+1,1)$
 whose Levi subgroup $L$ is a real form of $SO(n,\bC) \times SO(2,\bC)$ and thus is isomorphic to $ SO(n,0) \times SO(1,1)\simeq SO(n) \times GL(1,{\mathbb{R}})$. 
Note that $SO(1,1)\simeq GL(1,{\mathbb{R}})$ is a disconnected group and so determines the character  $\delta $.

The pure inner real forms of $SO(n,\bC)$
 with a compact Cartan subgroup are $SO(n-2p,2p)$, $0 \leq p \leq \frac n2$.  
For $n$ even, 
 we  assume that the center of $SO(n-2p,2p)$ is not contained
 in the kernel of the discrete series representation, 
 see Proposition \ref{prop:161648} (6).

By \cite[p.~35]{A1},
 if $G$ is $SO(2m+2,1)$ or $SO(2m+1,1)$,
 then there are $2^m$ representations
 in the Vogan packet containing a tempered representation
 $\overline{I}_{\delta}(\overline{V},\lambda)$
 and they are parametrized by characters
 of a finite group ${\mathcal{A}}_1 \simeq (\bZ/2\bZ)^m$.  
We write
\index{A}{VP@$VP(\cdot)$, Vogan packet|textbf}
 $VP(\overline{I}_\delta(\overline{V},\lambda))$
for this Vogan packet.

The representations in the Vogan packet 
 $VP(\overline{I}_\delta(\overline{V},\lambda))$
 can be described as  follows: 
 we call a real form $SO(\ell,k)$ of $SO(\ell + k,\bC)$ 
\index{B}{pureinnerform@pure inner form|textbf}
{\it{pure}}
 if $\ell$ is even and thus admits discrete series representations.
We consider parabolic subgroups of $SO(n-2p+1,1+2p)$ with Levi subgroups $L$,
 which are pure inner forms of
 $ SO(n) \times GL(1,{\mathbb{R}})$.  
Hence they are isomorphic to 
\[
   L \simeq SO(n-2p,2p) \times GL(1,{\mathbb{R}}).  
\]
The Vogan packet $VP(\overline{I}_\delta(\overline{V},\lambda ))$
 contains the {\bf tempered principal series representations}
 of $SO(n-2p+1,1+2p)$
 which have the same infinitesimal character
 as $\overline{I}_\delta(\overline{V},\lambda)$, 
 and which are induced from the outer tensor product of a discrete series representation of $SO(n-2p,2p)$, with the same infinitesimal character
 as $\overline{V}$  and a one-dimensional representation $\chi_\lambda$ of $GL(1,\mathbb{R})$, 
 \cite{V}.  

We use the same conventions
 for a Vogan packet $VP(\overline{J}_\varepsilon(\overline{W},\nu))$
 of the tempered principal series representation
 $\overline{J}_\varepsilon(\overline{W},\nu)$ of $\overline{G'}$.

\subsection{Vogan packets of discrete series representations with integral infinitesimal character of $SO(2m,1)$} 

We begin with the case
 $n=2m-1$.  
In this case $SO(n+1,1)=SO(2m,1)$
 has discrete series representations.  
We fix a set of positive roots 
$
   \Delta^+ \subset {\mathfrak{t}}_{\mathbb{C}}^{\ast}
$
 for the root system 
 $\Delta({\mathfrak{so}}(2m+1,\bC),{\mathfrak{t}}_{\mathbb{C}})$
 and denote by $\rho$ half the sum of positive roots as before. 
Let $\eta $ be an integral infinitesimal character, 
 which is dominant with respect to  $\Delta  ^+$.
For $\ell+k=2m+1$, 
 we call a real form $SO(\ell,k)$ {\it{pure}}
 if $\ell$ is even.
The Vogan packet containing the discrete series representation with infinitesimal character $\eta$ is the disjoint union of discrete series representations with infinitesimal character $\eta$ of the pure inner forms. 
The cardinality of this packet is 
\[
   2^m=\sum_{\substack{0 \le \ell \le 2m \\ \ell:\text{even}}} 
({m \atop{\frac \ell 2}}).  
\]
There exists a finite group 
$
  {\mathcal  A}_2 \simeq (\bZ/2\bZ)^m$  whose characters parametrize the representations in the Vogan packet. For the discrete series representation with parameter $\chi \in \widehat{{\mathcal{A}}_2}$
 we write $\overline{\pi}(\chi). $ For more details
 see \cite{GP} or \cite{V}. 
If $\overline{\pi}$ is a discrete series representation of $SO(2m,1)$ we write $VP(\overline{\pi})$ for the Vogan packet containing $\overline{\pi}$.

\begin{example}
Suppose that $\overline{\pi }$ is a discrete series representation
 of $SO(2m,1)$
 with trivial infinitesimal character $\rho$. 
{\rm{
\begin{itemize}
\item[{\rm{(1)}}] The trivial one-dimensional representation ${\bf{1}}$
 of the inner form $SO(0,2m+1)$ is in $VP(\overline{\pi })$.

\item[{\rm{(2)}}] We can define similarly a Vogan packet $VP(\overline{\pi })$ containing $(SO(1,2m),\overline{\pi })$.  
\end{itemize}
}}
\end{example}

\subsection{Embedding the group 
$\overline{G'}=SO(n-2p,2p+1)$ into the group 
 $\overline{G}=SO(n-2p+1,2p+1)$} 

To formulate the Gross--Prasad conjecture
 we have to fix an embedding of $\overline{G'}$ into $\overline{G}$. 

\medskip

We observe:
\begin{enumerate}
\item
[(1)] The quasisplit forms of the odd special orthogonal group are 
$SO(m,m+1)$ and $SO(m+1,m)$. 
The 
\index{B}{pureinnerform@pure inner form}
 pure inner forms in the same class
 as $SO(m,m+1)$
 are $SO(m-2p, m+2p+1)$ 
 and those in the same class as $SO(m+1,m)$ are  $SO(m+1-2p, m+2p)$.

\item
[(2)]  
The quasisplit forms 
 of the even special orthogonal group are $SO(m,m)$, 
 $SO(m-1,m+1)$, and $SO(m+1,m-1)$. 
The pure inner forms are $SO(n-2p, n+2p)$ and  $SO(m+1-2p, m-1 -2p)$,
respectively,
 with $p \leq \frac m 2$.
 \end{enumerate}

So
\begin{enumerate}
\item if $n=2m$, then
 the orthogonal group $SO(2m+1,1)$ is a pure inner form of $SO(m+1,m+1)$
 if $m$ is even and of $SO(m+2,m)$ if $m$ is odd; 
\item if $n= 2m-1$, then the orthogonal group $SO(2m,1)$ is a pure inner form of $SO(m+1,m)$ if $m$ is odd and of $SO(m,m+1)$ if $m$ is even.
\end{enumerate}

\medskip
We  consider an  indefinite quadric form 
\[ 
Q_{n-2p+1,2p+1}(x)
=x^2_1+ \dots + x^2_{n-2p+1} -x^2_{n-2p+2}-  \dots -x^2_{n+2}
\]
 of signature $(n-2p +1, 2p+1)$.
We assume that $n-2p +1> 0$ and identify
 $SO(n-2p,2p+1)$ with the subgroup of $SO(n-2p+1,2p+1)$ which stabilizes 
 the basis vector $e_{n-2p+1}$.  
This  allows us to  identify the Levi subgroup of the maximal parabolic subgroup of $SO(n-2p,2p+1)$ with the intersection of the corresponding maximal parabolic subgroup of $\overline{G}$. 
This embedding of $SO(n,1)$ into $SO(n+1,1)$  is conjugate to the one
 we consider in Section \ref{subsec:Xi}.  
We use this embedding in the  formulation of the Gross--Prasad conjectures.

For tempered principal series representations we 
 consider symmetry breaking operators, 
 namely, 
 $SO(n-2p,2p+1)$-homomorphisms from representations
 in 
 $VP(\overline{I}_\delta(\overline{V},\lambda))$
 to representations
 in $VP(\overline{J}_\varepsilon(\overline{W},\nu))$, 
 see Section \ref{subsec:IV.4}.

If the tempered representation
 of $\overline{G}$ or of $\overline{G'}$ is
 a discrete series representation, 
 we consider symmetry breaking from
 a Vogan packet of discrete series representations
 to a Vogan packet of tempered principal series representations,
 respectively from a Vogan packet
 of tempered principal series representations
 to a Vogan packet of discrete series representations
 (Section \ref{subsec:GPII}).

\subsection{The Gross--Prasad conjecture I: Tempered principal series  representations}
\label{subsec:IV.4}

By Theorem \ref{thm:tempVW},
 there is a nontrivial  symmetry breaking operator
 between the {\bf{tempered}} principal representations ${I}_\delta(V,\lambda)$
 of $G=O(n+1,1)$
 and ${J}_\varepsilon(W,\nu)$ of $G'=O(n,1)$ 
 if and only if $(\sigma,V) \in \widehat{O(n)}$
 and $(\tau,W) \in \widehat {O(n-1)}$ satisfy 
\[
  [V:W]=\dim_{\mathbb{C}}{\operatorname{Hom}}_{O(n-1)}
                          (V|_{O(n-1)},W) \not = 0.  
\]
An analogous result holds 
 for a pair of the {\it{special}} orthogonal groups 
 $(\overline G, \overline {G'})=(SO(n+1,1),SO(n,1))$.  
We set
\[
   [\overline{V}: \overline{W}]
   \equiv
   [\overline{V}|_{SO(n-1)}: \overline{W}]
  :=
  \dim_{\mathbb{C}}
  {\operatorname{Hom}}_{SO(n-1)}
  (\overline{V}|_{SO(n-1)},\overline{W}).  
\]
In Theorem \ref{thm:tempSO} in Appendix II 
 we prove:
\begin{theorem}
There is a nontrivial  symmetry breaking operator
 between the {\bf{tempered}} principal series representations
 $\overline{I}_\delta(\overline{V,}\lambda)$
 of $\overline{G}=SO(n+1,1)$
 and $\overline{J}_\varepsilon(\overline{W},\nu)$ of $\overline{G'}=O(n,1)$ 
 if and only if $(\overline \sigma,\overline{V}) \in \widehat{SO(n)}$
 and $(\overline \tau,\overline{W}) \in \widehat {SO(n-1)}$ satisfy 
\[
   [\overline{V}|_{SO(n-1)}: \overline{W}] \not = 0.    
\]
\end{theorem}

In their article B.~Gross and D.~Prasad presented a conjectured algorithm
 to determine the pair of representations in the Vogan packets
 $VP(\overline{I}_\delta(\overline{V},\lambda))$
 and  $VP(\overline{J}_\varepsilon(\overline{W},\nu))$
 with a nontrivial $SO(n,1)$-symmetry breaking operator. We  prove next that  the algorithm in fact predicts :
\[ 
   [\overline{V}|_{SO(n-1)}: \overline{W}] \not = 0 
\quad
\text{ if and only if }
\quad
   \mbox{Hom}_{\overline{G'}} 
  (\overline{I}_\delta(\overline{V},\lambda)|_{\overline{G'}},
   \overline{J}_\varepsilon(\overline{W},\nu)) \not = \{0\}.
\]

\medskip
\begin{observation}
A Levi subgroup $L$ with
$[L,L] = SO(r,s)$ of the maximal parabolic subgroup determines the class
 of pure inner forms of $SO(r+1,s+1)$. 
So for any algorithm to determine the pair $(SO(r+1,s), SO(r,s))$ of the groups in the Gross--Prasad conjectures it is enough to  determine the pair  of the Levi subgroups and their corresponding discrete series representations.
\end{observation}
\medskip

\noindent
{\bf First case:}
\enspace 
{\it Suppose that } $(\overline{G},\overline{G'}) = (SO(2m+1,1), SO(2m,1))$.
\\
Let $T_\bC$ be a torus in $SO(2m+2,\bC)  \times SO(2m+1,\bC)$, 
 and $X^{\ast}(T_\bC)$ the character group. 
Fix a basis 
\[
   X^*(T_{\mathbb{C}})
   = 
   \bZ e_1\oplus \bZ e_2 \oplus \dots \oplus \bZ {{e_{m+1}}} \oplus \bZ f_1\oplus \bZ f_2 \oplus \dots \oplus \bZ f_m  
\]
such that the standard root basis $\Delta_0$ is given by
\[
e_1-e_2, e_2-e_3,\dots,{{e_{m}-e_{m+1},e_{m}+e_{m+1}}}, f_1-f_2 ,f_2-f_3, \dots, f_{m-1}-f_m,f_m
\]
{if $m \ge 1$}. 

We fix $\delta, \varepsilon \in \{\pm\}$ as in Section \ref{subsec:Vpacket}.

Recall that all representations in a Vogan packet have
 the same Langlands parameter. 
We identify the Langlands parameter of the representations
 in the same Vogan packet as
\[
   (SO(2m+1,1)\times SO(2m,1), 
   \overline{I}_\delta(\overline{V}, \lambda) \boxtimes \overline{J}_\varepsilon(\overline{W},\nu))
\]
 for a pair $(\overline{V},\overline{W})$ of irreducible finite-dimensional representations with infinitesimal character
\begin{eqnarray*}
 \lefteqn{ ( v_1+m-1)e_1+(v_2+m-2)e_2
+ \dots 
+(v_m)e_m - (\lambda -m) e_{m+1}}
\\
&  & +(u_1+m-\frac 3 2 )f_1+
 (u_2+m- \frac 5 2 ) f_2 +\dots + (u_{m-1}+\frac 1 2)f_{m-1}
-(\nu - m+\frac 1 2 ) f_m, 
\end{eqnarray*}
see \eqref{eqn:ZGinfI}.  
Here $(v_1,v_2, \dots, v_m)$ is the highest weight
 of the $SO(2m)$-module $\overline{V}$,
 $(u_1,u_2,\dots, u_{m-1}) $ is the highest weight
 of the $SO(2m-1)$-module $\overline{W}$
 and the continuous parameter $\lambda-m$ and $\nu-m+\frac 12$
 are purely imaginary, 
 and thus $\overline I_{\delta}(\overline V, \lambda)$
 and $\overline J_{\varepsilon}(\overline W, \nu)$
 are (smooth) tempered principal series representations
 of $\overline G$ and $\overline{G'}$, 
 respectively.

As discussed before, 
 to determine the pair 
\[
     (SO(n-2p+1,2p+1),SO(n-2p,2p-1))
\]
 it suffices to solve this problem for the Levi subgroups. 
Hence it suffices to  consider the Langlands parameter
\begin{eqnarray*}
 \lefteqn{ ( v_1+m-2)e_1+(v_2+m-3)e_2+ \dots +(v_m)e_m }\\
&  & +(u_1+m-\frac 5 2 )f_1+
 (u_2+m- \frac 7 2 ) f_2 +\dots + (u_{m-1} + \frac 12) f_{m-1} .
\end{eqnarray*}

Let $\delta_i$ be the element which is $-1$ in the $i$-th factor
 of ${\mathcal{A}}_1$
 and equal to $1$ everywhere else, 
 and $\varepsilon_j$ the element which is $-1$
 in the $j$-th factor of ${\mathcal{A}}_2$ and $1$ everywhere else.  
Then the algorithm \cite[p.~993]{GP} determines
 $\chi_1 \in \widehat {{\mathcal{A}}_1}$
 and $\chi_2 \in \widehat {{\mathcal{A}}_2}$
 by 
\[
 \chi_1(\delta_i)= (-1)^{\#m-i+1>}
 \quad
 \mbox{ and }
 \quad 
 \chi_2(\varepsilon_j)=  (-1)^{\#m-j + \frac 1 2<}, 
\]
where $\#m-i+1>$ is the cardinality
 of the set 
\[
\{j: v_j+m-i> \text{the coefficients of $f_j$}\},
\]
 and $\# m-j+\frac 1 2 <$ is the cardinality
 of the set \[\{i:v_i +m-j-1+\frac 1 2 < \text{the coefficients of $e_i$}\}.  \]

If $\mbox{Hom}_{SO(n-1)} (\overline{V}|_{SO(n-1)},\overline{W}) \ne \{0\}$, 
 then $v_1 \leq u_1 \leq v_2 \leq \dots \leq u_{m-1} \leq |v_m|$.
Hence we deduce that both characters are  alternating characters
 if and only if 
$\mbox{Hom}_{SO(n-1)} (\overline{V}|_{SO(n-1)},  \overline{W}) \not = \{0\}$.  

\vskip 2pc
\noindent
{\bf Second case:}
\enspace 
{\it Suppose that } $(\overline{G},\overline{G'}) = (SO(2m,1), SO(2m-1,1))$.

We use the same arguments 
 for the pair 
\[
   (\overline G, \overline{G'})=(SO(2m,1), SO(2m-1,1)).
\]

\medskip
We normalize the quasisplit forms  by 
\begin{alignat*}{2}
 &SO(m+1,m) \times SO(m,m) \quad &&\mbox{if $m$ is even,} 
\\
 &SO(m,m+1) \times SO(m-1, m+1)  \quad &&\mbox{if $m$ is odd.} 
\end{alignat*}
Applying the formul{\ae} in \cite[(12.21)]{GP}, 
 we define the integers $p$ and $q$
 with $0 \leq p \leq m $ and $0 \leq q \leq m$
by
\[
 p = \# \{ i : \chi_1(    \delta_i  )  = (-1)^i  \}  \quad \mbox{ and }  
\quad  q  =  \# \{ j:  \chi_2(   \varepsilon_j     )  = (-1)^{m+j}  \}, 
\]
and we get the pure forms
\begin{alignat*}{2}
&SO(2m-2p+1, 2p) \times SO(2q ,2m-2q) \quad &&\mbox{if $m$ is even,}
\\
&SO(2p+1,2m-2p+1) \times SO(2m-2q, 2q+1) \quad &&\mbox{if $m$ is odd.}
\end{alignat*}
In our setting,  we get the pair of integers $(p,q)=(0,m)$ for $m$ even; $(p,q)=(m,0)$ for $m$ odd.
\medskip
Applying \cite[(12.22)]{GP} with correction by changing $n$ by $m$ 
 {\it{loc.~cit.}},
 we deduce
 that the alternating character $\chi$ defines
 the pure inner form 
\begin{equation*}
SO(2m+1,0) \times SO(2m,0)
\qquad
\mbox{for $m$ is even and odd.}
\end{equation*}  
Hence 
\[
   \overline G=SO(2m,1) \mbox{ and } \overline{G'}=SO(2m-1,1).   
\]
The only representation in $VP(\overline{I}_\delta(\overline{V,}\lambda))\times VP(\overline{J}_\varepsilon(\overline{W},\nu))$
 for this pair of pure inner forms  is 
\[
   \overline{I}_\delta(\overline{V},\lambda)
   \boxtimes 
   \overline{J}_\varepsilon(\overline{W},\nu). 
\]
If $\chi $ is not the alternating character,
 the calculation shows that  we obtain a different pair of groups.
Thus we can rephrase the conjecture by B.~Gross and D.~Prasad
 as follows:

\begin{conjecture}
[\bf Gross--Prasad conjecture I]
\label{conj:GPA}
Suppose that $\overline{I}_\delta(\overline{V},\lambda)\boxtimes \overline{J}_\varepsilon(\overline{W},\nu)$ are tempered principal series representations
 of $SO(n+1,1) \times SO(n,1)$.  
Then    
\[ 
    {\operatorname{Hom}}_{SO(n,1)}(\overline{I}_\delta(\overline{V},\lambda)
\boxtimes \overline{J}_\varepsilon(\overline{W},\nu),{\mathbb{C}})=\bC
\]
if and only if $\overline{V } \in \widehat{SO(n)}$ and $\overline{W} \in \widehat{SO(n-1)}$
 satisfies 
\[ 
   [\overline{V}|_{SO(n-1)}:\overline{W}] \ne 0.
\]
\end{conjecture}

\begin{theorem}
[see Theorem \ref{thm:tempSO}]
\label{theorem:GPtemp}
The Gross--Prasad conjecture I holds.  
\end{theorem}
We can deduce Theorem \ref{theorem:GPtemp} from the corresponding results
 (Theorem \ref{thm:tempVW}) for the orthogonal groups
 $O(n+1,1) \times O(n,1)$
 by using results
 about the reduction from $O(N,1)$
 to the special orthogonal group $SO(N,1)$.  
See the proof of Theorem \ref{thm:tempSO}
 in Section \ref{subsec:SOtemp}
 of Appendix II
 for details.

\subsection{The Gross--Prasad conjecture II: 
 Tempered representations with trivial infinitesimal character
$\rho$}
\label{subsec:GPII}

For completeness,
 we include the discussion of the Gross--Prasad conjectures
 for tempered representations with trivial infinitesimal character $\rho$
 which we also discussed in detail in \cite{sbonGP}.  
\medskip

We modify here the notation from \cite{sbonGP} by denoting the restriction of a representation $\Pi$ of $O(n+1,1)$ to the subgroup $SO(n+1,1)$ by $\overline{\Pi}$.

The Gross--Prasad conjecture I in the previous section 
treated the case
 where both $\overline \Pi$ and $\overline \pi$ are tempered 
 {\it{principal series representations}}
 of the group $G=SO(n+1,1)$
 and $G'=SO(n,1)$, 
 respectively.

Thus the remaining cases are
 when $\overline \Pi$ or $\overline \pi$ are 
 {\it{discrete series representations}}.  
We note 
 that both $\overline \Pi$ and $\overline \pi$ cannot be discrete series representations
 in our setting
 because $\overline G$ admit discrete series representations
 if and only if $n$ is odd and $\overline{G'}$ admit those 
 if and only if $n$ is even.  
Thus we discuss the Gross--Prasad conjecture
 in this case
 separately depending on the parity of $n$, 
 with the following notation.

Consider symmetry breaking operators
 for tempered representations
 with trivial infinitesimal character $\rho$
 of the group $SO(n+1,1)$
 for $n=2m$, $2m-1$, 
 and $2m-2$.  
We denote the corresponding representations
 by $\Pi$, $\pi$, and $\varpi$, 
 respectively,
 using the subscripts
 defined in Section \ref{subsec:PiSO} in Appendix II.  
We thus consider symmetry breaking from $SO(2m+1,1)$ to $SO(2m,1)$
 and further to $SO(2m-1,1)$:

\begin{equation*} 
\begin{tabular}{ccccc}
{$\overline{\Pi}_{m,(-1)^{m+1}}$}& $ \rightarrow   $   & $\overline{ \pi}_m     $      &    $   \rightarrow      $         &   {$\overline{ \varpi }_{m-1,(-1)^m}$}. \\
 &   &   &  &
\end{tabular}
\end{equation*}

Here $\overline{\Pi}_{m,(-1)^{m+1}}$ and $\overline{\varpi}_{m-1,(-1)^m}$ are tempered principal series representations
 which are nontrivial on the center of $SO(2m+1,1)$,  respectively
 $SO(2m-1,1)$, 
 and thus are genuine representations
 of the special orthogonal groups,
 see Proposition \ref{prop:161648} (6).  
Since $\overline{\pi}_{m,+} \simeq \overline{\pi}_{m,-}$ 
 as $SO(2m,1)$-modules,
 we simply write $\overline{\pi}_{m}$
 for $\overline{\pi}_{m,\pm}$, 
 which is a discrete series representation of $SO(2m,1)$. All representations have the trivial infinitesimal  character $\rho$.

\subsubsection{The Gross--Prasad conjecture II: 
Symmetry breaking from 
 $\overline{\Pi}_{m,(-1)^{m+1}}$
 to the discrete series representation $\overline \pi_m$}

We consider first the Vogan packet of tempered representations which contains the pair $(SO(2m+1,1) \times SO(2m,1), \overline{\Pi }_{m,\delta} \boxtimes \overline{\pi}_m )$ 
 or the Vogan packet which contains the pair
$(SO(1,1+2m)\times SO(1,2m), \overline{ \Pi}_{m,\delta} \boxtimes \overline{ \pi}_m )$. 
The representations in these packets are parametrized by characters of
\[
     {\mathcal{A}}_1 \times {\mathcal{A}}_2 
     \simeq 
     (\bZ/2\bZ)^{\color{black}{m}} \times (\bZ/2\bZ)^m \simeq (\bZ/2\bZ)^{2m} . 
\] 
We recall the algorithm proposed  by B.~Gross and D.~Prasad 
 which determines a pair 
$
     (\chi_1, \chi_2)
 \in 
     \widehat{{\mathcal{A}}_1}\times \widehat{{\mathcal{A}}_2}
$, 
 hence representations 
\[
   (\overline{\Pi}(\chi_1),\overline{\pi}(\chi_2 )) \in VP(\overline{\Pi}_{m,\delta}) \times VP(\overline{\pi}_m)
\]
 so that 
\[ 
   \mbox{Hom}_{\overline{G}(\chi_2)}(\overline{\Pi}(\chi_1)|_{\overline G(\chi_2)},\overline{\pi}(\chi_2) ) \not = \{0\},
\]
where $\overline{G}(\chi_2) $ is the pure inner form determined by $\chi_2$.

\medskip
Let $T_\bC$ be a torus in $SO(2m+2,\bC)  \times SO(2m+1,\bC)$,
 and $X^{\ast}(T_\bC)$ the character group. 
As before the standard root basis $\Delta_0$ is given by
\[
e_1-e_2, e_2-e_3,\dots,{{e_{m}-e_{m+1},e_{m}+e_{m+1}}}, f_1-f_2 ,f_2-f_3, \dots, f_{m-1}-f_m,f_m
\]
{if $m \ge 1$}.

We fix $\delta=(-1)^{m+1}$
 so that $\overline{\Pi}_{m,\delta}$ is a genuine representation
 of $SO(2m+1,1)$.  
We can identify the Langlands parameter of the Vogan packet containing 
\[
(SO(2m+1,1)\times SO(2m,1),\overline{ \Pi}_{m, {{\delta}}} \boxtimes \overline{ \pi}_m)
\]
with 
 \[ me_1+(m-1)e_2+ \dots +e_m+0e_{m+1} +(m-\frac 1 2 )f_1+(m- \frac 3 2 ) f_2 +\dots + \frac 1 2 f_m.\]

Let $\delta_i$ be the character in $ \widehat{{\mathcal{A}}_1}$ which is $-1$ in the $i$-th factor
 of ${\mathcal{A}}_1$
 and equal to $1$ everywhere else, 
 and $\varepsilon_j$ be the character  which is $-1$
 in the $j$-th factor of ${\mathcal{A}}_2$ and $1$ everywhere else.

Then the algorithm by B.~Gross and D.~Prasad 
\cite[p.~993]{GP} determines characters
 $\chi_1 \in \widehat {{\mathcal{A}}_1}$
 and $\chi_2 \in \widehat {{\mathcal{A}}_2}$
 by 
\[
 \chi_1(\delta_i)= (-1)^{\#m-i+1>}  \ \ \ \ 
 \mbox{and}  \ \ \ \  
 \chi_2(\varepsilon_j)=  (-1)^{\#m-j + \frac 1 2<}, 
\]
where $\#m-i+1>$ is the cardinality
 of the set 
\[
\{j: m-i+1> \text{the coefficients of $f_j$}\},
\]
 and $\# m-j+\frac 1 2 <$ is the cardinality
 of the set 
\[
\{i:m-j+\frac 1 2 < \text{the coefficients of $e_i$}\}.  
\]

As discussed before we normalize the quasisplit form by 
\begin{align*}
 SO(m+1,m+1) \times SO(m,m+1) \quad &\mbox{if $m$ is even,} 
\\
 SO(m+2,m) \times SO(m+1, m)  \quad &\mbox{if $m$ is odd.} 
\end{align*}
Applying the formul{\ae} in \cite[(12.21)]{GP}
 we define the integers $p$ and $q$
 with $0 \leq p \leq m $ and $0 \leq q \leq m$
by
\[
 p  =  \# \{ i : \chi_1(    \delta_i  )  = (-1)^i  \}  \quad \mbox{ and }  
\quad  q  =  \# \{ j:  \chi_2(   \varepsilon_j     )  = (-1)^{m+j}  \}
\]
and we get the pure forms
\begin{align}
SO(2m-2p+1, 2p+1) \times SO(2q ,2m-2q+1) \quad &\mbox{if $m$ is even,}
\\
SO(2p+1,2m-2p+1) \times SO(2m-2q, 2q+1) \quad &\mbox{if $m$ is odd.}
\end{align}
In our setting,  we get the pair of integers $(p,q)=(0,m)$ for $m$ even; $(p,q)=(m,0)$ for $m$ odd.
Applying \cite[(12.22)]{GP}
 with correction by changing $n$ by $m$ 
 {\it{loc.cit.}}, 
 we deduce that this character defines the pure inner form 
\[
 \text{$SO(2m+1,1) \times SO(2m,1)$ for $m$ even and odd.}
\]  

The only representation in $VP(\overline{\Pi}_{m,\delta} )\times VP(\overline{\pi}_m)$ for this pair of pure inner forms  is $\overline{\Pi}_{m,\delta} \boxtimes \overline{\pi}_m$.
Hence Theorem \ref{thm:170336} implies the Gross--Prasad conjecture
 in that case.

\subsubsection{The Gross--Prasad conjecture II: 
Symmetry breaking from the discrete series representation 
 $\pi_m$
 to $\varpi_{m-1,(-1)^m}$}

We now consider the Vogan packet of tempered representations containing the pair $(SO(2m,1) \times SO(2m-1,1), {\overline{\pi}}_{m} \boxtimes \overline{\varpi}_{m-1,(-1)^m} )$, 
 {\it{i.e.}}, 
 the Vogan packet 
\[
VP(\overline{\pi}_m\boxtimes \overline{\varpi}_{m-1,(-1)^m})\subset VP(\overline{\pi}_m) \times VP(\overline{\varpi}_{m-1,(-1)^m}).  
\] 
The packet $VP(\overline{\pi}_m) \times VP(\overline{\varpi }_{m,(-1)^m})$ is
 parametrized by  characters of the finite group  
\[
     {\mathcal{A}}_2 \times {\mathcal{A}}_3
 \simeq 
 (\bZ/2\bZ)^m \times (\bZ/2\bZ)^{m-1}\simeq (\bZ/2\bZ)^{2m-1}.
\]
Again 
 the algorithm by B.~Gross and D.~Prasad  determines a pair 
 $(\chi_2,\chi_3) \in \widehat{{\mathcal{A}}_2} \times \widehat{{\mathcal{A}}_3}$
 and hence representations 
\[
 (\overline{\pi}(\chi_2), \overline{\varpi}(\chi_3 )) \in VP(\overline{\pi}_{m})\times VP(\overline{\varpi}_{m-1,(-1)^m}  )\] so that 
\[ 
     \mbox{Hom}_{\overline{G}(\chi_3)}(\overline{\pi}(\chi_2)|_{\overline{G}(\chi_3)},\overline{\varpi}(\chi_3) ) \not = \{0\}, 
\]
where $\overline{G}(\chi_3) $ is the pure inner form determined by $\chi_3$.

\medskip

Let $T_\bC$ be a torus in $SO(2m+1,\bC)\times  SO(2m, \bC)$
 and $X^{\ast}(T_{\mathbb{C}})$ the character group.  
Fix a basis 
\[
  X^{\ast}(T_{\mathbb{C}})
  =
  {\mathbb{Z}}f_1 \oplus {\mathbb{Z}}f_2 \oplus \cdots
   \oplus {\mathbb{Z}}f_m
   \oplus {\mathbb{Z}}g_1 \oplus {\mathbb{Z}} g_2 \oplus
   \cdots
   \oplus {\mathbb{Z}} g_m 
\]
 such that the standard root basis
 $\Delta_0$ is given by 
\[
  f_1 - f_2, f_2 - f_3, \cdots, f_{m-1} - f_m, f_m,
  g_1-g_2, g_2-g_3, \cdots, g_{m-1}-g_m, g_{m-1}+g_m
\] 
for $m \ge 2$.  
Take $\varepsilon=(-1)^{m}$ as before.

We identify the Langlands parameter of the Vogan packet 
\[
VP(\overline{ \pi}_{m} ) \times VP(\overline{ \varpi }_{m,(-1)^m})
\]
 with 
 \[ (m-\frac 1 2 )f_1+(m- \frac 3 2 ) f_2 +\dots + \frac 1 2 f_m +(m-1)g_1+(m-2)g_2+ \dots +g_{m-1}+0g_{m} .\]

Again applying \cite[Prop.~12.18]{GP}
 we define characters $\chi_2 \in \widehat{{\mathcal{A}}_2} $,
 $\chi_3 \in \widehat{{\mathcal{A}}_3}$ as follows: 
Let $\varepsilon_j \in {\mathcal{A}}_2 \simeq (\bZ/2\bZ)^m$ be
 the element which is $-1$ in the $j$-th factor
 and equal to $1$ everywhere else
 as in Section \ref{subsec:IV.4};
 $\gamma_k \in {\mathcal{A}}_3 \simeq (\bZ/2\bZ)^{m-1}$ the element which is $-1$ in the $k$-th factor and $1$ everywhere else.
Then $\chi_2 \in \widehat{{\mathcal{A}}_2}$ and 
 $\chi_3 \in \widehat{{\mathcal{A}}_3}$ are determined by 
\[ 
 \chi_2(\varepsilon_j)= (-1)^{\#m- j +1/2<}
  \ \ \ \ \mbox{and} \ \ \ \ 
 \chi_3(\gamma_k)=  (-1)^{\#m-k >}, 
\]
where $\#m-j+\frac 1 2 <$ is the cardinality
 of the set
  \[\{k:m-j+\frac 1 2 < \text{the coefficients of $g_k$}\}, 
\]
 and $\# m-k >$ is the cardinality of the set
 \[\{j:m-k > \text{the coefficients of $f_j$}\}.\]

As discussed  we normalize the quasisplit form by 
\begin{alignat*}{2}
 &SO(m+1, m) \times SO(m+1,m-1) \quad &&\text{if $m$ is even,} 
\\
 &SO(m,m+1) \times SO(m , m)    \quad &&\text{if $m$ is odd.} 
\end{alignat*}

 We define the integers $p$ and $q$
 with $0 \leq p \leq m$ and $0 \leq q \leq m-1$
by
\[
 p  =  \# \{ j : \chi_2(    \varepsilon_j  )  = (-1)^j  \}  
\quad \mbox{ and }  
\quad  q  =  \# \{ k :  \chi_3(   \gamma _k      )  = (-1)^{m+k}  \}, 
\]
and we get
\begin{alignat*}{2}
&SO(2m-2p+1, 2p) \times SO(2q+1 ,2m-2q-1) 
\quad 
&&\mbox{if $m$ is even,}
\\
&SO(2p+1,2m-2p) \times SO(2m-2q-1, 2q+1) \quad 
&&\mbox{if $m$ is odd.}
\end{alignat*}

In our setting, the pair of integers 
 $(p,q)$ is given by $(p,q)=(m,0)$ for $m$ even;
 $(p,q)=(0,m-1)$ for $m$ odd.
We  deduce that this character defines the pure inner form 
\[
\text{
$SO(1,2m) \times SO(1,2m-1)$ for $m$ even and odd.  
}
\]

The only representation in $VP(\overline{\pi}_{m}) \times VP(\overline{\varpi}_{m-1,(-1)^m})$ with this pair of pure inner forms 
 is $(\overline{\pi}_{m}, \overline{\varpi}_{m-1,(-1)^m})$.

\medskip

In Chapter \ref{sec:SBOrho}, 
 we have determined
\[
   {\operatorname{Hom}}_{G'}(\Pi \boxtimes \pi, {\mathbb{C}})
\quad
\text{for all $\Pi \in {\operatorname{Irr}}(G)_{\rho}$
 and $\pi \in {\operatorname{Irr}}(G')_{\rho}$}, 
\]
see Theorems \ref{thm:SBOvanish} and \ref{thm:SBOone}
 and also Theorem \ref{thm:SBOBF}
 for orthogonal groups
\[
   G \times G'=O(n+1,1) \times O(n,1), 
\]
{}from which we deduce analogous results about
\[
   {\operatorname{Hom}}_{\overline{G'}}(\overline \Pi \boxtimes \overline \pi, {\mathbb{C}})
\quad
\text{for all $\overline\Pi \in {\operatorname{Irr}}(\overline G)_{\rho}$
 and $\overline \pi \in {\operatorname{Irr}}(\overline{G'})_{\rho}$},
\]
for the special orthogonal groups
\[
   \overline G \times \overline{G'} = SO(n+1,1) \times SO(n,1), 
\]
in Theorem \ref{thm:170336}.  
By the aforementioned argument,
 Theorem \ref{thm:170336} implies the following.  

\begin{theorem}
\label{thm:GPdisctemp}
The conjectures by B.~Gross and D.~Prasad \cite{GP}
 for tempered representations
 of special orthogonal groups $SO(n+1,1) \times SO(n,1)$ with trivial infinitesimal character $\rho$ hold.  
\end{theorem}

\begin{remark}
The Gross--Prasad conjectures concern tempered representations
 with trivial infinitesimal character $\rho$, 
 but one may expect similar results for unitary representations
 of orthogonal groups with integral infinitesimal character. 
Considering \lq\lq{Arthur--Vogan packets}\rq\rq\
 instead of the Vogan packets will include 
 other unitary representations
 which are of interest to number theory
 for example to the representation 
\index{A}{Aqlmd@$A_{\mathfrak{q}}(\lambda)$}
$A_{\mathfrak{q}}(\lambda)$.  
Low dimensional examples  and our results suggest
 that there exists  pairs of groups
 $\overline G\times \overline{G'}=SO(p+1,q) \times SO(p,q)$
 and of representations $ U_1 \boxtimes U_2$
 in this \lq\lq{Arthur--Vogan packet}\rq\rq\ 
 so that 
$
   \mbox{Hom}_{\overline{G'}}(U_1|_{\overline{G'}} \boxtimes U_2, \mathbb{C}) \not = \{0\}. 
$
The examples also suggest an algorithm
 to determine  pairs of groups
 and the pairs of representations with nontrivial multiplicity.
\end{remark}

\newpage
\section{Application II:
 Periods, distinguished representations and $(\mathfrak{g},K)$-cohomologies}
\label{sec:period}
\index{B}{gKcohomology@$(\mathfrak{g},K)$-cohomology}

Let $H$ be a subgroup of $G$.  
Following the terminology used in automorphic forms and the relative trace formula,
 we say that a smooth representation $U$ of $G$ is $H$-distinguished if
there exists 
 a nontrivial $H$-invariant linear functional 
\[
   F^{H}:U \rightarrow \bC, 
\]
 {\it{i.e.}}, 
 if $U$ has a nontrivial 
\index{B}{period@period}
$H$-period $F^{H}$. 
We consider first irreducible representations of $G$
 with infinitesimal character $\rho$ which are $H$-distinguished
 for the  pair $(G,H)=(O(n+1,1),O(m+1,1))$ 
 or for the pair $(G,H)= (O(n,1) \times O(m,1), O(m,1))$  with $m \le n$.  
We then discuss  a bilinear form on the $(\mathfrak{g},K)$-cohomology
 of the representations of $(O(n+1,1) \times O(n,1))$
 with infinitesimal character $\rho$
 which is induced by a symmetry breaking operator.

\subsection{Periods and $O(n,1)$-distinguished representations}
\subsubsection{Periods}
Let $\mathbb K$ be a number field, $\mathbb A$ its adels and 
 let $G_1\times G_2$ be a direct product
 of semisimple groups over a number field $\mathbb K$. 
We assume that $G_2 \subset G_1$. 
If the outer tensor product representation
 $\Pi_{\mathbb A} \boxtimes \pi_{\mathbb A}$ is an automorphic representation
 of the direct product group $G_1({\mathbb A}) \times G_2({\mathbb A})$, 
 then the $G_2$-period integral is defined as
\[
  \int_{G_2({\mathbb K} )\backslash G_2({\mathbb A} )}\Phi_1 (h) \phi_2(h) dh.
\]
Here $\Phi_1$ and $\phi_2$ are smooth vectors
 for the representation $\Pi_{\mathbb A} \boxtimes \pi_{\mathbb A}$. 
If $\Pi_{\mathbb A} \boxtimes \pi_{\mathbb A}$ is cuspidal, 
 then the integral converges and 
it defines a $G_2({\mathbb A})$-invariant linear functional
 on the smooth vectors of $\Pi_{\mathbb A} \boxtimes \pi_{\mathbb A}$.  
If this linear functional is not zero,
 then $\Pi_{\mathbb A}\boxtimes \pi_{\mathbb A}$
 is called $G_2$-{\it{distinguished}}. 
Conjecturally for certain pairs of groups
 the value of this integral is a multiple
 of the central value of an $L$-function, 
 see \cite{A1, II, IY}.

Often this period integral  factors into a product of local integrals. Following the global terminology
 we say that an admissible smooth representation $\Pi \boxtimes \pi$
 of the direct product group $G_1(\bR)\times G_2({\mathbb R})$ is $G_2({\mathbb R})$-{\it{distinguished}}
 if there is a nontrivial continuous linear functional 
\[F^{G_2({\mathbb R})}:\Pi \boxtimes \pi \rightarrow  {\mathbb C} \]
 which is invariant by $G_2({\mathbb{R}})$
 under the diagonal action.  
Here we recall Section \ref{subsec:BFSBO}
 for the topology on the tensor product.  
If $\Pi \boxtimes \pi$ is $G_2({\mathbb R})$-distinguished,
 we say that $F^{G_2({\mathbb R})}$ is a 
\index{B}{period@period|textbf}
{\it{period}} of $\Pi \boxtimes \pi$.  
We say that the period is nontrivial
 on a vector
 $\Phi \otimes \phi \in \Pi \boxtimes \pi$
 if $\Phi \otimes \phi$ is not in the kernel of $F^{G_2({\mathbb R})}$. 
If the period is nontrivial on a unit function $\Phi \otimes \phi$,
 we refer to its image as the value of the period on $\Phi \otimes \phi$.

\begin{remark} 
The integral 
\[  \int_{G_2(\mathbb R )}
\Phi (h) \phi (h) dh
\]
converges for some  smooth vectors of 
\index{B}{discreteseries@discrete series representation}
discrete series representations
 $\Pi \boxtimes \pi$
 for some symmetric pairs $(G_1(\mathbb R),G_2(\mathbb R))$.  
This was used by J.~Vargas \cite{Vargas}
 to determine  some subrepresentations
 in the restriction of some discrete series representations $\Pi$ of $G_1(\mathbb R)$
 to the subgroup $G_2(\mathbb R)$.  
\end{remark}

\medskip
We recall from Theorem \ref{thm:SBOBF}
 that the space of symmetry breaking operators 
\[ 
     \mbox{Hom}_{G_2(\bR)}(\Pi|_{G_2(\bR)}, \pi^{\vee}) 
\]
 and the space of $G_2(\bR)$-invariant continuous linear functionals
\[\mbox{Hom}_{G_2(\bR)}(\Pi \boxtimes \pi,\bC)\]
 are naturally isomorphic to each other. 
Thus, 
 instead of considering a $G_2(\bR)$-equivariant continuous linear functional
 defined by an integral, 
 we may use symmetry breaking operators
 to construct $G_2(\bR)$-invariant continuous  linear functionals. 
This technique allows us to obtain
$G_2(\bR)$-invariant continuous linear functionals
 not only for discrete series representations
 but also for nontempered representations.  
Thus we can determine 
 for the pair $(G,G')= (O(n+1,1), O(n,1))$
 the dimension of the space
 ${\operatorname{Hom}}_{G'}(\Pi \boxtimes \pi, {\mathbb{C}})$
for all $\Pi \in {\operatorname{Irr}}(G)_{\rho}$
 and $\pi \in {\operatorname{Irr}}(G')_{\rho}$
 as follows.

\begin{corollary}
\label{cor:bilin}
Suppose $0 \le i \le n+1$, 
 $0 \le j \le n$, 
 and $\delta$, $\varepsilon \in \{\pm\}$.  
Let $\Pi_{i,\delta}$ and $\pi_{j,\varepsilon}$ be
 irreducible admissible smooth representations of $G=O(n+1,1)$ and $G'=O(n,1)$,  respectively, 
 that have the trivial infinitesimal character $\rho$ as in \eqref{eqn:Pild}. 
Then the following three conditions 
 on $(i,j,\delta,\varepsilon)$ are equivalent:
\begin{enumerate}
\item[{\rm{(i)}}]
${\operatorname{Hom}}_{G'}( \Pi_{i,\delta} \boxtimes \pi_{j,\varepsilon}, \bC) \not = \{0\};$
\item[{\rm{(ii)}}]$\dim_{\mathbb{C}}{\operatorname{Hom}}_{G'}(\Pi_{i,\delta} \boxtimes \pi_{j,\varepsilon}, \bC) = 1;$
\item[{\rm{(iii)}}]$j \in \{i,i-1 \}$ and $\delta=\varepsilon$.  
\end{enumerate}
\end{corollary}
\begin{proof}
Owing to Theorem \ref{thm:SBOBF}, 
 this is a restatement
 of Theorems \ref{thm:SBOvanish} and \ref{thm:SBOone}.  
\end{proof}
\subsubsection{Distinguished representations}
Let $G$ be a reductive group, 
 and $H$ a reductive subgroup.  
We regard $H$ as a subgroup of the direct product group $G \times H$
 via the diagonal embedding $H \hookrightarrow G \times H$.  

\begin{definition}
Let $\psi$ be a one-dimensional representation of $H$. 
We say an admissible smooth representation  $\Pi$ of $G$ is
\index{B}{distinguishedHpsi@distinguished, $(H,\psi)$-|textbf}
$(H,\psi)$-{\it{distinguished}}
 if 
\[ \mbox{Hom}_{H}(\Pi \boxtimes \psi^{\vee}, {\mathbb{C}})
 \simeq \mbox{Hom}_{H}(\Pi|_H, \psi) \not = \{0\}.  
\]
If the character $\psi $ is trivial,
 we say $\Pi$ is 
\index{B}{distinguishedH@distinguished, $H$-|textbf}
$H$-{\it{distinguished}}.  
\end{definition}

In what follows, 
 we deal mainly with the pair
\[
   (G,H)=(O(n+1,1), O(m+1,1))
\quad
\text{for  $m \le n$.  }
\]

\begin{theorem}
\label{thm:period1}
Let $0 \le i \le n+1$.  
Then the representations $\Pi_{i,\delta}$ $(\delta \in \{\pm\})$
 of $G=O(n+1,1)$ are  $O(n+1-i,1)$-distinguished. 
\end{theorem}

The period is given by the composition 
 of the symmetry breaking operators
 that we constructed in Chapter \ref{sec:pfSBrho}
 with respect to the chain
 of subgroups
\begin{equation}
\label{eqn:seqOn}
  G=O(n+1,1) \supset O(n,1) \supset O(n-1,1) \supset \cdots \supset O(m+1,1)=H,
\end{equation}
 as we shall see in the proof
 in Section \ref{subsec:pfperiod}.  
Without loss of generality,
 we consider the case $\delta=+$, 
 and write simply $\Pi_i$ for $\Pi_{i,+}$.  
We recall from Theorem \ref{thm:LNM20} (3)
 that $\Pi_i \equiv \Pi_{i,+}$ has a 
\index{B}{minimalKtype@minimal $K$-type}
minimal $K$-type
 $\mub(i,+)=\Exterior^i({\mathbb{C}}^{n+1}) \boxtimes {\bf{1}}$.

Let $v\in \Exterior^i({\mathbb{C}}^{n+1})$ be the image of $1 \in {\mathbb{C}}$
 via the following successive inclusions:
\[
   \Exterior^{i}({\mathbb{C}}^{n+1}) \supset 
   \Exterior^{i-1}({\mathbb{C}}^{n}) \supset 
   \cdots \supset
   \Exterior^{i-l}({\mathbb{C}}^{n+1-l}) \supset 
   \cdots \supset 
  \Exterior^{0}({\mathbb{C}}^{n+1-i}) \simeq {\mathbb{C}} \ni 1, 
\]
 and we regard $v$ as an element
 of the minimal $K$-type $\mub(i,+)$
 of $\Pi_i$.  

\begin{theorem}
\label{thm:period2}
Let $\Pi_i$ be the irreducible representation of $G=O(n+1,1)$, 
 and $v$ be the normalized element
 of its minimal $K$-type as above.  
For $0 \le i \le n$, 
 the value $F(v)$ of the $O(n+1-i,1)$-period $F$ 
 on $v \in \Pi_i$ is
\[
  \frac{\pi^{\frac 1 4 i(2n-i-1)}}{((n-i)!)^{i-1}}
  \times
\begin{cases}  
  \frac{1}{(n-2i)!}\quad &\text{if $2i < n+1$}, 
\\
  (-1)^{n+1} (2i-n-1)!\quad & \text{if $2i \ge n+1$}.  
\end{cases}
\]
\end{theorem} 

\subsubsection{Symmetry breaking operators from $\Pi_{i,\delta}$
 to $\pi_{j,\delta}$
($j \in \{i-1,i\}$)}
\label{subsec:KKPipj}
Let $(G,G')=(O(n+1,1),O(n,1))$.  
We recall from Theorem \ref{thm:LNM20} (2)
 that 
\begin{align*}
{\operatorname{Irr}}(G)_{\rho} = 
&\{\Pi_{i,\delta} : 0 \le i \le n+1, \delta=\pm\}, 
\\
{\operatorname{Irr}}(G')_{\rho} = 
&\{\pi_{j,\varepsilon} : 0 \le j \le n, \varepsilon=\pm\}.  
\end{align*}
In Chapter \ref{sec:pfSBrho},
 we constructed nontrivial symmetry breaking operators
\[
  A_{i,j} \colon \Pi_{i,\delta} \to \pi_{j,\varepsilon}
\]
for $j \in \{i-1,i\}$
 and $\delta = \varepsilon$, 
 and investigated their $(K,K')$-spectrum
 for minimal $K$- and $K'$-types, 
\[
  (\mu,\mu')=(\mub(i,\delta),\mub(j,\delta)'), 
\]
see Proposition \ref{prop:Aidown} in the case $j=i-1$
 and Proposition \ref{prop:AiiAq} in the case $j=i$.

For the proofs of Theorems \ref{thm:period1} and \ref{thm:period2},
 we use these operators $A_{i,j}$ in the case $j=i-1$.  
For the study of the bilinear forms
 on $({\mathfrak{g}},K)$-cohomologies
 (see Section \ref{subsec:cohgkex} below), 
 we shall use them in the case $j=i$.

\subsection{Proofs of Theorems \ref{thm:period1} and \ref{thm:period2}}
\label{subsec:pfperiod}
We are ready to prove Theorems \ref{thm:period1} and \ref{thm:period2}
  by using Proposition \ref{prop:Aidown} successively.  
\begin{proof}
[Proof of Theorem \ref{thm:period1}]
Consider the chain \eqref{eqn:seqOn}
 of orthogonal subgroups 
 with $m=n-i$.  
For $1 \le \ell \le i$, 
 we denote by
\[
   A_{\ell, \ell-1}
   \colon
   \Pi_{i-\ell+1}^{O(n-\ell+2,1)} \to \Pi_{i-\ell}^{O(n-\ell+1,1)}
\]
 the symmetry breaking operator given in Proposition \ref{prop:Aidown}
 for the pair 
 $(O(n-\ell+2,1), O(n-\ell+1,1))$ of groups.  
Here \lq\lq{$\Pi_{i-\ell}^{O(n-\ell+1,1)}$}\rq\rq\ stands
 for the irreducible representation \lq\lq{$\Pi_{i-\ell,+}$}\rq\rq\
 of the group $O(n-\ell+1,1)$
 as given in Theorem \ref{thm:LNM20}, 
 by a little abuse of notation.  
Then the composition
\begin{equation}
\label{eqn:compfunct}
  F:= A_{1,0} \circ \cdots \circ A_{i-1,i-2} \circ A_{i,i-1}
\end{equation}
defines a nonzero $O(n+1-i,1)$-invariant functional
 on the irreducible representation $\Pi_i \equiv \Pi_{i,+}$
 of $G=O(n+1,1)$.  
\end{proof}

\begin{proof}
[Proof of Theorem \ref{thm:period2}]
The irreducible representation $\Pi_{i-\ell}^{O(n+1-\ell,1)}$, 
 namely, 
\lq\lq{$\Pi_{i-\ell,+}$}\rq\rq\
 of the group $O(n+1-\ell, 1)$
 has a minimal $K$-type
\[
   \mub(i-\ell,+)^{(\ell)}
   :=
   \Exterior^{i-\ell}({\mathbb{C}}^{n+1-\ell}) \boxtimes {\bf{1}}
   \in \widehat{O(n+1-\ell)}\times \widehat{O(1)}. 
\]
The $(K,K')$-spectrum of the symmetry breaking operator
 $A_{\ell, \ell-1}
   \colon
   \Pi_{i-\ell+1}^{O(n-\ell+2,1)} \to \Pi_{i-\ell}^{O(n-\ell+1,1)}
$
 for the minimal $K$-types
 $\mub (i-\ell+1,+)^{(\ell-1)} \hookleftarrow \mub (i-\ell,+)^{(\ell)}$
 is given by 
\[
\frac{\pi^{\frac {n-\ell}{2}}}{(n-i)!} \times 
\begin{cases}  
 n-2i+ \ell \qquad &\text{if $n \ne 2i -\ell$}, 
\\
  1 
  &\text{if $n = 2i -\ell$}, 
\end{cases}
\]
by Proposition \ref{prop:Aidown}.  
Applying this formula successively 
 to the sequence of minimal $K$-types:
\[
  \mub(i,+) \equiv \mub(i,+)^{(0)} \hookleftarrow \cdots \hookleftarrow
  \mub(i-\ell,+)^{(\ell)} \hookleftarrow \cdots 
  \hookleftarrow \mub(0,+)^{(i)}={\mathbb{C}}, 
\]
we get
\[
  F(v)=\prod_{\ell=1}^i \frac{\pi^{\frac{n-\ell}{2}} (n-2i+\ell)}{(n-i)!}
  =\frac{\pi^{\frac{1}{4}i(2n-i-1)}}{((n-i)!)^{i-1}(n-2i)!}
\]
 if $n > 2i-1$.

On the other hand, 
 if $n < 2i-1< 2n-1$, 
 then 
\begin{align*}
  F(v)
  =\,&
  \left(\prod_{\ell=1}^{2i-n-1} \frac{\pi^{\frac{n-\ell}{2}} (n-2i+\ell)}{(n-i)!}\right)
  \cdot
  \frac{\pi^{n-i}}{(n-i)!}
  \cdot
 \left(\prod_{\ell=2i-n+1}^{i}\frac{\pi^{\frac{n-\ell}{2}} (n-2i+\ell)}{(n-i)!}\right)
\\
=\,&
\frac{\pi^{\frac 1 4 i(2n-i-1)}}{((n-i)!)^i}
(((-1)^{2i-n-1}(2i-n-1)!)\cdot 1 \cdot ((n-i)!))
\\
=\,&
\frac{(-1)^{n+1} \pi^{\frac 1 4 i(2n-i-1)}(2i-n-1)!}{((n-i)!)^{i-1}}.  
\end{align*}
The cases $i= \frac{n+1}{2}$ ($n$: odd) or $i=n$ are treated separately, 
 and it turns out that the formula of $F(v)$ coincides
 with the one for $i < 2i-1< 2n-1$.  
Thus we have completed the proof of Theorem \ref{thm:period2}.  
\end{proof}

In the next theorem,
 we consider the pair
\[
 \text{$(G,H)=(O(n+1,1),O(m+1,1))$
 with $m \le n$}.  
\]
We write $\Pi_i^G$ ($0 \le i \le n+1$)
 for the irreducible representation $\Pi_{i,+}$
 of $G$
 (see \eqref{eqn:Pild}),
 and write $\pi_j^H$ for the irreducible representation 
 \lq\lq{$\Pi_{j,+}$}\rq\rq\
 of the subgroup $H$ 
 for $0 \le j \le m+1$.  
Theorem \ref{thm:171517} below generalizes Theorem \ref{thm:period1}, 
 which corresponds to the case $j=0$.  

\begin{theorem}
\label{thm:171517}
Let $0 \le i \le n+1$ and $0 \le j \le m+1$.  
\begin{enumerate}
\item[{\rm{(1)}}]
The outer tensor product representation $\Pi_i^G \boxtimes \pi_j^H$ 
 of the direct product group $G \times H$ has an $H$-period if $0 \le i -j \le n-m$.  
\item[{\rm{(2)}}]
The period constructed 
 by the composition of the symmetry breaking operators 
 via the sequence \eqref{eqn:seqOn} 
 is nontrivial on the minimal $K$-type.  
\end{enumerate}
\end{theorem}

\begin{proof}
[Proof of Theorem \ref{thm:171517}]
The proof is essentially the same  
 with the one for Theorem \ref{thm:period1}
 except
 that we use not only the surjective symmetry breaking operator
 $A_{i,i-1} \colon \Pi_{i,+} \to \pi_{i-1,+}$
 for the pair $(G,G')=(O(n+1,1),O(n,1))$
but also the one 
\[
  A_{i,i} \colon \Pi_{i,+} \to \pi_{i,+}
\]
for which the $(K,K')$-spectrum
 on minimal $K$-types $\mub(i,+)\hookleftarrow \mub(i,+)'$
 is nonzero by Proposition \ref{prop:AiiAq}.

Composing the symmetry breaking operators $A_{k,k-1}$ or $A_{k,k}$ successively
 to the sequence \eqref{eqn:seqOn} of orthogonal groups, 
 we get a nonzero continuous $H$-homomorphism
 $\Pi_i^G \to \pi_j^H$ 
 if $0 \le i-j \le n-m$.  
Then the first statement follows 
 because $\pi_j^H$ is self-dual.  
The second statement is clear by the construction
 and by the $(K,K')$-spectrum.   
\end{proof}

\subsection{Bilinear forms
 on $({\mathfrak{g}},K)$-cohomologies via symmetry breaking: General theory for nonvanishing}
For the rest of this chapter,
 we discuss $({\mathfrak{g}},K)$-cohomologies
 via symmetry breaking.  
In this section,
 we deal with a general setting
 where $G \supset G'$ is a pair of real reductive Lie groups.  
We shall define natural bilinear forms
 on $({\mathfrak{g}},K)$-cohomologies
 and $({\mathfrak{g}}',K')$-cohomologies
 via symmetry breaking $G \downarrow G'$, 
 and prove a nonvanishing result
 (Theorem \ref{thm:171556}) in the general setting generalizing a theorem of B.~Sun \cite{ S}.  
\subsubsection{Pull-back of $({\mathfrak{g}},K)$-cohomologies
 via symmetry breaking}
Let $G$ be a real reductive Lie group,
 and $K$ a maximal compact subgroup.  
We recall that the $({\mathfrak{g}},K)$-cohomology groups
 are the right derived functor of 
\[
   {\operatorname{Hom}}_{{\mathfrak{g}},K}({\mathbb{C}},\ast)
\]
{}from the category of $({\mathfrak{g}},K)$-modules.  
Suppose further 
 that $G'$ is a real reductive subgroup 
 such that $K':=K \cap G'$ is a maximal compact subgroup
 of $G'$.  
We write 
$
   {\mathfrak{g}}_{\mathbb{C}}
  ={\mathfrak{k}}_{\mathbb{C}}+{\mathfrak{p}}_{\mathbb{C}}
$ and 
$
   {\mathfrak{g}}_{\mathbb{C}}'
  ={\mathfrak{k}}_{\mathbb{C}}'+{\mathfrak{p}}_{\mathbb{C}}'
$
 for the complexifications of the corresponding Cartan decompositions.  
In what follows, 
we set 
\[
  d:=\dim G'/K'=\dim_{\mathbb{C}}{\mathfrak{p}}_{\mathbb{C}}'.  
\]
We shall use the Poincar{\'e} duality
 for the subgroup $G'$, 
 which may be disconnected.  
In order to deal with disconnected groups, 
 we consider the natural one-dimensional representation of $K'$
 defined by 
\begin{equation}
\label{eqn:Ksign}
   \chi \colon K' \to GL_{\mathbb{C}}
   (\Exterior^d {\mathfrak{p}}_{\mathbb{C}}')
   \simeq {\mathbb{C}}^{\times}.  
\end{equation}
The differential $d \chi$ is trivial 
 on the Lie algebra ${\mathfrak{k}}'$.  
We extend $\chi$ to a $({\mathfrak{g}}', K')$-module
 by letting ${\mathfrak{g}}'$ act trivially.  
Then we have 
\begin{equation}
\label{eqn:gKtop}
   H^d({\mathfrak{g}}',K';\chi) \simeq {\mathbb{C}}.  
\end{equation}

\begin{example}
\label{ex:chiOn}
For $G'=O(n,1)$, 
 the adjoint action of $K' \simeq O(n) \times O(1)$
 on ${\mathfrak{p}}_{\mathbb{C}}' \simeq {\mathbb{C}}^n$
 gives rise to the one-dimensional representation
\[
   \Exterior^n ({\mathfrak{p}}_{\mathbb{C}}')
   \simeq
   \Exterior^n ({\mathbb{C}}^n) \boxtimes (-1)^n.  
\]
Hence, 
 in terms of the one-dimensional character
\index{A}{1chipmpm@$\chi_{\pm\pm}$, one-dimensional representation of $O(n+1,1)$}
 $\chi_{a b}$ of $O(n,1)$
 defined in \eqref{eqn:chiab}, 
 the $({\mathfrak{g}}', K')$-module $\chi$
 defined in \eqref{eqn:Ksign} is isomorphic to $\chi_{-,(-1)^n}$.  
See also Example \ref{ex:gKchiOn} below.  
\end{example}

Now we recall the Poincar{\'e} duality
 for $({\mathfrak{g}}, K)$-cohomologies of $({\mathfrak{g}}, K)$-modules
 when $G$ is not necessarily connected:
\begin{lemma}
[Poincar{\'e} duality]
\label{lem:Poincare}
Let $\chi$ be the one-dimensional $({\mathfrak{g}}', K')$-module
 as in \eqref{eqn:Ksign}.  
Then for any irreducible $({\mathfrak{g}}', K')$-module $Y$, 
 there is a canonical perfect pairing
\begin{equation}
\label{eqn:Poincare}
   H^j({\mathfrak{g}}', K'; Y) \times H^{d-j}({\mathfrak{g}}', K'; Y^{\vee} \otimes \chi)
\to H^d({\mathfrak{g}}', K'; \chi) \simeq {\mathbb{C}}
\end{equation}
for all $j \in {\mathbb{N}}$.  
\end{lemma}

\begin{proof}
See  \cite[Cor.~3.6]{KV}
 (see also \cite[Chap.~I, Sect.~1]{BW}
 when $K$ is connected).  
\end{proof}

We use the terminology 
\index{B}{symmetrybreakingoperators@symmetry breaking operators|textbf}
 \lq\lq{symmetry breaking operator}\rq\rq\
 also in the category
 of $({\mathfrak{g}},K)$-modules, 
 when  we are given a pair $({\mathfrak{g}},K)$ and $({\mathfrak{g}}',K')$
 such that ${\mathfrak{g}} \supset {\mathfrak{g}}'$ and $K \supset K'$. 
We prove the following.  
\begin{proposition}
\label{prop:171543}
Let $X$ be a $({\mathfrak{g}}, K)$-module,
 $Y$ a $({\mathfrak{g}}', K')$-module,
 and $Y^{\vee}$ the contragredient $({\mathfrak{g}}', K')$-module
 of $Y$.  
Suppose $T\colon X \to Y$ is a $({\mathfrak{g}}', K')$-homomorphism, 
 where we regard the $({\mathfrak{g}}, K)$-module $X$
 as a $({\mathfrak{g}}', K')$ by restriction.  
Then the symmetry breaking operator $T$ induces a canonical homomorphism
\begin{equation}
\label{eqn:cohmap}
T_{\ast} \colon H^j({\mathfrak{g}}, K; X) \to H^j({\mathfrak{g}}', K'; Y)
\end{equation}
and a canonical bilinear form
\begin{equation}
\label{eqn:cohbi}
   B_T \colon 
   H^j({\mathfrak{g}}, K; X) \times H^{d-j}({\mathfrak{g}}', K'; Y^{\vee} \otimes \chi)
   \to 
  {\mathbb{C}}
\end{equation}
for all $j \in {\mathbb{N}}$.
\end{proposition}

\begin{proof}
The $({\mathfrak{g}}, K)$-module $X$ is viewed 
 as a $({\mathfrak{g}}', K')$-module
 by restriction.  
Then the map of pairs 
$
   ({\mathfrak{g}}', K') \hookrightarrow ({\mathfrak{g}}, K)
$
 induces natural homomorphisms
\[
   H^j({\mathfrak{g}}, K; X) \to H^j({\mathfrak{g}}', K'; X)
\quad
\text{for all $j \in {\mathbb{N}}$.}
\]

On the other hand,
 since $T \colon X \to Y$ is a $({\mathfrak{g}}', K')$-homomorphism,
 it induces natural homomorphisms
\[
   H^j({\mathfrak{g}}', K'; X) \to H^j({\mathfrak{g}}', K'; Y)
\quad
\text{for all $j \in {\mathbb{N}}$.}
\]
Composing these two maps,
 we get the homomorphisms
 \eqref{eqn:cohmap}.

In turn, 
 combining the morphism \eqref{eqn:cohmap}
 with 
\index{B}{Poincareduality@Poincar{\'e} duality}
 the Poincar{\'e} duality in \eqref{eqn:Poincare}
 in Lemma \ref{lem:Poincare}, 
 we get the bilinear map $B_T$
 as desired.  
\end{proof}

\subsubsection{Nonvanishing of pull-back of $({\mathfrak{g}}, K)$-cohomologies
 of $A_{{\mathfrak{q}}}$ via symmetry breaking}

Retain the setting where $(G,G')$ is a pair of real reductive Lie groups.  
In this subsection, 
 we discuss a nonvanishing result 
 for morphisms between $({\mathfrak{g}}, K)$-cohomologies
 and $({\mathfrak{g}}', K')$-cohomologies
 under certain assumption 
 on the 
\index{B}{KspectrumKprime@$(K,K')$-spectrum}
$(K,K')$-spectrum
 of the symmetry breaking operator,
 see Theorem \ref{thm:171556} and Remark \ref{rem:cohKK} below.

In order to formulate a nonvanishing theorem, 
 we begin with a setup 
 for finite-dimensional representations
 of compact Lie groups.  
Let $U$ be a $K$-module,
 $U'$ a $K'$-module,
 and $\varphi \colon U \to U'$ a $K'$-homomorphism.  
Via the inclusion map ${\mathfrak{p}}' \hookrightarrow {\mathfrak{p}}$, 
 the composition
 of the following two morphisms
\[
   {\operatorname{Hom}}_K(\Exterior^j {\mathfrak{p}}_{\mathbb{C}},U)
  \to 
   {\operatorname{Hom}}_{K'}(\Exterior^j {\mathfrak{p}}_{\mathbb{C}},U')
  \to
  {\operatorname{Hom}}_{K'}(\Exterior^j {\mathfrak{p}}_{\mathbb{C}}',U')
\]
 induces natural homomorphisms
\begin{equation}
\label{eqn:Upp}
  \varphi_{\ast} \colon 
  {\operatorname{Hom}}_K(\Exterior^j {\mathfrak{p}}_{\mathbb{C}},U)
  \to 
  {\operatorname{Hom}}_{K'}(\Exterior^j {\mathfrak{p}}_{\mathbb{C}}',U')
\end{equation}
for all $j \in {\mathbb{N}}$.  
\begin{definition}
\label{def:pnonvan}
A $K'$-homomorphism $\varphi$ is said to be ${\mathfrak{p}}$-{\it{nonvanishing
 at degree}} $j$
 if the induced morphism $\varphi_{\ast}$ in \eqref{eqn:Upp} is nonzero.  
\end{definition}

By a theorem of Vogan--Zuckerman \cite{VZ}
 every irreducible representation of $G$
 with nontrivial $({\mathfrak{g}},K)$-cohomology is equivalent
 to the representation, 
 to be denoted usually by $A_{\mathfrak{q}}$ 
 for some $\theta$-stable parabolic subalgebra ${\mathfrak {q}}$.  
Here $A_{\mathfrak{q}}$ is a $({\mathfrak {g}},K)$-module cohomologically
 induced from the trivial one-dimensional representation
 of the Levi subgroup $L=N_{G}({\mathfrak{q}}):=\{g \in G : {\operatorname{Ad}}(g){\mathfrak{q}}={\mathfrak{q}}\}$.  
Suppose ${\mathfrak{q}}={\mathfrak{l}}_{\mathbb{C}}+{\mathfrak{u}}$
 and ${\mathfrak{q}}'={\mathfrak{l}}_{\mathbb{C}}'+{\mathfrak{u}}'$
 be $\theta$-stable parabolic subalgebras
 of ${\mathfrak{g}}_{\mathbb{C}}$ and ${\mathfrak{g}}_{\mathbb{C}}'$, 
 respectively.  
In general, 
 we do not assume an inclusive relation
 of ${\mathfrak{q}}$ and ${\mathfrak{q}}'$.  
We shall work
 with a symmetry breaking operator $T \colon X \to Y$, 
 where $X$ is a $({\mathfrak{g}},K)$-module
 $A_{\mathfrak{q}}$ and $Y$ is a $({\mathfrak{g}}',K')$-module
 $A_{\mathfrak{q}'}$.  
We note that $Y$ contains a unique minimal $K'$-type,
 say $\mu'$.  
Let $Y'$ be the $K'$-submodule
 containing all the remaining $K'$-types in $Y$, 
 and 
\[
  {\operatorname{pr}} \colon Y \to \mu'
\] 
 be the first projection
 of the direct sum decomposition
 $Y= \mu' \oplus Y'$.  

\begin{theorem}
\label{thm:171556}
Let $T \colon X \to Y$ be a $({\mathfrak{g}}',K')$-homomorphism, 
 where $X$ is a $({\mathfrak{g}},K)$-module $A_{\mathfrak{q}}$
 and $Y$ is a $({\mathfrak{g}}',K')$-module $A_{\mathfrak{q}'}$.  
Let $U$ be the representation space
 of the minimal $K$-type $\mu$ in $X$, 
 and $U'$ that of the minimal $K'$-type $\mu'$ in $Y$.  
We define a $K'$-homomorphism by 
\begin{equation}
\label{eqn:phiT}
   \varphi_T := {\operatorname{pr}} \circ T|_{U}
   \colon U \to U'.  
\end{equation}
\begin{enumerate}
\item[{\rm{(1)}}]
If $\varphi_T$ is zero, 
 then the homomorphisms 
$
  T_{\ast} \colon 
  H^j({\mathfrak{g}},K;X) \to H^j({\mathfrak{g}}',K';Y)
$
 (see \eqref{eqn:cohmap})
 and the bilinear form $B_T$
 (see \eqref{eqn:cohbi})
 vanish for all degrees $j \in {\mathbb{N}}$.   
\item[{\rm{(2)}}]
If $\varphi_T$ is ${\mathfrak{p}}$-nonvanishing
 at degree $j$, 
 then $T_{\ast}$ and the bilinear forms $B_T$ are nonzero
  for this degree $j$.  
\end{enumerate}
\end{theorem}

\begin{proof}
[Proof of Theorem \ref{thm:171556}]
By Vogan--Zuckerman \cite[Cor.~3.7 and Prop.~3.2]{VZ}, 
 we have natural isomorphisms:
\begin{equation}
\label{eqn:gKp}
   {\operatorname{Hom}}_{K}(\Exterior^j {\mathfrak{p}}_{\mathbb{C}},U)
  \overset \sim \to
  {\operatorname{Hom}}_{K}(\Exterior^j {\mathfrak{p}}_{\mathbb{C}},X)
  \overset \sim \to 
  H^j ({\mathfrak{g}}, K; A_{\mathfrak{q}}).  
\end{equation}
By the definition \eqref{eqn:cohmap} of $T_{\ast}$
 in Proposition \ref{prop:171543}
 and $\varphi_{\ast}$ (see \eqref{eqn:Upp}), 
 the following diagram commutes:
\begin{alignat*}{5}
   &{\operatorname{Hom}}_{K}(\Exterior^j {\mathfrak{p}}_{\mathbb{C}},U)
   &&\overset \sim \to \,\,
   &&{\operatorname{Hom}}_{K}(\Exterior^j {\mathfrak{p}}_{\mathbb{C}},X)
   &&\overset \sim \to 
   &&H^j ({\mathfrak{g}}, K; X)
\\
   &\qquad (T|_U)_{\ast} \downarrow
   &&
   && \circlearrowright
   &&
   &&\qquad \downarrow T_{\ast}
\\
   &{\operatorname{Hom}}_{K'}(\Exterior^j {\mathfrak{p}}_{\mathbb{C}}',T(U))
  &&\subset
  &&{\operatorname{Hom}}_{K'}(\Exterior^j {\mathfrak{p}}_{\mathbb{C}}',Y)
  &&\overset \sim \to 
  &&H^j ({\mathfrak{g}}', K'; Y).  
\end{alignat*}

Since ${\operatorname{Hom}}_{K'}(\Exterior^j {\mathfrak{p}}_{\mathbb{C}}',Y') = \{0\}$
 for all $j$
 where $Y = \mu' \oplus Y'$ is the decomposition
 as a $K'$-module as before,
 we obtain the following commutative diagram
 by replacing $(T|_U)_{\ast}$ with $(\varphi_T)_{\ast}$:
\begin{alignat*}{5}
   &{\operatorname{Hom}}_{K}(\Exterior^j {\mathfrak{p}}_{\mathbb{C}},U)
   &&\overset \sim \to \,\,
   &&{\operatorname{Hom}}_{K}(\Exterior^j {\mathfrak{p}}_{\mathbb{C}},X)
   &&\overset\sim\to
   &&H^j ({\mathfrak{g}}, K; X)
\\
   &\qquad (\varphi_T)_{\ast} \downarrow
   &&
   && \circlearrowright
   &&
   &&\qquad \downarrow T_{\ast}
\\
   &{\operatorname{Hom}}_{K'}(\Exterior^j {\mathfrak{p}}_{\mathbb{C}}',U')
  &&\overset \sim \to
  &&{\operatorname{Hom}}_{K'}(\Exterior^j {\mathfrak{p}}_{\mathbb{C}}',Y)
  &&\overset\sim\to
  &&H^j ({\mathfrak{g}}', K'; Y).  
\end{alignat*}
Hence $T_{\ast}$ is a nonzero map
 if and only if $(\varphi_T)_{\ast}$ is nonzero.  
Since the bilinear map \eqref{eqn:Poincare} is a perfect pairing, 
 we conclude Theorem \ref{thm:171556}.  
\end{proof}

\begin{remark}
\label{rem:cohKK}
\begin{enumerate}
\item[{\rm{(1)}}]
The nonvanishing assumption of $\varphi_T$
 in the first statement of Theorem \ref{thm:171556}
 can be reformulated as the nonvanishing
 of the 
\index{B}{KspectrumKprime@$(K,K')$-spectrum}
 $(K,K')$-spectrum
 (see Section \ref{subsec:Kspec})
 of the symmetry breaking operator $T$
at $(\mu,\mu')$.  
\item[{\rm{(2)}}]
The verification of the ${\mathfrak{p}}$-vanishing assumption
 of $\varphi_T$ in the second statement of Theorem \ref{thm:171556} 
 reduces to a computation
 of finite-dimensional representations
 of compact Lie groups $K$ and $K'$.  
\item[{\rm{(3)}}]
If we set $R:=\dim_{\mathbb{C}}({\mathfrak{u}} \cap {\mathfrak{p}}_{\mathbb{C}})$
 and $R':=\dim_{\mathbb{C}}({\mathfrak{u}}'\cap {\mathfrak{p}}_{\mathbb{C}}')$,
 then the isomorphisms \cite[Cor.~3.7]{VZ} show
\begin{alignat*}{3}
&{\operatorname{Hom}}_{K}(\Exterior^j {\mathfrak{p}}_{\mathbb{C}},\mu)
   &&\simeq \,\,
   &&{\operatorname{Hom}}_{L \cap K}(\Exterior^{j-R} ({\mathfrak{l}}_{\mathbb{C}} \cap {\mathfrak{p}}_{\mathbb{C}}),{\mathbb{C}}), 
\\
  &{\operatorname{Hom}}_{K'}(\Exterior^j {\mathfrak{p}}_{\mathbb{C}}',\mu')
   &&\simeq
   &&{\operatorname{Hom}}_{L' \cap K'}(\Exterior^{j-R'}
   ({\mathfrak{l}}_{\mathbb{C}}' \cap {\mathfrak{p}}_{\mathbb{C}}'),{\mathbb{C}}).  
\end{alignat*}
\end{enumerate}
\end{remark}

\subsection{Nonvanishing bilinear forms on $(\mathfrak{g}, K)$-cohomologies
 via symmetry breaking for $(G,G')=(O(n+1,1),O(n,1))$}
\label{subsec:cohgkex}
\subsubsection{Nonvanishing theorem for $O(n+1,1) \downarrow O(n,1)$}
\label{subsec:cohbi}
In this section,
 we apply the general result
 (Theorem \ref{thm:171556})
 to the pair $(G,G')=(O(n+1,1),O(n,1))$.

In Proposition \ref{prop:gKq} in Appendix I, 
 we shall see
 that if $\Pi$ is an irreducible unitary representation
 of $G=O(n+1,1)$
 with $H^{\ast}({\mathfrak{g}},K;\Pi_K) \ne \{0\}$, 
 then the smooth representation $\Pi^{\infty}$ must be isomorphic
 to $\Pi_{\ell, \delta}$ defined in \eqref{eqn:Pild}
 for some $0 \le \ell \le n+1$ and $\delta \in \{\pm\}$.  
Thus, 
 we shall apply Theorem \ref{thm:171556}
 to the representations $\Pi_{\ell, \delta}$ of $G$
 and similar representations $\pi_{m,\varepsilon}$ 
 of the subgroup $G'=O(n,1)$.

In what follows,
 by abuse of notation,
 we denote an admissible  smooth representation and its underlying $({\mathfrak{g}},K)$-module by the same letter 
 when we discuss their $({\mathfrak{g}},K)$-cohomologies.

.
\begin{theorem}
\label{thm:gKOn}
Let $(G,G')=(O(n+1,1),O(n,1))$, 
 $0 \le i \le n$, 
 and $\delta \in \{ \pm \}$.  
Let $T:=A_{i,i}$ be the symmetry breaking operator 
 $\Pi_{i,\delta} \to \pi_{i,\delta}$
 given  in Proposition \ref{prop:AiiAq}.  
\begin{enumerate}
\item[{\rm{(1)}}]
$T$ induces bilinear forms
\[
   B_T \colon 
   H^j({\mathfrak{g}}, K; \Pi_{i,\delta}) 
   \times 
   H^{n-j}({\mathfrak{g}}', K'; \pi_{n-i,(-1)^n \delta})
   \to {\mathbb{C}}
\quad
\text{for all $j$.}
\]
\item[{\rm{(2)}}]
The bilinear form $B_T$ is nonzero
 if and only if $j=i$ and $\delta=(-1)^i$.  
\end{enumerate}
\end{theorem}

\begin{remark}
A similar theorem was proved 
 by B.~Sun \cite{S}
 for the $(\mathfrak{g},K)$-cohomology with nontrivial coefficients
 of a tempered representation
 of the  pair $(GL(n,\mathbb R), GL(n-1,\mathbb R))$. 
\end{remark}

We begin with the computation of the $({\mathfrak{g}},K)$-cohomologies
 of the irreducible representation $\Pi_{\ell, \delta}$ of $G=O(n+1,1)$.  
\begin{lemma}
\label{lem:172145}
Suppose $0 \le \ell \le n+1$, 
 $j \in {\mathbb{N}}$, 
 and $\delta \in \{\pm\}$. 
Then 
\[
   H^j({\mathfrak{g}},K; \Pi_{\ell, \delta})
  =
\begin{cases}
{\mathbb{C}} \quad &\text{if $j=\ell$ and $\delta=(-1)^{\ell}$, }
\\
\{0\} &\text{otherwise.}
\end{cases}
\]
\end{lemma}

In view of Theorem \ref{thm:LNM20} (4), 
 we have:
\begin{example}
\label{ex:gKchiOn}
For $G'=O(n,1)$, 
 we have $\pi_{n,(-1)^n} \simeq \chi_{-,(-1)^n}$ from 
 Theorem \ref{thm:LNM20} (4).  
In turn,
 the assertion 
$  
   H^{n}({\mathfrak{g}}', K'; \chi_{-,(-1)^n})
   \simeq
  {\mathbb{C}} 
$ from Lemma \ref{lem:172145} 
 corresponds to the equation \eqref{eqn:gKtop}
 by Example \ref{ex:chiOn}.  
\end{example}

By Proposition \ref{prop:161655} in Appendix I, 
 Lemma \ref{lem:172145} may be reformulated
 in terms of the cohomologically induced representations
\[
  (A_{{\mathfrak{q}}_i})_{a b}
  =
  A_{{\mathfrak{q}}_i} \otimes \chi_{a b}
  \simeq
  {\mathcal{R}}_{{\mathfrak{q}}_i}^{S_i} (\chi_{a b} \otimes {\mathbb{C}}_{\rho({\mathfrak{u}}_i)})
\]
 (see Section \ref{subsec:Aqgeneral} for notation)
 as follows:
\begin{lemma}
\label{lem:1721456}
Suppose $0 \le i \le [\frac{n+1}{2}]$ and $j \in {\mathbb{N}}$.  
Then we have
\begin{alignat*}{5}
H^j({\mathfrak{g}}, K;(A_{\mathfrak{q}_i})_{++})
=& {\mathbb{C}}\quad
&&\text{if $j=i$} 
&&\in 2{\mathbb{N}};
&&=\{0\}\quad
\text{otherwise,}
\\
H^j({\mathfrak{g}}, K;(A_{\mathfrak{q}_i})_{+-})
=& {\mathbb{C}}\quad
&&\text{if $j=i$}
&&\in 2{\mathbb{N}}+1;
&&=\{0\}\quad
\text{otherwise,}
\\
H^j({\mathfrak{g}}, K;(A_{\mathfrak{q}_i})_{-+})
=& {\mathbb{C}}\quad
&&\text{if $j=n+1-i$} 
&&\in 2{\mathbb{N}};
&&=\{0\}\quad
\text{otherwise,}
\\
H^j({\mathfrak{g}}, K;(A_{\mathfrak{q}_i})_{--})
=& {\mathbb{C}}\quad
&&\text{if $j=n+1-i$} 
&&\in 2{\mathbb{N}}+1;
&&=\{0\}\quad
\text{otherwise.}
\end{alignat*}
\end{lemma}

\begin{proof}
[Proof of Lemma \ref{lem:1721456}]
We recall from Theorem \ref{thm:LNM20} (3)
 (see also Proposition \ref{prop:161655} in Appendix I)
 that the irreducible $G$-module
$\Pi_{i,\delta}$
 contains $\mub(i,\delta) \simeq \Exterior^i ({\mathbb{C}}^{n+1}) \boxtimes \delta$
 as its minimal $K$-type.  
By \cite{VZ}, 
 we have then a natural isomorphism
\[
   {\operatorname{Hom}}_K
   (\Exterior^j {\mathfrak{p}}_{\mathbb{C}}, 
    \mub (i,\delta))
   \simeq
   H^j({\mathfrak{g}}, K; \Pi_{i,\delta}).  
\]

On the other hand, 
 the adjoint action of $K=O(n+1) \times O(1)$
 on ${\mathfrak{p}}_{\mathbb{C}} \simeq {\mathbb{C}}^{n+1}$
 gives rise to the $j$-th exterior tensor representation
\[
 \Exterior^j({\mathfrak{p}}_{\mathbb{C}}) 
 \simeq 
  \Exterior^j({\mathbb{C}}^{n+1}) \boxtimes (-1)^j.  
\]
Now the lemma follows.
\end{proof}

\begin{lemma}
\label{lem:phiT}
Let $\varphi_T$ be the $K'$-homomorphism
 defined in \eqref{eqn:phiT}
 for the symmetry breaking operator
 $T \colon \Pi_{i,\delta} \to \pi_{i,\delta}$
 in Theorem \ref{thm:gKOn}.  
Then $\varphi_T$ is ${\mathfrak{p}}$-nonvanishing
 at degree $j$ (Definition \ref{def:pnonvan})
 if and only if $j=i$ and $\delta=(-1)^i$.  
\end{lemma}
\begin{proof}
Similarly to the $G$-module $\Pi_{i,\delta}$, 
  the $G'$-module
$\pi_{i,\delta}$ contains 
$
   \mub(i,\delta)' \simeq \Exterior^i ({\mathbb{C}}^n) \boxtimes \delta
$
 as its minimal $K$-type.  
Then $\varphi_T$ in Theorem \ref{thm:171556} 
 amounts to a nonzero multiple
 of the projection 
 (see \eqref{eqn:Tii1}), 
\[
  \pr ii \colon 
  \Exterior^i({\mathbb{C}}^{n+1}) \boxtimes \delta 
  \to 
  \Exterior^i({\mathbb{C}}^{n}) \boxtimes \delta.  
\]
Then $(\varphi_T)_{\ast}$ is a nonzero multiple
 of the natural map from 
\begin{align*}
& {\operatorname{Hom}}_{O(n+1) \times O(1)}
  (\Exterior^j({\mathbb{C}}^{n+1})\boxtimes (-1)^j, \Exterior^i({\mathbb{C}}^{n+1})\boxtimes \delta)
\\
\intertext{to}
& {\operatorname{Hom}}_{O(n) \times O(1)}
  (\Exterior^j({\mathbb{C}}^{n})\boxtimes (-1)^j, \Exterior^i({\mathbb{C}}^{n})\boxtimes \delta)
\end{align*}
induced by the projection $\pr i i$.  
Now the lemma is clear.  
\end{proof}
We are ready to apply the general result (Theorem \ref{thm:171556})
 to prove Theorem \ref{thm:gKOn}.  

\begin{proof}
[Proof of Theorem \ref{thm:gKOn}]
By Example \ref{ex:chiOn}, 
 we have an isomorphism
 $\chi \simeq \chi_{-, (-1)^n}$ as $({\mathfrak{g}}',K')$-modules.  
Then it follows from Theorem \ref{thm:LNM20} (5) and (6) 
 that there are natural $G'$-isomorphisms:
\[
  \pi_{i,\delta}^{\vee} \otimes \chi_{-,(-1)^n}
  \simeq
  \pi_{i,\delta} \otimes \chi_{-,(-1)^n} \simeq \pi_{n-i,(-1)^n \delta}.  
\]
Thus Theorem \ref{thm:gKOn} (1) follows from Proposition \ref{prop:171543}.  
It then follows from Lemma \ref{lem:phiT}
 that Theorem \ref{thm:gKOn} (2) holds
 as a special case of Theorem \ref{thm:171556}.  
\end{proof}

In Proposition \ref{prop:161655}, 
 we shall see that the underlying $({\mathfrak{g}}',K')$-module 
 of $\pi_{n-i, (-1)^n\delta}$ is isomorphic to 
 $(A_{{\mathfrak{q}}_i'})_{-,(-1)^n\delta}$
 if $0 \le i \le [\frac n 2]$. 
The symmetry breaking operator 
$
   A_{i,i} \colon \Pi_{i,\delta} \to \pi_{i,\delta}$
 given in Proposition \ref{prop:AiiAq}
 induces a $(\mathfrak{g}', K')$-homomorphism $(A_{{\mathfrak{q}}_i})_{+,\delta} \to (A_{{\mathfrak{q}}_i'})_{+,\delta}$.

\begin{corollary}
\label{cor:12.14}
If  $0 \le 2i \le n$, 
 then the symmetry breaking operator $A_{i,i}\colon \Pi_{i,\delta} \to \pi_{i,\delta}$
 induces bilinear forms
\[
   H^j({\mathfrak{g}}, K; (A_{{\mathfrak{q}}_i})_{+,\delta})
  \times   
   H^{n-j}({\mathfrak{g}}', K'; (A_{{\mathfrak{q}}_i'})_{-,(-1)^n\delta})
  \to {\mathbb{C}}
\]
and linear maps
\[
  H^j({\mathfrak{g}}, K; (A_{{\mathfrak{q}}_i})_{+,\delta})
  \to  
   H^{j}({\mathfrak{g}}', K'; (A_{{\mathfrak{q}}_i'})_{+,\delta})
\]
for all $j$. 
They are nontrivial if and only if $j=i$ and $\delta=(-1)^i$.
\end{corollary}

Composing the symmetry breaking operators we deduce the following.  

\begin{corollary}
\label{cor:12.14.com}
If $0 \le 2i \le n$ and $H=O(n+1-i,1)$, 
 then the composition of the symmetry  breaking operators
 induces a linear map
\[
  H^{j}({\mathfrak{g}}, K; (A_{{\mathfrak{q}}_i})_{+,\delta})
  \to  
   H^{j}({\mathfrak{h}}, K\cap H';  (A_{{\mathfrak{q}}_i \cap \mathfrak{h}})_{+,\delta})
 \quad \text{for all } j. 
\]
It is nontrivial if and only if $j=i$ and $\delta=(-1)^{n+1-i}$.  
\end{corollary}

\begin{remark}
\label{rem:TongWang}
Y.~Tong and S.~P.~Wang \cite{TW} considered representations of $SO(n+1,1)$
 with nontrivial $({\mathfrak{g}},K)$-cohomology which are  $SO(n-i) \times SO(n+1-i,1)$-distinguished. Independently 
  S.~Kudla and J.~Millson \cite{KM2}
 considered representations
 of $O(n+1,1)$ with nontrivial $({\mathfrak{g}},K)$-cohomology
 which are  $O(n-i) \times O(n+1-i,1)$-distinguished. 
Since $O(n-i)$ commutes with $O(n-i+1,1)$,
 we have an action of $O(n-i)$
 on ${\operatorname{Hom}}_{O(n-i+1,1)}(\Pi_{i,\delta},{\mathbb C})$
 and ${\operatorname{Hom}}_{O(n-i+1,1)}(\Pi_{i,\delta},{\mathbb C})^{O(n-i)}$
 is isomorphic to ${\operatorname{Hom}}_{O(n-i)\times O(n-i+1,1)}(\Pi_{i,\delta},{\mathbb C})$.
 By results in \cite{KM2} this induces a nontrivial linear map
 on the $({\mathfrak{g}},K)$-cohomology.  
\end{remark}

\subsubsection{Special Cycles}
Geometric, topological and arithmetic properties of hyperbolic symmetric spaces $X_\Gamma=\Gamma\backslash O(n+1,1)/K$ for a discrete subgroup $\Gamma $ have been studied extensively using representation theoretic and geometric techniques. 
See for example \cite{BC, BMM} and references therein. 
If $X_\Gamma$ is compact, 
 then the Matsushima--Murakami formula 
 (\cite[Chap.~VII, Thm.~3.2]{BW})
shows
\[
   H^{\ast}(X_\Gamma, \bC) 
 \simeq \bigoplus_{\Pi \in \widehat{G}}  m(\Gamma, {\Pi})
   H^{\ast}({\mathfrak g},K; \Pi_K), 
\]
where
\index{A}{Ghat@$\widehat G$, unitary dual|textbf}
 $\widehat G$ is the set of equivalence classes 
 of irreducible unitary representations of $G$ 
({\it{i.e.}}, the {\it{unitary dual}} of $G$), 
 and we set for $\Pi \in \widehat G$
\[ 
  m(\Gamma,\Pi) := \dim_{\mathbb{C}} {\operatorname{Hom}}_G (\Pi,L^2(\Gamma\backslash G)). 
\]
By abuse of notation,
 we shall omit the subscript $K$ in the underlying $({\mathfrak g},K)$-module
 $\Pi_K$ of $\Pi$
 when we discuss its $({\mathfrak g},K)$-cohomologies.

In Proposition \ref{prop:gKq} in Appendix I,
 we shall show that every irreducible unitary representations with nontrivial $(\mathfrak{g},K)$-cohomology is isomorphic to a representations ${\Pi}_{i,\delta}$
 for some $i$ and $\delta \in \{\pm\}$, 
 see also Theorem \ref{thm:LNM20} (9). 
Thus
\[
   H^{\ast}(X_\Gamma, \bC) 
   = \bigoplus_{i,\delta} m(\Gamma, {\Pi}_{i,\delta})
   H^{\ast}({\mathfrak{g}},K; {\Pi}_{i,\delta}).  
\]

To obtain arithmetic information about the cohomology and the homology of $X_\Gamma $,
special cycles, 
{\it{i.e.,}} orbits of subgroups $H \subset G$ on $X_\Gamma$, 
 and their homology classes are frequently used. 
Suppose $0 \le 2i \le n+1$.  
We let 
\[
 G_i=O(n+1-i,1), 
\quad
 K_i :=K \cap G_i \simeq O(n+1-i) \times O(1), 
\]
 and $X_i$ be the Riemannian symmetric space $G_i/K_i$. 
Let $b_i :=n+1-i$, 
 the dimension of $X_i$.  
We set $\delta= (-1)^{n+1-i}$.  
By Corollary \ref{cor:12.14.com}  there  exists a nontrivial linear map  
$A^{n+1-i,n+1-i}:$
\[  
    H^{n+1-i}({\mathfrak{g}}, K; (A_{{\mathfrak{q}_{n+1-i}}})_{+,\delta})
  \to  
   H^{n+1-i}({\mathfrak{g}}_i, K_i; (A_{{\mathfrak{q}_{n+1-i}} \cap (\mathfrak{g}_i)_{\mathbb{C}}})_{+,\delta}). 
\]
Note that 
$(A_{{\mathfrak{q}}_{n+1-i} \cap ({\mathfrak{g}}_i)_{\mathbb{C}}})_{+,\delta}$ is one-dimensional
 and the image of $A^{n+1-i,n+1-i}$ is isomorphic to
\[
 {\operatorname{Hom}}_{K_i}
 (\Exterior^{n+1-i}({\mathfrak{p}}_{\mathbb{C}}\cap({\mathfrak{g}}_i)_{\mathbb{C}}),\chi_{+,\delta})
 \simeq  
 {\operatorname{Hom}}_{K_i}
 (\Exterior^{n+1-i}({{\mathbb C}^{n+1-i})\boxtimes {\bf{1}},
 {\mathbf 1}}). \]

Since the nonzero element of 
\[
   {\operatorname{Hom}}_{K_i}
   (\Exterior^{n+1-i}({{\mathbb C}^{n+1-i})\boxtimes {\bf{1}}, {\mathbf 1}})
\]
 gives a volume form on the symmetric space $X_i=G_i/K_i$, 
 this suggests
 that the homology classes defined
 by the orbits of $O(n+1-i,1)$
 for $0 \le 2i \leq n+1$ on $X_\Gamma$ are related
 to the contribution of $H^{n+1-i}({\mathfrak{g}},K;\Pi_{i,\delta})$
 to the cohomology of $X_\Gamma$.  
The work of S.~Kudla and J.~Millson confirms this.
We sketch their results   following the exposition in \cite{KM2, KM-I, KM-II}.

We have an embedding 
\[ 
     \iota_{X_i }: X_i \hookrightarrow X=G/K.   
\]
            
We fix an orientation of $X$ and $X_i$
 which is invariant under the connected component of $G$ respectively $G_i$.    Let $\mathbb A$ be the adels of the real number field $\mathbb K.$  
Then
\[ 
     X_{\mathbb A} = X \otimes G(\mathbb A_f) 
\]  
is the adelic  symmetric space. 
We set $G^+:= \prod SO_0(p,q)$ where we take the product over all real places
 of $\mathbb K $ and 
 $G^+(\mathbb K) := G(\mathbb K) \cap G^+G(\mathbb A_f)$.  
Then
 \[ H^*(G(\mathbb K) \backslash G(\mathbb A_f);\bC) = H^*({\mathfrak
 g},K; C^\infty(G(\mathbb Q)\backslash G(\mathbb A) ))\]
and 
\[
   H^*(G^+(\mathbb K)\backslash G(\mathbb A) /KK_f;\bC)
  =H^*(X_\mathbb A; \mathbb C)^{K_f}.  
\]
The cohomology here is the de Rham cohomology if  $X_\mathbb A$ is compact,
 otherwise the cohomology with compact support.

Following the exposition and notation in \cite[Sect.~2]{KM-II}
 we have an inclusion
\[
 \iota_{X_i} : \ X_i \times G_i(\mathbb A_f) \rightarrow X \times G({\mathbb A}_f)
\]
which is equivariant under the right action of $G({\mathbb A}_f)$. 
For $g \in G({\mathbb A}_f)$ we obtain a special cycle   
\[
     X_{i,g} = X_i G_i(\mathbb A_f) / (gK_f g^{-1} \cap G_i({\mathbb A}_f)).
\]
Consider the subspace $SX_i(X_{\mathbb A})$  spanned by special cycles
 in the homology group $H_i(X_{\mathbb A})$.

\medskip
We now assume that all but one factor of $G_\infty $ is compact and thus that $X_{\mathbb A}/K_f$ is compact.
Using the theta correspondence,
 S.~Kudla and J.~Millson show that there exist a subgroup $K_f$ and  nontrivial homomorphisms 
\[ 
\Psi \colon H^{b_i}({\mathfrak g},K; \Pi) 
            \rightarrow H^{b_i}(X_{\mathbb A}/K_f;\bC) \subset H^{b_i}(X_{\mathbb A})
\]
for some irreducible representation $\Pi$ of $G$.  

\medskip
Using integration, 
 S.~Kudla and J.~Millson \cite{KM2}, \cite[Thm.~7.1]{KM-II} prove the following:
\begin{theorem} 
\label{theorem:kudla+millson}
There exists a nontrivial pairing
\[ 
\Psi(H^{n+1-i}({\mathfrak g},K; \Pi) ) \times SX_i(X_{\mathbb A})  \rightarrow {\mathbb C}. 
\]
\end{theorem}

\begin{remark}
\begin{enumerate}
\item[{\rm{(1)}}]
As we see in Theorem \ref{thm:LNM20} (9), 
 Lemma \ref{lem:1721456}
 and Proposition \ref{prop:161655}, 
 the irreducible representation $\Pi$ of $G$
 with $H^{n+1-i}({\mathfrak{g}},K; \Pi) \ne \{0\}$
 must be of the form
\[
   \Pi \simeq \Pi_{n+1-i, (-1)^{n+1-i}}, 
\]
namely,
 $\Pi_K \simeq (A_{\mathfrak{q}_i})_{-,(-1)^{n+1-i}}$.  
\item[{\rm{(2)}}]
The nontrivial pairing in Theorem \ref{theorem:kudla+millson}
 defines an $O(n+1-i,1)$-invariant linear functional
 on the irreducible $G$-module $\Pi_{n+1-i, (-1)^{n+1-i}}$
 which is nontrivial on the minimal $K$-type. 
\end{enumerate}
\end{remark} 

\newpage
\section{A conjecture: Symmetry breaking for irreducible representations
 with regular integral infinitesimal character}
\label{sec:conjecture}

We conjecture that Theorems \ref{thm:SBOvanish} and \ref{thm:SBOone}
 hold  in more generality. 
We  will formalize and explain  this conjecture in this chapter more precisely and provide some supporting evidence.

As before we assume
that $G=O(n+1,1)$ and $G'=O(n,1)$.


\subsection{Hasse sequences and standard sequences
 of irreducible representations with regular integral infinitesimal character
 and their Langlands parameters} 
\label{subsec:Hasseps}
\index{B}{regularintegralinfinitesimalcharacter@regular integral infinitesimal character}

Before  stating the conjecture we define Hasse sequences and standard sequences  of irreducible representations, 
 and collect more information about the representations which occur in the Hasse and standard sequences.
In Chapter \ref{sec:aq} (Appendix I) we determine their $\theta$-stable parameters.

\subsubsection{Definition of Hasse sequence
 and standard sequence}
\begin{definitiontheorem}
[Hasse sequence]
\label{def:UHasse}
Let $n=2m$ or $2m-1$.  
For every irreducible finite-dimensional representation $F$
 of the group $G=O(n+1,1)$,
 there exists uniquely
 a sequence 
\begin{center}
\begin{tabular}{ccccccccccc}
& &$U_0$&, &\dots   & , & $U_{m-1} $ &, & $U_{m}$ 
\end{tabular}
\end{center}
of irreducible admissible smooth representations 
\index{A}{UHasse@$U_i(F)$, Hasse sequence starting with $F$}
 $U_i\equiv U_i(F)$ of $G$ such that 
\begin{enumerate}
\item  $U_0 \simeq F$; 
\item  consecutive representations are composition factors of a principal series representation;
\item
$U_i$ $(0 \le i \le m)$ are pairwise inequivalent as $G$-modules.  
\end{enumerate}

We refer to the sequence 
\begin{center}
\begin{tabular}{ccccccccccc}
& &$U_0$&, &\dots   & , & $U_{m-1} $ &, & $U_{m}$ 
\end{tabular}
\end{center}
 as the 
\index{B}{Hassesequence@Hasse sequence|textbf}
{\it {Hasse sequence}}
 of irreducible representations
 starting with the finite-dimensional representation $U_0=F$.  
We shall write $U_j(F)$ for $U_j$
 if we emphasize the sequence $\{U_j(F)\}$
 starts with $U_0=F$.  
\end{definitiontheorem}

\begin{proof}
[Sketch of the proof]  
D.~Collingwood \cite[Chap.~6]{C} computed embeddings
 of irreducible Harish-Chandra modules
 into principal series representations
 for all connected simple groups of real rank one, 
 which allowed him to define a diagrammatic description 
 of irreducible representations
 with regular integral infinitesimal character 
 of the connected group
 $G_0=SO_0(n+1,1)$.
For the disconnected group $G=O(n+1,1)$,
 we can determine similarly the composition factors
 of principal series representations,
 as in Theorems \ref{thm:171426} and \ref{thm:171425} below
 (see Sections \ref{subsec:SO}--\ref{subsec:PiSO} in Appendix II
 for the relationship 
 between irreducible representations of the disconnected group $G=O(n+1,1)$
 and those of a normal subgroup of finite index).  
To show the existence and the uniqueness
 of the Hasse sequence, 
 we note that there exists uniquely a principal series representation
 that contains a given irreducible finite-dimensional representation 
$F$ 
 as a subrepresentation.  
Then there exists only one irreducible composition factor
 other than $F$, 
 which is defined to be $U_1$.  
Repeating this procedure, 
 we can find irreducible representations $U_2$, $U_3$, $\cdots$, 
 whence the existence and the uniqueness of the Hasse sequence is shown for the disconnected group $G=O(n+1,1)$.  
\end{proof}

As we have seen in Theorem \ref{thm:LNM20} (1)
 when $F$ is the trivial one-dimensional representation ${\bf{1}}$, 
 the representations $U_i$ and $U_{i+1}$
 in this sequence have different signatures.  
The standard sequence (Definition \ref{def:Pii})
 starting with ${\bf{1}}$ 
 has given an adjustment
 for the different signatures.  
Extending this definition for the sequence
 starting 
 with an {\it{arbitrary}} irreducible finite-dimensional representation $F$, 
  we define the standard sequence of irreducible representations
 starting with $F$ as follows:

\begin{definition}
[standard sequence]
\label{def:Hasse}
If 
\begin{center}
\begin{tabular}{ccccccccccc}
& &\ $U_0$&, &\dots   & , & $U_{m-1} $ &, & $U_{m}$ 
\end{tabular}
\end{center}
is the Hasse sequence
 starting with an irreducible finite-dimensional representation 
 $F$ of $G$, 
then we refer to 
\begin{center}
\begin{tabular}{ccccccccccc}
& ${\Pi}_0:=U_0$&, &\dots   & , 
& ${\Pi}_{m-1}:=U_{m-1}\otimes (\chi_{+-})^{m-1} $ &, 
& $\Pi_m:=U_{m}\otimes (\chi_{+-})^m$ 
\end{tabular}
\end{center}
as the
\index{B}{standardsequence@standard sequence|textbf}
 {\it{standard sequence}}
 of irreducible representations 
\index{A}{1PiideltaF@$\Pi_i(F)$, standard sequence starting with $F$|textbf}
${\Pi}_i={\Pi}_i(F)$
 starting with ${\Pi}_0=U_0=F$, 
 where $\chi_{+-}$  
 is the one-dimensional representation of $G$
 defined in \eqref{eqn:chiab}.  
\end{definition}

\begin{remark}
\label{rem:Hassereg}
Clearly,
 any $U_j(F)$ in the Hasse sequence
 (or any $\Pi_j(F)$ in the standard sequence)
 starting with an irreducible finite-dimensional representation $F$
 of $G$ has a regular integral ${\mathfrak{Z}}_G({\mathfrak{g}})$-infinitesimal character
 (Definition \ref{def:intreg}).  
\end{remark}

The next proposition follows readily from the definition.  
\begin{proposition}
[tensor product with characters]
\label{prop:HStensor}
Let $F$ be an irreducible finite-dimensional representation of $G$, 
 and $\chi$ a one-dimensional representation of $G$.  
Then the representations in the Hasse sequences 
 (and in the standard sequence)
 starting with $F$ and $F \otimes \chi$ have the following relations:
\begin{alignat*}{3}
&\text{\rm{(Hasse sequence)}}
\qquad
&&U_i(F) \otimes \chi &&\simeq U_i(F \otimes \chi), 
\\
&\text{\rm{(standard sequence)}}
\qquad
&&
\Pi_i(F) \otimes \chi &&\simeq \Pi_i(F \otimes \chi).  
\end{alignat*}
\end{proposition}

The Hasse sequences and the standard sequences starting 
 with one-dimensional representations of $G$
 are described as follows.  
\begin{example}
\label{ex:171455}
We recall from Theorem \ref{thm:LNM20}
 that 
\index{A}{1Piidelta@$\Pi_{i,\delta}$, irreducible representations of $G$}
$\Pi_{\ell, \delta}$
 ($0 \le \ell \le n+1$, $\delta \in \{\pm\}$)
 are irreducible representations
 of $G=O(n+1,1)$
 with ${\mathfrak{Z}}_G({\mathfrak{g}})$-infinitesimal character 
\index{A}{1parho@$\rho_G$}
 $\rho_G$.  
Then for each one-dimensional representation
 $F \simeq \chi_{\pm\pm}$ of $G$
 (see \eqref{eqn:chiab}), 
 the Hasse sequence
 $U_i(F)$ ($0 \le i \le [\frac{n+1}{2}]$) 
 that starts with $U_0(F) \simeq F$, 
 and the standard sequence $\Pi_i(F):=U_i(F) \otimes (\chi_{+-})^i$
 are given as follows.  
\index{A}{1Piidelta@$\Pi_{i,\delta}$, irreducible representations of $G$}
\begin{alignat*}{3}
U_i({\bf{1}})&=\Pi_{i,(-1)^i},
&&
\Pi_i({\bf{1}})&&=\Pi_{i,+},
\\
U_i(\chi_{+-})&=\Pi_{i,(-1)^{i+1}},
&&
\Pi_i(\chi_{+-})&&=\Pi_{i,-},
\\
U_i(\chi_{-+})&=\Pi_{n+1-i,(-1)^{i}},
&&
\Pi_i(\chi_{-+})&&=\Pi_{n+1-i,+},
\\
U_i(\chi_{--})&=\Pi_{n+1-i,(-1)^{i+1}}, 
\qquad
&&
\Pi_i(\chi_{--})&&=\Pi_{n+1-i,-}.  
\end{alignat*}
\end{example}


\subsubsection{Existence of Hasse sequence}

In Section \ref{subsec:13.2}, 
 we formalize a conjecture
 about when 
\[
  \operatorname{Hom}_{G'}(\Pi|_{G'}, \pi ) \ne \{0\}
\]
for $\Pi \in \operatorname{Irr} (G)$
 and $\pi \in \operatorname{Irr} (G')$
 that have regular integral infinitesimal characters
 by using the standard sequence
 (Definition \ref{def:Hasse}).  
The formulation is based on the following theorem
 which asserts 
 that the converse statement to Remark \ref{rem:Hassereg}
 is also true.  

\begin{theorem}
\label{thm:Hassereg}
Any irreducible admissible representation of $G$
 of moderate growth
 with regular integral ${\mathfrak{Z}}_G({\mathfrak{g}})$-infinitesimal character is of the form 
 $U_j(F)$
 in the Hasse sequence
 for some $j$ $(0 \le j \le [\frac{n+1}{2}])$
 and for some irreducible finite-dimensional representation $F$ of $G$.

Similarly,
 any irreducible admissible representation of $G$
 of moderate growth 
 with regular ${\mathfrak{Z}}_G({\mathfrak{g}})$-infinitesimal character
 is of the form $\Pi_j(F')$
 in the standard sequence
 for some $j$ $(0 \le j \le [\frac{n+1}{2}])$
 and for some irreducible finite-dimensional representation $F'$ of $G$.  
\end{theorem}

The proof of Theorem \ref{thm:Hassereg} follows from the classification
 of $\operatorname{Irr}(G)$
 (Theorem \ref{thm:irrG} in Appendix I)
 and the Langlands parameter of the representations
 in the Hasse sequence below
 (see also Theorem \ref{thm:1808116}).  

\subsubsection{Langlands parameter of the representations
 in the Hasse sequence}
\label{subsec:LanHasse}
Let $F$  be an irreducible finite-dimensional representation of $G=O(n+1,1)$. 
 We now determine the Langlands parameter of the representations in the Hasse sequence  $\{ U_i(F) \}$ (and the standard sequence $\{ \Pi_i(F) \}$)  for $0 \le i \le [\frac {n+1}2]$  and their $K$-types.

We use the parametrization of the finite-dimensional representation
 of $O(n,1)$ introduced in Section \ref{subsec:fdimrep} in Appendix I.

We begin with the case
 where $F$ is obtained from an irreducible representation 
 of $O(n+2)$ of 
\index{B}{type1@type I, representation of $O(N)$}
 type I (Definition \ref{def:type})
 via the unitary trick.
The description of $U_i(F)$ and $\Pi_i(F)$
 for more general $F$ can be derived from this case
 by taking the tensor product with one-dimensional representations 
 $\chi_{\pm\pm}$
 of $G$, 
 see Theorem \ref{thm:1714107} below.  

\vskip 1pc
\par\noindent
{\bf{\underline{Case 1.}}}\enspace
$n=2m$ and $G=O(2m+1,1)$.  

For $F \in \widehat{O(n+2)}$ of type I, 
 we define $\sigma^{(i)} \equiv \sigma^{(i)}(F) \in \widehat{O(n)}$
 of type I for $0 \le i \le m = \frac n 2$
 as follows.  
We write $F= \Kirredrep {O(n+2, {\mathbb{C}})}{s}$
 with 
\index{A}{0tLambda@$\Lambda^+(N)$, dominant weight}
\[
   s=(s_{0},\cdots,s_{m},0^{m+1}) \in \Lambda^+(n+2)\equiv \Lambda^+(2m+2)
\]
 as in \eqref{eqn:CWOn}, 
 and regard it as an irreducible finite-dimensional representation
 of $G=O(n+1,1)$.  
We set
\[
   \sigma^{(i)} : = \Kirredrep {O(n)}{s^{(i)}} \in \widehat{O(n)}
\quad
  \text{for }\,\, 0 \le i \le m, 
\]
where
 $s^{(i)}\in\Lambda^+(n)\equiv \Lambda^+(2m)$ is given 
 for $0 \le i \le m$ as follows:
\begin{equation}
\label{eqn:sin=2m}
  s^{(i)}
  :=
  (s_{0}+1,\cdots,s_{i-1}+1,\widehat {s_{i}},s_{i+1},\cdots,s_{m},0^m). 
\end{equation}
It is convenient to introduce the {\em extended Hasse sequence}
 $\{U_i \equiv U_i(F)\}$
 ($0 \le i \le 2m+1$) by defining 
\index{A}{1chipmm@$\chi_{--}=\det$}
\begin{equation}
\label{eqn:Uiext1}
   U_i(F):=U_{n+1-i}(F) \otimes \chi_{--}
\quad
\text{for  $m+1 \le i \le n+1 = 2m+1$}.  
\end{equation}

\begin{theorem}
[$n=2m$]
\label{thm:171426}
Given an irreducible finite-dimensional representation
 $\Kirredrep{O(n+2,{\mathbb{C}})}{s}$of $G=O(n+1,1)$
 with 
\[
s =(s_0, s_1, \cdots, s_m,0,\cdots,0)\in \Lambda^+(n+2)
(=\Lambda^+(2m+2)), 
\]
 there exists uniquely an extended
\index{B}{Hassesequence@Hasse sequence}
Hasse sequence $U_0$, $U_1$, $\cdots$, $U_{2m+1}$
 starting with the irreducible finite-dimensional representation
 $U_0=\Kirredrep{O(n+2,{\mathbb{C}})}{s}$.  
Moreover,
 the extended Hasse sequence $U_0$, $\cdots$, $U_{2m+1}$
 satisfies the following properties.  
\begin{enumerate}
\item[{\rm{(1)}}]
There exist exact sequences of $G$-modules:
\begin{alignat*}{2}
& 0\to U_i \to I_{(-1)^{i-s_i}}(\sigma^{(i)}, i-s_i) \to U_{i+1} \to 0
\quad
&&(0 \le i \le m), 
\\
& 0\to U_i \to I_{(-1)^{n-i-s_{n-i}}}(\sigma^{(n-i)} \otimes \det, i+s_{n-i}) \to U_{i+1} \to 0
\quad
&&(m \le i \le 2m).   
\end{alignat*}

\item[{\rm{(2)}}]
The 
\index{B}{Ktypeformula@$K$-type formula}
$K$-type formula of the irreducible $G$-module $U_i$
 $(0 \le i \le m)$ is given by
\[
  \bigoplus_b \Kirredrep{O(n+1)}{b} \boxtimes 
  (-1)^{\sum_{k=0}^m(b_k-s_k)}, 
\]
where $b=(b_0, b_1, \cdots, b_m,0,\cdots,0)$
 runs over $\Lambda^+(n+1)\equiv \Lambda^+(2m+1)$
 subject to 
\begin{align*}
& b_{0} \ge s_{0}+1 \ge b_{1} \ge s_{1}+1 \ge \cdots \ge b_{i-1} \ge s_{i-1}+1, 
\\
& s_{i}\ge b_{i} \ge s_{i+1} \ge b_{i+1} \ge \cdots \ge s_{m}\ge b_{m} \ge 0, 
\\
&b_{m} \in \{0,1\}.  
\end{align*}
In particular, 
 the 
\index{B}{minimalKtype@minimal $K$-type|textbf}
 minimal $K$-type(s) of the $G$-module $U_i$ $(0 \le i\le m)$
 are given as follows:

for $s_m=0$, 
\begin{multline*}
  \Kirredrep{O(n+1)}{s^{(i)},0}\boxtimes (-1)^{i-s_i}
\\
  =
  \Kirredrep{O(n+1)}{s_{0}+1, \cdots, s_{i-1}+1, \widehat{s_i}, s_{i+1}, \cdots,s_m,0^{m+1}}\boxtimes (-1)^{i-s_i};
\end{multline*}

for $s_m >0$, 
\begin{equation*}
\Kirredrep{O(n+1)}{s^{(i)},0} \boxtimes (-1)^{i-s_i}
\,\,
\text{ and }
\,\,
(\Kirredrep{O(n+1)}{s^{(i)},0} \otimes \det)\boxtimes (-1)^{i-s_i+1}.  
\end{equation*}
\end{enumerate}
\end{theorem}

\begin{proof}
[Sketch of the proof]
\begin{enumerate}
\item[(1)]
By the translation principle,
 the first exact sequence follows from Theorem \ref{thm:LNM20} (1)
 which corresponds to the case $F \simeq {\bf{1}}$.  
Taking its dual, 
 we obtain another exact sequence
\[
  0 \to U_{i+1} \to I_{(-1)^{i-s_i}}(\sigma^{(i)}, n-i+s_i)
    \to U_i \to 0
\quad
  \text{for $0 \le i \le m$, }
\]
because $U_i$ is self-dual.  
Taking the tensor product with the one-dimensional representation $\chi_{--}$ of $G$, 
 we obtain by \eqref{eqn:Uiext1} and by Lemma \ref{lem:IVchi}
 another exact sequence of $G$-modules:
\[
  0 \to U_{n-i} \to I_{(-1)^{i-s_i}}(\sigma^{(i)}\otimes \det, n-i+s_i)
    \to U_{n+1-i} \to 0.  
\]
Replacing $i$ ($0 \le i \le m$)
 by $n-i$ ($m \le n-i \le 2m$), 
 we have shown the second exact sequence.  
\item[(2)]
The $K$-type formula of the irreducible finite-dimensional representation 
 $U_0=\Kirredrep {O(n+1,1)}{s}$ of $G$ is known by the classical branching law
 (see Fact \ref{fact:ONbranch}).  
Since the $K$-type formula of the principal series representation
 is given by the Frobenius reciprocity
 which we can compute by using Fact \ref{fact:ONbranch} again,
 the $K$-type formula of $U_{i+1}$ follows inductively from that of $U_i$
 by the exact sequence 
 in the first statement.  
\end{enumerate}
\end{proof}
See also Theorem \ref{thm:171471} in Appendix I
 for another description of the irreducible representation $U_i(F)$
 in terms of $\theta$-{\it{stable parameters}}.

\begin{remark}
\label{rem:Utemp}
When $i=m$ and $n=2m$, 
 $s^{(i)}$ is of the form
\[
  s^{(m)} =(s_0+1, \cdots, s_{m-1}+1, 0^m) \in \Lambda^+(2m), 
\]
and therefore
 the irreducible $O(n)$-module 
$\sigma^{(m)}=\Kirredrep{O(2m)}{s^{(m)}}$ is of type Y
 (Definition \ref{def:OSO}).  
Hence we have an isomorphism
\begin{equation}
\label{eqn:smdet}
    \sigma^{(m)} \simeq \sigma^{(m)} \otimes \det
\end{equation}
 as $O(2m)$-modules by Lemma \ref{lem:typeY}.  
We recall from Theorem \ref{thm:171426} (1)
 that there is an exact sequence of $G$-modules as follows:
\[
  0 \to U_m \to I_{(-1)^{m-s_m}}(\sigma^{(m)}, m-s_m) \to U_{m+1} \to 0.  
\]
Taking the tensor product with the character $\chi_{--} \simeq \det$, 
 we obtain from \eqref{eqn:Uiext1} and Lemma \ref{lem:IVchi}
 another exact sequence of $G$-modules:
\[
  0 \to U_{m+1} \to I_{(-1)^{m-s_m}}(\sigma^{(m)}\otimes \det, m-s_m)
 \to U_{m} \to 0.  
\]
By \eqref{eqn:smdet}, 
 the principal series representations 
\[
I_{(-1)^{m-s_m}}(\Kirredrep{O(2m)}{s^{(m)}}, m-s_m)
\simeq
I_{(-1)^{m-s_m}}(\Kirredrep{O(2m)}{s^{(m)}}\otimes \det, m-s_m)
\]
 split into a direct sum
 of two irreducible $G$-modules $U_m$ and $U_{m+1}$
 (see also Theorem \ref{thm:IndV} (3) in Appendix I).  
\end{remark}

\vskip 1pc
\noindent
{\bf{\underline{Case 2.}}}\enspace
$n=2m-1$ and $G=O(2m,1)$.  

For $F \in \widehat{O(n+2)}$ of type I, 
 we define $\sigma^{(i)} \equiv \sigma^{(i)}(F) \in \widehat{O(n)}$
 for $0 \le i \le m-1 = \frac 1 2 (n-1)$
 as follows.  
We write $F= \Kirredrep {O(n+2)}{s}$
 with 
\[
s=(s_{0},s_1,\cdots,s_{m-1},0^{m+1}) \in \Lambda^+(n+2)\equiv \Lambda^+(2m+1), 
\]
 as in \eqref{eqn:CWOn}.  
Then we define $s^{(i)}\in \Lambda^+(n)\equiv \Lambda^+(2m-1)$
 ($0 \le i \le m-1$) by
\begin{equation}
\label{eqn:sin=2m-1}
  s^{(i)}
  :=
  (s_{0}+1,\cdots,s_{i-1}+1,\widehat {s_{i}},s_{i+1},\cdots,s_{m-1},0^m), 
\end{equation}
and define irreducible finite-dimensional representations by 
\[
   \sigma^{(i)} : = \Kirredrep {O(n)}{s^{(i)}} \in \widehat{O(n)}
\quad
  \text{for }\,\, 0 \le i \le m-1.  
\]
For later purpose, we set
\[
s^{(m)}
  :=
  (s_{0}+1,\cdots,s_{m-2}+1,1,0^{m-1}) \in \Lambda^+(n). 
\]
Then there is an isomorphism as $O(n)$-modules: 
\[
  \sigma^{(m-1)} \otimes \det
  \simeq
  \Kirredrep{O(n)}{s^{(m)}}.  
\]
It is convenient to introduce the {\em extended Hasse sequence } $\{U_i \equiv U_i(F)\}$
 ($0 \le i \le 2m$) by defining for
 $m+1 \le i \le 2m$ 
\begin{equation}
\label{eqn:Uiext2}
  U_i(F):= U_{n+1-i}(F) \otimes \chi_{-+}.  
\end{equation}
Implicitly,
 the definition \eqref{eqn:Uiext2} includes a claim 
 that there is an isomorphism
 of discrete series representations
 (cf. Remark \ref{rem:Umdisc} below):
\index{A}{1chipmp@$\chi_{-+}$}
\begin{equation}
\label{eqn:UmF-+}
   U_m(F) \simeq U_m(F) \otimes \chi_{-+}
\end{equation}
when $G=O(n+1,1)$ with $n=2m-1$.

We note that the one-dimensional representations
 $\chi_{--}$ and $\chi_{-+}$ 
 in \eqref{eqn:Uiext1} and \eqref{eqn:Uiext2}
 are chosen differently 
 according to the parity of $n$.

The proof of the following theorem goes similarly to that of Theorem \ref{thm:171426}.  
\begin{theorem}
[$n=2m-1$]
\label{thm:171425}
Given an irreducible finite-dimensional representation
 $\Kirredrep {O(n+2,{\mathbb{C}})}{s}$
 of $G=O(n+1,1)$ with 
\[
s=(s_{0},s_1,\cdots,s_{m-1},0^{m+1}) \in \Lambda^+(n+2)
 (=\Lambda^+(2m+1)), 
\]
 there exists uniquely an extended 
\index{B}{Hassesequence@Hasse sequence}
Hasse sequence $U_0$, $U_1$, $\cdots$, $U_{2m}$
 of $G=O(2m,1)$
 starting with the irreducible finite-dimensional representation
 $U_0=\Kirredrep{O(n+2,{\mathbb{C}})}{s}$.  
Moreover, 
 the extended Hasse sequence $U_0$, $U_1$, $\cdots$, $U_{2m}$
 satisfies the following properties.  
\begin{enumerate}
\item[{\rm{(1)}}]
There exist exact sequences of $G$-modules:
\begin{alignat*}{2}
0 &\to U_i \to I_{(-1)^{i-s_i}}(\sigma^{(i)}, i-s_i)
   \to U_{i+1} \to 0
\qquad
&&(0 \le i \le m-1), 
\\
0 & \to U_i \to I_{(-1)^{n-i-s_{n-i}}}(\sigma^{(n-i)} \otimes \det, i+s_{n-i})
  \to U_{i+1} \to 0
\qquad
&&(m\le i \le 2m-1).  
\end{alignat*}

\item[{\rm{(2)}}]
The $K$-type formula of the irreducible $G$-module $U_i$
 $(0 \le i \le m)$ is given by
\[
  \bigoplus_b \Kirredrep{O(n+1)}{b} 
  \boxtimes 
  (-1)^{\sum_{k=0}^{m} b_k - \sum_{k=0}^{m-1}s_k}, 
\]
where $b=(b_0, b_1, \cdots, b_{m-1},0,\cdots,0)$
 runs over $\Lambda^+(n+1)\equiv \Lambda^+(2m)$
 subject to the following conditions:
\begin{align*}
& b_{0} \ge s_{0}+1 \ge b_{1} 
  \ge s_{1}+1 \ge \cdots \ge b_{i-1} \ge s_{i-1}+1, 
\\
& s_{i}\ge b_{i} \ge s_{i+1} \ge b_{i+1} \ge \cdots \ge s_{m-1}\ge b_{m-1}\ge 0.
\end{align*}
In particular,
 the minimal $K$-type of the $G$-module $U_i$ $(0 \le i \le m)$ is given by 
\begin{multline*}
  \Kirredrep{O(n+1)}{s^{(i)},0} \boxtimes (-1)^i
\\
  =
  \Kirredrep{O(n+1)}{s_0+1, \cdots, s_{i-1}+1, \widehat{s_i}, 
s_{i+1}, \cdots, s_{m-1}, 0^{m+1}} \boxtimes (-1)^{i-s_i}.  
\end{multline*}
\end{enumerate}
\end{theorem}

\begin{remark}
\label{rem:Umdisc}
$U_m$ is a discrete series representation
 of $G=O(2m,1)$.  
\end{remark}

See also Theorem \ref{thm:171471b} in Appendix I
 for another description of the irreducible representation
 $U_i(F)$ in terms of $\theta$-stable parameters.

By applying Proposition \ref{prop:HStensor} and Lemma \ref{lem:IVchi},
 we may unify the first statement of Theorems \ref{thm:171426}
 and \ref{thm:171425}
 as follows.  
\begin{theorem}
\label{thm:1714107}
Let $F$ be an irreducible finite-dimensional representation of $G=O(n+1,1)$
 of type I
 (see Definition \ref{def:typeone} in Appendix I), 
 and $a,b \in \{\pm\}$.  
Then for $F_{a, b}:= F \otimes \chi_{a b}$, 
 there exists uniquely a Hasse sequence
 $U_i(F_{a,b})$  $(0 \le i \le [\frac{n+1}{2}])$ 
 starting with $U_0(F_{a, b})=F_{a, b}$.  
Moreover, 
 the irreducible $G$-modules $U_i(F_{a,b})$ occur 
 in the following exact sequence of $G$-modules
\[
   0 \to U_i(F_{a, b}) \to I_{a b (-1)^{i-s_i}}(\sigma_a^{(i)}, i-s_i)
      \to U_{i+1}(F_{a, b}) \to 0
\]
for $0 \le i \le [\frac{n-1}{2}]$. 
Here $\sigma_a^{(i)}=\sigma^{(i)}$ if $a=+;$
$\sigma^{(i)} \otimes \det$ if $a=-$.  
\end{theorem}

\medskip
\noindent
\begin{remark}
\label{remark:171478}
By \eqref{eqn:314}, 
 we have linear bijections for all $i$, $j$:
\[
{\operatorname{Hom}}_{G'}
 (U_i(F)|_{G'}, U_j'(F'))
\simeq
{\operatorname{Hom}}_{G'}
 (U_{n+1-i}(F)|_{G'}, U_{n-j}'(F') \otimes \chi_{+-}).  
\]
\end{remark}

\medskip
\begin{remark}
Using the definition of the extended Hasse sequence we also define an extended standard sequence.
\end{remark}
\medskip

By abuse of notation we will from now on not distinguish between Hasse sequences and extended Hasse sequences and refer to both as Hasse sequences. A similar convention applies to standard sequences.

\medskip
The following observation will be used 
 in Section \ref{subsec:Evidence4}
 for the proof of Evidence E.4
 of Conjecture \ref{conj:GPver1} below.  

\begin{proposition}
\label{prop:172009}
Suppose $F$ and $F'$ are irreducible finite-dimensional representations
 of $G=O(n+1,1)$ and $G'=O(n,1)$, 
respectively, 
 such that ${\operatorname{Hom}}_{G'}(F|_{G'}, F') \ne \{0\}$.  
Suppose the principal series representations
 $I_{\delta}(V,\lambda)$ of $G$
 and $J_{\varepsilon}(W,\nu)$ of $G'$
 contain $F$ and $F'$, 
 respectively,
 as subrepresentations. 
Then the following hold.  
\begin{enumerate}
\item[{\rm{(1)}}]
$[V:W]=1;$
\item[{\rm{(2)}}]
\index{A}{1psi@$\Psising$,
          special parameter in ${\mathbb{C}}^2 \times \{\pm\}^2$}
$(\lambda, \nu, \delta, \varepsilon) \in \Psising$
 (see \eqref{eqn:singset}), 
namely,
 the quadruple $(\lambda, \nu, \delta, \varepsilon)$ does not satisfy
 the generic parameter condition \eqref{eqn:nlgen}.  
\end{enumerate}
\end{proposition}
\begin{proof}
For the proof,
 we use a description 
 of irreducible finite-dimensional representations
 of the disconnected group $G=O(n+1,1)$
 in Section \ref{subsec:fdimrep}
 of Appendix I.  
In particular,
 using Lemma \ref{lem:161612}, 
we may write 
\[
  F = \Kirredrep {O(n+1,1)}{\lambda_0, \cdots, \lambda_{[\frac n 2]}}_{a, b}
\]
for some $(\lambda_0, \cdots, \lambda_{[\frac n 2]})\in \Lambda^+([\frac n 2]+1)$ and $a, b \in \{\pm\}$.  
By the branching rule 
 for $O(n+1,1) \downarrow O(n,1)$
 (see Theorem \ref{thm:ON1branch}), 
 an irreducible summand $F'$ of $F|_{O(n,1)}$ is of the form
\[
  F' = \Kirredrep {O(n,1)}{\nu_0, \cdots, \nu_{[\frac {n-1} 2]}}_{a, b}
\]
for some $(\nu_0, \cdots, \nu_{[\frac {n-1} 2]})\in \Lambda^+([\frac {n+1} 2])$ such that 
\begin{alignat*}{2}
&\lambda_{0} \ge \nu_{0} \ge \lambda_{1} \ge \cdots \ge \nu_{[\frac{n-1}2]} \ge 0\quad
&&\text{for $n$ odd}, 
\\
&\lambda_{0} \ge \nu_{0} \ge \lambda_{1} \ge \cdots \ge \nu_{[\frac{n-1}2]} \ge \lambda_{[\frac n 2]}
\quad
&&\text{for $n$ even}.   
\end{alignat*}
We recall 
 that for every irreducible finite-dimensional representation $F$
 of a real reductive Lie group
 there exists only one principal series representation
 that contains $F$ 
 as a subrepresentation.  
By Theorem \ref{thm:1714107} with $i=0$, 
 the unique parameter $(V, \delta, \lambda)$ is given
 by 
\[
\text{$V=\Kirredrep{O(n)}{\lambda_1, \cdots, \lambda_{[\frac n 2]}}$
 ($\otimes \det$ if $a=-$), 
 $\lambda=-\lambda_0$ and $\delta=a b (-1)^{-\lambda_0}$.}
\]  
Likewise, 
 the unique parameter $(W, \varepsilon, \nu)$ for $F'$
 is given by 
\[
 \text{$W=\Kirredrep{O(n-1)}{\nu_1, \cdots, \nu_{[\frac {n-1} 2]}}$
 ($\otimes \det$ if $a=-$), 
 $\nu=-\nu_0$, 
 and $\varepsilon = a b (-1)^{-\nu_0}$.  }
\]
Hence $[V:W]\ne 0$, 
 or equivalently,
 $[V:W]=1$ 
 by the branching rule 
 for $O(n) \downarrow O(n-1)$.  
Moreover, 
 $\delta \varepsilon = (-1)^{\nu_0 + \lambda_0}$
 and $\nu-\lambda= \lambda_0 -\nu_0 \in {\mathbb{N}}$.  
Hence the generic parameter condition \eqref{eqn:nlgen} fails, 
 or equivalently,
$(\lambda, \nu, \delta, \varepsilon) \in \Psising$. 
 
\end{proof}

\subsection{The Conjecture}
\label{subsec:13.2}
We propose a conjecture about when 
\[
{\operatorname{Hom}}_{G'}(\Pi|_{G'},\pi)={\mathbb{C}}
\]
where $\Pi \in {\operatorname{Irr}}(G)$ and $\pi \in {\operatorname{Irr}}(G')$ have regular integral infinitesimal characters
 (Definition \ref{def:intreg}).  
We give two formulations of the conjecture, 
 see Conjectures \ref{conj:GPver1} and \ref{conj:V2} below.  
Supporting evidence is given in Section \ref{subsec:SuppConj}.  

\subsubsection{Conjecture: Version 1}
We begin with a formulation of the conjecture
 in terms of a standard sequence
 (Definition-Theorem \ref{def:UHasse})
 of irreducible representations $\Pi_i$ of $G=O(n+1,1)$
 and that of irreducible representations $\pi_j$ of the subgroup $G'=O(n,1)$.  
We note that both $\Pi_i$ and $\pi_j$ have
 regular integral infinitesimal characters
 because both $F:=\Pi_0$ and $F':=\pi_0$ are
 irreducible {\it{finite-dimensional}} representations of $G$
 and $G'$, respectively. 

\begin{conjecture}
\label{conj:GPver1}
Let $F$ be an irreducible finite-dimensional representations
 of $G=O(n+1,1)$,
 and 
\index{A}{1PiideltaF@$\Pi_i(F)$, standard sequence starting with $F$}
$\{\Pi_i(F)\}$ be the standard sequence
 starting at $\Pi_0(F)=F$.  
Let $F'$ be an irreducible finite-dimensional representation of the subgroup $G'=O(n,1)$, 
 and $\{\pi_j(F')\}$ the standard sequence starting 
 at $\pi_0(F')=F'$.  
Assume that 
\[ 
   {\operatorname{Hom}}_{G'}(F|_{G'},F') \not = \{0\}.  
\]
Then the symmetry breaking for representations $\Pi_i(F)$, $\pi_j(F')$
 in the standard sequences  is represented graphically
 in Diagrams \ref{tab:SBodd} and \ref{tab:SBeven}. 
In the first row are representations of G, in the second row are representations of $G'$. 
Symmetry breaking operators are represented by arrows, 
 namely,
 there exist nonzero symmetry breaking operators
 if and only if there are arrows in the diagram.  
\end{conjecture}
\medskip

\begin{figure}[htp]
\color{black}
\caption{Symmetry breaking for $O(2m+1,1) \downarrow O(2m,1)$}
\begin{center}
\begin{tabular}{@{}c@{~}c@{~}c@{~}c@{~}c@{~}c@{~}c@{~}c@{~}c@{~}c@{}}
$\Pi_0(F)$
& 
&$\Pi_1(F)$
&  
&\dots
&
&$\Pi_{m-1}(F)$
& 
&$\Pi_{m}(F)$ 
\\
$\downarrow$ 
&$\swarrow$
& $\downarrow$
& $\swarrow$ 
& 
& $\swarrow$ 
& $ \downarrow $
&  $\swarrow $  
&  $\downarrow$ 
\\
$\pi_0(F')$& &$\pi_1(F')$
& 
&\dots 
&
& $\pi_{m-1}(F')$ 
& 
& $\pi_{m}(F')$ 
\end{tabular}
\end{center}
\label{tab:SBodd}
\end{figure}%

\medskip
\begin{figure}[htp]
\color{black}
\caption{Symmetry breaking for $O(2m+2,1) \downarrow O(2m+1,1)$ }
\begin{center}
\begin{tabular}{@{}c@{~}c@{~}c@{~}c@{~}c@{~}c@{~}c@{~}c@{~}c@{~}c@{~}c@{}}
$\Pi_0(F)$& &$\Pi_1(F)$ & & \dots &  & $\Pi_{m-1}(F)$& & $\Pi_{m}(F)$ & & $\Pi_{m+1}(F)$\\
$\downarrow$ &$\swarrow$& $\downarrow $& $\swarrow$ &  & $\swarrow$ & $ \downarrow $& $\swarrow $ & $\downarrow$ & $\swarrow$\\
$\pi_0(F')$& &$\pi_1(F')$& &\dots  & & $\pi_{m-1}(F')$ & & $\pi_{m}(F')$ 
\end{tabular}
\end{center}
\label{tab:SBeven}
\end{figure}%

\begin{remark}
Instead of using standard sequences to state the conjecture
 it may be also useful to rephrase it using extended  Hasse sequences.
\end{remark}

\subsubsection{Conjecture: Version 2}
We rephrase the conjecture using $\theta$-stable parameters, 
 which will be introduced
 in Section \ref{subsec:thetapara}
 of Appendix I, 
 and restate Conjecture \ref{conj:GPver1} as an algorithm in this notation.

In Theorems \ref{thm:171471} and \ref{thm:171471b} of Appendix I, 
 we shall give the 
\index{B}{1thetastableparameter@$\theta$-stable parameter}
 $\theta $-stable parameters
 of the representations of the standard sequence starting
 with an irreducible finite-dimensional representation $F$
 summarized as follows.

\medskip
\noindent
\label{Hasse sequence} 
\begin{enumerate}
\item
Suppose that $n=2m$. 
Let 
\[
   F=
  \Kirredrep{O(2m+1,1)} {\mu}_{a,b}
  = \Kirredrep{O(2m+1,1)} {\mu} \otimes \chi_{a b}
\]
 for $\mu \in \Lambda^+(m+1)$ and $a,b \in \{\pm\}$
 be an irreducible finite-dimensional representation
 of $O(2m+1,1)$, 
 see Section \ref{subsec:fdimrep} in Appendix I.  
Its $\theta$-stable parameter is 
\[
(~||~\mu_1,\mu_2,\dots, \mu_{m},\mu_{m+1})_{a,b}
\]
 and we have the $\theta$-stable parameters of the representations in the  standard sequence (written in column).

\begin{eqnarray*}
\Pi_0(F)=&  (~||~\mu_1,\mu_2,\dots, \mu_{m},\mu_{m+1})_{a,b} & \\
\Pi_1(F)=& (\mu_1~|| ~\mu_2,\dots, \mu_{m}, \mu_{m+1})_{a,b} & \\
\vdots\qquad &  \vdots   & \\
\Pi_m(F)=&  (\mu_1,\mu_2,\dots, \mu_{m}~||~ \mu_{m+1})_{a,b}.  & 
\end{eqnarray*}
\item
Suppose that $n=2m+1$.  
Let 
\[
   F=\Kirredrep{O(2m+2,1)}{\mu}_{a,b} 
    = \Kirredrep{O(n+1,1)}{\mu} \otimes \chi_{ab} 
\]
   {for } $\mu \in \Lambda^+(m+1) \mbox{ and }
   a, b \in \{\pm\} $
 be an irreducible finite-dimensional representation of $O(2m+2,1)$.  
Its $\theta$-stable parameter is \[(||~\mu_1,\mu_2,\dots, \mu_{m},\mu_{m+1})_{a,b}\] 
and we have the $\theta$-stable parameters of the representations in standard sequence (written in column).

\begin{eqnarray*}
\Pi_0(F)&=  (||~\mu_1,\mu_2,\dots, \mu_{m+1})_{a,b} & 
\\
\Pi_1(F)&= (\mu_1~|| ~\mu_2,\dots, \mu_{m+1})_{a,b} & 
\\
\vdots\hphantom{mii} &  \vdots  &
\\
\Pi_{m+1}(F)&=  (\mu_1,\mu_2,\dots, \mu_{m+1}~|| )_{a,b}. & 
\end{eqnarray*}

\end{enumerate}

\medskip

We refer to the finite-dimensional representation
$\Pi_0(F)=F$
 as the {\it{starting representation}} of the standard sequence  
and to the tempered representation $\Pi_m(F)$
 (when $n=2m$)
 or the discrete series representation $\Pi_{m+1}(F)$
 (when $n=2m+1$)
 as the {\it{last representation}}
 of the standard sequence
 (see Remarks \ref{rem:Utemp} and \ref{rem:Umdisc}).

\begin{conjecture}
\label{conj:V2}
Let $\Kirredrep G \mu_{a,b}$ be an irreducible finite-dimensional representation of $G=O(n+1,1)$, 
 and $\Kirredrep {G'} \nu_{a,b}$ be an irreducible finite-dimensional representation of the subgroup $G'=O(n,1)$, 
 where $\mu \in \Lambda^+([\frac{n+2}{2}])$, 
 $\nu \in \Lambda^+([\frac{n+1}{2}])$, 
 and $a,b,c,d \in \{\pm\}$, 
 see \eqref{eqn:Fn1ab} and \eqref{eqn:ON1isom} in Appendix I.  
Assume that 
\begin{equation}
\label{eqn:fdabcd}
  {\operatorname{Hom}}_{G'} (\Kirredrep{G}{\mu}_{a,b}|_{G'}, \Kirredrep{G'}{\nu}_{c,d})\not = \{0\}.
\end{equation}

In (1) and (2) below, 
 nontrivial symmetry breaking operators are represented by arrows connecting the $\theta$-stable parameters of the representations.

\begin{enumerate}
\item
[{\rm{(1)}}]  
Suppose that $n=2m$.  
Then $\mu =(\mu_1, \cdots, \mu_{m+1})\in \Lambda^+(m+1)$
 and $\nu =(\nu_1, \cdots, \nu_{m})\in \Lambda^+(m)$.  
Then two  representations in the standard sequences have a nontrivial symmetry breaking operator if and only if the $\theta$-stable parameters of the representations satisfy one of the following conditions.

\begin{eqnarray*}
& (\mu_{1}, \dots , \mu_i ~|| ~\mu_{i+1} , \dots , \mu_{m+1})_{a,b}& \\
& \Downarrow &  \\
 &       (\nu_1 , \dots  ,\nu_i ~ ||~ \nu_{i+1}, \dots ,\nu_{m}  )_{c,d}
 \end{eqnarray*}
 \begin{center}
                      or
  \end{center}                    
\begin{eqnarray*}
& (\mu_{1}, \dots , \mu_i ~|| ~\mu_{i+1} , \dots , \mu_{m+1})_{a,b}&  
\\
& \Downarrow &  \\
&     (\nu_1, \dots  ,\nu_{i-1} ~ ||~ \nu_i, \nu_{i+1}, \dots ,\nu_{m})_{c,d} &
\end{eqnarray*}        

\item[{\rm{(2)}}]  
Suppose that $n=2m+1$.  
Then two infinite-dimensional representations in the standard sequences
 have a nontrivial symmetry breaking operator
 if and only if the $\theta$-stable parameters of the representations satisfy one of the following conditions:
\begin{eqnarray*} 
 & (\mu_{1}, \dots , \mu_i ~|| ~\mu_{i+1} , \dots , \mu_{m+1})_{a,b} &  \\  
 & \Downarrow &    \\
 &  (\nu_1, \dots  ,\nu_i ~ ||~ \nu_{i+1}, \dots ,\nu_{m+1}  )_{c,d} &
 \end{eqnarray*}
 \begin{center}
or 
\end{center}
\begin{eqnarray*}
& (\mu_{1}, \dots , \mu_i ~|| ~\mu_{i+1} , \dots , \mu_{m+1})_{a,b} &\\
& \Downarrow &  \\
 &  (\nu_1,  \dots  ,\nu_{i-1}  ~ ||~ \nu_{i}, \dots ,\nu_{m+1}  )_{c,d} &
 \end{eqnarray*}
\end{enumerate}  
\end{conjecture}

\begin{remark}
\label{rem:V2}
See Theorem \ref{thm:ON1branch} in Appendix I
 for the condition on the parameters $\mu$, $\nu$, and $a$, $b$, $c$, $d$
 such that \eqref{eqn:fdabcd} holds.  
In particular,
 \eqref{eqn:fdabcd} implies either 
 $(a,b)=(c,d)$ or $(a,b)=(-c,-d)$.  
See also Lemma \ref{lem:fdeq} (2) 
 for the description of overlaps in the expressions
 of irreducible finite-dimensional representations
 of $O(N-1,1)$ 
 when $N$ is even.  
\end{remark}
\subsection{Supporting evidence}
\label{subsec:SuppConj}

In this section,
 we provide some evidence 
 supporting our conjecture.  
\begin{description} 
\item[{\bf E.1}] \
If $F \in {\operatorname{Irr}}(G)_{\rho}$
 and $F' \in {\operatorname{Irr}}(G')_{\rho}$, 
 the Conjecture \ref{conj:GPver1} is true.  
(Equivalently,
 if $\Kirredrep{O(n+1,1)}\mu_{+,+}$ and $\Kirredrep{O(n,1)} {\nu}_{+,+}$ are both the trivial one-dimensional representations, 
Conjecture \ref{conj:V2} is true.)
\item[{\bf E.2}]  \ Some vanishing results  for symmetry breaking operators.  
\item[{\bf E.3}] \ Our conjecture is consistent
 with the Gross--Prasad conjecture for {\it{tempered}} representations of the special orthogonal group.
\item[{\bf E.4}]
  \  There exists a nontrivial symmetry breaking operator $\Pi_1 \rightarrow \pi_1$.  
\label{subsec:vanishtemp}
\end{description}

\bigskip
\noindent
\subsubsection{Evidence  E.1}
This was proved in Theorems \ref{thm:SBOvanish} and \ref{thm:SBOone}.  

\bigskip
\noindent
\subsubsection{Evidence E.2} 
Detailed proofs of the following propositions  will be published in a sequel to this monograph.

\medskip
Recall from Definition-Theorem \ref{def:UHasse}
 that $U_i(F_{a,b})$ refers to the $i$-th term in the Hasse sequence
 starting with the finite-dimensional representation 
 $F_{a,b}=F \otimes \chi_{a b}$ of $G$ 
 and $U_j(F'_{c,d})$ to the $j$-th term in the Hasse sequence starting
 with the finite-dimensional representation 
$
   F_{c,d}'=F' \otimes \chi_{c d}
$ of $G'$.  
\begin{proposition}
\label{conj:SBOvanish}
Let $a,b,c,d \in \{\pm\}$, 
 $0 \le i \le [\frac{n+1}{2}]$
 and $0 \le j \le [\frac{n}{2}]$.  
Then 
\[
{\operatorname{Hom}}_{G'}
 (U_i(F_{a,b})|_{G'}, U_j(F_{c,d}'))=\{ 0 \}
\quad
\text{if $j \not =i-1$, $i$}.  
\] 
\end{proposition}

\medskip 
If one of the representations of $G=O(n+1,1)$ respectively of $G'=O(n,1)$
 is tempered then the following vanishing theorems hold.

\vskip 1pc
\noindent
$\bullet$\enspace Assume first 
$(G,G')=(O(2m,1),O(2m-1,1))$.

Let $s=(s_0, \cdots, s_{m-1},0^{m+1}) \in \Lambda^+(2m+1)$
 and $t=(t_0, \cdots, t_{m-1},0^m) \in \Lambda^+(2m)$
 satisfy $t \prec s$
 (see Definition \ref{def:Young} for the notation).

\begin{proposition}
\label{prop:171427}
Let $U_0$, $\cdots$, $U_m$, $U_{m+1}$ be the Hasse sequence of $G=O(2m,1)$
 with $U_0= \Kirredrep{O(2m+1,{\mathbb{C}})}{s}$, 
 and $U_0', \cdots, U_{m-1}'$ be that of $G'=O(2m-1,1)$
 with $U_0'= \Kirredrep{O(2m,{\mathbb{C}})}{t}$.  
Then 
\[
{\operatorname{Hom}}_{G'}
(U_m|_{G'}, U_j')=\{0\}
\quad
\text{if $0 \le j \le m-2$.}
\]
\end{proposition}

\vskip 1pc
\noindent
$\bullet$\enspace Assume now
$(G,G')=(O(2m+1,1),O(2m,1))$.

Let $s=(s_0, \cdots, s_{m},0^{m+1}) \in \Lambda^+(2m+2)$
 and $t=(t_0, \cdots, t_{m-1},0^{m+1}) \in \Lambda^+(2m+1)$
 satisfy $t \prec s$.

\begin{proposition}
\label{prop:171429}
Let $U_0, \cdots, U_m$ be the Hasse sequence of $G=O(2m+1,1)$
 with $U_0= \Kirredrep{O(2m+2,{\mathbb{C}})}{s}$, 
 and $U_0', \cdots, U_{m}'$ be that of $G'=O(2m,1)$
 with $U_0'= \Kirredrep{O(2m+1,{\mathbb{C}})}{t}$.  
Then 
\[
{\operatorname{Hom}}_{G'}
(U_i|_{G'}, U_m')=\{0\}
\quad
\text{if \quad $0 \le i \le m-1$.}
\]
\end{proposition}

\medskip

\begin{remark}
These propositions prove only part of the vanishing statement
 of symmetry breaking operators
 formulated in Conjecture \ref{conj:V2}.
\end{remark}

\bigskip

\noindent 
\subsubsection{Evidence E.3}

We use the notations and assumptions
 of the previous section, 
 and show that our conjecture is consistent
 with the original 
\index{B}{GrossPrasadconjecture@Gross--Prasad conjecture|textbf}
Gross--Prasad conjecture
 on {\it{tempered representations}}
 \cite{GP}.  
For simplicity,
 we treat here 
 only for $(G, G')=(O(n+1,1), O(n,1))$
 with $n=2m$. 
We shall see that a special case
 of Conjecture \ref{conj:V2}
 ({\it{i.e.}}, the conjecture for the last representation
 of the standard sequence)
 implies some results
 (see \eqref{eqn:GPm} below)
 that were predicted by the original conjecture of Gross and Prasad
 for tempered representations of special orthogonal groups.

Assume that the irreducible finite-dimensional representations
 $\Pi_0$ of $G$ and $\pi_0$ of $G'$
 are of type I 
 (Definition \ref{def:typeone})
 and that 
 $(\mu_1,\dots, \mu_{m}, \mu_{m+1})$ and $(\nu_1,\dots, \nu_{m})$ are their highest weights.

By the branching law 
 for {\it{finite-dimensional}} representations
 with respect to $G \supset G'$
 (see Theorem \ref{thm:ON1branch} in Appendix I), 
 the condition 
 \[\mbox{Hom}_{O(n,1)}(\Pi_{0}|_{G'},\pi_{0}) \not = \{0\}\]
is equivalent to 
\begin{equation}
\label{eqn:munuGP}
   \mu_1 \geq \nu_1 \geq \mu_2 \geq \dots \geq \nu_{m} \geq \mu_{m+1} \geq 0.  
\end{equation}
Let $U_m$ (resp.~ $\Pi_m=U_m \otimes (\chi_{+-})^m$)
 be the $m$-th term of the Hasse sequence 
 (resp. the standard sequence)
 starting with the irreducible finite-dimensional representation
 $\Pi_0 =U_0$
 (see Definitions \ref{def:UHasse} and \ref{def:Hasse}).  
Then we have a direct sum decomposition 
 of the principal series representation
\[
   I_{(-1)^{m-\mu_{m+1}}}(\Kirredrep{O(2m)}{\mu_1+1, \cdots,\mu_m+1,0^m}, m-\mu_{m+1})
 \simeq U_m \oplus (U_m \otimes \det)
\]
by Theorem \ref{thm:171426} (1) and Remark \ref{rem:Utemp}.  
Assume that $\Pi_m$ is tempered.  
Then $U_m$ is also tempered,
 and the continuous parameter of the principal series representation
 must lie on the unitary axis, 
 that is, 
 $m-\mu_{m+1} \in m + \sqrt{-1}{\mathbb{R}}$.  
Hence $\mu_{m+1}=0$.

Since $\mu_{m+1}=0$, 
the $\theta$-stable parameters of the tempered representations
 $\Pi_m$, $\Pi_m \otimes \det$ are given by 
\[ 
   (\mu_1,\dots, \mu_m || 0  )_{+,+}, 
\quad
   \Rq{\mu_1, \cdots,\mu_m}{0}{-,-}, 
\]
 whereas the $\theta$-stable parameter 
 of the discrete series representation
 of $G'=O(2m,1)$ is given by 
\[
     (\nu_1, \dots , \nu_m ~|| )_{+,+}.  
\]
In view of the $K$-type formula in Theorem \ref{thm:171426} (2), 
 we see 
\[
   U_m \not \simeq U_m \otimes \det
\]
 as $G$-modules,
 and thus $\Pi_m \not \simeq \Pi_m \otimes \det$.  
Therefore,
 the restriction of the principal series representation $\Pi_m$
 of $G=O(2m+1,1)$ to the subgroup $\overline G=SO(2m+1,1)$ is irreducible
 by Lemma \ref{lem:171523} (1) in Appendix II.  
We set 
\[
  \overline \Pi_m := \Pi_m|_{\overline G}, 
\]
 which is an irreducible tempered representation of $\overline G$.

We now consider representations
 of the subgroups $G'=O(2m,1)$
 and $\overline{G'}=SO(2m,1)$.  
We observe that there is at most one
 discrete series representation of $\overline{G'}=SO(n,1)$
 for each infinitesimal character 
(see Proposition \ref{prop:disc} in Appendix I). 
Therefore the restriction of the discrete series representation $\pi_m$
 of $G'=O(2m,1)$ to the subgroup $\overline{G'}=SO(2m,1)$
 is irreducible,
 which is denoted by $\overline{\pi_m}$.

With these notations,
 Proposition \ref{prop:LNM26} in Appendix II yields
 a natural linear isomorphism:
\[
   {\operatorname{Hom}}_{G'}(\Pi_m|_{G'}, \pi_m)
   \oplus
   {\operatorname{Hom}}_{G'}((\Pi_m\otimes \det)|_{G'}, \pi_m)
   \simeq
   {\operatorname{Hom}}_{\overline{G'}}(\overline \Pi_m|_{\overline{G'}}, \overline{\pi_m}).  
\]
Conjecture \ref{conj:V2} for the pair $(G,G')=(O(n+1,1),O(n,1))$
 is applied to this specific situation;
 the first term in the left-hand side 
 equals
 ${\mathbb{C}}$
 and the second term vanishes.  
Thus Conjecture \ref{conj:V2} in this case implies the following statement
 for the pair $(\overline G,\overline{G'})=(SO(n+1,1),SO(n,1))$
 of special orthogonal groups:
\begin{equation}
\label{eqn:GPm}
   \text{
   ${\operatorname{Hom}}_{\overline{G'}}
   (\overline \Pi_m|_{\overline{G'}}, \overline{\pi_m})={\mathbb{C}}$
   \quad
   \text{if $\mu_{m+1}=0$ and \eqref{eqn:munuGP} is satisfied.}
}
\end{equation}

We now assume that the representation $\Pi_m$ is nontrivial on the center.
This determines the Langlands parameters of  the Vogan packets
 $VP(\overline{\Pi}_m)$ and $VP(\overline{\pi}_m)$ of $\overline G$
 respectively $\overline{G'}$, 
 and we follow exactly the steps of the algorithm by Gross and Prasad outlined
 in Chapter \ref{sec:Gross-Prasad}. 
We conclude again that  the Gross--Prasad conjecture predicts
 that $\{\overline \Pi_m, \overline \pi_m\}$ is the unique pair
 of representation in $VP(\overline \Pi_m)\times VP(\overline \pi_m)$
 with a nontrivial symmetry breaking operator.

\bigskip
\noindent
\subsubsection{Evidence E.4}
\label{subsec:Evidence4}
We will prove the existence of a nontrivial symmetry breaking operator 
\[
   \Pi _1 \rightarrow \pi_1.  
\]

We first introduce graphs to encode information
 about the images and kernels of symmetry breaking operators
 between reducible principal series representations
 as well as information about the images of the subrepresentation
 under the symmetry breaking operators. 
 This will be helpful to visualize the composition of an symmetry breaking operator with a Knapp--Stein operator.
 
\medskip
\noindent {\bf Admissible graphs} \\
\index{B}{admissiblegraph@admissible graph|textbf} 
Consider the vertices of a square. 
We call the following set of directed  graphs {\it{admissible}}:

\medskip

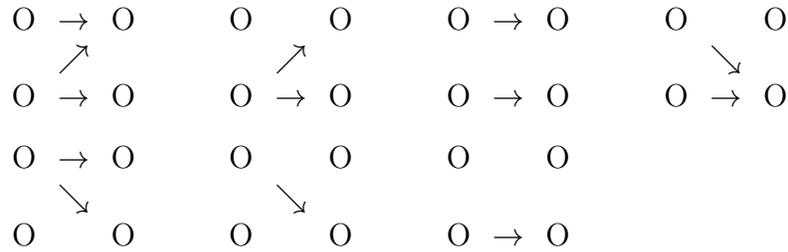
\begin{center}
\begin{tabular}{c@{\kern2em}c@{\kern2em}c@{\kern2em}c}
	\begin{tabular}{c@{\kern.7em}c@{\kern.7em}c} 
		\blackO   & $\rightarrow$  & \blackO \\
		&  $\nearrow$ & \\
		\blackO & $\rightarrow$ & \blackO
	\end{tabular} &
	\begin{tabular}{c@{\kern.7em}c@{\kern.7em}c} 
		\blackO  & &  \blackO \\
		& $\nearrow$ & \\
		\blackO & $\rightarrow$  & \blackO
	\end{tabular} &
	\begin{tabular}{c@{\kern.7em}c@{\kern.7em}c} 
		\blackO & $\rightarrow$  & \blackO \\
		&  &  \\
		\blackO & $\rightarrow$  & \blackO
	\end{tabular} &
	\begin{tabular}{c@{\kern.7em}c@{\kern.7em}c} 
		\blackO & &  \blackO \\
		&  $\searrow$ & \\
		\blackO & $\rightarrow$  & \blackO
	\end{tabular}
\\[2em]
	\begin{tabular}{c@{\kern.7em}c@{\kern.7em}c} 
		\blackO & $\rightarrow$  & \blackO \\
		&  $\searrow$ & \\
		\blackO & &  \blackO
	\end{tabular} &
	\begin{tabular}{c@{\kern.7em}c@{\kern.7em}c} 
		\blackO &   & \blackO \\
		&  $\searrow$  & \\
		\blackO & & \blackO
	\end{tabular} &
	\begin{tabular}{c@{\kern.7em}c@{\kern.7em}c} 
		\blackO &   & \blackO \\
		&  & \\
		\blackO & $\rightarrow$ & \blackO
	\end{tabular}
\end{tabular}
\end{center}

\noindent
and the zero graph without arrows:

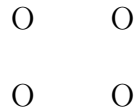
\begin{center}
	\begin{tabular}{c@{\kern.7em}c@{\kern.7em}c} 
		\blackO   &  & \blackO \\
             & \phantom{$\rightarrow$} & \\
            \blackO &  & \blackO
	\end{tabular}
\end{center}
Admissible graphs will encode
 information about the images
 and kernels of symmetry breaking operators.  
In the setting we shall use later,
 it is convenient to define the following equivalence relation 
among graphs,
 see Lemma \ref{lem:SBOadmgraph}.  
\begin{convention}
\label{conv:graph}
We identify two graphs ${\cal G}_1$ and ${\cal G}_2$
 if 
\[
   {\cal G}_1={\cal G}_2 \cup \{\ell\}
\]
 where $\ell$ is an arrow ending
 at the lower right vertex
 and ${\cal G}_2$ already contains an arrow 
 which starts from 
 the same vertex as $\ell$
 and which ends at the upper right vertex.  
\end{convention}
\begin{example}
\label{ex:samegraph}
The following graphs are pairwise equivalent.  
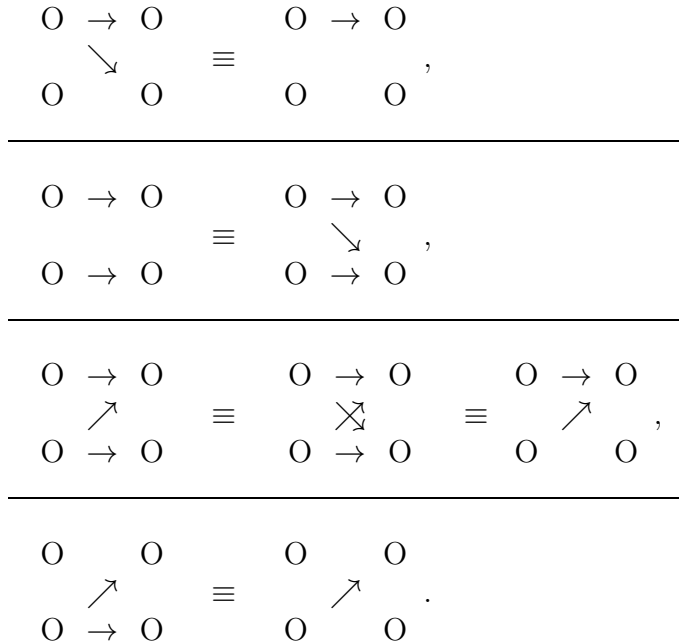
\begin{center}
\begin{tabular}{cccc}
	\begin{tabular}[m]{c@{\kern.7em}c@{\kern.7em}c}
		\blackO & $\to$ & \blackO \\
		 & $\searrow$ & \\ 
		\blackO &  & \blackO
	\end{tabular}& $\equiv$ & 
        \begin{tabular}[m]{c@{\kern.7em}c@{\kern.7em}c} 
		\blackO & $\to$ & \blackO \\
		~\\ 
		\blackO &  & \blackO
	\end{tabular}, & \\[2em]
\hline
\\
	\begin{tabular}[m]{c@{\kern.7em}c@{\kern.7em}c} 
		\blackO & $\to$ & \blackO \\
		~\\ 
		\blackO & $\to$ & \blackO
	\end{tabular} & $\equiv$ & 
	\begin{tabular}[m]{c@{\kern.7em}c@{\kern.7em}c}
		\blackO & $\to$ & \blackO \\
		 & $\searrow$ & \\ 
		\blackO & $\to$ & \blackO
	\end{tabular}, & \\[2em]
\hline
\\
	\begin{tabular}[m]{c@{\kern.7em}c@{\kern.7em}c}
		\blackO & $\to$ & \blackO \\
		 & $\nearrow$ & \\ 
		\blackO & $\to$ & \blackO
	\end{tabular}	& $\equiv$ & 
	\begin{tabular}[m]{c@{\kern.7em}c@{\kern.7em}c}
		\blackO & $\to$ & \blackO \\
		 & $\nearrow\kern-1em\searrow$ & \\ 
		\blackO & $\to$ & \blackO
	\end{tabular}
        & $\equiv$
          \begin{tabular}[m]{c@{\kern.7em}c@{\kern.7em}c}
		\blackO & $\to$ & \blackO \\
		 & $\nearrow$ & \\ 
		\blackO &  & \blackO
	\end{tabular},
\\[2em]
\hline
\\
        \begin{tabular}[m]{c@{\kern.7em}c@{\kern.7em}c}
		\blackO &  & \blackO \\
		 & $\nearrow$ & \\ 
		\blackO & $\to$ & \blackO
	\end{tabular}& $\equiv$ & 
        \begin{tabular}[m]{c@{\kern.7em}c@{\kern.7em}c}
		\blackO &  & \blackO \\
		 & $\nearrow$ & \\ 
		\blackO &  & \blackO
	\end{tabular}. &
\end{tabular}
\end{center}

\end{example}
We obtain a colored graph by coloring  the vertices
 of the graph by 4 different colors, each with a different color. We typically use the colors {\color{blue} blue } and {\color{red} red }
 for the vertices in the left column and {\color{green} green} and {\color{magenta} magenta } for the vertices in the right column.

\medskip
\noindent
{\bf Mutation of  admissible graphs} \\
We obtain a new colored graph ${\cal G}_2$  from a graph ${\cal G}_1$ 
 by 
\index{B}{mutation@mutation|textbf}
\lq\lq{mutation}\rq\rq.  
The rules of the mutation are given as follows.  
\begin{enumerate}
\item
[Rule 1.]  Consider the colored vertices on the right. 
Remove any arrow which ends at the lower right vertex. 
Interchange  the two colored  vertices on the right.  
The arrows which used to end at the upper right vertex now
 end at the lower right vertex.
\item
[Rule 2.] Consider the colored vertices on the left. 
Remove any arrow which starts at the upper left corner. 
Interchange the two colored vertices on the left. 
The arrows which used to start at the lower left vertex now start at the upper left vertex.
\item[Rule 3.] If the mutated graph ${\cal G}_2$ has no arrows, 
 {\it{i.e., }} ${\cal G}_2$ is the zero graph,  the mutation is not allowed.
\end{enumerate}

\medskip
We write {\bf R} for the mutation on the right column and {\bf L}  for the mutation on the left column.  
We sometimes refer to {\bf R} and {\bf L}
 as {\it{mutation rules}}.  

It is easy to see the following.  

\begin{lemma}
\begin{enumerate}
\item[{\rm{(1)}}] The mutated graph is again admissible.  
\item[{\rm{(2)}}] Mutation is well-defined for the equivalence relations
 given in Convention \ref{conv:graph}.  
\item[{\rm{(3)}}] Admissible graphs 
 for which no mutation is allowed
 do not have an arrow except for the one from the upper left vertex to the lower right vertex. 
\item[{\rm{(4)}}]  $\bf R \circ R $ and $\bf L \circ L$ are not allowed mutations. 
\item[{\rm{(5)}}]
 $\bf R \circ L= L \circ R$.
\end{enumerate}
\end{lemma}

\begin{definition}
[source and sink]
We call an admissible graph $\cal G$ a
\index{B}{source@source|textbf}
 {\it{source}} of a set of graphs
 if all other  graphs of the set are obtained through mutations of $\cal G$.
We call a graph $\cal G$ a
\index{B}{sink@sink|textbf}
 {\it{sink}} in a set of admissible graphs
 if neither $\bf R$ nor $\bf L$ is an allowed mutation of $\cal G$.
\end{definition}

Applying these rules,
 we obtain the following families of mutated graphs with one source. 
The source for the first,
 second, and third types is at the top right corner,
 applying {\bf R} changes the right column
 and applying {\bf L}  changes the left column.

\underline{First type} 

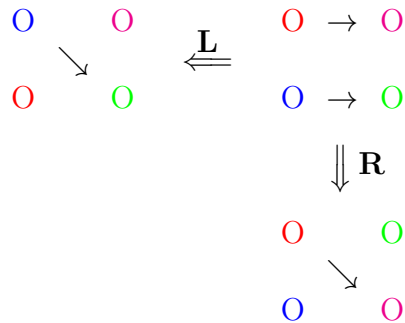
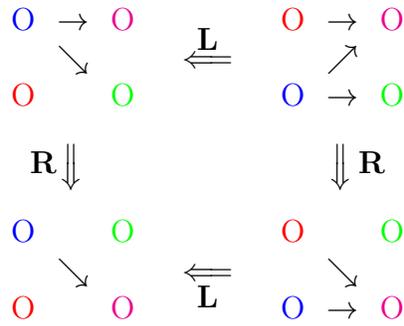
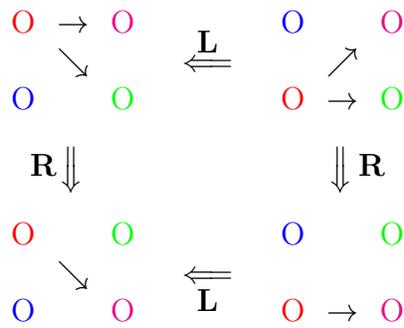
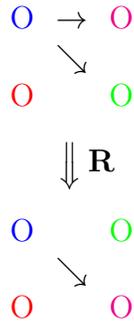
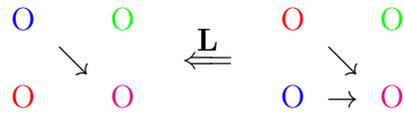
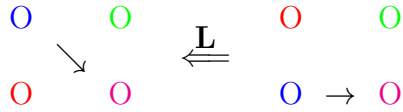
\def\Longdownarrow{\rotatebox[origin=c]{90}{$\Longleftarrow$}}
\def\Longuparrow{\rotatebox[origin=c]{90}{$\Longrightarrow$}}

\begin{center}
\begin{tabular}{ccc}
	\begin{tabular}{c@{\kern.7em}c@{\kern.7em}c} 
		\blueO &  & \magentaO \\
		 & $\searrow$  & \\
		 \redO & &  \greenO
	\end{tabular} & 	$\overset{\mbox{\textbf{L}}}{\Longleftarrow}$ &
	\begin{tabular}{c@{\kern.7em}c@{\kern.7em}c} 
		\redO  & $\rightarrow $ & \magentaO  \\
		 & & \\
		 \blueO & $\rightarrow$  & \greenO
	\end{tabular} \\[2em]
	& & \phantom{\textbf{R}}~\Longdownarrow~\textbf{R} \\[1em]
	& & 
	\begin{tabular}{c@{\kern.7em}c@{\kern.7em}c} 
		\redO &  & \greenO  \\
		& $\searrow$ & \\
		\blueO & & \magentaO
	\end{tabular}
\end{tabular}
\end{center}

\underline{Second type}

\begin{center}
\begin{tabular}{ccc}
	\begin{tabular}{c@{\kern.7em}c@{\kern.7em}c} 
		\blueO & $\rightarrow$ & \magentaO  \\
		& $\searrow$ & \\
		\redO & &  \greenO
	\end{tabular} 
	& $\overset{\mbox{\textbf{L}}}{\Longleftarrow}$ &
	\begin{tabular}{c@{\kern.7em}c@{\kern.7em}c} 
		\redO & $\rightarrow$ & \magentaO \\
		& $\nearrow$ & \\
		\blueO &  $\rightarrow$  & \greenO
	\end{tabular} \\[2em]
	\textbf{R}~\Longdownarrow~\phantom{\textbf{R}} & & \phantom{\textbf{R}}~\Longdownarrow~\textbf{R} \\[1em]
	\begin{tabular}{c@{\kern.7em}c@{\kern.7em}c} 
		\blueO & & \greenO \\
		&  $\searrow$  & \\
		\redO & & \magentaO
	\end{tabular} 
& $\underset{\mbox{\textbf{L}}}{\Longleftarrow}$ &
	\begin{tabular}{c@{\kern.7em}c@{\kern.7em}c} 
		\redO & & \greenO \\
		& $\searrow$ & \\
		\blueO & $\rightarrow$ & \magentaO
	\end{tabular} 
\end{tabular}
\end{center}

\underline{Third type}

\begin{center}
\begin{tabular}{ccc}
	\begin{tabular}{c@{\kern.7em}c@{\kern.7em}c} 
		\redO  & $\rightarrow$ & \magentaO  \\
		& $\searrow$ & \\
		\blueO & & \greenO
	\end{tabular}
	& $\overset{\mbox{\textbf{L}}}{\Longleftarrow}$ &
	\begin{tabular}{c@{\kern.7em}c@{\kern.7em}c} 
		\blueO  &  & \magentaO  \\
		& $\nearrow$ & \\
		\redO & $\rightarrow$ & \greenO 
	\end{tabular} 
\\[2em]
	\textbf{R}~\Longdownarrow~\phantom{\textbf{R}} & & \phantom{\textbf{R}}~\Longdownarrow~\textbf{R} \\[1em]	
\begin{tabular}{c@{\kern.7em}c@{\kern.7em}c} 
		\redO  &  & \greenO \\
		 &  $\searrow$  & \\
		 \blueO &  & \magentaO
	\end{tabular} 
& $\underset{\mbox{\textbf{L}}}{\Longleftarrow}$ &
	\begin{tabular}{c@{\kern.7em}c@{\kern.7em}c} 
		\blueO  &  &   \greenO \\
		&  & \\
		\redO & $\rightarrow$  & \magentaO
	\end{tabular}
\end{tabular}
\end{center}

\underline{Type A}

\begin{center}
\begin{tabular}{c}
	\begin{tabular}{c@{\kern.7em}c@{\kern.7em}c} 
		\blueO & $\rightarrow$ &  \magentaO  \\
		& $\searrow$ & \\
		\redO &  & \greenO
	\end{tabular}\\[2em]
 	\phantom{\textbf{R}}~\Longdownarrow~\textbf{R} \\[1em]
	\begin{tabular}{c@{\kern.7em}c@{\kern.7em}c} 
		\blueO  &  & \greenO \\
 		&  $\searrow$  & \\
 		\redO & &  \magentaO
 	\end{tabular} \\[2em]
\end{tabular}
\end{center}

\underline{Type B}

\begin{center}
\begin{tabular}{ccc}
	\begin{tabular}{c@{\kern.7em}c@{\kern.7em}c} 
		\blueO  &  &  \greenO  \\
		& $\searrow$  &    \\
		\redO &  & \magentaO
	\end{tabular}
& $\overset{\mbox{\textbf{L}}}{\Longleftarrow}$ &
	\begin{tabular}{c@{\kern.7em}c@{\kern.7em}c} 
		\redO  &  &\greenO  \\
		& $\searrow$ &  \\
		\blueO & $\rightarrow$ & \magentaO
	\end{tabular} 
\end{tabular}
\end{center}

\underline{Type C}

\begin{center}
\begin{tabular}{ccc}
	\begin{tabular}{c@{\kern.7em}c@{\kern.7em}c} 
		\blueO & & \greenO \\
		& $\searrow$ & \\
		\redO & & \magentaO     
	\end{tabular} 
& $\overset{\mbox{\textbf{L}}}{\Longleftarrow}$ &
	\begin{tabular}{c@{\kern.7em}c@{\kern.7em}c} 
		\redO &  & \greenO \\
		&   & \\
		\blueO & $\rightarrow$ & \magentaO
	\end{tabular}
\end{tabular}
\end{center}

\medskip

\noindent
This proves the following.  
 
\medskip

\begin{lemma}
\label{lem:graphs}
Let $\cal F$ be the family of admissible graphs
 that are obtained through mutations
 of a nonzero admissible graph.  
\begin{enumerate}
\item[{\rm{(1)}}] 
If $\cal F$ is not a singleton,
 it is one of the above six types.  
\item[{\rm{(2)}}] 
If $\cal F$ is a singleton,
 it is a coloring of the following graph. 

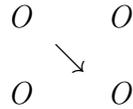
\begin{center}
	\begin{tabular}{c@{\kern.7em}c@{\kern.7em}c} 
		\blackO & & \blackO \\
		& $\searrow$ & \\
		\blackO & & \blackO
	\end{tabular}
\end{center}
\end{enumerate}
\end{lemma}

\medskip
{\bf From symmetry breaking operators to admissible graphs} \\

Assume that a principal series representation $I_{\delta}(V,\lambda)$ of $G$
 has exactly two composition factors $\Pi^1$ and $\Pi^2$, 
 which are not equivalent to each other.  
(The assumption is indeed satisfied for $G=O(n+1,1)$
 whenever $I_{\delta}(V,\lambda)$ is reducible.)
Thus there is an exact sequence of $G$-modules: 
\begin{equation}
  0  \to  {\color{red}\Pi^1}  \to  I_{\delta}(V,\lambda)
     \to {\color{blue}\Pi^2}  \to 0.  
\label{eqn:Inonsplit}
\end{equation}
Graphically,
 the irreducible inequivalent composition factors are represented by circles with different colors. The bottom circle represents the socle
 as follows.  

\begin{center}
	\begin{tabular}{c@{\kern2em}c@{\kern2em}c} 
		\blueO  & &\\
		& &\\
		\redO  & & 
	\end{tabular}
\end{center}
Later we shall assume in addition
 that the exact sequence \eqref{eqn:Inonsplit} does not split. 
(The assumption is satisfied
 if one of $\Pi^1$ or $\Pi^2$ is finite-dimensional.  
More generally,
 the assumption is indeed satisfied 
 for most of the pairs of the composition factors
 of the principal series representations of $G=O(n+1,1)$
 with regular integral infinitesimal characters, 
 see Theorem \ref{thm:LNM20} for example.)

\medskip
An analogous notation will be applied to principal series representations
 $J_{\varepsilon}(W,\nu)$
 of the subgroup $G'=O(n,1)$
 with two composition factors.  
Thus we represent the two composition factors
 of the reducible principal series representations
 $I_\delta(V, \lambda)$ and of $J_\varepsilon(W,\nu)$
 by four differently colored circles
 in a square; both the composition factors
 of a principal series representation are represented by circles vertically. 

We have the convention that the composition factors 
 of the representation $I_\delta(V,\lambda)$ of $G$
 are represented by the circles 
 on the left,
 those of $J_\varepsilon(W,\nu)$ of the subgroup on the right.  
Using this convention we get four squares with colored circles
 which  are obtained by changing the colors in each vertical column.

\bigskip

To a symmetry breaking operator 
\[{\mathbb B}^{V,W}_{\lambda,\nu}: I_\delta(V,\lambda) \rightarrow J_\varepsilon (W,\nu) \] 
we associate a graph
 which encodes information
 about the image and kernel of the symmetry breaking operator
 $\Bbb \lambda \nu {V,W}$
 as well as information about the image of the irreducible subrepresentation of the principal series representation
 $I_{\delta}(V,\lambda)$ of $G$
 under the symmetry breaking operator.  
We proceed as follows:
we obtain the arrows of the graph by considering the action of symmetry breaking operator ${\mathbb B}^{V,W}_{\lambda,\nu}$ on the composition factors. 
If no arrow starts at a circle, 
 then this means that the corresponding composition factor
 is in the kernel of the symmetry breaking operator. 
If no arrow ends at a circle,
 then this means
 that the $G'$-submodule
 of $J_{\varepsilon}(W,\nu)$ corresponding to the circle 
 is not in the image of the symmetry breaking operator. 
Then we have:
\begin{lemma}
\label{lem:SBOadmgraph}
Assume that both principal series representations
 $I_{\delta}(V,\lambda)$ and $J_{\varepsilon}(W,\nu)$
 have exactly two inequivalent composition factors
 with nontrivial extensions.  
Then with Convention \ref{conv:graph} the graph 
associated to our symmetry breaking operator
 $\Bbb \lambda \nu {V,W} \in {\operatorname{Hom}}_{G'}
(I_{\delta}(V,\lambda)|_{G'}, J_{\varepsilon}(W,\nu))$
 is an admissible graph.  
\end{lemma}
The proof of Lemma \ref{lem:SBOadmgraph} is straightforward.  
We illustrate it by examples as below.  
\begin{example}
[Graph of symmetry breaking operators]
\label{ex:SBOgraph}
{\rm{(1)}}\enspace
Suppose that the symmetry breaking operator is surjective and its restriction to the socle {\blueO} is also surjective. 
Then the associated graph is given by 
\begin{center}
\begin{tabular}{c@{\kern.7em}c@{\kern.7em}c} 
		\redO & $\rightarrow$ & \magentaO \\
		& $\nearrow\kern-1em\searrow$ & \\
		\blueO & $\rightarrow$ & \greenO
	\end{tabular}
\end{center}
by definition.  
With Convention \ref{conv:graph}, 
 we have 
\begin{center}
	\begin{tabular}{c@{\kern.7em}c@{\kern.7em}c} 
		\redO & $\rightarrow$ & \magentaO \\
		& $\nearrow$ & \\
		\blueO & $\rightarrow$ & \greenO
	\end{tabular}
$\equiv$
\begin{tabular}{c@{\kern.7em}c@{\kern.7em}c} 
		\redO & $\rightarrow$ & \magentaO \\
		& $\nearrow\kern-1em\searrow$ & \\
		\blueO & $\rightarrow$ & \greenO
	\end{tabular}, 
\end{center}
see Example \ref{ex:samegraph}.  
Then the graph in the left-hand side is admissible.  
\newline
{\rm{(2)}}\enspace
Suppose that the symmetry breaking operator is zero.  
Then it is depicted by the zero graph.  

\begin{center}
	\begin{tabular}{c@{\kern.7em}c@{\kern.7em}c} 
		\redO & & \magentaO \\
		& \phantom{$\rightarrow$} & \\
		\blueO & & \greenO
	\end{tabular}
\end{center}
\end{example}
To reduce the clutter in a digram representing a set of mutated graphs we often omit the zero graph, 
{\it{i.e.,}} the zero symmetry breaking operator. 

\medskip

We would like to encode information about a symmetry breaking operator
 and all its compositions with the Knapp--Stein operators at the same time. 
Composing symmetry breaking operators $\Bbb \lambda \nu {V,W}$
 with a Knapp--Stein intertwining operator 
\[
   \Ttbb \lambda {n-\lambda}V 
   \colon I_\delta(V,\lambda) \rightarrow I_\delta(V,n-\lambda)
\]
 for the group $G$ (see \eqref{eqn:KT}),  respectively 
\[
   \Ttbb \nu {n-1-\nu}W \colon J_\varepsilon(W,\nu) \rightarrow J_\varepsilon(W,n-1-\nu)
\]
 for the subgroup $G'$, 
 we obtain another symmetry breaking operator.  
If this new operator is not zero
 then it can be represented again by an admissible graph.~The graphs of these operators are arranged 
 compatible with our previous article \cite[Figs.~2.1--2.5]{sbon} where we draw $\nu$-value on the $x$-axis
 and the $\lambda$-value on the $y$-axis. 
We place the corresponding symmetry breaking operator in the corresponding quadrant.  
For example,
 if $\lambda \geq \frac{n}{2}$ and $\nu \geq \frac{n-1}{2}$, 
 then the parameters are arranged as 
\begin{alignat*}{2}
& ({n-1-\nu}, \lambda)
\qquad\qquad
&& (\nu,\lambda)
\\
&
&&
\\
& ({n-1-\nu},n-\lambda)
&& (\nu, {n-\lambda})
\end{alignat*}
 in the $(\nu,\lambda)$-plane, 
 and accordingly these symmetry breaking operators
 are arranged as follows.  

\begin{alignat*}{2}
& \Ttbb \nu {n-1-\nu}W \circ \Bbb \lambda \nu  {V,W}
&& \Bbb \lambda \nu {V,W}
\\
&
&&
\\
& \Ttbb \nu {n-1-\nu} W \circ \Bbb \lambda \nu {V,W} 
  \circ \Ttbb {n-\lambda} {\lambda}V
\qquad\qquad\qquad\quad\quad
&& \Bbb \lambda \nu {V,W} \circ \Ttbb {n-\lambda} {\lambda}V
\end{alignat*}

Accordingly,
 we shall consider four graphs
 of these four symmetry breaking operators.

By the definition of the mutation rule,
 we obtain:
\begin{lemma}
Assume that a principal series representation $I_{\delta}(V,\lambda)$ has 
 two irreducible composition factors $\Pi^1$ and $\Pi^2$
 with nonsplitting exact sequence
 \eqref{eqn:Inonsplit} 
 and that the Knapp--Stein operator
 $\Ttbb \lambda {n-\lambda} V \colon I_{\delta}(V,\lambda) \to 
 I_{\delta}(V,n-\lambda)$ is nonzero 
 but vanishes on the subrepresentation $\Pi^1$.  
Then the graph associated to 
 a symmetry breaking operator composed
 with $\Ttbb {n-\lambda}\lambda V$
 for the group $G$
 is obtained by using the mutation rule ${\mathbf {L}}$ for graphs.  
Similarly,
 the graph associated to a symmetry breaking operator
 composed with a nonzero Knapp--Stein operator
 $\Ttbb \nu{n-1-\nu}W \colon J_{\varepsilon}(W, \nu) 
  \to J_{\varepsilon}(W,n-1-\nu)$
 for the subgroup $G'$
 (with an analogous assumption on $J_{\varepsilon}(W,\nu)$)
 is obtained by using the mutation rule ${\mathbf{R}}$
 for graphs.
\end{lemma}

\begin{example}
In the Memoirs article \cite{sbon}
 we considered the case of two spherical principal series representations $I(\lambda)$ and $J(\nu)$ for integral parameters $i$, $j$.  
If 
\index{A}{Leven@$L_{\operatorname{even}}$}
$(-i,-j)\in L_{\operatorname{even}}$, 
 namely,
 if $i \ge j \ge 0$
 and $i \equiv j \mod 2$, 
 then the normalized regular symmetry breaking operator $I(-i) \rightarrow J(-j)$ is zero \cite[Thm.~8.1]{sbon}. 
The other symmetry breaking operators
 for spherical principal series representations
 with the same infinitesimal character are nonzero
 and we have functional equations with nonvanishing coefficients
 \cite[Thm.~8.5]{sbon}. 
Thus the family of mutated graphs associated to the regular symmetry breaking operators  is 
given as follows.

\medskip

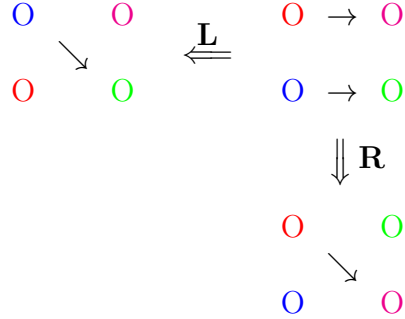
\begin{center}
\begin{tabular}{ccc}
	\begin{tabular}{c@{\kern.7em}c@{\kern.7em}c} 
		\blueO && \magentaO \\
		& $\searrow $ &\\
		\redO && \greenO
	\end{tabular} 
	& $\overset{\mbox{\textbf{L}}}{\Longleftarrow}$ & 
	\begin{tabular}{c@{\kern.7em}c@{\kern.7em}c}
		\redO  & $ \rightarrow $ &  \magentaO\\
		& \phantom{$\searrow$} & \\
		\blueO & $\rightarrow$ & \greenO
	\end{tabular} \\[2em]
	& & \phantom{\textbf{R}}~\Longdownarrow~\textbf{R}\\[1em]
	& & \begin{tabular}{c@{\kern.7em}c@{\kern.7em}c}
		\redO  &  & \greenO\\
		 & $\searrow$ & \\
		 \blueO &  & \magentaO
	\end{tabular}
\end{tabular}
\end{center}

We recall from \cite[Chap.~1]{sbon} 
 (or from Theorem \ref{thm:LNM20} in a more general setting)
 that both the $G$-module $I(-i)$ and the $G'$-module $J(-j)$ contain 
 irreducible finite-dimensional representations
 as their subrepresentations 
 ({\color{red}red} and {\color{magenta}magenta} circles)
 and irreducible infinite-dimensional representations
 $T(i)$ and $T(j)$
 ({\color{blue}blue} and {\color{green}green})
 as their quotients,
 respectively.  
The corresponding socle filtrations are given graphically
 as follows.
\[
  I(-i)=\begin{tabular}{c} 
		\blueO \\
		\\
		\redO
	\end{tabular} 
\qquad\qquad
  J(-j)=\begin{tabular}{c} 
		\greenO \\
		\\
		\magentaO
	\end{tabular}
\]
Note that, 
 under the assumption $i \ge j \ge 0$ and $i \equiv j \mod 2$, 
 we  have a nontrivial symmetry breaking operator
 between the two finite-dimensional representations 
 ({\color{red}red} and {\color{magenta}magenta} circles)
 and as well as between the nontrivial composition factors 
 $T(i) \rightarrow T(j)$ 
 ({\color{blue}blue} and {\color{green}green} circles), 
 see \cite[Thm~1.2 (1-a)]{sbon}. 
\end{example}

\begin{example}
\label{ex:317}
More generally in Corollary \ref{cor:IV1}
 we proved that 
\[
   \Atbb {\lambda_0}{\nu_0} {\gamma} {V,W}= 0 
\]
for negative integers $\lambda_0$, $\nu_0$ implies that 
\[
   \Atbb {n-\lambda_0}{n-1-\nu_0} {\gamma} {V,W} \not = 0.  
\]

Since $({n-1-\nu_0}, {n-\lambda_0}) \in {\mathbb{N}}^2$, 
 we may place the graph associated to the regular symmetry breaking operator 
$ 
   \Atbb {n-\lambda_0}{n-1-\nu_0} {\gamma} {V,W}
$
 in the NE corner
 according to the position 
 in the $(\nu,\lambda)$-plane
 as in \cite[Fig.~2.1, III.A or III.B]{sbon}.  

On the other hand,
 since $(\nu_0, \lambda_0) \in (-{\mathbb{N}})^2$, 
 we may place a zero graph associated to the zero operator
 $\Atbb {\lambda_0} {\nu_0} \gamma{V,W}$
 in the SW corner 
 according to the position
 in the $(\nu,\lambda)$-plane
 as in \cite[Fig.~2.1, I.A. or I.B.]{sbon}.  
\end{example}

\begin{example}
In the Memoirs article \cite[Thm.~11.1]{sbon}
 we prove that there is a differential symmetry breaking operator in the SW corner if the regular symmetry breaking operator is zero. 
To this operator and its composition with the Knapp--Stein operators the assigned graph is given as follows.  

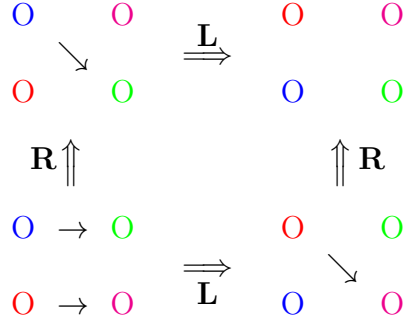
\begin{center}
\begin{tabular}{ccc}
	\begin{tabular}{c@{\kern.7em}c@{\kern.7em}c} 
		\blueO && \magentaO \\
		& $\searrow $ &\\
		\redO && \greenO
	\end{tabular} 
	& $\overset{\mbox{\textbf{L}}}{\Longrightarrow}$ & 
	\begin{tabular}{c@{\kern.7em}c@{\kern.7em}c}
		\redO  &  &  \magentaO\\
		& \phantom{$\searrow$} & \\
		\blueO &  & \greenO
	\end{tabular} \\[2em]
	\textbf{R}~\Longuparrow~\phantom{\textbf{R}} & & \phantom{\textbf{R}}~\Longuparrow~\textbf{R}\\[1em]
	\begin{tabular}{c@{\kern.7em}c@{\kern.7em}c}
		\blueO  &  $\rightarrow$ & \greenO\\
		 &  & \\
		 \redO &  $\rightarrow$ & \magentaO
	\end{tabular} 
	& $\underset{\mbox{\textbf{L}}}{\Longrightarrow}$ &
	\begin{tabular}{c@{\kern.7em}c@{\kern.7em}c}
		\redO  &  & \greenO\\
		 & $\searrow$ & \\
		 \blueO &  & \magentaO
	\end{tabular}
\end{tabular}
\end{center}
Note that the differential operator gives a source 
 in the mutation graphs in the SW corner in this setting.  
\end{example}

\bigskip

{\bf  Existence of a nontrivial symmetry breaking operators $\Pi_1 \rightarrow \pi_1$.} \\
\noindent
Recall that we assume that 
\[m(\Pi_0,\pi_0) =1\]
for the irreducible finite-dimensional representations
 $\Pi_0$ of $G$ and $\pi_0$ of the subgroup $G'$.
We consider now a pair of  reducible principal series representations  $I_\delta(V, \lambda)$ of $G$ and $J_\varepsilon (W,\nu)$ of $G'$
 with finite-dimensional composition factors $\Pi_0$, $\pi_0$, 
 respectively. 
\medskip
\begin{lemma}
 \label{lem:graphs}
Suppose that both \redO{} and \magentaO{} are representing irreducible finite-dimensional representations of $G$ and $G'$. 
We assume that \redO{} and \blueO{}
 respectively \magentaO{} and \greenO{} are representing the composition factors of a principal series representation of $G$, respectively $G'$. 
Then the following  graphs are  not  associated to a symmetry breaking operator.

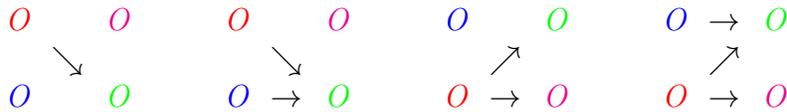
\begin{center}
\begin{tabular}{c@{\kern2em}c@{\kern2em}c@{\kern2em}c}
	\begin{tabular}{c@{\kern.7em}c@{\kern.7em}c}
		\redO  &  & \magentaO \\
		& $\searrow$ & \\
		\blueO &  & \greenO
	\end{tabular} &
	\begin{tabular}{c@{\kern.7em}c@{\kern.7em}c}
		\redO &  & \magentaO \\
		& $\searrow$  & \\
		\blueO & $\rightarrow$  & \greenO
	\end{tabular} &
	\begin{tabular}{c@{\kern.7em}c@{\kern.7em}c}
		\blueO & & \greenO \\
		& $\nearrow$ & \\
		\redO & $\rightarrow$ & \magentaO
	\end{tabular} &
	\begin{tabular}{c@{\kern.7em}c@{\kern.7em}c}
		\blueO & $\rightarrow$  & \greenO \\
		& $\nearrow$ & \\
		\redO & $\rightarrow$ & \magentaO
	\end{tabular}
\end{tabular}
\end{center}

\end{lemma}

\proof The representations \magentaO{} and \redO{} are
 finite-dimensional. 
The image of a finite-dimensional representation by a symmetry breaking operator is finite-dimensional.  \qed

\medskip

\begin{lemma}
\label{lem:fdgraph}
We keep Convention \ref{conv:graph} and the assumptions of Lemma \ref{lem:graphs}.  
\begin{enumerate}
\item[{\rm{(1)}}]
Suppose that \blueO{} and \greenO{} stand for 
 both irreducible subrepresentations of the principal series representations
 of $G$ and $G'$, respectively. 
The graph associated to a nontrivial symmetry breaking operator
 is one of the following.  

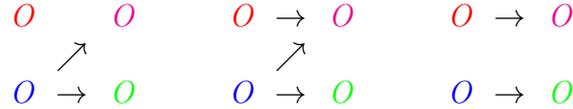
\begin{center}
\begin{tabular}{c@{\kern2em}c@{\kern2em}c@{\kern2em}c}
	\begin{tabular}{c@{\kern.7em}c@{\kern.7em}c}
		\redO &  & \magentaO \\
		& $\nearrow$ & \\
		\blueO & $\rightarrow$ & \greenO
	\end{tabular} & 
	\begin{tabular}{c@{\kern.7em}c@{\kern.7em}c}
		\redO & $\rightarrow$ &\magentaO \\
		& $\nearrow$ & \\
		\blueO & $\rightarrow$ & \greenO
	\end{tabular} & 
	\begin{tabular}{c@{\kern.7em}c@{\kern.7em}c}
		\redO & $\rightarrow$ & \magentaO \\
		& &\\
		\blueO & $\rightarrow$ & \greenO
	\end{tabular} 
\end{tabular}
\end{center}

\item[{\rm{(2)}}] 
Suppose that \redO{} and \magentaO{} stand
 for both irreducible finite-dimensional subrepresentations
 of the principal series representations. The   graph associated to a nontrivial symmetry breaking operator is one of the following.  

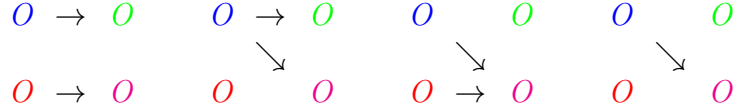
\begin{center}
\begin{tabular}{@{}c@{\kern1.4em}c@{\kern1.4em}c@{\kern1.4em}c@{\kern1.4em}c@{}}
	\begin{tabular}{c@{\kern.7em}c@{\kern.7em}c}
		\blueO & $\rightarrow$ & \greenO \\
		&  & \\
		\redO & $\rightarrow$ & \magentaO
	\end{tabular} &
	\begin{tabular}{c@{\kern.7em}c@{\kern.7em}c}
		\blueO & $\rightarrow$ & \greenO \\
		& $\searrow$ & \\
		\redO & & \magentaO
	\end{tabular} &
	\begin{tabular}{c@{\kern.7em}c@{\kern.7em}c}
		\blueO & & \greenO \\
		& $\searrow$ & \\
		\redO & $\rightarrow$ & \magentaO
	\end{tabular} &
	\begin{tabular}{c@{\kern.7em}c@{\kern.7em}c}
		\blueO & & \greenO \\
		& $\searrow$ &\\
		\redO & & \magentaO
	\end{tabular}
\end{tabular}
\end{center}
\end{enumerate}
\end{lemma}

\medskip

\medskip

Using the composition with the Knapp--Stein operators
 we obtain an action
 of the (little) Weyl group of $O(n+1,1) \times O(n,1)$
 on the continuous parameters of the symmetry breaking operators,
 hence on the symmetry breaking operators
 and also on their associated admissible graphs through the mutation rules.

\begin{example}
Let $\cal F$ be a family of mutated graphs
 such that the graph associated to the symmetry breaking operator 
 $\Atbb {n-\lambda_0}{n-1-\nu_0}{\gamma}{V,W}$
 is a source. 
If $\cal F$ is of first type,
 then the graph in the SE corner shows that there is a nontrivial symmetry breaking operator $\Pi_1 \rightarrow \pi_1$.
\end{example}
 
Using functional equations
 and the information about $(K,K')$-spectrum of regular symmetry breaking operators
 it is in some cases possible (see for example \cite{sbon})
 to show that the associated graph is of first type,
 but in general we do not have such explicit information about the regular symmetry breaking operators and so we have to proceed differently.

\bigskip

Suppose that $\Pi_0$ and $\pi_0$ are
 irreducible finite-dimensional subrepresentations
 of $I_\delta (V,\lambda)$ and $J_\varepsilon (W,\nu)$
 with ${\operatorname{Hom}}_{G'}(\Pi_0|_{G'}, \pi_0) \ne \{0\}$. 
By Proposition \ref{prop:172009}
 $[V:W] \not = 0$
 and
\index{A}{1psi@$\Psising$,
          special parameter in ${\mathbb{C}}^2 \times \{\pm\}^2$}
 $(\lambda, \nu, \delta, \varepsilon) \in \Psi_{\operatorname{sing}}$, 
 namely,
 the quadruple $(\lambda, \nu, \delta, \varepsilon)$
 does not satisfy the generic parameter condition  \eqref{eqn:nlgen}.
By Theorem \ref{thm:1716110} (see also Theorem \ref{thm:existDSBO} (1)), 
 there exists a nonzero differential symmetry breaking operator 
\[
   {\bf D}:I_\delta (V,\lambda) \rightarrow J_\varepsilon (W, \nu), 
\]
which we denote by $\bf D$. 
The image of $\bf D$ is infinite-dimensional
 by Theorem \ref{thm:imgDSBO}.  
Thus by Lemma \ref{lem:fdgraph} (2), 
 we obtain the following.

\begin{lemma}
The graph associated to $\bf D$ is one of the following.

\begin{center}
\begin{tabular}{c@{\kern2em}c@{\kern2em}c@{\kern2em}c}
	\begin{tabular}{c@{\kern.7em}c@{\kern.7em}c}
		\blueO & $\rightarrow$ & \greenO \\
		& & \\
		\redO & $\rightarrow$ & \magentaO
	\end{tabular} &
	\begin{tabular}{c@{\kern.7em}c@{\kern.7em}c}
		\blueO & $\rightarrow$ & \greenO \\
		& $\searrow$ &\\
		\redO & &\magentaO
	\end{tabular} &  &
\end{tabular}
\end{center}
\end{lemma}

\medskip
Mutating the graph of $\bf D$ by $\bf R$ we get the following.  

\begin{center}
	\begin{tabular}{c@{\kern.7em}c@{\kern.7em}c}
		\blueO &  & \magentaO \\
		& $\searrow$ & \\
		\redO & & \greenO
	\end{tabular}
\end{center}

Thus composing the differential symmetry breaking operator
 with a Knapp--Stein operator
 on the right
 we obtain a nontrivial symmetry breaking operator with this diagram
and thus a symmetry breaking operator  $U_1(F) \rightarrow  U_1(F')$.
We are ready to prove the following theorem,
 which gives evidence of our conjecture.  

\medskip

\begin{theorem}
Suppose that $F$ and $F'$ are irreducible finite-dimensional representations
 of $G$ and $G'$, 
 respectively.  
Let $\Pi_i$, $\pi_j$ be the standard sequences starting at $F$, $F'$, 
respectively.  
Then there exists a nontrivial symmetry breaking operator 
\[
   \Pi_1 \rightarrow \pi_1  
\]
 if ${\operatorname{Hom}}_{G'}(F|_{G'}, F')\ne \{0\}$.  
\end{theorem}
\begin{proof}
Recall from Definition \ref{def:Hasse}
 that  $\Pi_0= F$,  $\pi_0= F'$
 and $\Pi_1 = U_1(F) \otimes \chi_{+-}$,
 $\pi_1 = U_1(F') \otimes \chi_{+-}$
 and so 
\[ \mbox{Hom}_{G'}( \Pi_1|_{G'} , \pi_1)
\simeq \mbox{Hom}_{G'}(U_1(F)|_{G'}, U_1(F')).
\]
\end{proof}

\newpage
\section{Appendix I: Irreducible representations of $G=O(n+1,1)$, 
$\theta$-stable parameters,  and cohomological induction}
\label{sec:aq}
In Appendix I,
 we give a classification
 of irreducible admissible representations
 of $G=O(n+1,1)$
 in Theorem \ref{thm:irrG}.  
In particular, 
 we give a number of equivalent descriptions
 of irreducible representations
 with integral infinitesimal character
 (Definition \ref{def:intreg})
 by means of Langlands quotients
 (or subrepresentations), 
 coherent continuation starting at
\index{A}{1Piidelta@$\Pi_{i,\delta}$, irreducible representations of $G$}
 $\Pi_{i,\delta}$, 
 and cohomologically induced representations from 
 finite-dimensional representations of $\theta$-stable parabolic subalgebras,
 see Theorem \ref{thm:1808116}.  
Our results include a description
 of the following irreducible representations:
\begin{enumerate}
\item[$\bullet$]
\index{B}{Hassesequence@Hasse sequence}
\lq\lq{Hasse sequence}\rq\rq\ starting with arbitrary finite-dimensional 
irreducible representations
 (Theorems \ref{thm:171471} and \ref{thm:171471b});
\item[$\bullet$]
\index{B}{complementaryseries@complementary series representation}
complementary series representations
 with
\index{B}{singularintegralinfinitesimalcharacter@singular integral infinitesimal character}
 singular integral infinitesimal character
 (Theorem \ref{thm:compint}).  
\end{enumerate}
Since the Lorentz group $G=O(n+1,1)$ has four connected components, 
 we need a careful treatment
 even in dealing with finite-dimensional representations
 because not all of them extend holomophically 
 to $O(n+2,{\mathbb{C}})$.  
Thus Appendix I starts with irreducible finite-dimensional representations
 (Section \ref{subsec:fdimrep}), 
 and then discuss infinite-dimensional admissible representations
 for the rest of the chapter.

\subsection{Finite-dimensional representations of $O(N-1,1)$}
\label{subsec:fdimrep}
In this section we give a parametrization
 of irreducible finite-dimensional representations
 of the disconnected groups $O(N-1,1)$ and $O(N)$.   
The description here fits well with the $\theta$-stable parameters
 (Definition \ref{def:thetapara})
 for the Hasse sequence, 
 see Theorem \ref{thm:171471}.  
We note
 that the parametrization here 
 for irreducible finite-dimensional representations
 of $O(N)$ is different from what was defined in Section \ref{subsec:ONWeyl},  
 although the \lq\lq{dictionary}\rq\rq\
 is fairly simple,
 see Remark \ref{rem:FOn}.

There are two connected components in the compact Lie group $O(N)$.  
We recall from Definition \ref{def:type}
 that the set of equivalence classes
 of irreducible finite-dimensional representations
 of the orthogonal group $O(N)$
 can be divided into two types,
 namely,
 type I and II.  
On the other hand,
 there are four connected components in the noncompact Lie group $O(N-1,1)$, 
 and the division into two types 
 is not sufficient for the classification 
 of irreducible finite-dimensional representations of $O(N-1,1)$.  
We observe
 that some of the irreducible finite-dimensional representations
 of $O(N-1,1)$ cannot be extended to holomorphic representations
 of $O(N,{\mathbb{C}})$. 
For example,
 neither the one-dimensional representation $\chi_{+-}$ nor $\chi_{-+}$
 of $O(N-1,1)$ (see \eqref{eqn:chiab}) comes from a holomorphic character
 of $O(N,{\mathbb{C}})$.  
We shall use only representations
 of 
 \lq\lq{type I}\rq\rq\
 and tensoring them with four characters
 $\chi_{ab}$ ($a,b \in \{\pm\}$)
 to describe all irreducible finite-dimensional representations of $O(N-1,1)$.

First of all, 
 we recall from \eqref{eqn:Lambda}
 that 
\index{A}{0tLambda@$\Lambda^+(N)$, dominant weight}
$\Lambda^+(k)$ is the set of 
 $\lambda \in {\mathbb{Z}}^k$
 with $\lambda_1 \ge \lambda_2 \ge \cdots \ge \lambda_k
\ge 0$.

Let $N \ge 2$. 
For $\lambda \in \Lambda^+([\frac N 2])$, 
 we extend it to 
\begin{equation}
\label{eqn:lmdtilde}
   \widetilde \lambda 
 :=(\lambda_1, \cdots, \lambda_{[\frac N 2]}, 
    \underbrace{0,\cdots,0}_{[\frac{N+1}{2}]}) \in {\mathbb{Z}}^N, 
\end{equation}
 and define 
\index{A}{FONCl@$\Kirredrep{O(N,{\mathbb{C}})}{\lambda}$|textbf}
\begin{equation}
\label{eqn:Fholohw}
\Kirredrep {O(N,{\mathbb{C}})}{\lambda}_+
\equiv
\Kirredrep{O(N,{\mathbb{C}})}{\widetilde \lambda}, 
\end{equation}
to be the unique irreducible summand
 of $O(N,{\mathbb{C}})$ in the irreducible finite-dimensional representation 
 $\Kirredrep{G L(N,{\mathbb{C}})}{\widetilde \lambda}$
 of $GL(N,{\mathbb{C}})$
 that contains a highest weight vector
 corresponding to $\widetilde \lambda$, 
 see \eqref{eqn:CWOn}.  
Its restriction to the real forms $O(N)$ and $O(N-1,1)$
 will be denoted by $F^{O(N)}(\lambda)_+$ and $F^{O(N-1,1)}(\lambda)_{+,+}$, 
 respectively.  
Then the irreducible $O(N)$-module $\Kirredrep {O(N)}{\lambda}_+$ 
 is a representation of type I.  
We may summarize these notations
 as follows.  
\begin{equation}
\label{eqn:ONCreal}
   \Kirredrep {O(N)}{\lambda}_+
   \underset{{\operatorname{rest}}_{O(N)}}{\overset \sim \longleftarrow}
   \Kirredrep {O(N, {\mathbb{C}})}{\widetilde{\lambda}}
   \underset{{\operatorname{rest}}_{O(N-1,1)}}{\overset \sim \longrightarrow}
   \Kirredrep {O(N-1,1)}{\lambda}_{+,+}.  
\end{equation}
\begin{remark}
\label{rem:FOn}
With the notation as in \eqref{eqn:CWOn}, 
 we have
\[
   \Kirredrep {O(N)}{\lambda}_+
   \simeq
   \Kirredrep {O(N)}{\widetilde{\lambda}}
\]
for $\lambda \in \Lambda^+([\frac N 2])$.  
This is a general form 
 of representations of $O(N)$ of type I
 (Definition \ref{def:type}).  
Then other representations of $O(N)$, 
 {\it{i.e.}}, 
 representations of type II are obtained from the tensor product of those of type I
 with the one-dimensional representation, 
 $\det$, 
 as we recall now.  

Suppose $0 \le 2 \ell \le N$.  
If $\lambda \in \Lambda^+([\frac N 2])$ is of the form 
\[
   \lambda =(\lambda_1, \cdots, \lambda_{\ell}, 
             \underbrace{0, \cdots,0}_{[\frac N 2]-{\ell}})
\]
with $\lambda_{\ell} >0$, 
then by \eqref{eqn:type1to2}, 
 we have an isomorphism
 as representations of $O(N)$:
\[
   \Kirredrep {O(N)}{\lambda}_+ \otimes \det
   \simeq
   \Kirredrep {O(N)}{\underbrace{\lambda_1, \cdots,\lambda_{\ell}}_{\ell}, 
   \underbrace{1, \cdots,1}_{N-2{\ell}}, \underbrace{0, \cdots,0}_\ell}, 
\]
which is of type II
 if $N \ne 2\ell$.  
We shall denote this representation by $\Kirredrep {O(N)}{\lambda}_-$
 as \eqref{eqn:Fnpm} below.  
\end{remark}

Analogously, 
 $\Kirredrep {O(N-1,1)}{\lambda}_{+,+}$ is a general form
 of representations of the Lorentz group $O(N-1,1)$
 of type I in the following sense.  
\begin{definition}
[representation of type I for $O(N-1,1)$]
\label{def:typeone}
An irreducible
\newline
 finite-dimensional representation
 of $O(N-1,1)$ is said
 to be of {\it{type I}}
\index{B}{type1indef@type I, representation of $O(N-1,1)$|textbf}
 if it is obtained as the holomorphic continuation
of an irreducible representation
 of $O(N)$
 of type I
\index{B}{type1@type I, representation of $O(N)$}
 (see Definition \ref{def:type}).  
\end{definition}

We define for $\lambda \in \Lambda^+([\frac N2])$
\index{A}{FONpml@$\Kirredrep{O(N)}{\lambda}_{\pm}$|textbf}
\index{A}{FONpmpml@$\Kirredrep{O(N-1,1)}{\lambda}_{\pm\pm}$|textbf}
\begin{align}
  \Kirredrep {O(N)}{\lambda}_-
:=&
  \Kirredrep{O(N)}{\lambda}_+ \otimes \det, 
\label{eqn:Fnpm}
\\
  \Kirredrep {O(N-1,1)}{\lambda}_{a,b}
:=&
  \Kirredrep {O(N-1,1)}{\lambda}_{+,+} \otimes \chi_{ab}
\quad
(a,b \in \{\pm\}).  
\label{eqn:Fn1ab}
\end{align}
These are irreducible representations
 of $O(N)$ and $O(N-1,1)$, 
respectively.

With the notation \eqref{eqn:Fnpm} and \eqref{eqn:Fn1ab}, 
 irreducible finite-dimensional representations
 of $O(N)$ and of $O(N-1,1)$, 
 respectively,
 are described as follows:
\begin{lemma}
\label{lem:161612}
\begin{enumerate}
\item[{\rm{(1)}}]
Any irreducible finite-dimensional representation 
 of $O(N)$
 is of the form 
$
    \text{$F^{O(N)}(\lambda)_+$ or $F^{O(N)}(\lambda)_-$ for some $\lambda \in \Lambda^+([\frac N 2]).$}
$
\item[{\rm{(2)}}]
Suppose $N \ge 3$.  
Any irreducible finite-dimensional representation of $O(N-1,1)$
 is of the form
 $\Kirredrep {O(N-1,1)}{\lambda}_{a,b}$ 
for some $\lambda \in \Lambda^+([\frac N 2])$
 and $a, b \in \{\pm\}$.  
\index{A}{1chipmpm@$\chi_{\pm\pm}$, one-dimensional representation of $O(n+1,1)$}
\end{enumerate}
\end{lemma}

The point of Lemma \ref{lem:161612} (2) is 
 that an analogous statement of Weyl's unitary trick may fail
 for the disconnected group $O(N-1,1)$, 
 that is, 
 not all irreducible finite-dimensional representations
 of $O(N-1,1)$ cannot extend to holomorphic representations
 of $O(N, {\mathbb{C}})$.  
\begin{proof}
[Proof of Lemma \ref{lem:161612}]
(1)\enspace
This is a restatement of Weyl's description \eqref{eqn:CWOn}
 of $\widehat{O(N)}$.  
\par\noindent
(2)\enspace
Take any irreducible finite-dimensional representation
 $\sigma$ of $O(N-1,1)$.  
By the Frobenius reciprocity,
 $\sigma$ occurs as an irreducible summand of the induced representation 
 ${\operatorname{Ind}}_{SO_0(N-1,1)}^{O(N-1,1)}(\sigma|_{SO_0(N-1,1)})$.  
Since $N \ge 3$, 
 the fundamental group of $SO(N,{\mathbb{C}})/SO_0(N-1,1)$ is trivial
 because it is homotopic to 
 $SO(N)/SO(N-1) \simeq S^{N-1}$, 
 see \cite[Lem.~6.1]{xkpro89}.  
Hence the irreducible finite-dimensional representation $\tau$
 of $SO_0(N-1,1)$ extends to a holomorphic representation 
 of $SO(N,{\mathbb{C}})$, 
 which we shall denote by $\tau_{\mathbb{C}}$.

Let $\lambda \in \Lambda^+([\frac N 2])$ be the highest weight
 of the irreducible $SO(N,{\mathbb{C}})$-module $\tau_{\mathbb{C}}$.  
Then $\tau_{\mathbb{C}}$ occurs in the restriction 
 $\Kirredrep {O(N,{\mathbb{C}})} {\widetilde \lambda}|_{SO(N,{\mathbb{C}})}$, 
 and therefore the $SO_0(N-1,1)$-module $\tau$ occurs 
 in the restriction 
 ${\Kirredrep{O(N-1,1)}{\lambda}}_{+,+}|_{SO(N-1,1)}$.  
Hence $\sigma$ occurs as an irreducible summand of the induced representation
\begin{equation}
\label{eqn:IndSO0}
{\operatorname{Ind}}_{SO_0(N-1,1)}^{O(N-1,1)}
({\Kirredrep{O(N-1,1)}{\lambda}}_{+,+}|_{SO_0(N-1,1)}).  
\end{equation}
In light that ${\Kirredrep{O(N-1,1)}{\lambda}}_{+,+}$ is a representation
 of $O(N-1,1)$, 
 we can compute the induced representation \eqref{eqn:IndSO0}
 as follows.  
\begin{align*}
\text{\eqref{eqn:IndSO0}}
\simeq\,\, & {\Kirredrep{O(N-1,1)}{\lambda}}_{+,+} 
         \otimes {\operatorname{Ind}}_{SO_0(N-1,1)}^{O(N-1,1)}({\bf{1}})
\\
\simeq\,\, &{\Kirredrep{O(N-1,1)}{\lambda}}_{+,+} 
         \otimes (\bigoplus_{a, b \in \{\pm\}} \chi_{a b})
\\
\simeq\,\, & \bigoplus_{a,b \in \{\pm\}}{\Kirredrep{O(N-1,1)}{\lambda}}_{a,b}.  
\end{align*}
Thus Lemma \ref{lem:161612} is proved.  
\end{proof}
There are a few overlaps in the expressions \eqref{eqn:Fnpm} for 
$O(N)$-modules and \eqref{eqn:Fn1ab} for $O(N-1,1)$-modules.  
We give a necessary and sufficient condition
 for two expressions, 
 which give the same irreducible representation as follows.  
\begin{lemma}
\label{lem:fdeq}
\begin{enumerate}
\item[{\rm{(1)}}]
The following two conditions on $\lambda, \mu \in \Lambda^+([\frac N 2])$
 and $a, b \in \{\pm\}$ are equivalent:
\begin{enumerate}
\item[{\rm{(i)}}]
$\Kirredrep {O(N)}{\lambda}_a \simeq \Kirredrep {O(N)}{\mu}_b$
 as $O(N)$-modules;
\item[{\rm{(ii)}}]
\lq\lq{$\lambda = \mu$ and $a=b$}\rq\rq\
 or the following condition holds:
\begin{equation}
\label{eqn:ONisom}
\text{
$\lambda=\mu$, 
 $N$ is even, $\lambda_{\frac N 2} >0$, 
 and $a=-b$. 
}
\end{equation}
\end{enumerate}
\item[{\rm{(2)}}]
Suppose $N \ge 2$.  
Then the following two conditions on $\lambda, \mu \in \Lambda^+([\frac N 2])$
 and $a, b, c, d \in \{\pm\}$ are equivalent:
\begin{enumerate}
\item[{\rm{(i)}}]
$\Kirredrep {O(N-1,1)}{\lambda}_{a,b} \simeq \Kirredrep {O(N-1,1)}{\mu}_{c,d}$
 as $O(N-1,1)$-modules;
\item[{\rm{(ii)}}]
\lq\lq{$\lambda = \mu$ and $(a,b)=(c,d)$}\rq\rq\
 or the following condition holds:
\begin{equation}
\label{eqn:ON1isom}
\text{
$\lambda=\mu$, 
 $N$ is even, $\lambda_{\frac N 2}>0$, and $(a,b)=-(c,d)$. 
}
\end{equation}
\end{enumerate}
\end{enumerate}
\end{lemma}

\begin{proof}
(1)\enspace
The $O(N)$-isomorphism $\Kirredrep {O(N)}{\lambda}_{a} \simeq \Kirredrep {O(N)}{\mu}_{b}$ implies an obvious isomorphism
 $\Kirredrep {O(N)}{\lambda}_{a}|_{SO(N)} \simeq \Kirredrep {O(N)}{\mu}_{b}|_{SO(N)}$
 as $SO(N)$-modules,
 whence $\lambda=\mu$ by the classical branching law
 (Lemma \ref{lem:OSO})
 for the restriction $O(N) \downarrow SO(N)$.  
Then the equivalence (i) $\Leftrightarrow$ (ii) follows from the equivalence 
(i) $\Leftrightarrow$ (iii) in Lemma \ref{lem:branchII}.  
\par\noindent
(2)\enspace
Similarly to the proof for the first statement,
 we may and do assume $\lambda=\mu$
 by considering of the restriction $O(N-1,1) \downarrow SO(N-1,1)$.  
Then the proof of the equivalence (i) $\Leftrightarrow$ (ii) for $O(N-1,1)$
 reduces to the case for $O(N,1)$
 and the following lemma.  
\end{proof}

\begin{lemma}
\label{lem:holoON1}
Suppose $\sigma$ is an irreducible finite-dimensional representation of 
 $O(N-1,1)$.  
\begin{enumerate}
\item[{\rm{(1)}}]
Suppose $N \ge 2$.  
If $\sigma$ is extended to a holomorphic representation
 of $O(N,{\mathbb{C}})$, 
 then neither $\sigma \otimes \chi_{+-}$ nor $\sigma \otimes \chi_{-+}$
 can be extended to a holomorphic representation of $O(N,{\mathbb{C}})$.  
\item[{\rm{(2)}}]
Suppose $N \ge 3$.  
If $\sigma$ cannot be extended to a holomorphic representation
 of $O(N,{\mathbb{C}})$, 
 then both $\sigma \otimes \chi_{+-}$ and $\sigma \otimes \chi_{-+}$
 can be extended to a holomorphic representation of $O(N,{\mathbb{C}})$.  
\end{enumerate}
\end{lemma}

\begin{proof}
(1)\enspace
If $\sigma \otimes \chi_{ab}$ extends to a holomorphic representation
 of $O(N,{\mathbb{C}})$, 
 then so does the subrepresentation $\chi_{ab}$ in the tensor product
 $(\sigma \otimes \chi_{ab}) \otimes \sigma^{\vee}$, 
 where $\sigma^{\vee}$ stands for the contragredient representation of $\sigma$.  
Since $\chi_{ab}$ is the restriction of some holomorphic character of $O(N,{\mathbb{C}})$
 if and only if $(a,b)=(+,+)$ or $(-,-)$, 
 the first statement is proved.  
\par\noindent
(2)\enspace
As in the proof of Lemma \ref{lem:161612} (2), 
 we see that at least one element
 in $\{\sigma \otimes \chi_{ab}:a,b \in \{\pm\}\}$
 can be extended to a holomorphic representation
 of $O(N,{\mathbb{C}})$.  
Then the second statement follows from the first one.  
\end{proof}

\begin{example}
\label{ex:exteriorpm}
The natural action of $O(N)$
 on $i$-th exterior algebra
 $\Exterior^i({\mathbb{C}}^N)$ is given as 
\[
\Exterior^i({\mathbb{C}}^N)
\simeq
\begin{cases}
\Kirredrep {O(N)}{1^i, 0^{[\frac N 2]-i}}_{+}
\quad
&\text{if } i \le \frac N 2, 
\\
\Kirredrep {O(N)}{1^{N-i}, 0^{i-[\frac {N+1} 2]}}_{-}
\quad
&\text{if } i \ge \frac N 2, 
\end{cases}
\]
with the notation
 in this section,
 whereas the same representation was described as 
\[
    \Exterior^i({\mathbb{C}}^N)
    \simeq
    \Kirredrep {O(N)}{\underbrace{1,\cdots,1}_{i}, \underbrace{0,\cdots,0}_{N-i}}
\]
with the notation \eqref{eqn:CWOn}
 in Section \ref{subsec:repON}.  
\end{example}
As in the classical branching rule for $O(N) \downarrow O(N-1)$
 given in Fact \ref{fact:ONbranch}, 
 we give the irreducible decomposition of finite-dimensional representations
 of $O(N,1)$
 when restricted to the subgroup $O(N-1,1)$
 as follows:
\begin{theorem}
[branching rule for $O(N,1) \downarrow O(N-1,1)$]
\label{thm:ON1branch}
\index{B}{branchingruleON1@branching rule, for $O(N,1) \downarrow O(N-1,1)$}
Let $N \ge 2$.  
Suppose that $(\lambda_1, \cdots,\lambda_{[\frac{N+1}2]})
 \in \Lambda^+([\frac{N+1}2])$
 and $a,b \in \{\pm\}$.  
Then the irreducible finite-dimensional representation
 $\Kirredrep{O(N,1)}{\lambda_1, \cdots,\lambda_{[\frac{N+1}2]}}_{a, b}$
  of $O(N,1)$ decomposes
 into a multiplicity-free sum
 of irreducible representations
 of $O(N-1,1)$ as follows:
\begin{equation*}
   \Kirredrep{O(N,1)}{\lambda_1, \cdots,\lambda_{[\frac{N+1}2]}}_{a, b}|_{O(N-1,1)}
   \simeq
  \bigoplus
   \Kirredrep{O(N-1,1)}{\nu_1, \cdots, \nu_{[\frac{N}{2}]}}_{a, b}, 
\end{equation*}
where the summation is taken over
 $(\nu_1, \cdots, \nu_{[\frac{N}{2}]}) \in {\mathbb{Z}}^{[\frac N 2]}$
 subject to 
\begin{alignat*}{2}
&\lambda_{1} \ge \nu_{1} \ge \lambda_{2} \ge \cdots \ge \nu_{\frac N 2} \ge 0
\quad
&&\text{for $N$ even}, 
\\
&\lambda_{1} \ge \nu_{1} \ge \lambda_{2} \ge \cdots \ge \nu_{\frac {N-1} 2}
\ge \lambda_{\frac {N+1} 2}
\quad
&&\text{for $N$ odd}.   
\end{alignat*} 
\end{theorem}

\begin{proof}
The assertion follows 
 in the case $(a,b)=(+,+)$ from Fact \ref{fact:ONbranch}.  
The general case follows from the definition \eqref{eqn:Fn1ab}
 and from the observation 
 that the restriction $\chi_{a,b}|_{G'}$
 of the $G$-character $\chi_{a,b}$
 gives the same type of a character for $G'=O(N-1,1)$, 
 see \eqref{eqn:chiabrest}. 
\end{proof}

\subsection{Singular parameters for $V \in \widehat{O(n)}$: $S(V)$ and $S_Y(V)$}

In this section we prepare some notation 
 that describes the parameters
 of {\it{reducible}} principal series representations $I_{\delta}(V,\lambda)$
 of $G=O(n+1,1)$.

We recall from Lemma \ref{lem:IVchi}
 that both of the following subsets
\begin{align*}
&\{(\delta,V,\lambda)
:
\text{$I_{\delta}(V,\lambda)$ has regular integral ${\mathfrak{Z}}_G({\mathfrak{g}})$-infinitesimal character}
\}, 
\\
&\{(\delta,V,\lambda)
:
\text{$I_{\delta}(V,\lambda)$ is reducible}
\}
\end{align*}
 of $\{\pm\} \times \widehat {O(n)} \times {\mathbb{C}}$
 are preserved under the following transforms:
\begin{align*}
  (\delta,V,\lambda) &\mapsto (-\delta,V,\lambda), 
\\
    (\delta,V,\lambda) &\mapsto (\delta,V \otimes \det,\lambda).  
\end{align*}
Thus we omit 
 the signature $\delta$
 in our notation,
 and focus on the second and third components.  

\begin{definition}
\label{def:RIntRed}
We define two subsets
 of $\widehat {O(n)} \times {\mathbb{C}}$
 (actually, of $\widehat{O(n)} \times {\mathbb{Z}}$)
 by 
\index{A}{RzqSeducible@$\Reducible$ $(\subset \widehat {O(n)} \times {\mathbb{Z}})$|textbf}
\index{A}{RzqSeducibleIn@$\RInt$ $(\subset \widehat {O(n)} \times {\mathbb{Z}})$|textbf}
\begin{align}
\RInt
&:=\{(V,\lambda)
:
\text{$I_{\delta}(V,\lambda)$ has regular integral ${\mathfrak{Z}}_G({\mathfrak{g}})$-infinitesimal character}\}, 
\notag
\\
\label{eqn:reducible}
\Reducible
&:=\{(V,\lambda)
:
\text{$I_{\delta}(V,\lambda)$ is reducible}\}.  
\end{align}
\end{definition}

Both the  sets $\RInt$ and $\Reducible$ are preserved by the transformations
\begin{alignat*}{3}
&(V, \lambda) 
&&\mapsto 
&&(V \otimes \det, \lambda), 
\\
&(V, \lambda) 
&&\mapsto 
&&(V, n-\lambda).  
\end{alignat*}
This is clear for $\RInt$, 
 whereas the assertions for $\Reducible$ follows from
 the $G$-isomorphism $I_{\delta}(V,\lambda) \otimes \chi_{--} \simeq I_{\delta}(V\otimes \det,\lambda)$
 by Lemma \ref{lem:IVchi}
 and from the fact that $I_{\delta}(V,n-\lambda)$
 is isomorphic to the contragredient representation
 of $I_{\delta}(V,\lambda)$.  
We shall introduce two discrete sets 
 $S(V)$ and $S_Y(V)$ for $V \in \widehat {O(n)}$
 in Definition \ref{def:SVSYV} below, 
 and prove in Lemma \ref{lem:regint} and in Theorem \ref{thm:irrIV}
\begin{alignat*}{3}
  &\RInt
  &&=
  &&\{(V,\lambda) \in \widehat {O(n)} \times {\mathbb{Z}}
    :
    \lambda \not \in S(V) \}
\\
  &\hphantom{ii} \cup
  &&
  &&\hphantom{MMMMM}\cup
\\
  &\Reducible
  &&=
  &&\{(V,\lambda) \in \widehat {O(n)} \times {\mathbb{Z}}
    :
    \lambda \not \in S(V) \cup S_Y(V)\},  
\end{alignat*}
 see also Convention \ref{conv:SYsigma}.

\subsubsection{Infinitesimal character $r(V,\lambda)$ of $I_{\delta}(V,\lambda)$}
\label{subsec:muV}
Suppose that $V \in \widehat {O(n)}$ is given as 
\[
   V=\Kirredrep{O(n)}{\sigma}_{\varepsilon}
  \quad
  \text{for some }
  \sigma \in \Lambda^+\left([\dfrac n 2]\right)
  \text{ and }
  \varepsilon \in \{\pm\}
\]
with the notation
 as in Section \ref{subsec:fdimrep}.  
We define an element of ${\mathfrak {h}}_{\mathbb{C}}^{\ast}
\simeq {\mathbb{C}}^{[\frac n2]+1}$
by 
\index{A}{rVlmd@$r(V,\lambda)$|textbf}
\begin{equation}
\label{eqn:IVZG}
r(V,\lambda)
:=
(\sigma_{1}+\frac n 2-1,\sigma_{2}+\frac n 2-2, \cdots,\sigma_{[\frac n 2]}+
\frac n 2 - [\frac n 2],\lambda-\frac n 2).  
\end{equation}
The ordering in \eqref{eqn:IVZG} will play a crucial role
 in a combinatorial argument
 in later sections,
 whereas,
 up to the action of the Weyl group $W_G$, 
 $r(V,\lambda)$ gives the 
\index{A}{ZGg@${\mathfrak{Z}}_G({\mathfrak{g}})$}
\index{B}{infinitesimalcharacter@infinitesimal character}
\index{B}{ZGginfinitesimalcharacter@${\mathfrak{Z}}_G({\mathfrak{g}})$-infinitesimal character}
 ${\mathfrak{Z}}_G({\mathfrak{g}})$-infinitesimal character
 of the unnormalized induced representation
 $I_{\delta}(V, \lambda)$ of $G=O(n+1,1)$, 
 see \eqref{eqn:ZGinfI}.

\begin{example}
\label{ex:rhoi}
For $0 \le i \le n$, 
 we set $\ell := \min(i,n-i)$
 and 
\index{A}{1parhoi@$\rho^{(i)}$|textbf}
\begin{align*}
   \rho^{(i)}:=& r(\Exterior^i({\mathbb{C}}^n),i)  
\\
=&(\underbrace{\dfrac n 2, \dfrac n 2-1, \cdots, 
 \dfrac n 2-\ell+1}_\ell, 
 \underbrace{\dfrac n 2-\ell-1, \cdots, 
 \dfrac n 2-[\dfrac n 2]}_{[\frac n 2]-\ell},
 i-\frac n 2)
\\
=&\begin{cases}
(\underbrace{\dfrac n 2, \cdots, 
 \dfrac n 2-i+1}_i, 
 \underbrace{\dfrac n 2-i-1, \cdots, 
 \dfrac n 2-[\dfrac n 2]}_{[\frac n 2]-i},
 i-\frac n 2)
\quad
&\text{for $i \le [\dfrac n 2]$, }
\\[1em]
(\underbrace{\dfrac n 2, \cdots, 
 -\dfrac n 2+i+1}_{n-i}, 
 \underbrace{-\dfrac n 2+i-1, \cdots, 
 \dfrac n 2-[\dfrac n 2]}_{i-[\frac{n+1}2]},
 i-\dfrac n 2)
\quad
&\text{for $[\dfrac {n+1} 2] \le i$. }
\end{cases}
\end{align*}
Here are some elementary properties.  
\begin{enumerate}
\item[{\rm{(1)}}]
The following equations hold:
\begin{align}
\label{eqn:rhoi0}
\rho^{(i)}-\rho^{(0)}
=&
(\underbrace{1,\cdots,1}_{\ell},
 \underbrace{0,\cdots,0}_{[\frac n 2]-\ell},i)
\\
\notag
=&
\begin{cases}
(\underbrace{1,\cdots,1}_{i},
 \underbrace{0,\cdots,0}_{[\frac n 2]-i},
 i)
\quad
\text{for $0 \le i \le [\frac n 2]$, }
\\
(\underbrace{1,\cdots,1}_{n-i},
 \underbrace{0,\cdots,0}_{i-[\frac{n+1}2]},
 i)
\quad
\text{for $[\frac{n+1}2] \le i \le n$.  }
\end{cases}
\end{align}

\item[{\rm{(2)}}]
Let $r(V,\lambda)$ be defined as in \eqref{eqn:IVZG}.  
Then for any $i$ $(0 \le i \le n)$, 
 we have
\begin{align*}
r(V,\lambda)
 =&(\sigma_1, \cdots, \sigma_{[\frac n 2]}, \lambda) + \rho^{(0)}
\\
 =&(\sigma_1-1, \cdots, \sigma_{\ell}-1, \sigma_{\ell+1}, \cdots, \sigma_{[\frac n 2]}, \lambda-i) + \rho^{(i)}, 
\end{align*}
where we retain the notation $\ell = \min (i,n-i)$.  

\item[{\rm{(3)}}]
For all $i$ $(0 \le i \le n)$, 
\index{A}{1parho@$\rho_G$}
\begin{equation}
\label{eqn:rhoiW}
   \rho_G \equiv \rho^{(i)} \mod W_G.  
\end{equation}
\end{enumerate}
\end{example}

\subsubsection{Singular integral parameter:
$S(V)$ and $S_Y(V)$}
Retain the setting as in Section \ref{subsec:muV}.  
Let $G=O(n+1,1)$ and $m=[\frac n 2]$.  
Suppose $V \in \widehat{O(n)}$
 is given as $V = \Kirredrep{O(n)}{\sigma}_{\varepsilon}$
 with $\sigma=(\sigma_1, \cdots, \sigma_m)$ and $\varepsilon \in \{\pm\}$.  
Since $\sigma_1$, $\cdots$, $\sigma_{m} \in {\mathbb{Z}}$,  
 the following three conditions
 on $\lambda \in {\mathbb{C}}$ are equivalent:
\begin{enumerate}
\item[(i)]
The ${\mathfrak{Z}}_G({\mathfrak{g}})$-infinitesimal character
 of $I_{\delta}(V,\lambda)$ is integral 
 in the sense of Definition \ref{def:intreg};
\item[(ii)]
$\langle r(V,\lambda), \alpha^{\vee} \rangle \in {\mathbb{Z}}$ for any $\alpha \in \Delta({\mathfrak{g}}_{\mathbb{C}}, {\mathfrak{h}}_{\mathbb{C}})$;
\item[(iii)]
 $\lambda \in {\mathbb{Z}}$.  
\end{enumerate}

For each $V \in \widehat {O(n)}$, 
 we introduce a subset $S(V)$ in ${\mathbb{Z}}$
 (and a subset $S_Y(V)$ in ${\mathbb{Z}}$ for $V$ of type Y) 
 as follows.  
\begin{definition}
[$S(V)$ and $S_Y(V)$]
\label{def:SVSYV}
Let $m=[\frac n 2]$.  
For $V=\Kirredrep{O(n)}{\sigma}_{\varepsilon}$
 with $\sigma=(\sigma_1, \cdots, \sigma_m) \in \Lambda^+(m)$
 and $\varepsilon \in \{\pm\}$, 
 we define a finite subset of ${\mathbb{Z}}$ by 
\index{A}{Ssigma@$S(V)$|textbf}
\begin{equation}
\label{eqn:singint}
  S(V):=\{j-\sigma_j, n+\sigma_j-j: 1 \le j \le m\}.  
\end{equation}

When the irreducible $O(n)$-module $V$
 is of 
\index{B}{typeY@type Y, representation of ${O(N)}$\quad}
 type Y
 (see Definition \ref{def:OSO}), 
 namely, 
when $n$ is even $(=2m)$ and $\sigma_m > 0$, 
 we define also the following finite set
\index{A}{SsigmaY@$S_Y(V)$|textbf}
\begin{equation}
\label{eqn:SYsigma}
  S_Y(V):=\{\lambda \in {\mathbb{Z}} : 0 < |\lambda-m| <\sigma_m\}.  
\end{equation}
\end{definition}
We note that 
\[
  S(V) \cap S_Y(V) = \emptyset
\]
by definition.  
We shall sometimes adopt the following convention:
\begin{convention}
\label{conv:SYsigma}
When $V$ is of type X (see Definition \ref{def:OSO}), 
 we set
\[
   S_Y(V) = \emptyset.  
\]
\end{convention}

The definitions imply the following lemma.  
\begin{lemma}
\label{lem:regint}
The ${\mathfrak {Z}}_G({\mathfrak{g}})$-infinitesimal character
 of $I_{\delta}(V,\lambda)$ is regular integral 
 (see Definition \ref{def:intreg})
 if and only if $\lambda \in {\mathbb{Z}} \setminus S(V)$.  
Thus, 
 we have
\[
\RInt =\{(V,\lambda) \in \widehat{O(n)} \times {\mathbb{Z}}
:
 \lambda \not \in S(V) \}.  
\]
\end{lemma}
We refer to $S(V)$ 
as the set of 
 {\it{singular integral parameters}}.  
It should be noted
 that $I_{\delta}(V,\lambda)$ has
\index{B}{regularintegralinfinitesimalcharacter@regular integral infinitesimal character}
 regular integral infinitesimal character
 if $\lambda \in S_Y(V)$, 
 since $S_Y(V) \subset {\mathbb{Z}} \setminus S(V)$.

We shall see in Theorem \ref{thm:irrIV} below
 that the principal series representation $I_{\delta}(V, \lambda)$ is
 irreducible
 if and only if 
 $\lambda \in ({\mathbb{C}} \setminus {\mathbb{Z}}) \cup S(V) \cup S_Y(V)$.

We end this section with a lemma
 that will be used in Appendix III
 (Chapter \ref{sec:Translation})
 when we discuss translation functors.  

\begin{lemma}
\label{lem:1806115}
Let $V \in \widehat{O(n)}$
 and $\lambda \in {\mathbb{Z}} \setminus S(V)$.  
\begin{enumerate}
\item[{\rm{(1)}}]
Suppose $V$ is of type X (Definition \ref{def:OSO}).  
Then the $W_{\mathfrak{g}}$- and $W_G$-orbits
 through
 $r(V,\lambda) \in {\mathfrak{h}}_{\mathbb{C}}^{\ast} \simeq {\mathbb{C}}^{[\frac n 2]+1}$
 coincide:
\begin{equation}
\label{eqn:WGrV}
W_{\mathfrak{g}}\, r(V,\lambda)
=
W_G\, r(V,\lambda).  
\end{equation}
\item[{\rm{(2)}}]
Suppose $V$ is of type Y.  
Then \eqref{eqn:WGrV} holds
 if and only if $\lambda=\frac n 2$.  
\end{enumerate}
\end{lemma}

\begin{proof}
(1)\enspace 
The assertion is obvious
 when $n$ is odd
 because $W_{\mathfrak{g}}=W_G$
 in this case.  
Suppose $n$ is even, 
say, 
 $n=2m$.  
It is sufficient to show
 that $r(V,\lambda)$ contains zero
 in its entries.  
Since $V$ is of type X, 
 we have $\sigma_m=0$, 
 and therefore, 
 the $m$-th entry of $r(V,\lambda)$ amounts to 
 $\sigma_m + m-m=0$ by the definition \eqref{eqn:IVZG}.  
Thus the lemma is proved.  
\newline
(2)\enspace
Since $V$ is of type Y, 
 $n$ is even $(=2m)$ and $W_G \supsetneqq W_{\mathfrak{g}}$.  
Since $\lambda \not \in S(V)$, 
 $r(V,\lambda)$ is $W_{\mathfrak{g}}$-regular.  
Hence \eqref{eqn:WGrV} holds
 if and only if at least one of the entries in $r(V,\lambda)$ equals zero.  
Since $\sigma_1 \ge \sigma_2 \ge \cdots \ge \sigma_m >0$, 
 this happens only when the $(m+1)$-th entry
 of $r(V,\lambda)$ vanishes, 
 {\it{i.e.}}, $\lambda = \frac n 2 (=m)$.  
Hence Lemma \ref{lem:1806115} is proved.  
\end{proof}

\begin{remark}
\label{rem:zeroreg}
For $n=2m$ (even), 
 if $V$ is of type X
 or if $\lambda=m$, 
 then the ${\mathfrak {Z}}_G({\mathfrak{g}})$-infinitesimal character
 \eqref{eqn:IVZG} is regular 
 for 
\index{A}{Weylgroupg@$W_{\mathfrak{g}}$, Weyl group for ${\mathfrak{g}}_{\mathbb{C}}={\mathfrak{o}}(n+2,{\mathbb{C}})$}
 $W_{\mathfrak{g}}$
 in the sense of Definition \ref{def:intreg}, 
 but is 
\index{B}{singularintegralinfinitesimalcharacter@singular integral infinitesimal character}
\lq\lq{singular}\rq\rq\
 with respect to the Weyl group 
\index{A}{WeylgroupG@$W_G$, Weyl group for $G=O(n+1,1)$}
 $W_G$
 for the {\it{disconnected}} group 
 $G=O(n+1,1)$
 which is not in the Harish-Chandra class.  
\end{remark}

\subsection{Irreducibility condition of $I_{\delta}(V,\lambda)$}
\label{subsec:180595}
We are ready to state a necessary and sufficient condition
 for the principal series representation $I_{\delta}(V,\lambda)$
 of $G=O(n+1,1)$
 to be irreducible.  

We recall from \eqref{eqn:singint} and \eqref{eqn:SYsigma}
 the definitions of 
\index{A}{Ssigma@$S(V)$}
$S(V)$ and 
\index{A}{SsigmaY@$S_Y(V)$}
$S_Y(V)$, 
 respectively.  

\begin{theorem}
[irreducibility criterion of $I_{\delta}(V,\lambda)$]
\label{thm:irrIV}
Let $G=O(n+1,1)$, $\delta \in \{\pm\}$, 
 $V \in \widehat{O(n)}$,
 and $\lambda \in {\mathbb{C}}$.  
\begin{enumerate}
\item[{\rm{(1)}}]
If $\lambda \in {\mathbb{C}} \setminus {\mathbb{Z}}$, 
 then the principal series representation $I_{\delta} (V,\lambda)$ of $G$ 
 is irreducible.  
\item[{\rm{(2)}}]
Suppose $\lambda \in {\mathbb{Z}}$.  
Then $I_{\delta}(V,\lambda)$ is irreducible if and only if
\begin{alignat*}{2}
&\text{$\lambda \in S(V)$}
&&\text{when $V$ is of type X}, 
\\
&\text{$\lambda \in S(V) \cup S_Y(V)$}
\quad
&&\text{when $V$ is of type Y}.  
\end{alignat*}
\end{enumerate}
Thus $\Reducible$
 $($see \eqref{eqn:reducible}$)$
 is given by 
\begin{equation}
\label{eqn:Red}
  \Reducible
 =\{(V, \lambda) \in \widehat{O(n)} \times {\mathbb{Z}}
  : \lambda \not \in S(V) \cup S_Y(V)\}
\end{equation}
with Convention \ref{conv:SYsigma}.  
\end{theorem}
The proof of Theorem \ref{thm:irrIV} will be given 
in Section \ref{subsec:pfirrIV} in Appendix II
 by inspecting the restriction of $I_{\delta}(V,\lambda)$
 of $G=O(n+1,1)$
 to its subgroups $\overline G=SO(n+1,1)$ and $G_0=SO_0(n+1,1)$.  

\begin{example}
\label{ex:irrIilmd}
Let $0 \le i \le n$.  
The exterior tensor representation
 on $\Exterior^i({\mathbb{C}}^n)$ is of type X
 if and only if $n \ne 2i$
 (see Example \ref{ex:2.1}).  
A simple computation shows
\begin{alignat*}{2}
S(\Exterior^{(i)}({\mathbb{C}}^n))=&{\mathbb{Z}} 
   \setminus (\{i,n-i\} \cup (-{\mathbb{N}}_+) \cup (n+{\mathbb{N}}_+))
\qquad
&&\text{for $0 \le i \le n$}, 
\\
S_Y(\Exterior^{(m)}({\mathbb{C}}^n))=& \emptyset
\qquad
&&\text{for $n=2m$}, 
\end{alignat*}
 see also Example \ref{ex:isig}.  
Hence $I_{\delta}(i,\lambda)$ is reducible 
 if and only if
\[
   \lambda \in \{i,n-i\} \cup (-{\mathbb{N}}_+) \cup (n+{\mathbb{N}}_+)
\]
by Theorem \ref{thm:irrIV}.  
See Theorem \ref{thm:LNM20}
 for the socle filtration of $I_{\delta}(i,\lambda)$
 for $\lambda =i$ or $n-i$.  
\end{example}

For later purpose,
 we decompose $\Reducible$
 into two disjoint subsets as follows:
\begin{definition}
\label{def:Red12}
 We recall from Definition \ref{def:OSO}
 that any $V \in \widehat{O(n)}$ is either
 of type X or of type Y for
 $\widehat{O(n)}$.  
We set
\index{A}{RzqSeducibleI@$\RedI$|textbf}
\index{A}{RzqSeducibleII@$\RedJ$|textbf}
\begin{align*}
\RedI
:=&\{(V, \lambda) \in \Reducible
:
\text{$V$ is of type X or $\lambda = \frac n 2$}
\}, 
\\
\RedJ
:=&\{(V, \lambda) \in \Reducible
:
\text{$V$ is of type Y and $\lambda \ne \frac n 2$}
\}.  
\end{align*}
Then we have a disjoint union 
\[
\Reducible=\RedI \amalg \RedJ.  
\]
\end{definition}

\begin{remark}
\label{rem:Red12}
If $n$ is odd,
 then 
\[
  \RedJ = \emptyset
  \quad
  \text{and}
  \quad
  \Reducible =\RedI.  
\]
\end{remark}

\subsection{Subquotients of $I_{\delta}(V,\lambda)$
\label{subsec:subrep}}
By Theorem \ref{thm:irrIV}, 
 the principal series representation $I_{\delta}(V,\lambda)$
 of $G=O(n+1,1)$ is reducible
 {\it{i.e.}}, $(V,\lambda) \in \Reducible$
 if and only if
\begin{equation*}
   \lambda \in {\mathbb{Z}} \setminus (S(V) \cup S_Y(V))
\end{equation*}
 with Convention \ref{conv:SYsigma}.  
In this section,
 we explain the socle filtration of $I_{\delta}(V,\lambda)$.  
A number of different characterizations of the subquotients 
 will be given in later sections,
 see Theorem \ref{thm:1808116}
 for summary.  
We divide the arguments into the following two cases:
\par\noindent
Case 1.\enspace
$\lambda \ne \frac n 2$, 
 see Section \ref{subsec:Isub1};
\par\noindent
Case 2.\enspace
$\lambda = \frac n 2$, 
 see Section \ref{subsec:Isub3}.  

\subsubsection{Subquotients of 
 $I_{\delta}(V,\lambda)$ for $V$ of type X
\label{subsec:Isub1}}
\index{B}{typeX@type X, representation of ${O(N)}$\quad}

We begin with the case
 where $\lambda \ne \frac n 2$.  
This means that we treat the following cases:
\begin{enumerate}
\item[$\bullet$]
$V$ is of type X, 
 and $\lambda \in {\mathbb{Z}} \setminus S(V)$;
\item[$\bullet$]
$V$ is of type Y, 
 and $\lambda \in {\mathbb{Z}} \setminus (S(V) \cup S_Y(V)\cup \{\frac n 2\})$. \end{enumerate}

\begin{proposition}
\label{prop:Xred}
Let $G=O(n+1,1)$, 
 $V \in \widehat{O(n)}$, 
$\delta \in \{\pm\}$, 
 and $\lambda \in {\mathbb{Z}} \setminus (S(V) \cup S_Y(V))$.  
Assume further that $\lambda \ne \frac n 2$.  
Then there exists a unique proper submodule
 of the principal series representation
 $I_{\delta}(V,\lambda)$, 
 to be denoted by 
\index{A}{IdeltaVf@$I_{\delta}(V, \lambda)^{\flat}$|textbf}
$I_{\delta}(V,\lambda)^{\flat}$.  
In particular, 
 the quotient $G$-module
\index{A}{IdeltaVs@$I_{\delta}(V, \lambda)^{\sharp}$|textbf}
\[
  I_{\delta}(V,\lambda)^{\sharp}
  :=
  I_{\delta}(V,\lambda)
  /I_{\delta}(V,\lambda)^{\flat}
\]
is irreducible.  
\end{proposition}

The proof of Proposition \ref{prop:Xred} will be given 
 in Section \ref{subsec:Isub2} of  Appendix II.  

\begin{remark}
The $K$-type formul{\ae}
 and the minimal $K$-types
 of the irreducible $G$-modules $I_{\delta}(V,\lambda)^{\flat}$
 and $I_{\delta}(V,\lambda)^{\sharp}$
 will be given
 in Proposition \ref{prop:KXred}
 and Proposition \ref{prop:minKXred}, 
 respectively.  
\end{remark}

\subsubsection
{Subrepresentations of $I_{\delta}(V,\frac n 2)$
 for $V$ of type Y}
\label{subsec:Isub3}
\index{B}{typeY@type Y, representation of ${O(N)}$\quad}

Next we discuss the case:
\begin{enumerate}
\item[$\bullet$]
$V$ is of type Y and $\lambda = \frac n 2$.  
\end{enumerate}

In this case $I_{\delta}(V,\lambda)$ is the smooth representation
 of a tempered unitary representation.  

\begin{proposition}
[reducible tempered principal series representation]
\label{prop:IVmtemp}
Let $G=O(n+1,1)$ with $n=2m$, 
 $V\in \widehat{O(n)}$ be of type Y, 
 and $\delta \in \{\pm\}$.  
Then the principal series representation
 $I_{\delta}(V,m)$ of $G$
 is decomposed
 into the direct sum
 of two irreducible representations of $G$,
 to be written as:
\[
  I_{\delta}(V,m)
  \simeq
  I_{\delta}(V,m)^{\flat} 
  \oplus
  I_{\delta}(V,m)^{\sharp}.  
\]
If we express $V=\Kirredrep {O(n)}{\sigma}_{\varepsilon}$
 by $\sigma=(\sigma_1, \cdots, \sigma_m) \in \Lambda^+(m)$
 with $\sigma_m >0$
 and $\varepsilon \in \{\pm\}$, 
 then the irreducible $G$-modules
 $I_{\delta}(V,m)^{\flat}$ and $I_{\delta}(V,m)^{\sharp}$
 are characterized by their minimal $K$-types
 given respectively by the following:
\begin{align*}
 & \Kirredrep{O(n+1)}{\sigma_1, \cdots, \sigma_m}_{\varepsilon} \boxtimes \delta, 
\\
 & \Kirredrep{O(n+1)}{\sigma_1, \cdots, \sigma_m}_{-\varepsilon} \boxtimes (-\delta).  
\end{align*}
\end{proposition}
\begin{proof}
This is proved in Proposition \ref{prop:180572} (2)
 except for the assertion on the $K$-types.  
The last assertion on the minimal $K$-types follow from 
 the $K$-type formula
 of $I_{\delta}(V,m)^{\flat}$ and $I_{\delta}(V,m)^{\sharp}$
 in Proposition \ref{prop:KXred} (2).  
\end{proof}

\subsubsection{Socle filtration of $I_{\delta}(V,\lambda)$}
By Theorem \ref{thm:irrIV}
 together with Propositions \ref{prop:Xred} and \ref{prop:IVmtemp}, 
 we obtain the following:
\begin{corollary}
\label{cor:length2}
Let $G=O(n+1,1)$ for $n \ge 2$.  
Then the principal series representation
 $I_{\delta}(V,\lambda)$
 $(\delta \in \{\pm\}$, $V \in \widehat{O(n)},$ $\lambda \in {\mathbb{C}})$
 of $G$ 
 is either irreducible or of composition 
 series of length two.  
\end{corollary}

\subsection{Definition of the height $i(V,\lambda)$}
\label{subesec:defiVlmd}

In this section we introduce the 
\index{B}{height@height, of $(V,\lambda)$|textbf}
\lq\lq{height}\rq\rq
\[
  i \colon \RInt \to \{0,1,\ldots,n\}, 
\quad
  (V, \lambda) \mapsto i(V, \lambda)
\]
which plays an important role
 in the study of the principal series representation
 $I_{\delta}(V,\lambda)$ of $G$.  
For instance, 
 we shall see in Section \ref{subsec:HasseV}
 that the $K$-type formula
 for subquotients
 of $I_{\delta}(V,\lambda)$ 
 is described by using the height $i(V,\lambda)$
 when $(V,\lambda) \in \Reducible$
 (Definition \ref{def:RIntRed}).  
Moreover,
 we shall prove in Theorem \ref{thm:1807113}
 that the $G$-module $I_{\delta}(V,\lambda)$
 is obtained by the translation functor
 applied to the principal series representation $I_{\pm}(i, i)$
 with the trivial infinitesimal character $\rho_G$
 without \lq\lq{crossing the wall}\rq\rq\
 if we take $i$ to be the height $i(V,\lambda)$, 
 see Theorem \ref{thm:1807113}.  
We note that the group $G=O(n+1,1)$
 is not in the Harish-Chandra class 
 when $n$ is even, 
 and will discuss carefully
 a {\it{translation functor}}
 in Appendix III
 (Chapter \ref{sec:Translation}).

We recall from \eqref{eqn:IVZG}
 that 
\[
   r(V,\lambda)=(\sigma_1+\frac n 2-1, \sigma_2 + \frac n 2-2, \cdots, \sigma_m + \frac n 2-m, \lambda-\frac n 2), 
\]
where $m:=[\frac n 2]$.  
To specify the Weyl chamber for $W_{\mathfrak{g}}$
 that $r(V,\lambda)$ $\in (\frac 1 2 {\mathbb{Z}})^{m+1}$ belongs to, 
 we label the places
 where $\lambda-\frac n 2$ is located
 with respect to the following inequalities.  

\par\noindent
{\bf{\underline{Case 1.}}}\enspace
$n=2m$:
\[
   -\sigma_1-m+1< -\sigma_2-m+2
   < \cdots < -\sigma_m
   \le \sigma_m 
   < \cdots < \sigma_2 + m-2 < \sigma_1 + m-1;
\]
\par\noindent
{\bf{\underline{Case 2.}}}\enspace
$n=2m+1$:
\[
-\sigma_1-m+\frac 12< -\sigma_2-m+\frac 32
   < \cdots < -\sigma_{m}-\frac 1 2 < 0 
   < \sigma_{m}+ \frac 1 2 
   < \cdots < \sigma_1 + m-\frac 1 2.  
\]
Unifying these inequalities
 by adding $\frac n 2$
 to each term, 
 we may write as 
\[
  1-\sigma_1 < 2-\sigma_2 < \cdots< m - \sigma_m
  \le 
  \frac n 2
  \le \sigma_m + n-m
  < \cdots
  < \sigma_2 + n-2
  < \sigma_1 + n-1.  
\]
\begin{definition}
\label{def:RVi}
For $0 \le i \le n$, 
 we define the following subsets
\index{A}{RVi@$R(V;i)$|textbf}
 $R(V;i)$
 of ${\mathbb{Z}}$:
\begin{alignat*}{3}
&\{\lambda \in {\mathbb{Z}}: i - \sigma_{i} < \lambda < i+1-\sigma_{i+1}\}
\quad
&&\text{for $0 \le i < \frac {n-1} 2$}, 
&&
\\
&\{\lambda \in {\mathbb{Z}}: \frac{n-1}2 -\sigma_{\frac{n-1}2} < \lambda < \frac n 2\}
\quad
&&\text{for $i = \frac {n-1} 2$} \quad 
&&\text{($n$ odd)},  
\\
&\{\lambda \in {\mathbb{Z}}: \frac n 2-\sigma_{\frac n 2} < \lambda < \sigma_{\frac n 2}+\frac n 2\}
\quad
&&\text{for $i = \frac n 2$}
&&\text{($n$ even)}, 
\\
&\{\lambda \in {\mathbb{Z}}: \frac n 2 < \lambda < \sigma_{\frac {n-1} 2}+\frac {n+1} 2\}
\quad
&&\text{for $i = \frac {n+1} 2$}
&&\text{($n$ odd)}, 
\\
&\{\lambda \in {\mathbb{Z}}: \sigma_{n-i+1} +i-1< \lambda < \sigma_{n-i}+i\}
\quad
&&\text{for $\frac{n+1}2 < i \le n$}.  
&&
\end{alignat*}
\end{definition}
Here we regard $\sigma_0=\infty$.  
\begin{lemma}
\label{lem:ilsig}
Let $V \in \widehat{O(n)}$.  
We recall from \eqref{eqn:singint}
 that $S(V)$ is the set of singular integral parameters.  
\begin{enumerate}
\item[{\rm{(1)}}]
The set of regular integral parameters has the following disjoint decomposition:
\[
{\mathbb{Z}} \setminus S(V)
=
\coprod_{i=0}^n R(V;i). 
\]
In particular,
 there exists a map
\begin{equation}
\label{eqn:indexV}
   i(V, \cdot) \colon {\mathbb{Z}} \setminus S(V) \to \{0,1,\ldots,n\}
\end{equation}
such that $\lambda \in R(V;i(V, \lambda))$.  
\item[{\rm{(2)}}]
The set $S(V)$ is preserved
 by the transformations $\lambda \mapsto n-\lambda$
 and $V \mapsto V \otimes \det$, 
 and we have 
\begin{align*}
i(V,n-\lambda)\,=\,&n-i(V,\lambda)
\\
i(V \otimes \det ,\lambda)\,=\,&i(V,\lambda)
\end{align*}
for any $\lambda\in {\mathbb{Z}} \setminus S(V)$.    
\item[{\rm{(3)}}]
$R(V;\frac n 2) \ne \emptyset$
 if and only if $n$ is even and the irreducible $O(n)$-module $V$ 
 is of type Y.  
In this case,
 we have
\begin{equation}
\label{eqn:RSYV}
R(V;\frac n 2)=\{\frac n 2\} \cup S_Y(V)
\quad
\text{{\rm{(disjoint union).}}}
\end{equation}
\end{enumerate}
\end{lemma}

\begin{example}
\label{ex:isig}
Let $0 \le i \le n$.  
For the $i$-th exterior tensor representation
 $V=\Exterior^i({\mathbb{C}}^n)$ of $O(n)$, 
 we have
\[
   S(\Exterior^{(i)}({\mathbb{C}}^n))
   =
   {\mathbb{Z}} 
   \setminus 
   (\{i,n-i\}\cup (-{\mathbb{N}}_+) \cup (n+{\mathbb{N}}_+)).  
\]
Furthermore, 
 we see from Example \ref{ex:exteriorpm}
 that the set $R(V;j)$ is given as follows.  
\begin{enumerate}
\item[{\rm{(1)}}]
For $1 \le i \le n-1$, 
\[
  R(\Exterior^i({\mathbb{C}}^n);j)
  =
  \begin{cases}
  -{\mathbb{N}}_+
  \qquad
  &{\text{if $j=0$}}, 
\\
  \{j\}
  \qquad
  &{\text{if $j=i$ or $n-i$}}, 
\\
  n+{\mathbb{N}}_+
  \qquad
  &{\text{if $j=n$}}, 
\\
  \emptyset
  \qquad
  &{\text{otherwise}}.  
  \end{cases}
\]
\item[{\rm{(2)}}]
For $i =0$ or $n$, 
\[
  R(\Exterior^i({\mathbb{C}}^n);j)
  =
  \begin{cases}
  -{\mathbb{N}}
  \qquad
  &{\text{if $j=0$}}, 
\\
  n+{\mathbb{N}}
  \qquad
  &{\text{if $j=n$}}, 
\\
  \emptyset
  \qquad
  &{\text{otherwise}}.  
  \end{cases}
\]
\end{enumerate}

\end{example}

We recall from Definition \ref{def:RIntRed}
 that $\RInt$ is a subset of $\widehat {O(n)} \times {\mathbb{Z}}$.  
\begin{definition}
[height $i(V,\lambda)$]
\label{def:iVlmd}
By \eqref{eqn:indexV} in Lemma \ref{lem:ilsig}, 
 we define a map
\[
  i \colon \RInt \to \{0,1,\ldots,n\}, 
\]
see Lemma \ref{lem:regint}.  
We refer to 
\index{A}{iVlmd@$i(V,\lambda)$, height|textbf}
$i(V,\lambda)$ as the {\it{height}} of $(V,\lambda)$.  
We also refer it to as the 
\index{B}{heightIVlmd@height, of $I_{\delta}(V,\lambda)$|textbf}
height of the principal series representation
 $I_{\delta}(V,\lambda)$.  
\end{definition}

\begin{example}
\label{ex:iVlmdfig}
We illustrate the definition of the height
 $i(V,\lambda) \in \{0,1,\ldots,n\}$
 for $(V,\lambda) \in \RInt$
 by a graphic description 
 when $m(=[\frac n 2])=1$, 
 namely,
 when $n=2$ or $3$.  
In this case $G=O(n+1,1)$ is either $O(3,1)$ or $O(4,1)$, 
 and $V \in \widehat{O(n)}$ is given 
 by $V=\Kirredrep{O(n)}{\sigma_1}_{\varepsilon}$
 with $\sigma_1 \in {\mathbb{N}}$ and $\varepsilon \in \{\pm\}$.  
Then 
\begin{equation*}
\RInt \simeq
 \begin{cases}
\{(\sigma_1, \varepsilon,\lambda) \in {\mathbb{N}} \times \{\pm\}\times {\mathbb{Z}}:
\lambda-1 \ne \pm \sigma_1\}
\quad
 &\text{if $n=2$, }
\\
\{(\sigma_1, \varepsilon,\lambda) \in {\mathbb{N}} \times \{\pm\}\times {\mathbb{Z}}:
\lambda-2 \ne \pm \sigma_1, \lambda \ne 2\}
\quad
 &\text{if $n=3$.  }
\end{cases}
\end{equation*}
In the $(\sigma_1, \lambda)$-plane,
 the height $i(V,\lambda)$ is given as in Figure 14.1.


\begin{tabular}{cc}
\begin{tikzpicture}[scale=1] 
\draw [->] (0,-2) -- (0,5) node [left, inner sep=8] {$\lambda$}; 
\draw [->] (-.5,0) -- (3,0) node [below, inner sep=8] {$\sigma_1$}; 
\draw [dashed] (2.5,-1.5) -- (0,1) -- (2.5,3.5);
\fill [red] (0,0) circle [radius=.1];
\fill [red] (1,1) circle [radius=.1];
\fill [red] (0,2) circle [radius=.1];
\node at (1,3.5) {$i=2$};
\node at (2.5,1) {$i=1$};
\node at (1,-1.5) {$i=0$};
\end{tikzpicture}
&
\begin{tikzpicture}[scale=1] 
\draw [->] (0,-2) -- (0,5) node [left, inner sep=8] {$\lambda$}; 
\draw [->] (-.5,0) -- (3.5,0) node [below, inner sep=8] {$\sigma_1$}; 
\draw [dashed] (0,2) -- (2.5,4.5);
\draw [dashed] (0,1.5) -- (3,1.5);
\draw [dashed] (2.5,-1.5) -- (0,1);
\fill [red] (0,0) circle [radius=.1];
\fill [red] (1,1) circle [radius=.1];
\fill [red] (1,2) circle [radius=.1];
\fill [red] (0,3) circle [radius=.1];
\node at (1,4) {$i=3$};
\node at (2.5,2.5) {$i=2$};
\node at (2.5,.75) {$i=1$};
\node at (1,-1.5) {$i=0$};
\end{tikzpicture}
\\[1em]
$G=O(3,1)$ & $G=O(4,1)$
\\
$(n=2)$ & $(n=3)$
\end{tabular}
$$
\textrm{
Figure 14.1:
The height $i=i(V,\lambda)$ for $(V,\lambda) \in \RInt$ when $n=2,3$. 
}
$$
%

The red dots stand for 
 $(V,\lambda)=(\Exterior^j({\mathbb{C}}^n), j)$ 
 when $j=0,1,\ldots,n$. 
\end{example}

The case 
 where the height $i(V,\lambda)$ is equal to $\frac n 2$
 requires a special attention. 
\begin{lemma}
\label{lem:isig}
Let $m:=[\frac n 2]$.  
Suppose that $V =\Kirredrep {O(n)}{\sigma}_{\varepsilon}$ 
 with $\sigma \in \Lambda^+(m)$
 and $\varepsilon \in \{\pm\}$, 
 and $\lambda \in {\mathbb{Z}} \setminus S(V)$.  
\begin{enumerate}
\item[{\rm{(1)}}]
The height $i(V, \lambda)$ is equal to $\frac n 2$
 if and only if $n$ is even
 $(=2m)$
 and $\sigma_{m} > |\lambda -m|$.  
\item[{\rm{(2)}}]
If $\lambda \in S_Y(V)$
 (see \eqref{eqn:SYsigma}), 
 then $n$ is even $(=2m)$
 and $i(V, \lambda)=m$.  
\item[{\rm{(3)}}]
Suppose that $V$ is of type Y
 (Definition \ref{def:OSO}).  
Then,  
 for $(V, \lambda) \in \Reducible$, 
 the following two conditions are equivalent:
\begin{enumerate}
\item[{\rm{(i)}}]
$i(V,\lambda) = \frac n 2$;
\item[{\rm{(ii)}}]
$n$ is even 
 and $\lambda = \frac n 2$.  
\end{enumerate}
\end{enumerate}
\end{lemma}

\subsection{$K$-type formul{\ae} of irreducible $G$-modules
\label{subsec:Ktypes}}

In this section we provide explicit $K$-type formul{\ae}
 of irreducible representations of $G=O(n+1,1)$.  
The height $i(V,\lambda)$ plays a crucial role
 in describing the $K$-type formul{\ae}
 of irreducible subquotients
 of $I_{\delta}(V, \lambda)$, 
 see Proposition \ref{prop:KXred} (1).  

\subsubsection{$K$-type formula of $I_{\delta}(V, \lambda)$}
We begin with the $K$-type formula
 of the principal series representation
 $I_{\delta}(V, \lambda)$
 which generalizes Lemma \ref{lem:KtypeIi}
 for $I_{\delta}(i,i)$ 
 in the setting that $V=\Exterior^i({\mathbb{C}}^n)$.  

\begin{proposition}
[$K$-type formula of $I_{\delta}(V,\lambda)$]
\label{prop:KtypeIV}
\index{B}{Ktypeformula@$K$-type formula}
Let $G=O(n+1,1)$ and $m=[\frac n 2]$.  
Suppose that $V=\Kirredrep{O(n)}{\sigma}_{\varepsilon}$
 with $\sigma=(\sigma_1, \cdots,\sigma_m) \in \Lambda^+(m)$
 and $\varepsilon \in \{\pm\}$.  
\begin{enumerate}
\item[{\rm{(1)}}]
For $n=2m+1$, 
 the $K$-type formula 
 of the principal series representation 
 $I_{\delta}(V,\lambda)$ is given by
\[
  \bigoplus_{\mu} \Kirredrep{O(n+1)}{\mu_1, \cdots,\mu_{m+1}}_{\varepsilon}
  \boxtimes \delta(-1)^{\sum_{j=1}^{m+1} \mu_j - \sum_{j=1}^{m}\sigma_j},
\]
where $\mu=(\mu_1, \cdots,\mu_{m+1})$
 runs over $\Lambda^+(m+1)$ subject to 
\begin{equation}
\label{eqn:Kmuodd}
\mu_1 \ge \sigma_1 \ge \mu_2 \ge \sigma_2 \ge \cdots \ge \mu_m \ge \sigma_m 
\ge \mu_{m+1} \ge 0.  
\end{equation}

\item[{\rm{(2)}}]
For $n=2m$ and $V \in \widehat{O(n)}$ of type X
 (Definition \ref{def:OSO}), 
 the $K$-type formula 
 of $I_{\delta}(V,\lambda)$ is given by
\[
  \bigoplus_{\mu} \Kirredrep{O(n+1)}{\mu_1, \cdots,\mu_{m}}_{\varepsilon}
  \boxtimes \delta(-1)^{\sum_{j=1}^{m} \mu_j - \sum_{j=1}^{m}\sigma_j},
\]
where $\mu=(\mu_1, \cdots,\mu_{m})$
 runs over $\Lambda^+(m+1)$ subject to 
\begin{equation}
\label{eqn:Kmueven1}
\mu_1 \ge \sigma_1 \ge \cdots \ge \mu_m \ge \sigma_m 
(=0).  
\end{equation}

\item[{\rm{(3)}}]
For $n=2m$ and $V\in \widehat{O(n)}$ of type Y, 
 the $K$-type formula of $I_{\delta}(V,\lambda)$ is given by
\[
\bigoplus_{\kappa=\pm}
\bigoplus_{\mu} \Kirredrep{O(n+1)}{\mu_1, \cdots,\mu_{m}}_{\kappa\varepsilon}
  \boxtimes \kappa \delta(-1)^{\sum_{j=1}^{m} \mu_j - \sum_{j=1}^{m}\sigma_j},
\]
where $\mu=(\mu_1, \cdots,\mu_{m})$
 runs over $\Lambda^+(m)$ subject to 
\begin{equation}
\label{eqn:Kmueven2}
\mu_1 \ge \sigma_1 \ge \cdots \ge \mu_m \ge \sigma_m 
(>0).  
\end{equation}
\end{enumerate}
\end{proposition}

\begin{proof}
By the Frobenius reciprocity, 
Proposition \ref{prop:KtypeIV} follows from the classical
 branching rule for the restriction
 $O(n+1) \downarrow O(n)$, 
 see Fact \ref{fact:ONbranch}.  
\end{proof}

Since the principal series representation $I_{\delta}(V,\lambda)$
 of $G$ is multiplicity-free as a $K$-module,
 any subquotient of $I_{\delta}(V,\lambda)$
 can be characterized 
 by its $K$-types.  
In the next subsection, 
 we provide $K$-type formul{\ae}
 of subquotients of $I_{\delta}(V,\lambda)$
 based on Proposition \ref{prop:KtypeIV}.

\subsubsection{$K$-types of subquotients
 $I_{\delta}(V,\lambda)^{\flat}$ and $I_{\delta}(V,\lambda)^{\sharp}$}

We recall from \eqref{eqn:reducible}
 and Theorem \ref{thm:irrIV}
 that the following two conditions on $(V,\lambda) \in \widehat{O(n)} \times {\mathbb{C}}$ 
 are equivalent.  
\begin{enumerate}
\item[(i)]
$(V,\lambda) \in \Reducible$, 
 {\it{i.e.}}, the $G$-module $I_{\delta}(V,\lambda)$ is reducible; 
\item[(ii)]
$\lambda \in {\mathbb{Z}} \setminus (S(V) \cup S_Y(V))$.  
\end{enumerate}

We note that $\lambda = \frac n 2$ belongs
 to ${\mathbb{Z}} \setminus (S(V) \cup S_Y(V))$
 when $n$ is even.

In this section, 
  we describe the $K$-types of the subquotients
 $I_{\delta}(V,\lambda)^{\flat}$ and $I_{\delta}(V,\lambda)^{\sharp}$
 when the principal series representation
 $I_{\delta}(V,\lambda)$ is reducible, 
 {\it{i.e.}}, 
 when $(V,\lambda) \in \Reducible$, 
 see \eqref{eqn:Red}.

We shall see that the description 
 depends on the height $i(V,\lambda)$ 
 (Definition \ref{def:iVlmd})
 when $\lambda = \frac n 2$.  
To be more precise,
 let $m=[\frac n 2]$ and $V \in \widehat{O(n)}$.  
Suppose $\lambda \in {\mathbb{Z}} \setminus(S(V) \cup S_Y(V))$
 and we define $i$ to be the height
 $i(V,\lambda) \in \{0,1,\ldots,n\}$.  
We write $V=\Kirredrep{O(n)}{\sigma}_{\varepsilon}$
 with $\sigma=(\sigma_1, \cdots, \sigma_m)\in \Lambda^+(m)$
 and $\varepsilon \in \{\pm\}$ as before.  
We observe the following:
\begin{enumerate}
\item[$\bullet$]
if $i< \frac{n-1}{2}$, 
 then $1 \le i+1 \le m$
 and the condition $i-\sigma_i < \lambda< i+1 -\sigma_{i+1}$
 (Definition \ref{def:RVi})
 amounts to 
\begin{equation}
\label{eqn:barK1}
 \sigma_{i+1} \le i-\lambda
\quad\text{and}\quad 
 i-\lambda+1 \le \sigma_i;
\end{equation}
\item[$\bullet$]
if $i=\frac{n-1}{2}$, 
 then $n$ is odd $(=2m+1)$ 
 and we have
\begin{equation}
\label{eqn:barK3}
0 \le m-\lambda
\quad\text{and}\quad
m-\lambda+1 \le \sigma_m;
\end{equation}

\item[$\bullet$]
if $i=\frac{n+1}{2}$, 
 then $n$ is odd $(=2m+1)$ 
 and we have
\begin{equation}
\label{eqn:barK4}
0 \le \lambda-m-1
\quad\text{and}\quad
\lambda-m \le \sigma_m; 
\end{equation}
\item[$\bullet$]
if $\frac{n+1}2 < i$, 
 then $ 1 \le n-i+1 \le m$
 and the condition $\sigma_{n-i+1}+ i-1 < \lambda< \sigma_{n-i} +i$
 amounts to  
\begin{equation}
\label{eqn:barK2}
\sigma_{n-i+1} \le \lambda-i
\quad\text{and}\quad 
\lambda-i+1 \le \sigma_{n-i}.  
\end{equation}
\end{enumerate}

We recall that the principal series representation
$I_{\delta}(V,\lambda)$
 of $G=O(n+1,1)$
 is $K$-multiplicity-free,
 and its $K$-type formula 
 is given explicitly in Proposition \ref{prop:KtypeIV}.  
To describe the $K$-type formul{\ae}
 of subquotients of $I_{\delta}(V,\lambda)$, 
 we use the inequalities
 \eqref{eqn:barK1}--\eqref{eqn:barK2}
 in Proposition \ref{prop:KXred} (1) below.  

\begin{proposition}
[$K$-type formul{\ae} of subquotients]
\label{prop:KXred}
Suppose that $(V,\lambda) \in \Reducible$, 
 or equivalently, 
 $V \in \widehat{O(n)}$
 and $\lambda \in {\mathbb{Z}} \setminus(S(V) \cup S_Y(V))$, 
 see Theorem \ref{thm:irrIV}.  
Let 
\index{A}{iVlmd@$i(V,\lambda)$, height}
$i:=i(V,\lambda) \in \{0,1,\ldots,n\}$
 be the height of $(V,\lambda)$
as in Definition \ref{def:iVlmd}.  
\begin{enumerate}
\item[{\rm{(1)}}]
Suppose $\lambda \ne \frac n2$.  
In this case $i \ne \frac n 2$. 
Then the $K$-types of the submodule $I_{\delta}(V,\lambda)^{\flat}$
 and the quotient $I_{\delta}(V,\lambda)^{\sharp}$
 of $I_{\delta}(V,\lambda)$, 
 see Proposition \ref{prop:Xred}, 
 are subsets of the $K$-types of $I_{\delta}(V,\lambda)$
 (Proposition \ref{prop:KtypeIV})
 characterized by the following additional inequalities:
\par\noindent
$\bullet$\enspace
for $i \le \frac{n-1}2$, 
the condition $\sigma_{i+1} \le \mu_{i+1} \le \sigma_i$
 in \eqref{eqn:Kmuodd}--\eqref{eqn:Kmueven2}
 is divided as follows:
\begin{alignat*}{5}
(\sigma_{i+1} 
&\le)\,\,
&&\mu_{i+1} 
&& 
&&\le i-\lambda
\qquad
&&\text{for $I_{\delta}(V,\lambda)^{\flat}$}, 
\\
 i-\lambda+1
&\le\,\,
&&\mu_{i+1} 
&& (
&&\le \sigma_i)
\qquad
&&\text{for $I_{\delta}(V,\lambda)^{\sharp}$};
\\
\intertext{$\bullet$\enspace for $\frac{n+1}{2}\le i$, 
the condition $\sigma_{n-i+1} \le \mu_{n-i+1} \le \sigma_{n-i}$
 in \eqref{eqn:Kmuodd}--\eqref{eqn:Kmueven2}
 is divided as follows:}
 \lambda-i+1 
&\le\,\,
&&\mu_{n-i+1} 
&& (
&&\le \sigma_{n-i})
\qquad
&&\text{for $I_{\delta}(V,\lambda)^{\flat}$}, 
\\
(\sigma_{n-i+1} 
&\le)\,\,
&&\mu_{n-i+1} 
&& 
&&\le \lambda-i
\qquad
&&\text{for $I_{\delta}(V,\lambda)^{\sharp}$}.  
\end{alignat*}
Here we regard $\sigma_{m+1}=0$
 (this happens when $i=\frac{n \pm 1}{2}$).  
\item[{\rm{(2)}}]
Suppose $\lambda =\frac n 2$.  
In this case $n$ is even $(=2m)$ and $i=m$.  
Then the $K$-types
 of the submodules $I_{\delta}(V,\lambda)^{\flat}$
 and $I_{\delta}(V,\lambda)^{\sharp}$
 of the (tempered) principal series representation
 $I_{\delta}(V,\lambda)$, 
 see Proposition \ref{prop:IVmtemp}, 
 are given by 
\begin{alignat*}{2}
  &\bigoplus_{\mu}
  \Kirredrep{O(n+1)}{\mu_1, \cdots,\mu_{m}}_{\varepsilon}
  \boxtimes 
   \delta(-1)^{\sum_{j=1}^{m} \mu_j - \sum_{j=1}^{m}\sigma_j}
   \quad
&&\text{for $I_{\delta}(V,\lambda)^{\flat}$}, 
\\
  &\bigoplus_{\mu}
  \Kirredrep{O(n+1)}{\mu_1, \cdots,\mu_{m}}_{-\varepsilon}
  \boxtimes 
   \delta(-1)^{\sum_{j=1}^{m} \mu_j - \sum_{j=1}^{m}\sigma_j-1}
   \quad
&&\text{for $I_{\delta}(V,\lambda)^{\sharp}$}
\end{alignat*}
where $\mu=(\mu_1, \cdots,\mu_{m})$
 runs over $\Lambda^+(m)$
 subject to \eqref{eqn:Kmueven2}.  
\end{enumerate}
\end{proposition}

\begin{proof}
The $K_0$-types for all irreducible subquotients
 of principal series representations
 of the connected Lie group $G_0=SO_0(n+1,1)$
 were obtained in Hirai \cite{Hirai62}, from
 which analogous results
 for the group $\overline G=S O(n+1,1)$
 are easily shown.  
Our concern is with the group $G=O(n+1,1)$.  
Then the first assertion follows from Proposition \ref{prop:KtypeIV}
 on the $K$-type formula 
 of $I_{\delta}(V,\lambda)$
 and from the branching rule for the restriction
 $G \downarrow \overline G$
 in Propositions \ref{prop:20180906}
 and \ref{prop:1808104} in Appendix II.  
The second assertion follows from the branching rule
 of $I_{\delta}(V,\frac n 2)$
 for the restriction $G\downarrow \overline G$
 in Proposition \ref{prop:180572}.  
\end{proof}

\subsection{$(\delta, V, \lambda)
\rightsquigarrow (\delta^{\uparrow},V^{\uparrow}, \lambda^{\uparrow})$
 and $(\delta^{\downarrow},V^{\downarrow}, \lambda^{\downarrow})$}
\label{subsec:HasseV}

In this section,
 we introduce a correspondence
\begin{alignat*}{4}
&\delta \in \{\pm\}, 
&&V \in \widehat{O(n)}, 
&&
\text{and }
&&\lambda \in {\mathbb{Z}} \setminus(S(V) \cup S_Y(V))
\\
&
&&
&&
&&
\rotatebox{270}{$\rightsquigarrow$}
\\
&\delta^{\uparrow} \in \{\pm\},\,\,
&&V^{\uparrow} \in \widehat{O(n)},\,\, 
&&
\text{and }\,\,
&&
\lambda^{\uparrow} \in {\mathbb{Z}} \setminus(S(V^{\uparrow}) \cup S_Y(V^{\uparrow}))
\end{alignat*}
satisfying the following two properties (Proposition \ref{prop:Vlmdup}):
\begin{align*}
  i(V^{\uparrow}, \lambda^{\uparrow}) =\,& i(V,\lambda)+1
\\
  I_{\delta^{\uparrow}}(V^{\uparrow}, \lambda^{\uparrow})^{\flat} \simeq\,
& I_{\delta}(V,\lambda)^{\sharp}.  
\end{align*}
We retain the notation
 that $G=O(n+1,1)$
 and $m = [\frac n 2]$.  

\begin{definition}
\label{def:180909}
Suppose that 
 $V=\Kirredrep{O(n)}{\sigma}_{\varepsilon}$
 with $\sigma \in \Lambda^+(m)$
 and $\varepsilon \in \{\pm\}$, 
 and $\lambda \in {\mathbb{Z}} 
 \setminus S(V)$.  
Let $i:=i(V,\lambda) \in \{0,1,\ldots,n\}$
 be the height of $(V,\lambda)$
 as in Lemma \ref{lem:ilsig}.  
\begin{enumerate}
\item[{\rm{(1)}}]
We assume $0 \le i \le n-1$, 
 or equivalently, 
 $\lambda \le \sigma_1 +n-1$.  
We define
\index{A}{0deltaVlmd1@$(\delta,V,\lambda)^{\uparrow}
 =(\delta^{\uparrow},V^{\uparrow},\lambda^{\uparrow})$|textbf}
\begin{equation}
\label{eqn:upVlmd}
 (\delta,V,\lambda)^{\uparrow} \equiv (\delta^{\uparrow},V^{\uparrow},\lambda^{\uparrow})
\in {\mathbb{Z}}/2{\mathbb{Z}} \times \widehat{O(n)} \times {\mathbb{Z}}
\end{equation}
with $V^{\uparrow} := \Kirredrep{O(n)}{\sigma^{\uparrow}}_{\varepsilon}$ 
 as follows:

\begin{enumerate}
\item[$\bullet$]
For $\lambda < \frac n 2$, 
 we have $0 \le i < \frac n2$
 and set
\begin{align*}
\delta^{\uparrow} := & \delta(-1)^{i+1-\sigma_{i+1}-\lambda}, 
\\
\sigma^{\uparrow} := & (\sigma_{1}, \cdots, \sigma_{i}, i+1-\lambda, 
 \sigma_{i+2}, \cdots, \sigma_{m}), 
\\
\lambda^{\uparrow} := & i+1-\sigma_{i+1}.  
\end{align*}
\item[$\bullet$]
For $\frac n 2 \le \lambda \le \sigma_1 +n-1$, 
 we have $\frac n2 \le i \le n-1$
 and set
\begin{align*}
\delta^{\uparrow} := & \delta(-1)^{\lambda-\sigma_{n-i}-i}, 
\\
\sigma^{\uparrow} := & (\sigma_{1}, \cdots, \sigma_{n-i-1}, \lambda-i, 
 \sigma_{n-i+1}, \cdots, \sigma_{m}), 
\\
\lambda^{\uparrow} := & \sigma_{n-i} +i.  
\end{align*}
\end{enumerate}

\item[{\rm{(2)}}]
Conversely,
 for $1 \le i \le n$, 
 namely, 
 for $1-\sigma_1 \le \lambda$, 
 we define
\index{A}{0deltaVlmd2@$(\delta,V,\lambda)^{\downarrow}
 =(\delta^{\downarrow},V^{\downarrow},\lambda^{\downarrow})$|textbf}
\begin{equation}
\label{eqn:downVlmd}
  (\delta,V,\lambda)^{\downarrow}
  \equiv
  (\delta^{\downarrow},V^{\downarrow},\lambda^{\downarrow})
\end{equation}
 as the inverse of the correspondence
\[
   (\delta,V,\lambda) \mapsto (\delta,V,\lambda)^{\uparrow}.  
\]
\end{enumerate}
\end{definition}

A prototype for Definition \ref{def:180909}
 appeared implicitly
 in Theorem \ref{thm:LNM20}
 for the principal series representations $I_{\delta}(i,i)$
 having the trivial ${\mathfrak{Z}}_G({\mathfrak{g}})$-infinitesimal character $\rho_G$.  
We now explain this explicitly
 as an example for $(V,\lambda)=(\Exterior^i({\mathbb{C}}^n), i)$
 $(1 \le i \le n)$:  
\begin{example}
\label{ex:up}
For the exterior representations 
 $\Exterior^i({\mathbb{C}}^n)$
 of $O(n)$, 
 we have
\begin{alignat*}{2}
 (\delta,\Exterior^i({\mathbb{C}}^n), i)^{\uparrow}
 =&
 (-\delta,\Exterior^{i+1}({\mathbb{C}}^n), i+1)
 \quad
&&\text{for $0 \le i \le n-1$, }
\\
(\delta,\Exterior^i({\mathbb{C}}^n), i)^{\downarrow}
 =&
 (-\delta,\Exterior^{i-1}({\mathbb{C}}^n), i-1)
 \quad
&&\text{for $1 \le i \le n$.  }
\end{alignat*}
The proof follows directly from 
 the definition, 
 see Example \ref{ex:exteriorpm}.  
\end{example}

Here are basic properties of the correspondence
\[
  (\delta,V,\lambda) \mapsto
(\delta^{\uparrow},V^{\uparrow},\lambda^{\uparrow})
\quad
\text{or}
\quad
(\delta^{\downarrow}, V^{\downarrow}, \lambda^{\downarrow}).  
\]

\begin{proposition}
\label{prop:Vlmdup}
Suppose that $(V,\lambda) \in \RInt$, 
 {\it{i.e.}}, 
 $V \in \widehat{O(n)}$
 and $\lambda \in {\mathbb{Z}} \setminus S(V)$.  
In what follows, 
 we assume the height $i(V,\lambda)$ is not equal to $n$
 when we consider $(V^{\uparrow}, \lambda^{\uparrow})$, 
 and is nonzero 
 when we consider $(V^{\downarrow}, \lambda^{\downarrow})$.  
\begin{enumerate}
\item[{\rm{(1)}}]
$r(V^{\uparrow},\lambda^{\uparrow}), 
r(V^{\downarrow},\lambda^{\downarrow})
\in W_G\, r(V,\lambda)$, 
 see \eqref{eqn:IVZG}.  
In particular,
 $(V^{\uparrow},\lambda^{\uparrow})$, 
 $(V^{\downarrow},\lambda^{\downarrow}) \in \RInt$.  
\item[{\rm{(2)}}]
$i(V^{\uparrow}, \lambda^{\uparrow})-1=i(V, \lambda)
 =i(V^{\downarrow}, \lambda^{\downarrow})+1$.  
\item[{\rm{(3)}}]
$\delta^{\uparrow} (-1)^{\lambda^{\uparrow}}
=\delta (-1)^{\lambda}
=\delta^{\downarrow} (-1)^{\lambda^{\downarrow}}$. 
\item[{\rm{(4)}}]
$(V^{\uparrow}, \lambda^{\uparrow}),  
(V^{\downarrow}, \lambda^{\downarrow})\in \Reducible$, 
 if $(V, \lambda) \in \Reducible$,
 see \eqref{eqn:reducible}.  
\item[{\rm{(5)}}]
Suppose that $(V,\lambda) \in \Reducible$
 and $\lambda \ne \frac n 2$.  
Then the unique submodule of $I_{-\delta}(V^{\uparrow}, \lambda^{\uparrow})$
 is isomorphic to the unique quotient
 of $I_{\delta}(V, \lambda)$, 
 that is, 
 we have the following $G$-isomorphisms
 with the notation 
 as in Proposition \ref{prop:Xred}:
\begin{align*}
  I_{\delta^{\uparrow}}(V^{\uparrow}, \lambda^{\uparrow})^{\flat}
  \simeq&
  I_{\delta}(V, \lambda)^{\sharp}, 
\\
  I_{\delta^{\downarrow}}(V^{\downarrow}, \lambda^{\downarrow})^{\sharp}
  \simeq&
  I_{\delta}(V, \lambda)^{\flat}.  
\end{align*}
\end{enumerate}
\end{proposition}

With these notations,
 we give the formul\ae\ for the minimal $K$-types
 of the irreducible subquotients 
 $I_{\delta}(V,\lambda)^{\flat}$ and $I_{\delta}(V,\lambda)^{\sharp}$
 in $I_{\delta}(V,\lambda)$
 in the setting of Proposition \ref{prop:Xred}.  

\begin{proposition}
\label{prop:minKXred}
Let $G=O(n+1,1)$ and $m=[\frac n2]$.  
Suppose $V=\Kirredrep{O(n)}{\sigma}_{\varepsilon}$
 with $\sigma=(\sigma_1, \cdots,\sigma_m) \in \Lambda^+(m)$
 and $\varepsilon \in \{\pm\}$.  
Let $\delta \in \{\pm\}$ 
 and $\lambda \in {\mathbb{Z}} \setminus (S(V) \cup S_Y(V) \cup \{\frac n 2\})$. 
\begin{enumerate}
\item[{\rm{(1)}}]
The minimal $K$-types of $I_{\delta}(V,\lambda)^{\flat}$ for $\lambda < \frac n 2$
 and of $I_{\delta}(V,\lambda)^{\sharp}$
 for $\lambda>\frac n 2$ are given by
\begin{alignat*}{2}
&\Kirredrep{O(n+1)}{\sigma}_{\varepsilon}\boxtimes \delta
&&\text{for $n=2m$ and $\sigma_m=0$, }
\\
&\Kirredrep{O(n+1)}{\sigma}_{\varepsilon}\boxtimes \delta,
\quad
\Kirredrep{O(n+1)}{\sigma}_{-\varepsilon}\boxtimes (-\delta)
\qquad
&&\text{for $n=2m$ and $\sigma_m>0$, }
\\
&\Kirredrep{O(n+1)}{\sigma,0}_{\varepsilon}\boxtimes \delta
\qquad
&&\text{for $n=2m+1$.}
\end{alignat*}
\item[{\rm{(2)}}]
The minimal $K$-types of $I_{\delta}(V,\lambda)^{\sharp}$ for $\lambda < \frac n 2$
 and of $I_{\delta}(V,\lambda)^{\flat}$
 for $\lambda > \frac n 2$ are given by
\begin{alignat*}{2}
&\Kirredrep{O(n+1)}{\sigma^{\uparrow}}_{\varepsilon}\boxtimes \delta^{\uparrow}
&&\text{for $n=2m$ and $\sigma_m=0$, }
\\
&\Kirredrep{O(n+1)}{\sigma^{\uparrow}}_{\varepsilon}\boxtimes \delta^{\uparrow},
\quad
\Kirredrep{O(n+1)}{\sigma^{\uparrow}}_{-\varepsilon}\boxtimes (-\delta^{\uparrow})
\qquad
&&\text{for $n=2m$ and $\sigma_m>0$, }
\\
&\Kirredrep{O(n+1)}{\sigma^{\uparrow},0}_{\varepsilon}\boxtimes \delta^{\uparrow}
\qquad
&&\text{for $n=2m+1$.}
\end{alignat*}
\end{enumerate}
\end{proposition}

\subsection{Classification of irreducible admissible representations
 of $G=O(n+1,1)$}

Irreducible admissible representations
 of the connected group $G_0=SO_0(n+1,1)$
 were classified infinitesimally
 ({\it{i.e.,}} on the level of $({\mathfrak{g}}, K_0)$-modules)
 by Hirai \cite{Hirai62}, 
 see also Borel--Wallach \cite{BW}
 and Collingwood \cite[Chap.~5]{C}.  
However,
 we could not find in the literature
 a classification of irreducible admissible representations
 of the indefinite orthogonal group $G=O(n+1,1)$, 
 which is not in the Harish-Chandra class
 when $n$ is even.  
For the sake of completeness,
 we give an infinitesimal classification 
 of irreducible admissible representations of $G$, 
 or equivalently,
 give a classification
 of irreducible $({\mathfrak{g}}, K)$-modules
 in this section.  
Moreover we give three characterizations
 of the irreducible representations of $G$
 when they are neither principal series representations
 nor tempered representations, 
 see Theorem \ref{thm:1808116}.

\subsubsection{Characterizations of the irreducible subquotients
 $\Pi_{\delta}(V,\lambda)$}

We recall from Section \ref{subsec:subIii}
 the irreducible representations $\Pi_{\ell, \delta}$
 of $G$
 that have the trivial ${\mathfrak{Z}}_G({\mathfrak{g}})$-infinitesimal character $\rho_G$.  
Analogously to the notation $\Pi_{\ell, \delta}$
 in \eqref{eqn:Pild}
 for $\operatorname{Irr}(G)_{\rho}$, 
 we set
\index{A}{1PiideltaV@$\Pi_{\delta}(V,\lambda)$|textbf}
\begin{equation}
\label{eqn:PiVlmd}
  \Pi_{\delta}(V,\lambda):=I_{\delta}(V, \lambda)^{\flat}
\end{equation}
for $\delta \in \{\pm\}$
 and $(V,\lambda) \in \Reducible$.  
If $i(V,\lambda) \ne 0$, 
 then we have a $G$-isomorphism
\begin{equation}
\label{eqn:PiVlmddown}
  \Pi_{\delta}(V,\lambda) \simeq I_{\delta^{\downarrow}}(V^{\downarrow}, \lambda^{\downarrow})^{\sharp}, 
\end{equation}
where $(\delta^{\downarrow}, V^{\downarrow}, \lambda^{\downarrow})$
 is given in Definition \ref{def:180909}.  
We also have a $G$-isomorphism
\begin{equation}
\label{eqn:PiVlmddual}
  \Pi_{\delta}(V,\lambda) \simeq I_{\delta}(V, n-\lambda)^{\sharp}.  
\end{equation}
We have already discussed
 in Proposition \ref{prop:IVmtemp}
 irreducible subquotients
 of reducible tempered principal series representations
 $I_{\delta}(V, \lambda)$ 
under the assumption
 that $(V, \lambda)\in\Reducible$
 with $\lambda=\frac n 2$.  
This assumption implies that $n$ is even, 
 $V$ is of type Y
 and $\lambda = \frac n 2$.  
The next theorem discusses the remaining
 (and the important) case
 when the principal series representation $I_{\delta}(V, \lambda)$ 
 is reducible, 
 namely, 
 $(V,\lambda) \in \Reducible$
 with an additional condition $\lambda \ne \frac n 2$.  

\begin{theorem}
[characterizations of $\Pi_{\delta}(V, \lambda)$]
\label{thm:1808116}
Let $G=O(n+1,1)$, 
 and we set $m:=[\frac n 2]$.  
Assume that $(V, \lambda)\in\Reducible$.  
This means that $V \in \widehat{O(n)}$ and
\[
  \lambda \in 
\begin{cases}
  {\mathbb{Z}} \setminus S(V)
  \quad&\text{if $n=2m+1$, }
\\
{\mathbb{Z}} \setminus (S(V) \cup S_Y(V))
  \quad&\text{if $n=2m$, }
\end{cases}
\]
see Theorem \ref{thm:irrIV}.  
We further assume that $\lambda \ne \frac n 2$.  
\begin{enumerate}
\item[{\rm{(1)}}]
{\rm{(Langlands subrepresentation of principal series)}}\enspace
For $\delta \in \{\pm\}$, 
 $\Pi_{\delta}(V, \lambda)$ is the unique proper $G$-submodule
 of $I_{\delta}(V, \lambda)$.  

\item[{\rm{(2)}}]
{\rm{($\theta$-stable parameter)}}
\enspace
Let $i:=i(V,\lambda) \in \{0,1,\ldots,n\}$
 be the height of $(V,\lambda)$
 as in \eqref{eqn:indexV}.  
We write $V=\Kirredrep{O(n)}{\sigma}_{\varepsilon}$
 with $\sigma=(\sigma_1, \cdots,\sigma_m) \in \Lambda^+(m)$
 and $\varepsilon \in \{\pm\}$.   
Then the underlying $({\mathfrak {g}}, K)$-module of $\Pi_{\delta}(V, \lambda)$
 is given by means of $\theta$-stable parameter
 (see Section \ref{subsec:thetapara})
 as 
\[
\Pi_{\delta}(V, \lambda)_K
\simeq
\begin{cases}
\Rq{\sigma_1-1, \cdots, \sigma_i-1}{i-\lambda,\sigma_{i+1}, \cdots, \sigma_m}
{\varepsilon,\delta \varepsilon}
&\text{if $\lambda < \frac n 2$, }
\\
\Rq{\sigma_1-1, \cdots, \sigma_{n-i}-1, \lambda-i}{\sigma_{n-i+1}, \cdots, \sigma_m}
{\varepsilon,-\delta \varepsilon}
\,\,
&\text{if $\frac n 2< \lambda$.  }
\end{cases}
\]

\item[{\rm{(3)}}]
{\rm{(coherent family starting at $\Pi_{i,\delta} \in \operatorname{Irr}(G)_{\rho}$)}}\enspace
We set 
\begin{equation*}
   \text{$r(V, \lambda)\in {\mathbb{C}}^{m+1}$ $(\simeq {\mathfrak{h}}_{\mathbb{C}}^{\ast})$ 
         as in \eqref{eqn:IVZG}}.  
\end{equation*}
Denote by 
\index{A}{Pmu@$P_{\mu}$}
$P_{\mu}$ 
the projection to the primary component 
 with the generalized ${\mathfrak{Z}}_G({\mathfrak{g}})$-infinitesimal character
 $\mu \in {\mathfrak{h}}_{\mathbb{C}}^{\ast} \mod W_G$
 (see Section \ref{subsec:primary} in Appendix III).  
Let 
\index{A}{FVlmd@$F(V,\lambda)$}
$F(V,\lambda)$ be the irreducible finite-dimensional representation
 of $G=O(n+1,1)$, 
 which will be defined in Definition \ref{def:Fshift} 
 in Appendix III.  
Then there is a natural $G$-isomorphism:
\[
\Pi_{\delta}(V,\lambda)
\simeq
P_{r(V, \lambda)}(\Pi_{i,\delta} \otimes F(V,\lambda)). 
\]

\item[{\rm{(4)}}]
{\rm{(Hasse sequence and standard sequence starting at $F(V,\lambda)$)}}
\enspace
\index{B}{Hassesequence@Hasse sequence}
\index{B}{standardsequence@standard sequence}
Let 
\index{A}{1PiideltaF@$\Pi_i(F)$, standard sequence starting with $F$}
$\Pi_j(F)$
 $(j=0,1,\cdots,n)$
 be the standard sequence 
 starting with an irreducible finite-dimensional 
 representation $F$ of $G$
 (Definition \ref{def:Hasse}), 
 and $i=i(V,\lambda)$ the height of $(V,\lambda)$.  
Then there is a natural $G$-isomorphism:
\[
\Pi_{\delta}(V, \lambda)
\simeq
\Pi_{i}(F(V,\lambda)) \otimes \chi_{+\delta}. 
\]
\end{enumerate}
\end{theorem}

See Proposition \ref{prop:Xred} for (1), 
 Theorem \ref{thm:IndV} for (2), 
 Theorem \ref{thm:1807113} for (3)
 in Chapter \ref{sec:Translation}
 (Appendix III), 
 and Theorems \ref{thm:171471} and \ref{thm:171471b} for (4).

\subsubsection{Classification of $\operatorname{Irr}(G)$}
\label{subsec:IrrG}
We give an infinitesimal classification 
 of irreducible admissible representations of $G=O(n+1,1)$.  
One may reduce the proof to the case 
 of connected groups, 
 by inspecting the restriction 
 to the subgroups $\overline G=SO(n+1,1)$
 or $G_0=SO_0(n+1,1)$, 
 see Chapter \ref{sec:SOrest} (Appendix II).

\begin{theorem}
[classification of $\operatorname{Irr}(G)$]
\label{thm:irrG}
Irreducible admissible representations
 of moderate growth
 of $G=O(n+1,1)$ are listed as follows:
\index{A}{1PiideltaV@$\Pi_{\delta}(V,\lambda)$}
\begin{alignat*}{2}
&\bullet\enspace I_{\delta}(V, \lambda)
\qquad
&& \lambda \in ({\mathbb{C}} \setminus {\mathbb{Z}}) \cup S(V) \cup S_Y(V), 
\\
&\bullet\enspace 
\Pi_{\delta}(V, \lambda)
\qquad
&&\text{$\lambda \in {\mathbb{Z}} \setminus (S(V) \cup S_Y(V))$ and $\lambda \le \frac n 2$, }
\end{alignat*}
where $V \in \widehat{O(n)}$
 and $\delta \in \{\pm\}$.  
\end{theorem}
We note that there is an isomorphism
 of irreducible $G$-modules:
\[
  I_{\delta}(V, \lambda)\simeq I_{\delta}(V, n-\lambda)
\]
 when $\lambda \in ({\mathbb{C}} \setminus {\mathbb{Z}}) \cup S(V) \cup S_Y(V)$.

\subsection{$\theta$-stable parameters
 and cohomological parabolic induction}
\label{subsec:thetapara}

In this section
 we give a parametrization of irreducible subquotients
 of the principal series representations
\[
    I_{\delta}(V,\lambda)={\operatorname{Ind}}_P^G(V \otimes \delta \otimes {\mathbb{C}}_{\lambda})
\]
 of the group $G=O(n+1,1)$ in terms of cohomological parabolic induction.  

\subsubsection{Cohomological parabolic induction $A_{\mathfrak{q}}(\lambda)={\mathcal{R}}_{\mathfrak{q}}^S({\mathbb{C}}_{\lambda+\rho({\mathfrak{u}})})$}
\label{subsec:Aqgeneral}
We fix some notation of
\index{B}{cohomologicalparabolicinduction@cohomological parabolic induction|textbf}
 cohomological parabolic induction.  
A basic reference is Vogan \cite{Vogan81} and Knapp--Vogan \cite{KV}.  
We begin with a {\it{connected}} real reductive Lie group $G$.  
Let $K$ be a maximal compact subgroup, 
 and $\theta$ the corresponding Cartan involution.  
Given an element $X \in {\mathfrak{k}}$, 
 the complexified Lie algebra
 ${\mathfrak{g}}_{\mathbb{C}}={\operatorname{Lie}}(G) \otimes_{\mathbb{R}}
{\mathbb{C}}$ is decomposed into the eigenspaces
 of $\sqrt{-1}{\operatorname{ad}}(X)$, 
 and we write 
\[
   {\mathfrak{g}}_{\mathbb{C}}
   ={\mathfrak{u}}_- + {\mathfrak{l}}_{\mathbb{C}} + {\mathfrak{u}}
\]
 for the sum of the eigenspaces 
 with negative, zero, 
 and positive eigenvalues.  
Then ${\mathfrak{q}}:={\mathfrak{l}}_{\mathbb{C}}+{\mathfrak{u}}$
 is a $\theta$-stable parabolic subalgebra 
 with Levi subgroup 
\begin{equation}
\label{eqn:LeviLq}
   L =\{g \in G: {\operatorname{Ad}}(g) {\mathfrak{q}}={\mathfrak{q}}\}.  
\end{equation}
The homogeneous space $G/L$ is endowed
 with a $G$-invariant complex manifold structure 
 with holomorphic cotangent bundle $G \times_L {\mathfrak{u}}$.  
As an algebraic analogue of Dolbeault cohomology groups
 for $G$-equivariant holomorphic vector bundle over $G/L$, 
 Zuckerman introduced a cohomological parabolic induction functor
 ${\mathcal{R}}_{\mathfrak{q}}^j(\cdot \otimes {\mathbb{C}}
_{\rho({\mathfrak{u}})})$ ($j \in {\mathbb{N}}$) from the category
 of $({\mathfrak{l}}, L \cap K)$-modules
 to the category of $({\mathfrak{g}}, K)$-modules. 
We adopt here the normalization of the cohomological parabolic induction
 ${\mathcal{R}}_{{\mathfrak {q}}}^{j}$ from a $\theta$-stable parabolic subalgebra
 ${\mathfrak {q}}={\mathfrak {l}}_{\mathbb{C}}+{\mathfrak {u}}$
 so that the ${\mathfrak{Z}}({\mathfrak {g}})$-infinitesimal character
 of the $({\mathfrak{g}},K)$-module 
 ${\mathcal{R}}_{{\mathfrak {q}}}^{j}(F)$ equals
\[
\text{
 the ${\mathfrak{Z}}({\mathfrak {l}})$-infinitesimal character
 of the ${\mathfrak {l}}$-module $F$
}
\]
 modulo the Weyl group 
 via the Harish-Chandra isomorphism.  

We note that if $F'$ is an $({\mathfrak{l}}, L \cap K)$-module
 then $F:=F' \otimes {\mathbb{C}}_{\rho({\mathfrak{u}})}$
 may not be defined 
 as an $({\mathfrak{l}}, L \cap K)$-module, 
 but can be defined as a module of the metaplectic covering group of $L$.  
When $F$ satisfies a positivity condition called 
\index{B}{goodrange@good range}
 \lq\lq{good range of parameters}\rq\rq, 
 the cohomology ${\mathcal{R}}_{\mathfrak{q}}^j(F)$ concentrates
 on the degree 
\[
   S:=\dim_{\mathbb{C}}({\mathfrak{u}} \cap {\mathfrak{k}}_{\mathbb{C}}).  
\]
For a one-dimensional representation $F$, 
 we also use another convention
 \lq\lq{$A_{\mathfrak{q}}(\lambda)$}\rq\rq.  
Following the normalization of Vogan--Zuckerman \cite{VZ}, 
 we set
\[
  A_{\mathfrak{q}}(\lambda):={\mathcal{R}}_{\mathfrak{q}}^S({\mathbb{C}}
_{\lambda+\rho({\mathfrak{u}})})
\]
for a one-dimensional representation ${\mathbb{C}}_{\lambda}$
 of $L$.  
In particular,
 we set
\[
   A_{\mathfrak{q}}:=A_{\mathfrak{q}}(0)
                   ={\mathcal{R}}_{\mathfrak{q}}^S({\mathbb{C}}_{\rho({\mathfrak{u}})}), 
\]
 which is an irreducible $({\mathfrak{g}},K)$-module
 with the same ${\mathfrak{Z}}({\mathfrak{g}})$-infinitesimal character $\rho$
 as that of the trivial one-dimensional representation ${\bf{1}}$ of $G$.

Similar notation will be used for disconnected groups $G$.  
For a character $\chi$ of the component group $G/G_0$, 
 we have an isomorphism of $({\mathfrak{g}},K)$-modules:
\[
   (A_{\mathfrak{q}})_{\chi}:=A_{\mathfrak{q}} \otimes \chi
   \simeq {\mathcal{R}}_{\mathfrak{q}}^S(\chi\otimes {\mathbb{C}}_{\rho({\mathfrak{u}})}).  
\]

\subsubsection{$\theta$-stable parabolic subalgebra ${\mathfrak{q}}_i$
 for $G=O(n+1,1)$}
\label{subsec:qi}
We apply the general theory
 reviewed in Section \ref{subsec:Aqgeneral}
 to the group $G=O(n+1,1)$.  
For this, 
we set up some notation 
 for $\theta$-stable parabolic subalgebra ${\mathfrak{q}}_i$
 and ${\mathfrak{q}}_{\frac{n+1}{2}}^{\pm}$
 of ${\mathfrak{g}}_{\mathbb{C}} = {\operatorname{Lie}}(G) \otimes_{\mathbb{R}} {\mathbb{C}}
\simeq {\mathfrak{o}}(n+2,{\mathbb{C}})$
 as follows.

We take a Cartan subalgebra ${\mathfrak{t}}^c$ of ${\mathfrak{k}}$, 
 and extend it to a fundamental Cartan subalgebra
 ${\mathfrak{h}}={\mathfrak{t}}^c + {\mathfrak{a}}^c$.  
If $n$ is odd then ${\mathfrak{a}}^c=\{0\}$.  
Choose the standard coordinates
 $\{f_k: 1 \le k \le [\frac n 2]+1\}$
 on ${\mathfrak{h}}_{\mathbb{C}}^{\ast}$
 such that the root system of ${\mathfrak{g}}$
 and ${\mathfrak{k}}$ are given by 
\begin{align*}
\Delta({\mathfrak{g}}_{\mathbb{C}}, {\mathfrak{h}}_{\mathbb{C}})
=&\{\pm(f_i \pm f_j) : 1 \le i < j \le [\frac n2]+1\}
\\
&\left(\cup \{\pm f_{\ell} : 1 \le \ell \le [\frac n2]+1\}
 \quad(\text{$n$: odd})\right), 
\\
\Delta({\mathfrak{k}}_{\mathbb{C}}, {\mathfrak{t}}_{\mathbb{C}})
=&\{\pm(f_i \pm f_j) : 1 \le i < j \le [\frac {n+1}2] \}
\\
&\left(\cup \{\pm f_{\ell} : 1 \le \ell \le [\frac {n+1}2]\}
 \quad(\text{$n$: even})\right).  
\end{align*}

For $1 \le i \le [\frac {n+1}2]$, 
 we define elements of ${\mathfrak{t}}_{\mathbb{C}}^{\ast}$
 by 
\begin{align*}
\mu_i:=&\sum_{k=1}^i (\frac n 2+1-k) f_k, 
\\
\mu_i^-:=& \mu_i- (n +2-2i) f_i. 
\end{align*}
It is convenient to set 
 $\mu_0=\mu_0^-=0$.  
(We shall use $\mu_i^-$
 only when we consider the identity component group 
 $G_0=SO_0(n+1,1)$
 with $n$ odd and when $n+1=2i$
 for later arguments.)
Let $\langle \, , \, \rangle$
 be the standard bilinear form
 on ${\mathfrak {h}}_{\mathbb{C}}^{\ast} \simeq {\mathbb{C}}^{[\frac n2 ]+1}$. 

\begin{definition}
[$\theta$-stable parabolic subalgebra ${\mathfrak{q}}_i$]
\label{def:qi}
For $0 \le i \le [\frac {n+1} 2]$, 
 we define $\theta$-stable parabolic subalgebras 
\index{A}{qi@${\mathfrak{q}}_i$, $\theta$-stable parabolic subalgebra|textbf}
$
{\mathfrak {q}}_i\equiv {\mathfrak {q}}_i^+
$
$
=({\mathfrak {l}}_i)_{\mathbb{C}} +{\mathfrak {u}}_i
$
and 
$
{\mathfrak {q}}_i^-
=({\mathfrak {l}}_i)_{\mathbb{C}} +{\mathfrak {u}}_i^-
$
 in ${\mathfrak {g}}_{\mathbb{C}}=\operatorname{Lie}(G) \otimes_{\mathbb{R}} {\mathbb{C}}$
by the condition
 that ${\mathfrak {q}}_i$ and ${\mathfrak {q}}_i^-$
 contain the fundamental Cartan subalgebra 
 ${\mathfrak {h}}$ 
 and that their nilradicals ${\mathfrak {u}}_i$ and ${\mathfrak {u}}_i^-$
 are given respectively
 by 
\begin{align*}
\Delta({\mathfrak {u}}_i, {\mathfrak {h}}_{\mathbb{C}})
=&\{\alpha \in \Delta({\mathfrak {g}}_{\mathbb{C}}, {\mathfrak {h}}_{\mathbb{C}}): \langle \alpha, \mu_i \rangle >0\}, 
\\
\Delta({\mathfrak {u}}_i^-, {\mathfrak {h}}_{\mathbb{C}})
=&\{\alpha \in \Delta({\mathfrak {g}}_{\mathbb{C}}, {\mathfrak {h}}_{\mathbb{C}}): \langle \alpha, \mu_i^- \rangle >0\}.  
\end{align*}
Then the Levi subgroup of 
 ${\mathfrak {q}}={\mathfrak {q}}_i$ and ${\mathfrak{q}}_i^-$
 is given by 
\index{A}{LSi@$L_i=SO(2)^i \times O(n-2i+1,1)$, Levi subgroup of ${\mathfrak{q}}_i$|textbf}
\begin{equation}
\label{eqn:Li}
  L_i := N_G({\mathfrak{q}})
      \equiv \{g \in G:{\operatorname{Ad}}(g){\mathfrak{q}}={\mathfrak{q}}\}
      \simeq SO(2)^i \times O(n-2i+1,1).  
\end{equation}
We note that $L_i$ is not in the Harish-Chandra class
 if $n$ is even, 
 as is the case $G=O(n+1,1)$.  
\end{definition}

If we write $\rho({\mathfrak{u}}_i)$ and $\rho({\mathfrak{u}}_i^-)$
 for half the sum of roots in ${\mathfrak{u}}_i$ and ${\mathfrak{u}}_i^-$, 
respectively, 
 then 
\[
  \rho({\mathfrak{u}}_i)=\mu_i
\quad
\text{and}
\quad
  \rho({\mathfrak{u}}_i^-)=\mu_i^-.  
\]
We suppress the superscript $+$ for ${\mathfrak{q}}_i^+$
 except for the case $n+1=2i$.  
For later purpose,
 we compare the following three groups with the same Lie algebras:
\begin{alignat}{3}
\label{eqn:GGG}
 G_0 =&SO_0(n+1,1) &&\hookrightarrow \overline G=SO(n+1,1) &&\hookrightarrow G=O(n+1,1)
\intertext{with maximal compact subgroups}
 K_0 =& SO(n+1)   &&\hookrightarrow \overline K=O(n+1) &&\hookrightarrow K=O(n+1)\times O(1).  
\nonumber
\end{alignat}
\begin{lemma}
\label{lem:conjqi}
\begin{enumerate}
\item[{\rm{(1)}}]
A complete system of the $K_0$-conjugacy classes
 of $\theta$-stable parabolic subalgebras 
 of ${\mathfrak{g}}_{\mathbb{C}}$
 with Levi subgroup $L_i$
 \eqref{eqn:Li}
is given by 
\index{A}{qiapm@${\mathfrak{q}}_i^{\pm}$|textbf}
\begin{alignat*}{2}
  &\{{\mathfrak{q}}_i\}
  &&\text{for $0 \le i < [\frac {n+1}2]$}, 
\\
&\{{\mathfrak{q}}_{\frac {n+1}{2}}^+, {\mathfrak{q}}_{\frac{n+1}2}^-\}
\quad
  &&\text{for $i = \frac {n+1}2$ $($$n$:odd$)$}.  
\end{alignat*}
\item[{\rm{(2)}}]
The $\theta$-stable parabolic subalgebra ${\mathfrak{q}}_i$
 with the property \eqref{eqn:Li} is unique up to conjugation
 by the disconnected group $\overline K$
 (and therefore, also by $K$)
 for all $i$ $(0 \le i \le [\frac {n+1}2])$.  
\end{enumerate}
\end{lemma}

We also make the following two observations:
\begin{lemma}
\label{lem:Litorus}
$L_i$ is compact
 if and only if $n$ is odd and $2i=n+1$.  
In this case, 
 $L_i \simeq SO(2)^{\frac {n+1}2} \times O(1)$.  
\end{lemma}
\begin{lemma}
\label{lem:GoverLi}
The inclusion maps \eqref{eqn:GGG} induce
 the following inclusion and bijection:
\[
   G_0/N_{G_0}({\mathfrak{q}}_i)
   \hookrightarrow
   \overline G/N_{\overline G}({\mathfrak{q}}_i)
   \overset \sim \rightarrow
   G/N_{G}({\mathfrak{q}}_i)
   =
   G/L_i
\]
for all $i$ $(0 \le i \le [\frac{n+1}2])$.  
The first inclusion is bijective
 if $n+1 \ne 2i$.  
\end{lemma}
The second bijection is reflected by the irreducibility of the $G$-module
 $\Pi_{\ell, \delta}$ when restricted to the subgroup 
 $\overline G=SO(n+1,1)$, 
 see Proposition \ref{prop:161648} (1)
 in Appendix II.  

Lemmas \ref{lem:conjqi} and \ref{lem:Litorus} yield the following 
 (well-known) representation theoretic results:
\begin{proposition}
\label{prop:disc}
\begin{enumerate}
\item[{\rm{(1)}}]
$G$ (or $\overline G$, $G_0$) admits a discrete series representation
 if and only if $n$ is odd.  
\item[{\rm{(2)}}]
Suppose $n$ is odd.  
Then there exists only one discrete series representation of $\overline G$
 for each regular integral infinitesimal character;
there exist exactly two discrete series representations of $G$
 (also of $G_0$)
 for each regular integral infinitesimal character.  
\end{enumerate}
\end{proposition}
For $n=2m-1$ in the second statement of Proposition \ref{prop:disc}, 
 we note the following properties
 for the three groups $G \supset \overline G \supset G_0$: 
\begin{enumerate}
\item[$\bullet$]
$L_m \simeq SO(2)^m \times O(1)$
 has two connected components;
\item[$\bullet$]
$L_m \cap \overline G=L_m \cap G_0$ are connected;
\item[$\bullet$]
${\mathfrak{q}}_m^+$ and ${\mathfrak{q}}_m^-$ are not conjugate
 by $G_0$;
 they are conjugate by $\overline G$ or $G$.  
\end{enumerate}
See \cite[Thm.~3 (0)]{KMemoirs92} 
 for results in a more general setting
 of the indefinite orthogonal group $O(p,q)$.

For $\nu =(\nu_1, \cdots, \nu_i) \in {\mathbb{Z}}^i$, 
 $\mu \in \Lambda^+([\frac n 2]-i+1)$, 
 and $a, b \in \{\pm \}$, 
 we consider an irreducible finite-dimensional $L_i$-module
\[
  \Kirredrep {O(n-2i+1,1)}{\mu}_{a, b} \otimes {\mathbb{C}}_{\nu}
\]
and define an admissible smooth representation of $G$ of moderate growth, 
 to be denoted by 
\[
  \Rq {\nu_1, \cdots, \nu_i}{\mu_1, \cdots, \mu_{[\frac n 2]-i+1}}{a, b}, 
\]
 whose underlying $({\mathfrak{g}},K)$-module
 is given by the cohomological parabolic induction
\index{A}{RzqS@${\mathcal{R}}_{\mathfrak{q}}^S$, cohomological parabolic induction|textbf}
\begin{equation}
\label{eqn:Rqmunu}
  {\mathcal{R}}_{{\mathfrak {q}}_i}^{S_i}(\Kirredrep{O(n-2i+1,1)}{\mu}_{a, b} \otimes {\mathbb{C}}_{\nu+\rho({\mathfrak{u}}_i)})
\end{equation}
 of degree $S_i$,  
where we set
\begin{equation}
\label{eqn:cohSi}
   S_i:= \dim_{\mathbb{C}} ({\mathfrak {u}}_i \cap {\mathfrak {k}}_{\mathbb{C}})       =i(n-i).  
\end{equation}
We note
 that if $i=0$
 then $\Rq {}{\mu_1, \cdots, \mu_{[\frac n 2]+1}}{a, b}$ is
 finite-dimensional.  

\begin{definition}
[$\theta$-stable parameter]
\label{def:thetapara}
We call $\Rq{\nu_1, \cdots, \nu_i}{\mu_1, \cdots, \mu_{[\frac n 2]-i+1}}{a, b}$
 the 
\index{B}{1thetastableparameter@$\theta$-stable parameter|textbf}
{\it{$\theta$-stable parameter}}
 of the representation \eqref{eqn:Rqmunu}.  
\end{definition}

If the $\theta$-stable parameter 
 of a representation $\Pi$ of $G$ is given by 
\[
   \Rq{\nu_1, \cdots, \nu_i}{\mu_1, \cdots, \mu_{[\frac n 2]-i+1}}{a, b}, 
\]
 then that of $\Pi \otimes \chi_{c d}$ for $c, d \in \{\pm\}$ is given by
\begin{equation}
\label{eqn:thetatensor}
\Rq{\nu_1, \cdots, \nu_i}{\mu_1, \cdots, \mu_{[\frac n 2]-i+1}}{ac,bd}.
\end{equation}
The 
\index{B}{infinitesimalcharacter@infinitesimal character}
${\mathfrak{Z}}_G({\mathfrak {g}})$-infinitesimal character
 of $\Rq{\nu_1, \cdots, \nu_i}{\mu_1, \cdots, \mu_{[\frac n 2]-i+1}}{a, b}$
 is given by
\[
  (\nu_1, \cdots, \nu_i, \mu_1, \cdots, \mu_{[\frac n 2]-i+1})
+
  (\frac n 2, \frac n 2-1, \cdots, \frac n 2-[\frac n 2]).  
\]
In particular,
 the $G$-module 
\[
  (\underbrace{0, \cdots,0}_{i} || \underbrace{0, \cdots,0}_{[\frac n 2]-i+1})_{a, b}
\]
has the trivial infinitesimal character $\rho_G$.  
In this case
 we shall write 
\index{A}{Aqlmd@$A_{\mathfrak{q}}(\lambda)$|textbf}
\index{A}{Aqlmdab@$A_{{\mathfrak{q}}_i}$|textbf}
\index{A}{Aqlmdac@$(A_{{\mathfrak{q}}})_{\pm,\pm}$|textbf}
\begin{equation}
\label{eqn:AqRq}
   (A_{\mathfrak{q}_i})_{a, b}
   :={\mathcal{R}}_{\mathfrak{q}_i}^{S_i}(\chi_{ab}\otimes {\mathbb{C}}_{\rho({\mathfrak{u}}_i)})
\end{equation}
 for its underlying $({\mathfrak{g}},K)$-module, 
 see Proposition \ref{prop:161655} below.

Sometimes we suppress the subscript $+,+$
 and write simply $A_{\mathfrak{q}_i}$
 to denote the $({\mathfrak{g}},K)$-module $(A_{\mathfrak{q}_i})_{+,+}$.  
\begin{remark}
\label{rem:goodrange}
\begin{enumerate}
\item[(1)]
(good range)\enspace
The irreducible finite-dimensional representation 
 $\Kirredrep{O(n-2i+1,1)}{\mu}_{a, b} \otimes {\mathbb{C}}_{\nu+\rho({\mathfrak{u}})}$ 
 of the metaplectic cover of $L_i$ is 
 in the 
\index{B}{goodrange@good range}
 {\it{good range}} 
 with respect to the $\theta$-stable parabolic subalgebra
 ${\mathfrak{q}_i}$
 (see \cite[Def.~0.49]{KV} for the definition)
 if and only if
\[ \nu_1 \ge \nu_2 \ge \cdots \ge \nu_i \ge \mu_1.\]
In this case, 
 the $({\mathfrak{g}},K)$-module \eqref{eqn:Rqmunu}
 is nonzero and irreducible,
 and therefore
 $\Rq{\nu_1, \cdots,\nu_i}{\mu_1, \cdots, \mu_{[\frac n 2]-i+1}}{a, b}$
 is a nonzero irreducible $G$-module.  
For the description
 of the Hasse sequence
 (Theorem \ref{thm:IndV} below), 
 we need only the parameter in the good range.  
\item[(2)]
(weakly fair range)\enspace
If $\mu=(0,\cdots,0)$, 
then the $({\mathfrak{g}},K)$-module \eqref{eqn:Rqmunu} reduces to
\[
  A_{\mathfrak{q}_i}(\nu)_{a, b}
   :={\mathcal{R}}_{\mathfrak{q}_i}^{S_i}
     (\chi_{ab}\otimes {\mathbb{C}}_{\nu+\rho({\mathfrak{u}}_i)})  
\]
cohomologically induced from the one-dimensional representation
 $\chi_{ab}\otimes {\mathbb{C}}_{\nu+\rho({\mathfrak{u}})}$.  
We note that $\chi_{ab}\otimes {\mathbb{C}}_{\nu+\rho({\mathfrak{u}})}$
 is in the 
\index{B}{weaklyfairrange@weakly fair range}
{\it{weakly fair range}} 
 with respect to $\mathfrak{q}_i$
 (see \cite[Def.~0.52]{KV} for the definition)
 if and only if 
\begin{equation}
\label{eqn:wfair}
\nu_{1} + \frac{n}{2} 
\ge 
\nu_{2} + \frac{n}{2} -1
\ge 
\cdots
\ge
\nu_{i} + \frac{n}{2} -i+1
\ge 
0.  
\end{equation}
In this case the $({\mathfrak{g}}, K)$-module
 $A_{\mathfrak{q}_i}(\nu)_{a, b}$
 may or may not vanish.  
See \cite[Thm.~3]{KMemoirs92}
 for the conditions on $\nu \in {\mathbb{Z}}^i$
 in the weakly fair range
 that assure the nonvanishing 
 and the irreducibility
 of $A_{\mathfrak{q}_i}^{S_i}({\mathbb{C}}_{\nu})_{a, b}$.  
We shall see in Section \ref{subsec:Aqsing} 
 that the underlying $({\mathfrak{g}}, K)$-modules
 of singular complementary series representations
 are isomorphic to these modules.  
\end{enumerate}
\end{remark}

\subsubsection{Irreducible representations $\Pi_{\ell, \delta}$ and $(A_{\mathfrak{q}_i})_{\pm,\pm}$}

In this subsection,
 we give a description
 of the underlying $({\mathfrak{g}}, K)$-modules
 of the subquotients $\Pi_{\ell,\delta}$
 of the principal series representation
 of the disconnected group $G=O(n+1,1)$
 in terms of the cohomologically parabolic induced modules
 $(A_{\mathfrak{q}_i})_{\pm,\pm}$.  

We recall from \eqref{eqn:Pild} 
 the definition of the irreducible representations 
\index{A}{1Piidelta@$\Pi_{i,\delta}$, irreducible representations of $G$}
 $\Pi_{\ell,\delta}$
 $(0 \le \ell \le n+1$, $\delta = \pm)$
 of $G=O(n+1,1)$.  
The set 
\[
\{\Pi_{\ell,\delta}: 0 \le \ell \le n+1, \delta=\pm\}
\]
 exhausts irreducible admissible representations
 of moderate growth 
 having ${\mathfrak{Z}}_G({\mathfrak{g}})$-infinitesimal
 character $\rho_G$, 
 see Theorem \ref{thm:LNM20} (2).  
Their underlying $({\mathfrak{g}},K)$-modules
 $(\Pi_{\ell,\delta})_K$
 can be given by cohomologically parabolic induced modules as follows.  
\begin{proposition}
\label{prop:161655}
For $0 \le i \le [\frac {n+1}2]$, 
 let 
\index{A}{qi@${\mathfrak{q}}_i$, $\theta$-stable parabolic subalgebra}
${\mathfrak{q}}_i$ be the $\theta$-stable parabolic subalgebras
 with the Levi subgroup $L_i \simeq SO(2)^i \times O(n-2i+1,1)$
 as in Definition \ref{def:qi}.  
\begin{enumerate}
\item[{\rm{(1)}}]
The underlying $({\mathfrak {g}},K)$-modules
 of the irreducible $G$-modules $\Pi_{\ell, \delta}$
 $(0 \le \ell \le n+1$, $\delta \in \{\pm\})$
 are given by the cohomological parabolic induction
 as follows:
\index{A}{Aqlmd@$A_{\mathfrak{q}}(\lambda)$}
\begin{alignat*}{2}
  (\Pi_{i,+})_K \simeq& (A_{\mathfrak {q}_i})_{+,+} 
&& \supset \Exterior^i({\mathbb{C}}^{n+1}) \boxtimes {\bf{1}}, 
\\
   (\Pi_{i,-})_K \simeq & (A_{\mathfrak {q}_i})_{+,-}
&& \supset \Exterior^i({\mathbb{C}}^{n+1}) \boxtimes {\operatorname{sgn}}, 
\\
   (\Pi_{n+1-i,+})_K \simeq & (A_{\mathfrak {q}_i})_{-,+}
&& \supset \Exterior^{n+1-i}({\mathbb{C}}^{n+1}) \boxtimes {\bf{1}}, 
\\
  (\Pi_{n+1-i,-})_K \simeq & (A_{\mathfrak {q}_i})_{-,-}
&& \supset \Exterior^{n+1-i}({\mathbb{C}}^{n+1}) 
\boxtimes {\operatorname{sgn}}.  
\end{alignat*}
For later purpose, 
 we also indicated their 
\index{B}{minimalKtype@minimal $K$-type}
 minimal $K$-types
 in the right column
 (see Theorem \ref{thm:LNM20} (3)).  
\item[{\rm{(2)}}]
If $n$ is even or if $2i\ne n+1$, 
 then the four $({\mathfrak{g}},K)$-modules
 $(A_{\mathfrak{q}_i})_{a,b}$ $(a,b \in \{\pm\})$ are not isomorphic to each other.

If $2i=n+1$, 
 then there are isomorphisms
\[
   (A_{\mathfrak {q}_{\frac{n+1}{2}}})_{+,+}
 \simeq 
   (A_{\mathfrak {q}_{\frac{n+1}{2}}})_{-,+}
\quad
\text{and}
\quad 
   (A_{\mathfrak {q}_{\frac{n+1}{2}}})_{+,-}
\simeq
   (A_{\mathfrak {q}_{\frac{n+1}{2}}})_{-,-}
\]
 as $({\mathfrak {g}}, K)$-modules
 for the disconnected group $O(n+1,1)$.  
\end{enumerate}
\end{proposition}
Thus the left-hand sides of the formul\ae\ 
 in Proposition \ref{prop:161655} (1) have overlaps
 when $n$ is odd and $i=\frac{n+1}{2}$.  
In fact, 
the Levi part in this case is of the form 
 $L_{\frac{n+1}{2}} \simeq SO(2)^{\frac{n+1}{2}} \times O(0,1)$, 
 and $\chi_{-+} \simeq {\bf{1}}$
 and $\chi_{+-} \simeq \chi_{--}$
 as $O(0,1)$-modules.  

\subsubsection{Irreducible representations with nonzero $({\mathfrak{g}},K)$-cohomologies}
\label{subsec:gKnonzero}

In this section,
 we prove Theorem \ref{thm:LNM20} (9)
 on the classification 
 of irreducible unitary representations of $G=O(n+1,1)$
 with nonzero $({\mathfrak{g}},K)$-cohomologies.  
We have already seen in Lemma \ref{lem:172145}
 that $H^{\ast}({\mathfrak{g}},K;(\Pi_{\ell,\delta})_K) \ne \{0\}$
 for all $0 \le \ell \le n+1$
 and $\delta \in \{\pm\}$.  
Hence the proof of Theorem \ref{thm:LNM20} (9) will be completed
 by showing the following.  
\begin{proposition}
\label{prop:gKq}
Let $\Pi$ be an irreducible unitary representation
 of $G=O(n+1,1)$
 such that 
 $H^{\ast}({\mathfrak{g}},K;\Pi_K) \ne \{0\}$.  
Then the smooth representation $\Pi^{\infty}$ is isomorphic 
 to 
\index{A}{1Piidelta@$\Pi_{i,\delta}$, irreducible representations of $G$}
 $\Pi_{\ell,\delta}$
 (see \eqref{eqn:Pild})
 for some $0 \le \ell \le n+1$ and $\delta \in \{\pm\}$.  
\end{proposition}

\begin{proof}
We begin with representations of the identity component
 $G_0=SO_0(n+1,1)$.  
In this case,
 we write $A_{\mathfrak{q}}^0$
 by putting superscript 0
 to denote the $({\mathfrak{g}},K_0)$-module
 which is cohomologically induced from 
 the trivial one-dimensional representation
 of a $\theta$-stable parabolic subalgabra ${\mathfrak{q}}$.

By a theorem of Vogan and Zuckerman \cite{VZ},
 any irreducible unitary representation $\Pi^0$
 of $G_0$ 
 with $H^{\ast}({\mathfrak{g}},K_0;(\Pi^0)_{K_0}) \ne \{0\}$
 is of the form
 $(\Pi^0)_{K_0} \simeq A_{\mathfrak{q}}^0$
 for some $\theta$-stable parabolic subalgebra 
 ${\mathfrak{q}}$ in ${\mathfrak {g}}_{\mathbb{C}}$.  
We recall from Definition \ref{def:qi}
 that ${\mathfrak{q}}_i$
 ($0 \le i < \frac{n+1}{2}$)
 and ${\mathfrak{q}}_i^{\pm}$
 ($i = \frac{n+1}{2}$)
 are $\theta$-stable parabolic subalgebras
 such that the Levi subgroup 
 $N_{G_0}({\mathfrak{q}}_i)$
 (or $N_{G_0}({\mathfrak{q}}_i^{\pm})$)
 are isomorphic to 
 $SO(2)^i \times SO_0(n-2i+1,1)$.  
They exhaust all $\theta$-stable parabolic subalgebras 
 up to inner automorphisms
 and up to cocompact Levi factors,
 namely,
 there exists $0 \le i \le [\frac{n+1}{2}]$
 such that 
\[
 {\mathfrak{q}}_i \subset {\mathfrak{q}}
\quad
\text{and}
\quad
  \text{$N_{G_0}({\mathfrak{q}})/N_{G_0}({\mathfrak{q}}_i)$ is compact}
\]
 if we take a conjugation
 of ${\mathfrak{q}}$
 by an element of $G_0$.  
(For $i = \frac{n+1}{2}$, 
 ${\mathfrak{q}}_i$ is considered as either ${\mathfrak{q}}_i^+$
 or ${\mathfrak{q}}_i^-$.)
Then we have a $({\mathfrak{g}},K_0)$-isomorphism
\begin{equation}
\label{eqn:ZVconn}
   (\Pi^0)_{K_0}
   \simeq
   A_{\mathfrak{q}}^0
   \simeq
\begin{cases}
 A_{\mathfrak{q}_i}^0 \quad 
&\text{if $2i < n+1$,}
\\
  A_{\mathfrak{q}_i^+}^0 \quad \text{or}\quad A_{{\mathfrak{q}}_i^-}^0
  \quad 
& \text{if $2i=n+1$.}
\end{cases}
\end{equation}
Now we consider an irreducible unitary representation $\Pi$
 of the {\it{disconnected}} group $G=O(n+1,1)$
 such that 
 $H^{\ast}({\mathfrak{g}},K;\Pi_K) \ne \{0\}$.  
The assumption implies $H^{\ast}({\mathfrak{g}},K_0;\Pi_K) \ne \{0\}$, 
 and therefore there exists a $G_0$-irreducible submodule
 $\Pi^0$ of the restriction $\Pi|_{G_0}$
 such that 
 $H^{\ast}({\mathfrak{g}},K_0;(\Pi^0)_{K_0}) \ne \{0\}$.  
By the reciprocity,
 the underlying $({\mathfrak{g}},K)$-module $\Pi_K$ must be an irreducible summand
 in the induced representation
\[
   {\operatorname{ind}}_{{\mathfrak{g}},K_0}^{{\mathfrak{g}},K}
   ((\Pi^0)_{K_0}).  
\]
It follows from \eqref{eqn:ZVconn}
 and from Proposition \ref{prop:161655} (2)
 that
\[
{\operatorname{ind}}_{{\mathfrak{g}},K_0}^{{\mathfrak{g}},K}
   ((\Pi^0)_{K_0})
\simeq
\begin{cases}
\underset {a, b \in \{\pm\}}{\bigoplus} (A_{\mathfrak{q}_i})_{a, b}
&\text{if $2i < n+1$,}
\\
(A_{{\mathfrak{q}}_{\frac{n+1}{2}}})_{+,+}
\oplus
(A_{\mathfrak{q}_{\frac{n+1}2}})_{-,-}
\quad 
&\text{if $2i=n+1$.}
\end{cases}
\]
Thus Proposition \ref{prop:gKq} follows from Proposition \ref{prop:161655} (1). \end{proof}

\subsubsection{Description of subquotients 
 in $I_{\delta}(V,\lambda)$}
We use the $\theta$-stable parameter
 for the description
 of irreducible subquotients
 of the principal series representations 
 $I_{\delta}(V,\lambda)$
 of $G=O(n+1,1)$ with regular integral infinitesimal character.  
\begin{theorem}
\label{thm:IndV}
Suppose $V \in \widehat{O(n)}$ and $\lambda \in {\mathbb{Z}} \setminus S(V)$.  
Let $i := i(\lambda, V)$ be the height as in Lemma \ref{lem:ilsig}.  
We write $V=\Kirredrep{O(n)}{\sigma}_{\varepsilon}$
 with $\sigma =(\sigma_1, \cdots, \sigma_{[\frac n2]})\in \Lambda^+([\frac n 2])$ and $\varepsilon \in \{\pm\}$.  
Let $\delta \in \{\pm\}$.  
\begin{enumerate}
\item[{\rm{(1)}}]
Suppose $\lambda \ge \frac n 2$.  
Then $\frac n 2 \le i \le n$.  

If $i \ne \frac n 2$, 
then we have the following nonsplit exact sequence
 of $G$-modules of moderate growth:
\begin{align}
  0 \to& \Rq{\sigma_1-1, \cdots, \sigma_{n-i}-1, \lambda-i}{\sigma_{n-i+1}, \cdots, \sigma_{[\frac n 2]}}{\varepsilon , -\delta \varepsilon}
\notag
\\
\to& I_{\delta}(V, \lambda)
\notag
\\
\to& 
\Rq{\sigma_1-1, \cdots, \sigma_{n-i}-1}{\lambda-i, \sigma_{n-i+1}, \cdots, \sigma_{[\frac n 2]}}{\varepsilon , \delta \varepsilon}
\to 0.  
\label{eqn:161688new}
\end{align}

\item[{\rm{(2)}}]
Suppose $\lambda \le \frac n 2$.   
Then $0 \le i \le \frac n 2$.  

If $i \ne \frac n 2$, 
then we have the following nonsplit exact sequence
 of $G$-modules
 of moderate growth:
\begin{align}
  0 \to& \Rq{\sigma_1-1, \cdots, \sigma_i-1}{i-\lambda, \sigma_{i+1}, \cdots, \sigma_{[\frac n 2]}}{\varepsilon , \delta \varepsilon}
\notag
\\
\to& I_{\delta}(V, \lambda)
\notag
\\
\to&
\Rq{\sigma_1-1, \cdots, \sigma_i-1, i-\lambda}{\sigma_{i+1}, \cdots, \sigma_{[\frac n 2]}}{\varepsilon , -\delta \varepsilon}
\to 0.  
\label{eqn:161616new}
\end{align}

\item[{\rm{(3)}}]
Suppose $i=\frac n 2$, 
 or equivalently, 
 suppose that $n$ is even and $\sigma_{\frac n 2} >|\lambda-\frac n 2|$.

If $\lambda \ne \frac n 2$, 
then $\lambda \in S_Y(V)$
 (see \eqref{eqn:RSYV}).  
In this case,
 $I_{\delta}(V, \lambda)$ is irreducible
 and we have a $G$-isomorphism:
\[I_{\delta}(V, \lambda)
    \simeq 
\Rq{\sigma_1-1, \cdots, \sigma_{\frac n 2}-1}{|\lambda-\frac n2|}{a, b}
\]
whenever $a,b \in \{\pm\}$ satisfies $a b = \delta$.

If $\lambda = \frac n 2$,
 then $I_{\delta}(V, \lambda)$ splits
 into the direct sum of two irreducible representations of $G$:
\begin{equation}
\label{eqn:IVYdeco}
I_{\delta}(V, \lambda)
    \simeq 
\bigoplus_{a,b \in \{\pm\}, ab =\delta}
\Rq{\sigma_1-1, \cdots, \sigma_{\frac n 2}-1}{0}{a, b}.  
\end{equation}
\end{enumerate}
\end{theorem}

\begin{remark}
\label{rem:IVlisom}
In Theorem \ref{thm:IndV} (3), 
 we have a $G$-isomorphism 
\[
   I_{\delta}({\Kirredrep{O(n)}{\sigma}}_+, \lambda)
   \simeq
   I_{\delta}({\Kirredrep{O(n)}{\sigma}}_-, \lambda)
\quad
\text{ for each $\delta= \pm$.  }
\]
In fact,
 by Lemma \ref{lem:isig}, 
 $i(\lambda,V)=\frac n 2$ implies 
 that $V$ is of type Y, 
 hence there is an $O(N)$-isomorphism
$\Kirredrep {O(N)} {\sigma}_+ \simeq \Kirredrep {O(N)} {\sigma}_-$
 by Lemma \ref{lem:typeY}.

Moreover, 
 the restriction of each irreducible summand in \eqref{eqn:IVYdeco}
 to the special orthogonal group $SO(n+1,1)$ is irreducible
 (see Lemma \ref{lem:171523} (1) in Appendix II).  
\end{remark}

\subsubsection{Proof of Theorem \ref{thm:IndV}}

\begin{proof}
[Sketch of the proof of Theorem \ref{thm:IndV}]
If the ${\mathfrak{Z}}_G({\mathfrak{g}})$-infinitesimal character
 of the principal series representation
 $I_{\delta}(\Kirredrep {O(n)}{\sigma}_{\varepsilon}, \lambda)$ 
 is $\rho_G$, 
 then Theorem \ref{thm:IndV} is a reformulation of Theorem \ref{thm:LNM20}
 in terms of $\theta$-stable parameters.  
This is done in Proposition \ref{prop:IndV0} below.

The general case is derived from the above case by the translation principle, 
 see Theorems \ref{thm:181104} and \ref{thm:20180904}, 
 and also the argument there
 ({\it{e.g.}}, Lemma \ref{lem:Rqtensor}) in Appendix III.  
\end{proof}
Suppose $V=\Exterior^i({\mathbb{C}}^n)$.  
By Example \ref{ex:exteriorpm}, 
the principal series representation 
 $I_{\delta}(i,\lambda)={\operatorname{Ind}}_P^G(\Exterior^i({\mathbb{C}}^n) \otimes \delta \otimes {\mathbb{C}}_{\lambda})$
 is expressed as follows.  
\begin{lemma}
\label{lem:171276}
There are natural $G$-isomorphisms:
\begin{alignat*}{2}
I_{\delta}(\ell,\ell) 
\simeq & 
I_{\delta}(\Kirredrep {O(n)}{1^\ell, 0^{[\frac n 2]-\ell}}_+,\ell)
&& \text{if $\ell \le \frac n 2$}, 
\\
I_{\delta}(\ell,\ell) 
\simeq & 
I_{\delta}(\Kirredrep {O(n)}{1^{n-\ell}, 0^{\ell-[\frac {n+1} 2]}}_-,\ell)
\qquad
&& \text{if $\ell \ge \frac n 2$}.  
\end{alignat*}
\end{lemma}

\begin{proposition}
\label{prop:IndV0}
Suppose $0 \le \ell \le \frac n 2$.  
Then Theorem \ref{thm:IndV} holds for $\lambda=\ell$ 
 and $\sigma = (1^\ell, 0^{[\frac n 2]-\ell}) \in \Lambda^+([\frac n 2])$.  
\end{proposition}
\begin{proof}
By Theorem \ref{thm:LNM20} (1), 
 we have an exact sequence of $G$-modules
\[
   0 \to \Pi_{\ell,\delta} \to I_{\delta}(\ell,\ell) \to \Pi_{\ell+1,-\delta} \to 0, 
\]
which does not split as far as $\ell \ne \frac n 2$.  
By Proposition \ref{prop:161655}, 
 this yields an exact sequence of $({\mathfrak {g}}, K)$-modules:
\[
  0 \to (A_{{\mathfrak {q}}_\ell})_{+,\delta}
    \to I_{\delta}(\ell,\ell)_K
    \to (A_{{\mathfrak {q}}_{\ell+1}})_{+,-\delta}
    \to 0.  
\]
By Lemma \ref{lem:171276} and the definition of $\theta$-stable parameters,
 this exact sequence can be written as
\[
  0 \to 
  \Rq{0^\ell}{0^{[\frac n 2]-\ell+1}}{+,\delta}
  \to 
  I_{\delta}(\Kirredrep{O(n)}{\sigma}_+,\ell)
  \to
  \Rq{0^{\ell+1}}{0^{[\frac n 2]-\ell}}{+,-\delta}
  \to 0.  
\]
Since the height of $(\Kirredrep{O(n)}{\sigma}_+,\ell)=(\Exterior^{\ell}({\mathbb{C}}^n),\ell)$
 is given by $i(\Exterior^{\ell}({\mathbb{C}}^n),\ell)=\ell$.  
$i(\ell,\sigma)=\ell$ by Example \ref{ex:isig}, 
 we get Proposition \ref{prop:IndV0} from
 Lemma \ref{lem:IVchi} and \eqref{eqn:thetatensor}.  
\end{proof}

\subsection{Hasse sequence in terms of $\theta$-stable parameters}
\label{subsec:HasseAq}

This section gives a description of the Hasse sequence 
 (Definition-Theorem \ref{def:UHasse})
 and the standard sequence (Definition \ref{def:Hasse})
 in terms of $\theta$-stable parameters.

We set $m:=[\frac{n+1}{2}]$, 
 namely $n=2m-1$ or $2m$.  
Let $F$ be an irreducible finite-dimensional representation of $G=O(n+1,1)$,
 and $U_i \equiv U_i(F)$
 ($0 \le i \le [\frac{n+1}{2}]$) be the Hasse sequence with $U_0 \simeq F$.  
We write
\[
   F=\Kirredrep{O(n+1,1)}{s_0, \cdots, s_{[\frac n 2]}}_{a, b}
\]
as in Lemma \ref{lem:161612} (2).  
\begin{theorem}
\label{thm:171471}
Let $n=2m$
 and $0 \le i \le m$.  
\begin{enumerate}
\item[{\rm{(1)}}]
{\rm{(Hasse sequence)}}\enspace
$U_i(F) \simeq \Rq{s_0, \cdots, s_{i-1}}{s_i, \cdots, s_m}
{a, (-1)^{i-s_i} b}$.  
\item[{\rm{(2)}}]
{\rm{(standard sequence)}}\enspace
$U_i(F) \otimes \chi_{+-}^i \simeq 
\Rq{s_0, \cdots, s_{i-1}}{s_i, \cdots, s_m}{a, (-1)^{s_i}b}$. 
\end{enumerate}
\end{theorem}
\begin{proof}
(1)\enspace
We begin with the case $a=b=+$.  
Let $s:=(s_0,\cdots,s_m,0^{m+1}) \in \Lambda^+(2m+2)$.  
As in \eqref{eqn:sin=2m} of Section \ref{subsec:Hasseps}, 
 we define $s^{(\ell)} \in \Lambda^+(2m)$
 for $0 \le \ell \le m$.  
Then by Theorem \ref{thm:171426}, 
 there is an injective $G$-homomorphism
\[
  U_\ell(F) 
  \hookrightarrow 
  I_{(-1)^{\ell-s_{\ell}}}(\Kirredrep {O(n)}{s^{(\ell)}}, \ell-s_\ell).  
\]
The $O(n)$-module $\Kirredrep{O(n)}{s^{(\ell)}}$ is of type I
 (Definition \ref{def:type}), 
 and we have 
\[
  i(\Kirredrep{O(n)}{s^{(\ell)}}, \ell-s_\ell))=\ell
\]
 with the notation of Lemma \ref{lem:ilsig}.  

By Theorem \ref{thm:IndV}, 
 we get the theorem for $a=b=+$ case.  
The general case follows from the case
 $(a,b)=(+,+)$
 by the tensoring argument given in Proposition \ref{prop:HStensor}.  
\par\noindent
(2)\enspace
The second statement follows from Definition \ref{def:Hasse}
 and \eqref{eqn:thetatensor}.  
\end{proof}

The case $n$ odd is given similarly as follows.  
\begin{theorem}
\label{thm:171471b}
Let $n=2m-1$, 
 and $0 \le i \le m-1$.  
\begin{enumerate}
\item[{\rm{(1)}}]
{\rm{(Hasse sequence)}}\enspace
$U_i(F) \simeq \Rq{s_0, \cdots, s_{i-1}}{s_i, \cdots, s_{m-1}}
{a, (-1)^{i-s_i} b}$.  
\item[{\rm{(2)}}]
{\rm{(standard sequence)}}\enspace
$U_i(F) \otimes \chi_{+-}^i \simeq 
\Rq{s_0, \cdots, s_{i-1}}{s_i, \cdots, s_{m-1}}{a, (-1)^{s_i}b}$. 
\end{enumerate}
\end{theorem}
\begin{proof}
\begin{enumerate}
\item[(1)]
We begin with the case $a=b=+$.  
Let $s:=(s_0,\cdots,s_{m-1},0^{m+1}) \in \Lambda^+(2m+1)$.  
As in \eqref{eqn:sin=2m-1}, 
 we define $s^{(\ell)} \in \Lambda^+(2m-1)$
 for $0 \le \ell \le m-1$.  
Then by Theorem \ref{thm:171425}, 
\[
  U_\ell(F) 
  \subset 
  I_{(-1)^{\ell-s_{\ell}}}(\Kirredrep {O(n)}{s^{(\ell)}}, \ell-s_\ell).  
\]
The $O(n)$-module $\Kirredrep{O(n)}{s^{(\ell)}}$ is of type I, 
 and we obtain
\[
  i(\Kirredrep {O(n)}{s^{(\ell)}}, \ell-s_\ell)=\ell
\]
 with the notation of Lemma \ref{lem:ilsig}.  

By Theorem \ref{thm:IndV}, 
 we get the theorem for $a=b=+$ case.  
The general case follows from the case
 $(a,b)=(+,+)$
 by the tensoring argument given in Proposition \ref{prop:HStensor}.  
\item[(2)]
The second statement follows from Definition \ref{def:Hasse}
 and \eqref{eqn:thetatensor}.  
\end{enumerate}
\end{proof}

\subsection{Singular integral case}
\label{subsec:Aqsing}
\index{B}{singularintegralinfinitesimalcharacter@singular integral infinitesimal character}
We end this chapter 
 with cohomologically induced representations
 with {\it{singular}} parameter,
 and give a description 
 of complementary series representations
 with integral parameter
 (see Section \ref{subsec:singcomp})
 in terms of $\theta$-stable parameters.

For $0 \le r \le [\frac{n+1}{2}]$,
 we define ${\mathfrak {q}}_r = ({\mathfrak {l}}_r)_{\mathbb{C}}
 +{\mathfrak {u}}_r$
 to be the $\theta$-stable parabolic subalgebra
 with Levi subgroups
$
  L_r \simeq SO(2)^r \times O(n+1-2r, 1)
$
 in $G=O(n+1,1)$ as in Definition \ref{def:qi}.  
We set $S_r = r(n-r)$.

For $\nu=(\nu_1, \cdots, \nu_r) \in {\mathbb{Z}}^r \simeq (SO(2)^r)\hspace{1mm}
{\widehat{}}\hphantom{ii}$
 and $a, b \in \{\pm\}$, 
 we consider the underlying $({\mathfrak{g}},K)$-modules
 of the admissible smooth representations of $G$:
\[
  \Rq {\nu_1, \cdots, \nu_r}{\underbrace{0,\cdots,0}_{[\frac n 2]-r+1}}{a, b}, 
\]
 namely, 
the following $({\mathfrak{g}},K)$-modules
\[
  A_{\mathfrak{q}_r}(\nu)_{a, b}
   ={\mathcal{R}}_{\mathfrak{q}_r}^{S_r}
     (\chi_{ab} \otimes {\mathbb{C}}_{\nu+\rho({\mathfrak{u}})})
   \simeq 
   {\mathcal{R}}_{\mathfrak{q}_r}^{S_r}
   ({\mathbb{C}}_{\nu+\rho({\mathfrak{u}}_r)}) \otimes \chi_{ab}, 
\]
which are cohomologically induced from the one-dimensional representations
 ${\mathbb{C}}_{\nu} \boxtimes \chi_{a b}$
 of the Levi subgroup $L_r$, 
 see Remark \ref{rem:goodrange}
 for our normalization about \lq\lq{$\rho$-shift}\rq\rq.

Sometimes we suppress the subscript $+,+$
 and write simply $A_{{\mathfrak{q}}_{r}}(\nu)$
 for $A_{{\mathfrak{q}}_{r}}(\nu)_{+,+}$.

For a description of singular integral complementary series representations
 $I_{\delta}(i,s)$
 in terms of $\theta$-stable parameters,
 we need to treat the parameter $\nu$ outside the good range 
 (\cite[Def.~0.49]{KV})
 relative to the $\theta$-stable parabolic subalgebra ${\mathfrak {q}}_r$
 with $r = i+1$
 (see Theorem \ref{thm:compint} below), 
 for which the general theory 
 about the nonvanishing and irreducibility
 ({\it{e.g}}. \cite[Thm.~0.50]{KV})
 does not apply.  
For instance, 
the condition on the parameter $\nu$
 for which $A_{\mathfrak{q}_r}(\nu) \ne 0$
 is usually very complicated 
 when $\nu$ wanders outside the
\index{B}{goodrange@good range}
 good range.  
In our setting, 
 we use the following results from \cite{KMemoirs92}:
\begin{fact}
\label{fact:Memo92}
Let $0 \le r \le [\frac{n+1}{2}]$, 
 and $\mathfrak{q}_r$ be the $\theta$-stable parabolic subalgebra
 as defined in Definition \ref{def:qi}.  
Suppose that $\nu=(\nu_1, \cdots, \nu_r) \in {\mathbb{Z}}^r$
 satisfies the weakly fair condition \eqref{eqn:wfair}
 relative to ${\mathfrak{q}}_r$.  
Let $a, b \in \{\pm\}$.  
\begin{enumerate}
\item[{\rm{(1)}}]
The $G$-module $\Rq{\nu_1, \cdots, \nu_r}{0,\cdots,0}{a, b}$
 is nonzero
 if and only if
 $r=1$ or $\nu_{r-1} \ge -1$.  
\item[{\rm{(2)}}]
If the condition (1) is fulfilled, 
then $\Rq{\nu_1, \cdots, \nu_r}{0,\cdots,0}{a, b}$ is irreducible 
 and unitarizable.  
\end{enumerate}
\end{fact}
\begin{proof}
This is a special case
 of \cite[Thm.~3]{KMemoirs92}
 for the indefinite orthogonal group $O(p,q)$
 with $(p,q)=(n+1,1)$
 with the notation there.  
\end{proof}

Assume now $\nu_1 = \cdots = \nu_{r-1}=0$.  
Then the necessary and sufficient condition 
 for the parameter
 $\nu=(0,\cdots,0,\nu_r) \in {\mathbb{Z}}^r$ 
 to be in the weakly fair range 
 but outside the good range
 is given by 
\[
   \nu_r \in \{-1,-2, \cdots, r-1-[\frac n 2]\}.  
\]
In this case, 
 the $G$-module $\Rq{0,\cdots,0,\nu_r}{0,\cdots,0}{a, b}$ is nonzero,
 irreducible, 
 and unitarizable for $a,b \in \{\pm\}$
 as is seen in Fact \ref{fact:Memo92}. 
It turns out that these very parameters give rise to the complementary series representations
 with integral parameter stated in Section \ref{subsec:singcomp} as follows:
\begin{theorem}
\label{thm:compint}
Let $0 \le i \le [\frac n 2]-1$.  
Then the underlying $({\mathfrak{g}},K)$-modules
 of the
\index{B}{complementaryseries@complementary series representation}
 complementary series representations $I_\pm(i,s)$ and $I_\pm(n-i,s)$
 with integral parameter
 $s \in \{i+1,i+2, \cdots, [\frac n 2]\}$ are given by
\begin{align*}
  I_+(i,s)_K
  &\simeq
  A_{{\mathfrak{q}}_{i+1}}(0,\cdots,0,s-i)_{+,+};
\\
  I_-(i,s)_K
  &\simeq
  A_{{\mathfrak{q}}_{i+1}}(0,\cdots,0,s-i)_{+,-};
\\
  I_+(n-i,s)_K
  &\simeq
  A_{{\mathfrak{q}}_{i+1}}(0,\cdots,0,s-i)_{-,-};
\\
  I_-(n-i,s)_K
  &\simeq
  A_{{\mathfrak{q}}_{i+1}}(0,\cdots,0,s-i)_{-,+}.  
\end{align*}
Hence,
 their smooth globalizations are described
 by $\theta$-stable parameters as follows:
\begin{alignat*}{2}
   &I_+(i,s) 
   &&\simeq \Rq{\underbrace{0,\cdots,0,s-i}_{i+1}}{\underbrace{0,\cdots,0}_{[\frac n 2]-i}}{+,+};
\\
&I_-(i,s) 
   &&\simeq \Rq{0,\cdots,0,s-i}{0,\cdots,0}{+,-};
\\
&I_+(n-i,s) 
   &&\simeq \Rq{0,\cdots,0,s-i}{0,\cdots,0}{-,-};
\\
&I_-(n-i,s) 
   &&\simeq \Rq{0,\cdots,0,s-i}{0,\cdots,0}{-,+}.  
\end{alignat*} 
\end{theorem}

\newpage
\section{Appendix II: Restriction to $\overline G=SO(n+1,1)$}
\label{sec:SOrest}

So far we have been working with symmetry breaking for a pair of the orthogonal groups
 $(O(n+1,1), O(n,1))$.  
On the other hand, 
 the Gross--Prasad conjectures
 (Chapters \ref{sec:Gross-Prasad} and \ref{sec:conjecture})
 are formulated for special orthogonal groups
 rather than orthogonal groups.  
In this chapter, 
 we explain how to translate the results 
 for $(G, G')=(O(n+1,1), O(n,1))$
 to those for the pair of special orthogonal groups
\index{A}{GOindefinitespecial@$\overline{G}=SO(n+1,1)$|textbf}
 $(\overline{G},\overline{G'})=(SO(n+1,1), SO(n,1))$.  
A part of the results here 
({\it{e.g.}}, Theorem \ref{thm:SOmult}) was announced in \cite{sbonGP}.

In what follows, 
 we use a bar over representations of special orthogonal groups
 to distinguish them from those of orthogonal groups.

\subsection{Restriction of representations of $G=O(n+1,1)$ to $\overline G= SO(n+1,1)$}
\label{subsec:SO}

It is well-known that any irreducible admissible representation $\Pi$
 of a real reductive group $G$ is decomposed into the direct sum 
 of finitely many irreducible admissible representations of $\overline G$
if $\overline G$ is an open normal subgroup of $G$
 (see \cite[Chap.~II, Lem.~5.5]{BW}).  
In order to understand
 how the restriction $\Pi|_{\overline G}$ decomposes,
 we use the action
 of the quotient group $G/\overline G$
 on the ring ${\operatorname{End}}_{\overline G}(\Pi|_{\overline G}) =
 {\operatorname{Hom}}_{\overline G}(\Pi|_{\overline G}, \Pi|_{\overline G})$.

We apply this general observation to our setting 
where 
\[(G,\overline G)=(O(n+1,1), SO(n+1,1)).  
\]
In this case, 
 the quotient group $G/\overline G \simeq {\mathbb{Z}}/2 {\mathbb{Z}}$.  
With the notation \eqref{eqn:chiab} of the characters $\chi_{a b}$ of $G$, 
\index{A}{1chipmm@$\chi_{--}=\det$}
\[
   \{\chi_{++}, \chi_{--}\}=\{{\bf{1}}, \det\}
\]
 is the set
 of irreducible representations
 of $G=O(n+1,1)$
 which are trivial on $\overline G=S O(n+1,1)$.  
In other words,
 we have a direct sum decomposition as $G$-modules:
\[
   {\operatorname{Ind}}_{\overline G}^G \overline{\bf{1}}
  \simeq
   {\bf{1}} \oplus \det.  
\]
Then we have the following:
\begin{lemma}
\label{lem:LNM26}
Let $\Pi$ be a continuous representation of $G=O(n+1,1)$.  
Then there is a natural linear bijection:
\index{A}{1chipmpm@$\chi_{\pm\pm}$, one-dimensional representation of $O(n+1,1)$}
\[
  {\operatorname{End}}_{\overline{G}}(\Pi|_{\overline{G}})
\simeq
  \operatorname{Hom}_{G}(\Pi, \Pi)
  \oplus
  \operatorname{Hom}_{G}(\Pi, \Pi \otimes \det).  
\]
\end{lemma}

\begin{proof}
Clear from the following linear bijections:
\[
{\operatorname{End}}_{\overline G}(\Pi|_{\overline{G}})
\simeq 
{\operatorname{Hom}}_{G}(\Pi, {\operatorname{Ind}}_{\overline G}^{G}(\Pi|_{\overline{G}}))
\simeq 
{\operatorname{Hom}}_{G}(\Pi, \Pi \otimes {\operatorname{Ind}}_{\overline G}^{G}\overline{\bf{1}}).  
\]
\end{proof}

We examine the restriction
 of irreducible representations of $G$
 to the subgroup $\overline G$:

\begin{lemma}
\label{lem:171523}
Suppose that $\Pi$ is an irreducible admissible representation of $G=O(n+1,1)$.  
\begin{enumerate}
\item[{\rm{(1)}}]
If $\Pi \not \simeq \Pi \otimes \det$ as $G$-modules, 
 then the restriction $\Pi|_{\overline G}$ is irreducible.  
\item[{\rm{(2)}}]
If $\Pi \simeq \Pi \otimes \det$ as $G$-modules, 
then the restriction $\Pi|_{\overline G}$ is the direct sum 
of two irreducible admissible representations
of $\overline G$ that are not isomorphic to each other.  
\end{enumerate}
\end{lemma}
\begin{proof}
By Lemma \ref{lem:LNM26}, 
 we have 
\begin{align*}
\dim_{\mathbb{C}}{\operatorname{Hom}}_{\overline G}(\Pi|_{\overline G}, \Pi|_{\overline G})
=& \dim_{\mathbb{C}} {\operatorname{Hom}}_{G}(\Pi, \Pi) 
+ \dim_{\mathbb{C}}{\operatorname{Hom}}_{G}(\Pi, \Pi \otimes \det)
\\
=& 
\begin{cases}
1 \quad &\text{if $\Pi \not \simeq \Pi \otimes \det$, }
\\
2 \quad &\text{if $\Pi \simeq \Pi \otimes \det$.  }
\end{cases}
\end{align*}
Since the restriction $\Pi|_{\overline G}$ is the direct sum
 of irreducible admissible representations 
 of ${\overline G}$, 
 we may write the decomposition as 
\[
  \Pi|_{\overline G} \simeq \bigoplus_{j=1}^{N} m_j \overline{\Pi}_j,
\]
where $\overline{\Pi}_j$ are 
 (mutually inequivalent)
 irreducible admissible representations of $\overline G$
 and $m_j \in {\mathbb{N}}_+$
 denote the multiplicity of $\overline{\Pi}_j$ in $\Pi|_{\overline G}$
 for $1 \le j \le N$.  
By Schur's lemma, 
\[
  \dim_{\mathbb{C}} {\operatorname{End}}_{\overline{G}}
  (\Pi|_{\overline G}) = \sum_{j=1}^N m_j^2.  
\]
This is equal to 1 or 2 
 if and only if $N=m_1 =1$ or $N=2$ and $m_1 = m_2=1$, 
respectively.  
Hence we get  the conclusion.  
\end{proof}

\subsection{Restriction of principal series representation
 of $G=O(n+1,1)$ to $\overline G=SO(n+1,1)$ }
\label{subsec:psrest}
This section discusses the restriction of the principal series representation $I_{\delta}(V, \lambda)$
 of $G=O(n+1,1)$ to the normal subgroup $\overline G= SO(n+1,1)$ of index two.  
First of all, 
 we fix some notation 
 for principal series representations of $\overline G$.  
We set 
\index{A}{PLanglandsdecompbar@$\overline P = \overline M A N_+$|textbf}
$\overline P := P \cap \overline G$. 
Then $\overline P$ is a minimal parabolic subgroup of $\overline G$, 
 and its Langlands decomposition is given
 by $\overline P = \overline M A N_+$, 
 where 
\index{A}{M1primebar@$\overline M=SO(n) \times O(1)$}
\[
\overline M := M \cap \overline G
=
\{ 
\begin{pmatrix}
               \varepsilon & &
\\
                           & B & 
\\
                           & & \varepsilon
\end{pmatrix}
:
B \in SO(n), \varepsilon= \pm 1
\}
\simeq
SO(n) \times O(1)
\]
 is a subgroup of $M$ of index two.  
For an irreducible representation $(\overline\sigma, \overline V)$ of $SO(n)$, 
 $\delta \in \{\pm\}$, 
 and $\lambda \in {\mathbb{C}}$, 
 we denote by 
\index{A}{IdeltaVsbar@$\overline I_{\delta}(V, \lambda)$}
$\overline I_{\delta}(\overline V,\lambda)$
 the (unnormalized) induced representation
$
    \operatorname{Ind}_{\overline P}^{\overline G}
    (\overline V \otimes \delta \otimes {\mathbb{C}}_{\lambda})
$
 of $\overline G =SO(n+1,1)$.

Let us compare principal series representations of $G$
 regarded as $\overline G$-modules
 by restriction with principal series representations of $\overline G$.  
For this, 
 we suppose $V$ is an irreducible representation of $O(n)$, 
 $\delta \in \{\pm \}$, and $\lambda \in {\mathbb{C}}$, 
 and form a principal series representation $I_{\delta}(V,\lambda)$
 of $G=O(n+1,1)$.  
Then its restriction to the subgroup $\overline G=SO(n+1,1)$ is isomorphic to 
$
   \operatorname{Ind}_{\overline P}^{\overline G}
   (V|_{SO(n)} \otimes \delta \otimes {\mathbb{C}}_{\lambda})
$
 as a $\overline G$-module, 
 because the inclusion $\overline G \hookrightarrow G$
 induces an isomorphism
 $\overline G/\overline P \overset \sim \to G/P$.

Concerning the $SO(n)$-module $V|_{SO(n)}$, 
 we recall from Definition \ref{def:OSO} that $V \in \widehat {O(n)}$ is said to be 
\index{B}{typeX@type X, representation of ${O(N)}$\quad}
 of type X
or 
\index{B}{typeY@type Y, representation of ${O(N)}$\quad}
 of type Y
 according to whether $V$ is irreducible or reducible 
when restricted to $SO(n)$.  
In the latter case, 
$n$ is even (see Lemma \ref{lem:OSO}) 
 and $V$ is decomposed into the direct sum of two irreducible representations
 of $SO(n)$:
\begin{equation}
\label{eqn:Vpm}
   V= V^{(+)} \oplus V^{(-)}, 
\end{equation}
where $V^{(-)}$ is isomorphic to the contragredient representation 
 of $V^{(+)}$.  
Accordingly, 
 we have an isomorphism 
 as $\overline G$-modules:
\index{A}{IdeltaVsbarpm@$\overline I_{\delta}(V^{(\pm)}, \lambda)$|textbf}
\begin{equation}
\label{eqn:IVSOpm}
   I_{\delta}(V,\lambda)|_{\overline G}
   \simeq
   \begin{cases}
   \overline{I}_{\delta}(V,\lambda)
   & \text{if $V$ is of type X, }
\\
  \overline{I}_{\delta}(V^{(+)}, \lambda)
   \oplus
  \overline{I}_{\delta}(V^{(-)}, \lambda)
\quad
     & \text{if $V$ is of type Y.}  
   \end{cases}
\end{equation}
By using \eqref{eqn:IVSOpm}, 
 we obtain the structural results of the restriction 
 of the principal series representation ${I}_{\delta}(V,\lambda)$
 of $G=O(n+1,1)$
 to the subgroup $\overline G=SO(n+1,1)$
 and further to the identity component group 
\index{A}{GOindefinitespecialidentity@$G_0=SO_0(n+1,1)$, the identity component
 of $O(n+1,1)$\quad}
 $G_0=SO_0(n+1,1)$.

\subsubsection{Restriction $I_{\delta}(V,\lambda)|_{\overline G}$
 when $I_{\delta}(V,\lambda)$ is irreducible}
We begin with the case
 where $I_{\delta}(V,\lambda)$ is irreducible as a $G$-module.  
\begin{lemma}
\label{lem:170305}
Let $(\sigma,V) \in \widehat{O(n)}$, 
$\delta \in \{\pm \}$
 and $\lambda \in {\mathbb{C}}$.  
Suppose $I_{\delta}(V,\lambda)$ is irreducible 
 as a module of $G=O(n+1,1)$.  
\par\noindent
{\rm{(1)}}\enspace
Suppose $V$ is of type X.  
Then the following three conditions on $(\delta,V,\lambda)$
 are equivalent:
\begin{enumerate}
\item[{\rm{(i)}}]
$I_{\delta}(V,\lambda)$ is irreducible
 as a $G$-module;
\item[{\rm{(ii)}}]
The restriction $I_{\delta}(V,\lambda)|_{\overline G}$ is irreducible
 as a $\overline G$-module;
\item[{\rm{(iii)}}]
The restriction $I_{\delta}(V,\lambda)|_{G_0}$ is irreducible
 as a $G_0$-module.  
\end{enumerate}
{\rm{(2)}}\enspace
Suppose $V$ is of type Y.  
If $I_{\delta}(V,\lambda)$ is irreducible as a $G$-module, 
 then $I_{\delta}(V,\lambda)|_{\overline G}$ splits
 into the direct sum of two irreducible $\overline G$-modules
 that are not isomorphic to each other.  
In this case, 
 $n$ is even
 and we may write the irreducible decomposition
 of $V|_{SO(n)}$
 as in \eqref{eqn:Vpm}.  
Then there is a natural isomorphism 
\[
I_{\delta}(V,\lambda)|_{\overline G}
\simeq
 \overline{I}_{\delta}(V^{(+)},\lambda) \oplus \overline{I}_{\delta}(V^{(-)},\lambda)
\]
 as $\overline G$-modules.  
Moreover,
 both $\overline I_{\delta}(V^{(+)}, \lambda)$
 and $\overline I_{\delta}(V^{(-)}, \lambda)$ stay irreducible
 when restricted to $G_0$, 
 and they are not isomorphic to each other also as $G_0$-modules.  
\end{lemma}

\begin{proof}
We observe
 that the first factor of $M$ is isomorphic to $O(n)$, 
 whereas that of $M \cap \overline G$
 $(=\overline M)$
 and of $M \cap G_0$ is isomorphic to $S O(n)$.  
Since the crucial part is the restriction from 
 the Levi subgroup $M A$ of $G$
 to that of $\overline G$ or of $G_0$, 
 we focus on the restriction $G \downarrow \overline G$, 
 which involves the restriction of $V$
 with respect to the inclusion $O(n) \supset S O(n)$.  
The restriction $G \downarrow G_0$ can be analyzed
 similarly by using the four characters
 $\chi_{a b}$ ($a,b \in \{\pm\}$)
 instead of $\chi_{--}=\det$
 as in \cite[Chap.~2, Sect.~5]{KKP}.

{}From now on,
 we consider the restriction $G \downarrow \overline G$.  
We recall from Lemma \ref{lem:IVchi}
 the following isomorphism of $G$-modules:
\index{A}{1chipmm@$\chi_{--}=\det$}
\[I_{\delta}(V, \lambda) \otimes \chi_{--}
 \simeq
 I_{\delta}(V \otimes \det, \lambda).  
\]
\par\noindent
(1)\enspace
If $V$ is of type X, 
then $V \not \simeq V \otimes \det$ as $O(n)$-modules.  
In turn, 
 the $G$-modules $I_{\delta}(V,\lambda)$ and $I_{\delta}(V\otimes \det,\lambda)$
 are not isomorphic to each other, 
 because their $K$-structures are different
 by the Frobenius reciprocity
 and the branching rule for $O(n) \downarrow O(n-1)$
 (Fact \ref{fact:ONbranch}).  
Therefore, 
 $I_{\delta}(V,\lambda)|_{\overline G}$ is irreducible
 by Lemma \ref{lem:171523} (1).  
\par\noindent
(2)\enspace
If $V$ is of type Y, 
 then $V \otimes \det \simeq V$
 by Lemma \ref{lem:typeY}, 
 and therefore Lemma \ref{lem:171523} (2) concludes the first assertion.  
The remaining assertions are now clear.  
\end{proof}

\subsubsection{Restriction $I_{\delta}(V,\lambda)|_{\overline G}$
 when $V$ is of type Y}
We take a closer look at the case
 where $V \in \widehat {O(n)}$ is of type Y
 (Definition \ref{def:OSO}).  
This means that $n$ is even, 
 say $n=2m$, 
 and the representation $V$ is of the form
\[
   V = \Kirredrep{O(2m)}{\sigma_1, \cdots,\sigma_m}_{\varepsilon}
\]
with $\sigma_1 \ge \cdots \ge \sigma_m \ge 1$
 and $\varepsilon \in \{\pm\}$, 
 see Section \ref{subsec:fdimrep} for the notation.  
Then the restriction $V|_{S O(n)}$ decomposes
 as 
\[
  V|_{S O(n)}=V^{(+)} \oplus V^{(-)} 
\]
as in \eqref{eqn:Vpm}, 
 where the highest weights of the irreducible $SO(2m)$-modules
  $V^{(\pm)}$ are given by $(\sigma_1, \cdots,\sigma_{m-1}, \pm \sigma_m)$.  
We recall from Definition \ref{def:SVSYV}
 for the subsets $S(V)$ and $S_Y(V)$ of ${\mathbb{Z}}$.  

\begin{proposition}
\label{prop:irrIVpm}
Suppose $G=O(n+1,1)$ with $n=2m$
 and $(\sigma,V) \in \widehat {O(n)}$ is of type Y.  
Let $\delta \in \{\pm\}$.  
\begin{enumerate}
\item[{\rm{(1)}}]
The following four conditions on $\lambda  \in {\mathbb{C}}$ 
 are equivalent.  
\begin{enumerate}
\item[{\rm{(i)}}]
$\overline I_{\delta}(V^{(+)},\lambda)$ is reducible
 as a representation 
 of $\overline G=SO(n+1,1);$
\item[{\rm{(ii)}}]
$\overline I_{\delta}(V^{(-)},\lambda)$ is reducible
 as a $\overline G$-module$;$
\item[{\rm{(iii)}}]
$\pm(\lambda-m)\in {\mathbb{Z}}
 \setminus 
 (\{\sigma_j+m-j:j=1,\cdots,m\} \cup \{0,1,2,\cdots,\sigma_m-1\});$
\item[{\rm{(iv)}}]
$\lambda \in {\mathbb{Z}} \setminus
(S(V) \cup S_Y(V) \cup \{m\})$.  
\end{enumerate}
\item[{\rm{(2)}}]
Suppose that $\lambda$ satisfies
 one of (therefore any of)
 the above equivalent conditions.  
Then, 
 for $\varepsilon \in \{\pm\}$, 
 the principal series representation
 $\overline I_{\delta}(V^{(\varepsilon)}, \lambda)$ of $\overline G$
 has a unique $\overline G$-submodule, 
 to be denoted by $\overline I_{\delta}(V^{(\varepsilon)}, \lambda)^{\flat}$, 
 such that the quotient $\overline G$-module
\index{A}{IdeltaVsbarpmflat@$\overline I_{\delta}(V^{(\pm)}, \lambda)^{\flat}$|textbf}
\index{A}{IdeltaVsbarpms@$\overline I_{\delta}(V^{(\pm)}, \lambda)^{\sharp}$|textbf}
\[
   \overline I_{\delta}(V^{(\varepsilon)}, \lambda)^{\sharp}
  :=\overline I_{\delta}(V^{(\varepsilon)}, \lambda)
   /\overline I_{\delta}(V^{(\varepsilon)}, \lambda)^{\flat}
\]
 is irreducible.  
Moreover we have
\begin{align*}
   \overline I_{\delta}(V^{(+)}, \lambda)^{\flat}
 &\not\simeq 
 \overline I_{\delta}(V^{(-)}, \lambda)^{\flat}, 
\\
  \overline I_{\delta}(V^{(+)}, \lambda)^{\sharp}
  &\not\simeq 
  \overline I_{\delta}(V^{(-)}, \lambda)^{\sharp}
\end{align*}
as $\overline G$-modules.  
\end{enumerate}
\end{proposition}

\begin{proof}
Since $\overline G=SO(2m+1,1)$ is generated 
 by the identity component 
 $G_0=SO_0(2m+1,1)$
 and a central element $-I_{2m+2}$, 
 any irreducible $\overline G$-module remains
 irreducible 
 when restricted to the connected subgroup $G_0$.  
Then the equivalence (i) $\Leftrightarrow$ (iii)
(also (ii) $\Leftrightarrow$ (iii))
 and the last assertion 
 in Proposition \ref{prop:irrIVpm} follows from
 Hirai \cite{Hirai62}.  
See also Collingwood \cite[Lem.~4.4.1 and Thm.~5.2.1]{C}
 for the computation of $\tau$-invariants
 of irreducible representations
 and a graphic description of the socle filtrations
 of principal series representations.  
Finally the equivalence (iii) $\Leftrightarrow$ (iv)
 is immediate from the definitions
 \eqref{eqn:singint} and \eqref{eqn:SYsigma}
 of $S(V)$ and $S_Y(V)$, 
 respectively.

The last assertion about the $\overline G$-inequivalence follows from
 the Langlands theory \cite{La88}
 because $\operatorname{Re} \lambda \ne m$
 and $V^{(+)} \not \simeq V^{(-)}$
 as $SO(2m)$-modules.  
\end{proof}

In the following proposition,
 we treat the set
 of the parameters $\lambda$
 complementary to the one in Proposition \ref{prop:irrIVpm}.  
\begin{proposition}
\label{prop:180572}
Suppose $G=O(n+1,1)$ with $n=2m$
 and $V \in \widehat {O(n)}$ is of type Y.  
Let $\delta \in \{\pm\}$.  
Assume that $\overline I_{\delta}(V^{(\pm)},\lambda)$ are irreducible
 representations of $\overline G=SO(2m+1,1)$, 
 or equivalently, 
 assume that 
\[
   \lambda \in ({\mathbb{C}} \setminus {\mathbb{Z}}) \cup S(V) \cup S_Y(V) \cup \{m\}.  
\]
\begin{enumerate}
\item[{\rm{(1)}}]
The following two conditions on $\lambda \in {\mathbb{C}}$ 
 are equivalent:
\begin{enumerate}
\item[{\rm{(i)}}]
The two $\overline G$-modules $\overline I_{\delta}(V^{(+)},\lambda)$
 and $\overline I_{\delta}(V^{(-)},\lambda)$ are isomorphic
 to each other;
\item[{\rm{(ii)}}]
$\lambda =m$.  
\end{enumerate}

\item[{\rm{(2)}}]
If $\lambda=m$ then the principal series representation
 $I_{\delta}(V,\lambda)$ of $G$ is decomposed 
 into the direct sum of two irreducible representations
 of $G$.  
\item[{\rm{(3)}}]
If $\lambda \ne m$,  then $I_{\delta}(V,\lambda)$ is irreducible
 as a representation of $G$.  
\end{enumerate}
\end{proposition}

\begin{proof}
(1)\enspace
As in the proof of Proposition \ref{prop:irrIVpm} (2), 
 if $\operatorname{Re}\lambda \ne m$, 
 then the Langlands theory \cite{La88} implies
 $\overline I_{\delta}(V^{(+)},\lambda) \not \simeq \overline I_{\delta}(V^{(-)},\lambda)$
 because $V^{(+)} \not \simeq V^{(-)}$
 as $SO(2m)$-modules.

If $\operatorname{Re}\lambda = m$, 
 then $\overline I_{\delta}(V^{(\pm)},\lambda)$ are (smooth)
 irreducible tempered representations,
 and the equivalence
 (i) $\Leftrightarrow$ (ii) follows from Hirai \cite{Hirai62}.  
This would follow also from the general theory 
 of the \lq\lq{$R$-group}\rq\rq\ (Knapp--Zuckerman \cite{KZ}).  
\par\noindent
(2)\enspace
Since $\operatorname{Re}\lambda = m$ is the unitary axis
 of the principal series representation
 $I_{\delta}(V,\lambda)$ 
 in our normalization
 (Section \ref{subsec:smoothI}), 
 the $G$-module $I_{\delta}(V,\lambda)$ decomposes into the direct sum
 of irreducible $G$-modules,
 say, 
 $\Pi^{(1)}$, $\dots$, $\Pi^{(k)}$, 
 and then decomposes further into irreducible $\overline G$-modules
 when restricted to the subgroup 
 $\overline G=SO(2m+1,1)$.  
Therefore the cardinality $k$ of irreducible $G$-summands satisfies
 either $k=1$
 ({\it{i.e.}}, $I_{\delta}(V,\lambda)$ is $G$-irreducible)
 or $k=2$ because the summands $\overline I_{\delta}(V^{(\pm)},\lambda)$ 
 in \eqref{eqn:IVSOpm} are irreducible as $\overline G$-modules
 by assumption.  
Since $\overline I_{\delta}(V^{(+)},m)
 \simeq \overline I_{\delta}(V^{(-)},m)$
 by the first statement, 
 we conclude $k \ne 1$ 
 by Lemma \ref{lem:170305} (2).  
Thus the second statement is proved.  
\par\noindent
(3)\enspace
We prove 
 that $I_{\delta}(V,\lambda)$ is irreducible
 by {\it{reductio ad absurdum}}.  
Suppose there were an irreducible proper submodule $\Pi$
 of $I_{\delta}(V,\lambda)$.  
Then $\Pi$ would remain irreducible 
 when restricted to the subgroup $\overline G=SO(2m+1,1)$
 because the restriction $\Pi|_{\overline G}$ must be isomorphic 
 to one of the $\overline G$-irreducible summands
 $\overline I_{\delta}(V^{(\pm)},\lambda)$ 
 in \eqref{eqn:IVSOpm}.  
Then $\Pi \not \simeq \Pi \otimes \det$
 as $G$-modules
 by Lemma \ref{lem:LNM26}.  
Therefore the direct sum $\Pi \oplus (\Pi \otimes \det)$ would be a $G$-submodule
 of $I_{\delta}(V,\lambda)$
 because $I_{\delta}(V,\lambda) \simeq I_{\delta}(V,\lambda) \otimes \det$
 when $V$ is of type Y.  
In turn, 
 its restriction to the subgroup $\overline G$ would yield 
 an isomorphism 
 $\overline I_{\delta}(V^{(+)},\lambda) \simeq \overline I_{\delta}(V^{(-)},\lambda)$ of $\overline G$-modules,
 contradicting the statement (1)
 of the proposition.  
Hence $I_{\delta}(V,\lambda)$ must be irreducible.  
\end{proof}

Applying Propositions \ref{prop:irrIVpm} and \ref{prop:180572}
 to the middle exterior representation
 $\Exterior^m({\mathbb{C}}^{n})$ of $O(n)$
 when $n=2m$, 
 we obtain the following.
\begin{example}
\label{ex:Imm}
Let $G=O(n+1,1)$ with $n=2m$, 
 and $\delta \in \{\pm\}$.  
As in \eqref{eqn:ISOneven}, 
 we write $\overline I_{\delta}^{(\pm)}(m,\lambda)$
 for the $\overline G$-modules
 $\overline I_{\delta}(V^{(\pm)},\lambda)$
 when $V=\Exterior^m({\mathbb{C}}^{2m})$.  
\begin{enumerate}
\item[{\rm{(1)}}]
The $\overline G$-modules
 $\overline I_{\delta}^{(\pm)}(m,\lambda)$ are reducible
 if and only if $\lambda \in (-{\mathbb{N}}_+) \cup (n+{\mathbb{N}}_+)$.  
\item[{\rm{(2)}}]
The $G$-module $I_{\delta}(m,m)$ decomposes
 into a direct sum of two irreducible $G$-modules
 (see also Theorem \ref{thm:LNM20} (1)).  
\item[{\rm{(3)}}]
$I_{\delta}(m,\lambda)$ is irreducible
 if $\lambda \in {\mathbb{Z}}$ 
 satisfies $0 \le \lambda \le n$ $(=2m)$
 and $\lambda \ne m$.  
\end{enumerate}
\end{example}

We refer to Theorem \ref{thm:LNM20}
 (also to Example \ref{ex:irrIilmd})
 for the irreducibility condition
 of $I_{\delta}(i,\lambda)$
 for general $i$ $(0 \le i \le n)$;
to Theorem \ref{thm:irrIV}
 for that of $I_{\delta}(V,\lambda)$, 
 which will be proved in the next section.

\subsection{Proof of Theorem \ref{thm:irrIV}: 
 Irreducibility criterion of $I_{\delta}(V,\lambda)$}
\label{subsec:pfirrIV}

As an application of the results in the previous sections, 
 we give a proof of Theorem \ref{thm:irrIV}
 on the necessary and sufficient condition
 for the principal series representation 
 $I_{\delta}(V,\lambda)$
 of $G=O(n+1,1)$ to be irreducible.  
\begin{proof}
[Proof of Theorem \ref{thm:irrIV}]
Suppose first that $V$ is of type X
 (Definition \ref{def:OSO}). 
Then the restriction $V|_{SO(n)}$ is irreducible as an $SO(n)$-module,
 and $I_{\delta}(V,\lambda)$ is $G$-irreducible
 if and only if the restriction $I_{\delta}(V,\lambda)|_{G_0}$ is 
 $G_0$-irreducible 
 by Lemma \ref{lem:170305} (1).  
The latter condition was classified in Hirai \cite{Hirai62}, 
 which amounts to the condition that $\lambda \not \in {\mathbb{Z}}$
 or $\lambda \in S(V)$.  
Thus Theorem \ref{thm:irrIV} for $V$ of type X is proved.

Next suppose $V$ is of type Y.  
As in \eqref{eqn:Vpm}, 
 we write $V|_{SO(n)} \simeq V^{(+)} \oplus V^{(-)}$
 for the irreducible decomposition as $SO(n)$-modules.  
If $I_{\delta}(V,\lambda)$ is $G$-irreducible, 
then $\overline I_{\delta}(V^{(\pm)},\lambda)$ are $\overline G$-irreducible
 by Lemma \ref{lem:170305} (2).  
Then the condition (iv) in Proposition \ref{prop:irrIVpm} (1) implies
 that
\begin{equation}
\label{eqn:lmdSVm}
  \text{$\lambda \not \in {\mathbb{Z}}$ or $\lambda \in S(V) \cup S_Y(V) \cup \{m\}$}.  
\end{equation}
Conversely,
 under the condition \eqref{eqn:lmdSVm}, 
 Proposition \ref{prop:180572} tells
 that $I_{\delta}(V,\lambda)$ is irreducible
 if and only if $\lambda \not \in {\mathbb{Z}}$
or $\lambda \in S(V) \cup S_Y(V)$.  
Thus Theorem \ref{thm:irrIV} is proved
 also for $V$ of type Y.  
\end{proof}

\subsection{Socle filtration of $I_{\delta}(V,\lambda)$:
 Proof of Proposition \ref{prop:Xred}}
\label{subsec:Isub2}

In this section,
 we complete the proof of Proposition \ref{prop:Xred}
 about the socle filtration
 of the principal series representation $I_{\delta}(V,\lambda)$
 of $G=O(n+1,1)$
 when it is reducible and $\lambda \ne \frac n 2$
 by using the restriction to the subgroups 
 $\overline G=SO(n+1,1)$ or $G_0=SO(n+1,1)$.

We begin with the case 
 that $V \in \widehat{O(n)}$ is of type X
 (Definition \ref{def:OSO}).  

\begin{proof}[Proof of Proposition \ref{prop:Xred} when $V$ is of type X]
In this case,
 for any nonzero subquotient $\Pi$ 
 of the principal series representation $I_{\delta}(V,\lambda)$ of $G=O(n+1,1)$, 
 we have 
\[\Pi \not \simeq \Pi \otimes \det\]
as $G$-modules
 because their $K$-types are different
 by Proposition \ref{prop:KtypeIV}.  
In turn, 
 Lemma \ref{lem:171523} implies
 that $\Pi$ is irreducible as a $G$-module
 if and only if the restriction $\Pi|_{\overline G}$
 is irreducible.

For $n$ even, 
 the restriction $\Pi|_{G_0}$
 further to the identity component $G_0=SO_0(n+1,1)$ 
 is still irreducible
 because $\overline G=SO(n+1,1)$ is generated
 by $G_0$ 
 and a central element $-I_{n+2}$.  
Thus the assertion follows from the socle filtration
 of the principal series representation of $G_0$
 in Hirai \cite{Hirai62}.

For $n$ odd, 
 since the restriction $V|_{SO(n)}$ stays irreducible, 
 $I_{\delta}(V,\lambda)|_{G_0}$ is a principal series representation 
 of $G_0=SO_0(n+1,1)$.  
Therefore the restriction $\Pi|_{G_0}$ is a $G_0$-subquotient
 of a principal series representation
 of $G_0$, 
 of which the length of composition series is either 2 or 3
 by Hirai \cite{Hirai62}.  
Inspecting the $K$-structure of $I_{\delta}(V,\lambda)$ from 
 Proposition \ref{prop:KtypeIV} again 
 and the $K_0$-structure of subquotients
 of the principal series representation 
 of 
\index{A}{GOindefinitespecialidentity@$G_0=SO_0(n+1,1)$, the identity component
 of $O(n+1,1)$\quad}
$G_0=SO_0(n+1,1)$ in \cite{Hirai62}, 
 we see that the restriction $\Pi|_{G_0}$ is irreducible
 as a $G_0$-module
 if $\Pi$ is not (the smooth representation of)
 a discrete series representation, 
 whereas it is a sum 
 of two (holomorphic and anti-holomorphic) discrete series representations 
 of $G_0$
 if $\Pi$ is a discrete series representation.  
\end{proof}

Alternatively,
 one may reduce the proof for type X
 to the case $(V,\lambda)=(\Exterior^i({\mathbb{C}}^n), i)$
 by using the translation functor, 
 see Theorems \ref{thm:1807113} and \ref{thm:1808101} (1)
 in Appendix III.

As the above proof shows, 
 we obtain the restriction formula
 of irreducible subquotients
\index{A}{IdeltaVf@$I_{\delta}(V, \lambda)^{\flat}$}
 $I_{\delta}(V,\lambda)^{\flat}$ and 
\index{A}{IdeltaVs@$I_{\delta}(V, \lambda)^{\sharp}$}
$I_{\delta}(V,\lambda)^{\sharp}$
 of the $G$-module
 $I_{\delta}(V,\lambda)$
 (Proposition \ref{prop:Xred})
 to the normal subgroup $\overline G=SO(n+1,1)$ as follows.  
\begin{proposition}
\label{prop:20180906}
Suppose $V \in \widehat{O(n)}$ is of type X
 and $\lambda \in {\mathbb{Z}} \setminus S(V)$.  
Let $\delta \in \{\pm\}$.  
\begin{enumerate}
\item[{\rm{(1)}}]
The principal series representation 
 $\overline I_{\delta}(V,\lambda)$ of $\overline G$
 has a unique proper submodule,
 to be denoted by
\index{A}{IdeltaVsbarpmflat@$\overline I_{\delta}(V^{(\pm)}, \lambda)^{\flat}$|textbf}
 $\overline I_{\delta}(V,\lambda)^{\flat}$.  
In particular,
 the quotient $\overline G$-module
\index{A}{IdeltaVsbarpms@$\overline I_{\delta}(V^{(\pm)}, \lambda)^{\sharp}$|textbf}
$
   \overline I_{\delta}(V,\lambda)^{\sharp}
 :=\overline I_{\delta}(V,\lambda)/\overline I_{\delta}(V,\lambda)^{\flat}
$
 is irreducible.  
\item[{\rm{(2)}}]
The restriction of the irreducible $G$-modules
 $I_{\delta}(V,\lambda)^{\flat}$ and $I_{\delta}(V,\lambda)^{\sharp}$ 
 to the normal subgroup $\overline G$ is given by
\begin{align*}
 I_{\delta}(V,\lambda)^{\flat}|_{\overline G}
 & \simeq \overline I_{\delta}(V,\lambda)^{\flat}, 
\\
 I_{\delta}(V,\lambda)^{\sharp}|_{\overline G}
 & \simeq \overline I_{\delta}(V,\lambda)^{\sharp}.   
\end{align*}
\end{enumerate}
\end{proposition}
We end this section 
 with the restriction of $I_{\delta}(V,\lambda)^{\flat}$
 and $I_{\delta}(V,\lambda)^{\sharp}$
 to the subgroup $\overline G$
 when $V$ is of type Y:

\begin{proposition}
\label{prop:1808104}
Suppose $G=O(n+1,1)$ with $n=2m$.  
Assume that $V \in \widehat {O(n)}$
 is of type Y, 
 $\delta \in \{\pm\}$, 
 and $\lambda \in {\mathbb{Z}} \setminus (S(V) \cup S_Y(V) \cup \{m\})$.  
Then the restriction of $I_{\delta}(V,\lambda)^{\flat}$
 and $I_{\delta}(V,\lambda)^{\sharp}$
 to the normal subgroup $\overline G =S O(n+1,1)$
 decomposes into the direct sum 
 of two irreducible $\overline G$-modules: 
\begin{align*}
  I_{\delta}(V,\lambda)^{\flat}|_{\overline G}
  &\simeq
  \overline I_{\delta}(V^{(+)},\lambda)^{\flat}
  \oplus
  \overline I_{\delta}(V^{(-)},\lambda)^{\flat}, 
\\
  I_{\delta}(V,\lambda)^{\sharp}|_{\overline G}
  &\simeq
  \overline I_{\delta}(V^{(+)},\lambda)^{\sharp}
  \oplus
  \overline I_{\delta}(V^{(-)},\lambda)^{\sharp}, 
\end{align*}
where we recall from Proposition \ref{prop:irrIVpm}
 for the definition of the irreducible $\overline G$-modules
 $\overline I_{\delta}(V^{(\pm)},\lambda)^{\flat}$
 and $\overline I_{\delta}(V^{(\pm)},\lambda)^{\sharp}$.  
\end{proposition}

\begin{proof}
[Proof of Proposition \ref{prop:Xred} for $V$ of type Y]
By Lemma \ref{lem:171523} and by the structural results
 on $\overline G$-modules $\overline I_{\delta}(V^{(\pm)},\lambda)$
 in Proposition \ref{prop:irrIVpm}, 
 the proof is reduced to the following lemma. \end{proof}
 
\begin{lemma}
\label{lem:1808103}
Under the assumption of Proposition \ref{prop:1808104}, 
 any $G$-submodule $\Pi$ of $I_{\delta}(V,\lambda)$ satisfies
\begin{equation}
\label{eqn:Pidet}
  \Pi \simeq \Pi \otimes \det
\end{equation}
as $G$-modules.  
\end{lemma}

\begin{proof}
Since $V$ is of type Y, 
 $V \simeq V \otimes \det$
 as $O(n)$-modules,
 hence we have natural $G$-isomorphisms
\[
   I_{\delta}(V,\lambda) \simeq
   I_{\delta}(V,\lambda) \otimes \det
\]
by Lemma \ref{lem:IVchi}.  
We prove \eqref{eqn:Pidet}
 by {\it{reductio ad absurdum}}.  
Suppose that the $G$-module $\Pi$ were not isomorphic
 to $\Pi \otimes \det$.  
Then the direct sum representation $\Pi \oplus (\Pi \otimes \det)$
 would be a $G$-submodule
 of $I_{\delta}(V,\lambda)$.  
In turn, 
 the $\overline G$-module $\Pi|_{\overline G}$ would occur in 
$
I_{\delta}(V,\lambda)|_{\overline G}
\simeq \overline I_{\delta}(V^{(+)},\lambda)
\oplus \overline I_{\delta}(V^{(-)},\lambda)
$
at least twice. 
But this is impossible by Proposition \ref{prop:irrIVpm}.  
Thus Lemma \ref{lem:1808103} is proved.  
\end{proof}

\subsection{Restriction of $\Pi_{\ell,\delta}$
 to $SO(n+1,1)$}
\label{subsec:PiSO}
\index{A}{1Piidelta@$\Pi_{i,\delta}$, irreducible representations of $G$}
In this section 
 we treat the case 
 where $I_{\delta}(V, \lambda)$ is not irreducible
 as a $G$-module.  
We discuss the restriction of $G$-irreducible subquotients
 of $I_{\delta}(V, \lambda)$
 to the subgroup $\overline G=SO(n+1,1)$.

We focus on the case 
 when $(\sigma,V)$ is the exterior representation 
 on $V=\Exterior^i({\mathbb{C}}^{n})$.  
In particular,
 irreducible representations 
 that have
 the ${\mathfrak{Z}}_G({\mathfrak{g}})$-infinitesimal character $\rho$,  
 namely,
 the irreducible $G$-modules $\Pi_{\ell,\delta}$
 ($0 \le \ell \le n+1$, $\delta \in \{\pm\}$) arise as 
 $G$-irreducible subquotients
 of $I_{\delta}(V, \lambda)$.  
To be more precise,
 we recall from \eqref{eqn:Pild}
 that $\Pi_{\ell,\delta}$
 are the irreducible subrepresentations
 of $I_{\delta}(\ell,\ell)$
 for $0 \le \ell \le n$
 and coincidently those of $I_{-\delta}(\ell-1,\ell-1)$
 for $1 \le \ell \le n+1$.  
\begin{lemma}
\label{lem:Pilso}
For all $0 \le \ell \le n+1$
 and $\delta=\pm$, 
 the restriction of $\Pi_{\ell, \delta}$
 to the subgroup $\overline G = SO(n+1,1)$ stays irreducible.  
\end{lemma}
\begin{proof}
The restriction $\Pi_{\ell, \delta}|_{\overline G}$ is 
irreducible by the criterion in Lemma \ref{lem:171523} (1)
 because $\Pi_{\ell,\delta} \otimes \det \simeq \Pi_{n+1-\ell,-\delta}
\not \simeq \Pi_{\ell,\delta}$
 by Theorem \ref{thm:LNM20} (5).  
\end{proof}

We denote by $\overline{\Pi}_{\ell,\delta}$
 the restriction of the irreducible $G$-module ${\Pi}_{\ell,\delta}$
 ($0 \le \ell \le n+1$, $\delta=\pm$)
 to the subgroup ${\overline G} = SO(n+1,1)$.  
By a little abuse of notation,
 we write $\overline I _{\delta}(i, \lambda)$
 for the restriction of $I _{\delta}(i, \lambda)$,  
 to the subgroup $\overline G$.  
Then the $SO(n)$-isomorphism
 $\Exterior^i({\mathbb{C}}^n) \simeq \Exterior^{n-i}({\mathbb{C}}^n)$
 induces a $\overline G$-isomorphism
\[
   \overline I _{\delta}(i, \lambda)
   \simeq
   \overline I _{\delta}(n-i, \lambda).  
\]

Special attention is needed in the case
 when $n$ is even
 and $n=2i$.  
In this case,
 the $O(n)$-module $\Exterior^{i}({\mathbb{C}}^n)$ 
 is of type Y
 (see Example \ref{ex:2.1}), 
 and it splits into the direct sum
 of two irreducible $SO(n)$-modules:
\[
   \Exterior^{\frac n 2}({\mathbb{C}}^n) 
   \simeq 
   \Exterior^{\frac n 2}({\mathbb{C}}^n)^{(+)}
   \oplus
   \Exterior^{\frac n 2}({\mathbb{C}}^n)^{(-)}.  
\]
We set 
\[
   \overline I_{\delta}^{(\pm)}(\frac n 2, \lambda)
   :=
   {\operatorname{Ind}}_{\overline P}^{\overline G}
   (\Exterior^{\frac n 2}({\mathbb{C}}^n)^{(\pm)}
    \otimes \delta 
    \otimes {\mathbb{C}}_{\lambda}).  
\]
As in \eqref{eqn:IVSOpm}, 
 the restriction $I_{\delta}(\frac n 2, \lambda)|_{\overline G}$
 is the direct sum of two $\overline G$-modules:
\begin{equation}
\label{eqn:ISOneven}
   \overline I_{\delta}(\frac n 2, \lambda)
   \simeq
   \overline I_{\delta}^{(+)}(\frac n 2, \lambda)
   \oplus
   \overline{I}_{\delta}^{(-)}({\frac n 2}, \lambda)
\quad
\text{for all $\lambda \in {\mathbb{C}}$.}
\end{equation}

If $I_{\delta}(\frac n 2, \lambda)$ is $G$-irreducible, 
 then Lemma \ref{lem:170305} (2) tells
 that the representations $\overline {I}_{\delta}^{(\pm)}(\frac n 2, \lambda)$
 of the subgroup $\overline G$
 are irreducible,
 and that they are not isomorphic to each other.

On the other hand, 
 if $\lambda=i$ $(=\frac n 2)$, 
 then the principal series representation $I_{\delta}(i,\lambda)$ is not irreducible as a $G$-module 
 but splits into the direct sum
 of two irreducible $G$-modules
 (see Theorem \ref{thm:LNM20} (1)):
\[
   I_{\delta}(\frac n 2,\frac n 2) \simeq \Pi_{\frac n 2+1,-\delta} \oplus \Pi_{\frac n 2,\delta}, 
\]
which are not isomorphic to each other.  
Moreover,
 the tensor product with $\chi_{--}$
 switches $\Pi_{\frac n 2+1,-\delta}$ and $\Pi_{\frac n 2,\delta}$
 (Theorem \ref{thm:LNM20} (5)).  
Hence we have a 
$
   \overline G$-isomorphism $\overline \Pi_{\frac n 2+1,-\delta}
   \simeq \overline\Pi_{\frac n 2,\delta}
$,
 which are $\overline G$-irreducible by Lemma \ref{lem:171523} (1).  
Therefore, 
 for $n$ even, 
 we have the following isomorphisms as $\overline G$-modules:
\begin{equation}
\label{eqn:161733}
\overline{\Pi}_{\frac n 2, \delta}
\simeq
\overline{\Pi}_{\frac n 2+1, -\delta}
\simeq
\overline{I}_{\delta}^{(+)}({\frac n 2}, \frac n2)
\simeq 
\overline{I}_{\delta}^{(-)}({\frac n 2}, \frac n2)
\quad
\text{for $\delta = \pm$}.  
\end{equation}
Similarly to Theorem \ref{thm:LNM20} 
 about the $O(n+1,1)$-modules $\Pi_{\ell,\delta}$, 
 we summarize the properties of the restriction $\overline{\Pi}_{\ell,\delta}=\Pi_{\ell,\delta}|_{\overline G}$ as follows.  
\begin{proposition}
\label{prop:161648}
Let ${\overline G}=SO(n+1,1)$
 with $n \ge 1$.  
\begin{enumerate}
\item[{\rm{(1)}}]
$\overline{\Pi}_{\ell,\delta}$ is irreducible as a $\overline G$-module
 for all $0 \le \ell \le n+1$ and $\delta = \pm $.  
\item[{\rm{(2)}}]
$\overline{\Pi}_{\ell,\delta} \simeq \overline{\Pi}_{n+1-\ell,-\delta}$
 as $\overline G$-modules
 for all $0 \le \ell \le n+1$ and $\delta = \pm $.  
\item[{\rm{(3)}}]
Irreducible representations of $\overline G$
 with ${\mathfrak {Z}}({\mathfrak {g}})$-infinitesimal character $\rho_{\overline G}$
 can be classified as 
\begin{alignat*}{2}
&
 \{ 
 \overline{\Pi}_{\ell,\delta}: 0 \le \ell \le \frac{n-1}{2}, \delta = \pm 
 \}
 \cup
 \{ 
 \overline{\Pi}_{\frac{n+1}{2},+}
 \}
 \qquad
&&
 \text{if $n$ is odd}, 
\\
&
 \{ 
 \overline{\Pi}_{\ell,\delta}: 0 \le \ell \le \frac{n}{2}, \delta = \pm 
 \}
 \qquad
&&
 \text{if $n$ is even}.  
\end{alignat*}
\item[{\rm{(4)}}]
Every $\overline{\Pi}_{\ell,\delta}$
 $(0 \le \ell \le n+1, \delta=\pm)$
 is unitarizable.  
\end{enumerate}
In the next statement,
 we use the same symbol $\overline{\Pi}_{\ell,\delta}$
 to denote the irreducible unitary representation 
obtained by the Hilbert completion 
of $\overline \Pi_{\ell,\delta}$ with respect to a $\overline G$-invariant inner product.  
\begin{enumerate}
\item[{\rm{(5)}}]
For $n$ odd, 
 $\overline{\Pi}_{\frac{n+1}{2},+}$ is a discrete series representation
 of $\overline G=SO(n+1,1)$.  
For $n$ even, 
 $\overline{\Pi}_{\frac{n}{2},\delta}$ $(\delta=\pm)$
 are 
\index{B}{temperedrep@tempered representation}
tempered representations.  
All the other representations in the list (2)
 are nontempered representations of $\overline G$.  
\item[{\rm{(6)}}]
For $n$ even, 
 the center of $\overline G=SO(n+1,1)$ acts nontrivially
 on $\overline{\Pi}_{\ell,\delta}$
 if and only if $\delta=(-1)^{\ell+1}$.  
For $n$ odd, 
 the center of $\overline G$ is trivial, 
 and thus acts trivially
 on $\overline{\Pi}_{\ell,\delta}$
 for any $\ell$ and $\delta$.  
\end{enumerate}
\end{proposition}

In Proposition \ref{prop:161655}, 
 we gave a description of the underlying $({\mathfrak{g}},{K})$-module
 of the irreducible $G$-module $\Pi_{\ell,\delta}$ 
 in terms of cohomological parabolic induction.  
We end this section 
 with analogous results
for the irreducible $\overline G$-module 
 $\overline \Pi_{\ell,\delta}= \Pi_{\ell,\delta}|_{\overline G}$
 (see Proposition \ref{prop:161648} (1)).  
\begin{proposition}
\label{prop:SO161655}
For $0 \le i \le [\frac {n+1}2]$, 
 let 
\index{A}{qi@${\mathfrak{q}}_i$, $\theta$-stable parabolic subalgebra}
${\mathfrak{q}}_i$ be the $\theta$-stable parabolic subalgebras
 with the Levi subgroup $\overline{L_i} \simeq SO(2)^i \times SO(n-2i+1,1)$
 as in Definition \ref{def:qi}
 and write $S_i=i(n-i)$, 
 see \eqref{eqn:cohSi}.  
Then the underlying $({\mathfrak {g}},\overline{K})$-modules
 of the irreducible $\overline G$-modules $\overline{\Pi}_{\ell, \delta}$
 $(0 \le \ell \le n+1$, $\delta \in \{\pm\})$
 are given by the cohomological parabolic induction
 as follows:
\index{A}{Aqlmd@$A_{\mathfrak{q}}(\lambda)$}
\begin{alignat*}{4}
  (\overline{\Pi}_{i,+})_{\overline K} 
&\simeq
 (\overline{\Pi}_{n+1-i,-})_{\overline K}
 \simeq
&& {\mathcal{R}}_{\mathfrak {q}_i}^{S_i}({\mathbb{C}}_{\rho({\mathfrak{u}})})
&& \simeq (A_{\mathfrak{q}_i})_{+,+}|_{({\mathfrak {g}},\overline{K})} 
&&  \simeq (A_{\mathfrak{q}_i})_{-,-}|_{({\mathfrak {g}},\overline{K})}, 
\\
   (\overline{\Pi}_{i,-})_{\overline K} 
&\simeq 
 (\overline {\Pi}_{n+1-i,+})_{\overline K}
   \simeq 
   && {\mathcal{R}}_{\mathfrak {q}_i}^{S_i}({\mathbb{C}}_{\rho({\mathfrak{u}})} \otimes \chi_{+-})
&& \simeq (A_{\mathfrak{q}_i})_{+,-}|_{({\mathfrak {g}},\overline{K})} 
&&\simeq (A_{\mathfrak{q}_i})_{-,+}|_{({\mathfrak {g}},\overline{K})}.  
\end{alignat*}
\end{proposition}
We notice that the four characters $\chi_{\pm \pm}$
 of $O(n+1,1)$ induce the following isomorphisms $\chi_{--} \simeq {\bf{1}}$
 and $\chi_{+-} \simeq \chi_{-+}$
 when restricted to the last factor
 $SO(n-2i+1,1)$ of the Levi subgroup $\overline{L}_i$, 
 whence Proposition \ref{prop:SO161655} gives an alternative proof
 for the isomorphism
\[
  \overline{\Pi}_{i,\delta} \simeq \overline{\Pi}_{n+1-i,-\delta}
\]
 as $\overline G$-modules
 for $0 \le i \le n+1$ and $\delta = \pm$.

\subsection{Symmetry breaking for tempered principal series representations}
\label{subsec:SOtemp}

In this section, 
 we deduce a multiplicity-one theorem 
 for tempered principal series representations 
$
   \overline I_{\delta}(\overline V, \lambda)
$
 and 
$
   \overline J_{\varepsilon}(\overline W, \nu)
$
 of $\overline G =SO(n+1,1)$ and $\overline {G'} =SO(n,1)$, 
 respectively, from the corresponding result
 (see Theorem \ref{thm:tempVW})
 for the pair $(G, G') = (O(n+1,1), O(n,1))$.

In \cite[Chap.~2, Sect.~5]{KKP}, 
 a trick analogous to Lemma \ref{lem:LNM26} was used 
 to deduce symmetry breaking 
 for the pair $(G_0, G_0')=(SO_0(n+1,1), SO_0(n,1))$ from that 
 for the pair $(G, G')$
 by using an observation
 that $G_0$ and $G_0'$ are normal subgroups
 of $G$ and $G'$, 
 respectively
 (cf.~\cite[page 26]{KKP}).  
This is formulated in our setting as follows:
\begin{proposition}
\label{prop:LNM26}
Let $\Pi$ and $\pi$ be continuous representations
 of $G=O(n+1,1)$ and $G'=O(n,1)$, 
respectively.  
Let $(\overline G, \overline {G'})=(SO(n+1,1), SO(n,1))$.  
Then we have natural isomorphisms:
\begin{align*}
  \operatorname{Hom}_{\overline{G'}}(\Pi|_{\overline{G'}}, \pi|_{\overline{G'}})&\simeq
  \operatorname{Hom}_{G'}(\Pi|_{G'}, \pi)
  \oplus
  \operatorname{Hom}_{G'}(\Pi|_{G'}, \chi_{--}\otimes \pi)
\\
&\simeq
  \operatorname{Hom}_{G'}(\Pi|_{G'}, \pi)
  \oplus
  \operatorname{Hom}_{G'}((\Pi \otimes \chi_{--})|_{G'}, \pi).  
\end{align*}
\end{proposition}

For $\overline V \in \widehat{SO(n)}$
 and $W \in \widehat{SO(n-1)}$, 
 we set
\[
[\overline V|_{SO(n-1)}:\overline W]
:=\dim_{\mathbb{C}} {\operatorname{Hom}}_{SO(n-1)}(\overline V|_{SO(n-1)},W).  
\]
The main result of this section
 is the following.  
\begin{theorem}
[tempered principal series representation]
\label{thm:tempSO}
Let $\overline V \in \widehat {SO(n)}$, 
 $\overline W \in \widehat {SO(n-1)}$, 
 $\delta, \varepsilon \in \{\pm\}$, 
 and $(\lambda,\nu) \in (\sqrt{-1}{\mathbb{R}}+ \frac n 2, 
 \sqrt{-1}{\mathbb{R}}+\frac 1 2(n-1))$
 so that $\overline {I}_{\delta}(\overline V, \lambda)$
 and $\overline {J}_{\varepsilon}(\overline W, \nu)$ are 
 irreducible tempered principal series representations
 of $\overline G=SO(n+1,1)$ and $\overline {G'}=SO(n,1)$, 
 respectively.  
Then the following conditions are equivalent:
\begin{enumerate}
\item[{\rm{(i)}}]
$[\overline V|_{SO(n-1)}: \overline W] 
\ne 0$.  

\item[{\rm{(ii)}}]
$\operatorname{Hom}_{SO(n,1)}
 (\overline {I}_{\delta}(\overline V, \lambda)|_{SO(n,1)}, 
  \overline {J}_{\varepsilon}(\overline W, \nu)) 
\ne \{0\}.$

\item[{\rm{(iii)}}]
$\dim_{\mathbb{C}}\operatorname{Hom}_{SO(n,1)}
 (\overline {I}_{\delta}(\overline V, \lambda)|_{SO(n,1)}, 
  \overline {J}_{\varepsilon}(\overline W, \nu)) =1.$  
\end{enumerate}
\end{theorem}
For the proof, 
we use the following elementary lemma
 on branching rules of finite-dimensional representations of $O(n)$.  
\begin{lemma}
\label{lem:Xbranch}
Suppose $\sigma \in \widehat{O(n)}$ and $\tau \in \widehat{O(n-1)}$
 are of both type $X$
 (Definition \ref{def:OSO}).  
If $[\sigma|_{O(n-1)}:\tau]\ne0$, 
 then $[\sigma|_{O(n-1)}:\tau \otimes \det]=0$.  
\end{lemma}

\begin{proof}
[Proof of Lemma \ref{lem:Xbranch}]
Easy from Fact \ref{fact:ONbranch} and from the characterization in Lemma \ref{lem:OSO}
 of representations of type X
 by means of the Cartan--Weyl bijection \eqref{eqn:CWOn}.  
\end{proof}
\begin{proof}
[Proof of Theorem \ref{thm:tempSO}]
There exist unique $V \in \widehat{O(n)}$ and $W \in \widehat{O(n-1)}$
such that $[V|_{SO(n)}:\overline V]\ne0$
 and $[W|_{SO(n-1)}:\overline W]\ne0$.  
We divide the argument into the following three cases:

Case XX: Both $V$ and $W$ are of type X. 

Case XY: $V$ is of type X and $W$ is of type Y.  

Case YX: $V$ is of type Y and $W$ is of type X.

Then we have from \eqref{eqn:IVSOpm}
\[
{I}_{\delta}(V, \lambda)|_{\overline G}
\simeq
\begin{cases} 
  \overline {I}_{\delta}(\overline V, \lambda)
 &\text{if $V$ is of type X,}
\\
\overline{I}_{\delta}(\overline V, \lambda)
\oplus
\overline {I}_{\delta}(\overline V^{\vee}, \lambda)
\qquad
 &\text{if $V$ is of type Y,}
\end{cases}
\]
and similarly for the restriction $J_{\varepsilon}(W, \nu)|_{\overline {G'}}$.

By Proposition \ref{prop:LNM26}, 
 we have
\begin{equation*}
  {\operatorname{Hom}}_{\overline{G'}} 
  ({I}_{\delta}(V, \lambda)|_{\overline {G'}}, 
   {J}_{\varepsilon}(W, \nu)|_{\overline {G'}})
   \simeq
   \bigoplus_{\chi \in \{{\bf{1}},\det\}}
   {\operatorname{Hom}}_{G'}({{I}_{\delta}(V, \lambda)|_{G'}},
   J_{\varepsilon}(W, \nu) \otimes \chi).  
\end{equation*}
Applying the 
\index{B}{multiplicityonetheorem@multiplicity-one theorem}
multiplicity-one theorem (Theorem \ref{thm:tempVW})
 for tempered representations 
 of the pair $(G,G')=(O(n+1,1),O(n,1))$
 to the right-hand side, 
 we get the following multiplicity formula:
\begin{multline}
\label{eqn:VWdetdim}
     \dim_{\mathbb{C}} {\operatorname{Hom}}_{\overline{G'}} 
  ({I}_{\delta}(V, \lambda)|_{\overline {G'}}, 
   {J}_{\varepsilon}(W, \nu)|_{\overline {G'}})
\\
=[V|_{O(n-1)}:W]+[V|_{O(n-1)}:W \otimes \det].  
\end{multline}
The right-hand side of \eqref{eqn:VWdetdim} does not vanish
 if and only if 
 $[\overline V|_{SO(n-1)}:\overline W] \ne 0$.  
In this case,
 we have 
\begin{equation*}
\text{\eqref{eqn:VWdetdim}}
=\begin{cases}
  1 \qquad &\text{Case XX, }
\\
  2 \qquad &\text{Case XY or Case YX, }
  \end{cases}
\end{equation*}
by Lemmas \ref{lem:branchII} and \ref{lem:Xbranch}.  
Thus the conclusion holds in Case XX.

If $V$ is of type Y, 
 then the two $\overline G$-irreducible summands
 $\overline{I}_{\delta}(\overline V, \lambda)$
 and $\overline{I}_{\delta}(\overline V^{\vee}, \lambda)$
 in the restriction $I_{\delta}(V, \lambda)|_{G'}$ are switched 
 if we apply the outer automorphism of $\overline G$
 by an element $g_0:={\operatorname{diag}}(1,\cdots,1,-1,1) \in G=O(n+1,1)$.  
Since $g_0$ commutes with $\overline {G'}$, 
 we obtain an isomorphism
\[
  {\operatorname{Hom}}_{\overline{G'}} 
  (\overline{I}_{\delta}(\overline V, \lambda)|_{\overline {G'}}, 
   \overline{J}_{\varepsilon}(\overline W, \nu))
  \simeq
  {\operatorname{Hom}}_{\overline{G'}} 
  (\overline{I}_{\delta}(\overline V^{\vee}, \lambda)|_{\overline {G'}}, 
   \overline{J}_{\varepsilon}(\overline W, \nu)).
\]
Hence the conclusion holds
 for Case YX.

Similar argument holds for Case XY
 where $W$ is of type Y.  
Therefore Theorem \ref{thm:tempSO} is proved.  
\end{proof}

\subsection{Symmetry breaking from
 $\overline I_{\delta}(i,\lambda)$ to $\overline J_{\varepsilon}(j,\nu)$}
\label{subsec:SBOpsSO}

In this section, 
 we give a closed formula
 of the multiplicity
 for the restriction $\overline G \downarrow \overline G'$
 when $(\sigma,V)$ is the exterior tensor $\Exterior^i({\mathbb{C}}^n)$.  
For the admissible smooth representations
 ${\overline{I}}_{\delta}(i,\lambda)$ of ${\overline{G}}=SO(n+1,1)$
 and ${\overline{J}}_{\varepsilon}(j,\nu)$ of ${\overline{G'}}=SO(n,1)$, 
 we set
\[
   m(i,j)
   \equiv m({\overline{I}}_{\delta}(i,\lambda), {\overline{J}}_{\varepsilon}(j,\nu))
   :=
   \dim_{\mathbb{C}} \operatorname{Hom}_{{\overline{G'}}}
    ({\overline{I}}_{\delta}(i,\lambda)|_{{\overline{G'}}}, {\overline{J}}_{\varepsilon}(j,\nu)).  
\]
%

In order to state a closed formula
 for the multiplicity $m(i,j)$
 as a function of $(\lambda, \nu,\delta, \varepsilon)$, 
 we introduce the following subsets of ${\mathbb{Z}}^2 \times \{\pm 1\}$:
\index{A}{Leven@$L_{\operatorname{even}}$}
\index{A}{Lodd@$L_{\operatorname{odd}}$}
\begin{align*}
  L:=&\{(-i,-j, (-1)^{i+j}): (i,j) \in {\mathbb{Z}}^2, 0 \le j \le i \}
    = L_{\operatorname{even}} \cup L_{\operatorname{odd}}, 
\\
  L':=&\{(\lambda,\nu,\gamma) \in L: \nu \ne 0 \}.  
\end{align*}

In the theorem below, we shall see
\begin{alignat*}{2}
m(i,j) \in &\{ 1,2,4 \} \qquad
&&\text{if $j=i-1$ or $i$}, 
\\
m(i,j) \in &\{ 0,1,2 \} \qquad
&&\text{if $j=i-2$ or $i+1$}, 
\\
m(i,j) =& 0 \qquad
&&\text{otherwise}.  
\end{alignat*}
By Proposition \ref{prop:LNM26} and Lemma \ref{lem:LNM27}, 
 the multiplicity formula
 for $(\overline G, \overline {G'})$ is derived from 
 the one for $(G,G')$ 
 by using Proposition \ref{prop:LNM26}, 
 which amounts to  
\begin{multline*}
\operatorname{Hom}_{\overline{G'}}(\overline{I}_{\delta}(i,\lambda)|_{\overline{G'}}, \overline{J}_{\varepsilon}(j,\nu))
\\
\simeq
 \operatorname{Hom}_{G'}({I}_{\delta}(i,\lambda)|_{G'}, {J}_{\varepsilon}(j,\nu))
  \oplus
  \operatorname{Hom}_{G'}({I}_{\delta}(n-i,\lambda)|_{G'}, {J}_{\varepsilon}(j,\nu)).  
\end{multline*}
The right-hand side was computed 
 in Theorem \ref{thm:1.1}.  
Hence we get an explicit formula of the multiplicity
 for the restriction of nonunitary principal series representations
in this setting:

\begin{theorem}
[multiplicity formula]
\label{thm:SOmult}
Suppose $n \ge 3$, $0 \le i\le [\frac n2]$, $0 \le j \le [\frac{n-1}2]$,
 $\delta$, $\varepsilon \in \{\pm \}\equiv \{ \pm 1\}$, 
 and $\lambda,\nu \in {\mathbb{C}}$.

Then the multiplicity 
$
   m(i,j)=
   \dim_{\mathbb{C}} \operatorname{Hom}_{{\overline{G'}}}
    ({\overline{I}}_{\delta}(i,\lambda)|_{{\overline{G'}}}, 
     {\overline{J}}_{\varepsilon}(j,\nu))
$
 is given as follows.  
\begin{enumerate}
\item[{\rm{(1)}}]
Suppose $j=i$.  
\begin{enumerate}
\item[{\rm{(a)}}]
Case $i=0$.
\begin{equation*}
m(0,0)
=
\begin{cases}
2
\qquad
&\text{if }
(\lambda, \nu, \delta\varepsilon) \in L, 
\\
1
&
\text{otherwise.}
\end{cases}
\end{equation*}
\item[{\rm{(b)}}]
Case $1 \le i < \frac  n2-1$.  
\begin{equation*}
m(i, i)
=
\begin{cases}
2
\qquad
&\text{if }
(\lambda, \nu, \delta \varepsilon) \in L' \cup \{(i,i,+)\},  
\\
1 
&
\text{otherwise}.  
\end{cases}
\end{equation*}
\item[{\rm{(c)}}]
Case $i= \frac n 2 -1$ $($$n$: even$)$.   
\begin{equation*}
m(\frac n 2 -1, \frac n 2 -1)
=
\begin{cases}
2
\qquad
&\text{if }
(\lambda, \nu, \delta\varepsilon) \in L' \cup \{(i,i,+)\}\cup \{(i,i+1,-)\},   
\\
1
&
\text{otherwise}.  
\end{cases}
\end{equation*}
\item[{\rm{(d)}}]
Case $i= \frac {n-1} 2$ $($$n$: odd$)$.   
\begin{equation*}
m(\frac {n-1} 2,\frac {n-1} 2)
=
\begin{cases}
4
\qquad
&\text{if }
(\lambda, \nu, \delta\varepsilon) \in L' \cup \{(i,i,+)\}, 
\\
2
&
\text{otherwise}.  
\end{cases}
\end{equation*}
\end{enumerate}

\item[{\rm{(2)}}]
Suppose $j=i-1$.  
\begin{enumerate}
\item[{\rm{(a)}}]
Case $1 \le i < \frac{n-1}{2}$.
\begin{equation*}
m(i,i-1)
=
\begin{cases}
2
\qquad
&\text{if }
(\lambda, \nu, \delta\varepsilon) \in L' \cup \{(n-i,n-i,+)\}, 
\\
1
&
\text{otherwise.}
\end{cases}
\end{equation*}
\item[{\rm{(b)}}]
Case $i = \frac{n-1}{2}$ $($$n$: odd$)$.  
\begin{equation*}
m(\frac{n-1}{2},\frac{n-3}{2})
=
\begin{cases}
2
\qquad
&\text{if }
(\lambda, \nu, \delta \varepsilon) \in L',  
\\
2
\qquad
&\text{if }
(\lambda, \nu, \delta \varepsilon) \in \{(n-i,n-i,+)\}
                                       \cup \{(i,i+1,-)\} ,  
\\
1 
&
\text{otherwise}.  
\end{cases}
\end{equation*}
\item[{\rm{(c)}}]
Case $i= \frac n 2$ $($$n$: even$)$.   
\begin{equation*}
m(\frac n 2,\frac n 2-1)
=
\begin{cases}
4
\qquad
&\text{if }
(\lambda, \nu, \delta\varepsilon) \in L' \cup \{(n-i,n-i,+)\},   
\\
2
&
\text{otherwise}.  
\end{cases}
\end{equation*}
\end{enumerate}

\item[{\rm{(3)}}]
Suppose $j=i-2$.  
\begin{enumerate}
\item[{\rm{(a)}}]
Case $2 \le i < \frac{n}{2}$.
\begin{equation*}
m(i,i-2)
=
\begin{cases}
1
\qquad
&\text{if }
(\lambda, \nu, \delta\varepsilon) = (n-i,n-i+1,-), 
\\
0
&
\text{otherwise.}
\end{cases}
\end{equation*}
\item[{\rm{(b)}}]
Case $i = \frac{n}{2}$ $($$n$: even$)$.  
\begin{equation*}
m(\frac{n}{2},\frac{n}{2}-2)
=
\begin{cases}
2
\qquad
&\text{if }
(\lambda, \nu, \delta \varepsilon) =(\frac{n}{2},\frac{n}{2}+1,-),  
\\
0 
&
\text{otherwise}.  
\end{cases}
\end{equation*}
\end{enumerate}

\item[{\rm{(4)}}]
Suppose $j=i+1$.  
\begin{enumerate}
\item[{\rm{(a)}}]
Case $i=0$ and $n>3$.
\begin{equation*}
m(0,1)
=
\begin{cases}
1
\qquad
&\text{if }
\lambda \in - {\mathbb{N}}, \nu=1, \text{ and }\delta\varepsilon = (-1)^{\lambda+1}, 
\\
0
&
\text{otherwise.}
\end{cases}
\end{equation*}
\item[{\rm{(b)}}]
Case $1 \le i < \frac{n-3}{2}$.  
\begin{equation*}
m(i,i+1)
=
\begin{cases}
1
\qquad
&\text{if }
(\lambda, \nu, \delta \varepsilon) =(i,i+1,-),  
\\
0 
&
\text{otherwise}.  
\end{cases}
\end{equation*}
\item[{\rm{(c)}}]
Case $i=\frac{n-3}{2}$ and $n>3$, odd.  
\begin{equation*}
m(\frac{n-3}{2}, \frac{n-1}{2})
=
\begin{cases}
2
\qquad
&\text{if }
(\lambda, \nu, \delta \varepsilon) =(\frac{n-3}{2},\frac{n-1}{2},-),  
\\
0 
&
\text{otherwise}.  
\end{cases}
\end{equation*}
\item[{\rm{(d)}}]
Case $i=0$ and $n=3$.  
\begin{equation*}
m(0, 1)
=
\begin{cases}
2
\qquad
&\text{if 
$\lambda \in -{\mathbb{N}}$, $\nu=1$, 
 and $\delta \varepsilon = (-1)^{\lambda+1}$},  
\\
0 
&
\text{otherwise}.  
\end{cases}
\end{equation*}
\end{enumerate}

\item[{\rm{(5)}}]
Suppose $j \not\in \{i-2, i-1, i, i+1\}$.
Then $m(i,j)=0$ for all $\lambda, \nu, \delta, \varepsilon$.
\end{enumerate}

\end{theorem}

\begin{remark}
[multiplicity-one property]
\label{rem:mult1}
In \cite{SunZhu}
 it is proved that 
$$
   \dim_{\mathbb{C}}\operatorname{Hom}_{{\overline{G'}}}(\Pi|_{{\overline{G'}}},\pi) \le 1
$$
 for any irreducible admissible smooth representations
 $\Pi$ and $\pi$ of ${\overline{G}}=SO(n+1,1)$ and ${\overline{G'}}=SO(n,1)$, 
respectively.  
Thus Theorem \ref{thm:1.1} fits well 
 with their multiplicity-free results 
 for 
  $\lambda, \nu \in {\mathbb{C}}\setminus {\mathbb{Z}}$, 
 where ${\overline{I}}_{\delta}(i,\lambda)$ and ${\overline{J}}_{\varepsilon}(j,\nu)$
 are irreducible admissible representations
 of ${\overline{G}}$ and ${\overline{G'}}$, 
 respectively,
 except for the cases $n=2 i$ or $n= 2j+1$.  
In the case $n=2i$ or $n=2j+1$, 
 the multiplicity is counted twice 
 as we saw in \eqref{eqn:ISOneven} and \eqref{eqn:161733}, 
 and thus the statements (1-d), (2-c), (3-b), 
 and (4-c)
 in Theorem \ref{thm:1.1} fit again with \cite{SunZhu}.  
\end{remark}
\vskip 1pc
\begin{remark}
[generic multiplicity-two phenomenon]
\label{rem:mult2}
In addition to the subgroup ${\overline{G'}}=SO(n,1)$, 
 the Lorentz group $O(n,1)$ contains two other subgroups
 of index two, 
 that is, 
 $O^+(n,1)$ (containing orthochronous reflections)
 and  $O^-(n,1)$ (containing anti-orthochronous reflections)
 with terminology
 in relativistic space-time for $n=3$.  
Our results yield also the multiplicity formula for such pairs
 by using an analogous result to Proposition \ref{prop:LNM26},
and it turns out that a 
\index{B}{genericmultiplicityonetheorem@generic multiplicity-one theorem}
generic multiplicity-one statement
 fails 
 if we replace $({\overline{G}},{\overline{G'}})=(SO(n+1,1),SO(n,1))$
 by $(O^-(n+1,1),O^-(n,1))$.  
In fact, 
 the multiplicity $m(\Pi, \pi)$ is generically equal to 2
 for irreducible representations $\Pi$ and $\pi$
 of $O^-(n+1,1)$ and $O^-(n,1)$, 
 respectively, 
 as is expected by the general theory \cite{xKOfm, sbon}
 because there are two open orbits
 in $P'^{-} \backslash G^-/P^-$ in this case.  
\end{remark}

\subsection{Symmetry breaking
 between irreducible representations
 of $\overline G$ and $\overline{G'}$
 with trivial infinitesimal character $\rho$}
Similar to the notation $\overline{\Pi}_{i,\delta}$
 for the restriction of the irreducible representation
 ${\Pi}_{i,\delta}$ of $G=O(n+1,1)$ 
 to the special orthogonal group $\overline G=SO(n+1,1)$, 
 we denote by $\overline \pi_{j,\varepsilon}$
 the restriction of the irreducible representation $\pi_{j,\varepsilon}$
 ($0 \le j \le n$, $\varepsilon=\pm$)
 of $G'=O(n,1)$ to the special orthogonal group $\overline {G'}=SO(n,1)$.  
Then $\overline{\Pi}_{i,\delta}$ ($0 \le i \le n+1$, $\delta=\pm$)
 and $\overline{\pi}_{j,\varepsilon}$ ($0 \le j \le n$, $\varepsilon=\pm$)
 exhaust irreducible admissible smooth representations of $\overline G$
 and $\overline {G'}$ 
 having ${\mathfrak{Z}}({\mathfrak{g}})$-infinitesimal character 
 $\rho_{\overline G}$ 
 and ${\mathfrak{Z}}({\mathfrak{g}}')$-infinitesimal character 
 $\rho_{\overline {G'}}$
 respectively,
 by Lemma \ref{lem:Pilso}.

In this section,
 we deduce the formula
 of the multiplicity
\[ 
   \dim_{\mathbb{C}} {\operatorname{Hom}}_{\overline {G'}}
   (\overline \Pi_{i,\delta}|_{\overline {G'}}, \overline \pi_{j,\varepsilon})
\]
for the symmetry breaking for 
 $(\overline G, \overline {G'})=(SO(n+1,1), SO(n,1))$ from the one
 for $(G, G')=(O(n+1,1), O(n,1))$.  
 
In view of the $\overline G$-isomorphism
 $\overline{\Pi}_{\frac {n+1}2,+}\simeq \overline{\Pi}_{\frac {n+1}2,-}$
 for $n$ even
 and the $\overline {G'}$-isomorphism
 $\overline{\pi}_{\frac {n}2,+}\simeq \overline{\pi}_{\frac {n}2,-}$
 for $n$ odd, 
 we shall use the following convention
\begin{equation}
\label{eqn:halfsgn}
  \text{$+ \equiv -$ for $\delta$ if $n+1=2i$;
\quad
       $+ \equiv -$ for $\varepsilon$ if $n=2j$
}
\end{equation}
when we deal with the representations
 $\overline{\Pi}_{i,\delta}$
 ($0 \le i \le [\frac{n+1}{2}]$)
 and $\overline{\pi}_{j,\varepsilon}$
 ($0 \le j \le [\frac{n}{2}]$).

Owing to Proposition \ref{prop:LNM26}, 
 Theorem \ref{thm:LNM20} tells that 
\begin{equation*}
 \operatorname{Hom}_{\overline{G'}}(\overline{\Pi}_{i,\delta}|_{\overline{G'}}, \overline{\pi}_{j, \varepsilon})
\simeq
  \operatorname{Hom}_{G'}(\Pi_{i,\delta}|_{G'}, \pi_{j,\varepsilon})
  \oplus
  \operatorname{Hom}_{G'}(\Pi_{n+1-i,-\delta}|_{G'}, \pi_{j, \varepsilon}).
\end{equation*}
Applying Theorems \ref{thm:SBOvanish} and \ref{thm:SBOone}
 about symmetry breaking for the pair
 $(G,G')=(O(n+1,1), O(n,1))$
 to the right-hand side, 
 we determine the multiplicity
\[
   m(\overline \Pi, \overline \pi)
\quad
\text{for all $\overline \Pi \in {\operatorname{Irr}}(\overline G)_{\rho}$
      and $\overline \pi \in {\operatorname{Irr}}(\overline {G'})_{\rho}$}
\]
 for the pair $(\overline G, \overline{G'})=(SO(n+1,1), SO(n,1))$
 of special orthogonal groups as follows.  
\begin{theorem}
\label{thm:170336}
Suppose $0 \le i \le [\frac{n+1}{2}]$, 
 $0 \le j \le [\frac{n}{2}]$, 
 and $\delta, \varepsilon = \pm$
 with the convention \eqref{eqn:halfsgn}.  
then 
\[
\dim_{\mathbb{C}} {\operatorname{Hom}}_{\overline{G'}}
(\overline{\Pi}_{i,\delta}|_{\overline{G'}}, 
 \overline{\pi}_{j,\varepsilon})
=
\begin{cases}
1
\quad
&\text{if $\delta \equiv \varepsilon$ and $j \in \{i-1,i\}$},
\\
0
&\text{otherwise.}
\end{cases}
\]
\end{theorem}

\newpage
\section{Appendix III: A translation functor for $G=O(n+1,1)$}
\label{sec:Translation}

In this chapter,
 we discuss a translation functor 
 for the group $G=O(n+1,1)$, 
 which is not in the Harish-Chandra class
 if $n$ is even, 
 in the sense that $\operatorname{Ad}(G)$ is not contained 
 in the group $\operatorname{Int}({\mathfrak{g}}_{\mathbb{C}})$
 of inner automorphisms.  
Then the \lq\lq{Weyl group}\rq\rq\
 $W_G$ is larger
 than the group generated 
 by the reflections of simple roots.  
This causes some technical difficulties
 when we extend the idea of translation functor
 which is usually formulated 
 for reductive groups
 in the Harish-Chandra class or reductive Lie algebras, 
 see \cite{Jantzen, SpehVogan, Vogan81, Zuckerman} for instance.  

\subsection{Some features of translation functors
 for reductive groups
 that are not of Harish-Chandra class}
For $n$ even, 
say $n=2m$, 
 we write ${\mathfrak{h}}_{\mathbb{C}}$
 $(\simeq {\mathbb{C}}^{m+1})$
 for a Cartan subalgebra of ${\mathfrak{g}}_{\mathbb{C}}$.  
Then we recall from Section \ref{subsec:2.1.4}:
\begin{enumerate}
\item[$\bullet$]
\index{A}{Weylgroupg@$W_{\mathfrak{g}}$, Weyl group for ${\mathfrak{g}}_{\mathbb{C}}={\mathfrak{o}}(n+2,{\mathbb{C}})$}
the Weyl group $W_{\mathfrak{g}} \simeq {\mathfrak{S}}_{m+1} \ltimes ({\mathbb{Z}}/2{\mathbb{Z}})^{m}$ for the root system
 $\Delta({\mathfrak{g}}_{\mathbb{C}}, {\mathfrak{h}}_{\mathbb{C}})$
 is of index two in the Weyl group 
\index{A}{WeylgroupG@$W_G$, Weyl group for $G=O(n+1,1)$}
$W_G\simeq {\mathfrak{S}}_{m+1} \ltimes ({\mathbb{Z}}/2{\mathbb{Z}})^{m+1}$
 for the disconnected group $G$;
\item[$\bullet$]
the ${\mathfrak{Z}}_G({\mathfrak{g}})$-infinitesimal character
 for the irreducible admissible representation of $G$
 is parametrized by ${\mathfrak{h}}_{\mathbb{C}}^{\ast}/W_G$,
 but not by ${\mathfrak{h}}_{\mathbb{C}}^{\ast}/W_{\mathfrak{g}}$;
\item[$\bullet$]
$\rho_G=(m,\cdots,1,0)$ is not 
 \lq\lq{$W_G$-regular}\rq\rq, 
 although it is \lq\lq{$W_{\mathfrak{g}}$-regular}\rq\rq\
 (Definition \ref{def:intreg}).  
\end{enumerate}

We can still use the idea of a translation functor,
 but we need a careful treatment
 for disconnected groups $G$
 which are not in the Harish-Chandra class.  
In fact,
 differently from the usual setting for reductive Lie groups
 in the Harish-Chandra class, 
 we are faced with the following feature:
\begin{enumerate}
\item[$\bullet$]
translation from a
\index{B}{WGregular@$W_G$-regular}
 $W_G$-regular
 (in particular, 
\index{B}{Wgregular@$W_{\mathfrak{g}}$-regular}
$W_{\mathfrak{g}}$-regular) dominant parameter
 to the trivial infinitesimal character $\rho_G$ 
 does not necessarily preserve irreducibility, 
 see Theorem \ref{thm:1808101}. 
\end{enumerate}
This means that translation inside the same \lq\lq{$W_{\mathfrak{g}}$-regular Weyl chamber}\rq\rq\
 may involve a phenomenon 
 as if it were \lq\lq{translation from the wall 
 to regular parameter}\rq\rq, 
 {\it{cf.}},  \cite{SpehVogan}.

In what follows,
 we retain the terminology
 \lq\lq{regular}\rq\rq\ for $W_{\mathfrak{g}}$
 but not for $W_G$
 as in Definition \ref{def:intreg}
 (in particular, $\rho_G$ is regular in our sense), 
 whereas we need to use $W_G$
 (not $W_{\mathfrak{g}}$)
 in describing ${\mathfrak{Z}}_G({\mathfrak{g}})$-infinitesimal characters
 of $G$-modules.

\subsection{Translation functor for $G=O(n+1,1)$}

In this section 
we fix some notation
 for a 
\index{B}{translationfunctor@translation functor|textbf}
{\it{translation functor}}
 for the group $G=O(n+1,1)$.  
Usually, 
 a translation functor is defined in the category 
 of $({\mathfrak{g}},K)$-modules
 of finite length.  
However, 
 we also consider a translation functor
  in the category of {\bf{admissible representations
 of finite length of moderate growth}}.

\subsubsection{Primary decomposition of admissible smooth representations}
\label{subsec:primary}

Let $\Pi$ be an admissible smooth representation of $G$ of finite length.  
For $\mu \in {\mathfrak{h}}_{\mathbb{C}}^{\ast}/W_G$, 
 we define the $\mu$-primary component $P_{\mu}(\Pi)$ of $\Pi$ by 
\[
    P_{\mu}(\Pi) := \bigcup_{N >0} \bigcap_{z \in {\mathfrak{Z}}_G({\mathfrak{g}})}
   \operatorname{Ker}(z-\chi_{\mu}(z))^N, 
\]
where we recall the Harish-Chandra isomorphism \eqref{eqn:HCpara}
\index{A}{ZGg@${\mathfrak{Z}}_G({\mathfrak{g}})$}
\[
   \operatorname{Hom}_{\mathbb{C}\text{-alg}}({\mathfrak{Z}}_G({\mathfrak{g}}), {\mathbb{C}})
   \simeq {\mathfrak{h}}_{\mathbb{C}}^{\ast}/W_G, 
\quad
   \chi_{\mu} \leftrightarrow \mu.  
\]
Then $P_{\mu}(\Pi)$ is a $G$-module with generalized ${\mathfrak{Z}}_G({\mathfrak{g}})$-infinitesimal character $\mu$, 
 and $\Pi$ is decomposed into a direct sum 
 of finitely many primary components:
\[
  \Pi = \bigoplus_{\mu}P_{\mu}(\Pi)
\qquad
\text{(finite direct sum)}.  
\]
By abuse of notation,
 we use the letter 
\index{A}{Pmu@$P_{\mu}$|textbf}
$P_{\mu}$ 
 to denote the $G$-equivariant projection $\Pi \to P_{\mu}(\Pi)$
 with respect to the direct sum decomposition.

\subsubsection{Translation functor $\psi_{\mu}^{\mu+\tau}$
 for $G=O(n+1,1)$}
Let $G=O(n+1,1)$ and $m=[\frac n 2]$.  
We recall that $W_G \simeq {\mathfrak{S}}_{m+1} \ltimes ({\mathbb{Z}}/2{\mathbb{Z}})^{m+1}$
 acts on ${\mathfrak{h}}_{\mathbb{C}}^{\ast} \simeq {\mathbb{C}}^{m+1}$
 as a permutation group 
 and by switching the signatures of the standard coordinates.  
For $\tau \in {\mathbb{Z}}^{m+1}$, 
 we define 
\index{A}{1pataudom@$\tau_{\operatorname{dom}}$|textbf}
$\tau_{\operatorname{dom}}$ to be the unique element
 in $\Lambda^+(m+1)$
 (see \eqref{eqn:Lambda})
 in the $W_G$-orbit through $\tau$, 
{\it{i.e.}}, 
\begin{equation}
\label{eqn:Wdom}
  \tau_{\operatorname{dom}}= w \, \tau
  \qquad
  \text{for some $w \in W_G$.  }
\end{equation}
Let $\Kirredrep {O(n+1,1)}{\tau_{\operatorname{dom}}}_{+,+}$
 be the irreducible finite-dimensional representation of 
 $G=O(n+1,1)$ of type I 
 (Definition \ref{def:typeone})
 defined as in \eqref{eqn:ONCreal}.

\begin{definition}
[translation functor $\psi_{\mu}^{\mu+\tau}$]
\label{def:transG}
For $\mu \in {\mathbb{C}}^{m+1}$ and $\tau \in {\mathbb{Z}}^{m+1}$, 
 we define translation functor $\psi_{\mu}^{\mu+\tau}$ by 
\index{A}{1psinlmutau@$\psi_{\mu}^{\mu+\tau}$|textbf}
\begin{equation}
\label{eqn:translation}
   \psi_{\mu}^{\mu+\tau}(\Pi)
   :=
   P_{\mu+\tau}(P_{\mu}(\Pi) \otimes \Kirredrep{O(n+1,1)}{\tau_{\operatorname{dom}}}_{+,+}).  
\end{equation}
\end{definition}
Then $\psi_{\mu}^{\mu+\tau}$ is a covariant functor
 in the category of admissible smooth representations of $G$
 of finite length, 
 and also  
 in the category
 of $({\mathfrak{g}},K)$-modules of finite length.  
Clearly,
 we have 
\begin{equation}
\label{eqn:V731}
  \psi_{\mu}^{\mu+\tau}=\psi_{w\mu}^{w\mu+w\tau}
\qquad
\text{for all $w \in W_G$.}
\end{equation}

In defining the translation functor $\psi_{\mu}^{\mu+\tau}$
 in \eqref{eqn:translation}, 
 we have used only finite-dimensional representations
\index{B}{type1indef@type I, representation of $O(N-1,1)$}
 of type I
 (Definition \ref{def:typeone})
 of the disconnected group $G=O(n+1,1)$.  
We do not lose any generality
 because taking the tensor product
 with the one-dimensional characters
 $\chi_{a b}$
 ($a, b \in \{\pm\}$)
 yields the following isomorphism as $G$-modules:
\begin{equation}
\label{eqn:chipsi}
   \psi_{\mu}^{\mu+\tau}(\Pi) \otimes \chi_{a b}
   \simeq 
   P_{\mu+\tau}
   (P_{\mu}(\Pi) \otimes 
    \Kirredrep{O(n+1,1)}{\tau_{\operatorname{dom}}}_{a,b}).  
\end{equation}
We shall use a finite-dimensional representation
 $F(V,\lambda)$
 (Definition \ref{def:Fshift})
 which is not necessarily of type I
 in Theorems \ref{thm:181104} and \ref{thm:181107}, 
 which are a reformulation of the properties
 (Theorems \ref{thm:1807113} and \ref{thm:1808101}, respectively)
 of the translation functor
 \eqref{eqn:translation}
 via \eqref{eqn:chipsi}.

The translation functor $\psi_{\mu+\tau}^{\mu}$ is the adjoint functor
 of $\psi_{\mu}^{\mu+\tau}$.  
In our setting, 
 since $(-\tau)_{\operatorname{dom}}=\tau_{\operatorname{dom}}$, 
 the functor $\psi_{\mu+\tau}^{\mu}$
 takes the following form:
\[
  \psi_{\mu+\tau}^{\mu}(\Pi)
  =
  P_{\mu}
  (P_{\mu+\tau}(\Pi) \otimes 
   \Kirredrep{O(n+1,1)}{\tau_{\operatorname{dom}}}_{+,+}).  
\]

\subsubsection{The translation functor and the restriction $G \downarrow \overline G$}
We retain the notation of Appendix II, 
 and denote by $\overline G$ the subgroup 
 $SO(n+1,1)$ in $G=O(n+1,1)$.  
Then $\overline G=SO(n+1,1)$ is in the Harish-Chandra class
 for all $n$.  
For the group $\overline G$,
 we shall use the notation
 $\overline P_{\mu}$ and $\overline{\psi}_{\mu}^{\mu+\tau}$ 
instead of $P_{\mu}$ and $\psi_{\mu}^{\mu+\tau}$, 
 respectively.   
To be precise,
 for $\tau \in {\mathbb{Z}}^{m+1}$ 
 where $m=[\frac n 2]$, 
 we write  $\overline \tau_{\operatorname{dom}}$ for the unique element
 in the orbit $W_{\mathfrak{g}}\, \tau$
 which is dominant with respect to the positive system
 $\Delta^+({\mathfrak{g}}_{\mathbb{C}}, {\mathfrak{h}}_{\mathbb{C}})$.  
We denote by $\Kirredrep{SO(n+1,1)}{\overline \tau_{\operatorname{dom}}}_+$
 the irreducible representation of $\overline G=SO(n+1,1)$
 obtained by the restriction
 of the irreducible holomorphic representation of $SO(n+2,{\mathbb{C}})$
 having $\overline \tau_{\operatorname{dom}}$
 as its highest weight.  
For an admissible smooth representation $\overline \Pi$
 of $\overline G$ of finite length, 
 the translation functor $\overline \psi_{\mu}^{\mu+\tau}$
 is defined by
\begin{equation}
\label{eqn:transSO}
\overline \psi_{\mu}^{\mu+\tau}(\overline \Pi)
:=
\overline P_{\mu+\tau}
(\overline P_{\mu}(\overline \Pi) 
\otimes 
\Kirredrep{SO(n+1,1)}{\overline \tau_{\operatorname{dom}}}_+).  
\end{equation}
We collect some basic facts concerning the primary components
 for $G$-modules and $\overline G$-modules. 
The following lemma is readily shown by comparing \eqref{eqn:HCpara}
 of the Harish-Chandra isomorphisms
 for $G$ and $\overline G$.  
\begin{lemma}
\label{lem:primary}
Let $\Pi$ be an admissible smooth representation
of finite length of $G=O(n+1,1)$.  
We set $m:=[\frac n 2]$ as before.  
Suppose $\mu \in {\mathfrak{h}}_{\mathbb{C}}^{\ast} 
\simeq {\mathbb{C}}^{m+1}$.  
\begin{enumerate}
\item[{\rm{(1)}}]
If $n$ is odd or if $n$ is even
 and at least one of the entries $\mu_1$, $\cdots$, $\mu_{m+1}$
 is zero, 
then there is a natural isomorphism of $\overline G$-modules:
\[
   P_{\mu}(\Pi)|_{\overline G} \simeq \overline P_{\mu}(\Pi|_{\overline G}).  
\]
\item[{\rm{(2)}}]
If $n$ is even
 and all of $\mu_j$ are nonzero, 
 then we have a direct sum decomposition of a $\overline G$-module:
\[
   P_{\mu}(\Pi)|_{\overline G}=\overline P_{\mu}(\Pi|_{\overline G}) \oplus 
\overline P_{\mu'}(\Pi|_{\overline G}), 
\]
where we set $\mu':=(\mu_1, \cdots, \mu_m, -\mu_{m+1})$.  
\end{enumerate}
\end{lemma}

Now the following lemma is an immediate consequence
 of Lemma \ref{lem:primary}
 and of the definition of the translation functors
 $\psi_{\mu}^{\mu+\tau}$ and $\overline \psi_{\mu}^{\mu+\tau}$, 
 see \eqref{eqn:translation} and \eqref{eqn:transSO}. 
\begin{lemma}
\label{lem:Xtrans}
Let $G=O(n+1,1)$ and $\overline G=SO(n+1,1)$.  
Let $\Pi$ be an admissible smooth representation of $G$
 of finite length.  
\begin{enumerate}
\item[{\rm{(1)}}]
Suppose $n$ is odd.  
Then we have a canonical $\overline G$-isomorphism:
\begin{equation}
\label{eqn:restrans}
\psi_{\mu}^{\mu+\tau}(\Pi)|_{\overline G}
\simeq
\overline \psi_{\mu}^{\mu+\tau}(\Pi|_{\overline G}).  
\end{equation}
\item[{\rm{(2)}}]
Suppose that $n$ is even.  
If all of $\mu$, $\tau$ and $\mu+\tau$ contain 0 in their entries,
 then we have a canonical $\overline G$-isomorphism
 \eqref{eqn:restrans}.  
\end{enumerate}
\end{lemma}

\subsubsection{Some elementary properties
 of translation functor $\psi_{\mu}^{\mu+\tau}$}
Some of the properties of the translation functors remain true
 for the disconnected group $G=O(n+1,1)$.  
\begin{proposition}
\label{prop:transfunct}
Suppose $\mu \in {\mathfrak{h}}_{\mathbb{C}}^{\ast} (\simeq {\mathbb{C}}^{m+1})$
 and $\tau \in {\mathbb{Z}}^{m+1}$.  
\begin{enumerate}
\item[{\rm{(1)}}]
$\psi_{\mu}^{\mu+\tau}$ is a covariant exact functor.  
\item[{\rm{(2)}}]
Suppose $\mu$ and $\mu+\tau$ belong to the same Weyl chamber
 with respect to $W_{\mathfrak{g}}$.  
If $\mu+\tau$ is regular 
 (Definition \ref{def:intreg}), 
 then $\psi_{\mu}^{\mu+\tau}(\Pi)$ is nonzero 
 if $\Pi$ is nonzero.  
\end{enumerate}
\end{proposition}

\begin{proof}
(1)\enspace
The first statement follows directly from  the definition,
 see Zuckerman \cite{Zuckerman}.  
\par\noindent
(2)\enspace
By Lemma \ref{lem:primary} and the branching law from $G=O(n+1,1)$
 to the subgroup $\overline G=SO(n+1,1)$, 
 we have
\[
 \psi_{\mu}^{\mu+\tau}(\Pi)|_{\overline G}
  \supset 
  \overline \psi_{\mu}^{\mu+\tau}(\Pi|_{\overline G}).  
\]
Since ${\overline G}$ is in the Harish-Chandra class,
 $\overline \psi_{\mu}^{\mu+\tau}(\Pi|_{\overline G})$
 is nonzero under the assumption
 on $\mu$ and $\tau$.  
Hence $\psi_{\mu}^{\mu+\tau}(\Pi)$ is a nonzero $G$-module.  
\end{proof}

\begin{remark}
The regularity assumption 
 for $\mu + \tau$
 in Proposition \ref{prop:transfunct}
 is in the weaker sense
 ({\it{i.e.}}, $W_{\mathfrak{g}}$-regular), 
 and not in the stronger sense
 ({\it{i.e.}}, $W_G$-regular). 
\end{remark}

\subsection{Translation of principal series representation $I_{\delta}(V,\lambda)$}

We discuss how the translation functors affect 
 induced representations
 of $G=O(n+1,1)$.  
We recall that $G$ is not in the Harish-Chandra class when $n$ is even.  

\subsubsection{Main results: Translation of principal series representations}
\begin{theorem}
\label{thm:1807113}
Suppose $G=O(n+1,1)$
 and $(V,\lambda) \in \Reducible$, 
 see \eqref{eqn:reducible}, 
 or equivalently,
 $V \in \widehat {O(n)}$
 and $\lambda \in {\mathbb{Z}} \setminus (S(V) \cup S_Y(V))$, 
 see Theorem \ref{thm:irrIV}.  
Let $i:=i(V,\lambda) \in \{0,1,\ldots,n\}$ be the height of $(V,\lambda)$
 as in \eqref{eqn:indexV}, 
 and $r(V, \lambda)\in {\mathbb{Z}}^{m+1}$ as in \eqref{eqn:IVZG}.  
We write $V=\Kirredrep{O(n)}{\sigma}_{\varepsilon}$
 with $\sigma \in \Lambda^+(m)$ and $\varepsilon \in \{\pm\}$, 
 where $m:=[\frac n 2]$.  
We define a character $\chi$ of $G$ by
\index{A}{1chipmpmVlmd@$\chi(V,\lambda)$|textbf}
\begin{equation}
\label{eqn:sgnchi}
\chi \equiv \chi(V,\lambda):=
\begin{cases}
{\bf{1}}
&\text{if $\varepsilon(\frac n 2-i) \ge 0$,}
\\
\det
\qquad
&\text{if $\varepsilon(\frac n 2-i) < 0$.  }
\end{cases}
\end{equation}
Then there is a natural $G$-isomorphism:
\[
  \psi_{\rho^{(i)}}^{r(V,\lambda)}(I_{\delta}(i,i)) \otimes \chi
  \simeq 
  I_{(-1)^{\lambda-i}\delta}(V,\lambda).  
\]
\end{theorem}

\begin{remark}
\label{rem:FVchipm}
The conclusion of Theorem \ref{thm:1807113} does not change
 if we replace the definition \eqref{eqn:sgnchi}
 with $\chi=\det$
 when $i=\frac n 2$.  
In fact, 
$V$ is of type Y
 if the height $i(V,\lambda)$ equals $\frac n 2$, 
and thus $V \otimes \det \simeq V$
 as $O(n)$-modules
 (Lemma \ref{lem:typeY}).  
Then there is an isomorphism of $G$-modules
\[
   I_{\delta}(V,\lambda) \otimes \det \simeq I_{\delta}(V,\lambda)
\]
for any $\delta \in \{\pm\}$
 by Lemmas \ref{lem:IVchi} and \ref{lem:isig}.  
\end{remark}

The translation functor 
 $\psi_{r(V,\lambda)}^{\rho^{(i)}}$
 is the adjoint functor
 of $\psi_{\rho^{(i)}}^{r(V,\lambda)}$.  
Even when the infinitesimal character
 of $I_{\delta}(V,\lambda)$
 is $W_G$-regular
 (in particular,
 $W_{\mathfrak{g}}$-regular)
 (Definition \ref{def:intreg}), 
 the translation functor $\psi_{r(V,\lambda)}^{\rho^{(i)}}$
 does not always preserve {\it{irreducibility}}
 if $G$ is not of Harish-Chandra class as in the following theorem.  

\begin{theorem}
\label{thm:1808101}
Retain the setting and notation of Theorem \ref{thm:1807113}.  
In particular,
 we recall that 
\index{A}{RzqSeducible@$\Reducible$ $(\subset \widehat {O(n)} \times {\mathbb{Z}})$}
$(V,\lambda)\in \Reducible$, 
\index{A}{iVlmd@$i(V,\lambda)$, height}
 $i=i(V,\lambda)$ is the height of $(V,\lambda)$, 
 and $\chi\equiv \chi(V,\lambda)$, 
 see \eqref{eqn:sgnchi}.  
\begin{enumerate}
\item[{\rm{(1)}}]
If 
\index{A}{RzqSeducibleI@$\RedI$}
$(V,\lambda)\in \RedI$
 (Definition \ref{def:Red12}), 
 {\it{i.e.}}, 
 if $V$ is of type X
 (in particular, 
 if $n$ is odd)
 or if $\lambda=\frac n 2$, 
 then there is a natural $G$-isomorphism:
\[
  \psi_{r(V,\lambda)}^{\rho^{(i)}}(I_{(-1)^{\lambda-i}\delta}(V,\lambda))
  \otimes \chi
  \simeq 
  I_{\delta}(i,i).  
\]
\item[{\rm{(2)}}]
If 
\index{A}{RzqSeducibleII@$\RedJ$}
$(V,\lambda)\in \RedJ$, 
 {\it{i.e.}}, 
 if $V$ is of type Y 
 and $\lambda \ne \frac n 2$, 
 then $n$ is even, 
 $i \ne \frac n 2$
 and there is a natural $G$-isomorphism:
\[
  \psi_{r(V,\lambda)}^{\rho^{(i)}}(I_{(-1)^{\lambda-i}\delta}(V,\lambda))
  \otimes \chi
  \simeq 
  I_{\delta}(i,i) \oplus I_{\delta}(n-i,i).  
\]
\end{enumerate}
\end{theorem}
In Section \ref{sec:FVlmd}, 
 we introduce an irreducible finite-dimensional representation
 $F(V,\lambda)$
 by taking the tensor product
 of $\Kirredrep{O(n+1,1)}{\tau_{\operatorname{dom}}}_{+,+}$
 with an appropriate one-dimensional character of $G$, 
 see Definition \ref{def:Fshift}.  
Then, 
 by using $F(V,\lambda)$, 
 Theorems \ref{thm:1807113} and \ref{thm:1808101}
 can be reformulated in a simpler form about signatures, 
 see Theorems \ref{thm:181104} and \ref{thm:181107}.

\subsubsection{Strategy of the proof
 for Theorems \ref{thm:1807113} and \ref{thm:1808101}}
\label{subsec:sketchtrans}
If $n$ is odd,
 then $G=\langle SO(n+1,1), -I_{2n+2} \rangle$ is in the Harish-Chandra class, 
 and therefore Theorems \ref{thm:1807113} and \ref{thm:1808101}
 are a special case of the general theory, 
 see \cite[Chap.~7]{Vogan81} for instance.  
Moreover,
 the translation functor
 behaves
 as we expect from the general theory 
 for reductive groups in the Harish-Chandra class
 when it is applied to the induced representation
 $I_{\delta}(V,\lambda)$
 if $(V,\lambda) \in \RedI$, 
 see Theorem \ref{thm:1808101} (1).  
We note that $\Reducible=\RedI$ and $\RedJ=\emptyset$
 if $n$ is odd
 (Remark \ref{rem:Red12}).

On the other hand,
 its behavior is somewhat different
 if $(V,\lambda) \in \RedJ$, 
 see Theorem \ref{thm:1808101} (2)
 and Proposition \ref{prop:transrest2}
 for instance.  
Main technical complications arise from 
 the fact that we need the primary decomposition 
 for the generalized ${\mathfrak{Z}}_G({\mathfrak{g}})$-infinitesimal characters
 parametrized by ${\mathfrak{h}}_{\mathbb{C}}^{\ast}/W_G$
 where $W_G$ is larger than the group 
 generated by the reflections
 of simple roots
 if $n$ is even, 
 for which $G=O(n+1,1)$ is not in the Harish-Chandra class.

Our strategy is to use partly 
 the relation 
 of translation functors for $G=O(n+1,1)$
 and the subgroup 
 $\overline G=SO(n+1,1)$
 which is in the Harish-Chandra class.

Theorem \ref{thm:1807113} is proved in Section \ref{subsec:pftran1}
 as a consequence
 of the following two propositions.  

\begin{proposition}
\label{prop:IVsub}
Suppose that $(V,\lambda) \in \Reducible$.  
Retain the notation as in Theorem \ref{thm:1807113}.  
Then the $G$-module 
 $\psi_{\rho^{(i)}}^{r(V,\lambda)}(I_{\delta}(i,i))\otimes \chi$ contains
 $I_{(-1)^{\lambda-i}\delta}(V,\lambda)$ as a subquotient.   
Equivalently, the $G$-module $P_{r(V,\lambda)}(I_{\delta}(i,i) \otimes F(V,\lambda))$, 
 see Definition \ref{def:Fshift} below, 
 contains $I_{\delta}(V, \lambda)$ as a subquotient.  
\end{proposition}

We recall from \eqref{eqn:sgnchi}
 that the character 
\index{A}{1chipmpmVlmd@$\chi(V,\lambda)$}
$\chi \equiv \chi(V,\lambda)$
 is trivial 
 when restricted to the subgroup $\overline G=SO(n+1,1)$.  

\begin{proposition}
\label{prop:transrest}
Suppose that $(V,\lambda)\in \Reducible$.  
Retain the notation as in Theorem \ref{thm:1807113}.  
Then $\psi_{\rho^{(i)}}^{r(V,\lambda)}(I_{\delta}(i,i))|_{\overline G}$
 is isomorphic to 
 $I_{(-1)^{\lambda-i}\delta}(V,\lambda)|_{\overline G}$
 as a $\overline G$-module.  
\end{proposition}

Similarly, 
 Theorem \ref{thm:1808101} is proved 
in Section \ref{subsec:pf168}
 by using analogous results, 
 namely,
 Propositions \ref{prop:Iisub} and \ref{prop:PIFrest}
 in Section \ref{subsec:pf168}.

\subsubsection{Basic lemmas for the translation functor}
\label{subsec:tensorcomb}
We use the following well-known lemma, 
 which holds 
 without the assumption that $G$ is of Harish-Chandra class.  
\begin{lemma}
\label{lem:Ftensor}
Let $F$ be a finite-dimensional representation
 of $G$, 
 $V \in \widehat{O(n)}$, 
 $\delta \in \{\pm\}$, 
 and $\lambda \in {\mathbb{C}}$.  
Then there is a $G$-stable filtration
\[
  \{0\} = I_0 \subset I_1 \subset \cdots \subset I_k 
 = I_{\delta}(V, \lambda) \otimes F
\]
such that
\[
 I_j /I_{j-1} \simeq \operatorname{Ind}_P^G(V_{\lambda,\delta} \otimes F^{(j)})
\quad
  (1 \le j \le k)
\]
where $F^{(j)}$ is a $P$-module
 such that the unipotent radical $N_+$ acts trivially
 and that $F^{(j)}|_{M A}$ is isomorphic to a subrepresentation
 of the restriction $F|_{M A}$ to the Levi subgroup $M A$.  
\end{lemma}
For the sake of completeness, 
we give a proof.  
\begin{proof}
Take a $P$-stable filtration
\[
   \{0\} = F_0 \subset F_1 \subset \cdots \subset F_k = F
\]
 such that the unipotent radical $N_+$ of $P$ acts trivially
 on 
\[
  F^{(j)}:=F_j/F_{j-1}
\quad (1 \le j \le k).  
\]

As in \eqref{eqn:Vlmddelta}, 
 we denote by 
\index{A}{Vln@$V_{\lambda,\delta}=V \otimes \delta \otimes {\mathbb{C}}_{\lambda}$, representation of $P$}
$V_{\lambda,\delta}$
 the irreducible $P$-module
 which is an extension of the $M A$-module
 $V \boxtimes \delta \boxtimes {\mathbb{C}}_{\lambda}$
 with trivial $N_+$ action.  
We define $G$-modules $I_j$ ($0 \le j \le k$)
 by 
\[
  I_j:=\operatorname{Ind}_P^G(V_{\lambda,\delta} \otimes F_j|_P).  
\]
Then there is a natural filtration of $G$-modules
\[
 0 =
 I_0
\subset
I_1
\subset
\cdots
\subset
I_k=\operatorname{Ind}_P^G(V_{\lambda,\delta} \otimes F|_P)
\]
such that 
\[
I_j/I_{j-1}
\simeq
\operatorname{Ind}_P^G(V_{\lambda,\delta} \otimes (F_j/F_{j-1}))
\]
as $G$-modules.  
Since the finite-dimensional representation $F$ of $G$
 is completely reducible
 when viewed as a representation of the Levi subgroup $M A$, 
 the $M A$-module 
$
   F^{(j)}=F_j/F_{j-1}
$
 is isomorphic
 to a subrepresentation of the restriction $F|_{M A}$.  
Now Lemma \ref{lem:Ftensor} follows from the following $G$-isomorphism:
\[
 \operatorname{Ind}_P^G(V_{\lambda,\delta} \otimes F|_P)
 \simeq
 \operatorname{Ind}_P^G(V_{\lambda,\delta}) \otimes F.  
\]
\end{proof}
Similarly to Lemma \ref{lem:Ftensor}, 
 we have the following lemma
 for cohomological parabolic induction.  
Retain the notation
 as in Section \ref{subsec:Aqgeneral}.  
\begin{lemma}
\label{lem:Rqtensor}
Suppose that ${\mathfrak{q}}={\mathfrak{l}}_{\mathbb{C}}+{\mathfrak{u}}$
 is a $\theta$-stable parabolic subalgebra
 of ${\mathfrak{g}}_{\mathbb{C}}$
 with Levi subgroup $L$, 
 see \eqref{eqn:LeviLq}, 
 and that $W$ a finite-dimensional $({\mathfrak{l}}, L \cap K)$-module.  
Let $F$ be a finite-dimensional representation of $G$, 
 and 
\[
  \{0\}= F_0 \subset F_1 \subset \cdots \subset F_k =F
\]
 a $({\mathfrak{q}},L)$-stable filtration 
 such that the nilpotent radical ${\mathfrak{u}}$
 acts trivially on $F^{(j)}:=F_j/F_{j-1}$.  
Then there is a natural spectral sequence
\index{A}{RzqS@${\mathcal{R}}_{\mathfrak{q}}^S$, cohomological parabolic induction}
\[
  {\mathcal{R}}_{\mathfrak{q}}^p(W \otimes F^{(j)} \otimes {\mathbb{C}}_{\rho({\mathfrak{u}})})
  \Rightarrow
  {\mathcal{R}}_{\mathfrak{q}}^p(W \otimes F \otimes {\mathbb{C}}_{\rho({\mathfrak{u}})})
  \simeq 
  {\mathcal{R}}_{\mathfrak{q}}^p(W \otimes {\mathbb{C}}_{\rho({\mathfrak{u}})})
  \otimes F
\]
 as $({\mathfrak{g}},K)$-modules.  
\end{lemma}

The proof is similar to the case 
 where $G$ is in the Harish-Chandra class, 
 see \cite[Lem.~7.23]{Vogan81}.  

\vskip 1pc
By the definition \eqref{eqn:translation}
 of the translation functor $\psi_{\mu}^{\mu+\tau}$, 
 we need to estimate
 possible ${\mathfrak{Z}}_G({\mathfrak{g}})$-infinitesimal characters
 of $\operatorname{Ind}_P^G(V_{\lambda,\delta} \otimes F^{(j)})$
 in Lemma \ref{lem:Ftensor}
 or that of ${\mathcal{R}}_{\mathfrak{q}}^p(W \otimes F^{(j)} \otimes {\mathbb{C}}_{\rho({\mathfrak{u}})})$
 in Lemma \ref{lem:Rqtensor}.

In order to deal with reductive groups
that are not in the Harish-Chandra class, 
 we use the following lemma 
 which is formulated in a slightly stronger form
 than \cite[Lem.~7.2.18]{Vogan81}, 
but has the same proof.  

\begin{lemma}
\label{lem:V7218}
Let ${\mathfrak{h}}_{\mathbb{C}}$ be a Cartan subalgebra
 of a complex semisimple Lie algebra ${\mathfrak{g}}_{\mathbb{C}}$, 
 $W_{\mathfrak{g}}$ the Weyl group 
 of the root system $\Delta({\mathfrak{g}}_{\mathbb{C}}, {\mathfrak{h}}_{\mathbb{C}})$, 
 $\Delta^+({\mathfrak{g}}_{\mathbb{C}}, {\mathfrak{h}}_{\mathbb{C}})$ 
 a positive system,
 $\langle \, , \, \rangle$ 
 a $W_{\mathfrak{g}}$-invariant inner product 
 on ${\mathfrak{h}}_{\mathbb{R}}^{\ast}:=\operatorname{Span}_{\mathbb{R}} \Delta ({\mathfrak{g}}_{\mathbb{C}}, {\mathfrak{h}}_{\mathbb{C}})$, 
 and $||\,\cdot\,||$ its norm.

Suppose that $\nu$ and $\tau \in {\mathfrak{h}}_{\mathbb{R}}^{\ast}$
 satisfy
\begin{alignat*}{2}
\langle \nu, \alpha^{\vee} \rangle &\in {\mathbb{N}}_+
\qquad
&&({}^{\forall} \alpha \in \Delta^+({\mathfrak{g}}_{\mathbb{C}}, {\mathfrak{h}}_{\mathbb{C}})), 
\\
\langle \nu+\tau, \alpha^{\vee} \rangle &\in {\mathbb{N}}
\qquad
&&({}^{\forall} \alpha \in \Delta^+({\mathfrak{g}}_{\mathbb{C}}, {\mathfrak{h}}_{\mathbb{C}})).  
\end{alignat*}

If $\gamma \in {\mathfrak{h}}_{\mathbb{C}}^{\ast}$
 satisfies the following two conditions:
\begin{alignat}{2}
\nu + \gamma =& w (\nu+\tau)
\qquad
&&\text{for some $w  \in W_{\mathfrak{g}}$}, 
\label{eqn:V81a}
\\
|| \gamma || \le& ||\tau||, 
\qquad
&&
\label{eqn:V81b}
\end{alignat}
then $\gamma = \tau$.  
\end{lemma}

\begin{remark}
In \cite[Lem.~7.2.18]{Vogan81}, 
$\gamma$ is assumed to be a weight
 occurring in the irreducible finite-dimensional representation
 of $G$ 
 (in the Harish-Chandra class)
 with extremal weight $\tau$
 instead of our assumption \eqref{eqn:V81b}.  
\end{remark}

\subsection{Definition of an irreducible finite-dimensional
\\
 representation $F(V,\lambda)$ of $G=O(n+1,1)$
\label{sec:FVlmd}}
For $(V,\lambda)\in \RInt$, 
 {\it{i.e.}}, 
  for $V \in \widehat {O(n)}$
 and $\lambda \in {\mathbb{Z}} \setminus S(V)$, 
 we defined in Chapter \ref{sec:aq}
\index{B}{ZGginfinitesimalcharacter@${\mathfrak{Z}}_G({\mathfrak{g}})$-infinitesimal character}
\begin{alignat*}{2}
& i(V,\lambda) \in \{0,1,\ldots,n\}, 
\quad
&&\text{height of $(V,\lambda)$
 (Definition \ref{def:iVlmd})}, 
\\
& r(V,\lambda) \in {\mathbb{C}}^{[\frac n 2]+1}, 
&&\text{giving the ${\mathfrak{Z}}_G({\mathfrak{g}})$-infinitesimal character 
 of $I_{\delta}(V,\lambda)$, }
\\
&
&&\text{see \eqref{eqn:IVZG}.}
\end{alignat*}
In this section 
 we introduce an irreducible finite-dimensional representation
 $F(V,\lambda)$ of $G=O(n+1,1)$
 which contains important information
 on signatures.  

\subsubsection{Definition of $\sigma^{(i)}(\lambda)$
 and $\widehat{\sigma^{(i)}}$}
We begin with some combinatorial notation.  
\begin{definition}
\label{def:sigmailmd}
Let $m:=[\frac n 2]$.  
For $1 \le i \le n$, 
 $\sigma=(\sigma_1, \cdots,\sigma_m) \in \Lambda^+(m)$, 
 and $\lambda \in {\mathbb{Z}}$, 
 we define $\sigma^{(i)}(\lambda) \in {\mathbb{Z}}^{m+1}$
 as follows.  
\par\noindent
{\bf{\underline{Case 1.}}}\enspace
$n=2m$
\begin{equation*}
\sigma^{(i)}(\lambda):=
\begin{cases}
(\sigma_{1}-1,\cdots,\sigma_{i}-1,i-\lambda,\sigma_{i+1},\cdots,\sigma_{m})
\quad
&\text{for $0 \le i \le m-1$,}
\\
(\sigma_{1}-1,\cdots,\sigma_{m}-1,|\lambda-m|)
\quad
&\text{for $i=m$,}
\\
(\sigma_{1}-1,\cdots,\sigma_{n-i}-1,\lambda-i,\sigma_{n-i+1},\cdots,\sigma_{m})
\quad
&\text{for $m+1 \le i \le n$.  }
\end{cases}
\end{equation*}

\par\noindent
{\bf{\underline{Case 2.}}}\enspace
$n=2m+1$
\begin{equation*}
\sigma^{(i)}(\lambda):=
\begin{cases}
(\sigma_{1}-1,\cdots,\sigma_{i}-1,i-\lambda,\sigma_{i+1},\cdots,\sigma_{m})
\quad
&\text{for $0 \le i \le m$,}
\\
(\sigma_{1}-1,\cdots,\sigma_{n-i}-1,\lambda-i,\sigma_{n-i+1}, \cdots,\sigma_m)
\quad
&\text{for $m+1 \le i \le n$.  }
\end{cases}
\end{equation*}
Moreover we define $\widehat{\sigma^{(i)}} \in {\mathbb{Z}}^m$
 to be the vector obtained by removing the 
 $\min(i+1, n-i+1)$-th component from $\sigma^{(i)}(\lambda)\in {\mathbb{Z}}^{m+1}$.  
\end{definition}

\par\noindent
{\bf{\underline{Case 1.}}}\enspace
$n=2m$
\begin{equation*}
\widehat{\sigma^{(i)}}:=
\begin{cases}
(\sigma_{1}-1,\cdots,\sigma_{i}-1,\sigma_{i+1},\cdots,\sigma_{m})
\quad
&\text{for $0 \le i \le m-1$,}
\\
(\sigma_{1}-1,\cdots,\sigma_{m}-1)
\quad
&\text{for $i=m$,}
\\
(\sigma_{1}-1,\cdots,\sigma_{n-i}-1,\sigma_{n-i+1},\cdots,\sigma_{m})
\quad
&\text{for $m+1 \le i \le n$.  }
\end{cases}
\end{equation*}

\par\noindent
{\bf{\underline{Case 2.}}}\enspace
$n=2m+1$
\begin{equation*}
\widehat{\sigma^{(i)}}:=
\begin{cases}
(\sigma_{1}-1,\cdots,\sigma_{i}-1,\sigma_{i+1},\cdots,\sigma_{m})
\quad
&\text{for $0 \le i \le m$,}
\\
(\sigma_{1}-1,\cdots,\sigma_{n-i}-1,\sigma_{n-i+1}, \cdots,\sigma_m)
\quad
&\text{for $m+1 \le i \le n$.  }
\end{cases}
\end{equation*}

\begin{definitionlemma}
\label{deflem:180694}
Let $m:=[\frac n 2]$.  
For $(V, \lambda) \in \RInt$, 
 {\it{i.e.}}, 
 for $V \in \widehat{O(n)}$
 and $\lambda \in {\mathbb{Z}} \setminus S(V)$, 
 we write $V=\Kirredrep{O(n)}{\sigma}_{\varepsilon}$
 with $\sigma \in \Lambda^+(m)$
 and $\varepsilon \in \{\pm\}$.  
We set
\begin{equation}
\label{eqn:sigmalmd}
\sigma(\lambda) := \sigma^{(i)}(\lambda), 
\end{equation}
where $i:=i(V,\lambda) \in \{0,1,\ldots,n\}$
 is the height of $(V,\lambda)$ as in \eqref{eqn:indexV}.  
Then we have
\[
    \sigma(\lambda) \in \Lambda^+(m+1).  
\]
\end{definitionlemma}
\begin{proof}
Suppose $n=2m$ (even).  
Let $\lambda \in {\mathbb{Z}}$.  
By the definition of $R(V;i)$ (Definition \ref{def:RVi}), 
 we have the following equivalences:
\vskip 1pc
\par\noindent
$\bullet$\enspace for $0 \le i \le m-1$, 
\begin{align*}
\lambda \in R(V;i)
\,\,\Leftrightarrow\,\, & \sigma_i-i> -\lambda > \sigma_{i+1}-i-1
\\
\,\,\Leftrightarrow\,\, & \sigma_i-1 \ge i -\lambda \ge  \sigma_{i+1};
\intertext{$\bullet$\enspace for $i=m$, }
\lambda \in R(V;m)
\,\,\Leftrightarrow\,\, & -\sigma_m <\lambda-m< \sigma_{m}
\\
\,\,\Leftrightarrow\,\, & \sigma_m-1 \ge |\lambda -m|;
\intertext{$\bullet$\enspace for $m+1 \le i \le n$, }
\lambda \in R(V;i)
\,\,\Leftrightarrow\,\, & \sigma_{n-i+1}-1 <\lambda-i< \sigma_{n-i}
\\
\,\,\Leftrightarrow\,\, & \sigma_{n-i}-1 \ge \lambda -i \ge \sigma_{n-i+1}.  
\end{align*}
Thus in all cases $\sigma^{(i)}(\lambda) \in \Lambda^+(m+1)$.

The proof for $n$ odd is similar.  
\end{proof}

\subsubsection{Definition of a finite-dimensional representation $F(V,\lambda)$
 of $G$}
We are ready to define a finite-dimensional representation, 
 to be denoted by $F(V,\lambda)$, 
 for $(V,\lambda)\in \RInt$.  

\begin{definition}
[a finite-dimensional representation $F(V,\lambda)$]
\label{def:Fshift}
Suppose that  $(V,\lambda) \in \RInt$,
{\it{i.e.}},
 $V \in \widehat{O(n)}$ and $\lambda \in {\mathbb{Z}} \setminus S(V)
$.  
We write $V=\Kirredrep{O(n)}{\sigma}_{\varepsilon}$
 with $\sigma \in \Lambda^+(m)$
 and $\varepsilon \in \{\pm\}$
 where $m:=[\frac n 2]$.  
We set $i:=i(V,\lambda)$, 
 the height of $(V,\lambda)$
 as in \eqref{eqn:indexV}, 
 and 
\index{A}{1pasigmalmd@$\sigma(\lambda)$ $(=\sigma^{(i)}(\lambda))$|textbf}
$
\sigma(\lambda) \in \Lambda^+(m+1)
$
as in Definition-Lemma \ref{deflem:180694}.

We define an irreducible finite-dimensional representation
 $F(V, \lambda)$ of $G=O(n+1,1)$
 as follows:
\begin{enumerate}
\item[$\bullet$]
for $V$ of type Y
 and $\lambda=\frac n 2 (=m)$, 
\begin{align*}
F(V,\lambda):=& \Kirredrep{O(n+1,1)}{\sigma(\lambda)}_{+,+}
\\
             =& \Kirredrep{O(n+1,1)}{\sigma_1-1, \cdots,\sigma_m-1,0}_{+,+};
\end{align*}
\item[$\bullet$]
for $V$ of type X
 or $\lambda \ne \frac n 2$, 
\index{A}{FVlmd@$F(V,\lambda)$|textbf}
\begin{equation}
\label{eqn:FVlmd}
   F(V,\lambda)
:=
\begin{cases}
\Kirredrep {O(n+1,1)}
           {\sigma(\lambda)}_{\varepsilon,(-1)^{\lambda-i}\varepsilon}
\quad
&\text{if $i \le \frac n 2$}, 
\\
\Kirredrep {O(n+1,1)}
           {\sigma(\lambda)}_{-\varepsilon,(-1)^{\lambda-i-1}\varepsilon}
\quad
&\text{if $i > \frac n 2$}, 
\end{cases}
\end{equation}
see \eqref{eqn:Fn1ab} for notation.  
\end{enumerate}

By using the character $\chi \equiv \chi(V,\lambda)$ of $G$
 as defined in \eqref{eqn:sgnchi}, 
 we obtain a unified expression
\begin{equation}
\label{eqn:FVchi}
F(V,\lambda) \simeq F(\sigma(\lambda))_{+,(-1)^{\lambda-i}} \otimes \chi.
\end{equation}
\end{definition}

\begin{remark}
\label{rem:FVlmd}
We note
 that \eqref{eqn:FVlmd} is well-defined.  
In fact, 
 if $V$ is of type Y (Definition \ref{def:OSO}), 
 then $\varepsilon$ is not uniquely determined 
 because there are two expressions for $V$:
\[
  V \simeq \Kirredrep{O(n)}{\sigma}_+ \simeq \Kirredrep{O(n)}{\sigma}_-, 
\]
 see Lemma \ref{lem:fdeq} (1).  
On the other hand, 
 the $(m+1)$-th component of $\sigma(\lambda)$ does not vanish
 except for the case $i=\lambda=m$
 by Definition \ref{def:sigmailmd}.  
Hence we obtain an isomorphism
 of $O(n+1,1)$-modules:
\[
   \Kirredrep {O(n+1,1)}{\sigma(\lambda)}_{a,b}
   \simeq
   \Kirredrep {O(n+1,1)}{\sigma(\lambda)}_{-a,-b}
\]
for any $a, b \in \{\pm\}$ by Lemma \ref{lem:fdeq} (2).  
\end{remark}

By Definition \ref{def:typeone}, 
 the following lemma is clear.  
\begin{lemma}
\label{lem:FVlmddet}
Suppose that $V$ is of type X 
 or $\lambda \ne \frac n 2$.  
Then there is a natural isomorphism
 of $O(n+1,1)$-modules:
\[
  F(V \otimes \det, \lambda)
\simeq
 F(V, \lambda) \otimes \det.  
\]
\end{lemma}
\begin{lemma}
\label{lem:FVlmdY}
The following two conditions on $(V,\lambda)\in \RInt$
 ({\it{i.e.}}, 
 $V \in \widehat{O(n)}$
 and $\lambda \in {\mathbb{Z}} \setminus S(V)$)
 are equivalent:
\begin{enumerate}
\item[{\rm{(i)}}]
$F(V, \lambda) \otimes \det \simeq F(V, \lambda)$
 as $G$-modules;
\item[{\rm{(ii)}}]
$V$ is of type Y (Definition \ref{def:OSO})
 and $\lambda \ne \frac n 2$.  
\end{enumerate}
In particular,
 for $(V,\lambda)\in \Reducible$, 
 (i) holds 
 if and only if $(V,\lambda)\in \RedJ$
 (Definition \ref{def:Red12}).  
\end{lemma}

\begin{proof}
Any of the conditions (i) or (ii) implies
 that $n$ is even, 
 say,  $n=2m$.  
Let us verify (ii) $\Rightarrow$ (i).  
If we write $V = \Kirredrep{O(n)}{\sigma}_{\varepsilon}$
 for some $\sigma=(\sigma_1, \cdots, \sigma_m) \in \Lambda^+(m)$
 and $\varepsilon\in \{\pm\}$, 
 then $\sigma_m \ne 0$
 because $V$ is of type Y.  
On the other hand,
 the height $i:=i(V,\lambda)$ is not equal to $m$
 because $\lambda \ne m$, 
 hence the $(m+1)$-th component of $\sigma^{(i)}(\lambda)$ equals 
 $\sigma_m (\ne 0)$
 by Definition \ref{def:sigmailmd}.  
Thus there is a natural $G$-isomorphism
 $F(V, \lambda) \otimes \det \simeq F(V, \lambda)$.  
The converse implication is similarly verified.  
\end{proof}
\begin{example}
\label{ex:FVexterior}
Let $(V,\lambda)=(\Exterior^\ell({\mathbb{C}}^n), \ell)$
 for $\ell=0,1,\cdots,n$.  
We set $m=[\frac n 2]$ as usual.  
Then
\[
\text{
 $i(V,\lambda)=\ell,\,\,
 \sigma(\lambda)=0$
 $(\in {\mathbb{Z}}^{m+1}),$ 
 and
 $\widehat{\sigma^{(i)}}=0 \in {\mathbb{Z}}^m$.
}  
\]
Moreover, 
we have an isomorphism of $G$-modules:
\[
   F(V,\lambda)\simeq {\bf{1}}
\quad
  \text{for $0 \le \ell \le n$}.
\]
\end{example}

\subsubsection{Reformulation of Theorems \ref{thm:1807113} and \ref{thm:1808101}}

By using the finite-dimensional representation $F(V,\lambda)$
 of $G=O(n+1,1)$
 (Definition \ref{def:Fshift}), 
 Theorems \ref{thm:1807113} and \ref{thm:1808101}
 may be reformulated in simpler forms, 
respectively, 
 as follows.  
\begin{theorem}
\label{thm:181104}
For $(V,\lambda) \in \Reducible$
 (Definition \ref{def:RIntRed}), 
 we set 
\index{A}{iVlmd@$i(V,\lambda)$, height}
$i:=i(V,\lambda)$, 
 the height of $(V,\lambda)$
 as in \eqref{eqn:indexV}.  
Then there is a natural $G$-isomorphism:
\index{A}{FVlmd@$F(V,\lambda)$}
\[
  P_{r(V,\lambda)}(I_{\delta}(i,i) \otimes F(V,\lambda))
  \simeq 
  I_{\delta}(V, \lambda).  
\]
\end{theorem}

\begin{theorem}
\label{thm:181107}
Suppose $(V,\lambda) \in \Reducible$.  
Retain the notation as in Theorem \ref{thm:181104}.  
\begin{enumerate}
\item[{\rm{(1)}}]
If $(V,\lambda) \in \RedI$
(Definition \ref{def:Red12}), 
 then there is a natural $G$-isomorphism:
\[
  P_{r(V,\lambda)}(I_{\delta}(V, \lambda) \otimes F(V,\lambda))
  \simeq 
  I_{\delta}(i, i).  
\]
\item[{\rm{(2)}}]
If $(V,\lambda) \in \RedJ$, 
 then there is a natural $G$-isomorphism:
\[
  P_{r(V,\lambda)}(I_{\delta}(V, \lambda) \otimes F(V,\lambda))
  \simeq 
  I_{\delta}(i, i) \oplus I_{\delta}(n-i, i).  
\]
\end{enumerate}
\end{theorem}

\subsubsection{Translation of irreducible representations
 $\Pi_{\ell,\delta}$}

We recall from \eqref{eqn:Pild}
 that $\Pi_{\ell,\delta}$ ($0 \le \ell \le n+1$, $\delta \in \{\pm\}$)
 are irreducible admissible smooth representations
 of $G$
 with trivial infinitesimal character $\rho_G$,
 and from \eqref{eqn:PiVlmd}
 that $\Pi_{\delta}(V,\lambda)$ is an irreducible admissible smooth representation of $G$
 with ${\mathfrak{Z}}_G({\mathfrak{g}})$-infinitesimal character
 $r(V,\lambda) \mod W_G$.  
We also recall that $\rho^{(i)} \equiv \rho_G \mod W_G$
for all $0 \le i \le n$.  
In this section,
 we determine the action 
 of translation functor 
 $\psi_{\rho^{(i)}}^{r(V,\lambda)}$ 
 on irreducible representations.  

\begin{theorem}
\label{thm:20180904}
Suppose that $(V,\lambda) \in \Reducible$.  
Let $i:=i(V,\lambda)$ be the height of $(V,\lambda)$, 
 and $F(V,\lambda)$ be the irreducible finite-dimensional representation
 of $G$
 (Definition \ref{def:Fshift}).  
Then there is a natural $G$-isomorphism:
\[
  P_{r(V,\lambda)} (\Pi_{i,\delta} \otimes F(V,\lambda))
  \simeq
  \Pi_{\delta}(V,\lambda).  
\]
\end{theorem}

\begin{proof}
Since the translation functor is a covariant exact functor
 (Proposition \ref{prop:transfunct} (1)), 
 the exact sequence of $G$-modules
\[
  0 \to \Pi_{i,\delta}\to I_{\delta}(i,i) \to \Pi_{i+1,-\delta} \to 0
\]
(Theorem \ref{thm:LNM20} (1)) yields 
 an exact sequence of $G$-modules
\[
  0 \to P_{r(V,\lambda)}(\Pi_{i,\delta} \otimes F(V,\lambda))
    \to I_{\delta}(V,\lambda)
    \to P_{r(V,\lambda)}(\Pi_{i+1,-\delta}\otimes F(V,\lambda))
    \to 0, 
\]
 where we have used Theorem \ref{thm:181104}
 for the middle term.  
Since the first and third terms do not vanish 
 by Proposition \ref{prop:transfunct} (2), 
 we conclude the following isomorphisms of $G$-modules:
\begin{align*}
I_{\delta}(V,\lambda)^{\flat}
&\simeq P_{r(V,\lambda)}(\Pi_{i,\delta} \otimes F(V,\lambda)), 
\\
I_{\delta}(V,\lambda)^{\sharp}
&\simeq P_{r(V,\lambda)}(\Pi_{i+1,-\delta} \otimes F(V,\lambda))
\end{align*}
because $I_{\delta}(V,\lambda)$ has composition series
 of length two
 (Corollary \ref{cor:length2}).  
Hence Theorem \ref{thm:20180904} follows from the definition 
 \eqref{eqn:PiVlmd} of $\Pi_{\delta}(V,\lambda)$.  
\end{proof}

\subsubsection{Proof of Theorems \ref{thm:181104} and \ref{thm:181107}}

In this subsection,
 we explain that Theorem \ref{thm:1807113}
 is equivalent to Theorem \ref{thm:181104};
 Theorem \ref{thm:1808101} is equivalent to Theorem \ref{thm:181107}.  

For this 
 we begin with the following lemma which clarifies
 some combinatorial meaning
 of the height $i(V,\lambda) \in \{0,1,\ldots,n\}$
 and the dominant integral weight $\sigma(\lambda) \in \Lambda^+(m+1)$
 in Definition \ref{def:Fshift}.  
Here we recall $m=[\frac n 2]$.  
\begin{lemma}
\label{lem:transvec}
Suppose $V=\Kirredrep{O(n)}{\sigma}_{\varepsilon}$
 with $\sigma \in \Lambda^+(m)$
 and $\varepsilon \in \{\pm\}$.  
For $0 \le i \le n$ and $\lambda \in {\mathbb{Z}}$, 
 we set
\index{A}{1patauiVlmd@$\tau^{(i)}(V,\lambda)$|textbf}
\begin{equation}
\label{eqn:tauiVlmd}
\tau^{(i)}(V,\lambda) :=r(V,\lambda)-\rho^{(i)} \in (\frac 1 2{\mathbb{Z}})^{m+1}, 
\end{equation}
see \eqref{eqn:IVZG} and Example \ref{ex:rhoi}
 for the notation.  
\begin{enumerate}
\item[{\rm{(1)}}]
Then $\tau^{(i)}(V,\lambda) \in {\mathbb{Z}}^{m+1}$
 is given by 
\[
\begin{cases}
(\sigma_{1}-1,\cdots,\sigma_{i}-1,\sigma_{i+1},\cdots,\sigma_{m}, \lambda-i)
\quad
&\text{for $0 \le i \le m$, }
\\
(\sigma_{1}-1,\cdots,\sigma_{n-i}-1,\sigma_{n-i+1}, \cdots,\sigma_m, \lambda-i)
\quad
&\text{for $m+1 \le i \le n$.  }
\end{cases}
\]
\item[{\rm{(2)}}]
Assume that $\lambda \in {\mathbb{Z}} \setminus S(V)$, 
 and we take $i$ to be the height $i(V, \lambda)$ of $(V,\lambda)$
 as in \eqref{eqn:indexV}.  
Then, 
\[
  {\textsl{$r(V,\lambda)$ and $\rho^{(i)}$ belong to the same Weyl chamber
 for $W_{\mathfrak{g}}$.   }}
\]
\item[{\rm{(3)}}]
Let $\sigma(\lambda)$ be as defined in Definition \ref{def:Fshift}.  
Then we have
\begin{equation}
\label{eqn:tauisigma}
   \tau^{(i)}(V,\lambda)_{\operatorname{dom}}
   =
  \sigma(\lambda).  
\end{equation}
\end{enumerate}
\end{lemma}

\begin{proof}
(1)\enspace
 Clear from the definition \eqref{eqn:IVZG}
 of $r(V,\lambda)$ and $\rho^{(i)}$.  
(2)\enspace
The assertion is verified by inspecting
 the definition \eqref{eqn:indexV} of the height $i(V,\lambda)$.  (3)\enspace
The statement follows from Definition-Lemma \ref{deflem:180694}.  
\end{proof}

Now we determine the action of the translation functor
 $\psi_{\rho^{(i)}}^{r(V,\lambda)}$.  
We recall that the principal series representation
 $I_{\delta}(i,i)$ $(0 \le i \le n)$
 has the trivial ${\mathfrak{Z}}_G({\mathfrak{g}})$-infinitesimal character,
 which is $W_{\mathfrak{g}}$-regular
 but not always $W_G$-regular.  
We apply the translation functor \eqref{eqn:translation} to $I_{\delta}(i,i)$
 for an appropriate choice of $i$.  

\begin{proposition}
\label{prop:transIii}
Let $m=[\frac n 2]$.  
Suppose $G=O(n+1,1)$, 
 $\delta \in \{\pm\}$, 
 $V=\Kirredrep{O(n)}{\sigma}_{\varepsilon}$
 with $\sigma \in \Lambda^+(m)$ and $\varepsilon \in \{\pm\}$, 
 and $\lambda \in {\mathbb{Z}} \setminus (S(V) \cup S_Y(V))$.  
Let $i:=i(V,\lambda) \in \{0,1,\ldots,n\}$ be as in \eqref{eqn:indexV}. 
We define $r(V, \lambda)\in {\mathbb{C}}^{m+1}$ as in \eqref{eqn:IVZG}
 and $\sigma(\lambda)\in \Lambda^+(m+1)$.  
Then we have 
\[
  \psi_{\rho^{(i)}}^{r(V,\lambda)}(I_{\delta}(i,i)) 
= P_{r(V,\lambda)}(I_{\delta}(i,i)  
  \otimes 
  \Kirredrep {O(n+1,1)}{\sigma(\lambda)}_{+,+}).  
\]
\end{proposition}
\begin{proof}
Since $I_{\delta}(i,i)$ has the trivial ${\mathfrak{Z}}_G({\mathfrak{g}})$-infinitesimal character,
 $P_{\rho^{(i)}}(I_{\delta}(i,i))=I_{\delta}(i,i)$
 by \eqref{eqn:rhoiW}.  
Since $r(V,\lambda)=\rho^{(i)} + \tau^{(i)}(V,\lambda)$
 by \eqref{eqn:tauiVlmd}, 
 and since $\sigma(\lambda)=\tau^{(i)}(V,\lambda)_{\operatorname{dom}}$ by
 \eqref{eqn:tauisigma}, 
 the definition of the translation functor shows 
\[
  \psi_{\rho^{(i)}}^{\rho^{(i)} + \tau^{(i)}(V,\lambda)}(I_{\delta}(i,i)) 
= P_{r(V,\lambda)}(I_{\delta}(i,i) \otimes \Kirredrep{O(n+1,1)}{\sigma(\lambda)}_{+,+}).
\]

Thus Proposition \ref{prop:transIii} is proved.  
\end{proof}

It follows from Proposition \ref{prop:transIii}
 and from the definition of $F(V,\lambda)$
(Definition \ref{def:Fshift})
 that Theorem \ref{thm:1807113}
 is equivalent to Theorem \ref{thm:181104}
 and Theorem \ref{thm:1808101} is equivalent to Theorem \ref{thm:181107}.

\subsection{Proof of Proposition \ref{prop:IVsub}}
\label{subsec:pftrans}
In this section we complete the proof of Proposition \ref{prop:IVsub}.  
By Lemma \ref{lem:Ftensor}, 
 the proof reduces to some branching laws
 for the restriction
 of finite-dimensional representations
 of $G=O(n+1,1)$
 with respect to $M A \simeq O(n) \times SO(1,1)$
 and to the study of their tensor product representations, 
 see Proposition \ref{prop:180863}.

\subsubsection{Irreducible summands
 for $O(n+2) \downarrow O(n) \times O(2)$
 and for tensor product representations}
Before working with Proposition \ref{prop:180863}
 in the noncompact setting, 
 we first discuss analogous branching rules 
 for the restriction with respect to a pair
 of compact groups $O(n+2) \supset O(n) \times O(2)$:
\begin{lemma}
[$O(n+2)\downarrow O(n) \times O(2)$]
\label{lem:1810122}
Let $\mu=(\mu_1, \cdots,\mu_{m+1}) \in \Lambda^+(m+1)$, 
 where $m:=[\frac n 2]$ as before.  
For $1\le k \le m+1$, 
 we set
\[
\mu_{(k)}':=(\mu_1, \cdots,\mu_{k-1}, \widehat{\mu_k}, \mu_{k+1}, \cdots,\mu_{m+1})
 \in \Lambda^+(m).  
\]
Then the $O(n+2)$-module $\Kirredrep{O(n+2)}{\mu}_+$
 (see \eqref{eqn:ONCreal}) contains the 
 $(O(n) \times O(2))$-module
\[
   \bigoplus_{k=1}^{m+1} \Kirredrep{O(n)}{\mu_{(k)}'}_+
   \boxtimes \Kirredrep{O(2)}{\mu_k}_+
\]
 when restricted to the subgroup $O(n) \times O(2)$.  
\end{lemma}

\begin{proof}
Take a Cartan subalgebra ${\mathfrak{h}}_{\mathbb{C}}$
 of ${\mathfrak{g l}}(n+2,{\mathbb{C}})$
 such that ${\mathfrak{h}}_{\mathbb{C}} \cap {\mathfrak{o}}(n+2,{\mathbb{C}})$
 is a Cartan subalgebra of ${\mathfrak{o}}(n+2,{\mathbb{C}})$.  
We identify ${\mathfrak{h}}_{\mathbb{C}}^{\ast}$
 with ${\mathbb{C}}^{n+1}$ 
 via the standard basis $\{f_j\}$
 as before,
 and choose a positive system $\Delta^+({\mathfrak{gl }}(n+2,{\mathbb{C}}), {\mathfrak{h}}_{\mathbb{C}})=\{f_i-f_j: 1 \le i < j \le n+2\}$.  
Then 
\[
\widetilde \mu:=(\mu_1, \cdots,\mu_{m+1}, 0^{n+1-m})\in \Lambda^+(n+2)
\]
 is a dominant integral 
 with respect to the positive system. 
Let $v_{\widetilde \mu}$ be a (nonzero) highest weight vector
 of the irreducible representation 
 $(\tau, \Kirredrep{U(n+2)}{\widetilde \mu})$
 of the unitary group $U(n+2)$.  
By definition,
 the $O(n+2)$-module
 $\Kirredrep{O(n+2)}{\mu}_+$, 
 see \eqref{eqn:ONCreal}, 
 is the unique irreducible $O(n+2)$-summand
 of $\Kirredrep{U(n+2)}{\widetilde \mu}$
 containing the highest weight vector $v_{\widetilde \mu}$.  
We now take a closer look at the $U(n+2)$-module
 $\Kirredrep{U(n+2)}{\widetilde \mu}$.  
Fix $1 \le k \le m+1$.  
Iterating the classical branching rule 
 for $U(N) \supset U(N-1) \times U(1)$
 for $N=n+2$, $n+1$,  
 we see that the restriction
 $\Kirredrep{U(n+2)}{\widetilde \mu}|_{U(n) \times U(2)}$
 contains
\[
   \Kirredrep{U(n)}{\widetilde {\mu_{(k)}'}} \boxtimes W
\]
 as an irreducible summand,
 where
\[
  \widetilde{\mu_{(k)}'}:=(\mu_1, \cdots,\widehat{\mu_k}, \cdots,\mu_{m+1}, 0^{n-m}) \in \Lambda^+(n)
\]
and $W$ is an irreducible representation of $U(2)$
 which has a weight $(\mu_k,0)$.  
Since all the weights of an irreducible finite-dimensional representation
 are contained in the convex hull
 of the Weyl group orbit
 through the highest weight, 
 we conclude 
 that $(\mu_k,0)$ is actually the highest weight of the $U(2)$-module $W$.  
Hence the $(U(n) \times U(2))$-module
\[
   \Kirredrep{U(n)}{\widetilde {\mu_{(k)}'}}
   \boxtimes
   \Kirredrep{U(2)}{\mu_k,0}
\]
occurs as an irreducible summand
 of the $U(n+2)$-module 
 $\Kirredrep{U(n+2)}{\widetilde \mu}$.  
We now consider the following diagram
 of subgroups of $U(n+2)$, 
 and investigate the restriction of the $U(n+2)$-module
 $\Kirredrep{U(n+2)}{\widetilde \mu}$.  
\begin{alignat*}{3}
 &U(n+2) &&\supset\,\, &&U(n) \times U(2)
\\
 &\hphantom{MM} \cup && &&\hphantom{MM}\cup
\\
 &O(n+2) &&\supset\,\, &&O(n) \times O(2)
\end{alignat*}
By our choice of the Cartan subalgebra
 ${\mathfrak{h}}_{\mathbb{C}}$, 
 we observe
 that there exists $w_k \in O(n+2)$
 such that $\operatorname{Ad}(w_k){\mathfrak{h}}_{\mathbb{C}}
={\mathfrak{h}}_{\mathbb{C}}$
 and 
\[
  w_k \widetilde \mu
=(\mu_1, \cdots,\widehat{\mu_k}, \cdots,\mu_{m+1}, 0^{n-m},\mu_k,0) \in 
{\mathbb{Z}}^{n+2}, 
\]
where we write $w_k \widetilde \mu$ simply
 for the contragredient action
 of $\operatorname{Ad}(w_k)$ on $\widetilde \mu \in {\mathfrak{h}}_{\mathbb{C}}^{\ast}$
 $(\simeq {\mathbb{C}}^{n+2})$.  
In particular,
 the $O(n+2)$-submodule
 $\Kirredrep{O(n+2)}{\mu}_+$
 of the restriction $\Kirredrep{U(n+2)}{\widetilde \mu}|_{O(n+2)}$
 contains the weight vector
 $v_{w_k \widetilde \mu}:=\tau(w_k) v_{\widetilde \mu}$
 for the weight $w_k \widetilde \mu$.  
Since $w_k \widetilde \mu$ is an extremal weight, 
 the weight vector in $\Kirredrep{U(n+2)}{\widetilde \mu}|_{O(n+2)}$
 is unique up to scalar multiplication.  
Hence $v_{w_k \widetilde \mu}$ is contained also in the submodule
 $\Kirredrep{U(n)}{\widetilde {\mu_{(k)}'}} \boxtimes \Kirredrep{U(2)}{\mu_{k},0}$.  
Thus we conclude that the irreducible $O(n+2)$-module
 $\Kirredrep{O(n+2)}{\mu}_+$ contains
\[
   \Kirredrep{O(n)}{\mu_{(k)}'} \boxtimes \Kirredrep{O(2)}{\mu_k}
\]
 as an $(O(n) \times O(2))$-summand 
 when restricted to the subgroup $O(n) \times O(2)$ of $O(n+2)$.  
\end{proof}

Let $m=[\frac n 2]$
 as before.  
Let $V=\Kirredrep{O(n)}{\sigma}_{\varepsilon}$
 with $\sigma \in \Lambda^+(m)$
 and $\varepsilon \in \{\pm\}$.  
Suppose $\lambda \in {\mathbb{Z}} \setminus S(V)$.  
We recall from Definition-Lemma \ref{deflem:180694}
 that $\sigma(\lambda) \in \Lambda^+(m+1)$.  

\begin{lemma}
\label{lem:180862}
For $0 \le i \le n$, 
 the following $(O(n) \times O(2))$-module
\[
  (\Exterior^i({\mathbb{C}}^n) \boxtimes {\bf{1}})
  \otimes 
  \Kirredrep{O(n+2)}{\sigma(\lambda)}_{\varepsilon}|_{O(n) \times O(2)}
\]
contains
\begin{alignat*}{2}
  &V \boxtimes
  \Kirredrep{O(2)}{|\lambda-i|}_+
\qquad
  &&\text{if $i \le \frac n 2$}
\\
  &(V \otimes \det) \boxtimes \Kirredrep{O(2)}{|\lambda-i|}_+
\qquad
 &&\text{if $i \ge \frac n 2$}
\end{alignat*}
as an irreducible summand.  
\end{lemma}
We note that $V \simeq V \otimes \det$
 as $O(n)$-modules
 if $i=\frac n 2$
 by Lemmas \ref{lem:typeY} and \ref{lem:isig}.  
\begin{proof}
It suffices to prove Lemma \ref{lem:180862}
 for $\varepsilon =+$
 by using a similar argument
 to \eqref{eqn:314}
 for the pair $(O(n+2), O(n) \times O(2))$
 and for $\chi=\det$.  
Then Lemma \ref{lem:180862} is derived from the following two branching laws of compact Lie groups.  
\begin{enumerate}
\item[$\bullet$]
$O(n+2) \downarrow O(n) \times O(2)$:
\par
By Lemma \ref{lem:1810122}, 
 the $O(n+2)$-module $\Kirredrep{O(n+2)}{\sigma(\lambda)}_+$ contains
\[
  \Kirredrep{O(n)}{\widehat{\sigma^{(i)}}}_+
  \boxtimes
  \Kirredrep{O(2)}{|\lambda-i|}
\]
as an irreducible summand
when restricted to the subgroup $O(n) \times O(2)$, 
 see Definition \ref{def:sigmailmd} for the notation $\widehat{\sigma^{(i)}}$. 
\item[$\bullet$]
Tensor product for $O(n)$:
\par
The tensor product representation
\[
  \Exterior^i({\mathbb{C}}^n) 
  \otimes 
  \Kirredrep{O(n)}{\widehat{\sigma^{(i)}}}_+
\]
contains 
\begin{alignat*}{2}
V&\simeq  \Kirredrep {O(n)}{\sigma}_{+}
&&\text{if $i \le \frac n 2$}, 
\\
V\otimes \det &\simeq  \Kirredrep {O(n)}{\sigma}_{-}
\qquad
&&\text{if $i \ge \frac n 2$}
\end{alignat*}
 as an irreducible component.  
\end{enumerate}
\end{proof}

\subsubsection{Irreducible summand for the restriction
 $G \downarrow M A$
 and for tensor product representations}

We recall that the Levi subgroup $M A$
 of the parabolic subgroup $P$
 in $G=O(n+1,1)$ is expressed as
\[
   M A \simeq  O(n) \times S O(1,1)
       \simeq O(n) \times {\mathbb{Z}} / 2 {\mathbb{Z}} \times {\mathbb{R}}.
\]  
The goal of the subsection is to prove
 the following proposition.  
\begin{proposition}
[tensor product and the restriction $O(n+1,1) \downarrow M A$]
\label{prop:180863}
~~~\newline
Suppose that $(V,\lambda) \in \Reducible$, 
 {\it{i.e.}}, 
 $V \in \widehat{O(n)}$
 and $\lambda \in {\mathbb{Z}} \setminus (S(V) \cup S_Y(V))$.  
Let $i=i(V,\lambda)$ be the height of $(V,\lambda)$
 (see \eqref{eqn:indexV}), 
 and $F(V, \lambda)$ be the irreducible $O(n+1,1)$-module
 as in Definition \ref{def:Fshift}.  
Then the $M A$-module
\index{A}{FVlmd@$F(V,\lambda)$}
\begin{equation}
\label{eqn:itensorF}
  (\Exterior^i({\mathbb{C}}^n) \boxtimes \delta \boxtimes {\mathbb{C}}_i)
  \otimes
  F(V, \lambda)|_{M A}
\end{equation}
contains
\[
V \boxtimes \delta \boxtimes {\mathbb{C}}_{\lambda}
\]
 as an irreducible component.  
\end{proposition}
In what follows, 
 we use a mixture of notations
 in describing irreducible finite-dimensional representations
 (see Sections \ref{subsec:repON} and \ref{subsec:fdimrep}).  
To be precise, 
 we shall use:
\begin{enumerate}
\item[$\bullet$]
$\Lambda^+(O(n+2))$ $(\subset {\mathbb{Z}}^{n+2})$, 
 see \eqref{eqn:CWOn}, 
 to denote irreducible holomorphic finite-dimensional representations
 of the {\it{complex}} Lie group $O(n+1,{\mathbb{C}})$
 as in Section \ref{subsec:repON};
\item[$\bullet$]
$\Lambda^+(m+1)$ $(\subset {\mathbb{Z}}^{m+1})$
 and signatures to denote irreducible
 finite-dimensional representations
 of the {\it{real}} groups $O(n+2)$ and $O(n+1,1)$
 where $m:=[\frac n 2]$, 
 as in Section \ref{subsec:fdimrep}.

See \eqref{eqn:ONCreal} 
 for the relationship 
 among these representations.  
\end{enumerate}

\begin{proof}
[Proof of Proposition \ref{prop:180863}]
We write $V=\Kirredrep {O(n)}{\sigma}_{\varepsilon}$ as before
 where $\sigma \in \Lambda^+(m)$, 
 $\varepsilon \in \{\pm\}$, 
 and $m=[\frac n 2]$.  
By Weyl's unitary trick for the disconnected group $O(n+1,1)$, 
 see \eqref{eqn:ONCreal}, 
 the restrictions of the holomorphic representation
 $\Kirredrep {O(n+2,{\mathbb{C}})}{\sigma(\lambda), 0^{n+1-m}}$
 to the subgroups $O(n+2)$ and $O(n+1,1)$
 are given respectively by
\begin{align*}
&\Kirredrep {O(n+2)}{\sigma(\lambda)}_{+}, 
\\
&\Kirredrep {O(n+1,1)}{\sigma(\lambda)}_{+,+}.  
\end{align*}
Then Lemma \ref{lem:180862} implies
 that the holomorphic $(O(n,{\mathbb{C}}) \times O(2,{\mathbb{C}}))$-representation
\[
 (\Exterior^i({\mathbb{C}}^n) \boxtimes {\bf{1}})
 \otimes
 \Kirredrep{O(n+2,{\mathbb{C}})}{\sigma(\lambda), 0^{n+1-m}}
 |_{O(n,{\mathbb{C}}) \times O(2,{\mathbb{C}})}
\]
contains
\begin{alignat*}{2}
   &\Kirredrep{O(n,{\mathbb{C}})}{\sigma,0^{n-m}}
   \boxtimes
   \Kirredrep{O(2,{\mathbb{C}})}{|\lambda-i|,0}
\quad
   &&\text{if $i \le \frac n 2$, }
\\
   &(\Kirredrep{O(n,{\mathbb{C}})}{\sigma,0^{n-m}} \otimes \det)
   \boxtimes
   \Kirredrep{O(2,{\mathbb{C}})}{|\lambda-i|,0}
\quad
   &&\text{if $i \ge \frac n 2$}
\end{alignat*}
as an irreducible summand.  
Because the restriction of the first factor
 to compact real from $O(n)$ is isomorphic
 to $\Kirredrep{O(n)}{\sigma}_+$ or $\Kirredrep{O(n)}{\sigma}_-$
 according to whether $i \le \frac n 2$ or $i \ge \frac n 2$.  
Taking the restriction to another real form 
 $O(n) \times O(1,1)$ of $O(n,{\mathbb{C}}) \times O(2,{\mathbb{C}})$, 
 we set that the $(O(n) \times O(1,1))$-module
\[
   (\Exterior^i({\mathbb{C}}^n) \boxtimes {\bf{1}})
   \otimes
   \Kirredrep {O(n+1,1)}{\sigma(\lambda)}_{+,+}
   |_{O(n) \times O(1,1)}
\]
contains
\begin{alignat*}{2}
   &\Kirredrep{O(n)}{\sigma}_+
   \boxtimes
   \Kirredrep{O(2,{\mathbb{C}})}{|\lambda-i|,0}|_{O(1,1)}
\quad
   &&\text{if $i \le \frac n 2$}, 
\\
   &\Kirredrep{O(n)}{\sigma}_- 
   \boxtimes
   \Kirredrep{O(2,{\mathbb{C}})}{|\lambda-i|,0}|_{O(1,1)}
\quad
   &&\text{if $i \ge \frac n 2$}
\end{alignat*}
as an irreducible summand.

Since $V = \Kirredrep{O(n)}{\sigma}_{\varepsilon}$, 
 the definition of $F(V,\lambda)$
 (Definition \ref{def:Fshift})
 implies that the $M A$-module
\[
  (\Exterior^i({\mathbb{C}}^n) \boxtimes {\bf{1}}) \otimes F(V,\lambda)|_{M A}
\]
 contains 
\[
   V \boxtimes 
  (\Kirredrep {O(2,{\mathbb{C}})}{|\lambda-i|,0}|_{S O(1,1)} 
   \otimes 
   \chi_{\varepsilon,(-1)^{\lambda-i}\varepsilon}|_{S O(1,1)})
\]
as an $M A$-module.  
Here we have used that $M A \simeq O(n) \times S O (1,1)$
 and that $\chi_{a,b}|_{S O(1,1)} \simeq \chi_{-a,-b}|_{S O(1,1)}$.  
Hence Proposition \ref{prop:180863} is derived from the following lemma
 on the restriction $O(2,{\mathbb{C}}) \downarrow S O(1,1)$.  
\end{proof}

Let ${\mathbb{C}}_k$ denote the holomorphic character 
 of $S O(2,{\mathbb{C}})$ on ${\mathbb{C}} e^{i k \theta}$.  

\begin{lemma}
[$O(2,{\mathbb{C}}) \downarrow S O(1,1)$]
\label{lem:SOtworest}
\[
  \Kirredrep{O(2,{\mathbb{C}})}{k,0}|_{{\mathbb{Z}} / 2 {\mathbb{Z}} \times {\mathbb{R}}}
  \simeq
\begin{cases}
(-1)^k \boxtimes ({\mathbb{C}}_k \oplus {\mathbb{C}}_{-k}) 
\quad
&\text{for $k \in {\mathbb{N}}_+$}, 
\\
{\bf{1}} \boxtimes {\bf{1}}
\quad
&\text{for $k =0$, }
\end{cases}
\]
where we identify $S O(1,1) \simeq \{\pm I_2\} \times S O_0(1,1)$
 with ${\mathbb{Z}} / 2 {\mathbb{Z}} \times {\mathbb{R}}$.  
In particular, 
the $SO(1,1)$-module
\begin{align*}
&\Kirredrep {O(2,{\mathbb{C}})}{|\lambda-i|,0}|_{S O(1,1)} 
   \otimes 
   \chi_{\varepsilon,(-1)^{\lambda-i}\varepsilon}|_{S O(1,1)}
\\
\simeq
&\Kirredrep {O(2,{\mathbb{C}})}{|\lambda-i|,0}|_{S O(1,1)} 
   \otimes 
   \chi_{-\varepsilon,(-1)^{\lambda-i-1}\varepsilon}|_{S O(1,1)}
\end{align*}
contains
\[
  {\bf{1}} \boxtimes {\mathbb{C}}_{\lambda -i}
\]
 as an irreducible summand.  
\end{lemma}

\begin{proof}
For $k \in {\mathbb{N}}_+$, 
 the holomorphic representation $\Kirredrep{O(2,{\mathbb{C}})}{k,0}$
 is a two-dimensional representation 
 of $O(2,{\mathbb{C}})$, 
 which is isomorphic to $\operatorname{Ind}_{S O(2,{\mathbb{C}})}^{O(2,{\mathbb{C}})}({\mathbb{C}} e^{i k \theta})$.  
Its restriction to the connected subgroup $S O(2,{\mathbb{C}})$
 decomposes into a sum of two characters
 of $S O(2,{\mathbb{C}})$:  
\[
   \Kirredrep{O(2,{\mathbb{C}})}{k,0}|_{S O(2,{\mathbb{C}})}
  \simeq
  {\mathbb{C}} e^{i k \theta} \oplus {\mathbb{C}} e^{-i k \theta}, 
\]
on which the central element $-I_2$ acts
 as the scalar multiplication of $(-1)^k=(-1)^{-k}$.  
Since $S O(1,1)$ is generated by the central element
 $-I_2$ 
 and the identity component $S O_0(1,1)$, 
 Lemma \ref{lem:SOtworest} follows.  
\end{proof}

\subsubsection{Proof of Proposition \ref{prop:IVsub}}
\label{subsec:1653}
\begin{proof}
[Proof of Proposition \ref{prop:IVsub}]
Let $F(V,\lambda)$ be the finite-dimensional representation of $G$
 as in Definition \ref{def:Fshift}.  
Filter $F(V,\lambda)$ as in Lemma \ref{lem:Ftensor}.  
We may assume in addition 
 that each $F^{(j)} = F(V,\lambda)_j/F(V,\lambda)_{j-1}$ is irreducible
 as an $M A$-module.  
Then by Proposition \ref{prop:180863}, 
 $I_{\delta}(V,\lambda)$ occurs as a subquotient
 of the $G$-module
 $P_{r(V,\lambda)}(I_{\delta}(i,i) \otimes F(V,\lambda))$.  
Hence the second assertion 
 of Proposition \ref{prop:IVsub} is shown.  
By Proposition \ref{prop:transIii}, 
 the first assertion follows.  
\end{proof}

\subsection{Proof of Theorem \ref{thm:1807113}}
\label{subsec:pftran1}

In this section
 we complete the proof of Theorem \ref{thm:1807113}
 and also its reformulation Theorem \ref{thm:181104}.  
By Proposition \ref{prop:IVsub}, 
 it suffices to show Proposition \ref{prop:transrest}
 is an isomorphism
 in the level of $\overline G$-modules
 instead of the isomorphism
 in Theorem \ref{thm:1807113} as $G$-modules.

We divide the argument 
according to the decomposition 
\[
 \Reducible=\RedI \amalg \RedJ, 
\]
where we recall from Definition \ref{def:Red12}:
\begin{enumerate}
\item[$\bullet$]
$(V,\lambda) \in \RedI$, 
 if $V$ is of type X or $\lambda = \frac n 2$;
\item[$\bullet$]
$(V,\lambda) \in \RedJ$, 
 if $V$ is of type Y and $\lambda =\frac n 2$.  
\end{enumerate}

As we shall show in the proof of Proposition \ref{prop:transrest} below, 
 the following assertion holds with the notation therein.  
\begin{proposition}
\label{prop:transrest2}
There  is a natural isomorphism, 
 as $\overline G$-modules
\begin{equation*}
  \psi_{\rho^{(i)}}^{r(V,\lambda)}(I_{\delta}(i,i))|_{\overline G}
\simeq
\begin{cases}
  \overline \psi_{\rho^{(i)}}^{r(V,\lambda)}
 (I_{\delta}(i,i)|_{\overline G})
\quad
&\text{if $(V,\lambda)\in \RedI$}, 
\\
\bigoplus_{\xi=\pm}
\overline \psi_{\rho^{(i)}}^{r(V^{(\xi)},\lambda)}
 (I_{\delta}(i,i)|_{\overline G})
\quad
&\text{if $(V,\lambda)\in \RedJ$}.  
\end{cases}
\end{equation*}
\end{proposition}

\subsubsection{Case: $(V,\lambda) \in \RedI$}
\index{A}{RzqSeducibleI@$\RedI$}
In this subsection,
 we discuss the case 
 where $V$ is of type X
 or $\lambda=\frac n 2$.  
\begin{proof}
[Proof of Proposition \ref{prop:transrest} for $(V,\lambda) \in \RedI$.]
If $n$ is odd, 
 then Proposition \ref{prop:transrest}
 follows from Lemma \ref{lem:Xtrans} (1).

Hereafter we assume $n$ is even, 
 say $n=2m$.  
We claim 
 that Proposition \ref{prop:transrest} follows from
 Lemma \ref{lem:Xtrans} (2)
 if $V$ is of type X
 (Definition \ref{def:OSO})
 or $\lambda=m$.  
To see this, 
 it is enough 
 to verify 
 that all of $\rho^{(i)}$, 
 $r(V,\lambda)$, 
 and $\tau^{(i)}(V,\lambda)=r(V,\lambda)-\rho^{(i)}$, 
 see \eqref{eqn:tauiVlmd}, 
 contain 0
 in their entries.  
This is automatically true
 for $\rho^{(i)}$ as $\rho^{(i)} \in W_G \rho^G$
 (Example \ref{ex:rhoi} (3))
 and $n$ is even.  
For $r(V,\lambda)$, 
 one sees from \eqref{eqn:IVZG} 
 that the $m$-th component vanishes
 if $V$ is of type X
 and the $(m+1)$-th component vanishes
 if $\lambda=m$.  
For $\tau^{(i)}(V,\lambda)$, 
 one see from the formula 
 of $\tau^{(i)}(V,\lambda)$
in Lemma \ref{lem:transvec}
 that an analogous assertion holds
 because $\lambda=m$ $(=\frac n 2)$ implies
 that the height $i(V,\lambda)$ equals $m$
 by Definition \ref{def:iVlmd}.  
Hence Proposition \ref{prop:transrest}
 for $(V,\lambda) \in \RedI$ is shown.  
\end{proof}

\subsubsection{Case: $(V,\lambda) \in \RedJ$}
\index{A}{RzqSeducibleII@$\RedJ$}

In this subsection, 
 we discuss the case 
 where $V$ is of type Y and $\lambda \ne \frac n2$.  
In this case, 
 $n$ is even ($=2m$), 
 $i:=i(V,\lambda)=m$,  
 and the restriction of $V$ to $SO(n)$
 is a sum of two irreducible representations of $SO(n)$:
\[
  V=V^{(+)} \oplus V^{(-)}, 
\]
as in \eqref{eqn:Vpm}.  
We extend the definition \eqref{eqn:IVZG}
 of $r(V,\lambda)$ 
 to irreducible representations $V^{(\pm)}$ of $SO(n)$
 with $n=2m$ by 
\[
  r(V^{(\pm)},\lambda):=(\sigma_1+m-1,\cdots, \sigma_{m-1}+1, \pm \sigma_m, \lambda-m) \in {\mathbb{Z}}^{m+1}.  
\]
Then $r(V^{(\pm)}, \lambda)$ viewed as an element of ${\mathfrak{h}}_{\mathbb{C}}^{\ast}/W_{{\mathfrak{g}}}$ 
 is the ${\mathfrak{Z}}({\mathfrak{g}})$-infinitesimal character
 of the principal series representation
 $\overline I_{\delta}(V^{(\pm)},\lambda)$
 of $\overline G=SO(n+1,1)$.

As in \eqref{eqn:tauiVlmd}, 
 we set
\[
   \tau^{(i)}(V^{(\pm)},\lambda):=r(V^{(\pm)},\lambda)-\rho^{(i)}.  
\]
Inspecting the definition \eqref{eqn:indexV}
 of the height $i:=i(V,\lambda)$, 
 we see that both $r(V^{(\pm)},\lambda)$ and $\rho^{(i)}$
 belong to the same Weyl chamber
 with respect to the Weyl group $W_{\mathfrak{g}}$
 ({\it{not}} $W_G$)
 as in Lemma \ref{lem:transvec}.

By Lemma \ref{lem:FVlmdY}, 
 the irreducible finite-dimensional $G$-module
 $F(V, \lambda)$ decomposes into a direct sum
 of two irreducible $\overline G$-modules, 
 which we may write as
\[
  F(V, \lambda)|_{\overline G}=\overline F(V^{(+)}, \lambda)
  \oplus \overline F(V^{(-)}, \lambda).  
\]
To be precise, 
 we set $\sigma^{(+)} (\lambda):=\sigma(\lambda)$
 (Definition \ref{def:Fshift}), 
 and define $\sigma^{(-)}(\lambda)$
 by replacing the $(m+1)$-th component $\sigma_m$
 with $-\sigma_m$.  
For instance,
 if $\lambda<m$, 
 then the height $i=i(V,\lambda)$ is smaller than $m$
and
\begin{align*}
\sigma^{(+)}(\lambda) =&(\sigma_1-1,\cdots,\sigma_i-1,i-\lambda, \sigma_{i+1},\cdots, \sigma_{m-1}, \sigma_{m}), 
\\
\sigma^{(-)}(\lambda) =&(\sigma_1-1,\cdots,\sigma_i-1,i-\lambda, \sigma_{i+1},\cdots, \sigma_{m-1}, -\sigma_{m}).  
\end{align*}
Then $\overline F(V^{(\pm)}, \lambda)$ are
 the irreducible $\overline G$-modules
 such that 
\[
\overline F(V^{(\pm)}, \lambda) \otimes \chi_{+,(-1)^{\lambda-i}}|_{SO(n+1,1)}
\]
 extends to irreducible holomorphic finite-dimensional representations
 of the connected complex Lie group $SO(n+2,{\mathbb{C}})$
 with highest weights $\sigma^{(\pm)}(\lambda)$.

\begin{proof}
[Proof of Proposition \ref{prop:transrest} for $(V,\lambda) \in \RedJ$]
By the definition \eqref{eqn:translation}
 of the translation functor
 and by Lemma \ref{lem:primary}, 
 there is a natural $\overline G$-isomorphism:
\begin{multline}
\label{eqn:Ytrans}
   \psi_{\rho^{(i)}}^{r(V,\lambda)}(I_{\delta}(i,i))|_{\overline G} 
\\
  \simeq
  (\overline P_{r(V^{(+)},\lambda)}
   + 
   \overline P_{r(V^{(-)},\lambda)})
  (I_{\delta}(i,i)|_{\overline G}
  \otimes 
  (\overline F(V^{(+)},\lambda) \oplus \overline F(V^{(-)},\lambda))).  
\end{multline}
We claim for $\xi, \eta \in \{\pm\}$: 
\begin{alignat}{2}
\label{eqn:paraterm}
   \overline P_{r(V^{(\xi)},\lambda)}
  (I_{\delta}(i,i)|_{\overline G}
  \otimes 
 \overline F(V^{(\eta)},\lambda))
  =&
  \overline \psi_{\rho^{(i)}}^{r(V^{(\xi)}, \lambda)}(I_{\delta}(i,i)|_{\overline G})
\quad
&&\text{if $\xi \eta =+$}, 
\\
\label{eqn:crosszero}
   \overline P_{r(V^{(\xi)},\lambda)}
  (I_{\delta}(i,i)|_{\overline G}
  \otimes 
  \overline F(V^{(\eta)},\lambda))=&0
\quad
&&\text{if $\xi \eta =-$}.  
\end{alignat}
The first claim \eqref{eqn:paraterm}
 holds by definition \eqref{eqn:transSO}.  
To see the vanishing \eqref{eqn:crosszero}
 of the cross terms in \eqref{eqn:crosszero},
 suppose that 
\[
   \rho^{(i)} + \gamma = w (\rho^{(i)}+ \tau^{(i)}(V^{(\xi)}, \lambda))
\]
for some weight $\gamma$ in $F(V^{(\eta)}, \lambda)$
 and for some $w \in W_{\mathfrak{g}}$.  
Then we have 
\[
   ||\gamma|| \le ||\tau^{(i)}(V^{(\eta)}, \lambda)||=||\tau^{(i)}(V^{(\xi)}, \lambda)||.  
\]
Hence we can apply Lemma \ref{lem:V7218}
 and conclude
\[
  \gamma = \tau^{(i)}(V^{(\xi)}, \lambda).  
\]
By the vanishing \eqref{eqn:crosszero}
 of the cross terms in \eqref{eqn:Ytrans}, 
 we obtain the following $\overline G$-isomorphisms:
\begin{align*}
  \psi_{\rho^{(i)}}^{r(V,\lambda)}(I_{\delta}(i,i))|_{\overline G} 
  \simeq&
  \bigoplus_{\xi \in \{\pm\}}
  \overline \psi_{\rho^{(i)}}^{r(V^{(\xi)},\lambda)}(I_{\delta}(i,i)|_{\overline G}) \\
  \simeq& 
  \bigoplus_{\xi \in \{\pm\}} \overline I_{(-1)^{\lambda-i} \delta}(V^{(\xi)}, \lambda), 
\end{align*}
which is isomorphic to the restriction
 of the principal series representation
 $I_{(-1)^{\lambda-i} \delta}(V, \lambda)$ of $G$
 to the subgroup $\overline G$
 by \eqref{eqn:IVSOpm}.  
\end{proof}

\subsection{Proof of Theorem \ref{thm:1808101}}
\label{subsec:pf168}
In this section,
 we show Theorem \ref{thm:1808101}, 
 or its reformulation,
 Theorem \ref{thm:181107}.  
The proof is similar 
 to that of Theorem \ref{thm:1807113}, 
 hence we give only a sketch
 of the proof with focus on necessary changes.  
A part of the proof is carried out
 separately
 according to the decomposition 
\[
\Reducible = \RedI \amalg \RedJ
\quad
\text{(Definition \ref{def:Red12})}.  
\]

The following lemma is a counterpart 
 of Proposition \ref{prop:180863}.  
\begin{lemma}
[tensor product and $G \downarrow MA$]
\label{lem:1808101}
Suppose $(V, \lambda) \in \Reducible$.  
Let $i:=i(V,\lambda)$ be the height 
 of $(V,\lambda)$, 
 see \eqref{eqn:indexV}, 
 and $F(V,\lambda)$ be the irreducible
 finite-dimensional representation 
 of $G=O(n+1,1)$
 as in Definition \ref{def:Fshift}.  
Then the $M A$-module
\[
  (V \boxtimes \delta \boxtimes {\mathbb{C}}_{\lambda})
 \otimes 
 F(V, \lambda)|_{M A}
\]
contains
\begin{alignat*}{2}
& \Exterior^i({\mathbb{C}}^n) \boxtimes \delta \boxtimes {\mathbb{C}}_i
&&
\text{if $(V,\lambda) \in \RedI$, }
\\
& (\Exterior^i({\mathbb{C}}^n) \boxtimes \delta \boxtimes {\mathbb{C}}_i)
 \oplus 
(\Exterior^{n-i}({\mathbb{C}}^n) \boxtimes \delta \boxtimes {\mathbb{C}}_i)
\qquad
&&
\text{if $(V,\lambda) \in \RedJ$, }
\end{alignat*}
 as an irreducible component.  
\end{lemma}
\begin{proof}
The proof is similar to that of Proposition \ref{prop:180863}
 except that there is a $G$-isomorphism
 $F(V, \lambda) \otimes \det \simeq F(V, \lambda)$
 by Lemma \ref{lem:FVlmdY}
 if $(V,\lambda) \in \RedJ$.  
In this case, 
 the height $i=i(V,\lambda)$ is not equal to $\frac n 2$
 by Lemma \ref{lem:isig} (3).  
Thus both the $O(n)$-modules 
$
 \Exterior^i({\mathbb{C}}^n)
$
 and 
$
   \Exterior^{n-i}({\mathbb{C}}^n) 
   \simeq 
   \Exterior^{i}({\mathbb{C}}^n) \otimes \det
$
 occur simultaneously 
 in $V \otimes F(V,\lambda)|_{O(n)}$.  
\end{proof}

Theorem \ref{thm:181107}, 
 or equivalently, 
 Theorem \ref{thm:1808101}
 is deduced from the following two propositions.  

\begin{proposition}
\label{prop:Iisub}
 Suppose $(V,\lambda) \in \Reducible$.  
(Definition \ref{def:RIntRed}), 
equivalently,
 $V \in \widehat{O(n)}$
 and $\lambda \in {\mathbb{Z}}\setminus(S(V) \cup S_Y(V))$.  
Then the $G$-module
 $P_{r(V,\lambda)}(I_{\delta}(V,\lambda)
 \otimes F(V,\lambda))$ contains
\begin{alignat*}{2}
& I_{\delta}(i,i)
\quad
&&\text{for $(V,\lambda) \in \RedI$}, 
\\
& \text{$I_{\delta}(i,i)$ and $I_{\delta}(n-i,i)$}
\quad
&&\text{for $(V,\lambda) \in \RedJ$}, 
\end{alignat*}
as subquotients.  
\end{proposition}

\begin{proof}
As in the proof of Proposition \ref{prop:IVsub}
 in Section \ref{subsec:1653}, 
 Proposition \ref{prop:Iisub} follows readily from 
Lemma \ref{lem:Ftensor}
 by using Lemma \ref{lem:1808101}.  
\end{proof}

\begin{proposition}
\label{prop:PIFrest}
Suppose $(V,\lambda) \in \Reducible$, namely,
 $V \in \widehat{O(n)}$
 and $\lambda \in {\mathbb{Z}} \setminus (S(V) \cup S_Y(V))$.  
Then there is a natural isomorphism
 of $\overline G$-modules:
\begin{equation*}
P_{r(V,\lambda)}(I_{\delta}(V,\lambda)
 \otimes F(V,\lambda))|_{\overline G}
 \simeq
\begin{cases}
I_{\delta}(i,i)|_{\overline G}
\quad
&
\text{for $(V,\lambda) \in \RedI$}, 
\\
I_{\delta}(i,i)|_{\overline G} \oplus I_{\delta}(n-i,i)|_{\overline G}
\quad
&
\text{for $(V,\lambda) \in \RedJ$}.  
\end{cases}
\end{equation*}
\end{proposition}

\begin{proof}
The proof is similar to that of Proposition \ref{prop:transrest}, 
 again by showing the vanishing of the cross terms
 as in \eqref{eqn:crosszero}.  
\end{proof}


\printindex{A}{List of Symbols}
\printindex{B}{Index}

\end{document}